\documentclass[10pt,a4paper]{article}
\usepackage{amsfonts}
\usepackage{amssymb}
\usepackage{amsmath}
\usepackage{textcomp}
\usepackage{bbm}
\usepackage{stmaryrd}
\usepackage{bm}
\usepackage[unicode]{hyperref}
\usepackage[
top    = 3cm,
bottom = 3cm,
left   = 3cm,
right  = 4cm]{geometry}

\usepackage{tikz-cd}
\usepackage{graphicx}
\usepackage{theorem}
\newcommand{\df} {\,{\rm d}}

\def \MMM{\mathcal M}
\def \LL{\mathcal L}
\def \B{\mathcal B}
\def \F{\mathcal F}
\def \HHH{\mathcal H}

\def \J{\mathcal J}
\def \E{\mathcal E}
\def \NN{\mathbbm N}
\def \ZZ{\mathbbm Z}
\def \ZZ{\mathbbm Z}
\def \QQ{\mathbbm Q}
\def \RR{\mathbbm R}

\def \TT{\mathcal T}
\def \DD{\mathcal D}
\def \LL{\mathcal L}

\def \Y{\mathcal Y}
\def \X{\mathcal X}

\def \A{\mathcal A}

\def \EE{\mathbbm E}

{\theorembodyfont{\rmfamily} \newtheorem{ori}{Definition}[section]}
{\theorembodyfont{\rmfamily} \newtheorem{rem}{Remark}[section]}
{\theorembodyfont{\rmfamily} \newtheorem{example}{Example}[section]}
\newtheorem{theorema}{Theorem}[section]
\newtheorem{prop}{Proposition}[section]
\newtheorem{lemma}{Lemma}[section]
\newtheorem{cor}{Corollary}[section]
\newcommand{\ra}{\rightarrow}

\newcommand{\lra}{\longrightarrow}
\newcommand{\Lra}{\Longrightarrow}
\newcommand{\Llra}{\Longleftrightarrow}

\newcommand{\D}{\Delta}
\newcommand{\g}{\gamma}

\newcommand{\subs}{\subseteq}

\newcommand{\Lip}{{\rm{Lip}}}

\newcommand{\dv}{{\rm{div}}}

\newcommand{\Ima}{{\rm{Im}}}

\newcommand{\wt}{\widetilde}
\newcommand{\ee}{\varepsilon}

\newcommand{\phih}{\varphi}
\newcommand{\mx}{\vee}
\newcommand{\mn}{\wedge}
\newcommand{\x}{\times}

\newcommand{\W}{\Omega}
\newcommand{\s}{\sigma}
\newcommand{\ls}{\langle}
\newcommand{\rs}{\rangle}
\newcommand{\lls}{\llangle}
\newcommand{\rrs}{\rrangle}

\newcommand{\sm}{\setminus}
\newcommand{\bs}{\boldsymbol}

\newcommand{\bbar}{\overline}

\newcommand{\fr}{\frac{1}}
\newcommand{\PP}{\mathbbm{P}}

\newcommand{\T}{\mathbbm{T}}
\newcommand{\MM}{\mathbbm{M}}

\newcommand{\rN}{\|_{\X^1}}
\newcommand{\1}{\mathbbm{1}}

\newcommand{\w}{\omega}

\newcommand{\y}{\upsilon}
\newcommand{\pd}{\partial}

\makeatletter
\newcommand{\esssup}{\mathop{\operator@font ess\,sup}}
\makeatother
\makeatletter
\newcommand{\argmin}{\mathop{\operator@font argmin}}
\makeatother
\makeatletter
\newcommand{\Lim}{\mathop{\operator@font Lim}}
\makeatother
\makeatletter
\DeclareFontFamily{OMX}{MnSymbolE}{}
\DeclareSymbolFont{MnLargeSymbols}{OMX}{MnSymbolE}{m}{n}
\SetSymbolFont{MnLargeSymbols}{bold}{OMX}{MnSymbolE}{b}{n}
\DeclareFontShape{OMX}{MnSymbolE}{m}{n}{
	<-6>  MnSymbolE5
	<6-7>  MnSymbolE6
	<7-8>  MnSymbolE7
	<8-9>  MnSymbolE8
	<9-10> MnSymbolE9
	<10-12> MnSymbolE10
	<12->   MnSymbolE12
}{}
\DeclareFontShape{OMX}{MnSymbolE}{b}{n}{
	<-6>  MnSymbolE-Bold5
	<6-7>  MnSymbolE-Bold6
	<7-8>  MnSymbolE-Bold7
	<8-9>  MnSymbolE-Bold8
	<9-10> MnSymbolE-Bold9
	<10-12> MnSymbolE-Bold10
	<12->   MnSymbolE-Bold12
}{}

\let\llangle\@undefined
\let\rrangle\@undefined
\DeclareMathDelimiter{\llangle}{\mathopen}%
{MnLargeSymbols}{'164}{MnLargeSymbols}{'164}
\DeclareMathDelimiter{\rrangle}{\mathclose}%
{MnLargeSymbols}{'171}{MnLargeSymbols}{'171}
\makeatother

\begin{document}
		\title{Generalized Young measures and the hydrodynamic limit of condensing zero range processes} 
	\author{Michail Loulakis\thanks{
			School of Applied Mathematical and Physical Sciences, 
			National Technical University of Athens, 
			Polytechnioulopis Zographou, Athens 15780, Greece;
			and Institute of Applied and Computational Mathematics, 
			Foundation for Research and Technology - Hellas, Greece.
			Email: {\tt loulakis@math.ntua.gr}.
		} \ and Marios G.~Stamatakis\thanks{Institute of Applied and Computational Mathematics, 
			Foundation for Research and Technology - Hellas, Greece.
			Email: {\tt ma.ge.stams@gmail.com, mgstam@iacm.forth.gr}.
	}}
		
		\maketitle
		\begin{abstract}
	Condensing zero range processes (ZRPs) are stochastic interacting particle systems that exhibit phase separation with the emergence of a condensate. Standard approaches for deriving a hydrodynamic limit of the density fail in these models, and an effective macroscopic description has not been rigorously established, yet. In this article we prove that the limiting triple $(\pi,W,\s)$ of the empirical density, the empirical current, and the empirical jump rate of the ZRP satisfies the 
	continuity equation $\pd_t\pi=-\dv W$ in the sense of distributions. Here $(\pi_t)_{t\geq 0}$ is a $w^*$-continuous curve of finite non-negative measures on the torus $\T^d$, $\s_t\in H^1(\T^d)$ and $W_t=-\nabla\s_t$ is a vector-valued measure that is absolutely continuous with respect to the Lebesgue measure, for all almost all $t\ge 0$. In order to obtain a closed equation we propose a generalization of Young measures and we prove that for symmetric ZRPs on the torus, the hydrodynamic limit of  the density is a generalized Young-measure-valued weak solution $\bs\pi=(\bs\pi_t)_{t\geq 0}$ to a saturated filtration equation $\pd_t\bs\pi=\Delta\Phi(\bs\pi)$. Furthermore we prove a one-sided two-blocks estimate and we give an equivalent criterion for its validity. Assuming the validity of the two-blocks estimate one obtains the equation $\pd_t\pi=\Delta\Phi(\pi^{ac})$ for the empirical density, where $\pi=\pi^{ac}+\pi^\perp$ is the Radon-Nikodym decomposition.
	\bigskip
	
	\bigskip 
	
	
	\bigskip
\end{abstract}
	\tableofcontents
		
		\section{Introduction} 
		\indent Symmetric Zero Range Processes (ZRPs) are interacting particle systems on a lattice such that each particle jumps at an exponential rate $\mathfrak{g}(k)$ that depends only on the number $k$ of particles that occupy the same site of the lattice, through some function $\mathfrak{g}:\ZZ_+\lra\RR_+$ called the local jump rate. Particles that jump change position according to a symmetric transition probability $p$. In the study of their hydrodynamic limit it is customary to take as lattice the discrete torus $\T_N^d:=\{0,1,\ldots,N-1\}$ with periodic boundary conditions, so that the phase space of the ZRP is the space $\MM_N^d:=\ZZ_+^{\T_N^d}$ of configurations of particles, and for each $\eta\in\MM_N^d$ $\eta(x)$ is the number of particles at $x\in\T_N^d$. Then the empirical process $(\pi_t^N)_{t\geq 0}$ of the ZRP $(\eta_t)_{t\geq 0}$ is defined by
		$$\pi_t^N=\fr{N^d}\sum_{x\in\T_N^d}\eta_t(x)\delta_{\frac{x}{N}}\in\MMM_+(\T^d).$$
		 Here $\MMM_+(\T^d)$ is the space of non-negative Borel measures on $\T^d$. In the study of the hydrodynamic limit of the ZRP we are interested in proving the convergence of the diffusively rescaled laws of the empirical process $(\pi_t^N)_{t\geq 0}$ to a law on the Skorohod space $D(0,T;\MMM_+(\T^d))$ that is supported by trajectories that satisfy an evolutionary partial differential equation, the so called \emph{hydrodynamic limit}.
		
		Since their introduction by Spitzer in 1970, ZRPs have attracted a lot of attention, one reason being that for particular choices of local jump rate functions $\mathfrak{g}$ they exhibit phase transition phenomena, via the emergence of mass condensation at densities above a critical density $\rho_c$. So ZRPs can serve as a  simple prototype model for the study of condensation phenomena, \cite{Bialas1997a,Evans2000a,Grosskinsky2003a}. ZRPs that can exhibit condensation are called \emph{condensing}. Obtaining the hydrodynamic limit of condensing ZRPs in closed form is a difficult open problem since none of the existing methods of proving hydrodynamic limits applies due to the lack of a replacement lemma and the fact that expected hydrodynamic equation is not always well-posed.
		
			In this article, apart from the empirical density process $(\pi_t^N)_{t\geq 0}$ we consider the empirical current $(W_t^N)_{t\geq 0}$ and the empirical jump rate $(\s_t^N)_{t\geq 0}$ processes of the ZRP defined by 
		$$W_t^N=\fr{N^{d-1}}\sum_{j=1}^d\sum_{x\in\T_N^d}\big\{\mathfrak{g}(\eta(x))-\mathfrak{g}(\eta(x+e_j))\big\}\delta_{\frac{x}{N}}e_j,\quad\s_t^N=\fr{N^d}\sum_{x\in\T_N^d}\mathfrak{g}(\eta(x))\delta_{\frac{x}{N}}$$
		and we will prove that their laws are concentrated as $N\ra+\infty$ on paths $(\pi,W,\s)$ satisfying the continuity equation 
		\begin{equation}\label{ContEq}\begin{cases}
		\pd_t\pi+\dv W=0\\
		W=-\nabla\s\end{cases}\quad\mbox{in }(0,T)\x\T^d\end{equation} in the sense of distributions, where $\pd_t$ is the time-derivative, $\dv$ is the spatial divergence and $\nabla$ is the spatial gradient. Here for (almost) all $t\geq 0$ $\pi_t$ is finite non-negative measure, $W_t$ is a vector-valued measure absolutely continuous with respect to Lebesgue measure and $\s_t\in H^1(\T^d)$. The relative compactness of the empirical density process on the Skorohod space of paths of non-negative measures is well-known by \cite[Lemma 5.1.5]{Kipnis1999a}, which as noted there is valid even for condensing zero range processes. Here we prove the relative compactness of the empirical jump rate and current in the weak-star ($w^*$) topology of appropriate $L_{w^*}^\infty$-spaces of $w^*$-measurable Banach valued paths.
		
		More importantly, we give a first closed hydrodynamic equation for condensing ZRPs in terms of a notion of generalized Young measures, using only the extension of the one-block estimate to condensing ZRPs in~\cite{Stamatakis2014} and not the full replacement lemma. Ordinary Young-measures are not sufficient to yield a closed hydrodynamic equation, since they may lose track of the mass that is in condensed phase. For this reason we consider a generalization of Young measures that, loosely speaking, allows us to encode the mass in condensed phase on a separate coordinate, namely a measure $\mu\in\MMM_+(\T^d)$. More precisely we define the generalized Young functionals as elements of the dual of the Banach space 
		$$\bbar{C}_1(\T^d\x\RR_+):=\Big\{F\in C(\T^d\x\RR_+)\Bigm|\exists\bar{F}\in C(\T^d)\colon\;\lim_{\lambda\ra\infty}\sup_{u\in\T^d}\Big|\frac{F(u,\lambda)}{1+\lambda}-\bbar{F}(u)\Big|=0\Big\}$$
		equipped with norm $\|F\|_{\infty;1}:=\sup_{(u,\lambda)\in\T^d\x\RR_+}\frac{|F(u,\lambda)|}{1+\lambda}$. Viewing the empirical distribution of the ZRP as a generalized Young-functional via 
		$$\bs\pi^{N,\ell}_t:=\fr{N^d}\sum_{x\in\T_N^d}\delta_{\frac{x}{N}}\otimes\delta_{\eta_t^\ell(x)}$$
		where $\eta^\ell(x)=\fr{(2\ell+1)^d}\sum_{|y|\leq\ell}\eta(x+y)$ and using the generalization of the one-block estimate for condensing ZRPs with bounded jump rate proved in~\cite{Stamatakis2014} we obtain that all limit points as $N$ and then $\ell$ tend to infinity of the laws of the empirical density $\bs\pi^{N,\ell}$ of the ZRP are supported on trajectories that satisfy the closed hydrodynamic equation 
		\begin{equation}\label{YMHydrLim}\pd_t\bs\pi=\Delta\Phi(\bs\pi)\end{equation}
		in the sense of distributions. Again, the path-space for the empirical density of the ZRP in terms generalized Young measures will be an appropriate $L_{w^*}^\infty$-spaces. As we will see, generalized Young-functionals are represented by pairs $(\bs\rho,\rho^\perp)\in\MMM_{1}(\T^d\x\RR_+)\x\MMM(\T^d)$ where $\bs\rho$ is an ordinary Young-measure, referred to as the regular part of the generalized Young-functional, and $\rho^\perp$ is a non-negative Borel measure on $\T^d$, referred to as the singular part of the the generalized Young-functional $\bs\pi$. Such pairs act on maps $F\in\bbar{C}_1(\T^d\x\RR_+)$ via 
		$$\ls F,\bs\pi\rs=\int F(u,\lambda)d\bs\rho(u,\lambda)+\int\bbar{F}(u)d\rho^\perp(u),\quad\bs\pi=(\bs\rho,\rho^\perp).$$ Viewing the empirical distribution of the ZRP as a generalized Young-functional, any limiting point of the family of diffusively rescaled empirical processes of the ZRP is concentrated on trajectories $\bs\pi=(\bs\rho,\mu)$ such that $\bs\rho$ contains all mass at sites with local particle density $\leq M$ and $\mu$ the rest mass as $M\ra+\infty$. 
		
		Furthermore, the usefulness of generalized Young-functionals is not restricted in yielding a closed hydrodynamic equation for the ZRP. They are also a natural setting to study the two-blocks estimate in the case of condensing ZRPS. The two-blocks estimate is one of the two basic lemmas in the proof of the replacement lemma~\cite[Lemma 5.1.10]{Kipnis1999a}, the other being the one-block estimate. The one-block estimate was proved for condensing ZRPs in~\cite{Stamatakis2014} for bounded cylinder functions and is improved in this article by extending the class of admissible cylinder functions to the class of all asymptotically linear cylinder maps. The generalized Young-functionals allow us to separate the fluid from the solid phase and we are thus able to study how the two-blocks estimate may fail in condensing ZRPs. As we will see, in general the micro-block averages $\eta^{\ell}(x)$ underestimate the fluid phase compared to the macro-averages $\eta^{[N\ee]}$. Here $[N\ee]$ is the integer part of $N\ee$. Indeed, using the notion of generalized Young measures we are able to adapt the ``cut off of large densities"~Lemma~\cite[Lemma 5.4.2]{Kipnis1999a} used in the proof of the two-blocks estimate. Using this adaptation of cutting of the large densities, Lemma~\ref{TwoBlockCutLargeDens} herein, we are able to interchange micro-block averages by macro-block averages of truncated micro-block averages, which allows us to compare the micro and macro-block averages $\eta^\ell(x)$ and $\eta^{[N\ee]}$ as $N\uparrow\infty$, $\ee\downarrow 0$ and then $\ell\uparrow\infty$. This comparison result allows us to give an equivalent condition for the validity of the two-blocks estimate.
		
		 We note here that the cut-off Lemma~\ref{TwoBlockCutLargeDens} can not be proven by the argument used in the proof of the original ``cut-off lemma" in~\cite[Lemma 5.4.2]{Kipnis1999a} since condensing the equilibrium states of condensing ZRPs do not have full exponential moments. Thus our argument via the use of generalized Young measures seems to be necessary here.


Once the correct topologies have been chosen on the path-spaces of the empirical processes described above, obtaining their relative compactness, the continuity equation~\eqref{ContEq} and the hydrodynamic limit~\eqref{YMHydrLim} in terms of generalized Young-functionals is rather straightforward and relies on an application of Prokhorov's theorem and the portmanteu theorem. Since the $L_{w^*}^\infty$-spaces are not polish spaces as is usually the case, we collect in the appendix known results of functional analysis~\cite{PilarCembranos1997a} for a precise description of the $L^\infty_{w^*}$-spaces and results of topological measure theory~ \cite{Topsoe1970a,Smolyanov1976a} according to which the usual theory of convergence of probability measures on polish spaces remains valid in completely regular submetrizable spaces. In the case of the empirical current in particular, in order to obtain the relative compactness we view it as a first order distribution. Simplifying the relative compactness by the choice of topology comes with a price though, since the limiting empirical current need not be a vector-valued measure any more. Thus the additional regularity result $\s_t\in H^1(\T^d)$ a.s$.$ is required to conclude via the identity $W=-\nabla\s$ that any limiting point of the laws of the empirical current process is concentrated on paths of vector-valued measures.

\paragraph{Plan of the paper:} In Section~\ref{PrelimSect} we recall basic facts on condensing ZRPs and describe the various empirical processes that will be considered in article. In Section~\ref{MainResults} we state our main results regarding the hydrodynamic limit of condensing ZRPs. In Section~\ref{GYM} we study the notion of generalized Young measures and their decomposition in regular and singular parts. Section~\ref{ProofS} contains the proofs of our main results. More precisely,~in Section~\ref{RelCompYoungM} we collect relative compactness results for various empirical processes in terms of generalized Young measures and prove basic properties of their limiting laws. Section~\ref{OBESection} contains the proof of the one-block estimate (Theorem~\ref{OBETheorem} herein). Sections~\ref{ContinuityEquationSection} and~\ref{EnergySection} contain the proofs of the continuity equation (Theorem~\ref{Theorem1}) and the energy estimate (Theorem~\ref{TheoremEnergyEstimate}), respectively, and Section~\ref{ACWYE} contains the proof of the closed hydrodynamic equation $\pd_t\bs\pi_t=\Delta\Phi(\bs\pi)$ in terms of generalized Young measures (Theorem~\ref{ClosedWeakEquation}). In Section~\ref{TBCSection} we give the proof of the two-blocks comparison (Theorem~\ref{TBCTheorem}) and in Section~\ref{RLSection} we prove the one-sided replacement lemma (Theorem~\ref{RLTheorem}). Finally, for the convenience of the reader, in the appendix we collect the results from functional analysis and topological measure theory that will be used throughout the main text.
 	
\section{Preliminaries}\label{PrelimSect}
	\subsection{Zero Range Processes}
	We give in this section the definition of ZRPs. A standard reference for the material in this section is the textbook \cite{Kipnis1999a} and the article \cite{Grosskinsky2003a}. A {\it{local jump rate}} is a function $\mathfrak{g}:\ZZ_+\lra\RR_+$ such that
   \begin{subequations}
	\label{equations}
		\begin{align}
		\label{LocalJumpRateDef:a}
		&\mathfrak{g}(k)=0\quad\Longleftrightarrow\quad k=0\\
		\label{LocalJumpRateDef:b}
		&\|\mathfrak{g}'\|_\infty:=\sup\nolimits_{k\in\ZZ_+}|\mathfrak{g}(k+1)-\mathfrak{g}(k)|<+\infty,\mbox{ and}\\
		\label{LocalJumpRateDef:c}
		&\mbox{The limit }\varphi_c:=\lim\nolimits_{k\rightarrow\infty}\mathfrak{g}!(k)^{\fr{k}}>0\mbox{ exists.}\qquad\qquad \qquad\qquad\qquad 
	\end{align}
		\end{subequations}	
	where $\mathfrak{g}!(0):=1$ and $\mathfrak{g}!(k):=\mathfrak{g}(1)\cdot \mathfrak{g}(2)\cdot\dots\cdot \mathfrak{g}(k)$. Assumption~\eqref{LocalJumpRateDef:a} of local jump rates means that the rate at which particles leave a site is zero if and only if the site is empty, assumption~\eqref{LocalJumpRateDef:b} is necessary to define the ZRP on the infinite lattice $\MM_\infty^d:=\ZZ_+^{\ZZ^d}$. In the study of the hydrodynamic it is used to bound certain quantities by the total number of particles and can be relaxed the assumption that there exists $C>0$ such that $\mathfrak{g}(k)\leq (1+C)k$. Finally assumption~\eqref{LocalJumpRateDef:c} is mainly required for the equivalence of ensembles to hold. An {\it{elementary step distribution}} is a probability distribution $p\in\PP\ZZ^d$ (where for any polish space $M$ we denote by $\PP M$ the space of all Borel probability measures on $M$) such that its support $\{z\in\ZZ^d|p(z)>0\}$ is bounded and the markov kernel $p(x,y):=p(y-x)$
	is an irreducible random walk kernel.
	
	\indent Instead of defining the ZRP on the infinite lattice $\MM_\infty^d$ it is more convenient to consider ZRPs that evolve on the discrete $d$-dimensional tori $\T_N^d\cong(\,^\ZZ/_{N\ZZ})^d\cong\{0,1,\dots,N-1\}^d$, $N\in\NN$, and consider the limit as $N\ra\infty$. The state space of a ZRP evolving on $\T_N^d$ is the space of configurations $$\MM_N^d:=\ZZ_+^{\T_N^d}.$$ Elements of $\MM_N^d$ are usually denoted by $\eta$, in which case $\eta_x$ is the number of particles occupying the site $x\in\T_N^d$. We will denote by $\eta(x):\MM_N^d\lra\ZZ_+$, $x\in\T_N^d$, the natural projections. The \emph{(diffusively rescaled) symmetric nearest neighbour ZRP of local jump rate function} $\mathfrak{g}$ {\it{on the discrete torus}} $\T_N^d$ is the unique Markov jump process on the Skorohod path-space $D(\RR_+;\MM_N^d)$ with generator $L^N\colon D(L_N)\leq B(\MM_N^d)\lra B(\MM_N^d)$ given by the formula 
	$$L^Nf(\eta)=N^2\sum_{x,y\in\T_N^d}\big\{f(\eta^{x,y})-f(\eta)\big\}\mathfrak{g}\big(\eta_x\big)p(y-x),$$
	where $\eta^{x,y}$ is the configuration resulting from $\eta$ by moving a particle from $x$ to $y$ and $p$ is the nearest neighbour random walk kernel rescaled so as to have total ``probability" equal to $2d$, i.e.~$p(e_j)=p(-e_j)=1$, $j=1,\ldots,d$. The factor $N^2$ corresponds to the diffusive rescaling necessary, due to the symmetry of the kernel $p$, in order for the macroscopic profile to evolve. The rescaling of the kernel $p$ is made so that a coefficient $\fr{2d}$ that would otherwise appear in the hydrodynamic equation is set to $1$.
	We denote
	by $S^N_t\colon\MM_N^d\lra\PP\MM_N^d$ the transition semi-group of the ZRP. We will denote by $(P_N^\eta)_{\eta\in\MM_N^d}\subs\PP D(\RR;\MM_N^d)$ the diffusively rescaled Markov kernel defined by the generator $L_N$. Given a sequence of initial distributions $\{\mu_0^N\in\PP\MM_N^d\}_{N\in\NN}$ we will write 
	$$P_N^{\mu_0^N}:=\int_{\MM_N^d} P_N^\eta\df\mu_0^N(\eta)$$
	the law of the ZRP starting from $\mu_0^N$, and if the sequence of initial distribution is fixed we will simply write $P^N:=P_N^{\mu_0^N}$. 
	
	 The total number of particles is conserved by the stochastic dynamics of the ZRP. Furthermore, the assumption that the support of the elementary step distribution $p$ generates $\ZZ^d$ implies that all configurations with the same number of particles communicate and therefore the communication classes of the ZRP are exactly the hyperplanes 
	$$\MM_{N,K}^d:=\bigg\{\eta\in\MM_N^d\Big|\sum_{x\in\T_N^d}\eta(x)=K\bigg\},\quad K\in\ZZ_+,$$
	with a fixed number of particles, and for each $(N,K)\in\NN\x\ZZ_+$, there exists a unique equilibrium state $\nu_{N,K}^d\in\PP\MM_N^d$ concentrated on $\MM_{N,K}^d$. We will refer to the family $\{\nu_{N,K}\}_{(N,K)\in\NN\x\ZZ_+}$ as {\it{the canonical ensemble}} of the ZRP.
	
	We will consider $\MM_N^d$ embedded in $\MM_\infty^d$ via the periodic extension $\MM_N^d\ni\eta\mapsto\wt{\eta}\in\MM_\infty^d$ of configurations. This is defined via the pull-back of the natural projection $q_N\colon\ZZ^d\to \T_N^d$, $q_N(z)=z+N\ZZ^d$, i.e. 
	$$\wt{\eta}(z):=q_N^*\eta(z)=\eta(z+N\ZZ^d),\quad\forall\eta\in\MM_N^d.$$
	We will also consider the finite lattice $\T_N^d$ embedded in $\ZZ^d$ via the left inverse $j_N$ of the natural projection $q_N$ for which $j_N(\T_N^d)=\{-[\frac{N}{2}],\ldots,[\frac{N-1}{2}]\}^d$. For any $J\subs\ZZ^d$ we will write $\MM_J^d:=\ZZ_+^J$ so that with the identification $\T_N^d\subs\ZZ^d$ just described $\MM_N^d=\MM_{\T_N^d}^d$. Furthermore, for $\eta\in\MM_J^d$ we will write $|\eta|_{J,1}:=\sum_{x\in J}\eta(x)$ for the $\ell_1$-norm of the configuration $\eta$ and set $|\cdot|_{N,1}:=|\cdot|_{\T_N^d,1}$. A map $\Psi\colon\MM_\infty^d\to\RR$ is called a \emph{cylinder map} if it depends on a finite number of coordinates, i.e.~ if there exists a finite $J\subs\ZZ^d$ for which there exists map $\Psi_J\colon\MM_J^d\to\RR_+$ such that $\Psi=\Psi_J\circ\eta^J$. Here $\eta^J\colon\MM_\infty^d\to\MM_J^d$ denotes the natural projection. Such a set $J$ will be called a \emph{supporting set} for $\Psi$. If $J,K\subs\ZZ^d$ are supporting sets for the cylinder map $\Psi$ then the set $J\cap K$ is also a supporting set for $\Psi$ and thus for each cylndric map $\Psi$ there exists a unique minimal supporting set for $\Psi$ which will be called the support of $\Psi$ and will be denoted by $J_\Psi$, i.e. 
	$$J_\Psi=\bigcap\big\{J\subs\ZZ^d\bigm|J\mbox{ finite for which }\exists\Psi_J\colon\MM_J^d\to\RR\mbox{ such that }\Psi=\Psi\circ\eta^J\big\}.$$
	The number $k_\Psi:=\sharp J_\Psi$ of elements of the support of $\Psi$ is called \emph{the size of $\Psi$}. The cylinder map $\Psi=\Psi_J(\eta^J)$ is called \emph{sublinear} if in addition  
	\begin{equation}\label{SublinearCyl}\lim_{|\eta|_{J,1}\ra+\infty}\frac{\Psi_J(\eta)}{|\eta|_{J,1}}=0.\end{equation} This definition does not depend on the choice of $J\subs\ZZ^d$ and $\Psi_J\colon\MM_J^d\to\RR_+$ for which $h\Psi=\Psi_J(\eta_J)$. 
	
	We will use the following notation for spatial averages of a cylinder map $\Psi\colon\MM_N^d\to\RR$. For $\ell,L\in\ZZ_+$ we write
	\begin{equation}\label{MicroBlockCyl}
	\Psi^\ell(\eta):=\fr{(2\ell+1)^d}\sum_{|y|\leq\ell}\tau_y\Psi(\eta)
	\end{equation}
	for the block average of radius $\ell$ we will denote by
	\begin{equation}\label{DoubleBlock}
	\Psi^{\ell,L}=\Psi^{\ell,L}:=\fr{L_\star^d\ell_\star^d}\sum_{|y|\leq\ell}\sum_{|z|\leq L}\tau_{y+z}\Psi
	\end{equation}
	the double block average of radii $\ell$ and $L$. Thus for example $\eta^\ell(x)=\fr{(2\ell+1)^2}\sum_{|y|\leq\ell}\eta(x+y)$.
	
	 The function $Z\equiv Z_\mathfrak{g}:\RR_+\lra[1,\infty]$ defined by the power series $$Z(\varphi):=\sum_{k=0}^\infty\frac{\varphi^k}{\mathfrak{g}!(k)}
	$$ is called the {\it{normalising partition function associated to}} the local jump rate function $\mathfrak{g}$. The radius of convergence of $Z$ is 
	$\varphi_c=\liminf_{k\ra\infty}\sqrt[k]{\mathfrak{g}!(k)}$ and so assumption~\eqref{LocalJumpRateDef:c} of local jump rate functions guaranties that $Z$ has non-trivial domain of convergence. Obviously any partition function $Z:\RR_+\lra[1,+\infty]$ is $C^\infty$ on $[0,\phih_c)$ with all of its derivatives strictly positive there. By Abel's theorem on power-series, $Z$ and all of its derivatives are lower semi-continuous on $\RR_+$. For any $\phih\in\DD_Z:=\{\phih\in\RR_+|Z(\phih)<+\infty\}$, the product distribution $\bar{\nu}_{\phih}^N\equiv\bar{\nu}_{\phih,\mathfrak{g}}^N\in\PP\MM_N^d$ with common marginal $\bar{\nu}^1_{\phih}\in\PP\ZZ_+$ given by 
	$$
	\bar{\nu}^1_\varphi\{k\}=\fr{Z(\varphi)}\frac{\varphi^k}{\mathfrak{g}!(k)},\qquad k\in\ZZ_+
	$$
	is called the {\it{zero range product distribution on $\T_N^d$ of rate $g$ and fugacity}} $\varphi$.
	
	 Note that whenever $\phih_c\in\DD_Z$ the one-site zero range distribution $\bar{\nu}_{\phih_c}^1$ corresponding to the critical fugacity $\phih_c$ is defined. The zero range product distributions $\bar{\nu}_{\phih}^N\in\PP\MM_N^d$, $\phih\in \DD_Z$, are equilibrium distributions, i.e.~$\bar{\nu}_{\phih}^NL^N = 0$,
	and translation invariant, that is 
	$\tau_{x\sharp}\bar{\nu}_{\phih}^N:=\bar{\nu}_{\phih}^N\circ\tau_x^{-1}=\bar{\nu}_{\phih}^N$ for all $x\in\T_N^d$, where $\tau_x:\MM_N^d\lra\MM_N^d$ denotes the translation operator $(\tau_x\eta)_y=\eta_{x+y}$. In fact they are the only translation invariant equilibrium states of the ZRP that are also product measures. Let us note here that for any measurable map $f\colon M\to N$ and any measure $\mu$ on $M$ we will denote by $f_\sharp\mu:=\mu\circ f^{-1}$ the push forward measure of $\mu$ on $N$.
 
  As is well known, the zero range product distributions can be reparametrised by the density. The {\it{mean density function}} $R:\DD_Z\lra[0,+\infty]$ defined by 
\begin{equation}\label{MeanDensityMainFormula}\qquad
R(\phih)=E_{\bar{\nu}_{\phih}^N}[\eta(0)]=\int kd\bar{\nu}^1_\phih(k)=\frac{\phih Z'(\phih)}{Z(\phih)}
\end{equation} is continuous on $\DD_Z$, it
is obviously $C^\infty$ on $[0,\phih_c)$, and as shown in \cite{Kipnis1999a}, it is strictly increasing. Consequently, it's inverse $\Phi:=R^{-1}:R(\DD_Z)\lra\DD_Z$ is well defined. Of course $[0,\rho_c)\subs R(\DD_Z)\subs[0,\rho_c]$, where 
\begin{equation}\label{CriticalDensity}
\rho_c\equiv R(\phih_c):=\lim_{\phih\uparrow\phih_c}R(\phih)\in(0,\infty],
\end{equation}
and $\rho_c\in R(\DD_Z)$ iff $\phih_c\in\DD_Z$. The number $\rho_c$ defined in~\eqref{CriticalDensity} is called the critical density. We will say that a ZRP is a \emph{condensing ZRP} if $\phih_c<+\infty$ and we will say that a condensing ZRP is \emph{weakly condensing ZRP} if $\rho_c=+\infty$ and \emph{strongly condensing ZRP} if $\rho_c<+\infty$.  By reparametrising the zero-range distributions by the mean jump rate $\Phi$ we get for any $\rho\in R(\DD_Z)$ an equilibrium distribution $\nu^N_\rho$ of mean density $\rho$:
\begin{eqnarray}\label{GrandCanEns}\qquad
\nu^N_\rho:=\bar{\nu}^N_{\Phi(\rho)},\quad\rho\in R(\DD_Z).
\end{eqnarray}
	We will refer to the family defined in~\eqref{GrandCanEns} as the \emph{grand canonical ensemble} of the ZRP.
	
	 The various possibilities for the set $R(\DD_Z)$ are as follows. As is proved in \cite{Kipnis1999a}, whenever $\phih_c\notin\DD_Z$, that is whenever $\phih_c=+\infty$ or $\phih_c<+\infty$ and $Z(\phih_c)=+\infty$, we have that $\rho_c=+\infty$. In this case $R(\DD_Z)=\RR_+$, there is no equilibrium state $\bar{\nu}_\phih^1$ corresponding to the critical fugacity $\phih=\phih_c$ and the mean jump rate function $\Phi$ is defined on all of $\RR_+$. On the other hand if $\phih_c\in\DD_Z$ then $R(\DD_Z)=[0,\rho_c]$ and in this case, as is shown by~\eqref{MeanDensityMainFormula}, the critical density is finite if and only if $Z'(\phih_c)\equiv\sup_{\phih<\phih_c}Z '(\phih)<+\infty$. In particular, whenever $\rho_c<+\infty$ we have that $\phih_c\in\DD_Z$ and so the grand canonical ensemble contains the equilibrium distribution $\bar{\nu}_{\phih_c}^1$ with density equal to the critical density $\rho_c$. Note that in the special case that $\phih_c\in\DD_Z\sm\DD_{Z'}$ there exists an equilibrium state $\nu_\infty^N=\bar{\nu}^N_{\phih_c}$ corresponding to infinite density $\rho_c=+\infty$ and $R(\DD_Z)=[0,+\infty]$.
	
	\begin{example}{\rm{(The Evans Model)}}	As an example of a condensing ZRP in \cite{Evans2000a} Evans introduces ZRPs with local jump rate function 
	\begin{eqnarray}\label{EvansJumpRate}\qquad
	\mathfrak{g}_b(k)=\1_{\{k\geq 1\}}\Big(1+\frac{b}{k}\Big),\quad b\geq 0.
	\end{eqnarray}
	It is well known~(\cite{Grosskinsky2003a}) that $\phih_c=1$ for all $b\geq 0$ and that the corresponding ZRP is weakly condensing for $b\in[0,2]$ and strongly condensing for $b>2$ with critical density $\rho_c=\rho_c(b)=\fr{b-2}$. In fact for $b\in[0,1]$ we have that $\phih_c\notin\DD_Z$ and there is no equilibrium state with critical mean density $\rho_c=+\infty$ while for $b>1$ the critical equilibrium state $\bar{\nu}^1_{\phih_c}$ scales as $k\ra+\infty$ as a polynomial distribution of order $k^{-b}$. Thus more precisely $\phih_c\notin\DD_Z$ so that $R(\DD_Z)=[0,+\infty)$ iff $b\in[0,1]$, $\phih_c\in\DD_Z\sm\DD_{Z'}$ so that $R(\DD_Z)=[0,+\infty]$ iff $b\in(1,2]$ and finally for $b>2$ we have that $\phih_c\in\DD_R=\DD_{Z'}$ so that the first moment of the grand canonical distribution $\bar{\nu}_{\phih_c}^1$ is finite, thus leading to a finite critical density $\rho_c<\infty$ and $R(\DD_Z)=[0,\rho_c]\subs\RR_+$. For $b>3$ the critical equilibrium state $\nu_{\rho_c}^1$ has finite second order moments and $R_-'(\phih_c)<+\infty$ while $R_-'(\phih_c)=+\infty$ for $b\in[0,3]$. A precursor of the Evans model has been studied in \cite{Drouffe1998a}.
	\end{example}

		We note that the mean jump rate function $\Phi$ is Lipschitz with Lipschitz norm $\leq\|\mathfrak{g}'\|_\infty$ and is {\it{the mean jump rate function}} since for all $\rho\in R(\DD_R)$ we have that
	$$E_{\nu_\rho^N}[\mathfrak{g}(\eta(0))]=\int\mathfrak{g}(k)d\nu_\rho^1(k)=\fr{Z(\Phi(\rho))}\sum_{k=0}^\infty\mathfrak{g}(k)\frac{\Phi(\rho)^k}{\mathfrak{g}!(k)}= \Phi(\rho).$$
	More generally, for any cylinder map $\Psi\colon\MM_\infty^d\to\RR_+$ we define the \emph{(grand canoninical) homologue map $\wt{\Psi}\colon R(\DD_R)\to\RR$ of} $\Psi$ by 
	$$\wt{\Psi}(\rho)=\int\Psi\df\nu_\rho^\infty,\quad\rho\in R(\DD_R).$$
	With this definition we have that $\Phi=\wt{\mathfrak{g}(\eta(0))}$.

The logarithmic moment generating function $\Lambda_{\rho_*}:=\Lambda_{\nu^1_{\rho_*}}$ of $\nu^1_{\rho_*}$, $\rho_*\in(0,\rho_c)$, given by 
$$\Lambda_{\rho_*}(\theta)=\log\frac{Z\big(e^\theta\Phi(\rho_*)\big)}{Z\big(\Phi(\rho_*)\big)},$$ has proper domain $\mathcal{D}_{\Lambda_{\rho_*}}$ such that $(-\infty,b_{\rho_*})\subs\mathcal{D}_{\Lambda_{\rho_*}}\subs(-\infty,b_{\rho_*}]$, where $b_{\rho_*}:=\log\frac{\phi_c}{\Phi(\rho_*)}>0$. In particular when $\phih_c=+\infty$ then $\nu^1_{\rho_*}$ has full exponential moments for all $\rho_*\geq 0$, that is $\Lambda_{\rho_*}(\theta)=\int e^{\theta k}d\nu^1_{\rho_*}(k)<+\infty$ for all $\theta\in\RR$, $\rho_*\geq 0$. If $\phih_c<+\infty$ then $\nu^1_{\rho_*}$ has some exponential moments if $\rho_*<\rho_c$ while at the critical density $\rho=\rho_c$ we have that $b_{\rho_c}=0$ and $\nu^1_{\rho_c}$ does not have exponential moments.

	The phase transition in ZRPs with finite critical density has been described in \cite{Evans2000a} and proved rigorously in \cite[Theorem 1]{Grosskinsky2003a} as a continuous phase transition in the thermodynamic limit by using the relative entropy $\HHH(\cdot|\cdot)$ to count the distance between the canonical and grand canonical ensemble, which in general for any probability measures $\mu,\nu$ on a measurable space $(M,\F)$ is defined by 
	$$\HHH(\mu|\nu)=\begin{cases}
	\int\frac{d\mu}{d\nu}\log\frac{d\mu}{d\nu}d\nu\quad&{\rm{if }}\;\mu\ll\nu\\
	+\infty\quad&{\rm{otherwise}}
	\end{cases}.$$
	Here as usual the convention $0\log 0=\lim_{t\downarrow 0}t\log t=0$ is made. A useful inequality is the so called relative entropy inequality according to which for any bounded measurable $f\colon M\to\RR$
	\begin{equation}\label{RelEntrIneq}
		\int f\df\mu\leq\inf_{\theta>0}\fr{\theta}\Big\{\log\int e^{\theta f}\df\nu+\HHH(\mu|\nu)\Big\}
	\end{equation}
	 To be precise, the equivalence of ensembles states that if $\pi^L:\MM_N^d\lra\MM_L^d$, $N\geq L$, are the natural projections and we set $\nu_{N,K}^L:=\pi^L_\sharp \nu_{N,K}$, then for fixed $L\in\NN$, for all $\rho\geq 0$ it holds that 
	\begin{equation}\label{EoE}
		\lim_{\substack{N,K\ra+\infty\\K/N^d\ra\rho}}\HHH(\nu_{N,K}^L|\nu_{\rho\mn\rho_c}^L)=0.\end{equation} 
		In particular $\nu_{N,K}^L\lra\nu_{\rho\mn\rho_c}^L$ weakly as $N,K\ra\infty$ and $K/N^d\ra\rho$.

 An elegant application of this result has been recently given in \cite{Chleboun2014a}, where it is shown that for subcritical densities $\rho\leq\rho_c$ the equivalence of ensembles~\eqref{EoE} can be applied to yield weak convergence in duality with respect cylinder maps $\Psi\in L^{1+\ee}(\nu^\infty_\rho)$ for some $\ee>0$. As we will see in Lemma~\ref{WassersteinConvergenceEoE} this implies that for $\rho\leq\rho_c$ the weak convergence $\nu_{N,K}^L\lra\nu_{\rho}^L$ is in duality with respect all cylinder functions that have at most linear growth and that in the case $\rho>\rho_c$ the weak convergence $\nu_{N,K}^L\lra\nu_{\rho_c}^L$ can be strengthened to convergence in duality with respect to all sublinear cylinder maps $\Psi$. Of course this cannot strengthened to linear cylinder maps for $\rho>\rho_c$ since even for the linear cylinder function $\eta(0)$ 
	$$\int\eta(0)d\nu_{N,K}\lra\rho>\rho_c\qquad\mbox{as }\;N,K\ra\infty\;\;\mbox{and }\;K/N^d\ra\rho.$$
	In other words, at the thermodynamic limit we have a mean total loss of mass equal to $\rho-\rho_c$ at each site. As it has been proven, in many cases the excess mass of all the sites is concentrated on a single random site. We refer to \cite{Grosskinsky2003a,Armendariz2008a,Armendariz2009a} for a detailed description of the phase separation in the Evans model.\\
	\indent In particular the equivalence of ensembles yields via Lemma~\ref{WassersteinConvergenceEoE} that for any sublinear cylinder map $\Psi\colon\MM_\infty^d\to\RR_+$
	$$\lim_{\substack{N,K\ra\infty\\K/N^d\ra\infty}}\int \Psi d\nu_{N,K}=\int \Psi d\nu_{\rho\mn\rho_c}^\infty=\wt{\Psi}(\rho\mn\rho_c),$$ 
	for all $\rho\geq 0$. Thus for any sublinear map $\Psi\colon\MM_\infty^d\to\RR$ we define its extended homologue map $\bbar{\Psi}$ by extending $\wt{\Psi}$ on all of $\RR_+$ via
	\begin{equation}\label{ExtendedHomologueSublinear}\bbar{\Psi}(\rho)=\wt{\Psi}(\rho\mn\rho_c),\quad\mbox{ for all }\rho\geq 0.\end{equation} This extension has been considered in the particular case of the mean jump rate function $\Phi$ for bounded local jump rate functions $\mathfrak{g}$ in \cite{Grosskinsky2003a} and also in \cite{Harris2007a} which contains a heuristic discussion on the hydrodynamics of strongly assymetric ZRPs in the Eulerian scaling. It turns out~\cite{Stamatakis2014} that this choice of $\Phi$ is the right one in order to extend the one-block estimate to ZRPs with finite critical density. As we will see in this article the one-block estimate in condensing ZRPs holds in general for sublinear maps. In the case of weakly condensing ZRPs the one-block estimate holds for asymptotically linear cylinder maps, where a cylinder map $\Psi\colon\MM_\infty^d\to\RR$ is called {asymptotically linear} if there exists $a=(a_x)_{x\in J_\Psi}\in\RR^{J_\Psi}$ such 
	$$\lim_{|\eta|_{J_\Psi,1}\ra\infty}\Big|\frac{\Psi_{J_\Psi}(\eta)}{\ls a,\eta\rs}-1\Big|=0.$$ 
	If such $a\in\RR^{J_\Psi}$ exists then it is unique and it is denoted by $\nabla\Psi(\infty)=(\pd_x\Psi(\infty))_{x\in J_\Psi}$.
	Of course here $\ls a,\eta\rs=\sum_{x\in J_\Psi}a_x\eta(x)$. Furthermore, if we want to extend the one block estimate in strongly condensing ZRPs to asymptotically linear maps we have to define the extended homologue $\bbar{\Psi}$ of an asymptotically linear cylinder map $\Psi\colon\MM_\infty^d\to\RR$ by 
	\begin{equation}\label{ExtendedHomologueAsymptlinear}\bbar{\Psi}(\rho)=\wt{\Psi}(\rho\mn\rho_c)+\ls\nabla\Psi(\infty),\1_J\rs(\rho-\rho_c)^+,\quad\rho\geq 0.\end{equation}
Of course here $J_\Psi\subs\ZZ^d$ is the support of the cylinder map $\Psi$ and $\ls\nabla\Psi(\infty),\1_{J_\Psi}\rs=\sum_{x\in J_\Psi}\pd_x\Psi(\infty)$. We note that in the case of weakly condensing ZRPs and for sublinear cylinder maps $\Psi$ in the case of strongly condensing ZRPs formula~\eqref{ExtendedHomologueAsymptlinear} reduces to~\eqref{ExtendedHomologueSublinear}.
	
		So far, the hydrodynamic limit of ZRPs has only been proven under the assumption that $\phih_c=+\infty$ for $L^2$ initial profiles via the entropy method of Guo-Papanikolaou-Varadhan and in the case that $\phih_c\notin\DD_Z$ for $C^{2+\theta}$ initial profiles via the relative entropy method of H.T Yau, which both exclude ZRPs with finite critical density. The hydrodynamic limit was extended in~\cite{Stamatakis2014} to strongly condensing ZRPs with bounded jump rates for which the assumption $\phih_c\notin\DD_Z$ is not satisfied, but only in the case that we start the process from some $C^{2+\theta}$ strictly sub-critical initial profile $\rho_0$,~i.e. $\sup_{u\in\T^d}\rho_0(u)<\rho_c$.
	\subsection{The empirical processes} 
	In this section we briefly describe the various empirical processes that we will use to obtain information on the hydrodynamic behaviour of condensing ZRPs.
	\subsubsection{Empirical densities and the empirical jump rate}
	\emph{The empirical density} is the function $\pi^N\colon\MM_N^d\to\MMM_+(\T^d)$ given by $$\pi^N(\eta)=\fr{N^d}\sum_{x\in\T_N^d}\eta(x)\delta_{\frac{x}{N}}$$ and by a slight abuse of notation we continue to denote by $\pi^N$ the \emph{empirical density process} $\pi^N\colon D(0,T;\MM_N^d)\to D(0,T;\MMM_+(\T^d))$ induced on the Skorohod spaces by $\pi^N(\eta)(t):=\pi^N_{\eta_t}$. Since $\MM_N^d$ has the discrete topology the induced map $\pi^N$ on the Skorohod spaces is continuous regardless of the topology considered on $\MMM_+(\T^d)$. Here $\MMM_+(\T^d)$ denotes the set of all non-negative finite Borel measures equipped with the weak topology i.e.~the $w^*$-topology is inherits as a subspace of $C(\T^d)^*$. Even though the $w^*$ topology is never metrizable, the restriction of the $w^*$-topology of $\MMM(\T^d)=C(\T^d)^*$ on the cone $\MMM_+(\T^d)$ of non-negative measures is metrizable by a complete metric, and thus is a polish space. Such a metric $d$ on $\MMM_+(\T^d)$ is defined \cite[Section 4.1]{Kipnis1999a} by choosing a dense family $\{f_k\}_{k=1}\subs C(\T^d)$ with $f_1\equiv 1$ and setting 
	$$d(\mu,\nu)=\sum_{k\in\NN}\fr{2^k}\frac{|\ls\mu,f_k\rs-\ls\nu,f_k\rs|}{1+|\ls\mu,f_k\rs-\ls\nu,f_k\rs|}.$$ In what follows, the Skorohod space $D\big(0,T;\MMM_+(\T^d)\big)$ is considered with respect to this metric on $\MMM_+(\T^d)$.
	
	\emph{The empirical jump rate} is the map $\s^N\colon\MM_N^d\to\MMM_+(\T^d)$ defined by 
	$$\s^N(\eta)=\fr{N^d}\sum_{x\in\T_N^d}\mathfrak{g}(\eta(x))\delta_{\frac{x}{N}}.$$ Since the empirical jump rate $\s^N$ is not a conserved quantity in order to obtain the relative compactness of $\s^N$ we have to consider a weaker topology than the Skorohod one for the path space of the empirical jump rate process. We do this by considering the empirical jump rate process as a random variable taking values on the dual space $L^1(0,T;C(\T^d))^*$ equipped with the $w^*$-topology and the corresponding Borel $\s$-algebra. Since $C(\T^d)^*\cong\MMM(\T^d)$ does not have the Radon-Nikodym property, the dual space $L^1(0,T;C(\T^d))^*$ is not isometric to the space $L^\infty(0,T;\MMM(\T^d))$ of strongly measurable maps. Following \cite{PilarCembranos1997a} will give a precise description of the dual $L_{w^*}^\infty(0,T;X^*)$ of $L^1(0,T;X)$ for any Banach space $X$ in the appendix. Since $\MM_N^d$ has the discrete topology the map $\s^N$ is continuous and thus the induced mapping $\s^N\colon D(0,T;\MM_N^d)\to D(0,T;\MMM(\T^d))$ on the Skorohod spaces is continuous. Here we consider the space $\MMM(\T^d)$ equipped with the total variation norm. By composing this induced mapping with the continuous injection from $D(0,T;\MMM(\T^d))$ to $L_{w^*}^\infty(0,T;\MMM(\T^d))$ given in Proposition~\ref{SkorohodToLEmbed}, we obtain the \emph{the empirical jump rate process} as the continuous random variable
	$$\s^N\colon D(0,T;\MM_N^d)\to L_{w^*}^\infty(0,T;\MMM(\T^d)).$$ 
	Here continuity is with respect to the $w^*$-topology on the target space $L_{w^*}^\infty(0,T;\MMM(\T^d))$ and thus $\s^N$ is a random variable with respect to the Borel $\s$-algebra of the $w^*$-topology of $L^\infty_{w^*}(0,T;\MMM(\T^d))\cong L^1(0,T;C(\T^d))^*$. With respect to the duality pairing $\lls\cdot,\cdot\rrs$ between $L^1(0,T;C(\T^d))$ and its dual space $L_{w^*}^\infty(0,T;\MMM(\T^d))$ the empirical jump rate process is given by
 $$\lls f,\s^N\rrs=\int_0^T\fr{N^d}\sum_{x\in\T_N^d}f_t\Big(\frac{x}{N}\Big)\mathfrak{g}(\eta_t(x))\df t,\quad f\in L^1(0,T;C(\T^d)).$$

  With a slight abuse of notation we can also view the empirical distribution $\pi^N$ as taking values in $L_{w^*}^\infty(0,T;\MMM_+(\T^d))$ via the natural continuous injection from $D(0,T;\MMM_+(\T^d))\to L_{w^*}^\infty(0,T;\MMM_+(\T^d))$ of Proposition~\ref{SkorohodToLEmbed}. Since the injection $i$ is continuous, relative compactness of the laws $\pi^N_\sharp P^N\in\PP D(0,T;\MMM_+(\T^d))$, $N\in\NN$ of the empirical density process in the Skorohod space implies the relative compactness of the laws of the empirical density process also when viewed as taking values in $L_{w^*}^\infty(0,T;\MMM_+(\T^d))$.
  
   	More generally, for any cylinder map $\Psi\colon\MM_\infty^d\to\RR$ we will denote by $\s^{N,\Psi}\colon D(0,T;\MM_N^d)\to L_{w^*}^\infty(0,T;\MMM(\T^d))$ the empirical distribution process defined by duality via 
  \begin{equation}\label{GenEmpProc}
  \lls f, \s^{N,\Psi}\rrs=\int_0^T\fr{N^d}\sum_{x\in\T_N^d}f_t\Big(\frac{x}{N}\Big)\tau_x\Psi\df t.\end{equation}
  With this notation, $\pi^N=\s^{N,\eta(0)}$ and $\s^N=\s^{N,\mathfrak{g}(\eta(0))}$. With similar reasoning as in the definition of the empirical jump rate process this is a continuous map and thus a well defined random variable.
 

	\subsubsection{The empirical current}
	 The {\it{current along the bond}} $(x,y)\in\T_N^d\x\T_N^d$ for the ZRP in the discrete torus $\T_N^d$ is the function $W^N_{x,y}:\MM_N^d\lra\RR$ given by 
	\begin{equation}
	W^N_{x,y}(\eta)=L^N(\eta,\eta^{x,y})-L^N(\eta,\eta^{y,x})=[\mathfrak{g}(\eta_x)-\mathfrak{g}(\eta_y)]p(y-x)
	\end{equation}
	for all $\eta\in\MM_N^d$. The \emph{empirical current map} is the function $W^N\colon\MM_N^d\to \MMM(\T^d;\RR^d)$ defined by
	$$W^N=\fr{N^{d-1}}\sum_{i=1}^d\Big(\sum_{x\in\T_N^d}W^N_{x,x+e_i}\delta_{\frac{x}{N}}\Big)\cdot e_i=:\sum_{i=1}^dW^{N,i}\cdot e_i.$$ 
	An initial idea is to regard the empirical current process as a random variable $W^N\colon D(0,T;\MM_N^d)\mapsto L^1(0,T;C(\T^d;\RR^d))^*$ where the target space is considered equipped with its $w^*$-topology. However, proving the relative compactness of the empirical current in this $w^*$-convergence turns out to be difficult. We note that the empirical current has deterministically zero total current, that is $W^N(\T^d)\equiv 0$ on $\MM_N^d$. As such we can regard the empirical current map as taking values on the target space $\MMM_0(\T^d;\RR^d)$ of $\RR^d$-valued measures $W$ with zero total mass $W(\T^d)=0\in\RR^d$ and consider $\MMM_0(\T^d;\RR^d)$ as a subspace of the dual of \emph{the space $\X^1(\T^d):=C^1(\T^d;\RR^d)/\RR^d$ of $C^1$ vector fields $G\colon\T^d\to\RR^d$ modulo constants} equipped with the norm 
	$$\|G\rN:=\sup_{u\in\T^d}|DG(u)|=\big\||DG|\big\|_{C(\T^d)}=\|DG\|_{C(\T^d;\RR^{d\x d})}$$ where $DG$ is the derivative of $G$ and for a matrix $A\in\RR^d\x\RR^d$ we denote by $|A|^2:=\sum_{i,j=1}^da_{ij}^2$ its Frobenius norm. The space $\X^1(\T^d)$ is a separable Banach space since it is by definition isometric to the closed subspace $$\{DG|G\in\X^1(\T^d)\}\leq C(\T^d;\RR^{d\x d})$$ of the separable Banach space $C(\T^d;\RR^{d\x d})\cong C(\T^d)^{d^2}$.
	
	We will view the linear space $\MMM_0(\T^d;\RR^d)$ as a subspace of $\X^1(\T^d)^*$ via the natural injection $I\colon\MMM_0(\T^d;\RR^d)\hookrightarrow\X^1(\T^d)^*$ defined by $$I_W(G):=\int\ls G,\df W\rs\equiv\sum_{i=1}^d\int G^i\df W^i,\quad G=(G^i)_{i=1}^d\in\X^1(\T^d).$$ In this way we will identify each $W\in\MMM_0(\T^d;\RR^d)$ with $T_W\in\X^1(\T^d)^*$ and write $W(G)=I_W(G)$ for $G\in\X^1(\T^d)$ and the norm of a current $W\in\MMM_0(\T^d;\RR^d)\leq\X^1(\T^d)^*$ is given by 
	$$\|W\|_{\X^1(\T^d)^*}=\sup_{\|G\rN\leq 1}\int\ls G,\df W\rs.$$ 
	Via this embedding we consider the empirical current map as taking values in $\X^1(\T^d)^*$ i.e.
	$$W^N(G)=\fr{N^{d-1}}\sum_{j=1}^d\sum_{x\in\T_N^d}G^i\Big(\frac{x}{N}\Big)W_{x,x+e_i},\quad G=(G^i)_{i=1}^d\in\X^1(\T^d).$$
	Then the map $W^N\colon\MM_N^d\to\X^1(\T^d)^*$ defined above induces the continuous map $W^N\colon D(0,T;\MM_N^d)\to D(0,T;\X^1(\T^d)^*)$ which in turn by composing with the continuous injection $D(0,T;\X^1(\T^d)^*)\hookrightarrow L_{w^*}^\infty(0,T;\X^1(\T^d)^*)$ given in Proposition~\ref{SkorohodToLEmbed} yields \emph{the empirical current process} $$W^N\colon D(0,T;\MM_N^d)\to L_{w^*}^\infty(0,T;\X^1(\T^d)^*)$$ as a continuous random variable, where again the target space is equipped with its $w^*$-topology. Via the duality $L_{w^*}^\infty(0,T;\X^1(\T^d)^*)\cong L^1(0,T;\X^1(\T^d))^*$ and the corresponding pairing $\lls\cdot,\cdot\rrs$ the empirical current process is given by 
	$$\lls G,W^N\rrs=\int_0^T\fr{N^{d-1}}\sum_{i=1}^d\sum_{x\in\T_N^d}G^i_t\Big(\frac{x}{N}\Big)W_{x,x+e_i}(\eta_t)\df t,\quad G\in L^1(0,T;\X^1(\T^d)).$$

    \subsubsection{The micro and macro empirical densities}\label{YE} In order to give a closed hydrodynamic equation for the ZRP and in the study of the replacement lemma it will be useful to model the empirical distribution of the ZRP as a ``Young-measure". Since the ZRP takes non-negative values the corresponding empirical ``Young-measures" will be measures on $\T^d\x\RR_+$. In a particle system that takes real values or a particle system with $m$ species, $m\in\NN$, the corresponding empirical measures would be measures on $\T^d\x\RR$ or $\T^d\x\RR^m$ respectively.
    
    A measure $\bs\rho\in\PP(\T^d\x\RR_+)$ with marginal on $\T^d$ equal to the Lebesgue measure is called an (ordinary) Young measure and the space of Young measures is denoted by $\Y(\T^d)$. Via the disintegration theorem~\cite[Theorem 5.3.1]{AGS2000a} to each Young measure $\bs\rho\in\Y(\T^d)$ there corresponds a uniquely determined Lebesgue-a.s.~defined Borel family of probability measures $(\bs\rho^u)_{u\in\T^d}$ and $\bs\rho$ is recovered by the integral $\bs\rho=\int\delta_u\otimes\bs\rho^u\df u$, i.e. 
    $$\int F(u,\lambda)\df\bs\rho(u,\lambda)=\int_{\T^d}\int_{\RR_+}F(u,\lambda)\df\bs\rho^u(\lambda)\df u,\quad\forall F\in B(\T^d\x\RR_+).$$ By a slight abuse of language we will often refer to the elements of the space $\MMM(\T^d\x\RR_+)$ of (signed) Borel measures on $\T^d\x\RR_+$ with finite total variation as Young measures on $\T^d\x\RR_+$.
    
    A Young-measure $\bs\rho\in\Y(\T^d)$ is said to have finite $r$-th moments if $\int\lambda^r\df\bs\rho(u,\lambda)<+\infty$ and the space of all Young-measures with finite $r$-th moment will be denoted by $\Y_r(\T^d)$. We interpret the $1$-st moment as the mass of a Young-measure and for each $r\geq 1$ and $\mathfrak{m}>0$ we will denote by $\Y_{r,\mathfrak{m}}(\T^d)$ the space of all Young-measures $\bs\rho\in\Y_r(\T^d)$ with total mass 
    $$\int_{\T^d\x\RR_+}\lambda\df\bs\rho(u,\lambda)=\int_{\T^d}\int_{\RR_+}\lambda\bs{\rho}^u(\lambda)\df u=\mathfrak{m}.$$

    In order to obtain a closed hydrodynamic equation for condensing ZRPs relying only on the one-block estimate a first idea is to view the empirical distribution of the ZRP as an element of the space $\MMM_1(\T^d\x\RR_+)$ via the micro-empirical density map $\bs\rho^{N,\ell}\colon\MM_N^d\to\MMM_1(\T^d\x\RR_+)$ given by 
    $$\bs\rho^{N,\ell}:=\fr{N^d}\sum_{x\in\T_N^d}\delta_{\frac{x}{N}}\otimes\delta_{\eta^\ell(x)},$$ i.e.~as the process $\bs\rho^{N,\ell}\colon D(0,T;\MM_N^d)\to L_{w^*}^\infty(0,T;\MMM_1(\T^d\x\RR_+))$ given by 
    \begin{equation}\lls F,\bs\rho^{N,\ell}\rrs:=\int_0^T\fr{N^d}\sum_{x\in\T_N^d}F_t\Big(\frac{x}{N},\eta^\ell_t(x)\Big)\df t,\quad F\in L^1(0,T;C_1(\T^d\x\RR_+)).
    \end{equation} 
    
    However, as it turns out, ordinary Young-measures with finite $r$-th moments as described in Section~\ref{OYMSect} are not sufficient for this purpose since they may lose track of the mass in condensed phase in the macroscopic limit, i.e.~as $N\ra+\infty$ and then $\ell\ra+\infty$ To overcome this difficulty and to be able to take into account the mass that is lost by the Young measures we define a notion of generalized Young-measures. This generalized notion is based on a duality result in Section~\ref{OYMSect}, according to which the subspace $\MMM_r(\T^d\x\RR_+)$ of $\MMM(\T^d\x\RR_+)$ of Young-measures with finite $r$-th moment equipped with the norm 
    $$\|\bs\rho\|_{TV,r}:=\int_{\T^d\x\RR_+}(1+\lambda^r)\df|\bs\rho|(u,\lambda),$$
    where $|\bs\rho|$ is the variation of $\bs\rho$, is isometric to the dual of the Banach space $(C_r(\T^d\x\RR_+),\|\cdot\|_{\infty,r})$ of all continuous maps $F\colon\T^d\x\RR_+\to\RR_+$ for which the map 
    $$\T^d\x\RR_+\ni(u,\lambda)\to\frac{F(u,\lambda)}{1+\lambda^r}$$ belongs in the space $C_0(\T^d\x\RR_+)$ of maps that vanish at infinity, where the norm $\|\cdot\|_{\infty,r}$ is defined by \begin{equation}\label{rinftNorm}\|F\|_{\infty,r}:=\sup_{(u,\lambda)\in\T^d\x\RR_+}\frac{|F(u,\lambda)|}{1+\lambda^r}.\end{equation} This duality allows us to consider the space $\MMM_r(\T^d\x\RR_+)$ equipped with the corresponding $w^*$-topology. This $w^*$-topology is not strong enough to assure that no $r$-th moment is lost when taking limits since it can happen that $\bs{\rho}_n\lra\bs{\rho}$ in the $w^*$- topology while 
    $$\int\lambda^r\df\bs{\rho}(\lambda)<\liminf_{n\ra+\infty}\int\lambda^r\df\bs{\rho}_n.$$ 
    
    
    In view of this duality result we denote by $B_rC(\T^d\x\RR_+)$ the space of all $F\in C(\T^d\x\RR_+)$ such that $\|F\|_{\infty,r}<+\infty$ where $\|\cdot\|_{\infty,r}$ is the norm given in~\eqref{rinftNorm} and we define the space
    \begin{equation}\label{barcr}\bbar{C}_r(\T^d\x\RR^d):=\Big\{F\in B_rC(\T^d\x\RR_+)\Bigm|\exists f\in C(\T^d)\mbox{ s.t.~}\Big\|\frac{F(\cdot,\lambda)}{1+\lambda^r}-f(\cdot)\Big\|_\infty\stackrel{\lambda\ra+\infty}\lra 0\Big\}.\end{equation}
    Thus $\bbar{C}_r(\T^d\x\RR_+)$ consists of all $F\in B_rC(\T^d\x\RR_+)$ such that $\|F\|_{\infty,r}<+\infty$ and the recession function 
    \begin{equation}\label{Recession}
    R_rF(u)\equiv RF(u):=\lim_{\lambda\ra+\infty}\frac{F(u,\lambda)}{\lambda^r},\quad F\in\bbar{C}_r(\T^d\x\RR_+)
    \end{equation} is well-defined for all $u\in\T_N^d$, with the limit as $\lambda\ra+\infty$ being uniform over all $u\in\T^d$. As shown in~Proposition~\ref{BarCrBanach} the space $\bbar{C}_r(\T^d\x\RR_+)$ is a closed subspace of $B_rC(\T^d\x\RR_+)$ and thus a Banach space when equipped with the norm $\|\cdot\|_{\infty,r}$. We define the space of \emph{generalized Young-functionals (of order $r$ when $r\neq 1$)} as 
    $$\bbar{\MMM}_r(\T^d\x\RR_+):=\bbar{C}_r(\T^d\x\RR_+)^*$$
    equipped with the dual norm 
    $$\|\bs\pi\|_{TV,r}:=\sup_{\substack{F\in\bbar{C}_r(\T^d\x\RR_+)\\\|F\|_{\infty,r}\leq 1}}\bs\pi(F).$$
    The \emph{recession operator} $R\equiv R_r\colon\bbar{C}_r(\T^d\x\RR_+)\to C(\T^d)$ defined in~\eqref{Recession} is a linear contraction since the limit is assumed uniform over all $u\in\T^d$ and $C_r(\T^d\x\RR_+)=\ker R=R^{-1}(\{0\})$ is a closed subspace of $\bbar{C}_r(\T^d\x\RR_+)$. We will denote by $j\colon C_r(\T^d\x\RR_+)\hookrightarrow \bbar{C}_r(\T^d\x\RR_+)$
    the natural injection. Then the adjoint of $j$ gives a $w^*$-continuous surjective contraction 
    $$j^*\colon\bbar{\MMM}_r(\T^d\x\RR_+)\to\MMM_r(\T^d\x\RR_+).$$
    
     A generalized Young-functional $\bs{\pi}$ is called \emph{non-negative} if it is a positive functional i.e.~if $\bs{\pi}(F)\geq 0$ for all non-negative $F\in\bbar{C}_{r,+}(\T^d\x\RR_+)$. The space of all non-negative generalized Young-functionals on $\T^d\x\RR_+$ is denoted by $\bbar{\MMM}_{r,+}(\T^d\x\RR_+)$. The space of all non-negative Young-functionals $\bs\pi\in\bbar{\MMM}_{r,+}(\T^d\x\RR_+)$ such that the restriction $j^*(\bs{\pi})\equiv\bs{\pi}|_{C_r(\T^d\x\RR_+)}\in C_r(\T^d\x\RR_+)^*$ is via the Riesz isomorphism~Theorem~\ref{YoungRiesz} a probability measure $\bs{\rho}\in\PP_r(\T^d\x\RR_+)$ is denoted by $\bbar{\PP}_r(\T^d\x\RR_+)$. If $\bs\pi\in\bbar{\PP}_r(\T^d\x\RR_+)$ and $m\in\PP\T^d$ then $\bs{\pi}$ is called a \emph{generalized $m$-Young measure} if $U_\sharp\bs\rho=m$ where $U\colon\T^d\x\RR_+\to\T^d$ is the projection on the first coordinate, and the space of all generalized $m$-Young measures is denoted by $T_m\bbar{\PP}_r(\T^d\x\RR_+)$. In the case that $m=\LL_{\T^d}$ is the Lebesgue measure on $\T^d$ then elements of $\bbar{\Y}_r(\T^d):=T_{\LL_{\T^d}}\bbar{\PP}_r(\T^d\x\RR_+)$ are called \emph{generalized Young measures}. Finally, with $\Lambda\in\bbar{C_r}(\T^d\x\RR_+)$ denoting the projection on the second coordinate, we say that a non-negative generalized Young-functional $\bs\pi\in\bbar{\MMM}_{r,+}(\T^d\x\RR_+)$, $r\geq 1$, has \emph{total mass} $\mathfrak{m}>0$ if $\bs{\pi}(\Lambda)=\mathfrak{m}$. The space of all non-negative generalized Young-functionals with total mass $\mathfrak{m}>0$ will be denoted by $\bbar{\MMM}_{r,\mathfrak{m}}(\T^d\x\RR_+)$ and its subspace consisting of generalized Young measures will be denoted as $\bbar{\Y}_{r,\mathfrak{m}}(\T^d)=T_{\LL_{\T^d}}\bbar{\PP}_{r,\mathfrak{m}}(\T^d\x\RR_+)$.
   
  Note that any element $\bs\rho\in\MMM_r(\T^d\x\RR_+)$ can integrate any measurable map $F\colon\T^d\x\RR_+\to\RR$ with at most $r$-th growth at infinity, i.e.~$\|F\|_{\infty,r}<+\infty$ and thus an extension operator $E\colon\MMM_r(\T^d\x\RR_+)\to\bbar{\MMM}_r(\T^d\x\RR_+)$ is defined via 
  $$\ls F,E(\bs\rho)\rs=\int F\df\bs\rho,\quad F\in\bbar{C}_r(\T^d\x\RR_+).$$
   We will regard $\MMM_r(\T^d\x\RR_+)$ as a subspace of $\bbar{\MMM}_r(\T^d\x\RR_+)$ via the extension operator $E$. Also for all $r\geq 0$ and all $\mathfrak{m}>0$ we have $E(\Y_{r}(\T^d))\subs\bbar{\Y}_{r}(\T^d)$ and $E(\Y_{r,\mathfrak{m}}(\T^d))\subs\bbar{\Y}_{r,\mathfrak{m}}(\T^d)$ and thus we will regard $\Y_r(\T^d)$ and $\Y_{r,\mathfrak{m}}(\T^d)$ as subspaces of $\bbar{\Y}_{r}(\T^d)$ and $\bbar{\Y}_{r,\mathfrak{m}}(\T^d)$ respectively via the extension operator $E$.

    As we will see in Section~\ref{GYM}, where we will study generalized Young measures in more detail, any generalized Young-functional $\bs\pi\in\bbar{\MMM}_r(\T^d\x\RR_+)$ can be represented uniquely by a pair $(\bs\rho_{\bs\pi},\rho^\perp_{\bs\pi})\in\MMM_r(\T^d\x\RR_+)\x\MMM(\T^d)$. The action of $\bs\pi$ on functions $F\in\bbar{C}_r(\T^d\x\RR_+)$ is recovered from the pair $(\bs\rho_{\bs\pi},\rho^\perp_{\bs\pi})$ by the formula 
   $$\ls F,\bs\pi\rs=\int_{\T^d\x\RR_+}F(u,\lambda)\df\bs\rho_{\bs\pi}(u,\lambda)+\int_{\T^d} RF(u)\df\rho^\perp_{\bs\pi}(u)$$
   according to which $\bs\pi=E(\bs\rho_{\bs\pi})+R^*(\rho^\perp_{\bs\pi})$ with $R^*\colon \MMM(\T^d)\to\bbar{\MMM}_r(\T^d\x\RR_+)$ being the adjoint of the recession operator $R$. Thus as we will see any Young measure $\bs\pi$ can be written uniquely as the sum $\bs\pi=\widehat{\bs\pi}+\bs\pi^\perp$ where $\widehat{\bs\pi}:=E(\bs\rho_{\bs\pi})\in\bbar{\MMM}_r(\T^d\x\RR_+)$ is the extension of a uniquely determined Young measure $\bs\rho_{\bs\pi}$ and $\bs\pi^\perp:=R^*(\rho^\perp_{\bs\pi})\in\bbar{\MMM}_r(\T^d\x\RR_+)$ is a generalized Young measure that acts on maps $F\in\bbar{C}_r(\T^d\x\RR_+)$ only through their recession function $RF$ via the integration $\int RF\df\rho^\perp_{\bs\pi}$ for a uniquely determined Borel measure $\rho^\perp_{\bs\pi}\in\MMM(\T^d)$. We will refer to $\widehat{\bs\pi}$ as the regular part of $\bs\pi$ and to $\bs\rho_{\bs\pi}$ as the (ordinary) Young measure $\bs\rho_{\bs\pi}$ representing the regular part of $\bs\pi$, and we will refer to $\bs\pi^\perp$ as the singular part of $\bs\pi$ and to $\rho_{\bs\pi}^\perp$ as the measure representing the singular part of $\bs\pi$. We will sometimes write $\bs\pi=(\bs\rho,\rho^\perp)$ to denote the fact that the generalized Young-functional $\bs\pi\in\bbar{\MMM}_r(\T^d\x\RR_+)$ is represented by the pair $(\bs\rho,\rho^\perp)\in\MMM_r(\T^d\x\RR_+)\x\MMM(\T^d)$ via the relation $\bs\pi=R(\bs\rho)+R^*(\rho^\perp)$.
   
   This representation is also valid on the level of path-measures. Indeed as we will see in Section~\ref{GYM} the extension operator $E$ is bounded and $w^*$-Baire measurable (see the end of Section~\ref{Appendix1} in the appendix for the notion of $w^*$-Baire measurability considered here) and therefore induces by Proposition~\ref{InducedOpProp} a $w^*$-Baire measurable (and thus also $w^*$-measurable) operator, still denoted by $E$, on the corresponding $L_{w^*}^\infty$-spaces. Also the recession operator $R\colon\bbar{C}_r(\T^d\x\RR_+)\to C(\T^d)$ induces an operator, still denoted by $R$, on the corresponding $L^1$-Bochner spaces and its adjoint $R^*\colon L_{w^*}^\infty(0,T;\MMM(\T^d))\to L_{w^*}^\infty(0,T;\bbar{\MMM}_r(\T^d\x\RR_+))$ is $w^*$-continuous. Then any generalized Young path-measure $\bs\pi\in L_{w^*}^\infty(0,T;\bbar{\MMM}_1(\T^d\x\RR_+))$ is represented uniquely by a pair $(\bs\rho_{\bs\pi},\rho^\perp_{\bs\pi})\in L_{w^*}^\infty(0,T;\MMM_1(\T^d\x\RR_+))\x L_{w^*}^\infty(0,T;\MMM(\T^d))$ via the decomposition $\bs\pi=E(\bs\rho_{\bs\pi})+R^*(\rho^\perp_{\bs\pi})$ and this decomposition is $w^*$-Baire measurable. The action of $\bs\pi$ on test functions $F\in L^1(0,T;\bbar{C}_1(\T^d\x\RR_+))$ is recovered from the pair $(\bs\rho_{\bs\pi},\rho^\perp_{\bs\pi})$ via 
   $$\lls F,\bs\pi\rrs=\int_0^T\int_{\T^d\x\RR_+}F_t(u,\lambda)\df\bs\rho_{\bs\pi_t}(u,\lambda)\df t+\int_0^T\int_{\T^d}RF_t(u)\df\rho^\perp_{\bs\pi_t}(u)\df t.$$
   
   Using the notion of generalized Young-functionals we define \emph{the micro-empirical density map $\bs\pi^{N,\ell}\colon\MM_N^d\to\bbar{\MMM}_1(\T^d\x\RR_+)$, $N\in\NN$, $\ell\in\ZZ_+$ of the ZRP} by the formula 
   \begin{equation}\label{Memp}\ls F,\bs\pi^{N,\ell}\rs=\ls F,E(\bs\rho^{N,\ell})\rs=\fr{N^d}\sum_{x\in\T_N^d}F\Big(\frac{x}{N},\eta^\ell(x)\Big),\quad F\in\bbar{C}_1(\T^d\x\RR_+).
\end{equation}
   Since $\MM_N^d$ has the discrete topology the map $\bs\pi^{N,\ell}$ is continuous. Thus the induced map 
   \[\bs\pi^{N,\ell}\colon D(0,T;\MM_N^d)\to D(0,T;\bbar{\MMM}_1(\T^d\x\RR_+))\]
   is continuous where here $\bbar{\MMM}_1(\T^d\x\RR_+)$ is considered with the dual norm $\|\cdot\|_{TV;1}$. Thus by composing this map with the continuous injection 
   \[D(0,T;\bbar{\MMM}_1(\T^d\x\RR_+))\to L^\infty_{w^*}(0,T;\bbar{\MMM}_1(\T^d\x\RR_+))\]
   given by Proposition~\ref{SkorohodToLEmbed}, we obtain that the \emph{the micro-empirical process} 
   \[\bs\pi^{N,\ell}\colon D(0,T;\MM_N^d)\to L_{w^*}^\infty(0,T;\bbar{\PP}_1(\T^d\x\RR_+))\subs L_{w^*}^\infty(0,T;\bbar{\MMM}_1(\T^d\x\RR_+))\] 
   is continuous and thus a random variable. We define the \emph{macro-empirical process $\bs\pi^{N,\ee}$, $N\in\NN$, $\ee>0$ of the ZRP} by
   \[\bs\pi^{N,\ee}:=\bs\pi^{N,[N\ee]}\colon D(0,T;\MM_N^d)\to L_{w^*}^\infty(0,T;\bbar{\PP}_1(\T^d\x\RR_+)).\] 
   
   Finally, for any $M>0$ we will also consider the \emph{$M$-modified micro-empirical density process} $\bs\pi^{N,\ell;M}\colon D(0,T;\MM_N^d)\to L_{w^*}^\infty(0,T;\bbar{\PP}_1(\T^d\x\RR_+))$, $\ell\in\ZZ_+$, defined by 
   \begin{equation}\label{MModMicrEmpDens}
   \lls F,\bs\pi^{N,\ell;M}\rrs=\int_0^T\fr{N^d}\sum_{x\in\T_N^d}\Big\{F_t\Big(\frac{x}{N},\eta_t^\ell(x)\mn M\Big)+RF_t\Big(\frac{x}{N}\Big)(\eta^\ell_t(x)-M)^+\Big\}\df t
   \end{equation}
   in duality with respect to test functions $F\in L^1(0,T;\bbar{C}_1(\T^d\x\RR_+))$ and the \emph{$M$-modified macro-empirical density process} $\bs\pi^{N,\ee;M}:=\bs\pi^{N,[N\ee];M}$, $\ee>0$. By the same reasoning as all the other $L_{w^*}^\infty$-valued process defined in this section these are also well defined continuous random variables.

  Generalized Young measures are related to Borel measures via the barycentric projection map $B\colon\bbar{\MMM}_1(\T^d\x\RR_+)\to\MMM(\T^d)$ that is defined by 
  $$B(\bs\pi)(f)=\bs\pi\big(\Lambda f(U)\big),\quad f\in C(\T^d)$$
  where $U\colon\T^d\x\RR_+\to\T^d$ and $\Lambda\colon\T^d\x\RR_+\to\RR_+$ are the natural projections. Since $\Lambda f(U)\in\bbar{C}_1(\T^d\x\RR_+)$ for all $f\in C(\T^d)$ it is obvious $B$ is $w^*$-continuous. Thus by Proposition~\ref{AdjointInduced} it induces a $w^*$-continuous operator $B\colon L_{w^*}^\infty(0,T;\bbar{\MMM}_1(\T^d\x\RR_+))\to L_{w^*}^\infty(0,T;\MMM(\T^d))$ on the corresponding spaces of $L_{w^*}^\infty$-path-measures denoted by the same symbol $B$. More generally for any map $\Psi\in\bbar{C}_1(\RR_+)$, where $\bbar{C}_1(\RR_+)$ is the space of all continuous maps $\Psi\colon\RR_+\to\RR$ such that the limit $\Psi'(\infty):=\lim_{\lambda\ra+\infty}\frac{\Psi(\lambda)}{\lambda}\in\RR$ exists, we consider the $\Psi$-projection $B_\Psi\colon\bbar{\MMM}_1(\T^d\x\RR_+)\to\MMM(\T^d)$ given by 
  \begin{equation}\label{PsiProj}B_\Psi(\bs\pi)(f)=\bs\pi\big(\Psi(\Lambda)f(U)\big),\quad\bs\pi\in\MMM_1(\T^d\x\RR_+).\end{equation} Since $B_\Psi$ is $w^*$-continuous it also induces a $w^*$-continuous operator on the corresponding $L_{w^*}^\infty$-spaces. As we will see in section~\ref{Barycentric} if $\bs\pi\in T_m\bbar{\MMM}_{1,+}(\T^d\x\RR_+)$ and $\bs\pi=E(\bs\rho)+R^*(\rho^\perp)$ for some $\bs\rho\in\MMM_{1,+}(\T^d\x\RR_+)$, $\rho^\perp\in\MMM_+(\T^d)$ then $B_\Psi(\bs\pi)=\Psi(\bs\rho)\df m+\Psi'(\infty)\rho^\perp$ where $\Psi(\bs\rho)$ is the $m$-a.e.~defined map
  \begin{equation}\label{YoungComp}\Psi(\bs\rho)(u)=\int\Psi(\lambda)\df\bs\rho^u(\lambda).
  \end{equation}
with $(\bs\rho^u)_{u\in\T^d}$ being the $m$-a.e.~defined disintegration of $\bs\rho$ with respect to its first marginal $m=U\sharp\bs\rho$, i.e.~$\bs\rho=\int_{\T^d}\delta_u\otimes\bs\rho^u\df m(u)$, whose existence is guaranteed by~\cite[Theorem 5.3.1]{AGS2000a}. In particular if $\Psi$ is sublinear, i.e.~$\Psi_\infty=0$ then $B_\Psi(\bs\pi)$ depends only on the regular part $\widehat{\bs\pi}=E(\bs\rho)$ of $\bs\pi$ and is a measure absolutely continuous with respect to $m$.

\subsection{Hydrodynamic limits} In order to obtain the hydrodynamic limit of condensing ZRPs one has to prove when starting from a sequence of initial distributions $\mu_0^N\in\PP\MM_N^d$  \emph{associated to some macroscopic profile} $\mu_0\in\MMM_+(\T^d)$, i.e.
\begin{equation}\lim_{N\ra+\infty}\mu_0^N\big\{\big|\ls G,\pi^N-\mu_0\rs\big|>\delta\big\}=0,\quad\forall G\in C(\T^d),\;\delta>0,\end{equation}
 one would ideally like to prove that the laws $\big(\pi^N_t)_{0\leq t\leq T}\big)_\sharp P^N$ of the empirical process of the diffusively rescaled ZRP starting from $\{\mu_0^N\}$ converge weakly in an appropriate topology to a Dirac measure $\delta_\pi$ supported on the unique solution $\pi\in[0,T]\to\MMM_+(\T^d)$ (in some appropriate sense) of the saturated filtration equation 
 \begin{equation}\label{HL}
 \begin{cases}\pd_t\pi=\Delta\Phi(\pi)\\
 \pi_0=\mu_0\end{cases}
 \end{equation}
 where $\Phi\equiv\bbar{\Phi}$ is the extended mean jump rate function of the ZRP. In general the standard approach towards this aim is to prove that the laws $Q^N:=((\pi_t^N)_{0\leq t\leq T})_\sharp P^N\in\PP\X$ of the ZRP are relatively compact with respect to an appropriate path space $\X$ and then show that any subsequential limit $Q$ of $\{Q^N\}_{N\in\NN}$ is supported by a set of curves $\pi\in\X$ satisfying an evolutionary equation, e.g.~\eqref{HL}, which is then called \emph{the hydrodynamic equation}. If then the hydrodynamic equation satisfies uniqueness of solutions it follows that any subsequential limit point of the sequence $Q^N$ starting from $\{\mu_0^N\}$ is supported by the unique solution $\pi$ of~\eqref{HL} and thus the whole sequence converges to $\delta_\pi\in\PP\X$.
 
Giving a precise meaning to the notion of solutions of the hydrodynamic equation and choosing the path space $\X$ of solutions in the program above is part of the problem. In the particular case of condensing ZRPs equation~\eqref{HL} is ill-behaved. Namely, for weakly condensing ZRPs it is not uniformly parabolic for unbounded profiles since $$\lim_{\rho\ra+\infty}\Phi'(\rho)=0$$ while for strongly condensing ZRPs $\bbar{\Phi}'(\rho)=0$ on $(\rho_c,+\infty)$ with $\Phi$ being possibly non-differentiable at $\rho_c$. For example in the Evans model the extended mean jump rate $\bbar{\Phi}$ is differentiable for $b\in[0,3]$ while $\bbar{\Phi}_-'(\rho_c)>0$ for $b>3$.

In this article we will also examine the set of subsequential limit points of the sequence of laws of the joint empirical processes 
$$(\pi^N,W^N,\s^N)\colon D(0,T;\MM_N^d)\to D(0,T;\MMM_+(\T^d))\x L_{w^*}^\infty(0,T;\X^1(\T^d)^*)\x L^\infty_{w^*}(0,T;\MMM_+(\T^d))$$ and show that they are concentrated on trajectories $(\pi,W,\s)$ satisfying the continuity equation $\pd_t+\dv W=0$ with $W_t=-\nabla\s_t$ in $(0,T)\x\T^d$ in the sense of distributions. We will also examine the set of subsequential limit points of the sequence of laws of the micro-empirical distribution $\bs\pi^{N,\ell}$ on the space $L_{w^*}^\infty(0,T;\bbar{\MMM}_{1,+}(\T^d\x\RR_+))$ of generalized Young path-measures and show that any subsequential limit point is concentrated on solutions of the closed form equation~\eqref{HL} with respect to an appropriate sense in terms of generalized Young measures.

\subsubsection{The continuity equation}
 Let $(\pi_t)_{t\geq 0}\subs\MMM_+(\T^d)$ and $(W_t)_{t\geq 0}\subs\MMM(\T^d;\RR^d)$ be Borel families of measures. By~\cite[Lemma 4.1]{DNS2009a} if $\pi$ solves the continuity equation $\pd_t\pi=\dv W_t$ in the sense of distributions and $\bs\pi$, $W$ satisfy 
$$\int_0^T\pi_t(\T^d)\df t<+\infty\quad\mbox{and}\quad\int_0^T|W_t|(\T^d)\df t<+\infty$$ then there exists a weakly continuous curve $\wt{\pi}\colon[0,T]\to\MMM_+(\T^d)$ such that $\pi_t=\wt{\pi}_t$ for almost all $t\in[0,T]$ and for this curve $\wt{\pi}$ for all $G\in C^1([0,T]\x\T^d)$
\begin{equation}\label{ContEqContReprForm}\int G_t\df\wt{\pi}_t-\int G_s\df\wt{\pi}_s=\int_s^t\int_{\T^d}\pd_rG_r\df\wt{\pi}_r\df r+\int_s^t\int_{\T^d}\ls\nabla G_r,\df W_r\rs\df r.\end{equation}
Thus given $\mu_0\in\MMM_+(\T^d)$ one says that the Borel families $\pi=(\pi_t)_{0\leq t\leq T}\subs\MMM_+(\T^d)$ and $W=(W_t)_{0\leq t\leq T}\subs\MMM(\T^d;\RR^d)$ satisfy the initial value problem \begin{equation}\label{CEInitVal}
\begin{cases}\pd_t\pi_t+\dv W_t=0\\
\pi_0=\mu_0\end{cases}
\end{equation}  
if $\pi=(\pi_t)_{0\leq t\leq T}$ is a weakly continuous curve satisfying $\pd_t\pi+\dv W=0$ in the sense of distributions and $\pi_0=\mu_0$.

For our purposes the current $W=(W_t)_{0\leq t\leq T}$ in the continuity equation is has to be modelled as a $w^*$-measurable curve in $\X^1(\T^d)^*$. We say that \emph{a density/current pair $(\pi,W)\in L_{w^*}^\infty(0,T;\MMM_+(\T^d))\x L_{w^*}^\infty(0,T;\X^1(\T^d)^*)$ satisfies the continuity equation} $\pd_t\pi+\dv W=0$ if 
\begin{equation}\label{CEDistrDef}
\int_0^T\int_{\T^d}\pd_tG_t\df t+\int_0^T\ls\nabla G_t,W_t\rs\df t=0,\quad\forall G\in C_c^{1,2}((0,T)\x\T^d).
\end{equation}
Here we require $G$ to be twice continuously differentiable in space so that the curve $(\nabla G_t)_{0\leq t\leq T}$ defines an element of $L^\infty(0,T;\X^1(\T^d))$. As we will see in Proposition~\ref{ContEqBasicDensRegularity} the existence of a $w^*$-continuous representative $\wt{\pi}\in C_{w^*}(0,T;\MMM_+(\T^d))$ for which~\eqref{ContEqContReprForm} holds is also valid in the case that the current $W=(W_t)_{0\leq t\leq T}$ is more generally modelled as an element of the space $L^1_{w^*}(0,T;\X^1(\T^d)^*)$, but for all $G\in C^{1,2}([0,T]\x\T^d)$ in this case. Thus we say that $(\pi,W)\in L_{w^*}^1(0,T;\MMM_+(\T^d))\x L^1_{w^*}(0,T;\X^1(\T^d)^*)$ satisfies the initial value problem~\eqref{CEInitVal} if the continuous representative $\wt{\pi}$ in the a.s.~equality class of $\pi$ satisfies $\wt{\pi}_0=\mu_0$. In particular if $\pi\in L_{w^*}^1(0,T;\MMM_+(\T^d))$ solves the initial value problem~\eqref{CEInitVal} then by applying~\eqref{ContEqContReprForm} for the constant map $G\equiv1$ it follows that $\wt{\pi}_t(\T^d)=\wt{\pi}_0(\T^d)=\mu_0(\T^d)$ for all $0\leq t\leq T$.

We will say that a triple $(\pi,W,\s)\in L_{w^*}^\infty(0,T;\MMM_+(\T^d))\x L_{w^*}^\infty(0,T;\X_1(\T^d)^*)\x L_{w^*}^\infty(0,T;\MMM_+(\T^d))$ satisfies the continuity equation 
\begin{equation}\label{ContEqDivGrad}\begin{cases}\pd_t\pi=-\dv W\\
W=-\nabla\s\end{cases}\quad\mbox{ in }(0,T)\x\T^d
\end{equation}
in the sense of distributions if $\pd_t\pi+\dv W=0$ holds in the sense of~\eqref{CEDistrDef} and 
\begin{equation}\label{CurrentGrad}\int_0^T\ls F,W_t\rs\df t=\int_0^T\ls\dv F_t,\s_t\rs\df t,\quad \forall\;F\in L^1(0,T;\X^1(\T^d)).
\end{equation}
This is stronger than requiring that the equation $\pd_t\pi=-\dv W=\Delta\s$ holds in the sense if distributions since this would be equivalent to requiring that~\eqref{CurrentGrad} holds for maps $F\in L^1(0,T;\X^1(\T^d))$ that are spatial gradients of $C^2$-functions i.e.~$F_t=\nabla f_t$ for some $f_t\in C^2(\T^d)$ for almost  all $t\in[0,T]$.
\subsubsection{Closed form equations}
A rigorous interpretation of equation~\eqref{HL} allowing for measure-valued solutions has been given in \cite{Fornaro2012a} in dimension $d=1$ by interpreting equation~\eqref{HL} as a gradient flow in the quadratic Wasserstein space $\MMM_{+,\mathfrak{m}}(\T^d)$ of measures with a fixed total mass $\mathfrak{m}=\mu_0(\T^d)$. The gradient flow formulation result in~\cite{Fornaro2012a} allows then the writers to obtain uniqueness and existence of solutions in the sense of distributions to problem~\eqref{HL} by interpreting it as the problem 
\begin{equation}\label{HL2}
\begin{cases}\pd_t\pi=\Delta\Phi(\rho),\\
\pi_0=\mu_0\end{cases}
\end{equation}
where $\pi=\rho\df\LL_{\T^d}+\pi^\perp$, $\pi^\perp\perp\LL_{\T^d}$ is the Radon-Nikodym decomposition of $\pi$ with respect to Lebesgue measure $\LL_{\T^d}$ on the torus. Uniqueness of the weak solutions of problem~\eqref{HL2} is obtained in~\cite{Fornaro2012a} in the class of weakly continuous measure-valued curves with finite kinetic energy that take values in the space ${C}\MMM_+(\T^d)$ of continuous measures defined as $$\mathcal{C}\MMM_+(\T^d)=\big\{\pi=\rho\df u+\pi^\perp\bigm|\rho\in C(\T^d;[0,\rho_c]),\;\LL_{\T^d}\{\rho=\rho_c\}=0,\;\pi^\perp\{\rho<\rho_c\}=0\big\}.$$
Namely for any $\mu_0\in\MMM_+(\T^d)$ there exists a unique weak solution $\pi=\rho\df u+\pi^\perp\colon\RR_+\to\MMM_+(\T^d)$ to problem~\eqref{HL2} in the sense that $\Phi(\rho)\in W^{1,1}(\T^d)$ and 
$$\int_0^\infty\int\pd_tf_t\df\pi_t\df t=\int_0^\infty\int\ls\nabla f_t(u),\nabla\Phi(\rho_t(u))\rs\df u\df t,\quad\forall\;f\in C_c^1((0,T)\x\T^d)$$
such that 
\begin{itemize}\item[(a)] The curve $(\pi_t)_{t\geq 0}$ is weakly continuous in $\MMM_+(\T^d)$ and $\pi_0=\mu_0$,
	\item[(b)] $\pi_t\in\mathcal{C}\MMM_+(\T^d)$ for Lebesgue almost all $t>0$, and 
	\item[(c)] $\int_{T_0}^{T_1}\mathcal{J}(\pi_t)\df t<+\infty$ for all $0<T_0<T_1<+\infty$
\end{itemize}
where $\mathcal{J}\colon\MMM_+(\T^d)\to[0,+\infty]$ is the generalized Fisher dissipation functional  \[\mathcal{J}(\pi)=\begin{cases}\int_{\{\rho>0\}}\frac{|\nabla\Phi(\rho(u))|^2}{\rho(u)}\df u&\mbox{if }\pi=\rho\df u+\pi^\perp\in \mathcal{C}\MMM_+(\T^d)\mbox{ and }\Phi\circ\rho\in W^{1,1}(\T^d)\\
+\infty&\mbox{otherwise}\end{cases}\]
The assumptions in~\cite{Fornaro2012a} require that $\Phi$ is $C^1$ with $\Phi'(\rho)>0$ for all $\rho\in\RR_+$, which restricts the applicability of this gradient flow formulation to weakly condensing ZRPs.

In this article we consider a weaker notion of solutions via the notion of generalized Young measures. First for any (ordinary) Lebesgue-Young path-measure $\bs\rho\in T_{\LL_{\T^d}}\PP_1(\T^d\x\RR_+)$ and any map $\Phi\in C_1(\RR_+)$ we define the composition $\Phi(\bs\rho)\colon\T^d\to\RR_+$ by the formula
\begin{equation}\label{YoungCompose}
\Phi(\bs\rho)(u):=\int\Phi(\lambda)\df\bs\rho^u(\lambda),\quad\mbox{for }\LL_{\T^d}\mbox{-almost all }u\in\T^d,
\end{equation}
where $(\bs\rho^u)_{u\in\T^d}$ is the $\LL_{\T^d}$-a.s.~defined disintegration of $\bs\rho$ with respect to its first marginal $U_\sharp \bs\rho=\LL_{\T^d}$. Given a Lipschitz map $\Phi$ such that \begin{equation}\label{PhiSlope}\Phi'(\infty):=\lim_{\lambda\ra+\infty}\frac{\Phi(\lambda)}{\lambda}=0\end{equation} we say that a map $\bs\pi=(\bs\rho,\mu)\in L_{w^*}^\infty(0,T;\bbar{\Y}_1(\T^d)$ is a \emph{generalized Young-measure valued weak solution of the problem} 
\begin{equation}\label{YMHL}\pd_t\bs\pi=\Delta\Phi(\bs\pi)\end{equation}
if 
\begin{equation}\label{YMVZRPPDE}
\begin{cases}
\Phi(\bs\rho_t)\in W^{1,1}(\T^d)\mbox{ for almost all }t\in[0,T],\\
\int_0^T\bs\pi_t\big(\Lambda\pd_tf(U)\big)\df t=\int_0^T\ls\nabla f_t(u),\nabla\Phi(\bs\rho_t)(u)\rs\df u\df t,\quad\forall\;f\in C_c^1((0,T)\x\T^d)
\end{cases}.\end{equation}
We note that by~\eqref{PhiSlope} any weak solution $\bs\pi=(\bs\rho,\mu)\in L_{w^*}^\infty(0,T;\bbar{\Y}_1(\T^d))$ to problem~\eqref{YMHL} is also a \emph{mild solution} in the sense that
\begin{equation}\label{Mild}
\int_0^T\bs\pi_t\big(\Lambda\pd_tf(U)+\Phi(\Lambda)\Delta f(U)\big)\df t=0,\quad\forall\;f\in C_c^{1,2}((0,T)\x\T^d).
\end{equation}

\subsubsection{Subsequential limit sets}\label{SubLims} For any sequence $\{A_n\}_{n=1}^\infty$ of subsets $A_n\subs M$ of a submetrizable topological space $M$ and any sequence $\{q_n\}_{n=1}^\infty\subs M$ of points we set 
\[\Lim_{n\ra+\infty}A_n=\bigcap_{n=1}^\infty\bbar{\bigcup_{k=n}^\infty A_k},\quad\Lim_{n\ra+\infty}q_n:=\Lim_{n\ra+\infty}\{q_n\}.\]
If the union $\bigcup_{n\in\NN}A_n$ is relatively compact in $M$ then $\Lim_{n\ra+\infty}A_n$ is non-empty and consists of all subsequential limits of points of the sets $A_n$, i.e. 
\[\emptyset\neq \Lim_{n\ra+\infty}A_n=\Big\{q\in M\Bigm|\exists\;q_{k_n}\in A_{k_n},\;k_n<k_{n+1},\;n\in\NN:\;\lim_{n\ra+\infty}q_{k_n}=q\Big\}\]
and $\Lim_{n\ra+\infty}A_n$ is compact.

In exhibiting the continuity equation as a hydrodynamic limit of the laws of the joint process $(\pi^N,W^N,\s^N)$ we will show that the sequence 
$$Q^N:=(\pi^N,W^N,\s^N)_\sharp P^N\in \PP\big(L_{w^*}^\infty(0,T;\MMM_+(\T^d))\x L_{w^*}^\infty(0,T;\X^1(\T^d)^*)\x L_{w^*}^\infty(0,T;\MMM_+(\T^d))\big)$$
is relatively compact and that any limit point $Q\in\Lim_{N\ra+\infty}Q^N$ is concentrated on a measurable set of trajectories $(\pi,W,\s)$ that satisfy equation~\eqref{ContEqDivGrad}.

In proving the closed-form generalized Young-measure valued equation~\eqref{YMHL} in the hydrodynamic limit we will first show that the double sequence of laws 
\[\bs{Q}^{N,\ell}:=\bs\pi^{N,\ell}_\sharp P^N\equiv\big((\bs\pi^{N,\ell}_t)_{0\leq t\leq T}\big)_\sharp P^N\subs\PP L_{w^*}^\infty(0,T;\bbar{\MMM}_{1,+}(\T^d\x\RR_+))\]
is relatively compact and then that any subsequential limit 
\[\bs{Q}\in\Lim_{\ell\ra+\infty}\Big(\Lim_{N\ra+\infty}\bs{Q}^{N,\ell}\Big)=:\Lim_{\ell,N\uparrow\infty}\bs{Q}^{N,\ell}\] of $\{\bs{Q}^{N,\ell}\}$ as $N\ra+\infty$ and then $\ell\ra+\infty$ is concentrated on generalized Young-measure valued weak solutions of~\eqref{YMHL}. Equation~\eqref{YMHL} in terms of generalized Young measures is the first closed form equation given for condensing ZRPs and it relies only on the one-block estimate and not the full replacement lemma.

A property that will be used often is the following. If $\{Q^N\}\subs\PP M$ is a relatively compact sequence of probability measures on the completely regular submetrizable space $M$ and $\mathcal{Q}^\infty:=\Lim_{N\ra+\infty}Q^N$ then for all $f\in BC(M)$ $$\limsup_{N\ra+\infty}\int f\df Q^N=\max_{Q^\infty\in\mathcal{Q}^\infty}\int f\df Q^\infty.$$
A similar property holds for multi-parametric families of probability measures. For example, if $\{\bs{Q}^{N,\ell}\}_{N,\ell}\subs\PP M$ is a relatively compact double sequence of probability measures and $\bs{\mathcal{Q}}^{\infty,\infty}:=\Lim_{\ell,N\uparrow\infty}\bs{Q}^{N,\ell}$ then for all $f\in BC(M)$ $$\limsup_{\ell,N\uparrow\infty}\int f\df\bs{Q}^{N,\ell}=\max_{\bs{Q}^{\infty,\infty}\in\bs{\mathcal{Q}}^{\infty,\infty}}\int f\df\bs{Q}^{\infty,\infty}.$$ 

We note finally that if $\{\mathcal{Q}^N\}_{N\in\NN}$ is a sequence of families $\mathcal{Q}^N\subs\PP M$ of probability measures in a completely regular submetrizable space $M$ and $f\colon M\to N$ is a continuous map from $M$ to the completely regular submetrizable space $N$ then 
\begin{equation}\label{SetPushForwContMap}
f_\sharp\Big(\Lim_{N\ra+\infty}\mathcal{Q}^N\Big)=\Lim_{N\ra+\infty}f_\sharp\mathcal{Q}^N
\end{equation}
where for any  family $\mathcal{Q}\subs\PP M$ we set $f_\sharp\mathcal{Q}=\big\{f_\sharp Q\bigm|Q\in\mathcal{Q}\big\}$.
\subsubsection{Assumptions on the initial distributions}

A sequence $\{\mu_0^N\in\PP\MM_N^d\}$ of initial distributions satisfies \emph{the $O(N^d)$-entropy assumption} if there exist constant $\rho_*\in (0,\rho_c)$ such that 
\begin{eqnarray}\label{LeqThanONdEntr}\qquad
C(\rho_*):=\sup_{N\in\NN}\fr{N^d}H(\mu_0^N|\nu^N_{\rho_*})<+\infty.
\end{eqnarray}
Since $\nu_{\rho_*}^1$ has some exponential moments if $\rho_*<\rho_c$ it follows by the relative entropy inequality~\eqref{RelEntrIneq} that if~\eqref{LeqThanONdEntr} holds for some $\rho_*\in(0,\rho_c)$ then it also holds for all $\rho\in(0,\rho_c)\cap\RR$.

We will note here that in condensing ZRPs the $O(N^d)$-entropy assumption does not prohibit us from starting the ZRP from a sequence of initial states $\mu_0^N\in\PP\MM_N^d$ associated to a macroscopic profile $\mu_0\in\MMM_+(\T^d)$ having a condensate at a point $u\in\T^d$. For example consider the measure $\mu_0=\rho_0\df\LL_{\T^d}+\alpha\delta_u\in\MMM_+(\T^d)$ where $\rho_0\colon\T^d\to\RR$ is a measurable bounded and a.e.~continuous continuous function and $\alpha\geq0$. In the case that there is no condensate, i.e.~$a=0$ then the sequence $\{\nu_{\rho_0(\cdot)}^N\}_{N\in\NN}$ of product measures with slowly varying profile $\rho_0\in B(\T^d)$ satisfies the $O(N^d)$-entropy assumption since
	\begin{equation}\label{MicroscopEntropGoToMarcoscopEntrBeta}
	\lim_{N\ra+\infty}\fr{N^d}H(\nu_{\rho_0(\cdot)}^N|\nu^N_{\rho_*})=\int\Lambda_{\rho_*}^*\big(\rho_0(u)\mn\rho_c\big)\df u<+\infty
	\end{equation}
	for all $\rho_*\in(0,\rho_c)$. Here $\Lambda_{\rho_*}^*\colon\RR\to\RR_+$ is the Legendre transform of $\Lambda_{\rho_*}$ given by
	\begin{equation}\label{RateFunc}\Lambda^*_{\rho_*}(\rho)=\begin{cases}
	\rho\log\frac{\Phi(\rho\mn\rho_c)}{\Phi(\rho_*)}-\log\frac{Z(\Phi(\rho\mn\rho_c))}{Z(\Phi(\rho_*))},\quad&\rho\geq 0\\
	+\infty,\quad&\rho<0\end{cases}.\end{equation} Indeed, 
and $H(\nu_{\rho}^1|\nu_{\rho_*}^1)=\Lambda_{\rho_*}^*(\rho)$ for all $\rho\in[0,\rho_c]\cap\RR$. Therefore
$$\fr{N^d}H(\nu_{\rho_0(\cdot)}^N|\nu_{\rho_*}^N)=
\fr{N^d}\sum_{x\in\T_N^d}H(\nu_{\rho_0(x/N)\mn\rho_c}^1|\nu_{\rho_*}^1)=
\int_{\T^d}\Lambda_{\rho_*}^*\Big(\rho_0\Big(\frac{[Nu]}{N}\Big)\mn\rho_c\Big)\df u.$$
The function $\Lambda_{\nu_{\rho_*}}^*$ is always finite and smooth on all of $\RR_+$, and therefore since we assume the profile $\rho_0$ to be bounded and almost surely continuous, the required limit in~\eqref{MicroscopEntropGoToMarcoscopEntrBeta} follows by the bounded convergence theorem.

According to he following example the $O(N^d)$-entropy assumption is satisfied even by initial distributions that can have a condensate at some macroscopic point $u\in\T^d$. 

\begin{example} Let $\{\nu_{\rho_0(\cdot);u,\alpha}^N\}_{N\in\NN}$ be the sequence of product measures with slowly varying parameter associated to some bounded and a.s$.$ continuous profile $\rho_0\in B(\T^d)$ and a Dirac mass $\alpha>0$ at $x\in\T^d$, i.e$.$ 
	$$\nu_{\rho_0;u,\alpha}^N=\delta_{[\alpha N^d]}\otimes\bigotimes_{x\in\T_N^d\sm\{[Nu]\}}\nu_{\rho_0(\frac{x}{N})}=:\delta_{[\alpha N^d]}\otimes\nu_{\rho_0(\cdot)}^{N\sm u}\in\PP(\ZZ_+\x\ZZ_+^{\T_N^d\sm\{[Nu]\}})\cong\PP\MM_N^d,$$
	Then $\{\nu_{\rho_0;u,\alpha}^N\}_{N\in\NN}$ is associated to the measure $\mu_0=\alpha\delta_u+\rho_0\df\LL_{\T^d}\in\MMM_+(\T^d)$ and 
	\begin{equation}\label{MicroscopEntropGoToMarcoscopEntrBetaDiracCondens} 
	\lim_{N\ra+\infty}\fr{N^d}H(\nu_{\rho_0(\cdot);u,\alpha}^N|\nu_{\rho_*}^N)=
	\int_{\T^d}\Lambda_{\nu_{\rho_*}}^*\big(\rho_0(u)\mn\rho_c\big)du+
	\alpha\log\frac{\phi_c}{\Phi(\rho_*)}
	\end{equation}
	for all $\rho_*\in(0,\rho_c)$. In particular, whenever $\phi_c<+\infty$ the sequence $\{\nu_{\rho_0(\cdot);u,\alpha}^N\}$ is associated to the measure $\mu_0$ which has a condensate of mass $\alpha>0$ at $u\in\T^d$ and satisfies the $O(N^d)$-entropy assumption.
\end{example}\textbf{Proof} For all $N\in\NN$ we have that 
\begin{equation}\label{kleinmein}
H(\nu_{\rho_0(\cdot);u,\alpha}^N|\nu_{\rho_*}^N)=H(\delta_{[\alpha N^d]}|\nu_{\rho_*}^1)+
H(\nu_{\rho_0(\cdot)}^{N\sm u}|\nu_{\rho_*}^{\T_N^d\sm\{[Nu]\}}).
\end{equation}
By a simple computation $H(\delta_{[\alpha N^d]}|\nu_{\rho_*}^1)=-\log\nu_{\rho_*}([\rho N^d])$
for all $\rho\in[0,\rho_c]\cap\RR$ and therefore by property~\eqref{LocalJumpRateDef:c} of local jump rates
\begin{equation}\label{DiracCondensEntropWithRespToNoNCritInvarDist}
\lim_{N\ra+\infty}\fr{N^d}H(\delta_{[\rho N^d]}|\nu_{\rho_*}^1)=\lim_{N\ra+\infty}\frac{[\alpha N^d]}{N^d}\log\frac{\sqrt[{[\alpha N^d]}]{g!([\alpha N^d])}}{\Phi(\rho_*)}=\alpha\log\frac{\phi_c}{\Phi(\rho_*)}.
\end{equation}
Furthermore, we obviously have that 
\begin{equation*}H(\nu_{\rho_0(\cdot)}^{N\sm u}|\nu_{\rho^*}^{\T_N^d\sm\{[Nu]\}})=H(\nu_{\rho_0(\cdot)}^N|\nu^N_{\rho_*})
	-\Lambda_{\rho_*}^*\Big(\rho_0\Big(\frac{[Nu]}{N}\Big)\mn\rho_c\Big).
\end{equation*}
Since the profile $\rho_0:\T^d\lra\RR_+$ is assumed bounded and $\Lambda_{\rho_*}^*$ is continuous and finite on $\RR_+$ by~\eqref{RateFunc} we obviously have that $$\lim_{N\ra\infty}\fr{N^d}\Lambda_{\nu_{\rho_*}}^*\Big(\rho_0\Big(\frac{[Nu]}{N}\Big)\mn\rho_c\Big)=0$$ and therefore by~\eqref{DiracCondensEntropWithRespToNoNCritInvarDist},~\eqref{kleinmein} and~\eqref{MicroscopEntropGoToMarcoscopEntrBeta} it follows that
$$\lim_{N\ra+\infty}\fr{N^d}H(\nu_{\rho_0(\cdot);u,\alpha}^N|\nu_{\rho_*}^N)=
\int_{\T^d}\Lambda_{\rho_*}^*\big(\rho_0(u)\mn\rho_c\big)\df u+
\rho\log\frac{\phi_c}{\Phi(\rho_*)}$$
which is finite whenever $\phi_c<+\infty$.$\hfill\Box$

	\section{Main Results on condensing ZRPs}\label{MainResults}
	In all the following results we assume that the ZRP is weakly condensing i.e.~with finite critical fugacity $\phih_c<+\infty$ and that the laws $\{P^N\}$ of the (diffusively rescaled) ZRP start from a sequence of initial distributions $\mu_0^N\in\PP_1\MM_N^d$, $N\in\NN$, associated to a macroscopic profile $\mu_0\in\MMM_+(\T^d)$ and satisfying the $O(N^d)$-entropy assumption with constant $C_0<+\infty$ for some $\rho_*\in(0,\rho_c)$. For ZRPs with unbounded jump rate $\mathfrak{g}$ we furthermore assume that $\mu_0^N\in\PP_2\MM_N^d$, $N\in\NN$.
	
	The first result is the one-block estimate. In part (a) we generalize the one-block estimate for condensing ZRPs to unbounded cylinder functions $\Psi\colon\MM_N^d\to\RR$. In part (b) we reformulate the one-block estimate in terms of the joint law of the process $\s^{N,\Psi}$ defined in~\eqref{GenEmpProc}, the micro-empirical distribution $\bs\pi^{N,\ell}$ of the ZRP defined in~\eqref{Memp} and the $\bbar{\Psi}$-projection defined in~\eqref{PsiProj}. 
		\begin{theorema}\label{OBETheorem}{\rm{(One-block estimate for condensing ZRPs)}} \emph{(a)} Let $\bbar{\Psi}\colon\RR_+\to\RR$ be the extended homologue function defined in~\eqref{ExtendedHomologueAsymptlinear} of the asymptotically linear cylinder map $\Psi\colon\MM_\infty^d\to\RR$. Then for all $H\in L^1(0,T;C(\T^d))$ and all $\delta>0$  
		\begin{equation}\label{OBE}
		\limsup_{\ell,N\ra+\infty}\EE^N\Big|\int_0^T\fr{N^d}\sum_{x\in\T_N^d}H_t(x/N)\big[\tau_x\Psi^\ell(\eta_t)-\bbar{\Psi}\big(\eta_t^\ell(x)\big)\big]\df t\Big|=0.
		\end{equation}
	 \emph{(b)} Let $\Psi\colon\MM_\infty^d\to\RR$ be an asymptotically linear cylinder function. Then the extended homologue $\bbar{\Psi}\colon\RR_+\to\RR$ is asymptotically linear, i.e.~$\bbar{\Psi}\in\bbar{C}_1(\RR_+)$ and $\bbar{\Psi}'(\infty)=\ls\nabla\Psi(\infty,\1_J\rs$, the family of laws $\bbar{\bs{Q}}_\Psi^{N,\ell}:=(\s^{N,\Psi},\bs\pi^{N,\ell})_\sharp P^N$ is relatively compact and any limit point $$\bbar{\bs{Q}}_\Psi\in\Lim_{\ell,N\ra+\infty}\bbar{\bs{Q}}_\Psi^{N,\ell}\subs\PP\big(L_{w^*}^\infty(0,T;\MMM(\T^d))\x L_{w^*}^\infty(0,T;\bbar{\MMM}_1(\T^d\x\RR_+))\big)$$ is concentrated on the graph of the $\bbar{\Psi}$-projection, i.e.
			\begin{equation}\label{OBEYM}\bbar{\bs{Q}}_\Psi\big\{(\s,\bs\pi)\bigm|\s=B_{\bbar{\Psi}}(\bs\pi)\big\}=1\end{equation}
			where $B_{\bbar{\Psi}}\colon L_{w^*}^\infty(0,T;\bbar{\MMM}_1(\T^d\x\RR_+))\to L_{w^*}^\infty(0,T;\MMM(\T^d))$ is the $\bbar{\Psi}$-projection induced on the $L_{w^*}^\infty$-spaces of path-measures.
		\end{theorema}
	  
		\begin{theorema}\label{Theorem1}  We set $\W:=D(0,T;\MMM_+(\T^d))\x L_{w^*}^\infty(0,T;\X^1(\T^d)^*)\x L^\infty_{w^*}(0,T;\MMM_+(\T^d))$ and consider the image \begin{equation}\label{TripleLaw} Q^N:=(\pi^N,W^N,\s^N)_\sharp P^N\in\PP\W
		\end{equation} of the law $P^N$ of the diffusively rescaled ZRP starting from $\mu_0^N$ via the triple $(\pi^N,W^N,\s^N)$. Then the sequence $\{Q^N\}_{N\in\NN}\subs\PP\W$ is sequentially relatively compact in the weak topology of $\PP\W$. Furthermore, any limit point $Q^\infty$ of the sequence $\{Q^N\}$ is concentrated on trajectories $(\pi,W,\s)\in\W$ such that
		\begin{itemize}
		\item[{\rm{(a)}}] The continuity equation 
		\begin{equation}\label{CE}\begin{cases}\partial_t\pi+{\rm{div}}W_t=0\\
		W_t=-\nabla\s_t\end{cases}\quad\mbox{in }(0,T)\x\T^d
		\end{equation}
		holds in the sense of distributions.
		\item[{\rm{(b)}}] $\pi\in C(0,T;\MMM_+(\T^d))$, $\pi_0=\mu_0$ and $\pi_t(\T^d)=\mu_0(\T^d)$ for all $t\in[0,T]$.
		\item[{\rm{(c)}}] $\s_t\ll\LL_{\T^d}$, $\|\frac{\df\s_t}{\df\LL_{\T^d}}\|_{L^\infty(\T^d)}\leq\phih_c$ a.s$.$ for all $0\leq t\leq T$.
		\item[{\rm{(d)}}] $\s_t\ll\pi_t$, $\|\frac{\df\s_t}{\df\pi_t}\|_{L^\infty(\pi_t)}\leq\|\mathfrak{g}'\|_\infty$ a.s$.$ for all $0\leq t\leq T$, and 
		\item[{\rm{(e)}}] $\frac{\df\s_t}{\df\LL_{\T^d}}\in H^1(\T^d)$ and $W_t=-(\nabla\frac{\df\s_t}{\df\LL_{\T^d}})\df\LL_{\T^d}\in\MMM_0(\T^d;\RR^d)\leq\X^1(\T^d)^*$ for a.s.~all $0\leq t\leq T$.	
		\end{itemize}	
	\end{theorema}

	The proof of Theorem~\ref{Theorem1}(e) relies on the following regularity result which is worth stating in its own right. Let us note that for non-condensing ZRPs in which case $\phih_c=+\infty$ it is known by~\cite[Remark 5.1.8]{Kipnis1999a} that the first marginal of the law $Q^\infty$ is concentrated on trajectories $\pi\in D(0,T;\MMM_+(\T^d))$ such that $\pi_t\ll\LL_{\T^d}$ for all $0\leq t\leq T$. Thus by Theorem~\ref{Theorem1}(d) in case $\phih_c=+\infty$ in place Theorem~\ref{Theorem1}(c) it holds that $\s_t\ll\LL_{\T^d}$ and $\frac{\df\s_t}{\df\LL_{\T^d}}\leq\|\mathfrak{g}'\|_\infty\frac{\df\pi_t}{\df\LL_{\T^d}}$ for almost all $0\leq t\leq T$. In particular $\|\s\|_{TV;\infty}\leq\|\mathfrak{g}'\|_\infty\mu_0(\T^d)$ $Q_3$-a.s.~for all $\s\in L_{w^*}^\infty(0,T;\MMM_+(\T^d))$. Also, as will be clear from the proof, when $\phih_c=+\infty$ instead of Theorem~\ref{Theorem1}(e) it only holds that $\s_t\in W^{1,1}(\T^d)$ a.s.~for all $0\leq t\leq T$.
	\begin{theorema}\label{TheoremEnergyEstimate}
 Any subsequential limit point $Q_3$ of the sequence $\{Q_3^N\}$ of the third marginals on $L^\infty_{w^*}(0,T;\MMM_+(\T^d))$ of the laws $Q^N$ defined in~\eqref{TripleLaw} is concentrated on a $w^*$-measurable subspace of path-measures $\W_0\subs L_{w^*}^\infty(0,T;\MMM_{+,ac}(\T^d))\leq L_{w^*}^\infty(0,T;\MMM(\T^d))$ such that for all $\s\equiv\frac{\df\s}{\df\LL_{\T^d}}\in\W_0$ there exist $L^2((0,T)\x\T^d)$-functions ($L^1((0,T)\x\T^d)$ if the ZRP is non-condensing) denoted by $\pd_j\s$, $j=1,\dots,d$, satisfying 
		\begin{equation}\label{MeanTimeSob}
		\int_0^T\int_{\T^d}\pd_jH_t(u)\s(t,u)\df u\df t=-\int_0^T\int_{\T^d}H_t(x)\pd_j\s(t,u)\df u\df t
		\end{equation}
		and, setting $\nabla\s:=\sum_{j=1}^d\pd_j\s\cdot e_j$ for all $\s\in\W_0$, the energy estimate
		\begin{equation}\label{energyestimate}\int_0^T\int_{\T^d}\frac{\|\nabla\s(t,u)\|^2}{\s(t,u)}\df u\df t<+\infty.
		\end{equation}
		 holds. In particular $Q_3$ is concentrated in trajectories $\s\in L_{w^*}^\infty(0,T;\MMM_{+,ac}(\T^d))$ such that $\s_t\in H^1(\T^d)$ ($\s_t\in W^{1,1}(\T^d)$ if the ZRP is non-condensing) for almost all $0\leq t\leq T$.
	\end{theorema}

 \begin{theorema}\label{ClosedWeakEquation}
	Let $\bs{Q}^{N,\ell}:=\bs\pi^{N,\ell}_\sharp P^N$ be the law of the empirical generalized Young distribution of the ZRP and let 
	$$\bs{Q}\in\Lim_{\ell\ra+\infty}\Lim_{N\ra+\infty}\bs{Q}^{N,\ell}$$ be any limit point of $\{\bs{Q}^{N,\ell}\}$. Then $\bs{Q}$ is concentrated on generalized Young-measure-valued weak solutions $\bs\pi\in L_{w^*}^\infty(0,T;\bbar{\Y}_{1,\mathfrak{m}}(\T^d))$ of the non-linear diffusion equation
	$$\pd_t\bs\pi=\Delta\Phi(\bs\pi)$$
	in the sense of~\eqref{YMVZRPPDE} which (with the map $\Phi(\bs\rho_{\bs\pi_t})$ being defined as in~\eqref{YoungComp}) also satisfy the energy estimate
	\begin{equation}\label{PhiEnergEst}\int_0^T\int_{\T^d}\frac{\|\nabla\Phi(\bs\rho_{\bs\pi_t})(u)\|^2}{\Phi(\bs\rho_{\bs\pi_t})(u)}\df u\df t<+\infty.
	\end{equation} 
	
\end{theorema}

The next result concerns the two-blocks estimate for asymptotically linear cylinder maps $\Psi\colon\MM_\infty^d\to\RR$. Although we do not prove the full two-blocks estimate we prove a comparison property for the micro and macro-empirical density processes and a characterization of when the micro and macro-empirical processes $\bs\pi^{N,\ell}$ and $\bs\pi^{N,\ee}$ are interchangeable in the limit as $N\ra\infty$, $\ee\ra 0$ and then $\ell\ra\infty$. We will restrict attention to a subfamily $\{\bs\pi^{k_N^{(\ell)},m_\ell}\}_{(N,\ell)}$ of the micro-empirical processes along which the laws $\bs{Q}^{k_N^{(\ell)},m_\ell}=\bs\pi^{k_N^{(\ell)},m_\ell}_\sharp P^{k_N^{(\ell)}}$ converge. As a first we prove a truncated double-block estimate which in essence allows us to replace the mean number of particles $\eta^\ell(x)$ around a point $x\in\T_N^d$ with the truncated double-block average $(\eta^\ell(x)\mn M)^{[N\ee]}$. In terms of the \emph{micro-truncated double-block empirical density} $\bs\pi^{N,\ell;M;\ee}\colon D(0,T;\MM_N^d)\to L_{w^*}^\infty(0,T;\bbar{\MMM}_1(\T^d\x\RR_+))$ defined by 
\begin{equation}\label{DoubleBlockProc1}\lls F,\bs\pi^{N,\ell;M,\ee}\rrs= \int_0^T\fr{N^d}\sum_{x\in\T_N^d}\Big\{F_t\Big(\frac{x}{N},(\eta_t^\ell(x)\mn M)^{[N\ee]}\Big)+RF_t\Big(\frac{x}{N}\Big)(\eta_t^\ell(x)-M)^{+[N\ee]}\Big\}\df t\end{equation}
Furthermore, using the truncated double-block estimate we are able to compare the micro and macro-empirical densities $\bs\pi^{N,\ell}$ and $\bs\pi^{N,\ee}$ as $N\ra+\infty$, $\ee\ra 0$ and $\ell\ra+\infty$.
\begin{theorema}\label{TBCTheorem}\emph{(Two-blocks comparison)} Let $\{k_N^{(\ell)}\}_{N=1}^\infty$, $\ell\in\NN$, and $\{m_\ell\}_{\ell=1}^\infty$ be diverging sequences such that the subfamily of the micro-empirical laws $$\bs{Q}_*^{N,\ell}:=\bs{Q}^{k_N^{(\ell)},m_\ell}=\bs\pi^{k_N^{(\ell)},m_\ell}_\sharp P^{k_N^{(\ell)}}$$ converges weakly to a probability law $\bs{Q}_*^{\infty,\infty}\in\PP L_{w^*}^\infty(0,T;\bbar{\Y}_1(\T^d))$ as $N\uparrow\infty$ and then $\ell\ra\infty$. 
	\noindent{\rm{(a)}} {\rm{(Truncated double-block estimate)}} For all asymptotically linear maps $\Psi\in\bbar{C}_1(\RR_+)$ and all $G\in L^1(0,T;C(\T^d))$
	\begin{equation}\label{DBE}
	\lim_{M\ra+\infty}\limsup_{\ell\ra+\infty}\limsup_{\ee\ra0}\limsup_{N\ra+\infty}\EE^{k_N^{(\ell)}}\big|\lls G,B_\Psi(\bs\pi^{k_N^{(\ell)},m_\ell})-B_\Psi(\bs\pi^{k_N^{(\ell)},m_\ell;M;\ee})\rrs\big|=0.
	\end{equation}
	Consequently any subsequential limit point of the family of laws $$(\bs\pi^{k_N^{(\ell)},m_\ell},\bs\pi^{k_N^{(\ell)},m_\ell;M;\ee})_\sharp P^{k_N^{(\ell)}}\in\PP L_{w^*}^\infty(0,T;\bbar{\PP}_1(\T^d\x\RR_+))^2$$ as $N\ra+\infty$, $\ee\ra 0$, $\ell\ra+\infty$ and then $M\ra+\infty$ is concentrated on a measurable set of trajectory pairs $(\bs\pi^\infty,\bs\pi^0)\in L_{w^*}^\infty(0,T;\bbar{\Y}_1(\T^d))^2$ such that $B_\Psi(\bs\pi^\infty)=B_\Psi(\bs\pi^0)$ for all $\Psi\in\bbar{C}_1(\RR_+)$.\\		
\noindent{\rm{(b)}} {\rm{(Micro-macro block comparison)}} Any limit point of the family of laws
\begin{equation}\label{DBLAW}
\bbar{\bs{Q}}_*^{N,\ell,\ee}:=(\bs\pi^{k_N^{(\ell)},m_\ell},\bs\pi^{k_N^{(\ell)},\ee})_\sharp P^{k_N^{(\ell)}}\in\PP L_{w^*}^\infty(0,T;\bbar{\PP}_1(\T^d\x\RR_+))^2
\end{equation}
	as $N\ra+\infty$, $\ee\ra 0$ and $\ell\ra+\infty$ is concentrated on a measurable set of trajectory pairs $(\bs\pi^\infty,\bs\pi^0)\in L_{w^*}^\infty(0,T;\mathcal{\Y}_1(\T^d))^2$ such that for almost all $t\in[0,T]$ 
	\begin{itemize}
		\item[{\rm{(i)}}]  $B(\bs\pi_t^\infty)=B(\bs\pi_t^0)$ where $B\colon\bbar{\PP}_1(\T^d\x\RR_+)\to\MMM_+(\T^d)$ is the barycentric projection.
		\item[{\rm{(ii)}}]  The disintegrations $(\bs\rho^u_{\bs\pi_t^\infty})_{u\in\T^d}$, $(\bs\rho^u_{\bs\pi_t^0})_{u\in\T^d}$ satisfy $\bs\rho^u_{\bs\pi^\infty_t} \leq_{\rm{st}}\bs\rho^u_{\bs\pi_t^0}$ for Lebesgue a.s.~all $u\in\T^d$, i.e. for almost all $u\in\T^d$
		$$\int\Psi\df\bs\rho^u_{\bs\pi^\infty_t}\leq \int\Psi\df\bs\rho^u_{\bs\pi_t^0},\quad\mbox{for all non-decreasing maps }\Psi\in \bbar{C}_1(\RR_+).$$ 
		\item[{\rm{(iii)}}] $\rho^\perp_{\bs\pi^\infty_t}\geq\rho^\perp_{\bs\pi^0_t}$, i.e. $\rho^\perp_{\bs\pi^\infty_t}(f)\geq\rho^\perp_{\bs\pi^0_t}(f)$ for all $f\in C_+(\T^d)$.		
	\end{itemize}
\noindent{\rm{(c)}} Finally, the two-blocks estimate holds in the class $\bbar{C}_{1,\uparrow}(\RR_+)$ of non-decreasing maps $\Psi\in\bbar{C}_1(\RR_+)$ in the sense that the subsequential limit set $$\bbar{{\bs{\mathcal{Q}}}}_*^{\infty,\infty,0}:=\Lim_{\ell\uparrow,\ee\downarrow 0,N\uparrow\infty}\bbar{\bs{Q}}_*^{N,\ell,\ee}\subs\PP L_{w^*}^\infty(0,T;\bbar{\Y}_1(\T^d))$$ of the family defined in~\eqref{DBLAW} satisfies
\begin{equation}\label{TBE}\bbar{\bs{Q}}^*\big\{(\bs\pi^\infty,\bs\pi^0)\bigm|B_\Psi(\bs\pi^\infty)=B_\Psi(\bs\pi^0),\quad\forall\,\Psi\in\bbar{C}_{1,\uparrow}(\RR_+)\big\}=1,\quad\forall\bbar{\bs{Q}}_*\in\bbar{{\bs{\mathcal{Q}}}}_*^{\infty,\infty,0}\end{equation} if and only if for any subfamily $\{(k^{(m_\ell^{(1)})}_{k_N^{(1;\ell,i)}},m_{m_\ell^{(1)}},\ee_i^{(1;\ell)}\}$ of $\{(k_N^{(\ell)},m_\ell,\ee)\}$, where the sequence $\{m^{(1)}_\ell\}_{\ell\in\NN}$ is diverging, the sequences $\{\ee_i^{(1;\ell)}\}_{i=1}^\infty$ converge to $0$ for all $\ell\in\NN$ and $\{k_N^{(1;\ell,i)}\}_{N=1}^\infty$ is diverging for all $i,\ell\in\NN$, there exists a further subfamily 
\begin{equation}\label{Nightmare}
\bigg\{\bigg(k^{\big(m_{m_\ell^{(2)}}^{(1)}\big)}_{k_{k_N^{(2;\ell,i)}}^{\left(1;m_\ell^{(2)},\ee_i^{(2;\ell)}\right)}},m_{m_{m_\ell^{(2)}}^{(1)}},\ee_{\ee_i^{(2;\ell)}}^{\big(1;m_\ell^{(2)}\big)}\bigg)\bigg\}=:\{(\bar{k}_N^{(\ell,i)},\bar{m}_\ell,\bar{\ee}_i^{(\ell)})\}_{(N,\ell,i)}
\end{equation}
of $\{(k_N^{(\ell)},m_\ell,\ee)\}$ along which
\begin{equation}\label{ToProveReplacLemma}
\limsup_{M,\ell,i,N\uparrow\infty}
\EE^{\bar{k}_N^{(\ell,i)}}\int_0^T\fr{(\bar{k}_N^{(\ell,i)})^d}\sum_{x\in\T_{\bar{k}_N^{(\ell,i)}}^d}[(\eta_t^{\bar{m}_\ell}(x)-M)^+]^{[\bar{k}_N^{(\ell,i)}\bar{\ee}_i^{(\ell)}]}
\1_{[0,M]}(\eta_t^{\bar{m}_\ell}(x)^{[\bar{k}_N^{(\ell,i)}\bar{\ee}_i^{(\ell)}]})\df t=0.
\end{equation}
\end{theorema} 
\begin{rem}In case the two-blocks estimate as stated in~\eqref{TBE} holds then it also holds for all maps $\Psi\in\bbar{C}_1(\RR_+)$, i.e.
\begin{equation}\label{TBEVar}\bbar{\bs{Q}}^*\big\{(\bs\pi^\infty,\bs\pi^0)\bigm|B_\Psi(\bs\pi^\infty)=B_\Psi(\bs\pi^0),\quad\forall\,\Psi\in\bbar{C}_{1}(\RR_+)\big\}=1,\quad\forall\bbar{\bs{Q}}_*\in\bbar{{\bs{\mathcal{Q}}}}_*^{\infty,\infty,0}.\end{equation}
Indeed, since any function $\Psi\in\bbar{C}_1(\RR_+)\cap\Lip(\RR_+)$ can be written as the difference $\Psi=\Psi_1-\Psi_2$ of two non-decreasing maps $\Psi_1,\Psi_2\in\bbar{C}_{1,\uparrow}(\RR_+)$ and any map $\Psi\in\bbar{C}_1(\RR_+)$ can be approximated by the sequence of Lipschitz maps $\Psi_k(\lambda):=\Psi_{0,k}(\lambda)+\Psi'(\infty)\lambda$, $k\in\NN$, where $\Psi_{0,k}$ are the Moreau-Yosida approximations, given by~\eqref{MoreauYosida}, of the sublinear part $\Psi_0(\lambda):=\Psi(\lambda)-\Psi'(\infty)\lambda$ of $\Psi$,~\eqref{TBEVar} follows by the dominated convergence theorem.
\end{rem}

By the one-block estimate and the two-blocks comparison we obtain the following one-sided version of the replacement lemma in terms of the subsequential limit points of the family of joint laws 
\begin{equation}\label{RLLawsInitial}\bbar{\bs{Q}}_\Psi^N:=(\s^{N,\Psi},\pi^N)_\sharp P^N\in\PP\big(L_{w^*}^\infty(0,T;\MMM(\T^d))\x D(0,T;\MMM_+(\T^d)\big),\quad N\in\NN.
\end{equation}\vspace{-0.65cm}
\begin{theorema}\label{RLTheorem} {\rm{(Super-replacement lemma)}} Let $\bbar{\Psi}\colon\RR_+\to\RR$ be the extended homologue function defined in~\eqref{ExtendedHomologueSublinear} of a sublinear cylinder map $\Psi\colon\MM_\infty^d\to\RR_+$. Then:\\
\noindent\emph{(a)} The map 
	\begin{equation}\label{IPsi}
	D(0,T;\MMM_+(\T^d))\ni\pi\mapsto\big(\bbar{\Psi}(\pi^{ac}_t)\df\LL_{\T^d}\big)_{0\leq t\leq T}=:I_\Psi(\pi)\in L_{w^*}^\infty(0,T;\MMM_+(\T^d)),
	\end{equation} where $\pi=\pi^{ac}+\pi^\perp$ is the Radon Nikodym decomposition of $\pi$ with respect to Lebesgue measure, is measurable and if the extended homologue $\bbar{\Psi}$ is non-decreasing then for all subsequential limit points 
	\[\bbar{\bs{Q}}_\Psi^\infty\in\bbar{\bs{\mathcal{Q}}}_\Psi^\infty:=\Lim_{N\uparrow\infty}\bbar{\bs{Q}}_\Psi^N\]
	it holds that $\s^\Psi\leq\bbar{\Psi}(\pi^{ac})\df\LL_{\T^d}$ for $\bbar{\bs{Q}}_\Psi^\infty$-.a.s.~all $(\s^\Psi,\pi)\in L_{w^*}^\infty(0,T;\MMM_+(\T^d))\x D(0,T;\MMM_+(\T^d))$.

Furthermore, if each diverging sequence $\{k_N\}_{N=1}^\infty$ has a subsequence, still denoted by $\{k_N\}$, such that any subfamily $\{(k_N^{(\ell)},m_\ell)\}_{(N,\ell)}$ of $\{(k_N,\ell)\}_{(N,\ell)}$ has a subfamily, still denoted by $\{(k_N^{(\ell)},m_\ell)\}_{(N,\ell)}$, such that any subfamily $\{(k_N^{(\ell,i)},m_{\ell},\ee_i^{(\ell)})\}$ of $\{(k_N^{(\ell)},m_\ell,\ee)\}$ has a further subfamily $\{(\bar{k}_N^{(\ell,i)},\bar{m}_{\ell},\bar{\ee}_i^{(\ell)})\}$ along which~\eqref{ToProveReplacLemma} holds, then the full replacement lemma holds in the sense that for all sublinear cylinder maps $\Psi\colon\MM_\infty^d\to\RR_+$ 
\[\bbar{\bs{Q}}_\Psi^\infty\big\{\s^\Psi=\bbar{\Psi}(\pi^{ac})\df\LL_{\T^d}\big\}=1,\qquad\forall\bbar{\bs{Q}}_\Psi^\infty\in\bbar{\bs{\mathcal{Q}}}_\Psi^\infty.\]
\end{theorema}

By assuming that for a sequence $\{\mu_0^N\}$ of initial distributions of the ZRP the condition for the validity of the two-blocks estimate is satisfied, at least in some small time interval $[0,T_0]$, one obtains by the continuity equation~\eqref{CE} and the replacement lemma that all limit points of the laws $Q^N=\pi^N_\sharp\in\PP D(0,T;\MMM_+(\T^d))$ are concentrated on trajectories $\pi\in C(0,T;\MMM_+(\T^d))$ such that $\Phi(\pi^{ac})\equiv\Phi(\pi^{ac}\mn\rho_c)\in H^1(\T^d)$, the equation 
\[\pd_t\pi=\Delta\Phi(\pi^{ac}),\quad\pi=\pi^{ac}+\pi^\perp,\;\pi^{ac}\ll\LL_{\T^d},\;\pi^\perp\perp\LL_{\T^d}\]
holds in the sense of distributions and satisfy the energy estimate 
\[\int_0^T\int_{\{\pi_t^{ac}>0\}}\frac{\|\nabla\Phi(\pi_t^{ac}(u))\|^2}{\pi^{ac}_t(u)}\df u\df t<+\infty.\]
\section{Generalized Young measures}\label{GYM}
Our main goal is to study in more detail the space $\bbar{\MMM}_r(\T^d\x\RR_+):=\bbar{C}_r(\T^d\x\RR_+)^*$ of generalized Young-functionals. By definition they are the continuous linear functionals of the space $\bbar{C}_r(\T^d\x\RR_+)$, which is defined in~\eqref{barcr}. We recall also the recession operator $R_r\colon\bbar{C}_r(\T^d\x\RR_+)\to C(\T^d)$ which is defined in~\eqref{Recession}. Since the limit in the definition of the recession operator $R_r$ is assumed to be uniform it follows that $R_r$ is a contraction. Indeed for any $F\in\bbar{C}_r(\T^d\x\RR_+)$ we have that 
\[\|R_rF\|_\infty=\lim_{\lambda\ra+\infty}\|F^{(\lambda,r)}\|_\infty\leq\|F\|_{\infty,r}.\]
Here we denote by $F^{(\lambda,r)}\in C(\T^d)$ the map $F^{(\lambda,r)}:=\frac{F(\cdot,\lambda)}{1+\lambda^r}$. Furthermore $R_r$ is surjective since for any $f\in C(\T^d)$ we have that $R_rF=f$ where $F\in\bbar{C}_r(\T^d\x\RR_+)$ is given by $F(u,\lambda)=f(u)\lambda^r$ and its kernel is $\ker R_r=C_r(\T^d\x\RR_+)$. As is shown in Theorem~\ref{YoungRiesz} of Section~\ref{OYMSect}, by a simple application of the classic Riesz isomorphism $C_0(\T^d\x\RR_+)^*=\MMM(\T^d\x\RR_+)$
\begin{equation}\label{Riesz(r)}(C_r(\T^d\x\RR_+),\|\cdot\|_{\infty,r})^*\cong(\MMM_r(\T^d\x\RR_+),\|\cdot\|_{TV;r}).\end{equation} We will also identify the spaces $C(\T^d)^*$ and $\MMM(\T^d)$ via the classic Riesz isomorphism. Thus the adjoint $R_r^*$ of the recession operator yields a bounded and $w^*$-continuous operator $R_r^*\colon\MMM(\T^d)\to\bbar{\MMM}_r(\T^d\x\RR_+)$ via the formula $R^*(\mu)(F)=\int_{\T^d}RF\df\mu$. We introduce also the extension operator 
\[E\colon\MMM_r(\T^d\x\RR_+)\cong C_r(\T^d\x\RR_+)^*\to\bbar{\MMM}_r(\T^d\x\RR_+):=\bbar{C}_r(\T^d\x\RR_+)^*\]
defined by 
\begin{equation}\label{ExtensionOperator}
E(\bs\rho)(F)=\int_{\T^d\x\RR_+}F\df\bs\rho,\quad F\in\bbar{C}_r(\T^d\x\RR_+).
\end{equation}
This extension operator is well defined since $\bbar{C}_r(\T^d\x\RR_+)\subs \bigcap_{\bs\rho\in \MMM_r(\T^d\x\RR_+)}\LL^1(\bs\rho)$, where in the calligraphic $\mathcal{L}^1$-spaces we {\it{do not}} identify almost surely equal functions.

\begin{lemma}\label{InduceE} The extension operator $E\colon\MMM_r(\T^d\x\RR_+)\to\bbar{\MMM}_r(\T^d\x\RR_+)$ is the pointwise $w^*$-limit of a sequence of $w^*$-continuous operators and thus it is $w^*$-measurable. 
\end{lemma}\textbf{Proof} For each $M>0$ let $\Pi_M\colon\bbar{C}_r(\T^d\x\RR_+)\to C_r(\T^d\x\RR_+)$ be the linear operator defined by $\Pi_M(F)(u,\lambda)=F(u,\lambda\mn M)$. Then $\Pi_M$ is a contraction and its adjoint $E_M:=\Pi_M^*\colon\MMM_r(\T^d\x\RR_+)\to\bbar{\MMM}_r(\T^d\x\RR_+)$ is $w^*$-continuous. Thus it suffices to show that $E_M$ $w^*$-converges pointwise to $E$ as $M\ra+\infty$. Of course then $E$ will be $w^*$-measurable by Proposition~\ref{PointLimOp} in the appendix. So let $\bs\rho\in\MMM_r(\T^d\x\RR_+)$. For each $M>0$ and $F\in\bbar{C}_r(\T^d\x\RR_+)$ 
\[|E(\bs\rho)(F)-E_M(\bs\rho)(F)|\leq\int|F(u,\lambda)-F(u,\lambda\mn M)|\df|\bs\rho|(u,\lambda).\]
Obviously $F(u,\lambda)-F(u,\lambda\mn M)\lra 0$ as $M\ra+\infty$ and $|F(u,\lambda)-F(u,\lambda\mn M)|\leq 2\|F\|_{\infty,r}(1+\lambda^r)\in L^1(|\bs\rho|)$ so that an application of the dominated convergence theorem concludes the proof.$\hfill\Box$\\

 A generalized Young measure $\bs\pi$ is called \emph{regular} if $(E\circ j^*)(\bs\pi)=\bs\pi$ where $j^*$ is the adjoint of the natural inclusion $j\colon C_r(\T^d\x\RR_+)\to\bbar{C}_r(\T^d\x\RR_+)$ and it is called \emph{singular} if $j^*(\bs\pi)=0$. Thus a generalized Young measure $\bs\pi$ is regular if it is of the form $\bs\pi(F)=\int F\df\bs\rho$ for some $\bs\rho\in\MMM_r(\T^d\x\RR_+)$ and singular if it vanishes on all maps $F\in C_r(\T^d\x\RR_+)$, which as we will see implies that it is of the form $\bs\pi(F)=R^*(\mu)(F)=\int_{\T^d}RF\df\mu$ for some measure $\mu\in\MMM(\T^d)$. 

In the context of generalized Young measures we will denote by $U\colon\T^d\x\RR_+\to\T^d$ and $\Lambda\colon\T^d\x\RR_+\to\RR_+$ the natural projections. It is easy to see that for $r\in\RR_+$ and any continuous function $\Psi\colon\RR_+\to\RR_+$ such that 
$$\lim_{\lambda\ra+\infty}\frac{\Psi(\lambda)}{\lambda^r}=1,$$
the space $\bbar{C}_r(\T^d\x\RR_+)$ can be split as the direct sum
$$\bbar{C}_r(\T^d\x\RR_+)=C_r(\T^d\x\RR_+)\oplus\Psi(\Lambda)\cdot[C(\T^d)\circ U],$$
where $\Psi(\Lambda)\cdot[C(\T^d)\circ U]=\{\Psi(\Lambda)(f\circ U)\bigm|f\in C(\T^d)\}$. 

In Section~\ref{OYMSect} we prove the Riesz isomoprhism~\eqref{Riesz(r)} for ordinary Young measures. In Section~\ref{YMProof} we prove the representation Theorem~\ref{YoungMeasuresTheorem} below which yields the decomposition of a generalized Young measure $\bs\pi$ into a regular and a singular part. Then we adapt these results to the level of path-measures in Section~\ref{GYMP} and finally in Section~\ref{Barycentric} we describe the barycentric projection $B\colon\bbar{\MMM}_1(\T^d\x\RR_+)\to\MMM(\T^d)$.

\begin{theorema}\label{YoungMeasuresTheorem} {\rm{(a)}} Let the product space $\MMM_r(\T^d\x\RR_+)\x\MMM(\T^d)$ be equipped with the norm
	\begin{equation}\label{MixNormDef}
	\|(\bs{\rho},\mu)\|_{TV;r}:=\sup_{F\in\bbar{C}_r(\T^d\x\RR_+)\sm\{0\}}\frac{\int Fd\bs{\rho}+\int RFd\mu}{\|F\|_{\infty,r}}.
	\end{equation}
	The norm $\|\cdot\|_{TV;r}$ satisfies 
	\begin{equation}\label{MixNorm}\max\{\|\bs{\rho}\|_{TV;r}+\|\mu\|_{TV}\}\leq\|(\bs{\rho},\mu)\|_{TV,r}\leq\|\bs{\rho}\|_{TV;r}+\|\mu\|_{TV}\end{equation}
	and thus is equivalent to all the product norms on $\MMM_r(\T^d\x\RR_+)\x\MMM(\T^d)$.\\
	\noindent{\rm{(b)}} There is a unique isometry $I=(I^1,I^2)\colon\bbar{\MMM}_r(\T^d\x\RR_+)\to\MMM_r(\T^d\x\RR_+)\x\MMM(\T^d)$ such that
\begin{equation}\label{CharProp}\bs\pi(F)=\int_{\T^d\x\RR_+}F\df(E\circ I^1)(\bs\pi)+\int_{\T^d}RF\df I^2(\bs\pi),\quad\mbox{ for all }F\in\bbar{C}_r(\T^d\x\RR_+).
\end{equation}
Its inverse $J=I^{-1}$ is given by $J(\bs\rho,\mu)=E(\bs\rho)+R^*(\mu)$, i.e.
\begin{equation}\label{IInverse}
J(\bs\rho,\mu)(F)=\big(E(\bs\rho)+R^*(\mu)\big)(F)=\int F\df\bs\rho+\int RF\df\mu,\quad F\in\bbar{C}_r(\T^d\x\RR_+).
\end{equation}
\noindent{\rm{(c)}} The first coordinate $I^1$ of $I$ is the restriction operator $j^*$, i.e.~the adjoint of the natural inclusion $j\colon C_r(\T^d\x\RR_+)\hookrightarrow\bbar{C}_r(\T^d\x\RR_+)$, the second coordinate $I^2$ is given by the formula
\begin{equation}\label{I2AsymptForm}
I^2(\bs\pi)(f)=\lim_{M\ra+\infty}\bs\pi\big((\Lambda-M)^+\Lambda^{r-1}f(U)\big).
\end{equation}
\noindent{\rm{(d)}} The isometry $I$ is positive i.e.~it satisfies
\begin{equation}\label{ImageOfPositive}I\big(\bbar{\MMM}_{r,+}(\T^d\x\RR_+)\big)=\MMM_{r,+}(\T^d\x\RR_+)\x\MMM_+(\T^d).\end{equation}
\noindent{\rm{(e)}} The first coordinate $I^1=j^*$ of the isometry $I$ is $w^*$-continuous and the restriction of the second coordinate $I^2$ on $\bbar{\MMM}_{r,+}(\T^d\x\RR_+)$ is positively upper $w^*$-semicontinuous in the sense that for any net $\{\bs\pi_\alpha\}_{\alpha\in\A}\subs\bbar{\MMM}_{r,+}(\T^d\x\RR_+)$ converging to $\bs\pi\in\bbar{\MMM}_{r,+}(\T^d\x\RR_+)$ in the $w^*$-topology of $\bbar{\MMM}_r(\T^d\x\RR_+)$ and any non-negative map $f\in C_+(\T^d)$  
	\begin{equation}\label{IsomSemiCont}\limsup_\alpha I^2(\bs\pi_\alpha)(f)\leq I^2(\bs\pi)(f).\end{equation}
\end{theorema}

Before proceeding with the proof of Theorem~\ref{YoungMeasuresTheorem} let us see how it can easily be rephrased to yield a uniquely defined decomposition of a generalized Young measure into a regular and singular part. 
\begin{cor}\label{YMDecomp} (a) For any $\bs\pi\in\bbar{\MMM}_r(\T^d\x\RR_+)$ there exists a uniquely determined decomposition $\bs\pi=\widehat{\bs\pi}+\bs\pi^\perp$ of $\bs\pi$ where $\widehat{\bs\pi}\in\bbar{\MMM}_r(\T^d\x\RR_+)$ is regular and $\bs\pi^\perp\in\bbar{\MMM}_r(\T^d\x\RR_+)$ is singular. This decomposition satisfies
\begin{equation}\label{DecompNormIneq}
\max\{\|\widehat{\bs\pi}\|_{TV,r},\|\bs\pi^\perp\|_{TV,r}\}\leq\|\bs\pi\|_{TV,r}\leq\|\widehat{\bs\pi}\|_{TV,r}+\|\bs\pi^\perp\|_{TV,r}\end{equation}
and $\bs{\pi}$ is non-negative if and only if both $\widehat{\bs\pi}$ and $\bs\pi^\perp$ are non-negative.\\
\noindent(b) The operators $\widehat{D},D^\perp\colon\bbar{\MMM}_r(\T^d\x\RR_+)\to\bbar{\MMM}_r(\T^d\x\RR_+)$ defined by $\widehat{D}(\bs\pi)=\widehat{\bs\pi}$ and $D^\perp(\bs\pi)=\bs\pi^\perp$ are linear, bounded and $w^*$-Baire.\\
(c) The restriction $\widehat{D}|_{\bbar{\MMM}_{R,+}(\T^d\x\RR_+)}$ of $\widehat{D}$ on $\bbar{\MMM}_r(\T^d\x\RR_+)$ is positively $w^*$-lower semicontinuous and the restriction $D^\perp|_{\bbar{\MMM}_{r,+}(\T^d\x\RR_+)}$ is positively $w^*$-upper semicontinuous, i.e.~for any map $F\in\bbar{C}_{r,+}(\T^d\x\RR_+)$ and any net $\{\bs\pi_\alpha\}_{\alpha\in\A}\subs\bbar{\MMM}_{r,+}(\T^d\x\RR_+)$ converging to some $\bs\pi\in\bbar{\MMM}_{r,+}(\T^d\x\RR_+)$ in the $w^*$-topology
	\begin{equation}\label{DecompSemiCont}\liminf_{\alpha}\widehat{\bs\pi_\alpha}(F)\geq\widehat{\bs\pi}(F)\quad\mbox{and}\quad\limsup_{\alpha}\bs\pi_\alpha^\perp(F)\leq\bs\pi^\perp(F).
	\end{equation}
\end{cor}\textbf{Proof} (a) Let $\bs\pi\in\bbar{\MMM}_r(\T^d\x\RR_+)$, let $I=(I^1,I^2)$ be the isometry of Theorem~\ref{YoungMeasuresTheorem} and set $\widehat{\bs\pi}:=E(I^1(\bs\pi))$ and $\bs\pi^\perp:=R^*(I^2(\bs\pi))$ where $E$ and $R$ are the recession operators. Then since $E$ is a right inverse of the restriction operator $j^*$, i.e.~$j^*\circ E=\mathbbm{id}_{\MMM_r(\T^d\x\RR_+)}$ it obviously holds that 
$$(E\circ j^*)(\widehat{\bs\pi})=E\circ (j^*\circ E)(I^1(\bs\pi))=E(I^1(\bs\pi))=\widehat{\bs\pi}$$ so that $\widehat{\bs\pi}$ is regular. Since $R\circ j\equiv 0$ and thus $j^*\circ R^*\equiv 0$ it also holds that $j^*(\bs\pi^\perp)=j^*(R^*(\bs\pi))=0$. Thus $\bs\pi^\perp$ is singular and by the formula of $J:=I^{-1}$
$$E(I^1(\bs\pi))+R^*(I^2(\bs\pi))=J(I(\bs\pi))=\bs\pi$$
so that $(\widehat{\bs\pi},\bs\pi^{\perp})$ is a decomposition of $\bs\pi$ as the sum of a regular and a singular generalized Young measure. The decomposition $\bs\pi=\widehat{\bs\pi}+\bs\pi^\perp$ is unique since if $\bs\pi=\widehat{\bs\pi}_1+\bs\pi^\perp_1$ is another decomposition of $\bs\pi$ as a sum of a regular and singular generalized Young measure then $\bs\pi^\perp-\bs\pi^\perp_1$ is singular and thus $(\widehat{\bs\pi}-\widehat{\bs\pi}_1)(F)=(\bs\pi^\perp_1-\bs\pi^\perp)(F)=0$
for any $F\in C_r(\T^d\x\RR_+)$. Thus $j^*(\widehat{\bs\pi}-\widehat{\bs\pi}_1)=0$ and therefore since $\widehat{\bs\pi}$, $\widehat{\bs\pi}_1$ are regular  
$$\widehat{\bs\pi}-\widehat{\bs\pi}_1=E(j^*(\widehat{\bs\pi}-\widehat{\bs\pi}_1))=0,$$
which proves that the decomposition of generalized Young measures as a sum of regular and singular generalized Young measures is unique. Since $E$ is norm-preserving and, as we will see in the proof of Theorem~\ref{YoungMeasuresTheorem}, the adjoint $R^*$ of the recession operator is also norm-preserving, inequalities~\eqref{DecompNormIneq} follow by (a) of Theorem~\ref{YoungMeasuresTheorem}. The fact that $\bs\pi$ is non-negative if and only if $\widehat{\bs\pi}$ and $\bs\pi^\perp$ follows by (d) of the same Theorem. For the proof of (b) we note that $\widehat{D}=E\circ j^*$ is $w^*$-Baire according to Proposition~\ref{ContAnalyticComp} as the composition of a $w^*$-continuous operator $j^*$ with the $w^*$-Baire operator $E$. Consequently $D^\perp=\mathbbm{id}_{\bbar{\MMM}_r(\T^d\x\RR_+)}-\widehat{D}$ is also $w^*$-Baire. Finally the semicontinuity property~\eqref{DecompSemiCont} follows by Theorem~\ref{YoungMeasuresTheorem}(e) since $\bs\pi^\perp=R^*\circ I^2(\bs\pi)$, $R^*$ is $w^*$-continuous and $\bs\pi=\widehat{\bs\pi}+\bs\pi^\perp$. $\hfill\Box$\\

As it is evident by the proof of this corollary, the regular and singular decomposition operators ae given by the explicit relations 
\begin{equation}\label{RegularSingularDecompFormulas}
\widehat{D}=E\circ j^*,\qquad D^\perp=R^*\circ I^2,
\end{equation}
where $E\colon\MMM_1(\T^d\x\RR_+)\to\bbar{\MMM}_1(\T^d\x\RR_+)$ is the natural extension operator defined in~\eqref{ExtensionOperator}, $j\colon C_1(\T^d\x\RR_+)\hookrightarrow\bbar{C}_1(\T^d\x\RR_+)$ is the natural subspace inclusion, $R^*\colon\MMM(\T^d)\to\bbar{\MMM}_1(\T^d\x\RR_+)$ is the adjoint of the recession operator $R$ defined in~\eqref{Recession} and $I^2$ is the second coordinate of the isometry $I$ of Theorem~\ref{YoungMeasuresTheorem}.
	\subsection{A Riesz representation theorem for Young measures}\label{OYMSect}
Let $\Lambda\colon\T^d\x\RR_+\to\RR_+$ be the projection $\Lambda(u,\lambda)=\lambda$ on the second coordinate. Likewise we will denote by $U\colon\T^d\x\RR_+\to\T^d$ the projection on the first coordinate. Given $r\in(0,\infty)$ we denote by $B_r(\T^d\x\RR)$ the set of all measurable real-valued functions on $\T^d\x\RR_+$ with bounded polynomial growth of order $r$, i.e$.$ $$B_r(\T^d\x\RR_+):=\Big\{F\in\LL^0(\T^d\x\RR_+)\Bigm||F|\leq C(1+\Lambda^r)\;\mbox{for\;some}\;C\geq 0\Big\}.$$ 
Here $\LL^0(\T^d\x\RR_+)$ is the space of all measurable maps $F\colon\T^d\x\RR_+\to\RR$. By convention $B_0(\T^d\x\RR_+)=B(\T^d\x\RR_+)$ and for $r=+\infty$ we define $B_\infty(\T^d\x\RR_+)$ the set of all measurable maps that map bounded sets to bounded sets. As a shorthand we set $B_rC^k(\T^d\x\RR_+):=B_r(\T^d\x\RR_+)\cap C^k(\T^d\x\RR_+)$, $k\in\ZZ_+$.

Also, we will denote by $\MMM_r(\T^d\x\RR_+)$ the space of all finite Borel signed measures on $\T^d\x\RR_+$ with finite $r$-th moments, i.e$.$ 
$$\MMM_r(\T^d\x\RR_+):=\Big\{\bs\rho\in\MMM(\T^d\x\RR_+)\Bigm|\big\|\Lambda\big\|_{L^p(|\bs\rho|)}<\infty\Big\}.$$ Here for $r\in(0,1)$ we set $\|\Lambda\|_{L^r(|\bs\rho|)}:=\int\Lambda^r\df|\bs\rho|$. Again we set $\MMM_0(\T^d\x\RR_+)=\MMM(\T^d\x\RR_+)$. Note that $\MMM_\infty(\T^d\x\RR_+)$ is the space of all finite signed measures with bounded support. Finally, we set $\MMM_{r,+}(\T^d\x\RR):=\MMM_+(\T^d\x\RR_+)\cap \MMM_r(\T^d\x\RR)$ the set of non-negative measures with finite $r$-th moments and $\PP_r(\T^d\x\RR_+):=\PP(\T^d\x\RR_+)\cap\MMM_r(\T^d\x\RR_+)$ the set of probability measures with finite $r$-th moments. Then for all $r\in[0,+\infty]$
\begin{equation*}
B_r(\T^d\x\RR_+)=\bigcap_{\bs\rho\in\PP_r(\T^d\x\RR_+)}\mathcal{L}^1(\bs\rho).
\end{equation*}
\begin{prop}\label{Basic} Let $r\in(0,\infty)$. Then\\
	\noindent\emph{(a)} The space $B_r(\T^d\x\RR_+)$ becomes a Banach spaces with the norm 
	$$\|F\|_{\infty,r}:=\Big\|\frac{|F|}{1+\Lambda^r}\Big\|_\infty=\sup_{(u,\lambda)\in\T^d\x\RR_+}\frac{|F(u,\lambda)|}{1+\lambda^r},\quad\;F\in B_r(\T^d\x\RR_+).$$ 
	\noindent\emph{(b)} Convergence with respect to $\|\cdot\|_{\infty,r}$ implies uniform convergence in bounded subsets of $\T^d\x\RR_+$ (but not the converse). Therefore the subspace $B_rC(\T^d\x\RR_+)\leq B_r(\T^d\x\RR_+)$ is closed, and thus a Banach space when equipped with the restriction of $\|\cdot\|_{\infty,r}$.\\
	\noindent\emph{(c)} The space $\MMM_r(\T^d\x\RR_+)$ becomes a Banach space when equipped with the norm $$\|\bs\rho\|_{TV,r}=\|\bs\rho\|_{TV}+\|\Lambda\|^{r\mx 1}_{L^r(|\bs\rho|)}=\|(1+\Lambda^r)d\bs\rho\|_{TV},\;\bs\rho\in\MMM_r(\T^d\x\RR_+).$$
	\noindent\emph{(d)} The bilinear map $\ls\cdot,\cdot\rs:B_r(\T^d\x\RR_+)\x\MMM_r(\T^d\x\RR_+)\lra\RR$ given by $\ls F,\bs\rho\rs=\int F\df\bs\rho$
	satisfies 
	\begin{equation}\label{dualitybound}\qquad
	|\ls F,\bs\rho\rs|\leq\|F\|_{\infty,r}\|\bs\rho\|_{TV,r}
	\end{equation} for all $(F,\bs\rho)\in B_r(\T^d\x\RR_+)\x\MMM_r(\T^d\x\RR_+)$ and is a strongly non-degenerate dual pairing, i.e.~it induces the linear isometric inclusions 
	\begin{align*}
	&B_r(\T^d\x\RR_+)\ni F\mapsto\ls F,\cdot\rs\in\MMM_r(\T^d\x\RR_+)^*,\\
	&\MMM_r(\T^d\x\RR_+)\ni\bs\rho\mapsto\ls\cdot,\bs\rho\rs\in B_r(\T^d\x\RR_+)^*.\end{align*}
	\noindent\emph{(e)} The pairing $\ls\cdot,\cdot\rs$ is also a strongly non-degenerate dual pairing between the spaces $B_rC(\T^d\x\RR_+)$ and $\MMM_r(\T^d\x\RR_+)$. 
\end{prop} 
\textbf{Proof} (a) Let $\{F_n\}\subs B_r(\T^d\x\RR_+)$ be a Cauchy sequence. Then the sequence $G_n:=\frac{F_n}{1+\Lambda^r}$ is a Cauchy sequence in $B(\T^d\x\RR_+)$ which is a Banach space. Thus there exists $G\in B(\T^d\x\RR_+)$ such that $\|G_n-G\|_\infty\lra 0$. But then $F:=G(1+\Lambda^r)\in B_r(\T^d\x\RR_+)$ and $\|F_n-F\|_{\infty,r}\lra 0$ which proves that $B_r(\T^d\x\RR_+)$ is Banach. To prove (b) we suppose that $F_n\lra 0$ in $B_r(\T^d\x\RR_+)$ and let $M>0$. Then 
$$\sup_{\T^d\x[0,M]}|F_n|=(1+M^r)\sup_{\T^d\x[0,M]}\frac{|F_n|}{1+M^r}\leq
(1+M^r)\|F_n\|_{\infty;p}\stackrel{n\ra+\infty}\lra 0.$$
Therefore convergence in $B_r(\T^d\x\RR_+)$ implies uniform convergence in bounded subsets of $\T^d\x\RR_+$. To see that the converse is false consider the family $\{\Lambda^q\}_{0<q\leq r}\subs B_p(\T^d\x\RR_+)$. It is obvious that $\Lambda^q\lra\Lambda^r$ uniformly in bounded subsets of $\T^d\x\RR_+$ as $q\ra r$, while $\|\Lambda^r-\Lambda^q\|_{\infty,p}\geq 1$ for all $q<r$. The claim (c) follows from the fact that $(\MMM(\T^d\x\RR_+),\|\cdot\|_{TV})$ is a Banach space, since the function $T_r:\MMM_r(\T^d\x\RR_+)\lra\MMM(\T^d\x\RR_+)$ defined by $T_r(\bs\rho)=(1+\Lambda^r)\df\bs\rho$ is a surjective linear isometry. For the proof of the remaining claims (d) and (e) we note that inequality~\eqref{dualitybound} is obvious and it readily implies that $\|\ls F,\cdot\rs\|_{\MMM_r(\T^d\x\RR_+)^*}\leq \|F\|_{\infty,r}$ and $\|\ls\cdot,\bs\rho\rs\|_{B_r(\T^d\x\RR_+)^*}\leq \|\bs\rho\|_{TV;r}$ for all $F\in B_r(\T^d\x\RR_+)$ and all $\bs\rho\in\MMM_r(\T^d\x\RR_+)$. On the other hand, we have that 
$$\|\ls F,\cdot\rs\|_{\MMM_r(\T^d\x\RR_+)^*}=\sup_{\bs\rho\neq 0}\frac{|\ls F,\bs\rho\rs|}{\|\bs\rho\|_{TV,r}}\geq\sup_{(u,\lambda)\in\T^d\x\RR_+}\frac{|\ls F,\delta_{(u,\lambda)}\rs|}{\|\delta_{(u,\lambda)}\|_{TV,r}}=\|F\|_{\infty,r},$$
for all $F\in B_r(\T^d\x\RR_+)$ and therefore $\|\ls F,\cdot\rs\|=\|F\|_{\infty,r}$. Since $\|\bs\rho\|_{B_rC(\T^d\x\RR_+)^*}\leq\|\bs\rho\|_{B_r(\T^d\x\RR_+)^*}$ in order to complete the proof of the proposition it remains to show that $\|\bs\rho\|_{TV;r}\leq\|\bs\rho\|_{B_rC(\T^d\x\RR_+)^*}$ for all $\bs\rho\in\MMM_r(\T^d\x\RR_+)$. So let $\bs\rho\in\MMM_r(\T^d\x\RR_+)$ and let $X=P\cup N$ be a Hahn decomposition of $\T^d\x\RR_+$ with respect to $\bs\rho$. Since $\T^d\x\RR_+$ is polish, the finite measures $\bs\rho^+=\bs\rho|_{P}$ and $\bs\rho^-=-\bs\rho|_{N}$ are regular and thus for every $n\in\NN$ there exist compact sets $K_P^n\subs P$, $K_N^n\subs N$ such that $\bs\rho^+(P\sm K_P^n)\mx\bs\rho^-(N\sm K_N^n)\leq\fr{n}$ for all $n\in\NN$. Of course we can assume that the sequences $\{K_P^n\}_{n\in\NN}$ and $\{K_N^n\}_{n\in\NN}$ are increasing and if we set $K_P^\infty:=\bigcup_{n\in\NN}K_P^n$, $K_N^\infty:=\bigcup_{n\in\NN}K_N^n$ we obviously have that $\bs\rho^+(P\sm K_P^\infty)=\bs\rho^-(N\sm K_N^\infty)=0$. Since $\T^d\x\RR_+$ is a metric space, there exist for every $n\in\NN$ functions $\phi_P^n,\phi_N^n\in BC(\T^d\x\RR_+)$ such that 
$\1_{K_P^n}\leq\phi_P^n\leq 1-\1_{K_N^n}$ and $\1_{K_N^n}\leq\phi_N^n\leq 1-\1_{K_P^n}$ for all $n\in\NN$. Obviously the sequences $\{\phi_P^n\}$, $\{\phi_N^n\}$ converge pointwise in $K_P^\infty\cup K_N^\infty$. In particular $$\lim_{n\ra\infty}\phi_P^n=\begin{cases}1,&\mbox{in}\;\,K_P^\infty\\
0,&\mbox{in}\;\,K_N^\infty\end{cases}\qquad\lim_{n\ra\infty}\phi_N^n=\begin{cases}0,&\mbox{in}\;\,K_P^\infty\\
1,&\mbox{in}\;\,K_N^\infty\end{cases}.$$
But $|\bs\rho|(K_P^\infty\cup K_N^\infty)=|\bs\rho|(\T^d\x\RR_+)$ and therefore
$$\lim_{n\ra\infty}\phi_P^n=\1_{K_P^\infty}=\1_P,\quad\lim_{n\ra\infty}\phi_N^n=\1_{K_N^\infty}=\1_N,\quad|\bs\rho|\mbox{-a.e.}.$$
Since $|(\phi_P^n-\phi_N^n)(1+\Lambda^p)|\leq (1+\Lambda^r)\in L^1(|\bs\rho|)$ we have by the dominated convergence theorem that 
\begin{equation*}
\lim_{n\ra\infty}\int(\phi_P^n-\phi_N^n)(1+\Lambda^r)\df\bs\rho=
\int(\1_P-\1_N)(1+\Lambda^r)\df\mu=\int(1+\Lambda^r)\df|\bs\rho|
=\|\bs\rho\|_{TV;r}.
\end{equation*}
Therefore since $-1\leq2\1_{K_P^n}-1\leq\phi_P^n-\phi_N^n\leq1-2\1_{K_N^n}\leq 1$,
\begin{equation*}
\|\ls\cdot,\bs\rho\rs\|_{B_rC(\T^d\x\RR_+)^*}=\sup_{\|F\|_{\infty,r}\leq 1}\ls F,\bs\rho\rs
\geq\lim_{n\ra\infty}\int(\phi_P^n-\phi_N^n)(1+\Lambda^p)\df\bs\rho=\|\bs\rho\|_{TV,r}
\end{equation*}
and the proof is complete.$\hfill\Box$\\

We will denote by $C_0(\T^d\x\RR_+)$ the subspace of $BC(\T^d\x\RR_+)$ consisting of functions that vanish at infinity, i.e.~$F\in C_0(\T^d\x\RR_+)$ if and only if $$\lim_{\lambda\ra+\infty}\sup_{u\in\T^d}|F(u,\lambda)|=0,$$ which is a separable closed subspace of $BC(\T^d\x\RR_+)$. As we will see, by applying the Riesz representation theorem according to which $C_0(\T^d\x\RR_+)^*=\MMM(\T^d\x\RR_+)$ it follows that $\MMM_r(\T^d\x\RR_+)$ is a dual space, with separable predual the space 
$$C_r(\T^d\x\RR_+):=\Big\{F\in B_rC(\T^d\x\RR_+)\Bigm|F/(1+\Lambda^r)\in C_0(\T^d\x\RR_+)\Big\}.$$
Since $C_0(\T^d\x\RR_+)$ it follows that $C_r(\T^d\x\RR_+)$ is also a closed separable subspace of $B_rC(\T^d\x\RR_+)$ and so $C_r(\T^d\x\RR_+)$ is a separable Banach space with the restriction of the norm $\|\cdot\|_{\infty,r}$.

\begin{prop}\label{YoungRiesz}
	For any $r\in[0,+\infty)$ the dual pairing $$\ls\cdot,\cdot\rs\colon B_rC(\T^d\x\RR_+)\x\MMM_r(\T^d\x\RR_+)\to\RR$$ induces a linear surjective isometry 
	$$I_r\colon(\MMM_r(\T^d\x\RR_+),\|\cdot\|_{TV;r})\to (C_r(\T^d\x\RR_+),\|\cdot\|_{\infty,r})^*$$
	via the formula $I_r(\bs\rho)(F)=\ls F,\bs\rho\rs$. In particular 
	$$\|\bs\rho\|_{TV;r}=\sup_{\substack{F\in C_r(\T^d\x\RR_+)\\F\neq 0}}\frac{\int Fd\bs\rho}{\|F\|_{\infty,r}},\quad\bs\rho\in\MM_r(\T^d\x\RR_+).$$
\end{prop}\textbf{Proof} According to our definitions, the case $r=0$ is the Riesz representation theorem. The case $r>0$ is a simple consequence of the Riesz representation theorem. Indeed, recall that we have denoted by $T_r\colon\MMM(\T^d\x\RR_+)\to\MMM_r(\T^d\x\RR_+)$ the surjective isometry defined by $T_r(\bs\rho)=(1+\Lambda^r)d\bs\rho$. Also, the operator $S_r\colon C_r(\T^d\x\RR_+)\to C_0(\T^d\x\RR_+)$ defined by $S_r(F)=\frac{F}{1+\Lambda^r}$ is a linear surjective isometry and its adjoint $S_r^*\colon C_0(\T^d\x\RR_+)^*\to C_0(\T^d\x\RR_+)^*$ is also an isometry. It is elementary to check that the operator $I_r$ makes the diagram 
\begin{equation}\begin{tikzcd}
\MMM_0(\T^d\x\RR_+)\arrow[r, "I_0"]\arrow["T_r",d]
& C_0(\T^d\x\RR_+)^*\arrow[d, "S_r^*" ] \\
\MMM_r(\T^d\x\RR_+)\arrow[r, "I_r"]
& C_r(\T^d\x\RR_+)^*
\end{tikzcd}
\end{equation}
commutative. Therefore $I_r=S_r^*\circ I_0\circ T_r^{-1}$ is a linear surjective isometry as the composition of three surjective isometries.$\hfill\Box$\\

The $w^*$-topology that the space $\MMM_r(\T^d\x\RR_+)$ inherits as the dual of the separable Banach space $C_r(\T^d\x\RR_+)$ will be called the $C_r$-\emph{topology}. The Wasserstein topology of order $r$ is the (metrizable) topology characterized by 
$$\lim_{n\ra+\infty}\bs\rho_n=\bs\rho\quad\Llra\quad\lim_{n\ra+\infty}\int F\df\bs\rho_n=\int F\df\bs\rho,\quad\forall\;F\in B_rC(\T^d\x\RR_+).$$ Since $B_qC(\T^d\x\RR_+)\subs C_r(\T^d\x\RR_+)$ for all $q<r$ the $C_r$-topology is obviously stronger than the Wasserstein topology of order $q$ for all $q<r$ but weaker than the Wasserstein topology of order $r$. In fact $\lim_{n\ra\infty}\bs\rho_n=\bs\rho$ in the $r$-th Wasserstein topology if and only if $\bs\rho_n\lra\bs\rho$ in the $C_r$-topology and $\lim_{n\ra\infty}\int\Lambda^r\df \bs\rho_n=\int\Lambda^r\df\bs\rho$. 
\begin{prop}\label{BarCrBanach} The space $\bbar{C}_r(\T^d\x\RR_+)$ is a closed subspace of $B_rC(\T^d\x\RR_+)$ and thus a Banach space when equipped with the norm $\|\cdot\|_{\infty,r}$.\end{prop}
\textbf{Proof} Let $\{F_n\}_{n=1}^\infty\subs \bbar{C_r}(\T^d\x\RR_+)$ be a sequence converging to $F\in B_rC(\T^d\x\RR_+)$, i.e.
\begin{equation}\label{1ForBan}
\lim_{n\ra+\infty}\|F_n-F\|_{\infty,r}=0.\end{equation} By the definition of the space $\bbar{C}_r(\T^d\x\RR_+)$ for each $n\in\NN$ there exists $\bbar{F}_n\in C(\T^d)$, i.e.~$\bbar{F}_n:=RF_n$ is the recession function of $F_n$, such that \begin{equation}\label{2ForBan}\lim_{\lambda\ra+\infty}\|F_n^{(\lambda,r)}-\bbar{F}_n\|_\infty=0,\quad F^{(\lambda,r)}(\cdot):=\frac{F(\cdot,\lambda)}{1+\lambda^r}\in C(\T^d),\;\lambda\geq 0.\end{equation} With this notation $\|F\|_{\infty,r}=\sup_{\lambda\geq 0}\|F^{(\lambda,r)}\|_\infty$ for any $F\in B_rC(\T^d\x\RR_+)$ and thus for all $n,m\in\NN$ 
\begin{align*}
\|\bbar{F}_n-\bbar{F}_m\|_\infty
&\leq\|\bbar{F}_n-F^{(\lambda,r)}_n\|_\infty+\|F_n-F_m\|_{\infty,r}+\|F_m^{(\lambda,r)}-\bbar{F}_m\|_\infty.
\end{align*}
Taking the limit as $\lambda\ra+\infty$ we obtain by~\eqref{2ForBan} that $\|\bbar{F}_n-\bbar{F}_m\|_\infty\leq\|F_n-F_m\|_{\infty,r}$ which shows that $\{\bbar{F}_n\}$ is Cauchy in $C(\T^d)$ since $\{F_n\}$ converges to $F$ in $B_rC(\T^d\x\RR_+)$. Therefore, since $C(\T^d)$ is Banach, there exists $\bbar{F}\in C(\T^d)$ such that $\lim_{n\ra+\infty}\|\bbar{F}_n-\bbar{F}\|_\infty=0$ and with $F$ being the limit of $\{F_n\}$ in $B_rC(\T^d\x\RR_+)$ we have to show that 
\begin{equation}\label{ToProveBan1}\lim_{\lambda\ra+\infty}\|F^{(\lambda,r)}-\bbar{F}\|_\infty=0.\end{equation}
But 
$$\lim_{\lambda\ra+\infty}\|F^{(\lambda,r)}-\bbar{F}\|_\infty=\lim_{\lambda\ra+\infty}\lim_{n\ra+\infty}\|F_n^{(\lambda,r)}-\bbar{F}_n\|_\infty$$
and thus in order for~\eqref{ToProveBan1} to hold the double limit as $n\ra+\infty$ and then $\lambda\ra+\infty$ above must be interchangeable. But this is indeed true since by the assumption that $\{F_n\}$ converges to $F$ in $B_rC(\T^d\x\RR_+)$ that the convergence $\lim_{n\ra+\infty}F_n^{(\lambda,r)}=F^{(\lambda,r)}$ in $C(\T^d)$ is in fact uniform over all large $\lambda$ and thus the double limit can be interchanged. Indeed, for all $n\in\NN$ 
\begin{align*}
\sup_{\lambda\geq 0}\|F_n^{(\lambda,r)}-F^{(\lambda,r)}\|_\infty=\sup_{\lambda\geq 0}\Big\|\frac{F_n(\cdot,\lambda)-F(\cdot,\lambda)}{1+\Lambda(\cdot,\lambda)^r}\Big\|_\infty=\|F_n-F\|_\infty
\end{align*}
and therefore
\begin{align*}
\|F^{(\lambda,r)}-\bbar{F}\|\leq\|F^{(\lambda,r)}-F_n^{(\lambda,r)}\|_\infty+\|F_n^{(\lambda,r)}-\bbar{F}\|_\infty\leq\|F_n-F\|_{\infty,r}+\|F_n^{(\lambda,r)}-\bbar{F}\|_\infty
\end{align*}
which shows that 
$$\limsup_{\lambda\ra+\infty}\|F^{(\lambda,r)}-\bbar{F}\|\leq\|F_n-F\|_{\infty,r}+\|\bbar{F}_n-\bbar{F}\|_\infty$$
and taking the limit as $n\ra+\infty$ we conclude the proof.$\hfill\Box$
%
\subsection{Proof of Theorem~\ref{YoungMeasuresTheorem}}\label{YMProof}

The representation of generalized Young measures via a pair of an ordinary Young measure and a Borel measure is based on the following two functional analytic lemmas. Before stating those lemmas let us recall that for any Banach subspace $C_0$ of a Banach space $\bbar{C}_0$ we denote by $$C_0^\perp:=\big\{\bs\pi\in\bbar{C}_0^*\bigm|\bs\pi(F)=0\mbox{ for all }F\in C_0\big\}$$ the annihilator of $C_0$ in $\bbar{C}_0^*$. If $j\colon C_0\hookrightarrow\bbar{C}_0$ is the natural inclusion and $j^*\colon\bbar{C}_0^*\to C_0^*$ its adjoint operator $j^*(\bs\pi)=\bs\pi|_{C_0}$, $\bs\pi\in\bbar{C}_0^*$ then $\ker j^*=C_0^\perp$. In the abdtract Lemma~\ref{SplittingLemma} it might be useful conceptually to have in mind the spaces $C_0=C_r(\T^d\x\RR_+)$, $\bbar{C}_0=\bbar{C}_r(\T^d\x\RR_+)$ and $E$ the extension operator defined in~\eqref{ExtensionOperator} and in Lemma~\ref{RecessionSplit} that follows it the recession operator $R\colon\bbar{C}_0\to C:=C(\T^d)$ with kernel $\ker R=C_0$.
\begin{lemma}\label{SplittingLemma}
	Let $C_0$ be a subspace of the Banach space $\bbar{C}_0$ and let $j\colon C_0\hookrightarrow \bbar{C}_0$ be the natural inclusion. Let $E\colon C_0^*\to \bbar{C}_0^*$  be a linear extension operator, i.e. $j^*\circ E=\mathbbm{id}_{C_0^*}$. Then the map 
	$$P_E:=\mathbbm{id}_{\bbar{C}_0^*}-E\circ j^*\colon\bbar{C}_0^*\to\bbar{C}_0^*$$
	is a linear projection on $\ker j^*= C_0^\perp$, i.e.~$P_E(\bbar{C}_0)\subs C_0^\perp$ and $P_E|_{C_0^\perp}=\mathbbm{id}_{C_0^\perp}$ and the map
	\begin{equation}\label{SumSplit}
	T_E=(j^*,P_E)\colon\bbar{C}_0^*\to C_0^*\x C_0^\perp
	\end{equation} is a linear isomorphism with inverse $S_E$ given by 
$$S_E(\bs\rho,\bs\pi^\perp)=E(\bs\rho)+\bs\pi^\perp,\quad\forall(\bs\rho,\bs\pi^\perp)\in C_0^*\x C_0^\perp.$$
Furthermore, the norm $\|\cdot\|_0$ on $C_0^*\x C_0^\perp$ that makes $T_E$ an isometry, i.e. $$\|(\bs\rho,\bs\pi^\perp)\|_0:=\|S_E(\bs\rho,\bs\pi^\perp)\|_{\bbar{C}_0^*}=\sup_{\|F\|_{\bbar{C}_0}\neq 0}\frac{E(\bs\rho)(F)+\bs\pi^\perp(F)}{\|F\|_{\bbar{C}_0}},\quad(\bs\rho,\bs\pi^\perp)\in C_0^*\x C_0^\perp$$ satisfies 
$$\|\bs\rho\|_{C_0^*}\leq\|(\bs\rho,\bs\pi^\perp)\|_0\leq\|E(\bs\rho)\|_{\bbar{C}_0^*}+\|\bs\pi^\perp\|_{\bbar{C}_0^*}$$
and if $E$ is such that for all $F\in\bbar{C}_0$  
\begin{equation}\label{ForLowInftB}
\exists\{F_n\}\subs F+C_0:\; \|F_n\|_{\bbar{C}_0^*}\leq\|F\|_{\bbar{C}_0^*}\mbox{ and }\lim_{n\ra+\infty}E(\bs\rho)(F_n)=0,\quad\forall\bs\rho\in C_0^*
\end{equation} then also $\|(\bs\rho,\bs\pi^\perp)\|_0\geq\|\bs\pi^\perp\|_{\bbar{C}_0^*}$ and $P_E$ is a contraction.
\end{lemma}\textbf{Proof} The kernel $\ker j^*$ is always equal to the annihilator $C_0^\perp$ since $j^*(\bs\pi)=0$ if and only if $\bs\pi(F_0)=0$ for all $F\in C_0$.. The map $P_E$ is obviously linear and it is a projection on $\ker j^*$ since on one hand $$j^*\circ P_E=j^*-j^*\circ(E\circ j^*)=0$$ which implies that $\Ima P_E:=P_E(\bbar{C}_0^*)\subs\ker j^*$ and on the other hand for any $\bs\pi^\perp\in C_0^\perp$ we have $j^*(\bs\pi^\perp)=0$ and thus $P_E(\bs\pi^\perp)=\bs\pi^\perp-E(j^*(\bs\pi^\perp))=\bs\pi^\perp-E(0)=\bs\pi^\perp$ so that $P_E|_{C_0^\perp}=\mathbbm{id}_{C_0^\perp}$ and thus $P_E$ is indeed a linear projection on $\Ima P_E=\ker j^*=C_0^\perp$. Furthermore the projection $P_E$ induces the direct sum decomposition
\begin{equation}\label{GenYoungSplitAbs}
\bbar{C}_0^*=\ker P_E\oplus\Ima P_E=\Ima E\oplus\ker j^*=\Ima E\oplus C_0^\perp.
\end{equation}
Indeed, any $\bs\pi\in\bbar{C}_0^*$ can be written as the sum $\bs\pi=(\bs\pi-P_E(\bs\pi))+P_E(\bs\pi)$. Here obviously $P_E(\bs\pi)\in \Ima P_E=\ker j^*$ and $\bs\pi-P_E(\bs\pi)=E(j^*(\pi))\in\Ima E$. But by the assumption $j^*\circ E=\mathbbm{id}_{C_0^*}$ we have that 
$$P_E\circ E=E-(E\circ j^*)\circ E=E-E=0$$
and therefore $\Ima E\subs\ker P_E$. Of course it is obvious by the definition of $P_E$ that $\ker P_E\subs\Ima E$. This sum in~\eqref{GenYoungSplitAbs} is direct since if $\bs\pi\in \ker P_E\cap\Ima P_E$ then on one hand we have that $\bs\pi=E(j^*(\bs\pi))$ while on the other hand $\bs\pi\in\Ima P_E=\ker j^*$ and thus $\bs\pi=E(j^*(\bs\pi))=E(0)=0$.

Consequently the map $\wt{T}_E=(\mathbbm{id}_{\bbar{C}_0^*}-P_E,P_E)\colon\bbar{C}_0^*\to\Ima E\x C_0^\perp$ is a linear inverse of the summation operation $+\colon \Ima E\x C_0^\perp\to\bbar{C}_0^*$. Furthermore the map $E$ is necessarily injective as a right inverse and since it is a right inverse for $j^*$ in fact $\|E(\rho)\|_{\bbar{C}_0^*}\geq\|\bs\rho\|_{C_0^*}$ for all $\bs\rho\in C_0^*$. In particular $\Ima E\cong C_0^*$ with a linear isomorphism being the map $j^*|_{\Ima E}\colon\Ima E\to C_0^*$. Thus the map $j^*|_{\Ima E}\x\mathbbm{id}_{C_0^\perp}$ is a linear isomorphism and since $j^*\circ P_E=0$
$$(j^*|_{\Ima E}\x\mathbbm{id}_{C_0^\perp})\circ\wt{T}_E=\big(j^*(\mathbbm{id}_{\bbar{C}_0^*}-P_E),P_E\big)=(j^*,P_E)=T_E.$$
Thus the map $T_E$ defined in~\eqref{SumSplit} is a linear isomorphism as the composition of linear isomorphisms. Since $\wt{T}_E^{-1}=+$ and $j^*|_{\Ima E}^{-1}=E\colon C_0^*\to \Ima E$ it follows that the inverse of $T_E$ is given by 
$$S_E=T_E^{-1}=\wt{T}_E^{-1}\circ(j^*|_{\Ima E}\x\mathbbm{id}_{C_0^\perp})^{-1}=+\circ(E\x\mathbbm{id}_{C_0^\perp})=E+\mathbbm{id}_{C_0^\perp}$$
as claimed.

We prove finally the bounds for the norm $\|\cdot\|_0$ on $C_0^*\x C_0^\perp$. So let $(\bs\rho,\bs\pi^\perp)\in C_0^*\x C_0^\perp$. First, since $\|F\|_{C_0}=\|F\|_{\bbar{C}_0}$ and $\bs\pi^\perp(F)=0$ for all $F\in C_0$  
$$\|(\bs\rho,\bs\pi^\perp)\|_0=\sup_{F\in\bbar{C}_0\sm\{0\}}\frac{E(\bs\rho)(F)+\bs\pi^\perp(F)}{\|F\|_{\bbar{C}_0}}\geq \sup_{F\in C_0\sm\{0\}}\frac{\bs\rho(F)}{\|F\|_{C_0}}=\|\bs\rho\|_{C_0^*}.$$
On the other hand we obviously have that $\|(\bs\rho,\bs\pi^\perp)\|_0\leq\|E(\bs\rho)\|_{\bbar{C}_0^*}+\|\bs\pi^\perp\|_{\bbar{C}_0^*}$ and so it remains to show the inequality $\|(\bs\rho,\bs\pi^\perp)\|_0\geq\|\bs\pi^\perp\|_{\bbar{C}_0^*}$ under the assumption that $E$ satisfies property~\eqref{ForLowInftB} for all $F\in\bbar{C}_0$. So let $\ee>0$ and choose $F_\ee\in\bbar{C}_0$ such that $\|\bs\pi^\perp\|_{\bbar{C}_0^*}\leq\frac{\bs\pi^\perp(F_\ee)}{\|F_\ee\|_{\bbar{C}_0}}-\ee$. There exists then a sequence $\{F_n^\ee\}_{n=1}^\infty\subs F_\ee+C_0$ such that $\|F_n^\ee\|_{\bbar{C}_0}\leq\|F_\ee\|_{\bbar{C}_0}$ and $\lim_{n\ra+\infty}E(\bs\rho)(F_n^\ee)=0$ for all $\bs\rho\in C_0^*$. Since $\{F_n^\ee\}\subs F_\ee+C_0$ and $\bs\pi^\perp\in C_0^\perp$ we have that $\bs\pi^\perp(F_n^\ee)=\bs\pi^\perp(F_\ee)$ for all $n\in\NN$ and thus
\begin{align*}\|(\bs\rho,\bs\pi^\perp)\|_0&\geq\sup_{n\in\NN}\frac{E(\bs\rho)(F_n^\ee)+\bs\pi^\perp(F_n^\ee)}{\|F_n^\ee\|_{\bbar{C}_0}}\geq
\sup_{n\in\NN}\frac{E(\bs\rho)(F_n^\ee)+\bs\pi^\perp(F_\ee)}{\|F_\ee\|_{\bbar{C}_0}}\\
&\geq\limsup_{n\ra+\infty}\frac{E(\bs\rho)(F_n^\ee)+\bs\pi^\perp(F_\ee)}{\|F_\ee\|_{\bbar{C}_0}}=\frac{\bs\pi^\perp(F_\ee)}{\|F_\ee\|_{\bbar{C}_0}}\geq\|\bs\pi^\perp\|_{\bbar{C}_0^*}+\ee.
\end{align*}
Since this is true for any $\ee>0$ it follows that $\|(\bs\rho,\bs\pi^\perp)\|_0\geq\|\bs\pi^\perp\|_{\bbar{C}_0^*}$ as required. In particular for all $\bs\pi\in\bbar{C}_0^*$ we have $\|\bs\pi\|_{\bbar{C}_0^*}=\|(j^*(\bs\pi),P_E(\bs\pi))\|_0\geq\|P_E(\bs\pi)\|_{\bbar{C}_0^*}$
and thus $P_E$ is a contraction.$\hfill\Box$

\begin{lemma}\label{RecessionSplit} Let $R\colon\bbar{C}_0\to C$ be a bounded surjective of Banach spaces. We set $C_0:=\ker R$ and $j\colon C_0\hookrightarrow \bbar{C}_0$ the subspace inclusion. If $j^*\colon\bbar{C}_0^*\to C_0^*$ admits a norm-preserving right inverse $E\colon C_0^*\to\bbar{C}_0^*$ there exists a unique mapping $I=(I_1,I_2)\colon\bbar{C}_0^*\to C_0^*\x C^*$ such that 
	\begin{equation}\label{AbsCharProp}\bs\pi(F)=E(I_1(\bs\pi))(F)+I_2(\bs\pi)(RF),\quad\forall(F,\bs\pi)\in\bbar{C}_0\x\bbar{C}_0^*.\end{equation}
	This mapping is a linear isomorphism with inverse $J=I^{-1}$ given by the formula 
	$$J(\bs\rho,\mu)(F)=E(\bs\rho)(F)+R^*(\mu)(F)=E(\bs\rho)(F)+\mu(RF),\quad F\in\bbar{C}_0$$
	for all $(\bs\rho,\mu)\in C_0^*\x C^*$. Furthermore, if $E$ satisfies~\eqref{ForLowInftB} for all $F\in\bbar{C}_0$ and $R$ is a contraction such that 
	\begin{equation}\label{RecessBound}
	\mbox{for all }f\in C\mbox{ there exists }F\in R^{-1}(\{f\})\mbox{ such that }\|F\|_{\bbar{C}_0}\leq\|f\|_C
	\end{equation} then the adjoint $R^*$ is norm-preserving and the norm $\|\cdot\|_*$ on $C_0^*\x C^*$ that makes $I$ an isometry satisfies 
	$$\max\big\{\|\bs\rho\|_{C_0^*},\|\mu\|_{C^*}\big\}\leq\|(\bs\rho,\mu)\|_*\leq\|\bs\rho\|_{C_0^*}+\|\mu\|_{C^*},\quad(\bs\rho,\mu)\in C_0^*\x C^*$$
and is thus equivalent with all the product norms on $C_0^*\x C^*$.
\end{lemma}\textbf{Proof} If such a map $I$ exists then for all $F\in C_0=\ker R$ we have that $j^*(\pi)(F)=I_1(\bs\pi)$ and thus $I_1=j^*$. Then $I_2\colon\bbar{C}_0^*\to C^*$ satisfies $I_2(\bs\pi)(RF)=\bs\pi(F)-E(j^*(\bs\pi))(F)$ for all $F\in\bbar{C}_0$ which since $R$ is surjective characterizes $I_2$ uniquely. Thus if such a map $I$ exists it is unique. The existence of this map $I$ follows by Lemma~\ref{SplittingLemma} and the first isomorphism theorem of linear algebra. Indeed, Lemma~\ref{SplittingLemma} yields a linear isomorphism $T_E\colon\bbar{C}_0^*\to C_0^*\x C_0^\perp$. Furthermore $C_0^\perp$ is isomorphic to $(\,^{\bbar{C}_0}/_{C_0})^*$ with an isomorphism being given by the adjoint of the natural quotient map $[\cdot]_{C_0}\colon\bbar{C}_0\to\,^{\bbar{C}_0}/_{C_0}$. Indeed, since $[\cdot]_{C_0}$ is surjective its adjoint $[\cdot]^*_{C_0}\colon(\,^{\bbar{C}_0}/_{C_0})^*\to\bbar{C}_0^*$ is injective and $\Ima[\cdot]^*_{C_0}\leq C_0^\perp$. The map $Q\colon C_0^\perp\to(\,^{\bbar{C}_0}/_{C_0})^*$ given by the formula 
$Q(\bs\pi^\perp)(F+C_0)=\bs\pi^\perp(F)$ is well-defined on the specified domains since $\bs\pi^\perp\in C_0^\perp$ and $Q(\bs\pi^\perp)\colon\,^{\bbar{C}_0}/_{C_0}\to\RR$ is bounded with $\|Q(\bs\pi^\perp)\|_{(\,^{\bbar{C}_0}/_{C_0})^*}\leq\|\bs\pi^\perp\|_{\bbar{C}_0^*}$. The map $Q$ is the inverse of $[\cdot]^*_{C_0}\colon(\,^{\bbar{C}_0}/_{C_0})^*\to C_0^\perp$. Now since $R\colon\bbar{C}_0\to C$ is a bounded surjection with $\ker R=C_0$ by the first isomorphism theorem it induces a linear isomorphism $\,^R/_{C_0}\colon\,^{\bbar{C}_0}/_{C_0}\to C$ via $\,^R/_{C_0}(F+C_0)=R(F)$. The induced map $\,^R/_{C_0}$ is obviously bounded and by the inverse mapping theorem it has a continuous inverse. Consequently its adjoint $(\,^R/_{C_0})^*\colon C^*\to(\,^{\bbar{C}_0}/_{C_0})^*$ is a bi-Lipschitz linear isomorphism. Thus we can define the map $I\colon\bbar{C}_0^*\to C_0^*\x C^*$ as 
\begin{equation}\label{IIsomAbstractForm}
I:=\big(\mathbbm{id}_{C_0^*}\x((\,^R/_{C_0})^*)^{-1}\big)\circ(\mathbbm{id}_{C_0^*}\x Q)\circ T_E.
\end{equation}
Then $I$ is a linear isomorphism by definition. Since $T_E=(j^*,P_E)$ we obviously have that $I_1=j^*$ and $I$ satisfies~\eqref{AbsCharProp} since for all $\bs\pi\in\bbar{C}_0^*$ and all $F\in\bbar{C}_0$
\begin{align*}
E(I_1(\bs\pi))(F)+I_2(\bs\pi)(RF)&=E(j^*(\bs\pi))(F)+((\,^R/_{C_0})^{-1})^*\circ Q\circ P_E(\bs\pi)(RF)\\
&=E(j^*(\bs\pi))(F)+Q\circ P_E(\bs\pi)(F+C_0)\\
&=E(j^*(\bs\pi))(F)+P_E(\bs\pi)(F)=\bs\pi(F).
\end{align*}

The inverse of $I$ is the map \begin{equation}\label{J}
J=S_E\circ(\mathbbm{id}_{C_0^*}\x[\cdot]^*_{C_0})\circ\big(\mathbbm{id}_{C_0^*}\x(\,^R/_{C_0})^*\big)\colon C_0^*\x C^*\to\bbar{C}_0^*.
\end{equation} Since $[\cdot]^*_{C_0}\circ(\,^R/_{C_0})^*=\big((\,^R/_{C_0})\circ[\cdot]_{C_0}\big)^*=R^*\colon C^*\to C_0^\perp\leq\bbar{C}_0^*$ it follows that the inverse $J=I^{-1}$ is given for all $(\bs\rho,\mu)\in C_0^*\x C^*$ by the formula 

Now the norm $\|\cdot\|_*$ on $C_0^*\x C^*$ that makes $I$ an isometry is the norm that makes the map $\mathbbm{id}_{C_0^*}\x(R^*)^{-1}\colon C_0^*\x C_0^\perp\to C_0^*\x C^*$ an isometry where $C_0^*\x C_0^\perp$ is equipped with the norm $\|\cdot\|_0$ defined in Lemma~\ref{SplittingLemma}, so that 
$$\|(\bs\rho,\mu)\|_*=\|(\bs\rho, R^*\mu)\|_0=\sup_{F\in\bbar{C}_0}\frac{E(\bs\rho)(F)+\mu(RF)}{\|F\|_{\bbar{C}_0}}$$
and therefore if $E$ is norm-preserving and satifies~\eqref{ForLowInftB} for all $F\in\bbar{C}_0$ then 
$$\max\{\|\bs\rho\|_{C_0^*},\|R^*\mu\|_{\bbar{C}_0^*}\}\leq\|(\bs\rho,\mu)\|_*\leq\|E(\bs\rho)\|_{\bbar{C}_0^*}+\|R^*\mu\|_{\bbar{C}_0^*}=\|\bs\rho\|_{C_0^*}+\|R^*\mu\|_{\bbar{C}_0^*}.$$
If we assume now that $R\colon\bbar{C}_0\to C$ is a contraction then $\|R^*\mu\|_{\bbar{C}_0^*}\leq\|\mu\|_{C^*}$. If also $R$ is such that for all $f\in C$ there exists $F\in R^{-1}(\{f\})$ such that $\|F\|_{\bbar{C}_0}\leq\|f\|_C$ then given $\ee>0$ we can choose $f_\ee\in C$ such that $\|\mu\|_{C^*}\leq\frac{\mu(f_\ee)}{\|f_\ee\|_C}+\ee$ and then by choosing $F_\ee\in R^{-1}(f_\ee)$ with $\|F_\ee\|_{\bbar{C}_0}\leq\|f_\ee\|_C$ we obtain 
$$\|R^*\mu\|_{\bbar{C}_0^*}=\sup_{F\in\bbar{C}_0}\frac{\mu(RF)}{\|F\|_{\bbar{C}_0^*}}\geq\frac{\mu(RF_\ee)}{\|F_\ee\|_{\bbar{C}_0^*}}\geq\frac{\mu(f_\ee)}{\|f_\ee\|_C}\geq\|\mu\|_{C^*}-\ee$$
which since $\ee>0$ is arbitrary shows that $\wt{R}^*\colon C^*\to C_0^\perp\leq\bbar{C}_0^*$ is norm-preserving and completes the proof.$\hfill\Box$\\

We proceed now with the proof of Theorem~\ref{YoungMeasuresTheorem}. We start by proving parts (a) and (b) by applying Lemma~\ref{RecessionSplit} on the recession function $R\colon\bbar{C}_r(\T^d\x\RR_+)\to C(\T^d)$ for which $\ker R=C_r(\T^d\x\RR_+)$ and the extension operator $E\colon\MMM_r(\T^d\x\RR_+)\to\bbar{\MMM}_r(\T^d\x\RR_+)$ defined in~\eqref{ExtensionOperator} so that $\bbar{C}_0=\bbar{C}_r(\T^d\x\RR_+)$, $C_0=C_r(\T^d\x\RR_+)$ and $C=C(\T^d)$. We thus obtain the existence of a unique linear map $I=(I^1,I^2)\colon\bbar{\MMM}_r(\T^d\x\RR_+)\to\MMM_r(\T^d\x\RR_+)\x\MMM(\T^d)$ such that~\eqref{CharProp} holds, whose inverse $J$ is given by~\eqref{IInverse}. By definition the norm $\|\cdot\|_{TV,r}$ defined in~\eqref{MixNormDef} is the norm on $\MMM_r(\T^d\x\RR_+)\x\MMM(\T^d)$ that makes $I$ an isometry and thus in order to check that~\eqref{MixNorm} holds we have to verify that $E$ is norm-preserving and satisfies property~\eqref{ForLowInftB} for every $F\in\bbar{C}_r(\T^d\x\RR_+)$ and that the surjective contraction $R$ satisfies~\eqref{RecessBound}. By the Riesz isomorphism~\eqref{Riesz(r)} and Proposition~\ref{Basic} it follows that the extension operator $E$ is an isometric injection. To check that it satisfies property~\eqref{ForLowInftB} let $F\in\bbar{C}_r(\T^d\x\RR_+)$ and let $\Psi_n\colon\RR_+\to\RR_+$, $n\in\NN$, be given by the formula $\Psi_n(\lambda)=(\lambda-n)^+\lambda^{r-1}$. Then $$F_n:=\Psi_n(\Lambda)\cdot RF(U)=F+(\Psi_n(\Lambda)\cdot RF(U)-F)\in F+C_r(\T^d\x\RR_+),$$
the norms of the maps $F_n$ satisfy 
$$\|F_n\|_{\infty,r}=\sup_{\lambda\geq 0}\frac{(\lambda-n)^+\lambda^{r-1}}{1+\lambda^r}\|RF\|_\infty=\|RF\|_\infty\leq\|F\|_{\infty,r}$$
and by the dominated convergence theorem, for any $\bs\rho\in\MMM_r(\T^d\x\RR_+)$ 
\begin{equation}\label{ForFixed}
\lim_{n\ra+\infty}\int F_n\df\bs{\rho}=\lim_{n\ra+\infty}\int\Psi_n(\lambda)\cdot RF(u)\df\bs{\rho}(u,\lambda)=0.
\end{equation}
Thus $E$ satisfies property~\eqref{ForLowInftB}. Similarly, for any $f\in C(\T^d)$ the map $F:=(1+\Lambda^r)f(U)$ belongs in $R^{-1}(\{f\})$ and $\|F\|_{\infty,r}=\|f\|_\infty$ so that $R$ satisfies~\eqref{RecessBound} and $R^*$ is norm-preserving. This completes the proof of statements (a) and (b) of Theorem~\ref{YoungMeasuresTheorem}.

For the proof of (c) we start by noting that $I^1=j^*$ by~\eqref{IIsomAbstractForm} and thus we only have to obtain the formula~\eqref{I2AsymptForm} for the second coordinate $I^2$. Since $\Lambda^rf(U)\in R^{-1}(\{f\})\subs\bbar{C}_r(\T^d\x\RR_+)$ we have that $I^2(\bs\pi)(f)=\bs\pi(\Lambda^rf(U))-E\circ I^1(\bs\pi)\big(\Lambda^rf(U)\big)$. But since $I^1(\bs\pi)\in\MMM_r(\T^d\x\RR_+)$ is a measure with finite $r$-th moments $(\Lambda\mn M)\Lambda^{r-1}|f(U)|\leq\Lambda^r |f(U)|\in L^1(\bs\rho)$ for all $M>0$ and thus by the dominated convergence theorem 
\begin{equation}\label{RegApprox}
\lim_{M\ra+\infty}\int(\Lambda\mn M)\Lambda^{r-1}f(U)\df I^1(\bs\pi)=\int\Lambda^rf(U)\df I^1(\bs\pi)=E\circ I^1(\bs\pi)(\Lambda^r f(U)).
\end{equation}
But $\bs\pi\big((\Lambda\mn M)\Lambda^{r-1}f(U)\big)=\int(\Lambda\mn M)\Lambda^{r-1}f(U)\df I^1(\bs\pi)$ since $I^1=j^*$ is the restriction operator and $(\Lambda\mn M)\Lambda^{r-1}f(U)\in C_r(\T^d\x\RR_+)$, and therefore by the formula of $I^2(\bs\pi)$ we obtain~\eqref{I2AsymptForm}.

Next we prove (d) i.e~that $I$ is positive. So let $\bs{\pi}\in\bbar{C_1}(\T^d\x\RR_+)$ be positive. Then obviously $I^1(\bs{\pi})$ is positive and we have to show that $I^2(\bs{\pi})$ is also positive. We note that if $f\geq 0$ we have that $(\Lambda\mn M)f(U)\leq\Lambda f(U)$ for every $M>0$ and thus since $(\Lambda\mn M)f(U)\in C_1(\T^d\x\RR_+)$ and $\bs{\pi}$ is assumed positive we have that 
$$\int(\Lambda\mn M)f(U)\df I^1(\bs\pi)=\bs{\pi}(\Lambda\mn M f(U))\leq\bs{\pi}(\Lambda f(U)).$$
Taking the limit as $M\ra+\infty$ it follows by the monotone convergence theorem that 
$$E\circ I^1(\bs\pi)\big(\Lambda f(U)\big)=\int\Lambda f(U)\df I^1(\bs{\pi})\leq\bs{\pi}(\Lambda f(U)).$$
Therefore by~\eqref{CharProp} we have that
$$\int fd I^2(\bs{\pi})=\bs{\pi}\big(\Lambda f(U)\big)-\int\Lambda f(U)d I^1(\bs{\pi})\geq 0,$$
which proves that $I$ is positive. Since $J=I^{-1}$ is obviously positive equality~\eqref{ImageOfPositive} follows.

Finally we prove (e). Since $I^1=j^*$ is the adjoint of a bounded operator it is $w^*$-continuous. We thus have only to prove the positive upper $w^*$-semicontinuity of $I^2$ in~\eqref{IsomSemiCont}. Note that for all non-negative $\bs\pi\in\bbar{\MMM}_{r,+}(\T^d\x\RR_+)$ and $f\in C_+(\T^d)$ it follows by~\eqref{I2AsymptForm} that
$$I^2(\bs\pi)(f)=\inf_{M>0}\bs\pi\big((\Lambda-M)^+\Lambda^{r-1}f(U)\big).$$
Since $(\Lambda-M)^+\Lambda^{r-1}f(U)\in\bbar{C}_r(\T^d\x\RR_+)\leq\bbar{\MMM}_r(\T^d\x\RR_+)^*$ the functional $\ell_{f,M}\colon\bbar{\MMM}_r(\T^d\x\RR_+)\to\RR$ given by $\ell_{f,M}(\bs\pi)=\bs\pi\big((\Lambda-M)^+\Lambda^{r-1}f(U)\big)$ is $w^*$-continuous and therefore the map $\bbar{\MMM}_{r,+}(\T^d\x\RR_+)\ni\bs\pi\mapsto I^2(\bs\pi)(f)$ is upper $w^*$-semicontinuous as the infimum of the $w^*$-continuous linear functionals $\ell_{f,M}$ over $M>0$. This completes the proof of Theorem~\ref{YoungMeasuresTheorem}.$\hfill\Box$

\subsection{Generalized Young path-measures}\label{GYMP}
Our next goal is to lift the results of the previous section to the level of \emph{generalized Young path-measures}, which are elements 
$$\bs{\pi}=(\bs{\pi}_t)_{0\leq t\leq T}\in L_{w^*}^\infty\big(0,T;\bbar{\MMM}_r(\T^d\x\RR_+)\big)\cong L^1\big(0,T;\bbar{C}_r(\T^d\x\RR_+)\big)^*.$$ To do so we will apply Lemma~\ref{RecessionSplit} on the maps induced by the recession operator $R$ and the extension operator $E$. We will use the following notation for the norms of the $L^1$-Bochner spaces and the $L_{w^*}^\infty$-spaces:
\begin{align*}
&\|F\|_{\infty,r;1}:=\|F\|_{L^1(0,T;\bbar{C}_r(\T^d\x\RR_+))},\quad \|f\|_{\infty;1}:=\|f\|_{L^1(0,T;C(\T^d))}\\
&\|\bs\pi\|_{TV,r;\infty}:=\|\bs\pi\|_{L^\infty_{w^*}(0,T;\bbar{\MMM}_r(\T^d\x\RR_+))},\quad\|\mu\|_{TV;\infty}:=\|\mu\|_{L^\infty_{w^*}(0,T;\MMM(\T^d))}.
\end{align*}	
Since $L^1(0,T;C_r(\T^d\x\RR_+))$ is embedded in $L^1(0,T;\bbar{C}_r(\T^d\x\RR_+))$ via the operator $\bar{j}$ induced by the subspace inclusion $j\colon C_r(\T^d\x\RR_+)\to\bbar{C}_r(\T^d\x\RR_+)$ via $\bar{j}(F)(t)=j(F_t)$ for almost all $t\in[0,T]$ we will also use the symbol $\|\cdot\|_{\infty,r;1}$ for the norm of $L^1(0,T;C_r(\T^d\x\RR_+))$.


By Lemma~\ref{InduceE} the extension operator $E\colon\MMM_r(\T^d\x\RR_+)\to\bbar{\MMM}_r(\T^d\x\RR_+)$ is $w^*$-Baire norm preserving injection and by Proposition~\ref{InducedOpProp} it induces a $w^*$-Baire norm-preserving operator $$\bar{E}\colon L_{w^*}^\infty(0,T;\MMM_r(\T^d\x\RR_+))\to L_{w^*}^\infty(0,T;\bbar{\MMM}_r(\T^d\x\RR_+))$$ via the formula $\bar{E}(\bs\rho)(t)=E(\bs\rho_t)$ for almost all $t\in[0,T]$. We will view $L_{w^*}^\infty(0,T;\MMM_r(\T^d\x\RR_+))$ as a subspace of $L_{w^*}^\infty(0,T;\bbar{\MMM}_r(\T^d\x\RR_+))$ via the injection $\bar{E}$ and thus we will also write $\|\cdot\|_{TV,r;\infty}$ for the norm of $L_{w^*}^\infty(0,T;\MMM_r(\T^d\x\RR_+))$. Furthermore, the recession operator $R\colon\bbar{C}_r(\T^d\x\RR_+)\to C(\T^d)$ induces an operator $\bar{R}\colon L^1(0,T;\bbar{C}_r(\T^d\x\RR_+))\to L^1(0,T;C(\T^d))$ on the $L^1$-Bochner spaces. Since the recession operator has a bounded right inverse, for example the map $T\colon C(\T^d)\to\bbar{C}_r(\T^d\x\RR_+)$ given by $T(f)=\Lambda^rf(U)$, it follows that $\bar{R}$ is surjective with bounded right inverse the map $\bar{T}$.
\begin{theorema}\label{YoungPathMeasuresTheorem}(a) We set  
	$$L_{w^*}^\infty(0,T;\MMM_r(\T^d\x\RR_+)\x\MMM(\T^d)):=L_{w^*}^\infty(0,T;\MMM_r(\T^d\x\RR_+))\x L_{w^*}^\infty(0,T;\MMM(\T^d))$$ and on the product space $L_{w^*}^\infty(0,T;\MMM_r(\T^d\x\RR_+)\x\MMM(\T^d))$ we define the norm 
	\begin{equation}\label{LstarMixNorm}
	\|(\bs\rho,\mu)\|_{TV,r;\infty}:=\sup_{F\in L^1(0,T;\bbar{C}_r(\T^d\x\RR_+))\sm\{0\}}\frac{|\lls F,\bar{E}(\bs\rho)\rrs+\lls\bar{R}F,\mu\rrs|}{\|F\|_{\infty,r;1}}.
	\end{equation}
	The norm $\|\cdot\|_{TV,r;\infty}$ defined in~\eqref{LstarMixNorm} satisfies 
	$$\max\{\|\bs\rho\|_{TV,r;\infty},\|\mu\|_{TV;\infty}\}\leq\|(\bs\rho,\mu)\|_{TV,r;\infty}\leq\|\bs\rho\|_{TV,r;\infty}+\|\mu\|_{TV;\infty}$$
	and is thus equivalent to all the product norms on the space $L_{w^*}^\infty(0,T;\MMM_r(\T^d\x\RR_+)\x\MMM(\T^d))$.\\
	\noindent(b) There is a unique isometry 
	\begin{equation}\label{BarIIsom}
	\bar{I}=(\bar{I}^1,\bar{I}^2)\colon L_{w^*}^\infty(0,T;\bbar{\MMM}_r(\T^d\x\RR_+))\to L_{w^*}^\infty(0,T;\MMM_r(\T^d\x\RR_+)\x\MMM(\T^d))
	\end{equation} such that
	\begin{equation}\label{CharPathProp}\lls F,\bs\pi\rrs=\lls F,\bar{E}(\bar{I}^1(\bs\pi))\rrs+\lls\bar{R}F,\bar{I}^2(\bs\pi)\rrs,\quad\mbox{ for all }F\in L^1(0,T;\bbar{C}_r(\T^d\x\RR_+)).
	\end{equation} 
	\noindent(c) The first coordinate $\bar{I}^1$ of $I$ is the restriction operator $\bar{j}^*$, i.e.~the adjoint of the natural inclusion $\bar{j}\colon L^1(0,T;C_r(\T^d\x\RR_+))\hookrightarrow L^1(0,T;\bbar{C}_r(\T^d\x\RR_+))$ and the second coordinate $\bar{I}^2$ is given by the formula
	\begin{equation}\label{I2AsymptPathForm}
	\bar{I}^2(\bs\pi)(f)=\lim_{M\ra+\infty}\bs\pi\big((\Lambda-M)^+\Lambda^{r-1}f(U)\big),\quad f\in L^1(0,T;C(\T^d)).
	\end{equation}
	\noindent(d) The isometry $\bar{I}$ satisfies $\bar{I}(\bs\pi)(t)=I(\bs\pi_t)$ for almost all $t\in[0,T]$ and it is positive i.e.
	\begin{equation}\label{ImageOfPositivePath}\bar{I}\big(L_{w^*}^\infty(0,T;\bbar{\MMM}_{r,+}(\T^d\x\RR_+))\big)=L_{w^*}^\infty\big(0,T;\MMM_{r,+}(\T^d\x\RR_+)\x\MMM_+(\T^d)\big).\end{equation}
	\noindent(e) The restriction of the isometry $\bar{I}$ on $L_{w^*}^\infty(0,T;\bbar{\MMM}_{r,+}(\T^d\x\RR_+))$ is positively $w^*$-semicontinuous in the sense that if the net $\{\bs\pi_\alpha\}_{\alpha\in\A}\subs L_{w^*}^\infty(0,T;\bbar{\MMM}_{r,+}(\T^d\x\RR_+))$ converges to $\bs\pi\in L_{w^*}^\infty(0,T;\bbar{\MMM}_{r,+}(\T^d\x\RR_+))$ in the $w^*$-topology and $F\in L^1(0,T;\bbar{C}_{r,+}(\T^d\x\RR_+))$, $f\in L^1(0,T;C_+(\T^d))$ are non-negative then 
	\begin{equation}\label{IsomSemiContPaths}\liminf_\alpha\bar{I}^1(\bs\pi_\alpha)(F)\geq\bar{I}^1(\bs\pi)(F)\quad\mbox{and}\quad
	\limsup_\alpha\bar{I}^2(\bs\pi_\alpha)(f)\leq\bar{I}^2(\bs\pi)(f).\end{equation}
\end{theorema}\textbf{Proof} The proof is similar to the proof of Theorem~\ref{YoungMeasuresTheorem}. Here we apply Lemma~\ref{RecessionSplit} on the induced recession function $\bar{R}\colon L^1(0,T;\bbar{C}_r(\T^d\x\RR_+))\to L^1(0,T;C(\T^d))$ for which $\ker\bar{R}= L^1(0,T;C_r(\T^d\x\RR_+))$ and the induced extension operator $\bar{E}\colon L_{w^*}^\infty(0,T;\MMM_r(\T^d\x\RR_+))\to L_{w^*}^\infty(0,T;\bbar{\MMM}_r(\T^d\x\RR_+))$ defined in~\eqref{ExtensionOperator}, which is a right inverse to the induced injection $\bar{j}^*=\bbar{j^*}\colon L_{w^*}^\infty(0,T;\bbar{\MMM}_r(\T^d\x\RR_+))\to L_{w^*}^\infty(0,T;\MMM_r(\T^d\x\RR_+))$. We thus obtain the existence of a unique linear map $$\bar{I}=(\bar{I}^1,\bar{I}^2)\colon L_{w^*}^\infty(0,T;\bbar{\MMM}_r(\T^d\x\RR_+))\to L_{w^*}^\infty(0,T;\MMM_r(\T^d\x\RR_+)\x\MMM(\T^d))$$ satisfying~\eqref{CharPathProp}. By definition the norm $\|\cdot\|_{TV,r;\infty}$ defined in~\eqref{LstarMixNorm} is the norm on $L_{w^*}^\infty(0,T;\MMM_r(\T^d\x\RR_+)\x\MMM(\T^d))$ that makes $\bar{I}$ an isometry. Since $\bar{E}$ is norm-preserving in order to check that~\eqref{LstarMixNorm} holds, by Lemma~\ref{InduceE} it suffices to verify that $\bar{E}$ satisfies property~\eqref{ForLowInftB} for every $F\in L^1(0,T;\bbar{C}_r(\T^d\x\RR_+))$ and that the surjective contraction $\bar{R}$ satisfies~\eqref{RecessBound}. To check that $\bar{E}$ satisfies property~\eqref{ForLowInftB} let $F\in L^1(0,T;\bbar{C}_r(\T^d\x\RR_+))$ and let $\Psi_n\colon\RR_+\to\RR_+$, $n\in\NN$, be given by $\Psi_n(\lambda)=(\lambda-n)^+\lambda^{r-1}$. Then $$F_n:=\Psi_n(\Lambda)\cdot\bar{R}F(U)=F+(\Psi_n(\Lambda)\cdot \bar{R}F(U)-F)\in F+L^1(0,T;C_r(\T^d\x\RR_+)),$$
the norms of the maps $F_n$ satisfy 
$$\|F_n\|_{\infty,r;1}=\int_0^T\|F_{n,t}\|_{\infty,r}\df t=\sup_{\lambda\geq 0}\frac{(\lambda-n)^+\lambda^{r-1}}{1+\lambda^r}\int_0^T\|RF_t\|_\infty\df t=\|RF\|_{\infty;1}\leq\|F\|_{\infty,r;1}.$$
By a double application of the dominated convergence theorem 
\begin{equation}\label{ForFixedPath}
\lim_{n\ra+\infty}\int_0^T\int F_{n,t}\df\bs{\rho}_t\df t=\lim_{n\ra+\infty}\int_0^T\int\Psi_n(\lambda)\cdot RF_t(u)\df\bs{\rho}_t(u,\lambda)\df t=0
\end{equation}
for any $\bs\rho\in L_{w^*}^\infty(0,T;\MMM_r(\T^d\x\RR_+))$. Thus $\bar{E}$ satisfies property~\eqref{ForLowInftB}. Similarly, for any $f\in L^1(0,T;C(\T^d))$ the map $F:=(1+\Lambda^r)f(U)$ belongs in $\bar{R}^{-1}(\{f\})$ and $\|F\|_{\infty,r;1}=\|f\|_{\infty;1}$ so that $\bar{R}$ satisfies~\eqref{RecessBound}. This completes the proof of statements (a) and (b) of Theorem~\ref{YoungPathMeasuresTheorem}.

Claim (c) follows similarly to (c) of Theorem~\ref{YoungMeasuresTheorem} by a double application of the dominated convergence theorem. The first claim of (d) follows by~\eqref{CharPathProp} since the Lebesgue differentiation theorem implies that for almost all $t\in[0,T]$, for any $F\in\bbar{C}_r(\T^d\x\RR_+)$ it holds $\ls F,\bs\pi_t\rs=\ls F,E(\bs\rho_{\bs\pi_t})\rs+\ls RF,\mu_{\bs\pi_t}\rs$ and the second claim then follows by claim (d) of Theorem~\ref{YoungMeasuresTheorem}. The proof of (e) is also similar to the proof of (e) of Theorem~\ref{YoungMeasuresTheorem}. 
$\hfill\Box$\\


We will say that $\bs\pi\in L_{w^*}^\infty(0,T;\bbar{\MMM}_r(\T^d\x\RR_+))$ is a \emph{regular Young path-measure} if $\bar{E}\circ(\bar{j})^*(\bs\pi)=\bs\pi$ and \emph{singular} if $(\bar{j})^*(\bs\pi)=0$. We note that $\bar{E}\circ (\bar{j})^*=\bar{E}\circ\bbar{j^*}=\bbar{E\circ j^*}$ by Propositions~\ref{AdjointInduced} and~\ref{InducedOpProp} and thus $\bs\pi=(\bs\pi_t)_{0\leq t\leq T}$ is regular if and only if $\bs\pi_t$ is regular for almost all $t\in[0,T]$. Likewise $\bs\pi$ is singular if and only if $\bs\pi_t$ is singular for almost all $t\in[0,T]$.
\begin{cor}\label{PathMeasDecomp} (a) For any $\bs\pi\in L_{w^*}^\infty(0,T;\bbar{\MMM}_r(\T^d\x\RR_+))$ there exists a uniquely determined decomposition $\bs\pi=\widehat{\bs\pi}+\bs\pi^\perp$ of with $\widehat{\bs\pi}\in L_{w^*}^\infty(0,T;\bbar{\MMM}_r(\T^d\x\RR_+))$ being regular and $\bs\pi^\perp\in L_{w^*}(0,T;\bbar{\MMM}_r(\T^d\x\RR_+))$. Furthermore
	\begin{equation}\label{DecompNormIneqPath}
	\max\{\|\widehat{\bs\pi}\|_{TV,r;\infty},\|\bs\pi^\perp\|_{TV,r;\infty}\}\leq\|\bs\pi\|_{TV,r;\infty}\leq\|\widehat{\bs\pi}\|_{TV,r;\infty}+\|\bs\pi^\perp\|_{TV,r;\infty},\end{equation}
	and $\bs{\pi}\in L_{w^*}^\infty(0,T;\bbar{\MMM}_{r,+}(\T^d\x\RR_+))$ if and only if both $\widehat{\bs\pi}$ and $\bs\pi^\perp$ belong to $L_{w^*}^\infty(0,T;\bbar{\MMM}_{r,+}(\T^d\x\RR_+))$.\\
	\noindent(b) The operators $\widehat{D}_{\infty}$, $D^\perp_{\infty}$ on $L_{w^*}^\infty(0,T;\bbar{\MMM}_r(\T^d\x\RR_+))$ defined by $\widehat{D}_\infty(\bs\pi)=\widehat{\bs\pi}$ and $D^\perp_\infty(\bs\pi)=\bs\pi^\perp$ coincide with the operators $\bbar{\widehat{D}}$ and $\bbar{D^\perp}$ induced on the space $L_{w^*}^\infty(0,T;\bbar{\MMM}_r(\T^d\x\RR_+))$ by the maps $\widehat{D}$ and $D^\perp$, i.e.
	$$\widehat{D}_\infty(\bs\pi)=\big(\widehat{D}(\bs\pi_t)\big)_{t\in[0,T]},\;D^{\perp}_\infty(\bs\pi)=\big(D^\perp(\bs\pi_t)\big)_{t\in[0,T]}\quad\mbox{in }L_{w^*}^\infty(0,T;\bbar{\MMM}_r(\T^d\x\RR_+)).$$
	Consequently the maps $\widehat{D}_\infty$ and $D^\perp_\infty$ are pointwise $w^*$-limits of $w^*$-continuous operators. Thus they are $w^*$-Baire, and thus $w^*$-measurable.\\	
	(c) The restriction of $\widehat{D}_\infty$ on $L_{w^*}^\infty(0,T;\bbar{\MMM}_{r,+}(\T^d\x\RR_+))$ is positively $w^*$-lower semicontinuous and the restriction $D^\perp_\infty$ on $L_{w^*}^\infty(0,T;\bbar{\MMM}_{r,+}(\T^d\x\RR_+))$ is positively $w^*$-upper semicontinuous, i.e.~for any map $F\in\bbar{C}_{r,+}(\T^d\x\RR_+)$ and any net $\{\bs\pi_\alpha\}_{\alpha\in\A}\subs\bbar{\MMM}_{r,+}(\T^d\x\RR_+)$ converging to some $\bs\pi\in\bbar{\MMM}_{r,+}(\T^d\x\RR_+)$ in the $w^*$-topology
	\begin{equation}\label{DecompSemiContPath}\liminf_{\alpha}\widehat{\bs\pi_\alpha}(F)\geq\widehat{\bs\pi}(F)\quad\mbox{and}\quad\limsup_{\alpha}\bs\pi_\alpha^\perp(F)\leq\bs\pi^\perp(F).
	\end{equation}
\end{cor}\textbf{Proof} The proof of (a) and (c) is similar to the proof of Corollary~\ref{YMDecomp} where here we define $\widehat{\bs\pi}:=\bar{E}\circ(\bar{j})^*(\bs\pi)$ and $\bs\pi^\perp:=(\bar{R})^*\circ\bar{I}_2$ with $\bar{I}=(\bar{I}_1,\bar{I}_2)$ being the isometry of Theorem~\ref{YoungPathMeasuresTheorem}. For the proof of (b) we note that by Propositions~\ref{AdjointInduced} and~\ref{InducedOpProp} we have that $\widehat{D}_\infty=\bar{E}\circ(\bar{j})^*=\bbar{E\circ j^*}=\bbar{\widehat{D}}$ and $$D^\perp_\infty=\mathbbm{id}_{L_{w^*}^\infty(0,T;\bbar{\MMM}_r(\T^d\x\RR_+))}-\widehat{D}_\infty=\bbar{\mathbbm{id}_{\bbar{\MMM}_r(\T^d\x\RR_+)}}-\bbar{\widehat{D}}=\bbar{\mathbbm{id}_{\bbar{\MMM}_r(\T^d\x\RR_+)}-\widehat{D}}=\bbar{D^\perp}.$$ The fact that $\widehat{D}_\infty=\bbar{\widehat{D}}$ and $D^\perp_\infty=\bbar{D^\perp}$ implies by Proposition~\ref{InducedOpProp} that the maps $\widehat{D}_\infty$, $D^\perp_\infty$ are $w^*$-analytically measurable and thus they are also $w^*$-measurable. In fact, by Lemma~\ref{InduceE} we have that $E$ is the pointwise $w^*$-limit of the sequence of $w^*$-continuous operators $\Pi_M^*\colon\MMM_r(\T^d\x\RR_+)\to\bbar{\MMM}_t(\T^d\x\RR_+)$ where $\Pi_M\colon\bbar{C}_r(\T^d\x\RR_+)\to C_r(\T^d\x\RR_+)$ is the operator defined by $\Pi_M(F)(u,\lambda)=F(u,\lambda\mn M)$. Thus by Proposition~\ref{LinftInducedwStarOpCont} the operator $\widehat{D}_\infty$ is the pointwise $w^*$-limit of the sequence of the $w^*$-continuous operators $\bbar{\Pi_M^*\circ j^*}=(j\circ\Pi_M)^*\colon L_{w^*}^\infty(0,T;\MMM_r(\T^d\x\RR_+))\to L_{w^*}(0,T;\bbar{\MMM}_r(\T^d\x\RR_+))$, $M\in\NN$. Likewise, by the formula of $I^2$ we have that $D^\perp$ is the pointwise $w^*$-limit of the sequence of $w^*$-continuous operators $R^*\circ T_M^*=(T_M\circ R)^*\colon\bbar{\MMM}_r(\T^d\x\RR_+)\to\bbar{\MMM}_r(\T^d\x\RR_+)$ where here $T_M\colon C(\T^d)\to\bbar{C}_r(\T^d\x\RR_+)$ is the linear operator $T_Mf=(\Lambda-M)^+\Lambda^{r-1}f(U)$ and thus by Proposition~\ref{LinftInducedwStarOpCont} the operator $D^\perp_\infty$ is the pointwise $w^*$-limit of the sequence of operators $\bbar{R^*\circ T_M^*}$. $\hfill\Box$\\

The $M$-modified micro empirical density $\bs\pi^{N,\ell;M}$ of the ZRP defined in~\eqref{MModMicrEmpDens} is decomposed as the sum $\bs\pi^{N,\ell;M}=\widehat{\bs\pi}^{N,\ell;M}+\bs\pi^{\perp,N,\ell;M}$ of regular and singular generalized Young-functional valued process. Here $$\widehat{\bs\pi}^{N,\ell;M}\colon D(0,T;\MM_N^d)\to L_{w^*}^\infty(0,T;\bbar{\MMM}_1(\T^d\x\RR_+))$$ is the \emph{$M$-truncated empirical distribution of the ZRP} defined via duality by 
\begin{equation}\label{TruncMicroEmp}
\lls F,\widehat{\bs\pi}^{N,\ell;M}\rrs=\int_0^T\fr{N^d}\sum_{x\in\T^d_N}F_t\Big(\frac{x}{N},\eta_t^\ell(x)\mn M\Big)\df t,\quad F\in L^1(0,T;\bbar{C}_1(\T^d\x\RR_+))
\end{equation}
and $\bs\pi^{\perp,N,\ell;M}=\bar{R}^*\circ\rho^{\perp,N,\ell;M}$ where $\rho^{N,\ell;M,\perp}\colon D(0,T;\MM_N^d)\to L_{w^*}^\infty(0,T;\MMM_+(\T^d))$ is the $M$-excess empirical density defined by  
\[\lls f,\rho^{\perp,N,\ell;M}\rrs=\int_0^T\fr{N^d}\sum_{x\in\T_N^d}f_t\Big(\frac{x}{N}\Big)(\eta_t^\ell(x)-M)^+\df t\]
and $R^*$ is the adjoint of the recession operator $R\colon L^1(0,T;\bbar{C}_1(\T^d\x\RR_+))\to L^1(0,T;C_1(\T^d\x\RR_+))$.


\begin{prop}\label{PositiveRepresentingParts} The sets
	\begin{align*}
	&V_1=L_{w^*}^\infty(0,T;\bbar{\MMM}_{r,+}(\T^d\x\RR_+))\\
	&V_2=L_{w^*}^\infty(0,T;\bbar{\MMM}_{r,\mathfrak{m}}(\T^d\x\RR_+)),\;\mathfrak{m}>0\\
	&V_3=L_{w^*}^\infty(0,T;\bbar{\PP}_r(\T^d\x\RR_+))\\
	&V_4=L_{w^*}^\infty(0,T;T_m\bbar{\PP}_r(\T^d\x\RR_+)),\;m\in\MMM_+(\T^d)
	\end{align*}	 
 are ($w^*$-)closed subspaces of $L^\infty_{w^*}(0,T;\bbar{\MMM}_r(\T^d\x\RR_+))$. In particular $L_{w^*}^\infty(0,T;\Y_{1,\mathfrak{m}}(\T^d))$ is a $w^*$-closed subspace of $L^\infty_{w^*}(0,T;\bbar{\MMM}_r(\T^d\x\RR_+))$.
\end{prop}\textbf{Proof} The space $\bbar{\MMM}_{r,+}(\T^d\x\RR_+)$ is a $w^*$-closed subspace of $\bbar{\MMM}_r(\T^d\x\RR_+)$ as the dual cone in the sense of~\eqref{PosCone} of the positive cone $\bbar{C}_{r,+}(\T^d\x\RR_+)$ and thus by Proposition~\ref{PosCon} it follows that $L_{w^*}^\infty(0,T;\MMM_{r,+}(\T^d\x\RR_+))$ is a closed subspace of $L_{w^*}^\infty(0,T;\MMM_r(\T^d\x\RR_+))$. 

We check next that $V_2$ is a closed subspace. So let $\{\bs\pi_\alpha\}_{\alpha\in\A}\subs V_2$ be a net converging to some $\bs\pi\in L^\infty_{w^*}(0,T;\bbar{\MMM}_{r,+}(\T^d\x\RR_+))$. Then for all measurable $E\subs[0,T]$ with strictly positive Lebesgue measure 
$$\mathfrak{m}=\fr{\LL_{\T^d}(E)}\int_E\ls\Lambda,\bs\pi_{\alpha,s}\rs\df s=\Big\lls\fr{\LL_{\T^d}(E)}\1_{E}\Lambda,\bs\pi_\alpha\Big\rrs$$
and thus since $\{\bs\pi_\alpha\}_{\alpha\in\A}$ converges to $\bs\pi$,
$$\fr{\LL_{\T^d}(E)}\int_E\ls\Lambda,\bs\pi_s\rs\df s=\Big\lls\fr{\LL_{\T^d}(E)}\1_E\Lambda,\bs\pi\Big\rrs=\lim_\alpha\Big\lls\fr{\LL_{\T^d}(E)}\1_E\Lambda,\bs\pi_\alpha\Big\rrs=\mathfrak{m}.$$
Since this holds for all measurable $E\subs[0,T]$ with $\LL_{\T^d}(E)>0$ it follows that $\bs\pi_t(\Lambda)=\mathfrak{m}$ for almost all $t\in[0,T]$, i.e.~$\bs\pi\in L_{w^*}^\infty(0,T;\bbar{\MMM}_{r,\mathfrak{m}}(\T^d\x\RR_+))$. The fact that $V_3$ is also closed follows similarly.

We prove finally that $V_4$ is closed. The natural projection $U\colon\T^d\x\RR_+\to\T^d$ induces the pull-back $U^\sharp\colon C(\T^d)\to C_1(\T^d\x\RR_+)$ via $U^\sharp(f)=f(U)$, which is obviously a linear contraction. Its adjoint $(U^\sharp)^*\colon\MMM_1(\T^d\x\RR_+)\to\MMM(\T^d)$ is exactly the push-forward operator $U_\sharp$ of measures. Thus $U_\sharp\colon\MMM_1(\T^d\x\RR_+)\to\MMM(\T^d)$ is $w^*$-continuous and induces by Corollary~\ref{AdjointInduced} the linear operator 
$$(\bbar{U^\sharp})^*\equiv\bbar{U_\sharp}\colon L_{w^*}^\infty(0,T;\MMM_1(\T^d\x\RR_+))\to L_{w^*}^\infty(0,T;\MMM(\T^d)).$$
This is $w^*$-continuous as the adjoint of the induced operator $\bbar{U^\sharp}\colon L^1(0,T;C(\T^d))\to L^1(0,T;C_1(\T^d\x\RR_+))$ on the Bochner $L^1$-spaces. As we have seen the natural injection $j\colon C_1(\T^d\x\RR_+)\hookrightarrow \bbar{C}_1(\T^d\x\RR_+)$ induces the projection $\bbar{j^*}\equiv(\bar{j})^*\colon L_{w^*}^\infty(0,T;\bbar{\MMM}_1(\T^d\x\RR_+))\to L_{w^*}^\infty(0,T;\MMM_1(\T^d\x\RR_+))$ which is also $w^*$-continuous. Thus the composition $$\bbar{U_\sharp}\circ\bbar{j^*}\colon L_{w^*}^\infty(0,T;\bbar{\PP}_1(\T^d\x\RR_+))\to L_{w^*}^\infty(0,T;\MMM_+(\T^d))$$
is $w^*$-continuous and 
$$L_{w^*}^\infty(0,T;T_m\bbar{\PP}_1(\T^d\x\RR_+))=\big(\bbar{U_\sharp}\circ\bbar{j^*}\big)^{-1}(\{c_m\})$$
where $c_m\in L_{w^*}^\infty(0,T;\MMM_+(\T^d))$ is the path almost everywhere equal to $m\in\MMM_+(\T^d)$. Therefore $L_{w^*}^\infty(0,T;T_m\bbar{\PP}_1(\T^d\x\RR_+))$ is closed in $L_{w^*}^\infty(0,T;\bbar{\PP}_1(\T^d\x\RR_+))$ as the inverse image of a closed set via a continuous map. Since $L_{w^*}^\infty(0,T;\bbar{\PP}_1(\T^d\x\RR_+))$ is closed in $L_{w^*}(0,T;\bbar{\MMM}_1(\T^d\x\RR_+))$ it follows that $L_{w^*}^\infty(0,T;T_m\bbar{\PP}_1(\T^d\x\RR_+))$ is also closed in $L_{w^*}^\infty(0,T;\bbar{\MMM}_1(\T^d\x\RR_+))$. Therefore $V_4$ is closed since 
$V_4=\bbar{I}\big(L_{w^*}^\infty(0,T;T_m\bbar{\PP}_1(\T^d\x\RR_+))\big)$ and $\bar{I}$ is a homeomorphism. Finally $$L_{w^*}^\infty(0,T;\bbar{\Y}_{1,\mathfrak{m}}(\T^d))=L_{w^*}^\infty(0,T;\bbar{\PP}_{1,\mathfrak{m}}(\T^d\x\RR_+))\cap L_{w^*}^\infty(0,T;\bbar{\Y}_1(\T^d))$$ 
is $w^*$-closed as the intersection of $w^*$-closed subspaces.$\hfill\Box$\\


We close this section noting that since $\bbar{C_r}(\T^d\x\RR_+)$ is separable the space $L^1(0,T;\bbar{C}_r(\T^d\x\RR_+))$ is also separable by Proposition~\ref{ToApplyA2} and thus the results on topological measure theory of Section~\ref{Appendix1} apply for the space of probability measures on the dual space $L_{w^*}^\infty(0,T;\bbar{\MMM}_r(\T^d\x\RR_+))$ equipped with its $w^*$-topology.

\subsection{Barycentric projection}\label{Barycentric}
We consider first the barycentric projection $B_1\colon\MMM_1(\T^d\x\RR_+)\to\MMM(\T^d)$ ordinary Young measures defined by duality via
$$\ls f,B_1(\bs\rho)\rs=\ls f(U)\Lambda,\bs\rho\rs,\quad\forall f\in C(\T^d).$$
Each functional $B_1(\bs\rho)\in\MMM(\T^d)$ thus defined is indeed bounded with $\|B_1(\bs\rho)\|_{TV}\leq\int\Lambda\df|\bs\rho|$ and the map $B_1$ is a linear contraction. Note that $B_1$ is not $w^*$-continuous with respect to the topology of $\MMM_1(\T^d\x\RR_+)$.

According to the disintegration theorem~\cite[Theorem 5.3.1]{AGS2000a} for each $m\in\MMM_+(\T^d)$ and each $m$-Young measure $\bs\rho\in T_m\MMM_{1,+}(\T^d\x\RR_+)$ there exists a unique $m$-almost everywhere defined $w^*$-measurable family $(\bs\rho^u)_{u\in\T^d}\subs\PP_1\RR_+$ of probability measures such that
$$\bs\rho(F)=\iint F(u,\lambda)\df\bs\rho^u(\lambda)\df m(u),\quad\forall F\in B_1(\T^d\x\RR_+).$$  
Therefore if we define the \emph{barycentric density} $b(\bs\rho)$ of a Young-measure $\bs\rho\in\MMM_{1,+}(\T^d\x\RR_+)$ via its disintegration $(\bs\rho^u)_{u\in\T^d}\subs\PP_1\RR_+$ as the map $b(\bs\rho)\in L^1_+(U_\sharp\bs\rho)$ given by 
$$b(\bs\rho)(u)=\int\lambda\df\bs\rho^u(\lambda),\quad \mbox{for }U_\sharp\bs\rho\mbox{-almost all }u\in\T^d$$ then 
\begin{align*}
\ls f,B_1(\bs\rho)\rs=\int_{\T^d}f(u)\int\lambda\df\bs\rho^u(\lambda)\df U_\sharp\bs\rho(u)=\int_{\T^d}f(u)b(\bs\rho)(u)\df U_\sharp\bs\rho(u)
\end{align*}
and therefore $B_1(\bs\rho)\ll U_\sharp\bs\rho$ with $\frac{\df B_1(\bs\rho)}{\df U_\sharp\bs\rho}=b(\bs\rho)\in L^1_+(U_\sharp\bs\rho)$. The barycentric projection defined above yields for each $m\in\MMM_+(\T^d)$ a surjective linear operator $$B_1|_{T_m}\colon T_m\MMM_{1,+}(\T^d\x\RR_+)\to L^1_+(m)\df m:=\{f\df m|f\in L^1_+(m)\}\leq\MMM_+(\T^d).$$ A right inverse for $B_1|_{T_m}$ is given by $L^1_+(m)\df m\ni\rho\df m\mapsto\int\delta_u\otimes\delta_{\rho(u)}\df m(u)\in T_m\MMM_{1,+}(\T^d\x\RR_+)$.

The barycentric projection $B_1$ can be extended on the domain $\bbar{\MMM}_1(\T^d\x\RR_+)$ to a barycentric projection $B\colon\bbar{\MMM}_1(\T^d\x\RR_+)\to\MMM(\T^d)$. Namely, for each $\bs{\pi}\in\bbar{\MMM}_{1}(\T^d\x\RR_+)$ the \emph{barycentric projection} $\pi\in\MMM(\T^d)$ of $\bs\pi$ is the measure $\pi:=B(\bs{\pi})$ defined by
$$\pi(f):=B(\bs\pi)(f)=\bs{\pi}\big(\Lambda f(U)\big),\quad \forall\;f\in C(\T^d).$$
The map $B$ is the adjoint of the bounded linear injective contraction 
$$C(\T^d)\ni f\mapsto \Lambda f(U)\in\bbar{C}_1(\T^d\x\RR_+)$$
and as such it is surjective, bounded and $w^*$-continuous. Since $E$ is the $w^*$-pointwise limit of the operators $\Pi_M^*$ where $\Pi_M\colon \bbar{C}_1(\T^d\x\RR_+)\to C_1(\T^d\x\RR_+)$ is given by $\Pi_MF=F(U,\Lambda\mn M)$ and $B$ is $w^*$-continuous it follows that $B_1=B\circ E$ is the limit of the $w^*$-continuous operators $B_{1,M}=B\circ\Pi_M^*\colon\bbar{\MMM}_1(\T^d\x\RR_+)\to\MMM(\T^d)$ and thus it is $w^*$-Baire.

The restriction $B\colon\bbar{\MMM}_{1,+}(\T^d\x\RR_+)\to\MMM_+(\T^d)$ of the barycentric projection on non-negative Young measures remains surjective for the target space $\MMM_+(\T^d)$ when further restricted on $\bbar{\Y}_1(\T^d)$. Indeed, if $\mu\in\MMM_+(\T^d)$ has the Radon-Nikodym representation $\mu=\rho+\rho^\perp$ with $\rho\ll\LL_{\T^d}$ and $\rho^\perp\perp\LL_{\T^d}$ then for the generalized Young measure $\bs\pi:=E(\bs\rho)+R^*(\rho^\perp)$ where $\bs\rho=\int\delta_u\otimes\delta_{\rho(u)}\df u\in\PP_1(\T^d\x\RR_+)$ it holds that $B(\bs\pi)=\mu$. Indeed, by definition $\bs\pi(F)=\int F(u,\rho(u))\df u+\int RF(u)\df\rho^\perp(u)$ for all $F\in\bbar{C}_1(\T^d\x\RR_+)$ and therefore 
$$B(\bs\pi)(f)=\bs\pi(\Lambda f(U))=\int f(u)\rho(u)\df u+\int f(u)\df\rho^\perp(u)=\int f\df\mu.$$ 

Since $B_1=B\circ E$ it follows that $B\circ\widehat{D}=B\circ E\circ j^*= B_1\circ j^*$ and for any $\rho^\perp\in\MMM(\T^d)$ we have that 
\[B(R^*(\mu))(f)=R^*(\mu)(\Lambda f(U))=\mu\big(R(\Lambda f(U)\big)=\mu(f)\]
so that $B\circ R^*=\mathbbm{id}_{\MMM(\T^d)}$. Thus if $\bs\pi\in \bbar{\Y}_1(\T^d)$ and the decomposition $\bs\pi=\widehat{\bs\pi}+\bs\pi^\perp$ is represented by the pair $(\bs\rho_{\bs\pi},\rho_{\bs\pi}^\perp)\in\Y_1(\T^d)\x\MMM_+(\T^d)$, i.e.~$\widehat{\bs\pi}=E(\bs\rho_{\bs\pi})$ and $\bs\pi^\perp=R^*(\rho_{\bs\pi}^\perp)$ then the barycentric projection $B$ is given by $$B(\bs\pi)=B(\widehat{\bs\pi})+B(\bs\pi^\perp)=B_1(\bs\rho)+B(R^*(\rho^\perp))=b(\bs\rho)\df\LL_{\T^d}+\rho^\perp.$$
Since the measure $\rho_{\bs\pi}^\perp$ representing the singular part $\bs\pi^\perp$ of $\bs\pi$ is the measure $I^2(\bs\pi)$ where $I^2$ is the second coordinate of the isometry $I$ of Theorem~\ref{YoungMeasuresTheorem} we have the functional relation 
\begin{equation}\label{SingularRepresentativeBiaBerProj}
B\circ D^\perp(\bs\pi)=\rho^\perp_{\bs\pi}=I^2(\bs\pi),\quad\bs\pi\in\MMM_1(\T^d\x\RR_+),
\end{equation}
which together with the identity $B\circ R^*=\mathbbm{id}_{\MMM(\T^d)}$ shows that the restriction of $B$ on $\Ima D^\perp=D^\perp(\MMM(\T^d))$ is invertible with inverse the adjoint $R^*$ of the recession operator.

Since the barycentric projection $B\colon\bbar{\MMM}_1(\T^d\x\RR_+)\to\MMM(\T^d)$ is a $w^*$-continuous linear operator it induces a barycentric projection $\bar{B}\colon L_{w^*}^\infty(0,T;\bbar{\MMM}_1(\T^d\x\RR_+))\to L_{w^*}^\infty(0,T;\MMM(\T^d)))$ on the respective $L_{w^*}^\infty$-spaces such that $\bar{B}(\bs\pi)(t)=B(\bs\pi_t)$ for all $\bs\pi\in L_{w^*}^\infty(0,T;\bbar{\MMM}_1(\T^d\x\RR_+))$. This is of course $w^*$-continuous and restricts to a surjection $\bar{B}\colon L_{w^*}(0,T;\bbar{\Y}_1(\T^d))\to L_{w^*}^\infty(0,T;\MMM_+(\T^d))$. Since $B_1$ is bounded and $w^*$-analytic it follows that $B_1$ also induces an operator $\bar{B}_1$ on the corresponding $L_{w^*}^\infty$-spaces and $\bar{B}_1=\bar{B}\circ\bar{E}$.

More generally, for any $\Psi\in\bbar{C}_1(\RR_+)$ i.e.~such $\Psi\colon\RR_+\to\RR$ is continuous and the limit $\Psi'(\infty):=\lim_{\lambda\ra+\infty}\frac{\Psi(\lambda)}{\lambda}$ exists we will consider the projection $B_\Psi\colon\bbar{\MMM}_1(\T^d\x\RR_+)\to\MMM(\T^d)$ given by 
$$B_\Psi(\bs\pi)(f)=\bs\pi\big(\Psi(\Lambda)f(U)\big),\quad f\in C(\T^d)$$
and set $B_{1,\Psi}=B_\Psi\circ E$ the restriction of $B_\Psi$ on $\MMM_1(\T^d\x\RR_+)$ via the injection $E\colon\MMM_1(\T^d\x\RR_+)\to\bbar{\MMM}_1(\T^d\x\RR_+)$. Of course then $B=B_\Psi$ and $B_1=B_{1,\Psi}$ for $\Psi=\mathbbm{id}_{\RR_+}$. The map $B_\Psi$ is $w^*$-continuous for all $\Psi\in\bbar{C}_1(\RR_+)$ and $B_{1,\Psi}$ is the pointwise $w^*$-limit as $M\ra+\infty$ of the $w^*$-continuous operators $B_\Psi\circ\Pi_M^*$, where $\Pi_M$ is the operator defined in the proof of Lemma~\ref{InduceE} and thus it is $w^*$-Baire measurable by Proposition~\ref{PointLimOp}. Furthermore, if $\Psi\in C_1(\RR_+)$ i.e.~if $\Psi_\infty=0$, then $B_{1,\Psi}$ is also $w^*$-continuous.

 For any $\bs\rho\in\MMM_1(\T^d\x\RR_+)$ 
$$B_{1,\Psi}(\bs\rho)(f)=\bs\rho\big(\Psi(\Lambda)f(U)\big)=\int_{\T^d}f(u)\int_{\RR_+}\Psi(\lambda)\df\bs\rho^u(\lambda)\df U_\sharp\bs\rho(u)$$
where $(\bs\rho^u)_{u\in\T^d}$ is the $U_\sharp\bs\rho$-a.s.~defined disintegration of $\bs\rho$. Thus if for each $m\in\MMM_+(\T^d)$ we define the $\Psi$-density map $b_\Psi\colon T_m\MMM_{1,+}(\T^d\x\RR_+)\to L_+^1(m)\df m\leq\MMM(\T^d)$ via 
\begin{equation}\label{MapsDens}
b_\Psi(\bs\rho)(u)\equiv\Psi(\bs\rho)(u):=\int_{\RR_+}\Psi(\lambda)\df\bs\rho^u(\lambda),\quad m\mbox{-a.s.~}\forall u\in\T^d
\end{equation}
then $B_{1,\Psi}(\bs\rho)=b_\Psi(\bs\rho)\df m$. Also for any $\mu\in\MMM(\T^d)$ 
$$B_\Psi(R^*(\mu))(f)=R^*(\mu)\big(\Psi(\Lambda)f(U)\big)=\mu\big(R(\Psi(\Lambda)f(U))\big)=\Psi_\infty\mu(f)=\Psi_\infty\cdot B(R^*(\mu))$$
so that $B_\Psi\circ R^*=\Psi_\infty\cdot\mathbbm{id}_{\MMM(\T^d)}$. Therefore for all $\bs\pi\in T_m\bbar{\MMM}_1(\T^d\x\RR_+)$, if $\widehat{\bs\pi}=E(\bs\rho)$ and $\bs\pi^\perp=R^*(\rho^\perp)$, then 
\begin{equation}\label{PsiBaryc}
B_\Psi(\bs\pi)=B_\Psi\circ\widehat{D}(\bs\pi)+B_\Psi\circ D^\perp(\bs\pi)=b_\Psi(\bs\rho)\df m+\Psi'(\infty)\rho^\perp,
\end{equation}
which yields the functional relation~$B_\Psi=B_\Psi\circ\widehat{D}+\Psi'(\infty)\cdot B\circ D^\perp$.
\section{Proofs} \label{ProofS}
\subsection{Relative compactness of the empirical Young measures}\label{RelCompYoungM}

In this section we collect the results on the relative compactness of the laws of the micro and macro empirical density of the ZRP and their basic properties.

\begin{prop}\label{YoungRelComp} Let $P^N$ be the law of the diffusively rescaled ZRP on the Skorohod space $D(0,T;\MM_N^d)$, starting from a sequence $\{\mu_0^N\in\PP_1\MM_N^d\}_{N\in\NN}$ of initial distributions with total mass $\mathfrak{m}>0$ in probability. Then the family 
	\begin{equation}\label{MicroEmpLaw}
	\bs{Q}^{N,\ell}:=\bs\pi^{N,\ell}_\sharp P^N\in\PP L^\infty_{w^*}(0,T;\bbar{\PP}_1(\T^d\x\RR_+)),\quad N\in\NN,\;\ell\in\ZZ_+
	\end{equation}
	of the laws of the micro empirical density process of the ZRP is sequentially relatively compact. 
\end{prop}\textbf{Proof} Since $L_{w^*}^\infty(0,T;\bbar{\PP}_1(\T^d\x\RR_+))$ is a $w^*$-closed subspace of $L_{w^*}^\infty(0,T;\bbar{\MMM_1}(\T^d\x\RR_+))$ and the latter space is a completely regular submetrizable topological space, by the Prokhorov-Le Cam Theorem~\ref{ProkhorovLeCam} it suffices to show that $\{\bs{Q}^{N,\ell}\}$ is uniformly tight. By the Banach-Alaoglu theorem the ball $B_r:=\{\|\cdot\|_{TV,1;\infty}\leq r\}$ is $w^*$-compact in $L^\infty(\bbar{\MMM_1}(\T^d\x\RR_+))$ and therefore it suffices to show that 
$$\sup_{(N,\ell)\in\NN\x\ZZ_+}\bs{Q}^{N,\ell}(B_r^c)=\sup_{(N,\ell)\in\NN\x\ZZ_+}P^N\{\|\bs\pi^{N,\ell}\|_{TV,1;\infty}>r\}\stackrel{r\ra+\infty}\lra 0.$$
For any $F\in L^1(0,T;\bbar{C}_1(\T^d\x\RR_+))$
$$\lls F,\bs\pi^{N,\ell}\rrs\leq\int_0^T\frac{\|F_t\|_{\infty,1}}{N^d}\sum_{x\in\T_N^d}\big(1+\eta_t^\ell(x)\big)\df t
=\int_0^T\|F_t\|_{\infty,1}\big(1+\ls 1,\pi_t^N\rs\big)\df t$$
and therefore by the conservation of the total number of particles 
\begin{align}\label{TotalMassConservationBound}
\|\bs\pi^{N,\ell}\|_{TV,1;\infty}&=\sup_{\|F\|_{L^1(0,T;\bbar{C}_1(\T^d\x\RR_+))}\leq 1}\lls F,\bs\pi^{N,\ell}\rrs\nonumber\\
&\stackrel{P^N\mbox{-a.s.}}=\big(1+\ls 1,\pi_0^N\rs\big)\sup_{\|F\|_{L^1(0,T;\bbar{C}_1(\T^d\x\RR_+))}\leq 1}\int_0^T\|F_t\|_{\infty,1}\df t\nonumber\\
&\leq 1+\ls 1,\pi_0^N\rs
\end{align}
But by the bound~\eqref{TotalMassConservationBound} we have that for all $(N,\ell)\in\NN\x\ZZ_+$
$$P^N\{\|\bs\pi^{N,\ell}\|_{TV,1;\infty}>r\}\leq P^N\big\{1+\ls 1,\pi_0^N\rs>r\big\}=\mu_0^N\big\{\ls 1,\pi^N\rs>r-1\big\}$$
and therefore 
\begin{equation}\label{AlsoForEps}
\lim_{r\ra+\infty}\sup_{(N,\ell)\in\NN\x\ZZ_+}P^N\{\|\bs\pi^{N,\ell}\|_{TV,1;\infty}>r\}\leq\lim_{r\ra+\infty}\sup_{N\in\NN}\mu_0^N\big\{\ls 1,\pi^N\rs>r-1\big\}=0\end{equation}
where the last limit holds by the following Lemma~\ref{BoundedTotParticlInMacrLim} and the assumptions on the sequence $\{\mu_0^N\}$ of initial distributions.$\hfill\Box$\\

Since $\{\bs{Q}^{N,[N\ee]}|N\in\NN,\;\ee>0\}\subs\{\bs{Q}_{N,\ell}|N\in\NN,\;\ell\in\ZZ_+\}$ 
it is evident by~\eqref{AlsoForEps} that the family $\bs{Q}^{N,\ee}:=\bs{Q}^{N,[N\ee]}$, $N\in\NN$, $\ee>0$, is also relatively compact in $\PP L^\infty_{w^*}(0,T;\bbar{\PP}_1(\T^d\x\RR_+))$.
\begin{lemma}\label{BoundedTotParticlInMacrLim} Let $\{\mu_0^N\in\PP\MM_N^d\}$ be a sequence of initial distributions. If $\{\mu_0^N\}$ satisfies the $O(N^d)$-entropy assumption then 
	\begin{equation*}
	\sup_{N\in\NN}\int\ls\pi^N,1\rs d\mu_0^N<+\infty.
	\end{equation*}
If $\{\mu_0^N\}$ satisfies either the $O(N^d)$-entropy assumption or is associated to a macroscopic profile $\mu_0\in\MMM_+(\T^d)$ then 
\begin{equation}\label{TotalMassBound}
\lim_{M\ra+\infty}\sup_{N\in\NN}\mu_0^N\big\{\ls 1,\pi^N\rs>M\big\}=0.
\end{equation}
\end{lemma}\textbf{Proof} By the relative entropy inequality we have that
$$\int\ls\pi^N,1\rs d\mu_0^N\leq\fr{\theta N^d}\Big\{\log\int e^{\theta N^d\ls\pi^N,1\rs}d\nu_{\rho_*}^N+H(\mu_0^N|\nu_{\rho_*}^N)\Big\}$$ 
for all $\theta>0$ and all $N\in\NN$. But 
$$\int e^{\theta N^d\ls\pi^N,1\rs}d\nu_{\rho_*}=\int\prod_{x\in\T_N^d} e^{\theta\eta(x)}d\nu_{\rho_*}=M_{\nu^1_{\rho_*}}(\theta)^{N^d}$$
and therefore 
$$\int\ls\pi^N,1\rs d\mu_0^N\leq\fr{\theta}\Big\{\Lambda_{\rho_*}(\theta)+\fr{N^d}H(\mu_0^N|\nu_{\rho_*}^N)\Big\}$$ 
for all $\theta>0$ and all $N\in\NN$. It follows that 
$$\limsup_{N\ra+\infty}\int\ls\pi^N,1\rs d\mu_0^N\leq\frac{\Lambda_{\rho_*}(\theta)+K_*}{\theta}$$ for all $\theta>0$. But $\rho_*<\rho_c$ and thus $\nu^1_{\rho_*}$ has exponential moments, and therefore by choosing $\theta_*\in\DD_{\Lambda_{\rho_*}}\sm\{0\}$ in the inequality above we obtain~\eqref{BoundedTotParticlInMacrLim}.

If $\{\mu_0^N\}$ satisfies the $O(N^d)$-entropy assumption then~\eqref{TotalMassBound} follows by Chebyshev's inequality. In fact~\eqref{TotalMassBound} holds without the $O(N^d)$-entropy assumption as long as $\{\mu_0^N\}\subs\PP_1\MM_N^d$ has bounded total mass in probability in the sense that there exists $A>0$ such that $\lim_{N\ra+\infty}\mu_0^N\{\ls 1,\pi^N\rs>A\}=0$, and this holds with $A=\mu_0(\T^d)+\delta$, $\delta>0$, whenever $\{\mu_0^N\}$ is associated to a macroscopic profile $\mu_0\in\MMM_+(\T^d)$.$\hfill\Box$\\

We will denote by $\bs{\mathcal{Q}}^{\infty,\ell}$ the set of all subsequential limit points of $\{\bs{Q}^{N,\ell}\}_{N\in\NN}$ for each fixed $\ell\in\ZZ_+$, where $\bs{Q}^{N,\ell}$ is the law of the micro-empirical density of the ZRP, and we will denote by $\bs{\mathcal{Q}}^{\infty,\infty}$ the subsequential limit set $\Lim_{\ell\ra+\infty}\bs{\mathcal{Q}}^{\infty,\ell}$ so that
$$\bs{\mathcal{Q}}^{\infty,\infty}=\Lim_{\ell\ra+\infty}\bs{\mathcal{Q}}^{\infty,\ell}=\Lim_{\ell\ra+\infty}\Lim_{N\ra+\infty}\bs{Q}^{N,\ell}.$$
Likewise, we set 
$$\bs{\mathcal{Q}}^{\infty,0}:=\Lim_{\ee\ra0}\bs{Q}^{\infty,\ee},\quad\bs{Q}^{\infty,\ee}:=\Lim_{N\ra+\infty}\bs{Q}^{N,\ee},\quad\ee>0.$$
\vspace{-0.5cm}
\begin{prop}\label{Leb} For each $\ell\in\ZZ_+$, $\ee>0$, $$\bs{\mathcal{Q}}^{\infty,\ell}\cup\bs{\mathcal{Q}}^{\infty,\ee}\subs\PP L^\infty_{w^*}\big(0,T;T_{\LL_{\T^d}}\bbar{\PP}_{1}(\T^d\x\RR_+)\big)$$
	and if the sequence $\{\mu_0^N\}$ of initial distribution has total mass $\mathfrak{m}>0$
	then $$\bs{\mathcal{Q}}^{\infty,\ell}\cup\bs{\mathcal{Q}}^{\infty,\ee}\subs\PP L^\infty_{w^*}\big(0,T;\bbar{\Y}_{1,\mathfrak{m}}(\T^d)\big)$$
	In particular since $L_{w^*}^\infty\big(0,T;\bbar{\Y}_1(\T^d)\big)$ and $L_{w^*}^\infty(0,T;\bbar{\PP}_{1,\mathfrak{m}}(\T^d\x\RR_+)\big)$ are closed subspaces, the same inclusions hold for the set $\bs{\mathcal{Q}}^{\infty,\infty}\cup\bs{\mathcal{Q}}^{\infty,0}$.	
\end{prop}\textbf{Proof} We will show the claim for the set $\bs{\mathcal{Q}}^{\infty,\ee}$, the proof for the set $\bs{\mathcal{Q}}^{\infty,\ell}$ being similar. So let $\ee>0$ and let $\bs{Q}^{\infty,\ee}\in\bs{\mathcal{Q}}^{\infty,\ee}$ be a subsequential limit point. We note first that for any generalized Young path-measure $\bs\pi\in L^\infty_{w^*}(0,T;\bbar{\MMM}_{1,+}(\T^d\x\RR_+))$ it holds that 
$$\bs\pi\in L_{w^*}^\infty(0,T;\bbar{\Y}_1(\T^d))\;\Llra\;\lls H(U),\bs\pi\rrs=\int_0^T\int H_t(u)\df u\df t,\;\forall H\in L^1(0,T;C(\T^d)).$$
Therefore, since $L^1(0,T;C(\T^d))$ is separable, in order to prove the first claim it suffices to show that 
\begin{equation}\label{ToConcludeLebProof}\bs{Q}^{\infty,\ee}\Big\{\bs\pi\in L_{w^*}^\infty(0,T;\bbar{\PP}_{1}(\T^d\x\RR_+))\Bigm|\lls H(U),\bs\pi\rrs=\int_0^T\int H_t(u)\df u\df t\Big\}=1\end{equation}
for each $H\in L^1(0,T;C(\T^d))$. This follows by the portmanteau theorem. Indeed, the functional $\lls H(U),\cdot\rrs\colon L_{w^*}^\infty(0,T;\bbar{\MMM}_1(\T^d\x\RR_+))\to\RR$ is $w^*$-continuous and thus for each $\delta>0$ the set 
$$A_{H,\delta}:=\Big\{\bs\pi\in L_{w^*}^\infty(0,T;\bbar{\PP}_1(\T^d\x\RR_+))\Bigm|\Big|\lls H(U),\bs\pi\rrs-\int_0^T\int H_t(u)\df u\df t\Big|>\delta\Big\}$$ is open in $L_{w^*}^\infty(0,T;\bbar{\PP}_1(\T^d\x\RR_+))$. Thus if for each $\ee>0$ we pick a sequence $\{k_N^{(\ee)}\}\subs\NN$ such that $\bs{Q}^{\infty,\ee}=\lim_{N\ra+\infty} \bs{Q}^{k_N^{(\ee)},\ee}$ then by the portmanteau theorem 
\begin{align*}
\bs{Q}^{\infty,\ee}(A_{H,\delta})&\leq\liminf_{N\ra+\infty}\bs{Q}^{k_N^{(\ee)},\ee}(A_{H,\delta})\leq\limsup_{N\ra+\infty}\bs{Q}^{N,\ee}(A_{H,\delta})\\
&=\limsup_{N\ra+\infty}P^N\Big\{\Big|\int_0^T\big(\ls H_t(U),\bs\pi^{N,\ee}_t\rs-H_t(u)\big)\df u\df t\Big|>\delta\Big\}=0.\end{align*} The last limit inferior is indeed equal to $0$ since due to the fact that $H\in L^1(0,T;C(\T^d))$ it holds that 
$$\lim_{N\ra+\infty}\lls H(U),\bs\pi^{N,\ee}\rrs=\lim_{N\ra+\infty}\int_0^T\fr{N^d}\sum_{x\in\T_N}H_t\Big(\frac{x}{N}\Big)\df t=\int_0^T\int H_t(u)\df u\df t.$$
Since $\delta>0$ was arbitrary,~\eqref{ToConcludeLebProof} holds.

For the second claim we note that by the Lebesgue differentiation theorem for any generalized Young path-measure $\bs\pi\in L^\infty_{w^*}(0,T;\bbar{\MMM}_{1,+}(\T^d\x\RR_+))$
it holds that $$\bs\pi\in L_{w^*}^\infty(0,T;\bbar{\PP}_{1,\mathfrak{m}}(\T^d\x\RR_+))\;\Llra\;\lls f\Lambda,\bs\pi\rrs=\mathfrak{m}\int_0^Tf_t\df t,\quad\forall f\in L^1(0,T),$$
where $f\Lambda$ is the map given by $(f\Lambda)_t(u,\lambda)=\lambda f_t$ for $(t,u,\lambda)\in[0,T]\x\T^d\x\RR_+$. Therefore in order to prove the second claim we have to show  that for any $\bs{Q}^{\infty,\ee}\in\bs{\mathcal{Q}}^{\infty,\ee}$, $f\in L^1(0,T)$ and $\delta>0$ 
$$\bs{Q}^{\infty,\ee}\Big\{\Big|\lls f\Lambda,\bs\pi\rrs-\mathfrak{m}\int_0^Tf_t\df t\Big|>\delta\Big\}=0$$
But the set $B_{f,\delta}:=\big\{\big|\lls f\Lambda,\bs\pi\rrs-\mathfrak{m}\int_0^Tf_t\df t\big|>\delta\big\}$ is open in $L_{w^*}^\infty(0,T;\bbar{\MMM}_1(\T^d\x\RR_+))$ and therefore by the portmanteau theorem if $\{k_N^{(\ee)}\}_{N\in\NN}$ is a sequence such that $\lim_{N\ra+\infty}\bs{Q}^{k_N^{(\ee)},\ee}=\bs{Q}^{\infty,\ee}$ then 
\begin{align*}
\bs{Q}^{\infty,\ee}(B_{f,\delta})&\leq\liminf_{N\ra+\infty}\bs{Q}^{k_N^{(\ee)},\ee}(B_{f,\delta})\leq\limsup_{N\ra+\infty} P^N\Big\{\Big|\lls f\Lambda,\bs\pi^{N,\ee}\rrs-\mathfrak{m}\int_0^Tf_t\df t\Big|>\delta\Big\}\\
&=\limsup_{N\ra+\infty} P^N\Big\{\Big|\int_0^T\big(\ls\Lambda,\bs\pi^{N,\ee}_t\rs-\mathfrak{m}\big)f_t\df t\Big|>\delta\Big\}.
\end{align*}
Therefore, since $\ls \Lambda,\bs\pi^{N,\ee}\rs=\ls 1,\pi^N\rs$ for all $\ee>0$ it follows that
\begin{align*}
\bs{Q}^{\infty,\ee}(B_{f,\delta})&\leq\limsup_{N\ra+\infty}\mu_0^N\Big\{\big|\ls 1,\pi^N\rs-\mathfrak{m}\big|>\delta\Big(\int_0^T|f_t|\df t\Big)^{-1}\Big\}=0
\end{align*}
where the last limit superior is equal to $0$ due to the fact that $\{\mu_0^N\}$ has total mass $\mathfrak{m}$ in probability.$\hfill\Box$
\begin{prop}\label{ModRelComp} Let $P^N$ be the law of the diffusively rescaled ZRP on the Skorohod space $D(0,T;\MM_N^d)$, starting from a sequence $\{\mu_0^N\in\PP\MM_N^d\}_{N\in\NN}$ of initial distributions having total mass $\mathfrak{m}>0$ in probability and let $\bs{\pi}^{N,\ell;M}$ be the $M$-modified micro-empirical density process of the ZRP. Then the family of Borel probability measures 
	\begin{equation}\label{ModMicroEmpLaw}
	\bs{Q}^{N,\ell;M}:=\bs\pi^{N,\ell;M}_\sharp P^N\in\PP L^\infty(0,T;\bbar{\PP}_1(\T^d\x\RR_+)),\quad N\in\NN,\;\ell\in\ZZ_+,\;M>0
	\end{equation}
	is sequentially relatively compact. 
\end{prop}\textbf{Proof} Since the recession operator is a contraction, for all $F\in L^1(0,T;\bbar{C_1}(\T^d\x\RR_+))$
\begin{align*}
\lls F,\bs\pi^{N,\ell,M}\rrs&\leq\int_0^T\fr{N^d}\sum_{x\in\T_N^d}\Big\{\|F_t\|_{\infty,1}\big(1+\eta_t^\ell(x)\mn M\big)+\|RF_t\|_\infty(\eta_t^\ell(x)-M)^+\Big\}\df t\\
&\stackrel{P^N\mbox{-a.s.}}\leq\|F\|_{\infty,1;1}\big(1+\ls 1,\pi_0^N\rs\big)
\end{align*}
and therefore $\|\bs\pi^{N,\ell;M}\|_{TV,1;\infty}\leq1+\ls 1,\pi_0^N\rs$ for $P^N$-a.s.~all $\eta\in D(0,T;\MM_N^d)$ for all $N\in\NN$, $\ell\in\ZZ_+$, $M>0$. Therefore as in the proof of Proposition~\ref{YoungRelComp} it follows that 
$$\lim_{r\ra+\infty}\sup_{(N,\ell,M)\in\NN\x\ZZ_+\x(0,+\infty)}P^N\{\|\bs\pi^{N,\ell;M}\|_{TV,1;\infty}>r\}=0$$
which according to the Prokhorov-Le Cam and Banach-Alaoglu theorems proves the relative compactness of the family $\{\bs{Q}^{N,\ell;M}\}_{(N,\ell,M)}$. $\hfill\Box$\\

This proposition implies that the family $\bs{Q}^{N,\ee;M}:=\bs{Q}^{N,[N\ee];M}$, $N\in\NN$, $\ee>0$, $M>0$ is also relatively compact. The following simple lemma will be useful in comparing as $M\ra+\infty$ the micro and macro-empirical distributions $\bs\pi^{N,\ell}$ and $\bs\pi^{N,\ee}$ of the ZRP with their modified versions $\bs\pi^{N,\ell;M}$ and $\bs\pi^{N,\ee;M}:=\bs\pi^{N,[N\ee];M}$ defined in~\eqref{MModMicrEmpDens}.
\begin{lemma}\label{PathRecession} If for each $F\in L^1(0,T;\bbar{C}_r(\T^d\x\RR_+))$ we set $$F^{(\lambda,r)}=(F_t^{(\lambda,r)})_{0\leq t\leq T}=\Big(\frac{F_t(\cdot,\lambda)}{1+\lambda^r}\Big)_{0\leq t\leq T}\in L^1(0,T;C(\T^d))$$ 
	then \begin{equation}\label{RecessionInL1}
	\lim_{\lambda\ra+\infty}\big\|F_t^{(\lambda,r)}-RF_t\big\|_{\infty;1}\leq\lim_{M\ra+\infty}\int_0^T\sup_{\lambda>M}\big\|F_t^{(\lambda,r)}-\bar{R}_rF_t\big\|_\infty\df t=0.\end{equation}
\end{lemma}\textbf{Proof} The first inequality in~\eqref{RecessionInL1} is obvious. The right hand-side limit follows by the dominated convergence theorem. Indeed, if $F\in^1(\bbar{C_r}(\T^d\x\RR_+))$ then for almost all $0\leq t\leq T$ and all $\lambda\in\RR_+$
$$\|F_t^{(\lambda,r)}-R_rF_t\|_\infty\leq\|F_t^{(\lambda,r)}\|+\|R_rF_t\|_\infty\leq2\|F_t\|_{\infty,r},$$
and therefore the map $t\mapsto\sup_{M>0}\|F_t^{(\lambda,r)}-R_rF_t\|_\infty$ is dominated by the integrable function $2\|F_\cdot\|_{\infty,r}\in L^1(0,T)$. Thus, since by the definition of the recession operator 
$$\lim_{M\ra+\infty}\sup_{M>0}\big\|F^{(\lambda,r)}_t-R_rF_t\big\|_\infty=0$$ for almost all $0\leq t\leq T$, the limit~\eqref{RecessionInL1} follows by the dominated convergence theorem.$\hfill\Box$

\begin{prop}\label{FluidSolidSeparation} Let $P^N$ be the law of the diffusively rescaled ZRP on the Skorohod space $D(0,T;\MM_N^d)$ starting from a sequence $\{\mu_0^N\in\PP\MM_N^d\}_{N\in\NN}$ of initial distributions having total mass $\mathfrak{m}>0$. Then for any $F\in L^1(0,T;\bbar{C}_1(\T^d\x\RR_+))$
	$$\lim_{M\ra+\infty}\sup_{\ell\in\ZZ_+}\sup_{N\in\NN}\EE^N|\lls F,\bs\pi^{N,\ell}\rrs-\lls F,\bs\pi^{N,\ell;M}\rrs|=0.$$	
\end{prop}\textbf{Proof} For each $F\in L^1(0,T;\bbar{C}_1(\T^d\x\RR_+))$,
\begin{align*}
\lls F,\bs\pi^{N,\ell}\rrs-\lls F,\bs\pi^{N,\ell;M}\rrs&=\int_0^T\fr{N^d}\sum_{x\in\T_N^d}\Big\{F_t\Big(\frac{x}{N},\eta_t^\ell(x)\Big)-F_t\Big(\frac{x}{N},M\Big)\Big\}\1_{\{\eta_t^\ell(x)>M\}}dt\\
&\quad-\int_0^T\fr{N^d}\sum_{x\in\T_N^d}RF_t\Big(\frac{x}{N}\Big)(\eta_t^\ell(x)-M)^+\\
&=\int_0^T\fr{N^d}\sum_{x\in\T_N^d}\Big\{F_t\Big(\frac{x}{N},\eta_t^\ell(x)\Big)-RF_t\Big(\frac{x}{N}\Big)\eta_t^\ell(x)\Big\}\1_{\{\eta_t^\ell(x)>M\}}dt\\
&\quad-\int_0^T\fr{N^d}\sum_{x\in\T_N^d}\Big\{F_t\Big(\frac{x}{N},M\Big)-RF_t\Big(\frac{x}{N}\Big)M\Big\}\1_{\{\eta_t^\ell(x)>M\}}dt\\
&=\int_0^T\fr{N^d}\sum_{x\in\T_N^d}\Big\{\frac{F_t(\frac{x}{N},\eta_t^\ell(x))}{\eta^\ell_t(x)}-RF_t\Big(\frac{x}{N}\Big)\Big\}\eta_t^\ell(x)\1_{\{\eta_t^\ell(x)>M\}}dt\\
&\quad-\int_0^T\fr{N^d}\sum_{x\in\T_N^d}
\Big\{\frac{F_t(\frac{x}{N},M)}{M}-RF_t\Big(\frac{x}{N}\Big)\Big\}M\1_{\{\eta_t^\ell(x)>M\}}dt
\end{align*}
and therefore 
\begin{align*}
\big|\lls F,\pi^{N,\ell}\rrs-\lls F,\bs\pi^{N,\ell;M}\rrs\big|&\leq\int_0^T\sup_{\lambda\geq M}\Big\|\frac{F_t(\cdot,\lambda)}{\lambda}-RF_t\Big\|_\infty\fr{N^d}\sum_{x\in\T_N^d}\eta_t^\ell(x)\1_{\{\eta_t^\ell(x)>M\}}dt\\
&\quad+\int_0^T\Big\|\frac{F_t(\cdot,M)}{M}-RF_t\Big\|_\infty\fr{N^d}\sum_{x\in\T_N^d}
M\1_{\{\eta_t^\ell(x)>M\}}dt\\
&\leq 2\int_0^T\sup_{\lambda\geq M}\Big\|\frac{F_t(\cdot,\lambda)}{\lambda}-RF_t\Big\|_\infty\fr{N^d}\sum_{x\in\T_N^d}\eta_t^\ell(x)dt\\
&=2\int_0^T\sup_{\lambda\geq M}\Big\|\frac{F_t(\cdot,\lambda)}{\lambda}-RF_t\Big\|_\infty\frac{|\eta_t|_1}{N^d}dt.
\end{align*}
Consequently, by the conservation of the total number of particles, $P^N$-a.s.~it holds that 
$$\big|\lls F,\bs\pi^{N,\ell}\rrs-\lls F,\bs\pi^{N,\ell;M}\rrs\big|\leq 2\ls 1,\pi_0^N\rs\int_0^T\sup_{\lambda\geq M}\Big\|\frac{F_t(\cdot,\lambda)}{\lambda}-RF_t\Big\|_\infty dt=:2\ls 1,\pi^N_0\rs A_M(F).$$
By this inequality it follows that 
$$\EE^N|\lls F,\bs\pi^{N,\ell}\rrs-\lls F,\bs\pi^{N,\ell;M}\rrs|\leq2A_M(F)\EE^N\ls 1,\pi^N_0\rs=2A_M(F)\int\ls 1,\pi^N\rs\df\mu_0^N$$
for all $(N,\ell)\in\NN\x\ZZ_+$ and therefore 
\begin{equation}\label{ModIneq}
\sup_{(N,\ell)\in\NN\x\ZZ_+}\EE^N|\lls F,\bs\pi^{N,\ell}\rrs-\lls F,\bs\pi^{N,\ell;M}\rrs|
\leq 2A_M(F)\sup_{N\in\NN}\int\ls 1,\pi^N\rs\df\mu_0^N.
\end{equation}
But since $F\in L^1(0,T;\bbar{C}_1(\T^d\x\RR_+))$, by Lemma~\ref{PathRecession} we have that  
$$\lim_{M\ra+\infty}A_M(F)=\lim_{M\ra+\infty}\int_0^T\sup_{\lambda\geq M}\Big\|\frac{F_t(\cdot,\lambda)}{\lambda}-RF_t\Big\|_\infty dt=0$$ and therefore the claim follows by
Lemma~\ref{BoundedTotParticlInMacrLim} and the inequality~\eqref{ModIneq}.$\hfill\Box$\\

By inequality~\eqref{ModIneq} it is obvious that this last proposition is also valid for the macro-empirical density $\bs\pi^{N,\ee;M}$. Since the families $\{\bs{Q}^{N,\ell}\}_{(N,\ell)\in\NN\x\ZZ_+}$ and $\{\bs{Q}^{N,\ell;M}\}_{(N,\ell,M)\in\NN\x\ZZ_+^2}$ are sequentially relatively compact, the family of the joint laws 
$$\bbar{\bs{Q}}^{N,\ell;M}:=(\bs\pi^{N,\ell;M},\bs\pi^{N,\ell})_\sharp P^N\in\PP\big(L_{w^*}^\infty(0,T;\bbar{\PP}_1(\T^d\x\RR_+))^2\big)$$
is also sequentially relatively compact. Since $$\{\bbar{\bs{Q}}^{N,\ee;M}|N\in\NN,\;\ee>0\;M>0\}\subs\{\bbar{\bs{Q}}^{N,\ell;M}|N\in\NN,\;\ell\in\ZZ_+\;M>0\}$$ the family $\{\bbar{\bs{Q}}^{N,\ee;M}\}$ is also relatively compact. By Proposition~\ref{ModRelComp} it follows that any limit point of the family $\{\bbar{\bs{Q}}^{N,\ell;M}\}$ as $N$, $\ell$ and then $M$ tend to infinity is concentrated on the diagonal and the same is true for the family $\{\bbar{\bs{Q}}^{N,\ee;M}\}$. We state this more precisely as a proposition.

\begin{cor}\label{ModIsUnMod} Any subsequential limit point 
	$$\bbar{\bs{Q}}\in \Lim_{M\ra+\infty}\Lim_{\ell\ra+\infty}\Lim_{N\ra+\infty}\bbar{\bs{Q}}^{N,\ell;M}$$ 
	is concentrated on the diagonal, i.e.
	$$\bbar{\bs{Q}}\big\{(\bs\pi_1,\bs\pi_2)\in L_{w^*}^\infty(0,T;\bbar{\PP}_1(\T^d\x\RR_+))^2\bigm|\bs\pi_1=\bs\pi_2\big\}=1.$$
	The same is true for any subsequential limit point 
	$$\bbar{\bs{Q}}\in \Lim_{M\ra+\infty}\Lim_{\ee\ra 0}\Lim_{N\ra+\infty}\bbar{\bs{Q}}^{N,\ee;M}$$ 
\end{cor}\textbf{Proof} For any $F\in L^1(0,T;\bbar{C}_1(\T^d\x\RR_+))$ the map 
$I_F\colon L^\infty(0,T;\bbar{\PP}_1(\T^d\x\RR_+))^2\to\RR_+$ given by 
$$I_F(\bs\pi_1,\bs\pi_2)=\big|\lls F,\bs\pi_1\rrs-\lls F,\bs\pi_2\rrs\big|$$ is $w^*$-continuous and therefore for any $\delta>0$ the set 
$$A_{F,\delta}:=\big\{(\bs\pi_1,\bs\pi_2)\in L_{w^*}^\infty(0,T;\bbar{\PP}_1(\T^d\x\RR_+))^2\bigm||\lls F,\bs\pi_1\rrs-\lls F,\bs\pi_2\rrs|>\delta\big\}$$
is open in the product of the $w^*$-topologies on $L^\infty(0,T;\bbar{\PP}_1(\T^d\x\RR_+))^2$. Let now $\bbar{\bs{Q}}$ be a subsequential limit point of $\bbar{\bs{Q}}^{N,\ell;M}$ as $N$, $\ell$ and then $M$ tend to infinity. Then by the portmanteau theorem, Chebyshev's inequality and Proposition~\ref{FluidSolidSeparation} 
\begin{align}\label{CombinedToLimsup}
\bbar{\bs{Q}}(A_{F,\delta})&\leq\limsup_{M\ra+\infty}\limsup_{\ell\ra+\infty}\limsup_{N\ra+\infty}\bbar{\bs{Q}}^{N,\ell;M}(A_{F,\delta})\nonumber\\
&\leq\limsup_{M\ra+\infty}\sup_{\ell\in\ZZ_+}\sup_{N\in\NN}P^N\big\{|\lls F,\bs\pi^{M,\ell;M}\rrs-\lls F,\bs\pi^{N,\ell}\rrs|>\delta\big\}=0.
\end{align} 
Now, since $L^1(0,T;\bbar{C}_1(\T^d\x\RR_+))$ is separable we can choose a dense sequence $\{F_n\}_{n\in\NN}$ and then 
$$\bbar{\bs{Q}}\big\{(\bs{\pi}_1,\bs{\pi}_2)\bigm|\bs{\pi}_1\neq\bs{\pi}_2\big\}=\bbar{\bs{Q}}\Big(\bigcup_{n,k\in\NN}\Big\{|\lls F_n,\bs\pi_1\rrs-\lls F_n,\bs\pi_2\rrs|>\fr{k}\Big\}\Big)=0$$
which concludes the proof. Let us verify the first inequality in~\eqref{CombinedToLimsup}. By the definition of $\bbar{\bs{Q}}$ there exist then a sequence $\{M_n\}_{n\in\NN}\subs\RR_+$ converging to infinity and
$$\bbar{\bs{Q}}^{\infty,\infty;M_n}\in\Lim_{\ell\ra+\infty}\Lim_{N\ra+\infty}\bbar{\bs{Q}}^{N,\ell;M_n},\quad\forall n\in\NN$$ such that $\lim_{n\ra+\infty}\bbar{\bs{Q}}^{\infty,\infty;M_n}=\bbar{\bs{Q}}$. Then by the portmanteau theorem 
\begin{equation}\label{MLim1}\bbar{\bs{Q}}(A_{F,\delta})\leq\liminf_{n\ra+\infty}\bbar{\bs{Q}}^{\infty,\infty;M_n}(A_{F,\delta}).
\end{equation}
Next, by the definition of $\bs{Q}^{\infty,\infty;M_n}$, $n\in\NN$, there exists for each $n\in\NN$ a sequence $\{m_\ell^{(n)}\}_{\ell\in\NN}\subs\ZZ_+$ and 
$$\bbar{\bs{Q}}^{\infty,m_\ell^{(n)};M_n}\in\Lim_{N\ra+\infty}\bbar{\bs{Q}}^{N,m_\ell^{(n)};M_n},\quad\ell\in\NN$$
such that $\lim_{\ell\ra+\infty}\bbar{\bs{Q}}^{\infty,m_\ell^{(n)};M_n}=\bbar{\bs{Q}}^{\infty,\infty;M_n}$ for all $n\in\NN$ and thus by the portmanteau theorem 
\begin{equation}\label{lLim1}\bbar{\bs{Q}}^{\infty,\infty;M_n}(A_{F,\delta})\leq\liminf_{\ell\ra+\infty}\bbar{\bs{Q}}^{\infty,m_\ell^{(n)};M_n}(A_{F,\delta}).\end{equation}
Finally, by the definition of $\bbar{\bs{Q}}^{\infty,m_\ell^{(n)};M_n}$, for any $(\ell,n)\in\NN^2$ there exists a sequence $\{k_N^{(\ell,n)}\}_{N\in\NN}$ such that 
$$\lim_{N\ra+\infty}\bbar{\bs{Q}}^{k_N^{(\ell,n)},m_\ell^{(n)};M_n}=\bbar{\bs{Q}}^{\infty,m_\ell^{(n)};M_n},\quad\forall(\ell,n)\in\NN^2$$
and therefore by the portmanteau theorem again 
\begin{equation}\label{NLim1}\bbar{\bs{Q}}^{\infty,m_\ell^{(n)};M_n}(A_{F,\delta})\leq\liminf_{N\ra+\infty}\bbar{\bs{Q}}^{k_N^{(\ell,n)},m_\ell^{(n)};M_n}(A_{F,\delta})\leq
\limsup_{N\ra+\infty}\bbar{\bs{Q}}^{N,m_\ell^{(n)};M_n}(A_{F,\delta}).\end{equation}
By combining~\eqref{MLim1},~\eqref{lLim1} and~\eqref{NLim1} we thus obtain 
\begin{align*}
\bbar{\bs{Q}}(A_{F,\delta})&\leq\liminf_{n\ra+\infty}\bbar{\bs{Q}}^{\infty,\infty;M_n}(A_{F,\delta})\leq\liminf_{n\ra+\infty}\liminf_{\ell\ra+\infty}\bbar{\bs{Q}}^{\infty,m_\ell^{(n)};M_n}(A_{F,\delta})\\
&\leq\liminf_{n\ra+\infty}\liminf_{\ell\ra+\infty}\limsup_{N\ra+\infty}\bbar{\bs{Q}}^{N,m_\ell^{(n)};M_n}(A_{F,\delta})\leq\limsup_{M\ra+\infty}\limsup_{\ell\ra+\infty}\limsup_{N\ra+\infty}\bbar{\bs{Q}}^{N,\ell;M}(A_{F,\delta})
\end{align*}
which proves the first inequality in~\eqref{CombinedToLimsup} and concludes the proof.$\hfill\Box$

\begin{cor}\label{AcrossAlreadyConv} Let $\{m_\ell\}_{\ell=0}^\infty$ and $\{k_N^{(\ell)}\}_{N=1}^\infty$, $\ell\in\ZZ_+$ be sequences of integers increasing to infinity such that the limit
	\begin{equation}\label{LimElementApproxSubseqN}\lim_{N\ra+\infty}\bs{Q}^{k_N^{(\ell)},m_\ell}=:\bs{Q}^{\infty,m_\ell}\end{equation}
	exists for all $\ell\in\ZZ_+$ and the limit 
	\begin{equation}\label{LimElementApproxSubseqell}\lim_{\ell\ra+\infty}\bs{Q}^{\infty,m_\ell}=:\bs{Q}^{\infty,\infty}\end{equation}
	also exists. Then if for each $M>0$ we set $\bs{Q}_*^{N,\ell;M}:=\bs{Q}^{k_N^{(\ell)},m_\ell;M}$,
	$$\bs{\mathcal{Q}}_*^{\infty,\infty;\infty}:=\Lim_{M\ra+\infty}\Lim_{\ell\ra+\infty}\Lim_{N\ra+\infty}\bs{Q}_*^{N,\ell;M}=\{\bs{Q}^{\infty,\infty}\}.$$
\end{cor}\textbf{Proof} Let $\bs{Q}_*\in\bs{\mathcal{Q}}^{\infty\infty;\infty}_*$. We will show that $\bs{Q}_*=\bs{Q}^{\infty,\infty}$. By the definition of the set $\bs{\mathcal{Q}}^{\infty\infty;\infty}_*$ there exists a sequence $\{M_n\}_{n=1}^\infty\subs(0,+\infty)$ diverging to $+\infty$ and $$\bs{Q}_*^{\infty,\infty;M_n}\in\bs{\mathcal{Q}}_*^{\infty,\infty;M_n}:=\Lim_{\ell\ra+\infty}\Lim_{N\ra+\infty}\bs{Q}_*^{N,\ell;M_n}$$ such that 
\begin{equation}\label{3MargSub}\bs{Q}_*=\lim_{n\ra+\infty}\bs{Q}_*^{\infty,\infty;M_n}.\end{equation} 
Then by the definition of the set $\bs{\mathcal{Q}}_*^{\infty,\infty;M_n}$ there exists for each $n\in\NN$ a strictly increasing sequence $\{m_\ell^{(n)}\}_{\ell=1}^\infty$ of natural numbers and $$\bs{Q}_*^{\infty,m_\ell^{(n)};M_n}\in\bs{\mathcal{Q}}_*^{\infty,m_\ell^{(n)};M_n}:=\Lim_{N\ra+\infty}\bs{Q}_*^{N,m_\ell^{(n)};M_n}$$ such that 
\begin{equation}\label{2MargSub}\bs{Q}_*^{\infty,\infty;M_n}=\lim_{\ell\ra+\infty}\bs{Q}_*^{\infty,m_\ell^{(n)};M_n}\end{equation} 
and finally by the definition of the set $\bs{\mathcal{Q}}_*^{\infty,m_\ell^{(n)};M_n}$ there exists for each $n\in\NN$ and $\ell\in\ZZ_+$ a strictly increasing sequence $\{k_N^{(n,\ell)}\}_{N\in\NN}$ of natural numbers such that 
\begin{equation}\label{1MargSub}\bs{Q}_*^{\infty,m_\ell^{(n)};M_n}=\lim_{N\ra+\infty}\bs{Q}_*^{k_N^{(n,\ell)},m_\ell^{(n)};M_n}.
\end{equation} 
Consequently
$$\bs{Q}_*=\lim_{n\ra+\infty}\lim_{\ell\ra+\infty}\lim_{N\ra+\infty}\bs{Q}_*^{k_N^{(n,\ell)},m_\ell^{(n)};M_n}=\lim_{n\ra+\infty}\lim_{\ell\ra+\infty}\lim_{N\ra+\infty}\bs{Q}^{k_{k_N^{(n,\ell)}}^{(m_\ell^{(n)})},m_{m_\ell^{(n)}};M_n}.$$

%
%
 We consider now the family of the joint laws 
$$\bbar{\bs{Q}}^{N,\ell;M}:=(\bs\pi^{N,\ell;M},\bs\pi^{N,\ell})_\sharp P^N\in\PP L_{w^*}^\infty(0,T;\bbar{\MMM}_1(\T^d\x\RR_+))^2$$
and set $\bbar{\bs{Q}}^{N,\ell;M}_*:=\bbar{\bs{Q}}^{k_N^{(\ell)},m_\ell;M}$ for all $N\in\NN$, $\ell\in\ZZ_+$, $M>0$. This family is relatively compact and thus for each $(n,\ell)\in\NN\x\ZZ_+$ there exists a further subsequence $\{k_{\theta_N}^{(n,\ell)}\}_N$ of $\{k_N^{(n,\ell)}\}_N$ such that the limit 
$$\bbar{\bs{Q}}_*^{\infty,m_\ell^{(n)};M_n}:=\lim_{N\ra+\infty}\bbar{\bs{Q}}_*^{k_{\theta_N}^{(n,\ell)},m_\ell^{(n)};M_n}\in\Lim_{N\ra+\infty}\bbar{\bs{Q}}_*^{N,m_\ell^{(n)};M_n}=:\bbar{\bs{\mathcal{Q}}}_*^{\infty,m_\ell^{(n)};M_n}$$
exists, for each $n\in\NN$ there exists a further subsequence $\{m_{j_\ell^{(n)}}^{(n)}\}_\ell$ of $\{m_\ell^{(n)}\}_\ell$ such that the limit 
$$\bbar{\bs{Q}}_*^{\infty,\infty;M_n}:=\lim_{\ell\ra+\infty}\bbar{\bs{Q}}_*^{\infty,m_{j_\ell^{(n)}}^{(n)};M_n}\in\Lim_{\ell\ra+\infty}\bbar{\bs{\mathcal{Q}}}_*^{\infty,m_\ell^{(n)};M_n}\subs\Lim_{\ell\ra+\infty}\bbar{\bs{\mathcal{Q}}}_*^{\infty,\ell;M_n}=:\bbar{\bs{\mathcal{Q}}}_*^{\infty,\infty;M_n}$$ exists and finally there exists a further subsequence $\{M_{i_n}\}_n$ of $\{M_n\}_n$ such that the limit 
$$\bbar{\bs{Q}}_*^{\infty,\infty;\infty}:=\lim_{n\ra+\infty}\bbar{\bs{Q}}_*^{\infty,\infty;M_{i_n}}\in\Lim_{M\ra+\infty}\bbar{\bs{\mathcal{Q}}}^{\infty,\infty;M}=:\bbar{\bs{\mathcal{Q}}}_*^{\infty,\infty;\infty}$$
exists. Then 
\begin{equation}\label{QStar}
\bbar{\bs{Q}}_*^{\infty,\infty;\infty}=\lim_{n\ra+\infty}\lim_{\ell\ra+\infty}\lim_{N\ra+\infty}\bbar{\bs{Q}}_*^{k_{\theta_N}^{(i_n,j_\ell^{(i_n)})},m_{j^{(i_n)}_\ell}^{(i_n)};M_{i_n}}.
\end{equation} Since obviously
$$\bbar{\bs{Q}}^{\infty,\infty;\infty}_*\in \Lim_{M\ra+\infty}\Lim_{\ell\ra+\infty}\Lim_{N\ra+\infty}\bbar{\bs{Q}}_*^{N,\ell;M}\subs\Lim_{M\ra+\infty}\Lim_{\ell\ra+\infty}\Lim_{N\ra+\infty}\bbar{\bs{Q}}^{N,\ell;M}$$ it follows by Corollary~\ref{ModIsUnMod} that $\bbar{\bs{Q}}^{\infty,\infty;\infty}_*$ is concentrated on the diagonal and thus its marginals coincide. But the first marginal of $\bbar{\bs{Q}}_*^{\infty,\infty;\infty}$ is $\bs{Q}_*$ and its second marginal is $\bs{Q}^{\infty,\infty}$. Indeed, the first marginal of $\bbar{\bs{Q}}^{N,\ell;M}_*$ is the measure $\bs\pi^1_\sharp\bbar{\bs{Q}}^{N,\ell;M}_*=\bs{Q}_*^{N,\ell;M}$ and therefore by~\eqref{QStar} and the continuity of $\bs\pi^1_\sharp$
\begin{align*}
\bs\pi^1_\sharp\bbar{\bs{Q}}_*^{\infty,\infty;\infty}&=\lim_{n\ra+\infty}\lim_{\ell\ra+\infty}\lim_{N\ra+\infty}\bs\pi^1_\sharp\bbar{\bs{Q}}_*^{k_{\theta_N}^{(i_n,j_\ell^{(i_n)})},m_{j^{(i_n)}_\ell}^{(i_n)};M_{i_n}} \\
&=\lim_{n\ra+\infty}\lim_{\ell\ra+\infty}\lim_{N\ra+\infty}\bs{Q}_*^{k_{\theta_N}^{(i_n,j_\ell^{(i_n)})},m_{j^{(i_n)}_\ell}^{(i_n)};M_{i_n}}.
\end{align*}
But for each $n\in\NN$, $\ell\in\ZZ_+$ the sequence $\{k_{\theta_N}^{(i_n,j_\ell^{(i_n)})}\}$ is a subsequence of $\{k_N^{(i_n,j_\ell^{(i_n)})}\}$
 and therefore by~\eqref{1MargSub} 
 $$\lim_{N\ra+\infty}\bs{Q}_*^{k_N^{(i_n,j_\ell^{(i_n)})},m_{j_\ell^{(i_\ell)}}^{(i_n)};M_{i_n}}=\bs{Q}_*^{\infty,m_{j_\ell^{(i_n)}}^{(i_n)};M_{i_n}}.$$
 Likewise, for each $n\in\NN$ the sequence $\{m_{j_\ell^{(n)}}^{(n)}\}_\ell$ is a subsequence of $\{m_\ell^{(n)}\}_\ell$ and thus by~\eqref{2MargSub} 
 $$\lim_{\ell\ra+\infty}\bs{Q}_*^{\infty,m_{j_\ell^{(i_n)}}^{(i_n)};M_{i_n}}=\bs{Q}_*^{\infty,\infty;M_{i_n}}$$
 and finally since $\{M_{i_n)}\}$ is a subsequence of $\{M_n\}$ it follows by~\eqref{3MargSub} that 
 $$\bs\pi^1_\sharp\bbar{\bs{Q}}_*^{\infty,\infty;\infty}=\lim_{n\ra+\infty}\bs{Q}_*^{\infty,\infty;M_{i_n}}=\bs{Q}_*.$$
  Similarly, the second marginal of $\bbar{\bs{Q}}_*^{N,\ell;M}$ is the measure $$\bs\pi^2_\sharp\bbar{\bs{Q}}_*^{N,\ell;M}=\bs{Q}^{k_N^{(\ell)},m_\ell}=:\bs{Q}_*^{N,\ell}$$ 
 and therefore by~\eqref{QStar}, the continuity of $\bs\pi^2_\sharp$ and the assumptions~\eqref{LimElementApproxSubseqN} and~\eqref{LimElementApproxSubseqell} 
 \begin{align*}
 \bs\pi^2_\sharp\bbar{\bs{Q}}_*^{\infty,\infty;\infty}&=\lim_{n\ra+\infty}\lim_{\ell\ra+\infty}\lim_{N\ra+\infty}\bs\pi^2_\sharp\bbar{\bs{Q}}_*^{k_{\theta_N}^{(i_n,j_\ell^{(i_n)})},m_{j^{(i_n)}_\ell}^{(i_n)};M_{i_n}} \\
 &=\lim_{n\ra+\infty}\lim_{\ell\ra+\infty}\lim_{N\ra+\infty}\bs{Q}_*^{k_{\theta_N}^{(i_n,j_\ell^{(i_n)})},m_{j^{(i_n)}_\ell}^{(i_n)}}=\lim_{n\ra+\infty}\lim_{\ell\ra+\infty}\bs{Q}_*^{\infty,m_{j^{(i_n)}_\ell}^{(i_n)}}=\bs{Q}^{\infty,\infty}
 \end{align*}
 and the proof is complete.$\hfill\Box$\\

The analogous result to Corollary~\ref{AcrossAlreadyConv} for the case of the family $\{\bbar{\bs{Q}}^{N,\ee;M}\}$. We leave its statement to the reader. More importantly, the $M$-modified empirical density process $\bs\pi^{N,\ell;M}$ yields the correct decomposition in regular and singular path-measures. We make this precise in the following proposition in the case of the family $\{\bbar{\bs{Q}}^{N,\ell;M}\}$ and leave the statement of the analogous result for the family $\{\bbar{\bs{Q}}^{N,\ee;M}\}$ to the reader.
\begin{prop}\label{LimitOfTruncatedLaw} Let $\widehat{D},D^\perp\colon L_{w^*}^\infty(0,T;\bbar{\MMM}_1(\T^d\x\RR_+))\to L_{w^*}^\infty(0,T;\bbar{\MMM}_1(\T^d\x\RR_+))$ be the regular and singular decomposition operators and let $\{k_N^{(\ell)}\}_{N=1}^\infty$, $\ell\in\NN$ and $\{m_\ell\}_{\ell=1}^\infty$ be diverging sequences such that the iterated limit
	$$\lim_{\ell\ra+\infty}\lim_{N\ra+\infty}\bs{Q}^{k_N^{(\ell)},m_\ell}=:\bs{Q}^\infty$$
	exists. Then
	$$\lim_{M\ra+\infty}\lim_{\ell\ra+\infty}\lim_{N\ra+\infty}\widehat{D}_\sharp\bs{Q}^{k_N^{(\ell)},m_\ell;M}=\widehat{D}_\sharp\bs{Q}^\infty$$
	and 
	$$\lim_{M\ra+\infty}\lim_{\ell\ra+\infty}\lim_{N\ra+\infty} D^\perp_\sharp\bs{Q}^{k_N^{(\ell)},m_\ell;M}=D^\perp_\sharp\bs{Q}^\infty.$$
	The same results also hold if we take $\ell=[N\ee]$ and take the limit as $N\ra+\infty$, $\ee\ra 0$ and then $M\ra+\infty$.
\end{prop}\textbf{Proof} We denote by $\Pi_M\colon\bbar{C}_1(\T^d\x\RR_+)\to C_1(\T^d\x\RR_+)$, $M>0$, the bounded linear operator defined by $\Pi_MF(u,\lambda)=F(u,\lambda\mn M)$. We also denote by $\Pi_M\colon L^1(0,T;\bbar{C}_1(\T^d\x\RR_+))\to L^1(0,T;C_1(\T^d\x\RR_+))$ the induced operator on the $L^1$-spaces. Then the adjoint $\Pi_M^*\colon L_{w^*}^\infty(0,T;\MMM_1(\T^d\x\RR_+))\to L_{w^*}^\infty(0,T;\bbar{\MMM}_1(\T^d\x\RR_+))$ is bounded and $w^*$-continuous and 
\begin{align*}
\lls\widehat{D}\circ \bs\pi^{N,\ell;M},F\rrs=\int_0^T\fr{N^d}\sum_{x\in\T_N^d}F_t\Big(\frac{x}{N},\eta_t^\ell(x)\mn M\Big)\df t=\lls\Pi_M^*\circ j^*\circ\bs\pi^{N,\ell},F\rrs
\end{align*}
for all $F\in L^1(0,T;\bbar{C}_1(\T^d\x\RR_+))$. Therefore 
\begin{equation}\label{Trunc} \widehat{D}\circ \bs\pi^{N,\ell;M}=\Pi_M^*\circ j^*\circ\bs\pi^{N,\ell}
\end{equation}
which yields 
\begin{equation}\label{ApproxFirstCoord}
\widehat{D}_\sharp\bs{Q}^{N,\ell;M}=(\Pi_M^*\circ j^*)_\sharp\bs{Q}^{N,\ell}
\end{equation}
for all $(N,\ell,M)\in\NN\x\ZZ_+^2$. Thus since $\Pi_M^*\circ j^*\colon L_{w^*}^\infty(0,T;\bbar{\MMM}_1(\T^d\x\RR_+))\to L_{w^*}^\infty(0,T;\bbar{\MMM}_1(\T^d\x\RR_+))$ is $(w^*,w^*)$-continuous we have by~\eqref{ApproxFirstCoord} and the assumption that 
$$\lim_{\ell\ra+\infty}\lim_{N\ra+\infty}\widehat{D}_\sharp\bs{Q}^{k_N^{(\ell)},m_\ell;M}=(\Pi_M^*\circ j^*)_\sharp \lim_{\ell\ra+\infty}\lim_{N\ra+\infty}\bs{Q}^{k_N^{(\ell)},m_\ell}=(\Pi_M^*\circ j^*)_\sharp\bs{Q}^\infty$$ and therefore in order to conclude the proof it suffices to prove that 
$$\lim_{M\ra+\infty}(\Pi_M^*\circ j^*)_\sharp\bs{Q}=\widehat{D}_\sharp\bs{Q},\quad\forall\;\bs{Q}\in\PP L_{w^*}^\infty(0,T;\bbar{\MMM}_1(\T^d\x\RR_+)).$$
But this is true since $\Pi_M^*\circ j^*$ $w^*$-converges pointwise to $\widehat{D}=E\circ j^*$ as we have seen in the proof of Corollary~\ref{PathMeasDecomp}, and this in turn implies that $(\Pi_M^*\circ j^*)_\sharp$ converges pointwise to $\widehat{D}_\sharp$ on $\PP L_{w^*}^\infty(0,T;\bbar{\MMM}_1(\T^d\x\RR_+))$. For the second limit, recalling that $T_M\colon C(\T^d)\to\bbar{C}_1(\T^d\x\RR_+)$ is the operator $T_Mf=(\Lambda-M)^+f(U)$,
\begin{equation}\label{Trunc2} D^\perp\circ \bs\pi^{N,\ell;M}=R^*\circ T_M^*\circ\bs\pi^{N,\ell}\end{equation} and as we have also seen in the proof of Corollary~\ref{PathMeasDecomp} the maps $R^*\circ T_M^*$ $w^*$-converge pointwise to $D^\perp$ as $M\ra+\infty$ and the claim follows as the first limit.$\hfill\Box$\\

In the course of the proof of the replacement lemma where we will compare the processes $\bs\pi^{N,\ell}$ and $\bs\pi^{N,[N\ee]}$ apart from interpolating between them with the micro-truncated double-block empirical density process $\bs\pi^{N,\ell;M;\ee}$ introduced in~\eqref{DoubleBlockProc1} we will also interpolate with the process \emph{macro-truncated double-block empirical density process} $\bs\pi^{N,\ell,\ee;M}\colon D(0,T;\MM_N^d)\to L_{w^*}^\infty(0,T;\bbar{\MMM}_1(\T^d\x\RR_+))$ defined via duality with test maps $F\in L^1(0,T;\bbar{C}_1(\T^d\x\RR_+))$ by the formula 
\begin{equation}\label{DoubleBlockProc2}\lls F,\bs\pi^{N,\ell,\ee;M}\rrs=\int_0^T\fr{N^d}\sum_{x\in\T_N^d}\Big\{F_t\Big(\frac{x}{N},\eta_t^{\ell,[N\ee]}(x)\mn M\Big)+RF_t\Big(\frac{x}{N}\Big)(\eta^{\ell,[N\ee]}(x)-M)^+\Big\}\df t,\end{equation}
where $\eta^{\ell,[N\ee]}(x)$ is the double-block average as defined in~\eqref{DoubleBlock}. We note that the regular part $\widehat{\bs\pi}^{N,\ell,\ee;M}$ of the process $\bs\pi^{N,\ell,\ee;M}$ is equal to $(\Pi_{M}^*\circ j^*)_\sharp \bs\pi^{N,\ell,\ee}$, where $$\bs\pi^{N,\ell,\ee}\colon D(0,T;\MM_N^d)\to L_{w^*}^\infty(\bbar{\MMM}_{1,+}(\T^d\x\RR_+))$$
is \emph{the double-block empirical process} defined by 
\begin{equation}\label{DoubleBlockEmpDens} 
\lls F,\bs\pi^{N,\ell,\ee}\rrs=\int_0^T\fr{N^d}\sum_{x\in\T_N^d}F\Big(\frac{x}{N},\eta_t^\ell(x)^{[N\ee]}\Big)\df t,\quad F\in L^1(0,T;\bbar{C}_1(\T^d\x\RR_+)).\end{equation}
\vspace{-0.5cm}\begin{prop}\label{DoubleBlockTruncRelComp}
Suppose that the ZRP starts from a sequence of initial profiles $\{\mu_0^N\}$ with total mass $\mathfrak{m}>0$ in probability. Then the families of the path-laws
	\begin{align}
	&\bs{Q}^{N,\ell;M;\ee}:=\bs\pi^{N,\ell;M;\ee}_\sharp P^N\in \PP L_{w^*}^\infty(0,T;\bbar{\PP}_1(\T^d\x\RR_+)),\label{MicroDTrunc}\\
	&\bs{Q}^{N,\ell,\ee;M}:=\bs\pi^{N,\ell,\ee;M}_\sharp P^N\in \PP L_{w^*}^\infty(0,T;\bbar{\PP}_1(\T^d\x\RR_+)),\label{MacroDTrunc}\\
	&\bs{Q}^{N,\ell,\ee}:=\bs\pi^{N,\ell,\ee}_\sharp P^N\in \PP L_{w^*}^\infty(0,T;\bbar{\PP}_1(\T^d\x\RR_+))\label{DoubleBlockLaw}
	\end{align}
	over all $N\in\NN$, $\ee>0$, $\ell\in\ZZ_+$ and $M>0$ are relatively compact and any limit point of each of these families as $N\ra+\infty$, $\ee\ra 0$, $\ell\ra+\infty$ and then $M\ra+\infty$ is concentrated on $L_{w^*}^\infty(0,T;\bbar{\Y}_{1,\mathfrak{m}}(\T^d))$.
\end{prop}\textbf{Proof} Since the recession operator $R\colon \bbar{C}_1(\T^d\x\RR_+)\to C(\T^d)$ is a contraction
\begin{align*}
\big|\lls F,\bs\pi^{N,\ell;M;\ee}\rrs\big|&\leq\int_0^T\frac{\|F_t\|_{\infty;1}}{N^d}\sum_{x\in\T_N^d}\big\{1+(\eta_t^\ell(x)\mn M)^{[N\ee]}\}\df t\\
&\quad+\int_0^T\frac{\|RF_t\|_\infty}{N^d}\sum_{x\in\T_N^d}(\eta_t^\ell(x)-M)^{+[N\ee]}\df t\\
&\leq\int_0^T\frac{\|F_t\|_{\infty;1}}{N^d}\sum_{x\in\T_N^d}\big\{1+\eta_t^{\ell,[N\ee]}(x)\}\df t=\int_0^T\|F_t\|_{\infty;1}(1+\ls 1,\pi_t^N\rs)\df t.
\end{align*}
Therefore by the conservation of the total number of particles
\begin{equation}\label{Exotic1}
\sup_{(\ell,\ee,M)\in\ZZ_+\x(0,\infty)^2}\|\bs\pi^{N,\ell;M;\ee}\|_{TV,1;\infty}\leq1+\ls 1,\pi_0^N\rs,\quad P^N\mbox{-a.s.},\;N\in\NN.
\end{equation}
and thus
$$\sup_{(N,\ell,M,\ee)\in\NN\x\ZZ_+\x(0,\infty)^2}P^N\big\{\|\bs\pi^{N,\ell;M;\ee}\|_{TV,1;\infty}>A\big\}\leq\sup_{N\in\NN}\mu_0^N\Big\{\ls 1,\pi^N\rs>A-1\big\}.$$
The claim then follows by Lemma~\ref{BoundedTotParticlInMacrLim}. The processes $\bs\pi^{N,\ell,\ee;M}$ and $\bs\pi^{N,\ell,\ee}$ also satisfy the bound~\eqref{Exotic1} and their laws are also relatively compact. Finally we note that the double-block empirical density map $\bs\pi^{N,\ell,\ee}\colon\MM_N^d\to\bbar{\MMM}_1(\T^d\x\RR_+)$ satisfies for all $H\in L^1(0,T;C(\T^d))$ the relations 
$$\ls H(U),\bs\pi^{N,\ell,\ee}\rs=\fr{N^d}\sum_{x\in\T_N^d}H\Big(\frac{x}{N}\Big)\quad\mbox{and}\quad\ls\Lambda,\bs\pi^{N,\ell,\ee}\rs=\ls 1,\pi^N\rs$$
and therefore the proof of the second claim follows similarly to the proof of Proposition~\eqref{Leb}.$\hfill\Box$\\

By similar computations to the ones yielding inequality~\eqref{ModIneq} we also obtain the inequality 
$$\sup_{(N,\ell,\ee)\in\NN\x\ZZ_+\x(0,+\infty)}\EE^N|\lls F,\bs\pi^{N,\ell,\ee}\rrs-\lls F,\bs\pi^{N,\ell,\ee;M}\rrs|
\leq 2A_M(F)\sup_{N\in\NN}\int\ls 1,\pi^N\rs\df\mu_0^N$$ 
where $A_M(F)=\int_0^T\sup_{\lambda\geq M}\big\|\frac{F_t(\cdot,\lambda)}{\lambda}-RF_t\big\|_\infty dt$. This inequality yields the analogous result to Proposition~\ref{FluidSolidSeparation} for the double-block empirical density $\bs\pi^{N,\ell,\ee}$, i.e.
\begin{equation}\label{DoubleBlockFluidSolidSeparation}
	\lim_{M\ra+\infty}\sup_{\ell\in\ZZ_+}\sup_{\ee>0}\sup_{N\in\NN}\EE^N|\lls F,\bs\pi^{N,\ell,\ee}\rrs-\lls F,\bs\pi^{N,\ell,\ee;M}\rrs|=0,\quad\forall F\in L^1(0,T;\bbar{C}_1(\T^d\x\RR_+)).
	\end{equation} 
	\vspace{-0.5cm}
 \begin{prop} Let $\bbar{\bs{Q}}^{N,\ell,\ee;M}:=(\bs\pi^{N,\ell,\ee},\bs\pi^{N,\ell,\ee;M})_\sharp P^N$, $N\in\NN$, $\ell\in\ZZ_+$, $\ee,M>0$, be the family joint laws of the double-block and the macro-truncated double-block empirical density process of the ZRP and let 
 	$$\bbar{\bs{Q}}\in\bbar{\bs{\mathcal{Q}}}^{\infty,\infty,0;\infty}:=\Lim_{M\ra+\infty}\Lim_{\ell\ra+\infty}\Lim_{\ee\ra 0}\Lim_{N\ra+\infty}\bbar{\bs{Q}}^{N,\ell,\ee;M}$$ 
 be a limit point of this family.  Then $\bbar{\bs{Q}}$ is concentrated on the diagonal and if $\{m_\ell\}_{\ell\in\NN}\subs\ZZ_+$ is a diverging sequence, $\{\theta_\ee^{(\ell)}\}_{\ee>0}$, $\ell\in\NN$, are maps converging to zero as $\ee\ra 0$ and $\{k_N^{(\ell,\ee)}\}_{N\in\NN}$, $(\ell,\ee)\in\NN\x(0,+\infty)$ are diverging sequences such that the iterated limit
 $$\lim_{\ell\ra+\infty}\lim_{\ee\ra 0}\lim_{N\ra+\infty}\bs{Q}^{k_N^{(\ell,\ee)},m_\ell,\theta_\ee^{(\ell)}}=:\bs{Q}^{\infty,\infty,0}$$
 exists then $$\Lim_{M\ra+\infty}\Lim_{\ell\ra+\infty}\Lim_{k\ra+\infty}\Lim_{N\ra+\infty}\bs{Q}^{k_N^{(\ell,\ee)},m_\ell,\theta_\ee^{(\ell)};M}=\{\bs{Q}^{\infty,\infty,0}\}.$$
 \end{prop}\textbf{Proof} The proof is similar to the proof of Corollaries~\ref{ModIsUnMod} and~\ref{AcrossAlreadyConv}.$\hfill\Box$\\
 
 Similarly to the relations~\eqref{Trunc} and~\eqref{Trunc2} the double block empirical process satisfies the relations
\begin{equation}\label{TruncDouble}\widehat{D}\circ \bs\pi^{N,\ell,\ee;M}=\Pi_M^*\circ j^*\circ\bs\pi^{N,\ell,\ee}\quad\mbox{and}\quad D^\perp\circ \bs\pi^{N,\ell,\ee;M}=R^*\circ T_M^*\circ\bs\pi^{N,\ell,\ee}\end{equation} and the analogous result to Proposition~\ref{LimitOfTruncatedLaw} holds for the double-block process.
  \begin{prop}\label{LimitOfTruncatedLawDoubleBlock} Let $\widehat{D},D^\perp\colon L_{w^*}^\infty(0,T;\bbar{\MMM}_1(\T^d\x\RR_+))\to L_{w^*}^\infty(0,T;\bbar{\MMM}_1(\T^d\x\RR_+))$ be the regular and singular decomposition operators and let $\{m_\ell\}_{\ell\in\NN}\subs\ZZ_+$ be a diverging sequence, for each $\ell\in\NN$ let $\{\theta_\ee^{(\ell)}\}_{\ee>0}$ be a map converging to zero as $\ee\ra 0$ and let $\{k_N^{(\ell,\ee)}\}_{N\in\NN}$, $(\ell,\ee)\in\NN\x(0,+\infty)$ be a diverging sequence as $N\ra+\infty$ such that the iterated limit
  	$$\lim_{\ell\ra+\infty}\lim_{\ee\ra 0}\lim_{N\ra+\infty}\bs{Q}^{k_N^{(\ell,\ee)},m_\ell,\theta_\ee^{(\ell)}}=:\bs{Q}^{\infty,\infty,0}$$
  	exists. Then
  	$$\lim_{M\ra+\infty}\lim_{\ell\ra+\infty}\lim_{\ee\ra 0}\lim_{N\ra+\infty}\widehat{D}_\sharp\bs{Q}^{k_N^{(\ell,\ee)},m_\ell,\theta_\ee^{(\ell)}}=\widehat{D}_\sharp\bs{Q}^{\infty,\infty,0}$$
  	and 
  	$$\lim_{M\ra+\infty}\lim_{\ell\ra+\infty}\lim_{\ee\ra 0}\lim_{N\ra+\infty} D^\perp_\sharp \bs{Q}^{k_N^{(\ell,\ee)},m_\ell,\theta_\ee^{(\ell)}}=D^\perp_\sharp\bs{Q}^{\infty,\infty,0}.$$
  	 \end{prop}\textbf{Proof} The claim follows by the same arguments in the proof of Proposition~\ref{LimitOfTruncatedLaw} by using the relations~\eqref{TruncDouble} in place of the relations~\eqref{Trunc} and~\eqref{Trunc2} and is omitted.$\hfill\Box$.
\subsection{The one-block estimate}\label{OBESection}
For any asymptotically linear cylinder map $\Psi\colon\MM_\infty^d\to\RR$ we set $A_\Psi^{N,\ell}:=\Psi^\ell-\bbar{\Psi}(\eta^\ell(0))$. Then in order to prove the one-block estimate we have to show that for all $H\in L^1(0,T;C(\T^d))$
\begin{equation}\label{ToProveOBE}
\limsup_{\ell,N\ra+\infty}\EE^N\Big|\int_0^T\fr{N^d}\sum_{x\in\T_N^d}H_t\Big(\frac{x}{N}\Big)\tau_xA_\Psi^{N,\ell}(\eta_t)\df t\Big|=0.\end{equation} 
We will first reduce the case of asymptotically linear cylinder maps to sublinear cylinder maps. So let $\Psi\colon\MM_\infty^d\to\RR$ be an asymptotically linear cylinder map with support $J=J_\Psi\subs\ZZ^d$ and gradient $\nabla\Psi(\infty)\in\RR^J$ at infinity. Then, denoting by $\eta_J\colon\MM_\infty^d\to\ZZ_+^J$ the natural projection, the function $\Psi_0=\Psi-\ls\nabla\Psi(\infty),\eta_J\rs$ is sublinear and $\Psi=\Psi_0+\ls\nabla\Psi(\infty),\eta_J\rs$. For each $\ell\in\ZZ_+$ 
$$\Psi^\ell=\fr{\ell_\star^d}\sum_{|y|\leq\ell}\tau_y\Psi=\Psi_0^\ell+\ls \nabla\Psi(\infty),\eta_J\rs^\ell=\Psi_0^\ell+\ls\nabla\Psi(\infty),\eta_J^\ell\rs$$
where $\eta_J^\ell$ is defined coordinate-wise, i.e.~$\eta^\ell_J=(\eta^\ell(x))_{x\in J}$. Furthermore for all $\rho\in[0,\rho_c]\cap\RR_+$ 
$$\wt{\Psi}(\rho)=\wt{\Psi_0}(\rho)+\int\ls\nabla\Psi(\infty),\eta_J\rs\df\nu_{\rho}^\infty=\wt{\Psi_0}(\rho)+\rho\ls\nabla\Psi(\infty),\1_J\rs$$
where $\1_J\cong(1,\ldots,1)\in\RR^J$ and thus 
\begin{align}\label{ExtAsymptViaSubl}
\bbar{\Psi}(\rho)&=\wt{\Psi}(\rho\mn\rho_c)+(\rho-\rho_c)^+\ls\nabla\Psi(\infty),\1_J\rs\nonumber\\
&=\wt{\Psi_0}(\rho\mn\rho_c)+(\rho\mn\rho_c)\ls\nabla\Psi(\infty),\1_J\rs+(\rho-\rho_c)^+\ls\nabla\Psi(\infty),\1_J\rs\nonumber\\
&=\bbar{\Psi_0}(\rho)+\rho\ls\nabla\Psi(\infty),\1_J\rs.
\end{align}

Therefore $A_\Psi^\ell= A_{\Psi_0}^\ell+\ls a,\eta_J^\ell\rs-\eta^\ell(0)\ls a,1_J\rs$ and thus in order to reduce the case of asymptotically linear cylinder functions to sublinear functions it suffices to show that  
\begin{equation}\label{ToProveOBEWCZRPAsymptLin}
\limsup_{\ell,N\ra+\infty}
\EE^N\Big|\int_0^T\fr{N^d}\sum_{x\in\T_N^d}H_t\Big(\frac{x}{N}\Big)\tau_x(\ls a,\eta_J^\ell\rs-\eta^\ell(0)\ls a,\1_J\rs)(\eta_t)\df t\Big|=0
\end{equation}
for all $H\in L^1(0,T;C(\T^d))$, $J\subs\ZZ^d$ and $a\in\RR^J$. Since
$$\big(\ls a,\eta_J^\ell\rs-\eta^\ell(0)\ls a,\1_J\rs\big)(\eta_t)=\sum_{z\in J}a_z\eta_t^\ell(z)-\eta_t^\ell(0)\sum_{z\in J}a_z$$
the time integrand in~\eqref{ToProveOBEWCZRPAsymptLin} is equal to 
\begin{align*}I_t^{N,\ell,H,a}(\eta)&:=\fr{N^d}\sum_{x\in\T_N^d}H_t\Big(\frac{x}{N}\Big)\sum_{z\in J}a_z\big\{\eta_t^\ell(x+z)-\eta_t^\ell(x)\big\}\\
&\;=\fr{N^d}\sum_{x\in\T_N^d}\sum_{z\in J}a_z\Big\{H_t\Big(\frac{x-z}{N}\Big)-H_t\Big(\frac{x}{N}\Big)\Big\}\eta_t^\ell(x).
\end{align*}

Since $H\in L^1(0,T;C(\T^d))$, given $\ee>0$, there exists a function $\bar{\delta}_\ee\in L^\infty(0,T)$ depending on $H$ such that $\bar{\delta}_\ee(t)>0$ for almost all $t\in[0,T]$ and 
$$d_{\T^d}(u,\y)<\bar{\delta}_\ee(t)\quad\Lra\quad |H_t(u)-H_t(\y)|<\ee,\quad\mbox{a.s.-}\forall\;t\in[0,T].$$ Furthermore since $\|H_\cdot\|_\infty\in L^1(0,T)$ for each $\ell\in\ZZ_+$ there exists $\delta_\ell>0$ such that 
$$\LL_{(0,T)}(E)<\delta_\ell\quad\Lra\quad\int_E\|H_t\|_\infty\df t<\ell^{-1}.$$ Since $\LL_{(0,T)}(\{\bar{\delta}_{\ell^{-1}}(\cdot)=0\})=0$ for all $\ell\in\ZZ_+$, for each $\ell\in\ZZ_+$ there exists $k_\ell\in\NN$ such that $\LL_{(0,T)}(\{\bar{\delta}_{\ell^{-1}}(\cdot)<\fr{k_\ell}\})<\delta_\ell$. Then for all $N,\ell$
\begin{align}\label{OnSmallSet}
\Big|\int_{\{\bar{\delta}_{\ell^{-1}}(\cdot)<k_\ell^{-1}\}}I_t^{N,\ell,H,a}\df t\Big|&\leq 2|a|_1\int_{\{\bar{\delta}_{\ell^{-1}}(\cdot)<k_\ell^{-1}\}}\|H_t\|_\infty\fr{N^d}\sum_{x\in\T_N^d}\eta^\ell_t(x)\df t\nonumber\\
&\stackrel{P^N\mbox{-a.s.}}\leq 2|a|_1\ls 1,\pi_0^N\rs\int_{\{\bar{\delta}_{\ell^{-1}}(\cdot)<k_\ell^{-1}\}}\|H_t\|_\infty\df t\nonumber\\
&\leq2|a|_1\ell^{-1}\ls 1,\pi_0^N\rs\LL_{(0,T)}(\{\bar{\delta}_{\ell^{-1}}(\cdot)<k_\ell^{-1}\}).
\end{align}
On the other hand, setting $b:=\sup_{z\in J}|z|<+\infty$ for each $\ell\in\ZZ_+$ we can choose $N_\ell\in\NN$ such that $b/N<k_\ell^{-1}$ for all $N\geq N_\ell$ and then for all $\ell\in\ZZ_+$, $N\geq N_\ell$ 
\begin{align}\label{OnNotSmallSet}\Big|\int_{\{\bar{\delta}_{\ell^{-1}}(\cdot)\geq k_\ell^{-1}\}}I_t^{N,\ell,H,a}\df t\Big|&\leq\int_{\{\bar{\delta}_{\ell^{-1}}(\cdot)\geq k_\ell^{-1}\}}\fr{N^d}\sum_{x\in\T_N^d}\sum_{z\in J}|a_z|\Big|H_t\Big(\frac{x-z}{N}\Big)-H_t\Big(\frac{x}{N}\Big)\Big|\eta_t^\ell(x)\df t\nonumber\\
&\leq|a|_1\ell^{-1}\int_{\{\bar{\delta}_{\ell^{-1}}(\cdot)\geq k_\ell^{-1}\}}\fr{N^d}\sum_{x\in\T_N^d}\eta_t^\ell(x)\df t\nonumber\\
&\stackrel{P^N\mbox{-a.s.}}=|a|_1\ell_\star^{-1}\ls\pi_0^N,1\rs\LL_{(0,T)}(\{\bar{\delta}_{\ell^{-1}}(\cdot)\geq k_\ell^{-1}\}).
\end{align}It follows by~\eqref{OnSmallSet} and~\eqref{OnNotSmallSet} that for each $\ell\in\ZZ_+$ and all $N\geq N_\ell$
$$\EE^N\Big|\int_0^TI_t^{N,\ell,H,a}\df t\Big|\leq 2|a|_1\ell_\star^{-1}\EE^N\ls\pi_0^N,1\rs=2|a|_1\ell^{-1}\int\ls 1,\pi^N\rs\df\mu_0^N.$$
Therefore by the $O(N^d)$-entropy assumption and Lemma~\ref{BoundedTotParticlInMacrLim}
$$\limsup_{\ell,N\ra+\infty}\EE^N\Big|\int_0^TI_t^{N,\ell,H,a}\df t\Big|=0.$$

Then, since by~\cite[Lemma 5.5.3]{Kipnis1999a} the one-block estimate holds for sublinear cylinder functions in the case of weakly condensing ZRPs, the one-block estimate follows for weakly condensing ZRPs and in order to prove the one-block estimate for condensing ZRPs it suffices to extend the one-block estimate for condensing ZRPs in~\cite{Stamatakis2014} to sublinear cylinder maps.

 We proceed now with the proof of the one-block estimate for sublinear cylinder maps in the case of condensing ZRPs. This is based on the following consequence of the equivalence of ensembles. In fact since we have reduced the one-block estimate to sublinear cylinder functions we need the next result only for sublinear maps. We state for the more general case of asymptotically linear cylinder maps $\Psi$ to elucidate the definition of the associated extended homologue $\bbar{\Psi}$.

\begin{lemma}\label{WassersteinConvergenceEoE}{\rm{(Super-critical equivalence of ensembles)}} Let $\Psi\colon\MM_\infty^d\to\RR$ be an asymptotically linear cylinder map with support $J_\Psi\subs\ZZ^d$ and slope $\nabla\Psi(\infty)\in\RR^{J_\Psi}$ at infinity. Then 
	$$\lim_{\substack{N,K\ra+\infty\\K/N^d\ra\rho}}\int \Psi d\nu_{N,K}=\int \Psi\df\nu_{\rho\mn\rho_c}^\infty+\ls\nabla\Psi(\infty),\1_{J_\Psi}\rs(\rho-\rho_c)^+=\bbar{\Psi}(\rho),\quad\rho\geq 0.$$
\end{lemma}\textbf{Proof} By~\cite{Grosskinsky2003a} we know that the claim holds when $\Psi$ is bounded. First we will show by a truncation argument that it also holds for sublinear cylinder functions. So let $\Psi$ be a cylindric sublinear map with support $J\subs\ZZ^d$ and let $\Psi_J\colon\ZZ_+^J\to\RR$ be such that $\Psi=\Psi_J(\eta_J)$ and 
$$\lim_{|\eta_J|_1\ra+\infty}\frac{|\Psi_J(\eta_J)|}{|\eta_J|_1}=0.$$ Since $\Psi_J$ is sublinear there exists a constant $C_\Psi<+\infty$ such that $|\Psi|\leq C_\Psi(1+|\eta_J|_1)$. Thus since $\nu_\rho^1$ has exponential moments for $\rho\in[0,\rho_c)$ and $\nu_{\rho_c}^1$ has finite moments whenever $\rho_c<+\infty$ it follows that $\Psi_J\in\cap_{\rho\in[0,\rho_c]\cap\RR_+}L^1(\ZZ_+^J,\nu_\rho^J)$.

For each $\delta>0$ we pick $K_\delta>0$ such that $|\Psi_J(\eta)|<\delta|\eta|_1$ for all $\eta\in\MM_\infty^d$ with $|\eta|_1>K_\delta$ and $K_\delta\ra+\infty$ as $\delta\ra 0$, and decompose $\Psi_J$ as $\Psi_J=\Psi^{\leq}_{J,\delta}+\Psi^{>}_{J,\delta}$ where 
$\Psi_{J,\delta}^{\leq}:=\Psi_J\1_{[0,K_\delta]}(|\cdot|_1)$ and $\hat{\Psi}_{J,\delta}^{>}:=\Psi_J\1_{(K_\delta,+\infty)}(|\cdot|_1)$. Thus $\Psi=\Psi_{J,\delta}^{\leq}(\eta_J)+\Psi^{>}_{J,\delta}(\eta_J)=:\Psi_\delta^{\leq}+\hat{\Psi}^{>}_\delta$ so that $|\Psi-\Psi_\delta^{\leq}|=|\Psi^{>}_\delta|<\delta|\eta^J|_1$, and all $N\in\NN$ large enough so that $J\subs\T_N^d$ 
$$-\delta(\sharp J)\frac{K}{N^d}+\int\Psi_\delta\df \nu_{N,K}\leq\int hd\nu_{N,K}\leq\int\Psi_\delta\df\nu_{N,K}+\delta(\sharp J)\frac{K}{N^d}.$$
By its definition the map $\Psi_{J,\delta}$ is bounded for each $\delta>0$ and therefore 
$$\lim_{\substack{N,K\ra+\infty\\{K/N^d}\ra\rho}}\int\Psi_\delta\df\nu_{N,K}=\int\Psi_\delta d\nu_{\rho\mn\rho_c}^1.$$ Thus for each fixed $\delta>0$ 
\begin{align*}
-\delta\sharp J\rho+\int\Psi_\delta\df\nu_{\rho\mn\rho_c}^1&\leq\liminf_{\substack{N,K\ra+\infty\\{K/N^d}\ra\rho}}\int\Psi\df\nu_{N,K}\leq\limsup_{\substack{N,K\ra+\infty\\{K/N^d}\ra\rho}}\int\Psi\df\nu_{N,K}\\
&\leq\int\Psi_\delta\df\nu_{\rho\mn\rho_c}^1+\delta\sharp J\rho
\end{align*}
and passing to the limit as $\delta\ra 0$ we obtain 
\begin{align}\label{Tr11}\liminf_{\delta\ra 0}\int\Psi_\delta\df\nu_{\rho\mn\rho_c}^\infty&\leq\liminf_{\substack{N,K\ra+\infty\\{K/N^d}\ra\rho}}\int hd\nu_{N,K}\leq\limsup_{\substack{N,K\ra+\infty\\K/N^d\ra\rho}}\int\Psi\df\nu_{N,K}\nonumber\\
&\leq\limsup_{\delta\ra 0}\int\Psi_\delta\df\nu_{\rho\mn\rho_c}^\infty.\end{align}
Now, $\int\Psi_\delta\df\nu_{\rho\mn\rho_c}^\infty=\int\Psi_{J,\delta}d\nu_{\rho\mn\rho_c}^J$
and $|\Psi_{J,\delta}|\leq |\Psi_J|\in L^1(\nu_{\rho\mn\rho_c}^J)$, and therefore since $\Psi_{J,\delta}\lra \Psi_J$ pointwise as $\delta\ra 0$ it follows by the dominated convergence theorem that $$\lim_{\delta\ra 0}\int\Psi_\delta\df\nu_{\rho\mn\rho_c}^\infty=\int\Psi\df\nu_{\rho\mn\rho_c}^\infty$$
The claim for sublinear cylinder functions follows by this limit and~\eqref{Tr11}.

We suppose next that the cylinder map $\Psi$ is asymptotically linear. Then with $J=J_\Psi$ being the support of $\Psi$ $$\Psi=\Psi_J(\eta^J)=(\Psi_J(\eta^J)-\ls\nabla\Psi(\infty),\eta^J\rs)+\ls\nabla\Psi(\infty),\eta^J\rs:=\Psi_{J,0}+\ls\nabla\Psi(\infty),\eta^J\rs$$
with $\Psi_{J,0}$ being a sublinear cylinder map. Therefore 
\begin{align*}
\lim_{\substack{N,K\ra+\infty\\K/N^d\ra\rho}}\int\Psi\df\nu_{N,K}&=\wt{\Psi_{J,0}}(\rho\mn\rho_c)+\lim_{\substack{N,K\ra+\infty\\K/N^d\ra\rho}}\int\ls\nabla\Psi_J(\infty),\eta^J\rs\df\nu_{N,K}\\
&=\wt{\Psi_{0,J}}(\rho\mn\rho_c)+\sum_{x\in J}\pd_x\Psi(\infty)\lim_{\substack{N,K\ra+\infty\\K/N^d\ra\rho}}\int\eta(x)\df\nu_{N,K}\\
&=\wt{\Psi}(\rho\mn\rho_c)-\ls\nabla\Psi(\infty),\1_J\rs(\rho\mn\rho_c)+\ls\nabla\Psi(\infty),\1_J\rs\rho\\
&=\wt{\Psi}(\rho\mn\rho_c)+\ls\nabla\Psi(\infty),\1_J\rs(\rho-\rho_c)^+=\bbar{\Psi}(\rho),
\end{align*}
which completes the proof of the lemma.$\hfill\Box$

\begin{lemma}{\rm{(Uniform $L^1$ Law of Large numbers)}} Let $\mathcal{P}$ be a family of probability measures on a measurable space $(\W,\F)$ and let $(X_i)_{i\in\NN}$ be a sequence of random variables that is i.i.d.~with respect to all $P\in\mathcal{P}$. We denote by $\EE_P$ the expectation with respect to $P\in\mathcal{P}$ and set $\mu_P:=\EE_PX_1$. If $X_1$ is $\mathcal{P}$-uniformly integrable, i.e.~$\lim_{M\ra+\infty}\sup_{P\in\mathcal{P}}\EE_P\big(|X_1|\1_{\{|X_1|>M\}}\big)=0$
	then
	$$\lim_{n\ra+\infty}\sup_{P\in\mathcal{P}}\EE_P\Big|\fr{n}\sum_{i=1}^nX_i-\mu_P\Big|=0.$$
\end{lemma}\textbf{Proof} The claim obviously holds if the $X_i$'s have uniformly bounded variance, i.e.~if
$$\sup_{P\in\mathcal{P}}\EE_P|X_1-\mu_P|^2<+\infty.$$ The claim then follows by a truncation argument. For each $M>0$ set $X_i^M:=X_i\1_{\{|X_i|\leq M\}}$, $\bar{X}_i^M:=X_i\1_{\{|X_i|>M\}}$ and $\mu_P^M:=\EE_PX_1^M$, $\bar{\mu}_P^M:=\EE_P\bar{X}_1^M$. Then 
$X_i=X_i^M+\bar{X}_i^M$ and $\mu_P=\mu_P+\bar{\mu}_P^M$ for each $M>0$ and therefore 
$$\EE_P\Big|\fr{n}\sum_{i=1}^nX_i-\mu_P\Big|\leq\EE_P\Big|\fr{n}\sum_{i=1}^nX_i^M-\mu_P^M\Big|+\EE_P\Big|\fr{n}\sum_{i=1}^n\bar{X}_i^M-\bar{\mu}_P^M\Big|.$$
The random variables $X_i^M$ have obviously uniformly bounded variance for each $M>0$ and thus 
\begin{align*}\limsup_{n\ra+\infty}\sup_{P\in\mathcal{P}}\EE_P\Big|\fr{n}\sum_{i=1}^nX_i-\mu_P\Big|&\leq\limsup_{n\ra+\infty}\sup_{P\in\mathcal{P}}\EE_P\Big|\fr{n}\sum_{i=1}^n\bar{X}_i^M-\bar{\mu}_P^M\Big|\\
&\leq 2\sup_{P\in\mathcal{P}}\EE_P\big(|X_1|\1_{|X_1|>M}\big)\stackrel{M\ra+\infty}\lra 0,\end{align*}
by the $\mathcal{P}$-uniform integrability of the $X_i$'s.$\hfill\Box$\\

Using these two lemmas we can extend the one-block estimate of condensing ZRPs to sublinear cylinder maps.  As shown in \cite{Kipnis1999a}[Section 5.2] whenever the sequence $\{\mu_0^N\}$ satisfies the $O(N^d)$-entropy assumption with respect to some equilibrium state $\{\nu_{\rho_*}^N\}$, $\rho_*\in(0,\rho_c)$ with constant $C_0=C_0(\rho_*)$ then entropy and normalized Dirichlet form of the density
$\bar{f}_T^N:={d\bar{\mu}_T^N}/{d\nu_{\rho_*}^N}$ of the time averaged law $\bar{\mu}_T^N:=\fr{T}\int_0^T\mu_t^N\df t$ with respect to $\nu_{\rho_*}^N$ satisfy the bounds 
\begin{equation}\label{EntropyDirichletFormBounds}
H_N(\bar{f}_T^N)\leq C(\rho_*)N^d\quad\mbox{and}\quad D_N(\bar{f}_T^N)\leq\frac{C(\rho_*)}{2T}N^{d-2}
\end{equation}
Here the supremum is taken among all densities $f\in L^1_+(\nu_{\bs{\rho}_*}^N)$, we have set $H_N(f):=\HHH(f\df\nu_{\rho_*}^N|\nu_{\rho_*}^N)$ for all $f\in L^1_+(\nu_{\rho_*}^N)$ and the renormalized Dirichlet form $D_N\colon L^1_+(\nu_{\rho_*}^N)\to[0,+\infty]$
is given by $D_N(f)=\mathfrak{D}_N(\sqrt{f})$ where
$\mathfrak{D}_N\colon L^2(\nu_{\bs{\rho}_*})\to[0,+\infty]$ is the Dirichlet form associated to the generator $L_N$,
\begin{equation*}
\mathfrak{D}_N(f):=-\int fL_Nf\df\nu_{\bs{\rho}_*}=\fr{2}\sum_{\eta\in\MM_N^d}\sum_{x,y\in\T_N^d}\big(f(\eta^{x,y})-f(\eta)\big)^2\mathfrak{g}(\eta(x))p(y-x)\nu_{\rho_*}^N(\eta).
\end{equation*} 
Via these estimates the one-block estimate is reduced to proving that 
\begin{equation}
\label{StaticOBE}
\limsup_{\ell\ra\infty}\limsup_{N\ra\infty}\sup_{\substack{H_N(f)\leq C_0N^d\\D_N(f)\leq 
		C_0N^{d-2}}}\int\fr{N^d}\sum_{x\in\T_N^d}\tau_xV_\Psi^\ell f\df\nu_{\rho_*}^N\leq 0,\quad\forall\;C_0>0,
\end{equation}
Since the cylinder map is assumed here to be sublinear one can follows that steps 1 to 5 in~\cite[Sect. 5.4.1]{Kipnis1999a}to further reduce the one-block estimate to showing that for all constants $C_1>0$,
\begin{align}
\label{step5end}
\limsup_{\ell\ra+\infty}\max_{K|K\leq(2\ell+1)^dC_1}\int V^\ell d\nu_{2\ell+1,K}=0,
\end{align}
where the canonical measure $\nu_{2\ell+1,K}$ is considered as a measure on $\MM_\infty^d$ by identifying the cube
$\Lambda_\ell^d:=\{x\in\ZZ^d\bigm||x|\leq\ell\}\subs\ZZ^d$ with $\T_{2\ell+1}^d$. By fixing a positive integer $k$ which will tend to infinity after taking the limit as
$\ell\ra+\infty$, and decomposing the cube $\Lambda_\ell^d$ in smaller cubes of side-length $2k+1$, the one-block estimate is
reduced to showing that
\begin{align}
\label{FinalLimInOBE}	
\lim_{k\ra\infty}\lim_{m\ra\infty}S(m,k)=0,
\end{align}
where $S(m,k)$ denotes the supremum
\begin{equation*}S(m,k):=\sup_{\substack{\ell\geq m\\K\leq(2\ell+1)^dC_1}}
\int\Big|\fr{(2k+1)^d}\sum_{|x|\leq k}\tau_x\Psi-
\bbar{\Psi}\Big(\frac{K}{(2\ell+1)^d}\Big)\Big|_1\df\nu_{2\ell+1,K}.
\end{equation*} 
For each fixed $(m,k)\in\NN\x\NN$, we pick a sequence $\{(\ell_n^{m,k},K^{m,k}_n)\}_{n\in\NN}$ such that $\ell_n^{m,k}\geq m$ and
$K_n^{m,k}\leq(2\ell_n^{m,k}+1)^dC_1$ for all $n\in\NN$ that achieves the supremum, i.e.~such that
\begin{equation*}	
S(m,k)=\lim_{n\ra\infty}\int\Big|\fr{(2k+1)^d}\sum_{|x|\leq k}\tau_x\Psi-
\bbar{\Psi}\Big(\frac{K^{m,k}_n}{(2\ell_n^{m,k}+1)^d}\Big)\Big|_1\df\nu_{2\ell_n^{m,k}+1,K^{m,k}_n}.
\end{equation*}
Since the sequence $\{r_n^{m,k}\}_{n\in\NN}$ defined by
\begin{equation*}
r_n^{m,k}:=\frac{K^{m,k}_n}{(2\ell_n^{m,k}+1)^d},\qquad n\in\NN,
\end{equation*} 
is contained in the interval $[0,C_1]$, for each
fixed $(m,k)\in\NN\x\NN$, we can pick a sequence $\{n_j\}_{j\in\NN} := \{n_j^{m,k}\}$ such that $r^{m,k}_{n_j}$ converges to
some $r^{m,k}\in [0,C_1]$ as $j\ra\infty$. Since we assume that $\Psi$ is sublinear, it follows by the
equivalence of ensembles that
\begin{equation*}
S(m,k)=\int\Big|\fr{(2k+1)^d}\sum_{|x|\leq k}\tau_x\Psi-
\bbar{\Psi}\big(r^{m,k}\big)\Big|_1\df\nu_{r^{m,k}\mn\rho_c}^\infty.
\end{equation*}
Furthermore $r^{m,k}$ is also contained in $[0,C_1]$ and thus we can choose for each $k\in\NN$ a
sequence $\{m_j\}_{j\in\NN}\equiv\{m_j^{(k)}\}$ such that $\{r^{m_j,k}\}_{j\in\NN}$ converges to some
$r^k\in[0,C_1]$. Then $\lim_{j\uparrow\infty}r^{m_j,k}\mn\rho_c=r^k\mn\rho_c$ and since the grand canonical ensemble $(\nu_{\rho\mn\rho_c}^\infty)_{\rho\geq 0}\subs\PP_1\MM_\infty^d$
is continuous in the Wasserstein topology of $1$-st order and $\Psi$ is sublinear
\begin{align*}
\lim_{m\ra\infty}S(m,k)=\int\Big|\fr{(2k+1)^d}\sum_{|x|\leq k}\tau_x\Psi-
\wt{\Psi}(r^k\mn\rho_c)\Big|\df\nu_{r^k\mn\rho_c}^\infty.
\end{align*}
Therefore
\begin{equation*}
\limsup_{k\ra+\infty}\lim_{m\ra+\infty}S(m,k)\leq
\limsup_{k\ra\infty}\sup_{\rho\in[0,\rho_c]}\int\Big|\fr{(2k+1)^d}
\sum_{|x|\leq k}\tau_x\Psi-
\wt{\Psi}(\rho)\Big|\df\nu_{\rho}^\infty.
\end{equation*}

Now let $\ell_\Psi$ be the smallest integer $\ell\in\ZZ_+$ such that $J_\Psi\subs\Lambda_\ell^d$ and set 
$$E^k_x:=x+(2\ell_\Psi+1)\ZZ^d\cap\Lambda_k^d,\quad x\in\Lambda_{\ell_\Psi}^d.$$
Then the sets $E_x^k$, $x\in\Lambda_{\ell_\Psi}^d$ are disjoint and the cube $\Lambda_k^d$ is equal to the disjoint union 
$$\Lambda_k^d=\bigsqcup_{x\in\Lambda_{\ell_\Psi}^d}E_x^k.$$
Since $k$ is an arbitrary parameter that is introduced by splitting a larger cube of radius $\ell$ into smaller ones of radius $k$ and $k$ tends to infinity after $\ell$ is sent to infinity, we can assume without loss fo generality that $2\ell_\Psi+1$ divides $2k+1$ so that 
$$\sharp E_x^k=\frac{(2k+1)^d}{(2\ell_\Psi+1)^d}\in\NN,\quad\forall x\in\Lambda_{\ell_\Psi}^d.$$ 
Then 
$$\fr{(2k+1)^d}
\sum_{|x|\leq k}\tau_x\Psi=\fr{(2\ell_\Psi+1)^d}\sum_{x\in\Lambda_{\ell_\Psi}^d}\fr{\sharp E_x^k}\sum_{y\in E_x^k}\tau_y\Psi$$
and thus for each $\rho\in[0,\rho_c]$
$$\fr{(2k+1)^d}
\sum_{|x|\leq k}\tau_x\Psi-
\wt{\Psi}(\rho)=\fr{(2\ell_\Psi+1)^d}\sum_{x\in\Lambda_{\ell_\Psi}^d}\fr{\sharp E_x^k}\sum_{y\in E_x^k}\big(\tau_y\Psi-\wt{\Psi}(\rho)\big).$$
Therefore
\begin{align*}
\sup_{\rho\in[0,\rho_c]}\int\Big|\fr{(2k+1)^d}&
\sum_{|x|\leq k}\tau_x\Psi-
\wt{\Psi}(\rho)\Big|\df\nu_{\rho}^\infty\\
&\leq\fr{(2\ell_\Psi+1)^d}\sum_{x\in\Lambda_{\ell_\Psi}^d}\sup_{\rho\in[0,\rho_c]}\int\Big|\fr{\sharp E_x^k}\sum_{y\in E_x^k}\tau_y\Psi-\wt{\Psi}(\rho)\Big|\df\nu_{\rho}^\infty\end{align*}
and thus in order to complete the proof of the one-block estimate it suffices to show that 
$$\lim_{k\ra+\infty}\sup_{\rho\in[0,\rho_c]}\int\Big|\fr{\sharp E_x^k}\sum_{y\in E_x^k}\tau_y\Psi-\wt{\Psi}(\rho)\Big|\df\nu_{\rho}^\infty=0$$
for each $x\in\Lambda_{\ell_\Psi}^d$.
So let $x\in\Lambda_{\ell_\Psi}^d$. The random variables $\{\tau_y\Psi\}_{y\in E_x^k}$ are independent and identically distributed and by the uniform $L^1$-law of large numbers it suffices to show that $\tau_x\Psi$ is $\{\nu_\rho^\infty\}_{\rho\in[0,\rho_c]}$-uniformly integrable. Since $\Psi$ is sublinear it has at most linear growth and thus there exists a constant $C\geq 0$ such that $\Psi\leq C(1+|\eta^{J_\Psi}|_1)$. Then for each $x\in E_x^k$, $\rho\in[0,\rho_c]$ and all $M>0$
\begin{align*}\int|\tau_x\Psi|\1_{\{|\tau_x\Psi|>M\}}\df\nu_\rho^\infty&
=\int\tau_x(|\Psi|\1_{\{|\Psi|>M\}})\df\nu_\rho^\infty=\int|\Psi|\1_{\{|\Psi|>M\}}\df\nu_\rho^\infty\\
&\leq C\int(1+|\eta^{J_\Psi}|_1)\1_{(M,\infty)}(1+|\eta^{J_\Psi}|_1)\df\nu_\rho^{\infty}(\eta).
\end{align*}
But the map $\Psi_M\colon\MM_\infty^d\to\RR$ given by $\Psi_M(\eta)=C(1+|\eta^{J_\Psi}|_1)\1_{(M,\infty)}(1+|\eta^{J_\Psi}|_1)$ is increasing, i.e.~$\Psi_M(\eta)\leq\Psi_M(\zeta)$ whenever $\eta(x)\leq\zeta(x)$ for all $x\in\ZZ^d$ and by \cite[Lemma 2.3.5]{Kipnis1999a} the grand canonical ensemble ensemble $\{\nu_\rho^1\}_{\rho\in[0,\rho_c]}$ is a stochastically increasing function of the parameter $\rho\in[0,\rho_c]$. Therefore
$$\int|\tau_x\Psi|\1_{\{|\tau_x\Psi|>M\}}\df\nu_\rho^\infty\leq\int\Psi_M\df\nu_\rho^\infty\leq\int\Psi_M\df\nu_{\rho_c}^\infty$$ 
for all $\rho\in[0,\rho_c]$ and thus 
$$\lim_{M\ra+\infty}\sup_{\rho\in[0,\rho_c]}\int|\tau_x\Psi|\1_{\{|\tau_x\Psi|>M\}}\df\nu_\rho^\infty\leq\lim_{M\ra+\infty}\int\Psi_M\df\nu_{\rho_c}^\infty=0.$$
This shows that we can apply the uniform $L^1$-law of large numbers and completes the proof of~\eqref{OBE}. 

We prove next claim (b). In order to prove that the family $\{\bbar{\bs{Q}}_\Psi^{N,\ell}\}_{(N,\ell)\in\NN\x\ZZ_+}$ is relatively compact it suffices to show that the sequences of its marginals $Q_{\Psi,1}^N:=\bs\pi^{N,\ell}_\sharp P^N$ and $\bs{Q}_{\Psi,2}^{N,\ell}:=\s^{N,\Psi}_\sharp P^N$, $N\in\NN$, are relatively compact. The relative compactness of $\{\bs{Q}_{\Psi,1}^{N,\ell}\}_{(N,\ell)\in\NN\x\ZZ_+}$ has been proved in Proposition~\ref{YoungRelComp}. The relative compactness of $\{Q_{\Psi,2}^N\}$ is proved in the following.
\begin{prop}\label{PsiEmpDistComp} Let $\Psi\colon\MM_\infty^d\to\RR$ be an asymptotically linear cylinder map. If the sequence of initial distributions has total mass bounded above in probability then the sequence of distributions $Q_{\Psi,2}^N:=\s^{N,\Psi}_\sharp P^N\in\PP L_{w^*}^\infty(0,T;\MMM(\T^d))$ is relatively compact.  
\end{prop}\textbf{Proof} By the Prokhorov-Le Cam and Banach-Alaoglu theorems it suffices to show that $$\lim_{M\ra+\infty}\sup_{N\in\NN}
Q_{\Psi,2}^N\big\{\|\s\|_{L_{w^*}^\infty(0,T;\MMM(\T^d))}>M\big\}=0.$$ 
Since in any case the cylinder map $\Psi$ satisfies $|\Psi|\leq C(1+|\eta_J|_1)$ for some constant $0<C<+\infty$ and some finite $J\subs\ZZ^d$ 
\begin{align*}
\|\s^{N,\Psi}\|_{L_{w^*}^\infty(0,T;\MMM(\T^d))}&=\sup_{\|f\|_{L^1(0,T;C(\T^d))}\leq 1}\lls f,\s^{N,\Psi}\rrs\\
&=\sup_{\|f\|\leq 1}\int_0^T\fr{N^d}\sum_{x\in\T_N^d}f_t\Big(\frac{x}{N}\Big)\tau_x\Psi(\eta_t)\df t\\
&\leq C\sup_{\|f\|\leq 1}\int_0^T\|f_t\|_\infty\fr{N^d}\sum_{x\in\T_N^d}\Big(1+\sum_{y\in J}\tau_x\eta_t(y)\Big)\df t\\
&\stackrel{P^N\mbox{-a.s.}}=C(1+\sharp J\ls 1,\pi_0^N\rs)\sup_{\|f\|\leq 1}\int_0^T\|f_t\|_\infty\df t\leq C(1+\sharp J\ls 1,\pi_0^N\rs).
\end{align*}
Therefore 
\begin{align*}
Q_{\Psi,2}^N\big\{\|\s\|_{L_{w^*}^\infty(0,T;\MMM(\T^d))}>M\big\}&=P^N\big\{\|\s^N\|_{L_{w^*}^\infty(0,T;\MMM(\T^d))}>M\big\}\\
&\leq\mu_0^N\Big\{\ls 1,\pi_0^N\rs>\fr{\sharp J}\Big(\frac{M}{C}-1\Big)\Big\}
\end{align*}
and thus the tightness of $\{Q_{\Psi,2}^N\}_{N\in\NN}$ follows by Lemma~\ref{BoundedTotParticlInMacrLim}.$\hfill\Box$\\

We show next that $\bbar{\Psi}\in\bbar{C}_1(\RR_+)$ with $\bbar{\Psi}'(\infty)=\ls\nabla\Psi(\infty),\1_J\rs$. If the ZRP is strongly condensing this claim is obvious and thus we can assume that $\rho_c=+\infty$. We will show first that if $\Psi$ is sublinear then $\bbar{\Psi}'(\infty)=0$. So let $\ee>0$. Since $\Psi$ is sublinear there exists $K_\ee>0$ such that $|\Psi(\eta)|<\ee|\eta_J|_1$ whenever $|\eta_J|_1>K_\ee$. Then for all $\rho>0$
$$\Big|\frac{\bbar{\Psi}(\rho)}{\rho}\Big|\leq\fr{\rho}\int|\Psi|\df\nu_\rho^\infty=\fr{\rho}\int_{\{|\eta_J|_1\leq K_\ee\}}|\Psi|\df\nu_\rho^\infty+\fr{\rho}\int_{\{|\eta_J|_1>K_\ee\}}|\Psi|\df\nu_\rho^\infty$$
and therefore since $|\Psi|$ is bounded on $\{|\eta_J|_1\leq K_\ee\}$
$$\limsup_{\rho\ra+\infty}\Big|\frac{\bbar{\Psi}(\rho)}{\rho}\Big|\leq\limsup_{\rho\ra+\infty}\Big(\fr{\rho}\int_{\{|\eta_j|_1>K_\ee\}}|\Psi|\df\nu_\rho^\infty\Big).$$
But $$\int_{\{|\eta_J|_1>K_\ee\}}|\Psi|\df\nu_\rho^\infty\leq\ee\int|\eta_J|_1\df\nu_\rho^\infty=\ee\sharp J\rho$$
and therefore $\limsup_{\rho\ra+\infty}\big|\frac{\bbar{\Psi}(\rho)}{\rho}\big|\leq\ee\sharp J$ which since $\ee>0$ was arbitrary proves that $\bbar{\Psi}'(\infty)=0$. Now if $\Psi$ is asymptotically linear it is of the form $\Psi=\Psi_0+\ls\nabla\Psi(\infty),\eta_J\rs$ for some sublinear cylinder map $\Psi_0$ and $\bbar{\Psi}(\rho)=\bbar{\Psi_0}(\rho)+\rho\ls\nabla\Psi(\infty),\1_J\rs$ and thus $\bbar{\Psi}'(\infty)=\bbar{\Psi_0}'(\infty)+\ls\nabla\Psi(\infty),\1_J\rs=\ls\nabla\Psi(\infty),\1_J\rs$.

We prove finally~\eqref{OBEYM}. Since $L^1(0,T;C(\T^d))$ is separable it suffices to show that 
$$\bbar{\bs{Q}}_\Psi\big\{\big|\lls G,\s-\bar{B}_{\bbar{\Psi}}(\bs\pi)\rrs\big|>\ee\big\}=0$$
for all $G\in L^1(0,T;C(\T^d))$ and $\ee>0$. Since $\bbar{\bs{Q}}_\Psi$ is a limit point of $\{\bbar{\bs{Q}}_\Psi^{N,\ell}\}_{(N,\ell)}$ as $N\ra+\infty$ and then $\ell\ra+\infty$ there exist a sequence $\{m_\ell\}_{\ell\in\NN}\subs\ZZ_+$ and diverging sequences $\{k_N^{(\ell)}\}_{N\in\NN}\subs\NN$, $\ell\in\ZZ_+$ such that 
$$\bbar{\bs{Q}}_\Psi=\lim_{\ell\ra+\infty}\lim_{N\ra+\infty}\bbar{\bs{Q}}_\Psi^{k_N^{(\ell)},m_\ell}.$$
Since the map $\bar{B}_{\bbar{\Psi}}$ is $w^*$-continuous the set $\{|\s-\bar{B}_{\bbar{\Psi}}(\bs\pi)|>\ee\}$ is open and thus by the portmanteau theorem 
\begin{align*}
\bbar{\bs{Q}}_\Psi\big\{\big|\lls G,\s-\bar{B}_{\bbar{\Psi}}(\bs\pi)\rrs\big|>\ee\big\}&\leq\liminf_{\ell\ra+\infty}\liminf_{N\ra+\infty}\bbar{\bs{Q}}_\Psi^{k_N^{(\ell)},m_\ell}\big\{\big|\lls G,\s-\bar{B}_{\bbar{\Psi}}(\bs\pi)\rrs\big|>\ee\big\}\nonumber\\
&\leq\limsup_{\ell\ra+\infty}\limsup_{N\ra+\infty}\bbar{\bs{Q}}_\Psi^{N,\ell}\big\{\big|\lls G,\s-\bar{B}_{\bbar{\Psi}}(\bs\pi)\rrs\big|>\ee\big\}\nonumber\\
&=\limsup_{\ell\ra+\infty}\limsup_{N\ra+\infty}P^{N,\ell}\big\{\big|\lls G,\s^{N,\Psi}-\bar{B}_{\bbar{\Psi}}(\bs\pi^{N,\ell})\rrs\big|>\ee\big\}.
\end{align*}
Since 
$$\lls G,\bbar{B}_{\bbar{\Psi}}(\bs\pi^{N,\ell})\rrs=\lls G(U)\bbar{\Psi}(\Lambda),\bs\pi^{N,\ell}\rrs=\int_0^T\fr{N^d}\sum_{x\in\T_N^d}G_t\Big(\frac{x}{N}\Big)\bbar{\Psi}(\eta^\ell_t(x))\df t$$
the claim follows by (a).$\hfill\Box$
\subsection{The continuity equation}\label{ContinuityEquationSection}
In order to prove the relative sequential compactness of the sequence $\{Q^N\}_{N\in\NN}$ it suffices show that each one of the sequences of its marginals on the spaces $D(0,T;\MMM_+(\T^d))$, $L_{w^*}^\infty(0,T;\X^1(\T^d)^*)$ and $L^\infty_{w^*}(0,T;\MMM(\T^d))$ is relatively compact. The proof is based on the analysis of the martingales associated to the ZRP via the martingale problem, in the spirit of the Guo-Papanikolaou-Varadhan approach~\cite{Guo1988a} to proving hydrodynamic limits. 

For any Banach space we denote by $C^1(0,T;X)$ the space of all continuous curves $F\colon[0,T]\to X$ such that there exists a continuous map $\pd F\colon[0,T]\to X$ such that 
$$\lim_{h\ra 0}\frac{\|F_{t+h}-F_t-h \pd F(t)\|}{h}=0.$$ For any initial distribution $\mu_0^N\in\PP_{2}\MM_N^d$ ($\mu^N\in\PP_1\MM_N^d$ if the jump rate $\mathfrak{g}$ is bounded) and any $G\in C^1(0,T;C(\T^d))$, the real process 
	$$A_t^{N,G}:=\ls G_t,\pi_t^N\rs-\ls G_0,\pi_0^N\rs-\int_0^t(\pd_s+L^N)\ls G_s,\pi^N\rs(\eta^N_s)ds
	,\quad t\geq 0,$$ defined on the filtered probability space $\big(D(\RR_+,\MM_N^d),(\F^N_t)_{t\geq 0},P^N\big)$ is a martingale, where $(\F_t^N)$ is the minimal right continuous filtration to which the ZRP is adapted and $P^N$ is the distribution of the ZRP starting from $\mu_0^N$.

By the definition of the generator $L^N$ of the ZRP for any function $G\in C(\T^d)$
\begin{align*}
L^N\ls G,\pi^N\rs(\eta)&=\fr{N^{d-2}}\sum_{j=1}^d\sum_{x\in\T_N^d}\Big[G\Big(\frac{x+e_j}{N}\Big)+G\Big(\frac{x-e_j}{N}\Big)-2G\Big(\frac{x}{N}\Big)\Big]\mathfrak{g}\big(\eta(x)\big)\\
&=\fr{N^{d-2}}\sum_{j=1}^d\sum_{x\in\T_N^d}\Big[G\Big(\frac{x+e_j}{N}\Big)-G\Big(\frac{x}{N}\Big)\Big]W_{x,x+e_j}(\eta)
\end{align*}
where $W_{x,x+e_j}(\eta)=\mathfrak{g}(\eta(x))-\mathfrak{g}(\eta(x+e_j))$ is the current along the bond $x,x+e_j$. Since $\pd \ls G,\pi^N\rs=\ls\pd G,\pi^N\rs$, the martingale $A^{N,G}$ can be written in more detail as 
\begin{align*}
A_t^{N,G}=\ls G_t,\pi_t^N\rs-\ls G_0,\pi_0^N\rs&-\int_0^t\ls\pd_sG_s,\pi_s^N\rs ds-\int_0^t\fr{N^{d}}
\sum_{x\in\T_N^d}\D^NG_s\Big(\frac{x}{N}\Big)\mathfrak{g}\big(\eta_s(x)\big)\df s\\
=\ls G_t,\pi_t^N\rs-\ls G_0,\pi_0^N\rs&-\int_0^t\ls\pd_sG_s,\pi_s^N\rs ds\\
&-\int_0^t\fr{N^{d-1}}\sum_{j=1}^d\sum_{x\in\T_N^d}\pd_{+j}^NG_s\Big(\frac{x}{N}\Big)W_{x,x+e_j}(\eta_s)\df s
\end{align*}
where for any function $G\colon\T^d\to\RR$ we denote by $\D^NG:\T^d\lra\RR$ the discrete Laplacian 
$$\D^NG(u):=N^2\sum_{j=1}^d\Big[G\Big(u+\frac{e_j}{N}\Big)+G\Big(u-\frac{e_j}{N}\Big)
-2G(u)\Big],\quad u\in\T^d$$
and by 
$$\pd_{+j}^NG(u):=N\Big[G\Big(u+\frac{e_j}{N}\Big)-G(u)\Big],\quad u\in\T^d$$
the discrete right $j$-th partial derivative. Therefore in terms of the empirical jump rate process $\s^N$ and the empirical current process $W^N$ we can write the martingale $A^{N,G}$ as 
\begin{align}\label{VelocityFieldBeta1}
A_t^{N,G}&=\ls G_t,\pi_t^N\rs-\ls G_0,\pi_0^N\rs-\int_0^t\Big[\ls\pd_sG_s,\pi_s^N\rs+\ls\D^NG_s,\s_s^N\rs\Big]\df s\\
&=\ls G_t,\pi_t^N\rs-\ls G_0,\pi_0^N\rs-\int_0^t\Big[\ls\pd_sG_s,\pi_s^N\rs+\ls\nabla_+^NG_s,W_s^N\rs\Big]\df s
\end{align}
where here $\nabla_+^NG(u):=\sum_{j=1}^N\pd_{+j}^NG(u)e_j$ is the discrete right gradient.

As long as the initial sequence $\{\mu_0^N\}$ of initial distributions satisfies $\mu_0^N\in\PP_2\MM_N^d$ ($\mu_0^N\in\PP_1\MM_N^d$ if $\mathfrak{g}$ is bounded) the martingale $A^{N,G}$ is integrable and its quadratic variation is given by 
	\begin{align}\label{MartOfZRProc}\ls A^{N,G}\rs_t&=
		\int_0^t\Big\{L^N(\ls G_s,\pi^N\rs^2)(\eta_s)-2\ls G_s,\pi_s^N\rs L^N\ls G_s,\pi^N\rs(\eta_s)\Big\}ds\nonumber\\
		&=\fr{N^{2d-2}}\int_0^t\sum_{x,y\in\T_N^d}\Big[G_s\Big(\frac{y}{N}\Big)-G_s\Big(\frac{x}{N}\Big)\Big]^2\mathfrak{g}\big(\eta_s(x)\big)p(x,y)ds.
	\end{align}
\vspace{-0.5cm}
\begin{lemma} The martingale $(A_t^{N,G})_{t\geq 0}$ defined in~\eqref{MartOfZRProc} is asymptotically negligible, i.e. for all $\delta>0$
	\begin{equation}\label{MartAsymptNegl}
	\lim_{N\ra+\infty}P^N\Big\{\sup_{0\leq t\leq T}|A_t^{N,G}|\geq\delta\Big\}=0.
	\end{equation}
\end{lemma}\textbf{Proof} Let $\delta>0$. By the Chebyshev and the Burkholder-Davis-Gundy inequality~\cite[Section IV.4]{RevuzYor1999} (which also holds for cadlag martingales for $p\geq 1$) there exists a constant $C<+\infty$ such that
$$P^N\Big\{\sup_{0\leq t\leq T}\big|A_t^{N,G}|\geq\delta\Big\}\leq\fr{\delta^2}\EE^N\Big(\sup_{0\leq t\leq T}|A_t^{N,G}|^2\Big)\leq\frac{C}{\delta^2}\EE^N\ls A^{N,G}\rs_T$$ where $\ls A^{N,G}\rs$ denotes the quadratic variation of the martingale $A^{N,G}$. But by the formula of the quadratic variation $\ls A^{N,G}\rs$, the mean value theorem and the conservation of particles, we have for all $0\leq s\leq t\leq T$
\begin{align}\label{QVForm}
\ls A^{N,G}\rs_t-\ls A^{N,G}\rs_s&=
\fr{N^{2d-2}}\int_s^t\sum_{x,y\in\T_N^d}\Big[G_r\Big(\frac{y}{N}\Big)-G_r\Big(\frac{x}{N}\Big)\Big]^2\mathfrak{g}\big(\eta_r(x)\big)p(x,y)\df r\nonumber\\
&\leq\frac{\|\nabla G\|_{C(\RR_+\x\T^d)}}{N^{2d}}\int_s^t\sum_{x,y\in\T_N^d}|x-y|^2\mathfrak{g}\big(\eta_r(x)\big)p(x,y)\df r\nonumber\\
&\leq\frac{\|\nabla G\|_{C(\RR_+\x\T^d)}\|\mathfrak{g}'\|_\infty}{N^{2d}}\int_s^t\sum_{x,z\in\T_N^d}|z|^2p(z)\eta_r(x)\df r\nonumber\\
&=\frac{2d\|\nabla G\|_{C(\RR_+\x\T^d)}\|\mathfrak{g}'\|_\infty}{N^{d}}\int_s^t\ls 1,\pi^N_r\rs\df r\nonumber\\
&\stackrel{P^N-a.s.}=\frac{2d\|\nabla G\|_{C(\RR_+\x\T^d)}\|\mathfrak{g}'\|_\infty(t-s)}{N^{d}}\ls 1,\pi^N_0\rs.
\end{align}
It follows that 
$$P^N\Big\{\sup_{0\leq t\leq T}\big|A_t^{N,G}|\geq\delta\Big\}\leq\frac{2d\|\nabla G\|_{C(\RR_+\x\T^d)}\|\mathfrak{g}'\|_\infty T}{\delta^2N^{d}}\int\ls 1,\pi^N\rs d\mu_0^N,$$
where $\{\mu_0^N\}$ is the sequence of initial distributions of the ZRP, and taking the limit as $N\ra+\infty$ in the inequality above, it follows by Lemma~\ref{BoundedTotParticlInMacrLim} we obtain the asymptotic negligibility of the martingale $(A_t^{N,G})_{t\geq 0}$.$\hfill\Box$\\

Using Taylor's theorem it follows that for $C^3$ functions we can replace the discrete Laplacian $\D^N$ and the discrete gradient $\nabla_+^N$ by their continuous analogues. More precisely there exists a constant $C=C(G,d,g)\geq 0$ such that 
\begin{equation}\label{HeatContEqatthediscretoruses}
\big|V^{1,N,G}_t-A_t^{N,G}\big|\mx\big|V^{2,N,G}_t-A_t^{N,G}\big|\leq\frac{C}{N}\int_0^t\ls 1,\pi^N_s\rs\df s
\end{equation}
for all $t\geq 0$, where 
\begin{equation}\label{Vjump}
V_t^{1,N,G}:=\ls G_t,\pi_t^N\rs-\ls G_0,\pi_0^N\rs-\int_0^t\ls\pd_sG_s,\pi_s^N\rs+\ls\D G_s,\s_s^N\rs ds
\end{equation}
and 
\begin{equation}\label{Vcurr}
V_t^{2,N,G}:=\ls G_t,\pi_t^N\rs-\ls G_0,\pi_0^N\rs-\int_0^t\ls\pd_sG_s,\pi_s^N\rs+\ls\nabla G_s,W_s^N\rs ds.
\end{equation}
. 

We turn now to the proof of the relative compactness of the first marginal $Q_1^N:=[(\pi_t^N)_{0\leq t\leq T}]_\sharp P^N$ of the sequence $\{Q^N\}_{N\in\NN}$. As we know by the description of the relatively compact subsets of $\PP D(0,T;\MMM_+(\T^d))$ in order to prove that $\{Q_1^N\}$ is relatively compact it suffices to prove that for some countable subset $\{G_k|k\in\NN\}\subs C(\T^d)$ such that $G_1\equiv 1$, the sequence $$\ls G,\cdot\rs_\sharp Q_1^N\in\PP D(0,T;\RR),\quad N\in\NN$$ is relatively compact for all $k\in\NN$, where $\ls G,\cdot\rs\colon D(0,T;\MMM_+(\T^d))\lra D(0,T;\RR)$ is the mapping induced on the Skorohod spaces by the map $\MMM_+(\T^d)\ni\mu\mapsto\int G\df\mu$. In particular it suffices to prove that the sequence $\{\ls G,\cdot\rs_\sharp R_1^N\}$ is relatively compact for all $G\in C^\infty(\T^d)$.\\
\indent So let $G\in C^\infty(\T^d)$. In order to prove the relative compactness of $\{\ls G,\cdot\rs_\sharp Q_1^N\}_{N\in\NN}$ it suffices to prove that: (a) for all $t\in\RR_+$ 
$$\lim_{A\uparrow+\infty}\sup_{N\in\NN}\ls G,\cdot\rs_\sharp Q_1^N\big\{f\in D(\RR_+,\RR)\,\big|\,|f_t|>A\big\}=0$$
and (b) the condition of Aldous, i.e$.$ that for all $\ee,T>0$ we have 
$$\lim_{\delta\ra 0}\limsup_{N\ra+\infty}\sup_{\substack{\tau\in\mathfrak{T}^T(\F^{\RR}_+)\\ \theta\leq\delta}}\ls G,\cdot\rs_\sharp Q_1^N\Big\{f\in D(\RR_+,\RR)\,\Big|\,\big|f_{\tau(f)}-f_{[\tau(f)+\theta]\mn T}\big|>\ee\Big\}=0.$$
Here $\mathfrak{T}^T(\F^\RR_+)$ is the set of all stopping times $\tau\colon D(\RR_+;\RR)\lra[0,T]$ with respect to the continuation $\F^\RR=(\F_{t+}^{0,\RR})_{t\geq 0}$ of the natural filtration $(\F_t^{0,\RR})_{t\geq 0}$ in $D(\RR_+;\RR)$.\\
\noindent(a) Let $t\in\RR_+$. Of course we can assume that $\|G\|_\infty\neq 0$ or else we have nothing to prove, and for all $N\in\NN$ and all $A>0$ we have that 
\begin{align*}
	\ls G,\cdot\rs_\sharp Q_1^N\big\{f\in D(0,T;\RR)\,\big|\,|f_t|>A\big\}&=Q_1^N\big\{\pi\in D(\RR_+;\MMM_+(\T^d))\,\big|\,|\ls G,\pi_t\rs|>A\big\}\\
	&\leq Q_1^N\big\{\pi\in D(0,T;\MMM^d_+)\bigm|\ls\pi_t,1\rs>A/\|G\|_\infty\big\}\\
	&=P^N\big\{\eta\in D(0,T;\MM_N^d)\bigm|\ls 1,\pi^N_t\rs>A/\|G\|_\infty\big\}\\
	&=P^N\big\{\eta\in D(0,T;\MM_N^d)\bigm|\ls 1,\pi^N_0\rs>A/\|G\|_\infty\big\}\\
	&=\mu_0^N\big\{\ls\pi^N,1\rs>A/\|G\|_\infty\big\}.
\end{align*}
It follows that in order to prove (a) it suffices to show that 
\begin{equation}\label{fortotalmassofempiricbounded}
\lim_{A\uparrow+\infty}\sup_{N\in\NN}\mu_0^N\big\{\ls\pi^N,1\rs>A\big\}=0.
\end{equation}
But since the sequence of initial distributions $\{\mu_0^N\}$ is assumed to be associated to a macroscopic profile $\mu_0\in\MMM_+(\T^d)$ it \emph{has total mass $\mathfrak{m}:=\mu_0(\T^d)>0$ in probability} in the sense that 
	$$\lim_{N\ra+\infty}\mu_0^N\big\{\big|\ls 1,\pi^N\rs-\mathfrak{m}\big|>\delta\big\}=0$$
	for all $\delta>0$. Thus~\eqref{fortotalmassofempiricbounded} follows by the next Lemma~\ref{BoundedTotParticlInMacrLim}.

We prove now the Aldous condition (b). So let $\ee, T>0$ be fixed. As we know, given any continuous function $F\colon M\to N$ between polish spaces the induced mapping $\bar{F}\colon D(\RR_+;M)\lra D(\RR_+;N)$ is $(\F_t^M,\F_t^N)$-measurable for all $t\geq 0$, where $(\F^X)$ is the (right) continuation of the natural filtration $(\F_t^{0,X})_{t\geq 0}$ in $D(\RR_+;X)$, $X=M,N$, which shows that 
$$\mathfrak{T}^T(\F^N)\circ\bar{F}:=\big\{\tau\circ\bar{F}\,\big|\,\tau\in\mathfrak{T}^T(\F^N)\big\}\subs\mathfrak{T}^T(\F^M),$$
and we obviously have that $$\bar{F}(x)_{\tau(\bar{F}(x))}=F\big(x_{\tau\circ\bar{F}(x)}\big)\qquad\forall\;x\in D(\RR_+;M),\;\tau\in\mathfrak{T}^T(\F^N).$$ In our particular case we have that $\mathfrak{T}^T(\F^\RR)\circ\bar{F}_G\subs\mathfrak{T}^T(\F^{\MMM_+^d}_+)$ and if for each stopping time $\tau\in\mathfrak{T}^T(\F^\RR)$ we set $\tau_G:=\tau\circ\bar{F}_G$ then $\ls G,\mu\rs_{\tau(\ls G,\pi\rs)}=\ls G,\mu_{\tau_G(\pi)}\rs$ and so 
\begin{align*}
\ls G,\cdot\rs_\sharp Q_1^N\Big\{f\in& D(\RR_+;\RR)\,\Big|\,\big|f_{\tau(f)}-f_{[\tau(f)+\theta]\mn T}\big|>\ee\Big\}\\
&=Q_1^N\Big\{\pi\in D(\RR_+;\MMM_+(\T^d))\,\Big|\,\big|\ls G,\pi_{\tau_G(\pi)}\rs-\ls G,\pi_{[\tau_G(\pi)+\theta]\mn T}\rs\big|>\ee\Big\}
\end{align*}
for all $\theta>0$ and all $\tau\in\mathfrak{T}^T(\F^\RR)$. It follows that 
for all $\delta>0$ we have 
$$\sup_{\substack{\tau\in\mathfrak{T}^T(\F^\RR)\\ \theta\leq\delta}}\ls G,\cdot\rs_\sharp Q_1^N\Big\{\big|f_{\tau}-f_{[\tau+\theta]\mn T}\big|>\ee\Big\}
\leq\sup_{\substack{\tau\in\mathfrak{T}^T(\F^{\MMM_+(\T^d)})\\ \theta\leq\delta}}Q_1^N\Big\{\big|\ls G,\pi_{\tau}-\pi_{[\tau+\theta]\mn T}\rs\big|>\ee\Big\},$$
where of course in the inequality above, $f$ and $\pi$ are the canonical cadlag processes on the Skorohod spaces $D(\RR_+;\RR)$ and $D(\RR_+;\MMM_+(\T^d))$ respectively. With similar reasoning we get that 
$$\sup_{\substack{\tau\in\mathfrak{T}^T(\F^{\MMM_+(\T^d)})\\ \theta\leq\delta}}Q_1^N\Big\{\big|\ls G,\pi_{\tau}-\pi_{[\tau+\theta]\mn T}\rs\big|>\ee\Big\}
\leq\sup_{\substack{\tau\in\mathfrak{T}^T(\F^{\MM_N^d})\\ \theta\leq\delta}}P^N\Big\{\big|\ls G,\pi^N_{\tau}-\pi^N_{[\tau+\theta]\mn T}\rs\big|>\ee\Big\}$$ for all $\delta>0$, where here of course $\pi^N=(\pi^N_t)_{t\geq 0}$ is the empirical process.

Let now $A^{N,G}$ be the martingale defined in~\eqref{MartOfZRProc}. By~\eqref{HeatContEqatthediscretoruses} there exists a constant $C=C(G,d,g)\geq 0$ such that 
$$\Big|\ls G,\pi_t^N\rs-\ls G,\pi_s^N\rs-\fr{2}\int_s^t\ls\D G,\s_r^N\rs dr-(A_t^{N,G}-A_s^{N,G})\Big|\leq\frac{C}{N}\int_s^t\ls 1,\pi_r\rs\df fr,$$ for all $0\leq s\leq t$ and thus 
$$|\ls G,\pi^N_t-\pi^N_s\rs\big|\leq\big|A_t^{N,G}-A_s^{N,G}\big|+\fr{2}\int_s^t\big|\ls\D G,\s_r^N\rs\big|\df r+\frac{C}{N}\int_s^t\ls\pi^N_r,1\rs\df r.$$
But we obviously have that 
$$\big|\ls\D G,\s^N\rs\big|\leq\|\D G\|_\infty\ls 1,\s^N\rs\leq
\|\D G\|_\infty\|\mathfrak{g}'\|_\infty\ls 1,\pi^N\rs$$ and therefore, taking into account the conservation of the total number of particles by the dynamics of the ZRP, we can write that $$|\ls G,\pi^N_t-\pi^N_s\rs\big|\leq\big|A_t^{N,G}-A_s^{N,G}\big|+C_1\cdot(t-s)\ls\pi^N_0,1\rs$$ 
$P^N$-a.s$.$ for some constant $C_1\geq 0$, namely $C_1=C+\fr{2}\|\D G\|_\infty\|\mathfrak{g}'\|_\infty$. It follows that $$\big|\ls G,\pi^N_{[\tau+\theta]\mn T}-\pi^N_\tau\rs\big|\leq\big|A_{[\tau+\theta]\mn T}^{N,G}-A_\tau^{N,G}\big|+C_1\delta\ls\pi^N_0,1\rs$$
for all $\tau\in\mathfrak{T}^T(\F^{\MM_N^d})$ and all $0<\theta\leq\delta$, and therefore 
\begin{align*}
	\sup_{\substack{\tau\in\mathfrak{T}^T\\ \theta\leq\delta}}P^N\big\{\big|\ls G,\pi^N_{[\tau+\theta]\mn T}-\pi^N_\tau\rs\big|>\ee\big\}&\leq\sup_{\substack{\tau\in\mathfrak{T}^T\\ \theta\leq\delta}}
	P^N\Big\{\big|A_{[\tau+\theta]\mn T}^{N,G}-A_\tau^{N,G}\big|>\frac{\ee}{2}\Big\}\\
	&\quad+\mu_0^N\Big\{C_1\delta\ls\pi^N_0,1\rs>\frac{\ee}{2}\Big\}
\end{align*}
for all $\delta>0$. So since the term $\mu_0^N\{C_1\delta\ls\pi^N_0,1\rs>\ee/2\}$ converges to $0$ as $\delta\ra 0$ uniformly over $N$ by~\eqref{fortotalmassofempiricbounded}, in order to prove Aldous' criterion it remains to prove that 
$$\lim_{\delta\ra 0}\limsup_{N\ra+\infty}\sup_{\substack{\tau\in\mathfrak{T}^T(\F^{\MM_N^d})\\ \theta\leq\delta}}P^N\big\{\big|A_{[\tau+\theta]\mn T}^{N,G}-A_\tau^{N,G}\big|>\ee\big\}=0,$$
and by the Chebyshev-Markov inequality it suffices to prove that 
\begin{equation}\label{toendrelatcompactnessofempirtraject}
\lim_{\delta\ra 0}\limsup_{N\ra+\infty}\sup_{\substack{\tau\in\mathfrak{T}^T(\F^{\MM_N^d})\\ \theta\leq\delta}}\EE^N\big(A_{[\tau+\theta]\mn T}^{N,G}-A_\tau^{N,G}\big)^2=0.
\end{equation} By Doob's optional stopping theorem and inequality~\eqref{QVForm}, for any $\theta>0$ and any stopping time $\tau\in\mathfrak{T}^T(\F^{\MM_N^d})$
\begin{align*}
\EE^N\big(A_{[\tau+\theta]\mn T}^{N,G}-A_\tau^{N,G}\big)^2&=
	\EE^N\big(\ls A^{N,G}\rs_{[\tau+\theta]\mn T}-\ls A^{N,G}\rs_\tau\big)\\
	&\leq\frac{2d\|\nabla G\|_\infty^2\|\mathfrak{g}'\|_\infty\theta}{N^{d}}\int\ls\pi^N,1\rs\df\mu_0^N.
\end{align*}
It follows that 
$$\sup_{\substack{\tau\in\mathfrak{T}^T(\F^{\MM_N^d})\\ \theta\leq\delta}}\EE^N\big(A_{[\tau+\theta]\mn T}^{N,G}-A_\tau^{N,G}\big)^2\leq
\frac{2d\|\nabla G\|_u^2\|\mathfrak{g}'\|_\infty\delta}{N^{d}}\int\ls\pi^N,1\rs d\mu_0^N.$$
Therefore by the $O(N^d)$-entropy assumption and Lemma~\ref{BoundedTotParticlInMacrLim} the claim follows.

We proceed next to show that the sequence $Q_2^N:=W^N_\sharp P^N\in \PP L_{w^*}^\infty(0,T;\X^1(\T^d)^*)$, $N\in\NN$ is relatively compact. For this we will use the Prokhorov-Le Can theorem~\ref{ProkhorovLeCam} according to which it suffice to show that the family $\{W^N_\sharp P^N\}$ is uniformly tight. Since by the Banach-Alaoglu theorem the closed balls of $L_{w^*}^\infty(0,T;\X^1(\T^d)^*)$ are compact, it suffices to show that 
\begin{equation}\label{CurrentTightness}
\lim_{M\ra+\infty}\sup_{N\in\NN}Q_2^N\big\{\|W\|_{L_{w^*}^\infty(0,T;\X^1(\T^d)^*)}>M\big\}=0.
\end{equation} By the vector-valued mean value theorem for any $G\in\X^1(\T^d)$
$$|G(u)-G(\y)|_2\leq\|G\|_{\X^1}\cdot d_{\T^d}(u,\y),\quad u,\y\in\T^d.$$
and thus for the empirical current process $W^N\colon D(0,T;\MM_N^d)\to L_{w^*}^\infty(0,T;\X^1(\T^d)^*)$ 
\begin{align*}
\|W^N\|_{L_{w^*}^\infty}&=\sup_{\|G\|_{L^1(0,T;\X^1(\T^d))}\leq 1}|W^N(G)|\\
&=\sup_{\|G\|\leq 1}\Big|\int_0^T\fr{N^{d-1}}\sum_{\substack{x\in\T_N^d\\j=1,\ldots,d}}G_t^j\Big(\frac{x}{N}\Big)\{\mathfrak{g}(\eta_t(x))-\mathfrak{g}(\eta_t(x+e_j))\}\df t\Big|\\
&=\sup_{\|G\|\leq 1}\Big|\int_0^T\fr{N^{d-1}}\sum_{\substack{x\in\T_N^d\\j=1,\ldots,d}}\Big\{G_t^j\Big(\frac{x}{N}\Big)-G_t^j\Big(\frac{x-e_j}{N}\Big)\Big\}\mathfrak{g}(\eta_t(x))\df t\Big|\\
&\leq\sup_{\|G\|\leq 1}\int_0^T\fr{N^{d-1}}\sum_{x\in\T_N^d}\mathfrak{g}(\eta_t(x))\sum_{j=1}^d\Big|G_t^j\Big(\frac{x}{N}\Big)-G_t^j\Big(\frac{x-e_j}{N}\Big)\Big|\df t\\
&\leq d\sup_{\|G\|\leq 1}\int_0^T\fr{N^{d-1}}\sum_{x\in\T_N^d}\mathfrak{g}(\eta_t(x))\Big|G_t\Big(\frac{x}{N}\Big)-G_t\Big(\frac{x-e_j}{N}\Big)\Big|_2\df t\\
&\stackrel{P^N\mbox{-a.s.}}\leq d\|\mathfrak{g}'\|_\infty\ls 1,\pi_0^N\rs\sup_{\|G\|_{L^1(0,T;\X^1(\T^d))}\leq 1}\int_0^T\|G_t\|_{\X^1}\df t\\
&\leq d\|\mathfrak{g}'\|_\infty\ls 1,\pi_0^N\rs.
\end{align*}
Therefore 
\begin{align*}
Q_2^N\big\{\|W\|_{L_{w^*}^\infty(0,T;\X^1(\T^d)^*)}>M\big\}&=P^N\big\{\|W^N\|_{L_{w^*}^\infty(0,T;\X^1(\T^d)^*)}>M\big\}\\
&\leq P^N\Big\{\ls 1,\pi_0^N\rs>\frac{M}{d\|\mathfrak{g}'\|_\infty}\Big\}\\
&=\mu_0^N\Big\{\ls 1,\pi^N\rs>\frac{M}{d\|\mathfrak{g}'\|_\infty}\Big\}
\end{align*}
and thus~\eqref{CurrentTightness} follows by Lemma~\ref{BoundedTotParticlInMacrLim} since $\{\mu_0^N\}$ has total mass $\mu_0(\T^d)$.

Finally the relative compactness of the third marginal $Q_3^N:=\s^N_\sharp P^N\in\PP L_{w^*}^\infty(0,T;\MMM_+(\T^d))$ in $L_{w^*}^\infty(0,T;\MMM(\T^d))$ follows by Proposition~\ref{PsiEmpDistComp}. By definition the sequence $\{Q_3^N\}_{N\in\NN}$ is supported by the set $L_{w^*}^\infty(0,T;\MMM_+(\T^d))$, which according to Proposition~\ref{PosCon} is a closed subspace of $L_{w^*}(0,T;\MMM(\T^d))$. Therefore by the portmanteau theorem it follows that the sequence $\{Q_3^N\}$ is also relatively compact in $L_{w^*}^\infty(0,T;\MMM_+(\T^d))$, i.e.~any subsequential limit point $Q_3$ of $\{Q_3^N\}$ is supported by the set $L_{w^*}^\infty(0,T;\MMM_+(\T^d))$.


We will prove now properties (a) to (e). We start by proving (a), i.e.~that any limit point $Q$ of the sequence $\{Q^N\}$ is concentrated on trajectories $(\pi,W,\s)\in\W$ such that the continuity equation~\eqref{CE} holds. By the estimate~\eqref{HeatContEqatthediscretoruses} 
\begin{equation*}
\sup_{0\leq t\leq T}\big|V_t^{1,N,G}-A_t^{N,G}\big|\mx\big|V_t^{2,N,G}-A_t^{N,G}\big|\leq\frac{C}{N}\int_0^T\ls 1,\pi^N_t\rs\df t\stackrel{P^N\mbox{-a.s.}}=\frac{CT}{N}\ls 1,\pi^N_0\rs
\end{equation*}
and therefore by the asymptotic negligibility~\eqref{MartAsymptNegl} and Proposition~\ref{BoundedTotParticlInMacrLim} it follows that for all $G\in C_c^3(\RR_+\x\T^d)$ and all $\delta>0$
\begin{equation}\label{MicroCE}
\lim_{N\ra+\infty}P^N\Big\{\sup_{0\leq t\leq T}|V_t^{1,N,G}|\mx|V_t^{2,N,G}|\geq\delta\Big\}=0.
\end{equation}
Let us define for any $G\in C_c^3((0,T)\x\T^d)$ the maps $f^{i,G}\colon\W\to\RR$, $i=1,2,3$ by the formulas
$$
f^{1,G}(\w)=\int_0^T\ls\pd_tG_t,\pi_t\rs\df t,\quad
f^{2,G}(\w)=\int_0^T\ls\nabla G_t,W_t\rs\df t,\quad
f^{3,G}(\w)=\int_0^T\ls\D G_t,\s_t\rs\df t
$$ where $\w=(\pi,W,\s)$. Then with this notation in order to prove that the continuity equation holds in the hydrodynamic limit it suffices to show that 
\begin{equation}\label{JRCCE}Q\Big(\bigcap_{G\in C_c^3((0,T)\x\T^d)}\big\{f^{1,G}+f^{j,G}=0\big\}\Big)=1,\quad j=2,3.\end{equation}

By the separability of $C_c(M)$ when $M$ is a locally compact topological space it follows that $C_c^3((0,T)\x\T^d)$ is separable in the $C^2$-uniform norm 
$$\|G\|_{C^2}:=\|G\|_\infty+\|\nabla_{(t,u)} G\|_\infty+\|D_{(t,u)}^2G\|_\infty,$$
where the gradient $\nabla_{(t,u)}$ and the second derivative $D_{(t,u)}^2$ appearing in the definition of the $C^2$-uniform norm are with respect to both time and space. Thus there exists a countable family $\mathcal{G}:=\{G_k\}_{k=1}^\infty\subs C_c^3((0,T)\x\T^d)$ that is dense in $C_c^3((0,T)\x\T^d)$ in the uniform $C^2$-norm and then 
\begin{equation}\label{FromAllToOne}
\bigcap_{G\in C_c^3((0,T)\x\T^d)}\big\{f^{1,G}+f^{j,G}=0\big\}=\bigcap_{G\in\mathcal{G}}\big\{f^{1,G}+f^{j,G}=0\big\},\quad j=2,3.
\end{equation}
Indeed, if $\{G_n\}_{n=1}^\infty$ is a sequence in $\mathcal{G}$ that converges to $G\in C_c^3((0,T)\x\T^d)$ then 
$\pd G_n$, $\nabla G_n$ and $\D G_n$ converge uniformly on $(0,T)\x\T^d$ to $\pd G$, $\nabla G$ and $\D G$ respectively and therefore $\lim_{n\ra+\infty}f^{j,G_n}=f^{j,G}$ pointwise on $\W$ for $j=1,2,3$ which proves equality~\eqref{FromAllToOne}. Thus~\eqref{JRCCE} is reduced to showing that 
\begin{equation}\label{JRCCEOne}Q\big(\big\{f^{1,G}+f^{j,G}=0\big\}\big)=1,\quad\forall G\in C_c^3((0,T)\x\T^d)\mbox{ for }j=2,3.\end{equation}
By well known results on induced mappings on Skorohod spaces the map $f^{1,G}$ is continuous and for $G\in C_c^3((0,T)\x\T^d)$ the maps $\nabla G$ and $\D G$ induces elements in the spaces $L^1(0,T;\X^1(\T^d))$ and $L^1(0,T;C(\T^d))$ respectively so that the maps $f^{2,G}$ and $f^{3,G}$ are given by 
$$f^{2,G}(\pi,W,\s)=\lls \nabla G,W\rrs,\quad f^{3,G}(\pi,W,\s)=\lls \D G,\s\rrs$$
and are thus continuous. Therefore for any $\delta>0$ the set 
$$A_\delta:=\big\{|f^{1,G}+f^{2,G}|\mx|f^{1,G}+f^{3,G}|>\delta\big\}$$
is open in $\W$. Furthermore, for $G\in C_c^3((0,T)\x\T^d)$ the processes $V^{j,N,G}$, $j=1,2$, defined in~\eqref{Vjump} and~\eqref{Vcurr} respectively satisfy 
$$V^{j,N,G}=-(f^{1,G}+f^{j+1,G}).$$ Therefore if $\{Q^{k_N}\}$ is a subsequence of $\{Q^N\}_{N=1}^\infty$ converging to $Q$ then by the portmanteau theorem and~\eqref{MicroCE}
\begin{align*}
Q(A_\delta)&\leq\liminf_{N\ra+\infty}Q^{k_N}(A_\delta)=\liminf_{N\ra+\infty}P^N\big\{|V^{1,N,G}|\mx|V^{2,N,G}|>\delta\big\}=0.
\end{align*}
Now since this holds for any $\delta>0$ we obtain~\eqref{JRCCEOne} which proves that the equation $\pd_t\pi=-\dv W=\Delta\s$ holds in the hydrodynamic limit.

We prove next the second equation of~\eqref{CE}, i.e.~that $W=-\nabla\s$. For this we define the gradient operator $\nabla\colon\MMM(\T^d)\to\X^1(\T^d)^*$ by 
$$\ls F,\nabla\mu\rs=-\ls\dv F,\mu\rs,\quad F\in\X^1(\T^d).$$
This is $w^*$-continuous since $\nabla=-\dv^*$ where $\dv\colon\X^1(\T^d)\to C(\T^d)$ is the divergence operator and thus it induces a $w^*$-continuous gradient operator 
$\nabla\colon L^1(0,T;\X^1(\T^d))\to L^1(0,T;C(\T^d))$. We also define the family of discrete gradient operators 
$$\nabla_-^N\colon\MMM(\T^d)\to\MMM_0(\T^d;\RR^d)\leq\X^1(\T^d)^*$$
by the formula
$$\nabla_-^N\mu=N\sum_{j=1}^d(\tau_{-\frac{e_j}{N}\sharp}\mu-\mu)e_j.$$ It is easy to verify that $\nabla_-^N=-(\dv_+^N)^*$ where $\dv_+^N\colon\X^1(\T^d)\to C(\T^d)$ is the discrete divergence operator 
$$\dv_+^NF(u)=N\sum_{j=1}^d\Big(F_j\Big(u+\frac{e_j}{N}\Big)-F_j(u)\Big).$$
Consequently the gradient operators induce the $w^*$-continuous gradient operators \begin{equation}\label{MeasureGradient}\nabla_-^N:=-(\dv_\mp^N)^*\colon L_{w^*}^\infty(0,T;\MMM(\T^d))\to L_{w^*}^\infty(0,T;\X^1(\T^d)^*)\end{equation} on the $L_{w^*}^\infty$-spaces as the adjoints of the induced divergence operators $-\dv_+^N\colon L^1(0,T;\X^1(\T^d))\to L^1(0,T;C(\T^d))$ on the $L^1$-Bochner spaces. Then with this notation the empirical current and the empirical jump rate are related by the equality $W^N=-\nabla_-^N\circ\s^N$.
\begin{lemma} Suppose that the sequence of initial distributions total mass bounded above by $\bar{\mathfrak{m}}>0$ in probability. Then for any limit point $Q\in\Lim_{N\ra+\infty}(W^N,\s^N)_\sharp P^N$ 
	$$W=-\bar{\nabla}\s\quad Q\mbox{-a.s.~for all }(W,\s)\in L_{w^*}^\infty(0,T;\X^1(\T^d)^*)\x L^\infty(0,T;\MMM(\T^d))$$
	where $\bar{\nabla}\colon L_{w^*}^\infty(0,T;\MMM(\T^d))\to L_{w^*}^\infty(0,T;\X^1(\T^d)^*)$ is the induced gradient operator on the level of path-measures.
\end{lemma}\textbf{Proof} We have to prove that $Q\{W=-\nabla\s\}=1$. Since $L^1(0,T;\X^1(\T^d))$ is separable it suffices to show that $Q\{|\lls G,W+\nabla\s\rrs|>\ee\}=0$ for all $G\in L^1(0,T;\X^1(\T^d))$ and all $\ee>0$. Since the gradient operator is $w^*$-continuous, the set $\{|\lls G,W+\nabla\s\rrs|>\ee\}$ is open and therefore by the portmanteau theorem 
$$Q\{|\lls G,W+\nabla\s\rrs|>\ee\}\leq\limsup_{N\ra+\infty}P^N\{|\lls G,W^N+\nabla\s^N\rrs|>\ee\}.$$
By the equality $W^N=-\bar{\nabla}_-^N\s^N$,
\begin{align*}|\lls G,W^N+\bar{\nabla}\s^N\rrs|&=|\lls G,\bar{\nabla}\s^N-\bar{\nabla}_-^N\s^N\rrs|=|\lls \dv_+^NG-\dv G,\s^N\rrs|\\
&\leq\lls \|\dv_+^NG-\dv G\|_\infty,\s^N\rrs\leq\|\mathfrak{g}'\|_\infty\lls\|\dv_+^NG-\dv G\|_\infty,\pi^N\rrs\\
&\leq T\|\mathfrak{g}'\|_\infty\|\dv_+^NG-\dv G\|_{\infty,1}\ls 1,\pi_0^N\rs
\end{align*}

Now, for any $G\in\X^1(\T^d)$ by the fundamental theorem of calculus $\dv_+^NG=\int_0^1\dv G(\cdot+\frac{se_j}{N})\df s$ and, since $\T^d$ is compact, the map $\dv G$ is uniformly continuous and therefore $$\lim_{N\ra+\infty}\|\dv_+^NG-\dv G\|_\infty=0,\quad\forall G\in\X^1(\T^d).$$ In other words the sequence of operators $\dv_+^N\colon\X^1(\T^d)\to C(\T^d)$ converges strongly to $\dv$. This implies by Proposition~\ref{InducedOpPropL1StrongConv} that the induced operators on the corresponding $L^1$-spaces pointwise converge strongly to $\bbar{\dv}\colon L^1(0,T;\X^1(\T^d))\to L^1(0,T;C(\T^d))$ and thus the sequence $a_N(G):=\|\dv_+^NG-\dv G\|_{\infty,1}$ converges to $0$ as $N\ra+\infty$. Consequently
$$Q\{|\lls G,W+\nabla\s\rrs|>\ee\}\leq\limsup_{N\ra+\infty}\mu_0^N\Big\{\ls 1,\pi^N\rs>\frac{\ee}{T\|\mathfrak{g}'\|_\infty a_N(G)}\Big\}=0$$
by Lemma~\ref{BoundedTotParticlInMacrLim} since $\lim_{N\ra+\infty}a_N(G)=0$.$\hfill\Box$\\

We prove next (b), i.e.~that any subsequential limit point $Q$ of the sequence $\{Q^N\}$ is concentrated on trajectories $\pi\in C(0,T;\MMM_+(\T^d))$ such that $\pi_0=\mu_0$. So let $\{Q^{k_N}\}$ be a subsequence of $\{Q^N\}$ converging to $Q$. We show first that $\pi_0=\mu_0$. The evaluation mapping $e_t\colon D(0,T;\MMM_+(\T^d))\to\MMM_+(\T^d)$ given by $e_t(\pi)=\pi_t$ is continuous at each $\pi\in D(0,T;\MMM_+(\T^d))$ that is continuous at $t\in[0,T]$. In particular the evaluation $e_0:D(0,T;\MMM_+(\T^d))\lra\MMM_+(\T^d)$ is continuous and therefore for all $G\in C(\T^d)$ the composite mapping $\ls G,\cdot\rs\circ e_0:D(0,T;\MMM_+(\T^d))\lra\RR$ is continuous. Therefore, for all $G\in C(\T^d)$ and all $\ee>0$ we have by the portmanteau theorem that 
\begin{align*}
	Q\big\{|\ls G,\pi_0\rs-\ls G,\mu_0\rs|>\ee\big\}&\leq
	\liminf_{N\ra\infty}Q^N\big\{|\ls G,\pi_0\rs-\ls G,\mu_0\rs|>\ee\big\}\\
	&=\liminf_{N\ra\infty}\mu_0^N\big\{|\ls G,\pi^N\rs-\ls G,\mu_0\rs|>\ee\big\}=0,
\end{align*} since the sequence $\{\mu_0^N\}$ is associated to the macroscopic profile $\mu_0\in\MMM_+(\T^d)$. Since and this holds for all $\ee>0$ it follows that $Q\big\{\ls G,\pi_0\rs=\ls G,\mu_0\rs\big\}=1$ for all $G\in C(\T^d)$ and thus is $C(\T^d)$ separable we can choose a countable subset $D\subs C(\T^d)$ dense in $C(\T^d)$ in the uniform norm and then $$Q\{\pi_0=\mu_0\}=Q\bigg(\bigcap_{G\in D}\{|\ls G,\pi_0\rs-\ls G,\mu_0\rs|=0\}\bigg)=1.$$

The set $C(0,T;\MMM_+(\T^d))$ is a closed subspace of $D(0,T;\MMM_+(\T^d))$ in the Skorohod metric and thus measurable. The fact that $Q(C(0,T;\MMM_+(\T^d))=1$ follows by the existence of continuous representatives for solutions of the continuity equation.

\begin{prop}\label{ContEqBasicDensRegularity}
	Let $(\pi,W)\in L_{w^*}^\infty\big(0,T;\MMM_+(\T^d)\big)\x L^\infty_{w^*}\big(0,T;\X^1(\T^d)^*\big)$ be a density-current pair satisfying the continuity equation. Then there exists a weakly continuous curve $\wt{\pi}$ in the class of $\pi$ in $L_{w^*}^\infty\big(0,T;\MMM_+(\T^d)\big)$, and for this continuous representative for all $0\leq s\leq t\leq T$
	\begin{equation}\label{CECR}
	\int_{\T^d}G_t\df\wt{\pi}_t-\int_{\T^d}G_s\df\wt{\pi}_s=\int_s^t\Big(\int_{\T^d}\pd_rG_r\df\wt{\pi}_r+\ls\nabla G_r,W_r\rs\Big)\df r,\quad\forall\;G\in C^{1,2}([0,T]\x\T^d).\end{equation}
\end{prop}\textbf{Proof} We fix a function $\zeta\in C^2(\T^d)$ and let $\phi_\zeta:(0,T)\lra\RR$ denote the function defined a.s.~by 
$$\phi_\zeta(t)=\int_{\T^d}\zeta\df\pi_t.$$ Then $\phi_\zeta\in L^\infty(0,T)$ with norm $\|\phi_\zeta\|_\infty\leq\|\zeta\|_\infty\|\pi\|_{TV;\infty}<+\infty$ since $\pi\in L_{w^*}^\infty\big(0,T;\MMM_+(\T^d)\big)$.

Let now $G\in C_c^{1,2}((0,T)\x\T^d)$ be any function of the form $G(t,u)=f(t)\zeta(u)$ for some function $f\in C_c^1(0,T)$ and some function $\zeta\in C^2(\T^d)$. Then since the pair $(\pi,W)\in L^\infty(0,T;\MMM_+(\T^d))\x L^\infty_{w^*}(0,T;\X^1(\T^d)^*)$ satisfies the continuity equation, 
$$\int_0^Tf'(t)\varphi_\zeta(t)\df t=\int_0^Tf'(t)\int_M\zeta\df\wt{\pi}_t\df t=-\int_0^Tf(t)\ls\nabla\zeta,W_t\rs\df t.$$
Since the equality above holds for all $f\in C_c^1(0,T)$, the measurable function $\psi_\zeta\colon[0,T]\to\RR$ defined a.s.~by $\psi_\zeta(t)=\ls\nabla\zeta,W_t\rs$
is the weak derivative of the function $\phi_\zeta$. But since $W\in L_{w^*}^\infty(0,T;\X^1(\T^d)^*)$ the function $\psi_\zeta$ is in $L^\infty(0,T)$ with $\|\psi_\zeta\|_\infty\leq\|\nabla\zeta\|_{\X^1(\T^d)}\|W\|_{L_{w^*}^\infty(0,T;\X^1(\T^d)^*)}$. Therefore $\phi_\zeta\in W^{1,\infty}(0,T)$ with distributional derivative $\psi_\zeta$. Consequently, the equivalence class $\phi_\zeta$ contains a Lipschitz representative $\bar{\phi}_\zeta$ with Lipschitz constant
$$\|\bar{\phi}_\zeta\|_{{\Lip}}\leq\|\psi_\zeta\|_{L^\infty(0,T)}\leq\|\nabla\zeta\|_{\X^1(\T^d)}\|W\|_{L_{w^*}^\infty(0,T;\X^1(\T^d)^*)}.$$

Let now $\mathcal{Z}$ be a countable subset of $C^\infty(\T^d)$ that is dense in $C^2(\T^d)$ in the usual $C^2$-norm $\|\cdot\|_{C^2}$ given by $\|\zeta\|_{C^2}=\|\zeta\|_\infty+\|\nabla\zeta\|_\infty+\|D^2\zeta\|_\infty$ for $\zeta\in C^2(\T^d)$. Then $\mathcal{Z}$ is also dense with the uniform norm $\|\cdot\|_\infty$ in $C(\T^d)$ and we set 
$$I_\mathcal{Z}:=\bigcap_{\zeta\in\mathcal{Z}}\big\{t\in[0,T]\bigm|\phi_\zeta(t)=\bar{\phi}_\zeta(t)\big\}.$$
Then $I_\mathcal{Z}$ is of full Lebesgue measure in $[0,T]$. We continue to denote by $\pi\colon I_\mathcal{Z}\to\MMM_+(\T^d)$ the restriction of $\pi\in L_{w^*}^\infty(0,T;\MMM_+(\T^d))$ on $I_\mathcal{Z}$. Then since $\MMM_+(\T^d)\leq C(\T^d)^*$ and $C(\T^d)^*$ is naturally injected in $C^2(\T^d)^*$ via the restriction of domains, i.e.~via the adjoint $i^*$ of the natural inclusion $i\colon C^2(\T^d)\hookrightarrow C(\T^d)$, we can regard $\pi$ as the function $\hat{\pi}:= i^*\circ\hat{\pi}\colon I_\mathcal{Z}\to C^2(\T^d)^*$. As such the function $\hat{\pi}$ is Lipschitz, with Lipschitz constant $\leq\|W\|_{L_{w^*}^\infty(0,T;\X^1(\T^d)^*)}$. Indeed, for all $s,t\in I_\mathcal{Z}$ and all $\zeta\in\mathcal{Z}$ 
\begin{align*}|\hat{\pi}_t(\zeta)-\hat{\pi}_s(\zeta)|&=|\bar{\phi}_\zeta(t)-\bar{\phi}_\zeta(s)|
	\leq\|\hat{\varphi}_\zeta\|_{{\Lip}}|t-s|\leq\|W\|_{L_{w^*}^\infty(0,T;\X^1(\T^d)^*)}\|\nabla\zeta\|_{\Lip}|t-s|\\
	&\leq\|W\|_{L_{w^*}^\infty(0,T;\X^1(\T^d)^*)}\|\zeta\|_{C^2}|t-s|,
\end{align*}
which since $\mathcal{Z}$ is dense in $C^2(\T^d)$ in the $C^2$-norm $\|\cdot\|_{C^2}$ shows that 
$$\|\hat{\pi}_t-\hat{\pi}_s\|_{C^2(\T^d)^*}=
\sup_{\zeta\in\mathcal{Z}}\frac{|\hat{\pi}_t(\zeta)-\hat{\pi}_s(\zeta)|}{\|\zeta\|_{C^2}}\leq\|W\|_{L_{w^*}^\infty(0,T;\X^1(\T^d)^*)}|t-s|.$$
Therefore $\hat{\pi}:I_\mathcal{Z}\to C^2(\T^d)^*$ has a Lipschitz extension $\wt{\pi}\colon[0,T]\to C^2(\T^d)^*$ with the same Lipschitz constant.

\indent Now, since $\pi$ belongs in $L_{w^*}^\infty(0,T;\MMM_+(\T^d))$ by hypothesis, we can assume that $I_\mathcal{Z}$ has been chosen so that $\|\pi_t\|_{TV}=\pi_t(\T^d)\leq\|\pi\|_{\infty;TV}<+\infty$ for all $t\in I_\mathcal{Z}$. Thus since $\T^d$ is compact the family $\{\pi_t\}_{t\in I_\mathcal{Z}}$ is relatively compact in the weak topology of $\MMM_+(\T^d)$. It follows that the Lipschitz extension $\wt{\pi}\colon[0,T]\lra C^2(\T^d)^*$ takes values in $\MMM_+(\T^d)\equiv i^*(\MMM_+(\T^d))\leq C^2(\T^d)^*$ and is weakly continuous. Indeed, if $t\in[0,T]\sm I_\mathcal{Z}$ and $\{t_n\}_{n=1}^\infty\subs I_\mathcal{Z}$ is a sequence converging to $t$ then 
\begin{equation}\label{C2Karuga}
\lim_{n\ra+\infty}\ls f,\pi_{t_n}\rs=\ls f,\wt{\pi}_t\rs,\quad\forall f\in C^2(\T^d).
\end{equation}
But since $\{\pi_{t_n}\}$ is contained in the compact set of measures with total variation norm $\leq\|\pi\|_{TV;\infty}$ there exists a subsequence $\{t_{k_n}\}$ of $\{t_n\}$ and a measure $\mu_t\in\MMM_+(\T^d)$ such that
\begin{equation}\label{CKaruga}
\lim_{n\ra+\infty}\ls f,\pi_{t_{k_n}}\rs=\ls f,\mu_t\rs,\quad\forall\;f\in C(\T^d).
\end{equation}
By~\eqref{C2Karuga} and~\eqref{CKaruga} it follows that $\ls f,\wt{\pi}_t\rs=\ls f,\mu_t\rs$ for all $f\in C^2(\T^d)$ and therefore $\wt{\pi}_t=i^*(\mu_t)\in i^*(\MMM_+(\T^d))\cong\MMM_+(\T^d)$. Since any measure $\mu_t$ satisfying~\eqref{CKaruga} for any sequence $\{t_n\}$ converging to $t$ must necessarily coincide with $\wt{\pi}_t$ on $C^2$-functions it is unique and thus we can identify $\wt{\pi}_t$ with $\mu_t\in\MMM_+(\T^d)$. To see that the curve $\wt{\pi}\colon[0,T]\to\MMM_+(\T^d)$ thus defined is weakly continuous let $\{t_n\}\subs[0,T]$ be any sequence converging to $t$. We will show that any subsequence $\{t_{k_n}\}$ of $\{t_n\}$ has a further subsequence $\{t_{m_{k_n}}\}$ such that 
$$\lim_{n\ra+\infty}\ls f,\wt{\pi}_{t_{m_{k_n}}}\rs=\ls f,\wt{\pi}_t\rs,\quad\forall\; f\in C(\T^d).$$  
So let $\{t_{k_n}\}$ be a subsequence of $\{t_n\}$. Since  $\wt{\pi}_{t_{k_n}}$ is relatively compact in $\MMM_+(\T^d)$ there exists $\mu_t\in\MMM_+(\T^d)$ and a subsequence $\{t_{k_{m_n}}\}$ of $\{t_{k_n}\}$ such that $\lim_{n\ra+\infty}\ls f,\wt{\pi}_{t_{k_{m_n}}}\rs=\ls f,\mu_t\rs$ for all $f\in C(\T^d)$ and since $i^*\circ\wt{\pi}$ is Lipschitz $\lim_{n\ra+\infty}\ls f,\wt{\pi}_{t_{m_{k_n}}}\rs=\ls f,\wt{\pi}_t\rs$ for all $f\in C^2(\T^d)$. Therefore $\mu_t=\wt{\pi}_t$ and the curve $\wt{\pi}$ is continuous.

\indent We prove finally~\eqref{CECR}. So let $G\in C^{1,2}([0,T]\x\T^d)$, let $0\leq s<t\leq T$ and let $W$ be any curve in $L^\infty_{w^*}(0,T;\X^1(\T^d)^*)$. Let $f_\ee\in C_c^\infty\big((s,t);[0,1]\big)$, $\ee>0$, be such that $f_\ee\lra\1_{(s,t)}$ pointwise in $[0,T]$ and such that 
\begin{equation}\label{Der1Approx}\lim_{\ee\ra 0}\int_0^Tf_\ee '(r)h(r)\df r=h(s)-h(t),\quad\forall\;h\in C([0,T]).
\end{equation}
Then since the pair $(\pi,W)$ satisfies the continuity equation we have for all $\ee>0$ that 
\begin{align}\label{ToGetCECR}
	0&=\int_0^T\Big(\int_{\T^d}\pd_r[f_\ee(r)G_r]\df\wt{\pi}_r+\ls\nabla[f_\ee(r) G_r],W_r\rs\Big)\df r\nonumber\\
	&=\int_0^T\Big(\int_{\T^d}[f_\ee '(r)G_r+f_\ee(r)\pd_rG_r]\df\wt{\pi}_r+\ls f_\ee(r)\nabla G_r,W_r\rs\Big)\df r\nonumber\\
	&=\int_0^Tf_\ee '(r)\int_{\T^d}G_r\df\wt{\pi}_r
	+\int_0^Tf_\ee(r)\Big(\int_{\T^d}\pd_rG_r\df\wt{\pi}_r+\ls\nabla G_r,W_r\rs\Big)\df r.
\end{align}
Now, since $G\in C^{1,2}([0,T]\x\RR^d)$ the curve $[0,T]\ni t\mapsto G_t\in C(\T^d)$ is continuous with respect to the uniform norm in $C(\T^d)$ and therefore due to the weak continuity of $\wt{\pi}$, the function $[0,T]\ni t\mapsto\int G_t\df\wt{\pi}_t=\int G_t\df\wt{\pi}_t$ is continuous. Therefore taking the limit $\ee\ra 0$ and using (\ref{Der1Approx}) in the first term of the sum in~\eqref{ToGetCECR} and the bounded convergence theorem in the second term, we obtain~\eqref{CECR} as required.$\hfill\Box$\\

 Now, since $\pi\in D(0,T;\MMM_+(\T^d))$ and $\pi$ is almost everywhere equal to a weakly continuous path by Proposition~\ref{ContEqBasicDensRegularity}, it follows that $\pi\in C(0,T;\MMM_+(\T^d))$. Indeed, let $t\in[0,T]$ and let $E\subs[0,T]$ be a set of full measure in $[0,T]$ such that $\pi_s=\wt{\pi}_s$ for all $s\in E$. Since $I$ is of full measure in $[0,T]$ there exist sequences $\{s_n\}\subs I\cap[0,t]$ and $\{r_n\}\subs I\cap[t,T]$ such that $s_n\uparrow t$ and $r_n\downarrow t$. Then since $\pi$ is in $D(0,T;\MMM_+(\T^d))$ and $\wt{\pi}$ is continuous on one hand 
$$\pi_t=\lim_{n\ra+\infty}\pi_{s_n}=\lim_{n\ra+\infty}\wt{\pi}_{s_n}=\wt{\pi}_t$$
and on the other hand
$$\pi_{t-}:=\lim_{n\ra+\infty}\pi_{r_n}=\lim_{n\ra+\infty}\wt{\pi}_{r_n}=\wt{\pi}_t.$$
It follows that $\pi_{t-}=\pi_t$ and thus $\pi_{t-}=\pi_t=\wt{\pi}_t$ so that $\pi\in C(0,T;\MMM_+(\T^d))$.

Statement (c) follows by applying the following corollary of the one-block estimate to the cylinder map $\Psi=\mathfrak{g}(\eta(0))$.
\begin{cor} Assume that the sequence $\{\mu_0^N\}_{N\in\NN}$ of initial distributions is associated to a macroscopic profile $\mu_0\in\MMM_+(\T^d)$ and satisfies the $O(N^d)$-entropy assumption and let $\Psi\colon\MM_\infty^d\to\RR$ be a sublinear cylinder map. Then any limit point of the sequence of laws $\{\s^{N,\Psi}_\sharp P^N\}\subs\PP L_{w^*}^\infty(0,T;\MMM(\T^d))$ is concentrated on trajectories $\s\in L_{w^*}^\infty(0,T;\MMM(\T^d))$ such that $\s_t\ll\LL_{\T^d}$ and 
	$$\Big\|\frac{\df\s_t}{\df\LL_{\T^d}}\Big\|_\infty\leq\psi_c:=\sup_{\rho\in[0,\rho_c]\cap\RR}|\wt{\Psi}|(\rho)$$
	for almost all $t\in[0,T]$.
\end{cor}\textbf{Proof} We consider the joint laws $$\bs{Q}_\Psi^{N,\ell}:=(\bs\pi^{N,\ell},\s^{N,\Psi})_\sharp P^N\in\PP(L_{w^*}^\infty(0,T;\bbar{\MMM}_{1,+}(\T^d\x\RR_+)\x L_{w^*}^\infty(0,T;\MMM(\T^d)),\quad(N,\ell)\in\NN\x\ZZ_+.$$ By the statement of the one-block estimate in terms of generalized Young measures we know that the family $\{\bs{Q}_\Psi^{N,\ell}\}_{(N,\ell)}$ is relatively compact and any limit point $\bs{Q}_\Psi$ of $\{\bs{Q}_\Psi^{N,\ell}\}$ is concentrated on the graph of the $\bbar{\Psi}$-projection $\bar{B}_{\bbar{\Psi}}\colon L_{w^*}^\infty(0,T;\bbar{\MMM}_{1,+}(\T^d\x\RR_+))\to L_{w^*}^\infty(0,T;\MMM(\T^d))$, i.e. $$(\bar{B}_{\bbar{\Psi}}\x\mathbbm{id}_{L_{w^*}^\infty(0,T;\MMM(\T^d))})_\sharp\bs{Q}_{\Psi}(\Delta_{L_{w^*}^\infty(0,T;\MMM(\T^d))})=1$$ 
where $\Delta_{L_{w^*}^\infty(0,T;\MMM(\T^d))}:=\{(\mu,\mu)\bigm|\mu\in L_{w^*}^\infty(0,T;\MMM(\T^d))\}$
is the diagonal in $L_{w^*}^\infty(0,T;\MMM(\T^d))^2$
and $\mathbbm{id}_{L_{w^*}^\infty(0,T;\MMM(\T^d))}$ is the identity mapping on $L_{w^*}^\infty(0,T;\MMM(\T^d))$. Let now $Q^2_\Psi\in\Lim_{N\ra+\infty}\s^{N,\Psi}_\sharp P^N$. Since $L_{w^*}^\infty(0,T;\MMM_{ac}(\T^d))$ is a $w^*$-measurable subspace of $L_{w^*}^\infty(0,T;\MMM(\T^d))$ by Proposition~\ref{ACMEAS}, for the first claim it suffices to show that $Q_\Psi^2(L_{w^*}^\infty(0,T;\MMM_{ac}(\T^d))=1$. Also, the set $$L_{w^*}^\infty\big(0,T;B_{L^\infty(\nu)}(0,\psi_c)\big)=\Big\{\mu\in L_{w^*}^\infty(0,T;\MMM_{ac}(\T^d)\Bigm|\Big\|\Big\|\frac{\df\mu_\cdot}{\df\LL_{\T^d}}\Big\|_{L^\infty(\T^d)}\Big\|_{L^\infty(0,T)}\leq\psi_c\Big\}$$ is a $w^*$-measurable subspace of $L_{w^*}^\infty(0,T;\MMM(\T^d))$ by Lemma~\ref{InftyLipInfty} and thus for the second claim it suffices to show that if $\psi_c<+\infty$ then 
$Q_\Psi^2(L_{w^*}^\infty\big(0,T;B_{L^\infty(\T^d)}(0,\psi_c)\big)=1$.

Since $Q_\Psi^2\in\Lim_{N\ra+\infty}\s^{N,\Psi}_\sharp P^N$ there exists a law
$$\bs{Q}_\Psi\in\Lim_{\ell,N\ra+\infty}\bs{Q}_\Psi^{N,\ell}\subs\PP\big(L_{w^*}^\infty(0,T;\MMM(\T^d))\x L_{w^*}^\infty(0,T;\bbar{\MMM}_1(\T^d\x\RR_+))\big)$$
with second marginal $p^2_\sharp\bs{Q}_\Psi=Q_\Psi^2$ on $L_{w^*}^\infty(0,T;\MMM(\T^d))$ where here $p^2\colon L_{w^*}^\infty(0,T;\bbar{\MMM}_1(\T^d\x\RR_+))\x L_{w^*}^\infty(0,T;\MMM(\T^d))\to L_{w^*}^\infty(0,T;\MMM(\T^d))$ is the projection on the second coordinate. Let us also denote by $\bs{p}^1\colon L_{w^*}^\infty(0,T;\bbar{\MMM}_1(\T^d\x\RR_+))\x L_{w^*}^\infty(0,T;\MMM(\T^d))\to L_{w^*}^\infty(0,T;\bbar{\MMM}_1(\T^d\x\RR_+))$ the projection on the first coordinate. By Proposition~\ref{Leb} the law $Q^1:=\bs{p}^1_\sharp\bs{Q}_\Psi\in\Lim_{\ell,N\ra+\infty}\bs\pi^{N,\ell}_\sharp P^N$ is concentrated on the $w^*$-closed set $L_{w^*}^\infty(0,T;\bbar{\Y}_1(\T^d))$ and 
$$L_{w^*}^\infty(0,T;\bbar{\Y}_1(\T^d))\subs\bar{B}_{\bbar{\Psi}}^{-1}\big(L_{w^*}^\infty(0,T;\MMM_{ac}(\T^d))\big)$$
since $\bar{B}_{\bbar{\Psi}}(\bs\pi)=b_{\bbar{\Psi}}\circ j^*(\bs\pi)\df\LL_{\T^d}=b_{\bbar{\Psi}}(\bs\rho_{\bs\pi})\df\LL_{\T^d}$ by~\eqref{PsiBaryc} due to the fact that $\Psi$ is sublinear. Therefore
\begin{align*}
Q_\Psi^2\big(L_{w^*}^\infty(0,T;\MMM_{ac}(\T^d))\big)&=p^2_\sharp\bs{Q}_\Psi\big(L_{w^*}^\infty(0,T;\MMM_{ac}(\T^d))\big)=(\bar{B}_{\bbar{\Psi}})_\sharp Q^1\big(L_{w^*}^\infty(0,T;\MMM_{ac}(\T^d))\big)\\
&=Q^1\big(\bar{B}_{\bbar{\Psi}}^{-1}\big(L_{w^*}^\infty(0,T;\MMM_{ac}(\T^d))\big)\big)\\
&\geq Q^1\big(L_{w^*}^\infty(0,T;\bbar{\Y}_1(\T^d))\big)=1.
\end{align*}
For the second claim we note if $\psi_c<+\infty$ then we similarly have that 
\begin{align*}
Q_\Psi^2\big(L_{w^*}^\infty(0,T;B_{L^\infty(\T^d)}(0,\psi_c))\big)&=Q^1\big(\bar{B}_{\bbar{\Psi}}^{-1}\big(L_{w^*}^\infty(0,T;B_{L^\infty(\T^d)}(0,\psi_c))\big)\big)\\
&\geq Q^1\big(L_{w^*}^\infty(0,T;\bbar{\Y}_1(\T^d))\big)=1
\end{align*}
which completes the proof.$\hfill\Box$\\

We prove next statement (d). For all $G\in L^1(0,T;C_+(\T^d))$ $$\lls G,\s^N\rrs\leq\|\mathfrak{g}'\|_\infty\int_0^T\fr{N^d}\sum_{x\in\T_N^d}G_t\Big(\frac{x}{N}\Big)\eta_t(x)\df t=\|\mathfrak{g}'\|_\infty\lls G,\pi^N\rrs.$$
For any $G\in L^1(0,T;C(\T^d))$ the functional $J_G\colon\W\to\RR$ given by 
$$J_G(\w)=\|\mathfrak{g}'\|_\infty\lls G,i(\pi)\rrs-\lls G,\s\rrs,\quad\w=(\pi,W,\s)\in\W$$
is continuous, where here $i\colon D(0,T;\MMM_+(\T^d))\to L_{w^*}^\infty(0,T;\MMM(\T^d))$ is the continuous injection given in Proposition~\ref{SkorohodToLEmbed}. Therefore the set $\{J_G\geq 0\}$ is a closed set in $\W$ and thus by the portmanteau theorem for any subsequential limit point $Q=\lim_{N\ra\infty}Q^{k_N}$ of $\{Q^N\}$ 
\begin{align*}
Q\{J_G\geq 0\}&\geq\limsup_{N\ra+\infty}Q^{k_N}\big\{J_G\geq 0\big\}=\limsup_{N\ra\infty}P^{k_N}\big\{\|\mathfrak{g}'\|_\infty\lls G,i(\pi)\rrs\geq\lls G,\s\rrs\big\}=1.
\end{align*}
Since $L^1(0,T;C(\T^d))$ is separable it follows that 
$$Q\Big(\bigcap_{G\in L^1(0,T;C_+(\T^d))}\{J_G\geq 0\}\Big)=1.$$

Let $(\pi,W,\s)\in\bigcap_{G\in L^1(0,T;C_+(\T^d))}\{J_G\geq 0\}$. We will show that there exists a subset $E\subs[0,T]$ of full measure in $[0,T]$ such that for any $G\in C_+(\T^d)$ and any $t\in E$ it holds that $\ls G,\s_t\rs\leq\|\mathfrak{g}'\|_\infty\ls G,\pi_t\rs$. Ideed, for any $G\in C_+(\T^d)$ the maps $\ls G,\s_\cdot\rs$ and $\ls G,\pi_\cdot\rs$ are measurable and thus by Lebesgue's differentiation theorem there exists a measurable set $E_G\subs[0,T]$ of full measure in $[0,T]$ such that 
$$\lim_{\ee\ra 0}\fr{2\ee}\int_{t-\ee}^{t+\ee}\ls G,\s_s\rs\df s=\ls G,\s_t\rs,\quad \lim_{\ee\ra 0}\fr{2\ee}\int_{t-\ee}^{t+\ee}\ls G,\pi_s\rs\df s=\ls G,\pi_t\rs,\quad\forall t\in E_G.$$
For any $\ee>0$ the map $G_\ee(t)=\fr{2\ee}\1_{[t-\ee,t+\ee]}G$ is in $L^1(0,T;C_+(\T^d))$ and thus 
$$0\leq\|\mathfrak{g}'\|_\infty\lls G_\ee,\pi\rrs-\lls G_\ee,\s\rrs=\frac{\|\mathfrak{g}'\|_\infty}{2\ee}\int_{t-\ee}^{t+\ee}\ls G,\pi_s\rs\df s-\fr{2\ee}\int_{t-\ee}^{t+\ee}\ls G,\s_s\rs\df s.$$
Taking the limit as $\ee\ra 0$ it follows that 
$$0\leq\|\mathfrak{g}'\|_\infty h_G(t)-f_G(t)=\|\mathfrak{g}'\|_\infty\ls G,\pi_t\rs-\ls G,\s_t\rs,\quad\forall t\in E_G.$$
Since $C(\T^d)$ is separable the set $E:=\bigcap_{G\in C_+(\T^d)}E_G$ is of full measure and $\ls G,\s_t\rs\leq\|\mathfrak{g}'\|_\infty\ls G,\pi_t\rs$ for all $G\in C_+(\T^d)$ and all $t\in E$. According to~\eqref{LipMeasNormInTermsOfContFunct} this shows that $\Lip_{\nu_t}(\mu_t)<+\infty$ for all $t\in E$ and completes the proof of (d). The proof of (e) follows by the energy estimate of Theorem~\ref{TheoremEnergyEstimate} proved in the next section and the equality $W=-\nabla\s$ of equation~\eqref{CE}.$\hfill\Box$

\subsection{The energy estimate}\label{EnergySection}
In this section we prove the energy estimate Theorem~\ref{TheoremEnergyEstimate}.
The energy estimate is based on Lemma~\ref{LimitingDiffRateIsSobolevCoreLemma} below, whose proof is contained in the proof in Lemma 5.7.3 in~\cite{Kipnis1999a}. We state this lemma below for the convenience of the reader but before we do so let us fix some notation. For each $j=1,\dots,d$, $N\in\NN$, $\ee>0$ and each function $H\in C([0,T]\x\T^d)$ we define the function $V^N_{\ee,H}\equiv V^{N,j}_{\ee,H}:[0,T]\x\MM_N^d\lra\RR$ by the formula 
\begin{align*}
	V^N_{\ee,H}(t,\eta)&=\fr{N^{d-1}}\sum_{x\in\T_N^d}H_t\Big(\frac{x}{N}\Big)\frac{\mathfrak{g}\big(\eta(x)\big)-\mathfrak{g}\big(\eta(x+[N\ee]e_j)\big)}{[N\ee]}\\
	&\quad\quad\qquad-\frac{2}{N^d}\sum_{x\in\T_N^d}H_t\Big(\frac{x}{N}\Big)^2\fr{[N\ee]}\sum_{k=0}^{[N\ee]}\mathfrak{g}\big(\eta(x+ke_j)\big).
\end{align*}
The map $V^N_{\ee,H}$ induces a continuous map $\bar{V}^N_{\ee,H}\colon D(0,T;\MM_N^d)\to D(0,T;\RR)$ via the formula $\bar{V}^N_{\ee,H}(\eta)(t)=V^N_{\ee,H}(t,\eta_t)$ and since the integral map $\int_0^T\colon D(0,T;\RR)\to\RR$ that assigns to each path $f\in D(0,T;\RR)$ the number $\int_0^Tf(t)\df t$ is continuous, the map 
$$D(0,T;\MM_N^d)\ni\eta\mapsto \int_0^T \bar{V}^N_{\ee,H}(\eta)(t)\df t\in\RR$$
is also continuous and thus Borel measurable. Finally we denote by $(V_t^{N,\ee,H})_{0\leq t\leq T}$ the corresponding canonical process i.e.~$V_t^{N,\ee,H}(\eta)=\bar{V}^N_{\ee,H}(\eta)(t)$, $\eta\in D(0,T;\MM_N^d)$.
\begin{lemma}\label{LimitingDiffRateIsSobolevCoreLemma} Let $\{H^i\}_{i=1}^m\subs C^1([0,T]\x\T^d)$, $m\in\NN$, be a finite sequence of functions and let $\{\mu_0^N\in\PP_1\MM_N^d\}_{N\in\NN}$ be a sequence of initial distributions satisfying the $O(N^d)$-entropy assumption for some finite constant $C_0>0$. Then for all $\ee>0$ we have that 
	$$\limsup_{N\ra\infty}\EE^{\mu_0^N}\Big\{\max_{1\leq i\leq m}\int_0^TV^{N,\ee,H^i}_t\df t\Big\}\leq C_0.$$\end{lemma}
\begin{cor}\label{CorEn1}
	Let $\{H^i\}_{i=1}^m\subs C^{0,1}(I\x\T^d)$, $m\in\NN$, be a finite sequence of functions and let $$Q_3\in\Lim_{N\ra+\infty}Q_3^N\subs\PP L^\infty_{w^*}(0,T;\MMM_+(\T^d)),$$ where $Q_3^N:=\s^N_\sharp P^N$ is the law of the empirical jump rate process. Then 
	$$\int\Big\{\max_{1\leq i\leq m}\int_0^T\int_{\T^d}\big[\pd_jH^i_t(u)-H^i_t(u)^2\big]\df\s_t(u)\df t\Big\}\df Q_3(\s)\leq C_0.$$
\end{cor}\textbf{Proof} For each $H\in C(\T^d)$, $N\in\NN$ and $\ee>0$ we set $\pd_{-j}^{N,\ee}H(u):=\frac{N}{[N\ee]}\big(H(u)-H\big(u-\frac{[N\ee]}{N}e_j\big)\big)$ then for each $H\in C([0,T]\x\T^d)$ by a simple summation by parts we can write $V^{N,\ee,H}$ as 
\begin{align*}
	V_t^{N,\ee,H}&=\fr{N^{d-1}}\sum_{x\in\T_N^d}\frac{H_t(\frac{x}{N})-H_t\big(\frac{x-[N\ee]e_j}{N}\big)}{[N\ee]}\mathfrak{g}\big(\eta_t(x)\big)\\
	&\quad\quad\qquad-\frac{2}{N^d}\sum_{x\in\T_N^d}\fr{[N\ee]}\sum_{k=0}^{[N\ee]}H_t\Big(\frac{x-ke_j}{N}\Big)^2\mathfrak{g}\big(\eta_t(x)\big)\\
	&=\fr{N^d}\sum_{x\in\T_N^d}\Big\{\pd_{-j}^{N,\ee}H_t\Big(\frac{x}{N}\Big)-\frac{2}{[N\ee]}\sum_{k=0}^{[N\ee]}\tau_{-\frac{k}{N}e_j}H_t\Big(\frac{x}{N}\Big)^2\Big\}\mathfrak{g}\big(\eta_t(x)\big)\\
	&=\Big\ls\pd_{-j}^{N,\ee}H_t-\frac{2}{[N\ee]}\sum_{k=0}^{[N\ee]}\tau_{-\frac{k}{N}e_j}H_t^2,\s_t^N\Big\rs.
\end{align*}
Therefore if we denote by $A_j^{N,\ee}\colon C(\T^d)\to C(\T^d)$ the map 
$$A^{N,\ee}_j(H):=\pd_{-j}^{N,\ee}H-\frac{2}{[N\ee]}\sum_{k=0}^{[N\ee]}\tau_{-\frac{k}{N}e_j}H^2,\quad H\in C(\T^d)$$
then we can write $\int_0^TV^{N,j,\ee,H}_t\df t=\int_0^T\ls A^{N,\ee}_j(H_t),\s_t^N\rs\df t=\lls A^{N,\ee}_jH,\s^N\rrs$ and 
$$\EE^{\mu_0^N}\Big\{\max_{1\leq i\leq m}\int_0^TV_t^{N,j,\ee,H^i}\df t\Big\}
=\int\max_{1\leq i\leq m}\lls A_j^{N,\ee}(H^i),\s\rrs\df Q_3^N(\s)$$
where with a slight abuse of notation we write $A_j^{N,\ee}=\bar{A}^{N,\ee}_j$ also for the induced map on the $L^1$-Bochner spaces. We claim first that for each $\ee>0$ and each $H\in C(\T^d)$
\begin{equation}\label{AjNepsilonConv} 
\lim_{N\ra+\infty}A_j^{N,\ee}(H)=\frac{H-\tau_{-\ee e_j}H}{\ee}-\frac{2}{\ee}\int_0^\ee H(\cdot-se_j)^2\df s=:A_j^\ee(H)\end{equation}
uniformly over all $u\in\T^d$. The fact that 
$$\lim_{N\ra+\infty}\frac{N}{[N\ee]}\Big[H-\tau_{-\frac{[N\ee]}{N}e_j}H\Big]=\frac{H-\tau_{-\ee e_j}H}{\ee}=:\pd^\ee_{-j}H$$ uniformly on $\T^d$ is obvious and so we prove that 
\begin{align}\label{SecTermInfLim}
\lim_{N\ra+\infty}\fr{[N\ee]}\sum_{k=0}^{[N\ee]}\tau_{-\frac{k}{N}e_j}H^2=\fr{\ee}\int_0^\ee H(\cdot-se_j)^2\df s
\end{align} uniformly on $\T^d$. Indeed, for each $u\in\T^d$ 
$$
	\fr{[N\ee]}\sum_{k=0}^{[N\ee]}\tau_{-\frac{k}{N}e_j}H(u)^2
	=\frac{N}{[N\ee]}\sum_{k=0}^{[N\ee]}\int_{\frac{k}{N}}^{\frac{k+1}{N}}H\Big(u-\frac{k}{N}e_j\Big)^2\df s
	=\frac{N}{[N\ee]}\int_0^{\frac{[N\ee]+1}{N}}H\Big(u-\frac{[Ns]}{N}e_j\Big)^2\df s.
$$
So for all $u\in\T^d$ 
\begin{align*}\fr{[N\ee]}\sum_{k=0}^{[N\ee]}\tau_{-\frac{k}{N}e_j}H(u)^2-\fr{\ee}\int_0^\ee H(u-se_j)^2\df s&=\fr{\ee}\int_0^\ee \Big[H\Big(u-\frac{[Ns]}{N}e_j\Big)^2-H(u-se_j)^2\Big]\df s\\
	&\quad\qquad+\Big(\frac{N}{[N\ee]}-\fr{\ee}\Big)\int_0^\ee H\Big(u-\frac{[Ns]}{N}e_j\Big)^2\df s\\
	&\quad\qquad\quad+\frac{[N\ee]}{N}\int_\ee^{\frac{[N\ee]+1}{N}}H\Big(u-\frac{[Ns]}{N}e_j\Big)^2\df s
\end{align*}
and therefore 
\begin{align*}
	\Big|\fr{[N\ee]}\sum_{k=0}^{[N\ee]}\tau_{-\frac{k}{N}e_j}H(u)^2-\fr{\ee}\int_0^\ee H(u-se_j)^2\df s\Big|&\leq\frac{2\|H\|_\infty}{\ee}\int_0^\ee\Big|H\Big(u-\frac{[Ns]}{N}e_j\Big)-H(u-se_j)\Big|\df s\\
	&+\|H\|_\infty^2\Big(\frac{N\ee}{[N\ee]}-1\Big)+\frac{[N\ee]}{N}\|H\|_\infty^2\Big(\frac{[N\ee]+1}{N}-\ee\Big).
\end{align*}
The second and third term in the right hand side above are independent of the the variable $u\in\T^d$ and obviously converge to $0$ as $N\ra+\infty$. The first term also converges to $0$ as $N\ra+\infty$ uniformly over all $u\in\T^d$ since $H\in C(\T^d)$ is uniformly continuous by the compactness of $\T^d$ and thus the limit in~\eqref{AjNepsilonConv} holds uniformly on $\T^d$. Furthermore, for all $H\in C(\T^d)$ and all $\ee>0$ there exists large enough $N_{H,\ee}\in\NN$ 
such that \begin{equation}\label{forBcon2}\quad
\sup_{N\geq N_{H,\ee}}\|A_j^{N,\ee}(H)\|_\infty\leq
2\|H\|_\infty+\|H\|_\infty^2<+\infty.
\end{equation}
\indent Thus if $H\in L^2(0,T;C(\T^d))\subs L^1(0,T;C(\T^d))$ then $\lim_{N\ra\infty}\big\|A_j^{N,\ee}(H_t)-A^\ee_j(H_t)\big\|_\infty=0$ for all $t\in[0,T]$ and $$\sup_{N\geq N_{H,\ee}}\|A_j^{N,\ee}(H_\cdot)\|_\infty\leq
2\|H_\cdot\|_\infty+\|H_\cdot\|_\infty^2\in L^1(0,T).$$
Therefore by the dominated convergence theorem  
\begin{equation}\label{onedown}
\lim_{N\ra+\infty}\|A_j^{N,\ee}(H)-A_j^\ee(H)\|_{L^1(0,T;C(\T^d))}=0,\quad\forall\;H\in L^2(0,T;C(\T^d)).\end{equation}

We consider now a subsequence of $\{Q_3^{N}\}_{N\in\NN}$, which we continue to denote by $\{Q_3^{N}\}$, converging weakly to $Q$. Then, using the elementary inequality 
\begin{equation}\label{MaxDifferenceIneq}
\max_{1\leq i\leq m}a_i-\max_{1\leq i\leq m}b_i\leq\max_{1\leq i\leq m}(a_i-b_i)
\end{equation} which holds for all finite sequences $\{a_i\}_{i=1}^m$, $\{b_i\}_{i=1}^m$ of real numbers, we write 
\begin{align*}
	\int\max_{1\leq i\leq m}\lls A_j^\ee(H^i),\s\rrs\df Q_3^{N}(\s)&\leq\int\max_{1\leq i\leq m}\lls A_j^{N,\ee}(H^i),\s\rrs\df Q_3^{N}(\s)\\
	&\quad+\int\max_{1\leq i\leq m}\lls A_j^\ee(H^i)-A^{N,\ee}(H^i),\s\rrs\df Q_3^{N}(\s).
\end{align*}
The function $L^\infty_{w^*}(0,T;\MMM_+(\T^d))\ni\s\mapsto\max_{1\leq i\leq m}\lls\y_\ee(H^i),\s\rrs$ is continuous in the $w^*$-topology as the maximum of a finite number of continuous functionals and therefore since $Q_3^N$ converges weakly to $Q_3$ we have that 
$$\lim_{N\ra\infty}\int\max_{1\leq i\leq m}\lls A^\ee(H^i),\s\rrs\df Q_3^N(\s)=\int\max_{1\leq i\leq m}\lls A^\ee(H^i),\s\rrs\df Q_3(\s).$$
On the other hand, by Lemma \ref{LimitingDiffRateIsSobolevCoreLemma}  
$$\limsup_{N\ra\infty}\int\max_{1\leq i\leq m}\lls A^{N,\ee}(H^i),\s\rrs\df Q_3^N(\s)\leq C_0$$
and thus it follows that 
\begin{equation}\label{C0PlusLimSup}
\int\max_{1\leq i\leq m}\lls A_j^\ee(H^i),\s\rrs\df Q_3(\s)\leq C_0+\limsup_{N\ra\infty}\int\max_{1\leq i\leq m}\lls A_j^\ee(H^i)-A^{N,\ee}(H^i),\s\rrs\df Q_3^N(\s).
\end{equation}

Now, the limit superior in the right hand side above vanishes since 
$$\int\max_{1\leq i\leq m}\lls A_j^\ee(H^i)-A_j^{N,\ee}(H^i),\s\rrs\df Q_3^N(\s)\leq\EE^N\Big(\max_{1\leq i\leq m}|\lls A_j^\ee(H^i)-A_j^{N,\ee}(H^i),\s^N\rrs|\Big)$$
and for all $i=1,\dots,m$ we have that $P^N$-a.s$.$ 
\begin{align*}\big|\lls A_j^\ee(H^i)-A_j^{N,\ee}(H^i),\s^N\rrs\big|&\leq
	\int_0^T\big|\ls A_j^\ee(H_t^i)- A_j^{N,\ee}(H^i_t),\s_t^N\rs\big|\df t\\
	&\leq\|\mathfrak{g}'\|_\infty\int_0^T\|A_j^\ee(H^i_t)-A_j^{N,\ee}(H^i_t)\|_\infty\ls 1,\pi_t^N\rs\df t\\
	&=\|\mathfrak{g}'\|_\infty\ls 1,\pi_0^N\rs\int_0^T\|A_j^\ee(H^i_t)-A_j^{N,\ee}(H^i_t)\|_\infty\df t
\end{align*}
so that the integral in the limit superior in the right hand side of~\eqref{C0PlusLimSup} is bounded above by
$$\|\mathfrak{g}'\|_\infty\max_{1\leq i\leq m}\int_0^T\|A_j^\ee(H^i_t)-A_j^{N,\ee}(H^i_t)\|_\infty\df t\int\ls1,\pi^N\rs\df\mu_0^N$$
which tends to zero as $N\ra\infty$ by~\eqref{onedown} and Lemma~\ref{BoundedTotParticlInMacrLim}. Therefore 
\begin{equation}\label{almostthere}
\int\max_{1\leq i\leq m}\lls A_j^\ee(H^i),\s\rrs\df Q_3(\s)\leq C_0.
\end{equation}

We show finally that if $H\in L^1(0,T;C^1(\T^d))\cap L^2(0,T;C(\T^d))$ then 
\begin{equation}\label{ToCompleteCorEn1}
\lim_{\ee\ra+\infty}\int\max_{1\leq i\leq m}\lls A_j^\ee(H^i),\s\rrs\df Q_3(\s)=\int\max_{1\leq i\leq m}\lls \pd_jH_t-2H_t^2,\s\rrs\df Q_3(\s)
\end{equation}
which concludes the proof. We note first that for each function $H\in C^1(\T^d)$  
$$\lim_{\ee\ra 0}A^\ee(H)=\pd_jH-2H^2=:A(H)$$ uniformly in $\T^d$ since the function $\pd_jH-2H^2$ is uniformly continuous. Indeed, by the fundamental theorem of calculus and a simple change of variables 
$$A^\ee(H)(u)=\int_0^1\pd_jH(u-\ee se_j)\df s-\frac{2}{\ee}\int_0^\ee H(u-se_j)^2\df s=\int_0^1\big\{\pd_jH(u-\ee se_j)-2H(u-\ee se_j)^2\big\}\df s$$ and thus 
$$|A^\ee(H)(u)-A(H)(u)|\leq\int_0^1|\pd_jH(u-\ee se_j)-\pd_jH(u)+2H(u)^2-2H(u-\ee se_j)^2|\df s.$$
This shows that $\lim_{\ee\ra 0}\|A^\ee(H)-A(H)\|_\infty=0$ for all $H\in C^1(\T^d)$. Furthermore 
$$\sup_{\ee>0}\|A^\ee(H)-A((H)\|_\infty\leq2\|\pd_jH\|_\infty+4\|H\|_\infty^2$$
and therefore if $H\in L^1(0,T;C^1(\T^d))\cap L^2(0,T;C(\T^d))$ then the family of maps $\{\|A^\ee(H_\cdot)-A(H)\|_\infty\}_{\ee>0}$ is dominated by the function $2\|\pd_jH_\cdot\|_\infty+4\|H_\cdot\|^2_\infty\in L^1(0,T)$ and thus
\begin{equation}\label{L1BochAConv}
\lim_{\ee\ra 0}\|A^\ee(H)-A(H)\|_{L^1(0,T;C(\T^d))}=0,\quad\forall\;H\in L^1(0,T;C(\T^d)^1)\cap L^2(0,T;C(\T^d)).
\end{equation}
Consequently, if $\{H^i\}_{i=1}^m\subs L^1(0,T;C^1(\T^d))\cap L^2(0,T;C(\T^d))$ then by inequality~\eqref{MaxDifferenceIneq}
\begin{align}\label{B}
\Big|\int\Big\{\max_{1\leq i\leq m}\lls A_j^\ee(H^i),\s\rrs&-\max_{1\leq i\leq m}\lls A(H),\s\rrs\Big\}\df Q_3(\s)\Big|\leq\int\max_{1\leq i\leq m}|\lls A^\ee(H_i)-A(H^i),\s\rrs|\df Q_3(\s)\nonumber\\
&\leq\int \int_0^T\max_{1\leq i\leq m}\|A^\ee(H_t^i)-A(H_t^i)\|_\infty\|\s_t\|_{TV}\df t\df Q_3(\s)\nonumber\\
&\leq\max_{1\leq i\leq m}\|A^\ee(H^i)-A(H^i)\|_{L^1(0,T;C(\T^d))}\int\|\s\|_{TV;\infty}\df Q_3(\s).
\end{align}
But by the $w^*$-lower semicontinuity of the norm $\|\cdot\|_{TV;\infty}$ of $L_{w^*}^\infty(0,T;\MMM(\T^d))$ and Lemma~\ref{BoundedTotParticlInMacrLim}
\begin{align}\label{sIntegr}
\int\|\s\|_{TV;\infty}\df Q_3(\s)&\leq\liminf_{N\ra+\infty}\int\|\s^N\|_{TV;\infty}\df P^N\leq\|\mathfrak{g}'\|_\infty\limsup_{N\ra+\infty}\int\|\pi^N\|_{TV;\infty}\df P^N\nonumber\\
&=\|\mathfrak{g}'\|_\infty\limsup_{N\ra+\infty}\int\ls 1,\pi^N\rs\df\mu_0^N<+\infty.
\end{align}
Consequently~\eqref{ToCompleteCorEn1} follows by~\eqref{L1BochAConv},~\eqref{sIntegr} and~\eqref{B} and the proof of the corollary is complete.$\hfill\Box$\\

Let $K_0\colon L_{w^*}^\infty(0,T;\MMM_+(\T^d))\to[0,+\infty]$ be the map defined by 
\begin{equation}\label{KZeroMap}
K_0(\s)=\sup\limits_{H\in C^{0,1}([0,T]\x\T^d)}\int_0^T\int\big(\pd_jH_t(u)-2H_t(u)^2\big)\df\s_t(u)\df t.
\end{equation} Then $K_0$ is $w^*$-measurable and 
\begin{equation}\label{IntegrEnergBound}
\int K_0(\s)\df Q_3(\s)\leq C_0<+\infty.\end{equation} Indeed, if $\{H_i\}_{i\in\NN}\subs C^{0,1}([0,T]\x\T^d)$ is a sequence dense in $C^{0,1}([0,T]\x\T^d)$ in the usual $C^{0,1}$-uniform norm $\|H\|_{C^{0,1}}:=\|H\|_{C([0,T]\x\T^d)}+\|\nabla H\|_{C([0,T]\x\T^d)}$ then for all $\s\in L^\infty_{w^*}(0,T;L^\infty(\T^d))$ 
\begin{align*}
	K_0(\s)&=\sup_{i\in\NN}\int_0^T\int_{\T^d}\big[\pd_jH^i_t(u)-H^i_t(u)^2\big]\s_t(u)\df u\df t\\
&=\lim_{m\ra\infty}\max_{1\leq i\leq m}\int_0^T\int_{\T^d}\big[\pd_jH^i_t(u)-H^i_t(u)^2\big]\s_t(u)\df u\df t. 
\end{align*}
So the map $K_0$ is $w^*$-lower semicontinuous as the supremum of a family of $w^*$-lower semicontinuous functions and~\eqref{IntegrEnergBound} follows by the monotone convergence theorem and Corollary~\ref{CorEn1}.$\hfill\Box$\\

We are now ready to prove Theorem~\ref{TheoremEnergyEstimate}. Fix $j=1,\dots,d$. By~\eqref{IntegrEnergBound} for $Q_3$-a.s.~all paths $\s\in L^\infty(0,T;\MMM_+(\T^d))$ 
\begin{equation}\label{BE} K_0(\s)<+\infty.
\end{equation}
For each a path $\s\in\{K_0<+\infty\}$ we denote by $\bs\s\in\MMM_+([0,T]\x\T^d)$ the corresponding space-time measure characterized by 
$$\ls H,\bs\s\rs=\int_0^T\int_{\T^d}H(t,u)\df\s_t\df t,\quad H\in C([0,T]\x\T^d).$$
In other words $\bs\s:=i^*(\s)$ where $i\colon C([0,T]\x\T^d)\hookrightarrow L^1(0,T;C(\T^d))$ is the natural injection. For each $\s\in L_{w^*}^\infty(0,T;\MMM_+(\T^d))$ we define $L^2_\s$ to be the closure of $C^{0,1}([0,T]\x\T^d)$ in $L^2([0,T]\x\T^d,{\bs\s})$. Then $L^2_\s$ is a Hilbert space with the inner product 
$$\ls H,G\rs_{\bs\s}=\int_{[0,T]\x\T^d}H(t,u)G(t,u)\df\bs\s(t,u)=\int_0^T\int_{\T^d}H(t,u)G(t,u)\df\s_t(u)\df t.$$

Let now $\ell_\s^j\colon C^{0,1}([0,T]\x\T^d)\to\RR$, $j=1,\ldots,d$ denote the linear function given by the formula 
$$\ell_\s^j(H)=\int_0^T\int_{\T^d}\pd_jH(t,u)\df\s_t(u)\df t.$$
It follows from estimate~\eqref{BE} that 
$$a\ell_\s^j(H)-2a^2\|H\|_{L_\s^2}^2\leq K_0(\s)$$
for all $a\in\RR$ and all $H\in C^1([0,T]\x\T^d)$. The maximum over all $a\in\RR$ of the quantity in the left hand side of the inequality above is achieved at $a=\ell_\s^j(H)/\|2H\|_{L_\s^2}^2$ and therefore 
$$\frac{\ell_\s^j(H)^2}{8\|H\|_{L_\s^2}^2}=\frac{\ell_\s^j(H)^2}{4\|H\|_{L_\s^2}^2}-2\frac{\ell_\s^j(H)^2}{16\|H\|_{L_\s^2}^4}\|H\|_{L_\s^2}^2\leq K_0(\s)$$
for all $H\in C^1(I\x\T^d)$. It follows that $|\ell_\s^j(H)|\leq 2\sqrt{2K(\s)}\|H\|_{L_\s^2}$ for all $H\in C^1([0,T]\x\T^d)$ and thus $\ell_\s^j$ can be extended to a bounded linear function $\ell_\s^j\colon L_\s^2\lra\RR$ with norm $\|\ell_\s^j\|\leq 2\sqrt{2K_0(\s)}$.

By the Riesz representation theorem now there exists an $L_\s^2$-function, which we denote by $\pd_j(\log\s)$, such that 
\begin{equation}\label{endthis}
\ell_\s^j(H)=-\ls H,\pd_j(\log\s)\rs_{\bs\s}=-\int_0^T\int_{\T^d}H(t,u)\pd_j(\log\s)(t,u)\df\bs\s(t,u)
\end{equation}
for all $H\in C^{0,1}([0,T]\x\T^d)$. Since $\pd_j(\log\s)\in L_\s^2$ represents $\ell_\s^j$ via the Riesz representation theorem we have that 
\begin{equation}\label{enestinf}
\int_0^T\int_{\T^d}\big[\pd_j(\log\s)(t,u)\big]^2\df\bs\s(t,u)=\|\pd_j(\log\s)\|_{L_\s^2}^2=\|\ell_\s^j\|^2\leq 8K_0(\s)<+\infty.
\end{equation}
Now, by Theorem~\ref{Theorem1}(c) we know that $Q_3$ is concentrated on the $w^*$-measurable $L^\infty_{w^*}(0,T;\MMM_{+,ac}(\T^d))$ and thus for $Q_3$-almost all $\s\in\{K_0<+\infty\}$ we have that $\bs\s\ll\LL_{[0,T]\x\T^d}$ with density 
$$\s(t,u):=\frac{\df\bs\s}{\df\LL_{[0,T]\x\T^d}}(t,u)=\frac{\df\s_t}{\df\LL_{\T^d}}(u)\geq 0$$ for Lebesgue almost all $(t,u)\in[0,T]\x\T^d$. Since the Radon-Nikodym density $\s$ is Lebesgue almost surely equal to $0$ on any $\bs\s$-null set $E\subs[0,T]\x\T^d$ the Lebesgue almost sure equality class of the function $\pd_j\s:=\s\cdot\pd_j(\log\s)$ does not depend on the representative from the $\bs\s$-almost sure equality class of the function $\pd_j(\log\s)\in L^2_\s$. Furthermore the Lebesgue integral of $|\pd_j\s|$ is 
\begin{align*}
\int_{[0,T]\x\T^d}|\pd_j\s(t,u)|\df u\df t&=\int_{[0,T]\x\T^d}|\pd_j(\log\s)(t,u)|\df\bs\s(t,u)\\
&\leq\sqrt{\bs\s([0,T]\x\T^d)}\Big(\int_{[0,T]\x\T^d}|\pd_j(\log\s)|^2\df\bs\s(t,u)\Big)^{\fr{2}}\\
&\leq 2\sqrt{\bs\s([0,T]\x\T^d)2K_0(\s)}\leq 2\sqrt{2T\|\mathfrak{g}'\|_\infty\mu_0(\T^d)K_0(\s)}<+\infty
\end{align*}
and thus $\pd_j(\log\s)\in L^1([0,T]\x\T^d)$ for all $\s\in L_{w^*}^\infty(0,T;\MMM_{+,ac}(\T^d))\cap\{K_0<+\infty\}$. By~\eqref{endthis} $\pd_j\s$ satisfies property~\eqref{MeanTimeSob} and is thus the $L^1$-weak $j$-th spatial derivative of $\s$ in $[0,T]\x\T^d$. By the identity $\pd_j\s=\s\cdot\pd_j(\log\s)$ we have that 
$$\frac{|\nabla\s(t,u)|^2_2}{\s(t,u)}=\s(t,u)|\nabla(\log\s)(t,u)|_2^2$$ and therefore the energy estimate~\eqref{energyestimate} follows from~\eqref{enestinf}. Finally, if $\phih_c<+\infty$ then as we know from Theorem~\ref{Theorem1}(c) the law $Q_3$ is concentrated on the $w^*$-closed subspace of paths $\s\in L_{w^*}^\infty(0,T;\MMM_{+,ac})(\T^d))$ that satisfy $$\big\|\|\s_t\|_{L^\infty(\T^d)}\big\|_{L^\infty(0,T)}\leq\phih_c<+\infty$$ and thus the function $\pd_j\s:=\s\cdot\pd_j(\log\s)$ satisfies 
\begin{align*}
\int_{[0,T]\x\T^d}\big(\pd_j\s(t,u)\big)^2\df u\df t&=\int_0^T\int_{\T^d}\big[\pd_j(\log\s)(t,u)\big]^2\s(t,u)^2\df u\df t\leq\phih_c\|\pd_j(\log\s)\|_{L_\s^2}^2\\
&\leq 8\phih_c K_0(\s)<+\infty,
\end{align*}
and thus is in $L^2([0,T]\x\T^d)$. By~\eqref{endthis} $\pd_j\s$ satisfies property~\eqref{MeanTimeSob} and is thus the required $L^2$-weak $j$-th spatial derivative of $\s$ in $[0,T]\x\T^d$. 

It is now easy to see that $\s_t\in H^1(\T^d)$ ($\s_t\in W^1(\T^d)$ if $\phih_c=+\infty)$ a.s.~for all $t\in[0,T]$ $Q_3$-a.s.~for all $\s\in L_{w^*}^\infty(0,T;\MMM_+(\T^d))$ since 
$$L_{w^*}^\infty(0,T;\MMM_{+,ac}(\T^d))\cap\{K_0<+\infty\}\subs\big\{\s\bigm|\s_t\in W^1(\T^d)\mbox{ a.s.~}\forall\;t\in[0,T]\big\}$$ and the set in the left hand side above is of full $Q_3$-measure and likewise if $\phih_c<+\infty$ then 
$$\big\{\s\bigm|\big\|\|\s_t\|_{L^\infty(\T^d)}\big\|_{L^\infty(0,T)}\leq\phih_c\big\}\cap\{K_0<+\infty\}\subs\big\{\s\bigm|\s_t\in H^1(\T^d)\mbox{ a.s.~}\forall\;t\in[0,T]\big\}$$ and the set in the left hand side is of full $Q_3$-measure. By the previous paragraph we know that and $\s\in L^\infty_{w^*}(0,T;\MMM_{+,ac}(\T^d))\cap\{K_0<+\infty\}$ there exist functions $\pd_j\s\in L^1([0,T]\x\T^d)$, $j=1,\dots,d$, satisfying~\eqref{MeanTimeSob} for all $H\in C^{0,1}([0,T]\x\T^d)$. We will show that $\s_t\in W^{1,1}(\T^d)$ for almost all $t\in[0,T]$. For each $t\in[0,T]$ and $\ee>0$ we consider a sequence of smooth functions $\{f^n_{t,\ee}\}_{n\in\NN}$ defined on $[0,T]$ such that $f^n_{t,\ee}\leq\1_{[t-\ee,t+\ee]}$ for all $n\in\NN$ and $f^n_{t,\ee}\lra\1_{(t-\ee,t+\ee)}$ pointwise as $n\ra\infty$. Then for all functions $H\in C^{0,1}([0,T]\x\T^d)$ we have by~\eqref{MeanTimeSob} that 
$$\fr{2\ee}\int_0^T\int_{\T^d}f^n_{t,\ee}(s)\pd_jH_s(u)\s(s,u)\df u\df s=-\fr{2\ee}\int_0^T\int_{\T^d}f^n_{t,\ee}(s)H_s(u)\pd_j\s(s,u)\df u\df s.$$
Then taking the limit as $n\ra\infty$ in both sides of the inequality above, we get that  
$$\int_{t-\ee}^{t+\ee}\int_{\T^d}\pd_jH_s(u)\s(s,u)\df u\df s=-\int_{t-\ee}^{t+\ee}\int_{\T^d}H_s(x)\pd_j\s(s,x)\df u\df s.$$
Then taking the limit as $\ee\ra 0$ in both sides of the equality above, it follows by Lebesgue's differentiation theorem that for each $H\in C^{0,1}([0,T]\x\T^d)$,  \begin{equation}\label{atpointsob}
\int_{\T^d}\pd_jH_t(x)\s(t,u)\df u=-\int_{\T^d}H_t(x)\pd_j\s(t,x)\df u
\end{equation} for all $t\in E_H$, for some measurable set $E_H\subs[0,T]$ of full Lebesgue measure in $[0,T]$. Taking then a sequence $\{H^i\}_{i\in\NN}\subs C^{0,1}([0,T]\x\T^d)$ dense in $C^{0,1}([0,T]\x\T^d)$ in the $C^{0,1}$-uniform norm $\|\cdot\|_{C^{0,1}}$, we have that the set $E:=\bigcap_{i\in\NN}E_{H^i}$ is of full Lebesgue measure and for each $t\in E$ we have that~\eqref{atpointsob} holds for all $H\in C^{0,1}([0,T]\x\T^d)$. In particular, since $C^1(\T^d)$ can be considered as a subspace of $C^1([0,T]\x\T^d)$ it follows that 
$$\int_{\T^d}\pd_jH(u)\s(t,u)\df u=-\int_{\T^d}H(u)\pd_j\s(t,u)\df u,\quad\forall\;(t,H)\in E\x C^1(\T^d).$$ Consequently, $\s_t$ is weakly differentiable for almost all $t\in[0,T]$ with weak $j$-th partial derivative $\pd_j\s_t$. Finally, since $\pd_j\s\in L^1([0,T]\x\T^d)$ we have that 
$$\int_0^T\|\pd_j\s_t\|_{L^1(\T^d)}\df t=\|\pd_j\s\|_{L^1([0,T]\x\T^d)}<+\infty,$$
and therefore $\|\pd_j\s_t\|_{L^1(\T^d)}<+\infty$ for almost all $t\in[0,T]$. Consequently, $\s_t\in W^{1,1}(\T^d)$ for almost all $t\in[0,T]$. If $\phih_c<+\infty$ then $\pd_j\s\in L^2([0,T]\x\T^d)$ and thus $\s_t\in H^1(\T^d)$ for almost all $t\in[0,T]$ and the proof of Theorem~\ref{TheoremEnergyEstimate} is complete.$\hfill\Box$

 \subsection{A closed hydrodynamic equation}\label{ACWYE}
In this section we prove Theorem~\ref{ClosedWeakEquation}. We will first show that any limit point $\bs{Q}$ of the laws $\bs{Q}^{N,\ell}:=\bs\pi^{N,\ell}_\sharp P^N$ of the micro empirical density process of the ZRP is concentrated on mild Young measure-valued solutions of the hydrodynamic equation in the sense of~\eqref{Mild}. Then by the energy estimate of Theorem~\ref{TheoremEnergyEstimate} it will follow that it is in fact concentrated on weak solutions in the sense of~\eqref{YMVZRPPDE}. So let $\bs{Q}\in\Lim_{\ell\ra+\infty}\Lim_{N\ra+\infty}\bs{Q}^{N,\ell}$. For any test function $f\in C_c^{1,2}((0,T)\x\T^d)$ the map 
$$(0,T)\ni t\mapsto \Lambda\pd_tf_t(U)+\Phi(\Lambda)\Delta f_t(U)\in\bbar{C}_1(\T^d\x\RR_+)$$
is in $L^1(0,T;\bbar{C}_1(\T^d\x\RR_+))$. Here $\Phi=\bbar{\Phi}$ is the extended mean jump rate function of the ZRP. Therefore the functional $$\lls \Lambda\pd f(U)+\Phi(\Lambda)\Delta f(U),\cdot\rrs\colon L_{w^*}^\infty(0,T;\bbar{\MMM}_1(\T^d\x\RR_+))\to\RR$$ is $w^*$-continuous, and thus for each $\ee>0$ and each $f\in C_c^{1,2}((0,T)\x\T^d)$ the set 
$$A_{f,\ee}:=\big\{\bs\pi\in L_{w^*}^\infty(0,T;\bbar{\MMM}_1(\T^d\x\RR_+)\bigm||\lls\Lambda\pd f(U)+\Phi(\Lambda)\Delta f(U),\bs\pi\rrs|>3\ee\big\}$$
is open in $w^*$-topology of $L_{w^*}^\infty(0,T;\bbar{\MMM}_1(\T^d\x\RR_+))$. Consequently for any sequences $\{m_\ell\}_{\ell\in\ZZ_+}$ and $\{k_N^{(\ell)}\}_{N\in\NN}$, $\ell\in\ZZ_+$,  such that 
$$\bs{Q}=\lim_{\ell\ra+\infty}\lim_{N\ra+\infty}\bs{Q}^{k_N^{(\ell)},m_\ell}$$ we have by the portmanteau theorem that 
\begin{align*}\bs{Q}(A_{f,\ee})&\leq\liminf_{\ell\ra+\infty}\liminf_{N\ra+\infty}\bs{Q}^{k_N^{(\ell)},m_\ell}(A_{f,\ee})\leq\limsup_{\ell\ra+\infty}\limsup_{N\ra+\infty}\bs{Q}^{N,\ell}(A_{f,\ee})
\end{align*}
By the definition of $\bs{Q}^{N,\ell}$ 
\begin{align*}
\bs{Q}^{N,\ell}(A_{f,\ee})&=P^N\Big\{\Big|\int_0^T\bs\pi_t^{N,\ell}\big(\Lambda\pd_tf_t(U)+\Phi(\Lambda)\Delta f_t(U)\big)\df t\Big|>3\delta\Big\}\\
&=
P^N\Big\{\Big|\int_0^T\fr{N^d}\sum_{x\in\T_N^d}\Big[\pd_tf_t\Big(\frac{x}{N}\Big)\eta_t^{\ell}(x)+\Delta f_t\Big(\frac{x}{N}\Big)\Phi\big(\eta^\ell_t(x)\big)\Big]\df t\Big|>3\delta\Big\}\\
&=P^N\big\{|\lls\pd_tf_t,\pi^{N,\ell}\rrs+\lls\Delta f_t,B_\Phi(\bs\pi^{N,\ell})\rrs|>3\delta\big\}.\end{align*}
In the last equality $\pi^{N,\ell}:=B(\bs\pi^{N,\ell})$ is the barycentric projection of the micro-empirical density and $B_\Phi$ is the $\Phi$-projection defined as in~\eqref{PsiProj}. By adding and subtracting the terms $\ls\pd_tf_t,\pi^N\rs$ and $\ls\Delta f_t,\s^N\rs$ it follows that 
\begin{align}\label{ThreeTerms}
\bs{Q}^{N,\ell}(A_{f,\ee})&\leq P^N\big\{\big|\lls\pd f,\pi^{N,\ell}-\pi^N\rrs\big|>\ee\big\}+P^N\Big\{\big|\lls\pd f,\pi^N\rrs+\lls\Delta f,\s^N\rrs\big|>\ee\big\}\nonumber\\
&\quad+P^N\big\{\big|\lls\Delta f,\s^N-B_\Phi(\bs\pi^{N,\ell})\rrs\big|>\ee\big\}.
\end{align}

The second term in~\eqref{ThreeTerms} converges to $0$ as $N\ra+\infty$ by~\eqref{MicroCE} for $f\in C_c^3((0,T)\x\T^d)$ and the third term converges to $0$ for $f\in C_c^2((0,T)\x\T^d)$ by the one-block estimate. For the first term by a change of variables we have for all $N\in\NN$, $\ell\in\ZZ_+$ that
\begin{align}\label{Init11}\lls\pd f,\pi^{N,\ell}-\pi^N\rrs&=\int_0^T\fr{N^d}\sum_{x\in\T_N^d}\pd_tf_t\Big(\frac{x}{N}\Big)\big(\eta_t^{\ell}(x)-\eta_t(x)\big)\df t\nonumber\\
&=\int_0^T\fr{N^d}\sum_{x\in\T_N^d}\Big[\fr{\ell_\star^d}\sum_{|y|\leq\ell}\pd_tf_t\Big(\frac{x+y}{N}\Big)-\pd_tf_t\Big(\frac{x}{N}\Big)\Big]\eta_t(x)\df t\nonumber\\
&=\Big\lls\fr{\ell_\star^d}\sum_{|y|\leq\ell}\pd f\Big(\frac{\cdot+y}{N}\Big)-\pd f,\pi^N\Big\rrs.\end{align}
Since $\pd f\in L^1(0,T;C(\T^d))$ for each $\ee>0$ there exists a map $\bar{\delta}_\ee\in L_+^\infty(0,T)$ such that $\bar{\delta}_\ee(t)>0$ for almost all $t\in[0,T]$ and
$$d_{\T^d}(u,\y)<\bar{\delta}_\ee(t)\quad\Lra\quad|\pd_tf_t(u)-\pd_tf_t(\y)|<\ee$$
and since $\|\pd f\|_\infty\in L^1(0,T)$ for each $\ell\in\ZZ_+$ there exists $\delta_\ell>0$ such that 
$$\LL_{(0,T)}(E)<\delta_\ell\quad\Lra\quad\int_E\|\pd_tf_t\|_\infty\df t<\ell^{-1}.$$ 
Since $\LL_{(0,T)}(\{\bar{\delta}_{\ell^{-1}}=0\})=0$ for all $\ell\in\ZZ_+$, for each $\ell\in\ZZ_+$ there exists $k_\ell\in\NN$ such that $\LL_{(0,T)}(\{\bar{\delta}_{\ell^{-1}}<\fr{k_\ell}\})<\delta_\ell$. Then for all $(N,\ell)\in\NN\x\ZZ_+$
\begin{align}\label{OnSmallSet1}
\Big|\int_{\{\bar{\delta}_{\ell^{-1}}<\fr{k_\ell}\}}\Big\ls\fr{\ell_\star^d}\sum_{|y|\leq\ell}\pd_tf_t\Big(\frac{\cdot+y}{N}\Big)-\pd_tf_t,\pi^N_t\Big\rs\df t\Big|&\leq 2\int_{\{\bar{\delta}_{\ell^{-1}}<k_\ell^{-1}\}}\|\pd_tf_t\|_\infty\ls 1,\pi^N_t\rs\df t\nonumber\\
&\stackrel{P^N\mbox{-a.s.}}\leq 2\ls\pi_0^N,1\rs\int_{\{\bar{\delta}_\ell<k_\ell^{-1}\}}\|\pd_tf_t\|_\infty\df t\nonumber\\
&\leq2\ell^{-1}\ls\pi_0^N,1\rs\LL_{(0,T)}(\{\bar{\delta}_{\ell^{-1}}<k_\ell^{-1}\}).
\end{align}
On the other hand, for each $\ell\in\ZZ_+$ we can choose $N_\ell\in\NN$ such that $\ell/N<k_\ell^{-1}$ for all $N\geq N_\ell$ and then for all $\ell\in\ZZ_+$, $N\geq N_\ell$ 
\begin{align}\label{OnNotSmallSet1}\Big|\int_{\{\bar{\delta}_{\fr{\ell}}\geq\fr{k_\ell}\}}\Big\ls\fr{\ell_\star^d}\sum_{|y|\leq\ell}\pd_tf_t\Big(\frac{\cdot+y}{N}\Big)-\pd_tf_t,\pi^N_t\Big\rs\df t\Big|&\leq\int_{\{\bar{\delta}_{\fr{\ell}}\geq k_\ell^{-1}\}}\ell^{-1}\ls 1,\pi^N_t\rs\df t\nonumber\\
&\stackrel{P^N\mbox{-a.s.}}= \ell^{-1}\ls\pi_0^N,1\rs\LL_{(0,T)}(\{\bar{\delta}_{\fr{\ell}}\geq k_\ell^{-1}\}).
\end{align}It follows by~\eqref{Init11},~\eqref{OnSmallSet1} and~\eqref{OnNotSmallSet1} that for each $\ell\in\ZZ_+$ and all $N\geq N_\ell$
$$\lls\pd f,\pi^{N,\ell}-\pi^N\rrs\stackrel{P^N\mbox{-a.s.}}\leq 2 \ell^{-1}\ls\pi_0^N,1\rs$$
and therefore 
\begin{align*}
\limsup_{\ell,N\ra+\infty}P^N\big\{\big|\lls\pd f,\pi^{N,\ell}-\pi^N\rrs\big|>\ee\big\}&\leq\limsup_{\ell,N\ra+\infty}P^N\{2 \ell^{-1}\ls\pi_0^N,1\rs>\ee\}\\
&\leq\limsup_{\ell,N\ra+\infty}\mu_0^N\{2\ls\pi_0^N,1\rs>\ell\ee\}=0
\end{align*}	 
where the last limit follows by Lemma~\ref{BoundedTotParticlInMacrLim}.

We have thus shown that $\bs{Q}(A_{f,\ee})=0$ for all $C_c^3((0,T)\x\T^d)$. Since this holds for all $\ee>0$ it follows that 
$$\bs{Q}\big(\big\{\lls\Lambda\pd f(U)+\Phi(\Lambda)\Delta f(U),\bs\pi\rrs=0\}\big)=1,\quad\forall f\in C_c^3((0,T)\x\T^d).$$
Since there exists a countable family $\mathcal{G}\subs C_c^3((0,T)\x\T^d)$ that is dense in $C_c^2((0,T)\x\T^d)$ in the $C^2$-uniform norm it follows that 
$$\bs{Q}\Big(\bigcap_{f\in C_c^2((0,T)\x\T^d)}\big\{\lls\Lambda\pd f(U)+\Phi(\Lambda)\Delta f(U),\bs\pi\rrs=0\}\Big)=1.$$

We have thus shown that $\bs{Q}$ is concentrated on mild generalized Young measure-valued solutions to the hydrodynamic equation $\pd_t\bs\pi=\Delta\Phi(\bs\pi)$ in the sense of~\eqref{Mild}. By Proposition~\ref{Leb} we also know that $\bs{Q}$ is concentrated on $L_{w^*}^\infty(0,T;\bbar{\Y}_1(\T^d))$. Therefore if we set 
$$A:=\bigcap_{f\in C_c^2((0,T)\x\T^d)}\big\{\bs\pi\in L_{w^*}^\infty(0,T;\bbar{\Y}_1(\T^d))\bigm|\lls\Lambda\pd f(U)+\Phi(\Lambda)\Delta f(U),\bs\pi\rrs=0\big\}$$
we have that $\bs{Q}(A)=1$. Since $\Phi$ is sublinear, i.e.~$\lim_{\lambda\ra+\infty}\Phi(\lambda)/\lambda=0$, for any $\bs\pi\in A$ and any $f\in L^1(0,T;C(\T^d))$ $$\lls \Phi(\Lambda) f(U),\bs\pi\rrs=\int_0^T\ls\Phi(\Lambda)f_t(U),\bs\rho_{\bs\pi_t}\rs\df t=\int_0^T\int_{\T^d}f_t(u)\Phi(\bs\rho_{\bs\pi_t})(u)\df u\df t$$ 
where for any $\bs\rho\in\PP_1(\T^d\x\RR_+)$ the composition $\Phi(\bs\rho)\equiv b_\Phi(\bs\rho)\colon\T^d\to\RR$ is Lebesgue a.s.~defined as in~\eqref{YoungCompose}. Thus we can express the fact that  $\bs\pi\in A$ is a mild solution to the hydrodynamic equation as
\begin{equation}\label{MildForLeb}
\int_0^T\ls \pd_tf_t,\bs\pi_t\rs\df t=-\int_0^T\int_{\T^d}\Delta f_t(u)\Phi(\bs\rho_{\bs\pi_t})(u)\df u\df t,\quad\forall f\in C_c^2((0,T)\x\T^d).
\end{equation}

We will show now that $\bs{Q}$ is in fact concentrated on weak solutions in the sense of~\eqref{YMVZRPPDE} that satisfy the energy estimate~\eqref{PhiEnergEst}. For this we consider the joint laws 
$$\bs{Q}_{\mathfrak{g}(\eta(0))}^{N,\ell}:=(\s^N,\bs\pi^{N,\ell})_\sharp P^N\in L_{w^*}^\infty(0,T;\MMM_+(\T^d))\x L_{w^*}^\infty(0,T;\bbar{\MMM}_{1,+}(\T^d\x\RR_+))$$ 
and the $\Phi$-projection $B_\Phi\colon L_{w^*}^\infty(0,T;\bbar{\MMM}_{1,+}(\T^d\x\RR_+))\to L_{w^*}^\infty(0,T;\MMM_+(\T^d))$.
By the one-block estimate we know that for any limit point $\bs{Q}_{\mathfrak{g}(\eta(0))}$ of $\{\bs{Q}_{\mathfrak{g}(\eta(0))}^{N,\ell}\}$ 
\begin{equation}\label{JRequalPhiRho}
\bs{Q}_{\mathfrak{g}(\eta(0))}\big\{(\s,\bs\pi)\in L_{w^*}^\infty(0,T;\MMM_+(\T^d))\x L_{w^*}^\infty(0,T;\bbar{\MMM}_{1,+}(\T^d\x\RR_+))\bigm|\s=B_\Phi(\bs\pi)\big\}=1.
\end{equation}
Since $\Phi$ is sublinear
\begin{equation}\label{Equals}
B_\Phi(\bs\pi)=\Phi(\bs\rho_{\bs\pi})\df\LL_{\T^d}\quad\mbox{in }L_{w^*}^\infty(0,T;\MMM_+(\T^d)).
\end{equation}
Thus by Theorem~\ref{Theorem1}(e) and Theorem~\ref{TheoremEnergyEstimate} it follows that for $\bs{Q}$-a.s.~all $\bs\pi\in L_{w^*}^\infty(0,T;\bbar{\MMM}_{1,+}(\T^d\x\RR_+))$ it holds that $\Phi(\bs\rho_{\bs\pi_t})\in H^1(\T^d)$ for almost all $t\in[0,T]$ and 
$$\int_0^T\int_{\T^d}\frac{\|\nabla\Phi(\bs\rho_{\bs\pi_t})(u)\|^2}{\Phi(\bs\rho_{\bs\pi_t})(u)}\df u\df t<+\infty.$$ This proves the regularity $\Phi(\bs\rho_{\bs\pi_t})\in H^1(\T^d)$ for almost all $t\in[0,T]$ and the energy estimate~\eqref{PhiEnergEst}. By a standard mollification argument in the space variable and~\eqref{MildForLeb} that 
$$\int_0^T\ls \pd_tf_t,\bs\pi_t\rs\df t=\int_0^T\int_{\T^d}\ls\nabla f_t(u),\nabla\Phi(\bs\rho_{\bs\pi_t})(u)\rs\df u\df t,\quad\forall f\in C_c^1((0,T)\x\T^d).$$
Therefore any $\bs\pi\in A$ is a weak solution in the sense of~\eqref{YMVZRPPDE}.

 \subsection{Two-blocks Comparison}\label{TBCSection}
(a) In order to simplify the notation we will work with the quantity $\EE^N\big|\lls G,B_\Psi(\bs\pi^{N,\ell})-B_\Psi(\bs\pi^{N,\ell;M;\ee})\rrs|$ for the full set of parameters $N\in\NN$, $\ell\in\ZZ_+$, $\ee,M>0$ and restrict our attention to the subfamily of parameters $(k_N^{(\ell)},m_\ell,M,\ee)_{N,\ell,\ee,M}$ along which $\bs{Q}^{k_N^{(\ell)},m_\ell}:=\bs\pi^{k_N^{(\ell)},m_\ell}_\sharp P^{k_N^{(\ell)}}$ converges as $N\ra+\infty$ and then $\ell\ra+\infty$ to $\bs{Q}^{\infty,\infty}$ only when necessary.
	
Let $\Psi\in\bbar{C}_1(\RR_+)$. By Proposition~\ref{FluidSolidSeparation} for all $G\in L^1(0,T;C(\T^d))$ 
$$\lim_{M\ra+\infty}\limsup_{\ell,N\ra+\infty}\EE^N|\lls G,B_{\Psi}(\bs\pi^{N,\ell})-B_{\Psi}(\bs\pi^{N,\ell;M})\rrs|=0$$ and thus in order to prove the truncated double-block estimate~\eqref{DBE} it suffices to show that 
\begin{equation}\label{TDBE}
\lim_{M\ra+\infty}\limsup_{\ell\ra+\infty}\limsup_{\ee\ra 0}\limsup_{N\ra+\infty}\EE^{k_N^{(\ell)}}\big|\lls G,B_{\Psi}(\bs\pi^{k_N^{(\ell)},m_\ell;M})-B_{\Psi}(\bs\pi^{k_N^{(\ell)},m_\ell;M;\ee})\rrs\big|=0
\end{equation}
for all $\delta>0$ and $G\in L^1(0,T;C(\T^d))$.

We will reduce first the proof of~\eqref{TDBE} to the case where $\Psi$ is sublinear, i.e.~$\Psi\in C_1(\RR_+)$. Indeed, since $\Psi\in\bbar{C}_1(\RR_+)$ the limit $\Psi'(\infty)=\lim_{\lambda\ra\infty}\frac{\Psi(\lambda)}{\lambda}$ exists and the map $\Psi_0(\lambda)=\Psi(\lambda)-\Psi'(\infty)\lambda$ belongs in $C_1(\RR_+)$ and $\Psi(\lambda)=\Psi_0(\lambda)+\Psi'(\infty)\lambda$, $\lambda\geq 0$. Therefore
\begin{align*}
\ls G,B_{\Psi}(\bs\pi)\rs&=\ls G(U)\Psi(\Lambda),\bs\pi\rs=\ls G(U)\Psi_0(\Lambda),\bs\pi\rs+\Psi'(\infty)\ls G(U)\Lambda,\bs\pi\rs\\
&=\ls G,B_{\Psi_0}(\bs\pi)\rs+\Psi'(\infty)\ls G,B(\bs\pi)\rs.
\end{align*} 
Consequently for all $N\in\NN$, $\ell\in\ZZ_+$, $\ee,M>0$
\begin{align*}P^N\big\{\big|\lls G,B_{\Psi}(\bs\pi^{N,\ell;M})-B_{\Psi}(\bs\pi^{N,\ell;M;\ee})\rrs\big|>\delta\big\}&\leq
P^N\Big\{\big|\lls G,B_{\Psi_0}(\bs\pi^{N,\ell;M})-B_{\Psi_0}(\bs\pi^{N,\ell;M;\ee})\rrs\big|>\frac{\delta}{2}\Big\}\\
&+P^N\Big\{\big|\lls G,B(\bs\pi^{N,\ell;M})-B(\bs\pi^{N,\ell;M;\ee})\rrs\big|>\frac{\delta}{2\Psi'(\infty)}\Big\}.
\end{align*}
Thus if show that for any $\delta>0$ and any $G\in L^1(0,T;C(\T^d))$
\begin{equation}\label{TDBEL}
\lim_{M\ra+\infty}\limsup_{\ell\ra+\infty}\limsup_{\ee\ra 0}\limsup_{N\ra+\infty}P^N\big\{\big|\lls G,B(\bs\pi^{N,\ell;M})-B(\bs\pi^{N,\ell;M;\ee})\rrs\big|>\delta\big\}=0
\end{equation}
we reduce the proof of~\eqref{TDBE} to the case that $\Psi$ is sublinear. But this is elementary since by a change of variables 
\begin{align*}
\lls G,B(\bs\pi^{N,\ell;M})-B(\bs\pi^{N,\ell;M;\ee})\rrs=\int_0^T\fr{N^d}\sum_{x\in\T_N^d}\fr{[N\ee]_\star^d}\sum_{|y|\leq[N\ee]}\Big\{G_t\Big(\frac{x}{N}\Big)-G_t\Big(\frac{x-y}{N}\Big)\Big\}\eta_t^\ell(x)\df t
\end{align*}
and therefore by the conservation of the total number of particles $P^N$-a.s.~in $D(0,T;\MM_N^d)$
\begin{equation*}
|\lls G,B(\bs\pi^{N,\ell;M})-B(\bs\pi^{N,\ell;M;\ee})\rrs|\leq\ls 1,\pi_0^N\rs\int_0^T\w_{G_t}(2\ee)\df t,\end{equation*} where for any function $G\colon\T^d\lra\RR$ we denote by 
$$\w_G(\ee):=\sup_{\substack{u,\y\in\T^d\\|u-\y|<\ee}}|G(u)-G(\y)|$$ its modulus of continuity. Consequently 
$$P^N\Big\{\big|\lls G,B(\bs\pi^{N,\ell;M})-B(\bs\pi^{N,\ell;M;\ee})\rrs\big|\big\}
\leq\mu_0^N\Big\{\ls 1,\pi^N\rs\leq\frac{\delta}{\int_0^T\w_{G_t}(2\ee)\df t}\Big\}.$$
Since $G\in L^1(0,T;C(\T^d))$ it follows that $\lim_{\ee\ra 0}\int_0^T\w_{H_t}(2\ee)dt=0$ and therefore~\eqref{TDBEL} follows from Lemma~\ref{BoundedTotParticlInMacrLim}.

We will show next that in the proof of~\eqref{TDBE} we can further assume that the map $\Psi\in C_1(\RR_+)$ is Lipschitz. To prove this, we define for each $k\in\NN$ the Moreau-Yosida approximations $\Psi_k\colon\RR_+\to\RR$ of $\Psi$ by
\begin{equation}\label{MoreauYosida}
\Psi_k(\rho)=\inf_{\lambda\geq 0}\big\{\Psi(\lambda)+k|\lambda-\rho|\big\}.
\end{equation}
Then as is well known, each map $\Psi_k$ is Lipschitz with Lipschitz constant $\Lip(\Psi_k)\leq k$ and $\{\Psi_k\}_{k=1}^\infty$ increases pointwise to $\Psi$ pointwise. It is also easy to see that $\Psi_k$ is sublinear and $\|\Psi_k\|_{C_1(\RR_+)}\leq\|\Psi\|_{C_1(\RR_+)}$ for large $k$. Indeed, since $\Psi$ is sublinear there exists for each $\delta\in(0,1)$ a constant $C_\delta<+\infty$ such that $|\Psi(\lambda)|\leq C_\delta+\delta\lambda$ for all $\lambda\geq 0$ and therefore 
$$-C_\delta+\inf_{\lambda\geq 0}\big\{-\delta\lambda+k|\lambda-\rho|\big\}\leq\Psi_k(\rho)\leq C_\delta+\inf_{\lambda\geq 0}\big\{\delta\lambda+k|\lambda-\rho|\big\}.$$
It is elementary to check that the infima above are both obtained at $\lambda=\rho$ for all $k\in\NN$ and therefore $|\Psi_k(\rho)|\leq C_\delta+\delta\rho$ for all $\rho\geq 0$ and all $\delta\in(0,1)$, which proves that $\Psi_k$ is sublinear. Similarly one can check that $\|\Psi_k\|_{C_1(\RR_+)}\leq\|\Psi\|_{C_1(\RR_+)}$ for all $k>\|\Psi\|_{C_1(\RR_+)}$. Now, assuming that~\eqref{TDBE} holds for sublinear Lipschitz maps, it follows that for each $k\in\NN$ 
$$\lim_{M\ra+\infty}\limsup_{\ell\ra+\infty}\limsup_{\ee\ra 0}\limsup_{N\ra+\infty}\EE^{K_N^{(\ell)}}|\lls G,B_{\Psi_k}(\bs\pi^{k_N^{(\ell)},m_\ell;M}-\bs\pi^{k_N^{(\ell)},m_\ell;M;\ee})\rrs\big|=0$$
for all $G\in L^1(0,T;C(\T^d))$. In terms of the family of laws 
\begin{equation}\label{RLFamily1}
\bbar{\bs{Q}}_*^{N,\ell;M;\ee}:=(\bs{\pi}^{k_N^{(\ell)},m_\ell;M},\bs{\pi}^{k_N^{(\ell)},m_\ell;M;\ee})_\sharp P^N\in \PP\big(L_{w^*}^\infty(0,T;\bbar{\PP}_1(\T^d\x\RR_+))^2\big)
\end{equation}
this limit can be written as 
$$\lim_{M\ra+\infty}\limsup_{\ell\ra+\infty}\limsup_{\ee\ra 0}\limsup_{N\ra+\infty}\int|\lls G,B_{\Psi_k}(\bs\pi^\infty-\bs\pi^0)\rrs\big|\df\bbar{\bs{Q}}^{N,\ell;M;\ee}(\bs\pi^\infty,\bs\pi^0)=0.$$
Therefore if 
$$\bbar{\bs{\mathcal{Q}}}_*^{\infty,\infty;\infty;0}:=\Lim_{M\ra+\infty}\Lim_{\ell\ra+\infty}\Lim_{\ee\ra 0}\Lim_{N\ra+\infty}\bbar{\bs{Q}}_*^{N,\ell;M;\ee}$$
it follows by the discussion on subsequential limit sets in Section~\ref{SubLims} that 
\begin{equation}\label{k}
\max_{\bbar{\bs{Q}}\in\bbar{\bs{\mathcal{Q}}}_*^{\infty,\infty;\infty;0}}\int|\lls G,B_{\Psi_k}(\bs\pi^\infty-\bs\pi^0)\rrs\big|\df\bbar{\bs{Q}}(\bs\pi^1,\bs\pi^2)=0,\quad\forall k\in\ZZ_+
\end{equation} and we have to show that 
\begin{equation}\label{ForLip}\max_{\bbar{\bs{Q}}\in\bbar{\bs{\mathcal{Q}}}_*^{\infty,\infty;\infty;0}}\int|\lls G,B_{\Psi}(\bs\pi^\infty-\bs\pi^0)\rrs\big|\df\bbar{\bs{Q}}(\bs\pi^\infty,\bs\pi^0)=0.\end{equation} 
Now, if $\bbar{\bs{Q}}_{\Psi}\in\bbar{\bs{\mathcal{Q}}}_*^{\infty,\infty;\infty;0}$ is a maximizer in the maximum above, we have by~\eqref{k} that 
\begin{equation}\label{ToTakeLimForLip}
\int|\lls G,B_{\Psi_k}(\bs\pi^\infty-\bs\pi^0)\rrs\big|\df\bbar{\bs{Q}}_{\Psi}(\bs\pi^\infty,\bs\pi^0)=0,\quad\forall\in\ZZ_+.
\end{equation} 
Thus the claim follows by the dominated convergence, since as we will see the sequence $\{\mathcal{I}^G_k\}_{k\in\NN}$ of the functionals defined on $L_{w^*}^\infty(0,T;\bbar{\MMM}_{1,+}(\T^d\x\RR_+))^2$ by 
$$\mathcal{I}^G_k(\bs\pi^\infty,\bs\pi^0)=|\lls G,B_{\Psi_k}(\bs\pi^\infty-\bs\pi^0)\rrs|$$ converges $\bbar{\bs{Q}}_{\Psi}$-a.s.~pointwise to the functional $\mathcal{I}^G\colon L_{w^*}^\infty(0,T;\bbar{\MMM}_{1,+}(\T^d\x\RR_+))^2\to\RR_+$ and it is dominated by an $L^1(\bbar{\bs{Q}}_{\Psi})$-function. First, since $B_{\Psi_k}$ is linear, in order to prove this pointwise convergence it suffices to show that 
$$\lim_{k\ra+\infty}\lls G,B_{\Psi_k}(\bs\pi)\rrs=\lls G,B_{\Psi}(\bs\pi)\rrs$$ 
$\bs{Q}_{\Psi,1}$-a.s.~and $\bs{Q}_{\Psi,2}$-a.s.~for all $\bs\pi\in L_{w^*}^\infty(0,T;\bbar{\MMM}_{1,+}(\T^d\x\RR_+))$, where $\bs{Q}_{\Psi,i}$ is the $i$-th marginal of $\bbar{\bs{Q}}_{\Psi}$, $i=1,2$. Both marginals are supported on the set $L_{w^*}^\infty(0,T;\bbar{\Y}_{1,\leq\mathfrak{m}}(\T^d))$ of all generalized Young measures $\bs\pi\in L_{w^*}^\infty(0,T;\bbar{\Y}_1(\T^d))$ such that $\ls\Lambda,\bs\pi_t\rs\leq\mathfrak{m}$ for almost all $t\in[0,T]$ and thus in proving this limit we can assume that $\bs\pi\in L_{w^*}^\infty(0,T;\bbar{\Y}_{1,\leq\mathfrak{m}}(\T^d))$. Since $\{\Psi_k\}$ converges pointwise to $\Psi$ it obviously follows that $G(U)\Psi_k(\Lambda)$ converges pointwise to $G(U)\Psi(\Lambda)$ as $k\ra+\infty$ and since $\|\Psi_k\|_{\infty;1}\leq\|\Psi\|_{\infty;1}$ for large $k$, we have that $|G(U)\Psi_k(\Lambda)|\leq\|\Psi\|_{\infty;1}|G(U)|(1+\Lambda)$ for large $k$ and thus 
\begin{align}\label{ForDomin}
|G_t(U)\Psi_k(\Lambda)|\leq\|\Psi\|_{\infty;1}\ls|G_t(U)|(1+\Lambda),\bs\pi_t\rs\leq\|\Psi\|_{\infty;1}(1+\mathfrak{m})\|G_t\|_\infty
\end{align} for almost all $t\in[0,T]$. Therefore, since the maps $\Psi_k$, $\Psi$ are sublinear, by the bounded convergence theorem 
\begin{align}\label{ForDomin1}\lim_{k\ra+\infty}\ls G_t,B_{\Psi_k}(\bs\pi_t)\rs&=\lim_{k\ra+\infty}\ls G_t(U)\Psi_k(\lambda),\bs\pi_t\rs=\lim_{k\ra+\infty}\ls G_t(U)\Psi_k(\lambda),\bs\rho_{\bs\pi_t}\rs\nonumber\\
&=\ls G_t(U)\Psi(\lambda),\bs\rho_{\bs\pi_t}\rs=\ls G_t(U)\Psi(\lambda),\bs\pi_t\rs=\ls G_t,B_{\Psi}(\bs\pi_t)\rs\end{align} 
for almost all $t\in[0,T]$, where for each $\bs\pi\in\bbar{\MMM}_1(\T^d\x\RR_+))$ we denote by $\bs\rho_{\bs\pi}=j^*\bs\pi\in\Y_1(\T^d)$ the ordinary Young measure representing the regular part $\widehat{\bs\pi}$ of $\bs\pi$. Now, by~\eqref{ForDomin} the sequence $\{\ls G_\cdot,B_{\Psi_k}(\bs\pi_\cdot)\rs\}_{k\in\NN}$ is dominated for large $k\in\NN$ by the $L^1(0,T)$-function $\|\Psi\|_{\infty;1}(1+\mathfrak{m})\|G_\cdot\|_\infty$ and thus by~\eqref{ForDomin1} and the dominated convergence theorem $\lim_{k\ra+\infty}\lls G,B_{\Psi_k}(\bs\pi)\rrs=\lls G,B_{\Psi}(\bs\pi)\rrs$. This proves that $\mathcal{I}_k$ converges $\bbar{\bs{Q}}$-almost surely to $\mathcal{I}$ and thus in order to use the dominated convergence theorem in~\eqref{ToTakeLimForLip} to obtain~\eqref{ForLip} it remains to check that $\{\mathcal{I}_k\}$ is dominated by an $L^\infty(\bbar{\bs{Q}}_{\Psi})$-integrable function. But this is easy, since for all $\bs\pi\in L_{w^*}^\infty(0,T;\bbar{\Y}_{1,\leq\mathfrak{m}}(\T^d\x\RR_+))$ and all $G\in L^1(0,T;C(\T^d))$
\begin{align*}\mathcal{I}_k(\bs\pi^\infty,\bs\pi^0)&\leq|\lls G(U)\Psi_k(\Lambda),\bs\pi^\infty\rrs|+|\lls G(U)\Psi_k(\Lambda),\bs\pi^0\rrs|\\
&\leq\lls|G(U)\Psi_k(\Lambda)|,\bs\pi^\infty\rrs+\lls|G(U)\Psi_k(\Lambda)|,\bs\pi^0\rrs\\
&\leq 2\|\Psi\|_{\infty;1}\|G\|_{L^1(0,T;C(\T^d))}+\|\Psi\|_{\infty;1}\lls|G|,B(\bs\pi^\infty+\bs\pi^0)\rrs\\
&\leq 2\|\Psi\|_{\infty;1}\|G\|_{L^1(0,T;C(\T^d))}(1+\mathfrak{m})
\end{align*} 
where the last inequality holds $\bbar{\bs{Q}}_{\Psi}$-a.s.

So in what follows we will assume that $\Psi$ is a sublinear and Lipschitz cylinder map and we will prove~\eqref{TDBE}. We recall that the fact that $\Psi$ is sublinear implies that $B_{\Psi}=B_{\Psi}\circ\widehat{D}$. By considering the empirical process $\psi^{N,\ell;M;\ee}\colon D(0,T;\MM_N^d)\to L_{w^*}^\infty(0,T;\MMM(\T^d))$ defined by 
$$\lls G,\psi^{N,\ell;M;\ee}\rrs:=\int_0^T\fr{N^d}\sum_{x\in\T_N^d}\fr{[N\ee]_\star^d}\sum_{|z|\leq[N\ee]}G_t\Big(\frac{x}{N}\Big)\Psi\big(\eta_t^\ell(x+z)\mn M\big)\df t$$ the truncated double-block estimate~\eqref{TDBE} is split in proving the limits 
\begin{align}\label{TDBELim1}
\lim_{M\ra+\infty}\limsup_{\ell\ra+\infty}\limsup_{\ee\ra 0}\limsup_{N\ra+\infty}P^N\big\{\big|\lls G,B_{\Psi}(\bs\pi^{N,\ell;M})-\psi^{N,\ee,\ell;M}\rrs\df t\big|>\delta\big\}=0
\end{align}
and
\begin{align}\label{TDBELim2}
\lim_{M\ra+\infty}\limsup_{\ell\ra+\infty}\limsup_{\ee\ra 0}\limsup_{N\ra+\infty}P^N\big\{\big|\lls G,\psi^{k_N^{(\ell)},m_\ell;M;\ee}-B_{\Psi}(\bs\pi^{k_N^{(\ell)},m_\ell;M;\ee})\rrs\df t\big|>\delta\big\}=0
\end{align}
for all $G\in L^1(0,T;C(\T^d))$ and all $\delta>0$. For the quantity in~\eqref{TDBELim1}, by a change of variables we have that 
\begin{align*}
\lls G,B_{\Psi}(\bs\pi^{N,\ell;M})-\psi^{N,\ell;M;\ee}\rrs=\fr{N^d}\sum_{x\in\T_N^d}\fr{[N\ee]_\star^d}\sum_{|y|\leq[N\ee]}\Big\{G_t\Big(\frac{x}{N}\Big)-G_t\Big(\frac{x-y}{N}\Big)\Big\}\Psi\big(\eta_t^\ell(x)\mn M\big)
\end{align*}
and therefore $P^N$-a.s.~in $D(0,T;\MM_N^d)$ 
\begin{equation*}
|\lls G,B_{\Psi}(\bs\pi^{N,\ell;M})-\psi^{N,\ee,\ell;M}\rrs|\leq\|\Psi\|_{\infty;1}(1+\ls 1,\pi_0^N\rs)\int_0^T\w_{H_t}(2\ee)\df t
\end{equation*} where $\|\Psi\|_{\infty;1}:=\sup_{\lambda\geq 0}\frac{|\Psi(\lambda)|}{1+\lambda}<+\infty$. Thus the limit~\eqref{TDBELim1} is shown to vanish similarly to~\eqref{TDBEL}.

We prove next the limit~\eqref{TDBELim2}. By Chebyshev inequality for this term it suffices to show that 
$$\limsup_{M\uparrow\infty,\ell\uparrow\infty,\ee\downarrow 0,N\uparrow\infty}\EE^{k_N^{(\ell)}}\int_0^T|\ls G_t,\psi_t^{k_N^{(\ell)},m_\ell;M;\ee}-B_{\Psi}(\bs\pi_t^{k_N^{(\ell)},m_\ell;M;\ee})\rs|dt=0.$$
For all parameters $(N,\ell,\ee,M)\in\NN\x\ZZ_+\x(0,\infty)^2$ and each $t\in[0,T]$
\begin{equation*}\ls G_t,\psi_t^{N,\ee,\ell;M}-B_{\Psi}(\bs\pi_t^{N,\ell;M;\ee})\rs=\fr{N^d}\sum_{x\in\T_N^d}G_t\Big(\frac{x}{N}\Big)\big\{\Psi\big(\eta_t^\ell(x)\mn M\big)^{[N\ee]}-\Psi\big((\eta^\ell(x)\mn M)^{[N\ee]}\big)\big\}.\end{equation*}
Since the macroscopic averages appear inside that non-linear map $\Psi$ this term can not be dealt by an integration by parts. Since $\Psi$ is assumed Lipschitz we can estimate the absolute value of this term by 
\begin{align*}
|\ls G_t,\psi_t^{N,\ee,\ell;M}-B_{\Psi}(\bs\pi_t^{N,\ell;M;\ee})\rs|&\leq\frac{\|\Psi\|_\Lip\|G_t\|_\infty}{N^d[N\ee]_\star^d}\sum_{x\in\T_N^d}\sum_{|y|\leq[N\ee]}\big|\eta_t^\ell(x+y)\mn M-(\eta_t^\ell(x)\mn M)^{[N\ee]}\big|\\
&\leq\frac{\|\Psi\|_\Lip\|G_t\|_\infty}{N^d[N\ee]_\star^{2d}}\sum_{\substack{x\in\T_N^d\\|y|\mx|z|\leq[N\ee]}}\big|\eta_t^\ell(x+y)\mn M-\eta_t^\ell(x+z)\mn M\big|\\
&=\frac{\|\Psi\|_\Lip\|G_t\|_\infty}{N^d[N\ee]_\star^{2d}}\sum_{x\in\T_N^d}\sum_{\substack{|y|\mx|z|\leq[N\ee]\\2\ell<|y-z|}}
\Psi_M\big(\eta_t^\ell(x+y),\eta_t^\ell(x+z)\big)\\
&\quad+\frac{\|\Psi\|_\Lip\|G_t\|_\infty}{N^d[N\ee]_\star^{2d}}\sum_{x\in\T_N^d}\sum_{\substack{|y|\mx|z|\leq[N\ee]\\|y-z|\leq 2\ell}}\Psi_M\big(\eta_t^\ell(x+y),\eta_t^\ell(x+z)\big),
\end{align*}
 where $\Psi_M\colon\RR_+^2\to\RR$, $M>0$, is the map $\Psi_{M}(a,b):=|a\mn M-b\mn M|$, $a,b\in\RR_+$. The last term above is bounded above by 
$$\frac{\|\Psi\|_\Lip\|G_t\|_\infty}{N^d[N\ee]_\star^{2d}}\sum_{x\in\T_N^d}\sum_{\substack{|z|\leq[N\ee]\\|y-z|\leq 2\ell}}\Psi_M\big(\eta_t^\ell(x+y),\eta_t^\ell(x+z)\big)\leq M\|\Psi\|_\Lip\|G_t\|_\infty\frac{(2\ell)_\star^d}{[N\ee]_\star^{d}},$$
and the time integral of this last term vanishes as $N\ra+\infty$ since $G\in L^1(0,T;C(\T^d))$ for each $\ell\in\ZZ_+$ and $\ee,M>0$. Thus it suffices to show that 
\begin{align}\label{TBETr}
\limsup_{M\uparrow,\ell\uparrow\infty,\ee\downarrow 0,N\uparrow\infty}\EE^{k_N^{(\ell)}}\int_0^T\fr{(k_N^{(\ell)})^d[k_N^{(\ell)}\ee]_\star^{2d}}\sum_x\sum_{\substack{|y|\mx|z|\leq[k_N^{(\ell)}\ee]\\2m_\ell<|y-z|}}\Psi_{M}\big(\eta_t^{m_\ell}(x+y),\eta_t^{m_\ell}(x+z)\big)\df t=0,
\end{align}
where the sums are taken among all $x,y,z\in\T_{k_N^{(\ell)}}^d$.

The proof of~\eqref{TBETr} is similar to the proof of the two-blocks estimate in~\cite{Kipnis1999a}. The main new element in the proof is that the introduction of the truncating parameter $M>0$ allows us to cut off the large densities for the fluid phase as described in the following lemma. It is in this lemma that we need to restrict the limit superior
	 \begin{lemma}\label{TwoBlockCutLargeDens}{\rm{(Cutting off large densities for the fluid-phase)}} Suppose that the ZRP starts from a sequence of initial profiles $\mu_0^N\in\PP_1\MM_N^d$, $N\in\NN$, with total mass $\mathfrak{m}>0$ in probability and let $\{m_\ell\}_{\ell=1}^\infty$, $\{k_N^{(\ell)}\}_{N=1}^\infty$ be sequences such that $\bs{Q}^{k_N^{(\ell)},m_\ell}:=\bs\pi^{k_N^{(\ell)},m_\ell}_\sharp P^{k_N^{(\ell)}}$ converges to some $\bs{Q}_*^{\infty,\infty}\in\PP L_{w^*}^\infty(0,T;\bbar{\Y}_1(\T^d))$ as $N\ra+\infty$ and then $\ell\ra+\infty$. Then for any $T>0$
	 	$$\lim_{A\ra+\infty}\sup_{M>0}\limsup_{\ell,N\ra+\infty}\EE^{k_N^{(\ell)}}\int_0^T\fr{(k_N^{(\ell)})^d}\sum_{x\in\T_{k_N^{(\ell)}}^d}(\eta_t^{m_\ell}(x)\mn M)\1_{\{\eta_t^{m_\ell}(x)>A\}}\df t=0.$$
	 \end{lemma}
	 \textbf{Proof} The expectation in the conclusion of the lemma is an increasing function of $M>0$ and thus the supremum over $M>0$ is equal to the limit as $M\ra+\infty$. Thus since $A$ tends to infinity after $M$ has been sent to infinity, we can always assume that $M>A$. Now, for $M>A$ and all $\lambda\geq 0$ it obviously holds that 
	 $$(\lambda\mn M)\1_{(A,+\infty)}(\lambda)=(\lambda\mn M)\1_{(A,+\infty)}(\lambda\mn M)$$
	 and for the continuous map $\Psi_A(\lambda)=\lambda\cdot[(\lambda-A+1)^+\mn 1]$, $\lambda\geq 0$, we have that 
	 $$0\leq\lambda\1_{(A,+\infty)}(\lambda)\leq\Psi_A(\lambda)\leq\lambda,\quad\forall\lambda\geq 0.$$
	 Consequently in order to prove the lemma it suffices to show that 
	 $$\lim_{A\ra+\infty}\limsup_{M,\ell,N\ra+\infty}\EE^{k_N^{(\ell)}}\int_0^T\fr{(k_N^{(\ell)})^d}\sum_{x\in\T_{k_N^{(\ell)}}^d}\Psi_A(\eta_t^{m_\ell}(x)\mn M)\df t=0.$$
	 In terms of the regular part $\widehat{\bs\pi}^{N,\ell;M}$ of the $M$-modified micro-empirical density and equation~\eqref{Trunc} the expected value above can be written for all parameters $(N,\ell)$ as 
	 \begin{align*}
	 \EE^N\int_0^T\fr{N^d}\sum_{x\in\T_N^d}\Psi_A(\eta_t^\ell(x)\mn M)\df t&=\EE^N\lls\Psi_A(\Lambda),\widehat{\bs\pi}^{N,\ell;M}\rrs
	 =\EE^N\big\lls\Psi_A(\Lambda),\Pi_M^*\circ j^*\circ\bs{\pi}^{N,\ell}\big\rrs\\
	 &=\EE^N\big\lls j\circ\Pi_M\big(\Psi_A(\Lambda)\big),\bs{\pi}^{N,\ell}\big\rrs=\int\big\lls\Psi_A(\Lambda\mn M),\bs{\pi}\big\rrs\df\bs{Q}^{N,\ell}(\bs{\pi}).
	 \end{align*}
	 Since $\bs{Q}_*^{\infty,\infty}:=\lim_{\ell\ra+\infty}\lim_{N\ra+\infty}\bs{Q}^{k_N^{(\ell)},m_\ell}$ we have that 
	 	 \begin{equation}\label{ToOneInd}
	 \limsup_{\ell,N\ra+\infty}\EE^{k_N^{(\ell)}}\int_0^T\fr{(k_N^{(\ell)})^d}\sum_{x\in\T_{k_N^{(\ell)}}^d}\Psi_A(\eta_t^{m_\ell}(x)\mn M)\df t=\int\lls(j\circ\Pi_M)\Psi_A(\Lambda),\bs{\pi}\rrs\df\bs{Q}_*^{\infty,\infty}(\bs{\pi}).\end{equation}

	 Now, equation~\eqref{ToOneInd} reduces the proof of the lemma to showing that 
	 $$\lim_{A\ra+\infty}\limsup_{M\ra+\infty}\int\lls\Psi_A(\Lambda),\bs{\pi}\rrs\df(\Pi_M^*\circ j^*)_\sharp \bs{Q}^{\infty,\infty}_*(\bs{\pi})=0.$$
	 But in the proof of Corollary~\ref{PathMeasDecomp} we have seen that $\Pi_M^*\circ j^*$ $w^*$-converges pointwise to $\widehat{D}$ in $L_{w^*}^\infty(0,T;\bbar{\MMM}_{1,+}(\T^d\x\RR_+))$ and thus since the linear map $I_A(\cdot)\equiv\lls\Psi_A(\Lambda),\cdot\rrs\colon L_{w^*}^\infty$ is $w^*$-continuous,
	 $$\lim_{M\ra+\infty}\int\lls\Psi_A(\Lambda),\bs{\pi}\rrs\df(\Pi_M^*\circ j^*)_\sharp \bs{Q}_*^{\infty,\infty}(\bs{\pi})=
	 \int\lls\Psi_A(\Lambda),\bs{\pi}\rrs\df\widehat{D}_\sharp \bs{Q}_*^{\infty,\infty}(\bs{\pi}).$$
	Therefore in order to complete the proof of the lemma it suffices to show that 
	 \begin{equation}\label{ToProveLemmaTwoBlockCutLargeDens}
	 \lim_{A\ra+\infty}\int\lls\Psi_A(\Lambda),\bs{\pi}\rrs\df\widehat{D}_\sharp\bs{Q}_*^{\infty,\infty}(\bs{\pi})=0.
	 \end{equation}
	 
	 As we will see,~\eqref{ToProveLemmaTwoBlockCutLargeDens} follows by the dominated convergence theorem, which can be applied due to the fact that $\widehat{D}_\sharp\bs{Q}_*^{\infty,\infty}$ is concentrated on $L_{w^*}^\infty(0,T;\MMM_1(\T^d\x\RR_+))$ seen as a subspace of $L_{w^*}^\infty(0,T;\bbar{\MMM}_1(\T^d\x\RR_+))$. To apply the dominated convergence theorem we check first that the family $\{I_A\}_{A\geq 0}$ is dominated by an $L^1(\widehat{D}_\sharp\bs{Q}_*^{\infty,\infty})$-function and then that $\lim_{A\ra+\infty}I_A=0$ $\widehat{D}_\sharp{\bs{Q}}_*^{\infty,\infty}$-a.s.~pointwise.
	 
	 For the first claim, for all $\bs\pi\in L_{w^*}^\infty\big(0,T;\bbar{\Y}_{1,\mathfrak{m}}(\T^d)\big)$ we have $|I_A(\bs\pi)|=\lls\Psi_A(\Lambda),\bs\pi\rrs\leq\lls\Lambda,\bs\pi\rrs=T\mathfrak{m}$ for all $A>0$. Since $\bs{Q}_*^{\infty,\infty}$ is supported on generalized Young measures with total mass $\mathfrak{m}<+\infty$ and $\lls \Lambda,\widehat{\bs\pi}\rrs\leq\lls\Lambda,\bs\pi\rrs$ for all $\bs\pi\in L_{w^*}^\infty(0,T;\bbar{\MMM}_{1,+}(\T^d\x\RR_+))$ this shows that $\lls\Lambda,\cdot\rrs\in L^\infty(\widehat{D}_\sharp\bs{Q}_*^{\infty,\infty})$.

	 The fact that $\lim_{A\ra+\infty}I_A=0$ $\widehat{D}_\sharp\bs{Q}_*^{\infty,\infty}$-a.s.~pointwise follows by a double application of the dominated convergence theorem. Indeed, the space $$L^\infty_{w^*}(0,T;\Y_1(\T^d))=L_{w^*}^\infty(0,T;\bbar{\Y}_{1}(\T^d)\cap\ker D^\perp$$ is a $w^*$-measurable subspace of $L_{w^*}^\infty(0,T;\bbar{\Y}_1(\T^d))$ and since obviously 
	 $$\widehat{D}\big(L_{w^*}^\infty(0,T;\bbar{\Y}_1(\T^d))\big)\subs L_{w^*}^\infty(0,T;\Y_1(\T^d)),$$ it follows that $\widehat{D}_\sharp\bs{Q}_*^{\infty,\infty}$ is concentrated on the $w^*$-measurable subspace $L^\infty_{w^*}(0,T;\Y_1(\T^d))$. Furthermore since also $\bs{Q}_*^{\infty,\infty}(L_{w^*}^\infty(0,T;\bbar{\Y}_{1,\mathfrak{m}}(\T^d)))=1$ and $\widehat{D}(\bs\pi)\leq\bs\pi$ for all $\bs\pi\in\bbar{\Y}_1(\T^d)$ we obviously have that $\widehat{D}_\sharp\bs{Q}_*^{\infty,\infty}$ is concentrated on the measurable subspace 
	 $$\W_0:=\big\{\bs\rho\in L_{w^*}^\infty(0,T;\Y_1(\T^d))\bigm|\ls\Lambda,\bs\rho_t\rs\leq\mathfrak{m}\;\mbox{a.s.~}\forall\;t\in[0,T]\big\}.$$ 
	 Thus it suffices to show that $\lim_{A\ra+\infty}I_A(\bs\rho)=0$ for all all $\bs\rho\in\W_0$. So let $\bs\rho\in\W_0$. There exists then a Borel set $E_{\bs\rho}\subs[0,T]$ of full Lebesgue measure in $[0,T]$ such that $\bs\rho_t\in\Y_1(\T^d)$ and $\bs\rho_t(\Lambda)\leq\mathfrak{m}$ for all $t\in E_{\bs\pi}$. But then for all $t\in E_{\bs\rho}$ we have that $\Psi_A(\Lambda)\leq\Lambda\in L^1(\bs\rho_t)$ so that the family $\{\Psi_A(\Lambda)\}_{A\geq 0}$ is dominated by the $L^1(\bs\rho_t)$-function $\Lambda$. Since obviously $\Psi_A(\Lambda)\lra 0$ as $A\ra+\infty$ it follows by the dominated convergence theorem that 
	 $$\lim_{A\ra+\infty}\int_{\T^d\x\RR_+}\Psi_A(\Lambda)\df\bs\rho_t=0,\quad\forall\;t\in E_{\bs\rho}$$ and since 
	 $$\int_{\T^d\x\RR_+}\Psi_A(\Lambda)\df\bs\rho_t\leq\int\Lambda\df\bs\rho_t\leq\mathfrak{m},\quad\forall\;t\in E_{\bs\pi}$$ is follows by another application of the dominated convergence theorem that 
	 $$\lim_{A\ra+\infty}I_A(\bs\rho)=\lim_{A\ra+\infty}\int_0^T\int_{\T^d\x\RR_+}\Psi_A(\Lambda)\df\bs\rho_t\leq\int\Lambda\df\bs\rho_t\df t=0$$
	 and the proof of the lemma is complete.$\hfill\Box$
	 
	 \begin{rem} Since the linear functional $\lls\Psi_A(\Lambda\mn M),\cdot\rrs$ is $w^*$-continuous for each $M,A>0$
	  \begin{align}\label{ToOneInd1}
	 \limsup_{\ell,N\ra+\infty}\EE^N\int_0^T\fr{N^d}\sum_{x\in\T_N^d}\Psi_A(\eta_t^\ell(x)\mn M)\df t=\max_{\bs{Q}\in\bs{\mathcal{Q}}^{\infty,\infty}}\int\lls\Psi_A(\Lambda\mn M),\bs{\pi}\rrs\df\bs{Q}(\bs{\pi}),\end{align}
	 where $\bs{\mathcal{Q}}^{\infty,\infty}=\Lim_{\ell,N\ra+\infty}\bs{Q}^{N,\ell}$. In the proof of Lemma~\eqref{TwoBlockCutLargeDens} we have shown that 
	 $$\mathcal{I}_{A,M}(\bs{Q}):=\int\lls\Psi_A(\Lambda\mn M),\bs{\pi}\rrs\df\bs{Q}(\bs{\pi})$$
	 converges pointwise to $0$ as $M\ra+\infty$ and then $A\ra+\infty$ and thus one could wonder whether this pointwise convergence could be strengthened to $\Gamma$-convergence of the maps $-\mathcal{I}_{A,M}$ to the map 
	 $$-\mathcal{I}_A(\bs{Q}):=-\int\lls\Psi_A(\Lambda),\widehat{D}(\bs\pi)\rrs\df\bs{Q}(\bs\pi)$$as $M\ra+\infty$, and then the $\Gamma$-convergence of the maps $-\mathcal{I}_A$ to zero as $A\ra+\infty$, which would ensure that
	 $$\lim_{A,M\ra+\infty}\max_{\bs{Q}\in\bs{\mathcal{Q}}^{\infty,\infty}}\int\lls\Psi_A(\Lambda\mn M),\bs{\pi}\rrs\df\bs{Q}(\bs{\pi})=0.$$
	 This is not true for the weak topology on $\PP L_{w^*}^\infty(0,T;\bbar{\Y}_1(\T^d))$ induced by the $w^*$-topology. For example, note that the map $-\mathcal{I}_A$ is not $w^*$-lower semicontinuous and thus can not be a $\Gamma$-limit.$\hfill\Box$
	 \end{rem}
	 

	 To complete the proof of~\eqref{TDBE} it remains now to show that~\eqref{TBETr} holds. By introducing a parameter $A>0$ that will eventually be sent to $+\infty$ and writing the map $\Psi_M(a,b)=|a\mn M-b\mn M|$, $ a,b\in\RR_+$, as 
	 $$\Psi_M(a,b)=\Psi_{M}(a,b)\1_{[0,A]}(a\mx b)+\Psi_{M}(a,b)\1_{(A,\infty)}(a\mx b)=:\Psi_M^{\leq A}(a,b)+\Psi_M^{>A}(a,b),$$ it follows by the cut-off of large densities in Lemma~\ref{TwoBlockCutLargeDens} that in order to prove~\eqref{TBETr} it suffices to show that 
	 \begin{equation}\label{ToCompleteTDBE}
	 \limsup_{M\uparrow\infty,\ell\uparrow\infty,\ee\downarrow 0,N\uparrow\infty}\EE^N\int_0^T\fr{N^d[N\ee]_\star^{2d}}\sum_{x\in\T_N^d}\sum_{\substack{|y|\mx|z|\leq[N\ee]\\2\ell<|y-z|}}\tau_x\Psi_M^{\leq A}\big(\eta_t^\ell(y),\eta_t^\ell(z)\big)\df t=0,\quad\forall A>0,
	 \end{equation}
	 where $\tau_x\Psi_M^{\leq A}\big(\eta_t^\ell(y),\eta_t^\ell(z)\big):=\Psi_M^{\leq A}\big(\eta_t^\ell(x+y),\eta_t^\ell(x+z)\big)$. Indeed, if~\eqref{ToCompleteTDBE} holds then for every $A>0$ the iterated limit superior in~\eqref{TBETr} is bounded above by
	 $$\limsup_{M\uparrow\infty,\ell\uparrow\infty,\ee\downarrow 0,N\uparrow\infty}\EE^{k_N^{(\ell)}}\int_0^T\fr{(k_N^{(\ell)})^d[k_N^{(\ell)}\ee]_\star^{2d}}\sum_{x\in\T_N^d}\sum_{\substack{|y|\leq[k_N^{(\ell)}\ee]\\|z|\leq[k_N^{(\ell)}\ee]}}\Psi_M^{>A}\big(\eta_t^{m_\ell}(x+y),\eta_t^{m_\ell}(x+z)\big)\df t.$$
	 By the elementary inequality $|a-b|\1_{(A,\infty)}(a\mx b)\leq a\1_{(A,\infty)}(a)+b\1_{(A,\infty)}(b)$ which holds for all $a,b\geq 0$
	 if follows that for all $M>A>0$ and all $a,b\geq 0$
	 \begin{align*}
	 \Psi_M^{>A}(a,b)&=|a\mn M-b\mn M|\1_{(A,\infty)}(a\mx b)=|a\mn M-b\mn M|\1_{(A,\infty)}\big((a\mx b)\mn M\big)\\
	 &=|a\mn M-b\mn M|\1_{(A,\infty)}\big((a\mn M)\mx(b\mn M)\big)\\
	 &\leq(a\mn M)\1_{(A,\infty)}(a\mn M)+(b\mn M)\1_{(A,\infty)}(b\mn M)\\
	 &=(a\mn M)\1_{(A,\infty)}(a)+(b\mn M)\1_{(A,\infty)}(b).
	 \end{align*}
	 Consequently, if~\eqref{ToCompleteTDBE} holds then for every $A>0$ the iterated limit superior in~\eqref{TBETr} is bounded above by 
	 \begin{align*}
	 2\limsup_{M\uparrow\infty,\ell\uparrow\infty,\ee\downarrow 0,N\uparrow\infty}&\EE^{k_N^{(\ell)}}\int_0^T\fr{(k_N^{(\ell)})^d[k_N^{(\ell)}\ee]_\star^{d}}\sum_{x\in\T_N^d}\sum_{|y|\leq[k_N^{(\ell)}\ee]}\big(\eta_t^{m_\ell}(x+y)\mn M\big)\1_{\{\eta_t^{m_\ell}(x+y)>A\}}\df t\\
	 &= 2\limsup_{M\uparrow\infty,\ell\uparrow\infty,N\uparrow\infty}\EE^{k_N^{(\ell)}}\int_0^T\fr{(k_N^{(\ell)})^d}\sum_{x\in\T_N^d}\big(\eta_t^{m_\ell}(x)\mn M\big)\1_{\{\eta_t^{m_\ell}(x)>A\}}\df t,
	 \end{align*}
	 which converges to zero as $A\uparrow\infty$ by Lemma~\eqref{TwoBlockCutLargeDens}.
	 
	  But for every $A,M>0$ we have that 
	 $$\Psi_M^{\leq A}(a,b)=|a\mn M-b\mn M|\1_{[0,A]}(a\mx b)\leq|a-b|\1_{[0,A]}(a\mx b).$$
	 Consequently in order to complete the proof of the truncated double block estimate it remains to show that for all $A>0$ the term 
	 $$\EE^N\int_0^T\fr{N^d[N\ee]_\star^{2d}}\sum_{\substack{x\in\T_N^d\\|z|\leq[N\ee]}}\sum_{2\ell<|y-z|\leq 2[N\ee]}\big|\eta_t^\ell(x+y)-\eta_t^\ell(x+z)\big|\1_{\{\eta^\ell(x+y)\mx(x+z)\leq A\}}\df t$$
	 converges to zero as $N\ra+\infty$, $\ee\ra 0$ and then finally $\ell\ra+\infty$. By making the change of variables $x'=x+z$ and $y':=y-z$ the summation of all $|z|\leq[N\ee]$ disappears and this term becomes equal to 
	 $$\fr{[N\ee]_\star^d}\sum_{2\ell<|y|\leq2[N\ee]}\EE^N\int_0^T\fr{N^d}\sum_{x\in\T_N^d}\big|\eta_t^\ell(x+y)-\eta_t^\ell(x)\big|\1_{\{\eta_t^\ell(x+y)\mx\eta_t(x)\leq A\}}\df t .$$
	 Finally by replacing the average over $y\in\T_N^d$ with $2\ell<|y|\leq 2[N\ee]$ by the supremum of the summands and using the bounds~\eqref{EntropyDirichletFormBounds} on the entropy and Dirichlet form of the time averaged law $\bar{\mu}_T^N:=\fr{T}\int_0^T\mu_t^N\df t$ with respect to $\nu_{\rho_*}^N$, $\rho_*\in(0,\rho_c)$, this term is bounded above by 
	 $$\frac{(2[N\ee])_\star^d-(2\ell)_\star^d}{[N\ee]_\star^d}\sup_{\substack{H_N(f)\leq C_0N^d\\D_N(f)\leq C_0N^{d-2}}}
	 \sup_{2\ell<|y|\leq 2[N\ee]}\int\fr{N^d}\sum_{x\in\T_N^d}\big|\eta^\ell(x+y)-\eta^\ell(x)\big|\1_{\{\eta^\ell(x+y)\mx\eta(x)\leq A\}}f\df\nu_{\rho_*}^N.$$
	 Since $\lim_{N\uparrow\infty}\frac{(2[N\ee])_\star^d-(2\ell)_\star^d}{[N\ee]_\star^d}=2$ for all $\ell\in\ZZ_+$, $\ee>0$ the proof of the truncated double block estimate~\eqref{TDBE} is reduced to showing that
	 $$\limsup_{\ell\uparrow\infty,\ee\downarrow 0,N\uparrow\infty}\sup_{\substack{H_N(f)\leq C_0N^d\\D_N(f)\leq C_0N^{d-2}}}
	 \sup_{2\ell<|y|\leq 2[N\ee]}\int\fr{N^d}\sum_{x\in\T_N^d}\big|\eta^\ell(x+y)-\eta^\ell(x)\big|\1_{\{\eta^\ell(x+y)\mx\eta(x)\leq A\}}f\df\nu_{\rho_*}^N=0$$
	 for every $A>0$. Since the large densities have been cut, this term can now be handled as in the proof in~\cite[Section 5.5]{Kipnis1999a} and thus the proof of the truncated double block estimate~\eqref{TDBE} is complete. 
	 
	 We prove next the second claim of part (a). So we set $\bbar{\bs{Q}}_1^{N,\ell;M;\ee}:=(\bs\pi^{k_N^{(\ell)},m_\ell},\bs\pi^{k_N^{(\ell)},m_\ell;M;\ee})_\sharp P^{k_N^{(\ell)}}$ and let 
	 $$\bbar{\bs{Q}}_1\in\bbar{\bs{\mathcal{Q}}}_1^{\infty,\infty;0;\infty}:=\Lim_{M\ra+\infty}\Lim_{\ell\ra+\infty}\Lim_{\ee\ra 0}\Lim_{N\ra+\infty}\bbar{\bs{Q}}_1^{N,\ell;M;\ee}.$$
	 be a subsequential limit point of the family $\{\bbar{\bs{Q}}_1^{N,\ell;M;\ee}\}$. Then, denoting by $$\bs\pi^\infty,\bs\pi^0\colon L_{w^*}^\infty(0,T;\bbar{\PP}_1(\T^d\x\RR_+))^2\to L_{w^*}^\infty(0,T;\bbar{\PP}_1(\T^d\x\RR_+))$$ the natural projections on the first and second coordinate respectively, it follows by portmanteau theorem and~\eqref{DBE} that for all $G\in L^1(0,T;C(\T^d))$, all $\Psi\in\bbar{C}_1(\RR_+)$ and all $\delta>0$
\begin{align*}
\bbar{\bs{Q}}_1\big\{|\lls G,B_\Psi(\bs\pi^\infty-\bs\pi^0)\rrs|&>\delta\big\}\leq\limsup_{M\uparrow\infty,\ell\uparrow\infty,\ee\downarrow 0,N\uparrow \infty}\bbar{\bs{Q}}_1^{N,\ell;M;\ee}\big\{|\lls G,B_\Psi(\bs\pi^\infty-\bs\pi^0)\rrs|>\delta\big\}\\
&\leq\limsup_{M\uparrow\infty,\ell\uparrow\infty,\ee\downarrow 0,N\uparrow \infty}P^{k_N^{(\ell)}}\big\{|\lls G,B_\Psi(\bs\pi^{k_N^{(\ell)},m_\ell}-\bs\pi^{k_N^{(\ell)},m_\ell;M;\ee})\rrs|>\delta\big\}=0.
\end{align*}
	 Since this holds for all $\delta>0$ it follows that
	 $$\bbar{\bs{Q}}_1\big\{\lls G,B_\Psi(\bs\pi^\infty)\rrs=\lls G,B_\Psi(\bs\pi^0)\rrs\big\}=1$$
	 for all $G\in L^1(0,T;C(\T^d))$ and all $\Psi\in\bbar{C}_1(\RR_+)$. Since $L^1(0,T;C(\T^d))$ is separable it follows that 
	 $$\bbar{\bs{Q}}_1\big\{B_\Psi(\bs\pi^\infty)=B_\Psi(\bs\pi^0)\big\}=1,\quad\forall\Psi\in\bbar{C}_1(\RR_+).$$
	 
	 The space $\bbar{C}_1(\RR_+)$ is also separable. Indeed, $C_1(\RR_+)$ is separable since it is isometric to $C_0(\RR_+)$ and thus there exists a countable subset $D\subs C_1(\RR_+)$ dense in $C_1(\RR_+)$. Then the set $\bbar{D}$ of all maps $\Psi\in\bbar{C}_1(\RR_+)$ of the form $\Psi(\lambda)=\Psi_0(\lambda)+q\lambda$ for some $\Psi_0\in D$ and some $q\in\QQ$ is obviously countable and is we will check it is also dense in $\bbar{C}_1(\RR_+)$. Indeed, let $\Psi\in\bbar{C}_1(\RR_+)$. Then the map $\Psi_0(\lambda):=\Psi(\lambda)-\Psi'(\infty)\lambda$ is in $C_1(\RR_+)$ and thus there exists a sequence $\{\Psi_{0,n}\}\subs D$ such that $\lim_{n\ra+\infty}\|\Psi_{0,n}-\Psi_0\|_{\infty,1}=0$. Then if $\{q_n\}_{n=1}^\infty\subs\QQ$ is a sequence of rational numbers converging to $\Psi'(\infty)$ then the sequence of maps $\Psi_n(\lambda)=\Psi_{0,n}(\lambda)+q_n\lambda$ converges to $\Psi$ in $\bbar{C}_1(\RR_+)$ since 
	 $$\|\Psi_n-\Psi\|_{\infty,1}\leq\|\Psi_{0,n}-\Psi_0\|+|\Psi'(\infty)-q_n|.$$ Therefore the second claim of part (a) follows if we show that 
	 $$\bigcap_{\Psi\in\bbar{C}_1(\RR_+)}\big\{B_\Psi(\bs\pi^\infty)=B_\Psi(\bs\pi^0)\big\}=\bigcap_{\Psi\in\bbar{D}}\big\{B_\Psi(\bs\pi^\infty)=B_\Psi(\bs\pi^0)\big\}.$$ For this it suffices to show that if $\{\Psi_k\}\subs\bbar{C}_1(\RR_+)$ is a sequence converging in norm to $\Psi\in\bbar{C}_1(\RR_+)$ and $(\bs\pi^\infty,\bs\pi^0)\in L_{w^*}^\infty(0,T;\bbar{\Y}_1(\T^d))^2$ is a pair of trajectories such that $B_{\Psi_k}(\bs\pi^\infty)=B_{\Psi_k}(\bs\pi^0)$ for all $k\in\NN$, then $B_\Psi(\bs\pi^\infty)=B_\Psi(\bs\pi^0)$. But this is true since if $\lim_{k\ra+\infty}\|\Psi_k-\Psi\|_{\infty,1}=0$ then for all $\bs\pi\in L_{w^*}^\infty(0,T;\Y_1(\T^d))$
	 \begin{align*}
	 \|B_{\Psi_k}(\bs\pi)-B_\Psi(\bs\pi)\|_{TV,1;\infty}&\leq\sup_{\|G\|_{\infty,1;1}\leq 1}\lls |G(U)||\Psi_k(\Lambda)-\Psi(\Lambda)|,\bs\pi\rrs\\
	 &\leq\|\Psi_k-\Psi\|_{\infty,1}\sup_{\|G\|_{\infty,1;1}\leq 1}\lls |G(U)|(1+\Lambda),\bs\pi\rrs\\
	 &\leq\|\Psi_k-\Psi\|_{\infty,1}(1+\|B(\bs\pi)\|_{TV;\infty})\stackrel{k\ra+\infty}\lra 0.
	 \end{align*}
	 Therefore if $\{\Psi_k\}$ converges to $\Psi$ in $\bbar{C}_1(\RR_+)$ the sequence of operators $\{B_{\Psi_k}\}$ converges strongly to $B_\Psi$ and so if $(\bs\pi^\infty,\bs\pi^0)\in L_{w^*}^\infty(0,T;\bbar{\Y}_1(\T^d))$ is a pair of trajectories such that $B_{\Psi_k}(\bs\pi^\infty)=B_{\Psi_k}(\bs\pi^0)$ for all $k\in\NN$ then 
	 $$\|B_\Psi(\bs\pi^\infty)-B_\Psi(\bs\pi^0)\|_{TV;\infty}\leq\|B_\Psi(\bs\pi^\infty)-B_{\Psi_k}(\bs\pi^\infty)\|_{TV;\infty}+\|B_{\Psi_k}(\bs\pi^0)-B_\Psi(\bs\pi^0)\|_{TV;\infty}$$ for all $k\in\NN$ and taking the limit as $k\ra+\infty$ we conclude that $B_\Psi(\bs\pi^\infty)=B_\Psi(\bs\pi^0)$. This completes the proof of part (a) of Theorem~\ref{TBCTheorem}.
	 
	 \noindent(b) We start by proving (i), i.e.~that in the limit as $N\ra+\infty$, $\ee\ra 0$ and then $\ell\ra 0$ the laws the empirical processes $\bs\pi^{N,\ell}$ and $\bs\pi^{N,\ee}$ have that same barycentric projection. This point in the proof can be shown for the whole family $(\bs\pi^{N,\ell},\bs\pi^{N,\ee})$. So we set $\bbar{\bs{Q}}^{N,\ell,\ee}:=(\bs\pi^{N,\ell},\bs\pi^{N,\ee})_\sharp P^N$ and let $$\bbar{\bs{Q}}\in\bbar{\bs{\mathcal{Q}}}^{\infty,\infty,0}:=\Lim_{\ell\ra+\infty}\Lim_{\ee\ra 0}\Lim_{N\ra+\infty}\bbar{\bs{Q}}^{N,\ell,\ee}$$ be a subsequential limit of this family. We have to show that 
 $$\bbar{\bs{Q}}\big\{(\bs\pi^\infty,\bs\pi^0)\in L_{w^*}^\infty\big(0,T;\bbar{\Y}_1(\T^d)\big)^2\bigm|B(\bs\pi^\infty)=B(\bs\pi^0)\big\}=1.$$ Since $L^1(0,T;C(\T^d))$ is separable in order to show this it suffices to show that 
 $$\bbar{\bs{Q}}\big\{(\bs\pi^\infty,\bs\pi^0)\bigm|||\lls f,B(\bs\pi^\infty)-B(\bs\pi^0)\rrs|>\delta\big\}=0,\quad\forall f\in L^1(0,T;C(\T^d)),\;\delta>0.$$
  Indeed, the set $A_{f,\delta}:=\{|\lls f,B(\bs\pi^\infty)-B(\bs\pi^0)\rrs|>\delta\}$ is open and thus by the portmanteau theorem 
   \begin{align*}
  \bbar{\bs{Q}}(A_{f,\delta})&\leq\limsup_{\ell\ra+\infty}\limsup_{\ee\ra 0}\limsup_{N\ra+\infty}\bbar{\bs{Q}}^{N,\ell,\ee}\{|\lls f,B(\bs\pi^\infty)-B(\bs\pi^0)\rrs|>\delta\}\\
  &\leq\limsup_{\ell\ra+\infty}\limsup_{N\ra+\infty}P^N\{|\lls f,B(\bs\pi^{N,\ell})-\pi^N\rrs|>\delta/2\}\\
  &\quad+\limsup_{\ee\ra 0}\limsup_{N\ra+\infty}P^N\{|\lls f,B(\bs\pi^{N,\ee})-\pi^N\rrs|>\delta/2\}\\
  &\leq\limsup_{\ell,N\ra+\infty}P^N\Big\{\Big|\Big\lls \fr{\ell_\star^d}\sum\nolimits_{|y|\leq\ell}f\Big(\frac{\cdot-y}{N}\Big)-f,\pi^N\Big\rrs\Big|>\frac{\delta}{2}\Big\}\\
  &\quad+\limsup_{\ee\ra 0}\limsup_{N\ra+\infty}P^N\Big\{\Big|\Big\lls \fr{[N\ee]_\star^d}\sum\nolimits_{|y|\leq[N\ee]}f\Big(\frac{\cdot-y}{N}\Big)-f,\pi^N\Big\rrs\Big|>\frac{\delta}{2}\Big\}.
  \end{align*}
  Since $f\in L^1(0,T;C(\T^d))$ it follows that these two terms vanish.
  
   For the proof (ii) we have to consider subsequential limit points along the subfamily $\{\bbar{\bs{Q}}_*^{N,\ell,\ee}\}_{(N,\ell,\ee)}$ of $\bbar{\bs{Q}}^{N,\ell,\ee}:=(\bs\pi^{N,\ell},\bs\pi^{N,\ee})_\sharp P^N$ defined in ~\eqref{DBLAW} so that we will be able to apply the truncated double-block estimate of part (a). In order to prove (ii) it suffices to show that for any limit point $$\bbar{\bs{Q}}_*\in\bbar{\bs{\mathcal{Q}}}_*^{\infty,\infty,0}:=\Lim_{\ell\uparrow\infty}\Lim_{\ee\downarrow 0}\Lim_{N\uparrow\infty}\bbar{\bs{Q}}_*^{N,\ell,\ee}\subs\bbar{\bs{\mathcal{Q}}}^{\infty,\infty,0}$$ any non-decreasing map $\Psi\in\bbar{C}_1(\RR_+)$, any $G\in L^1(0,T;C_+(\T^d))$ and any $\delta>0$
   \begin{equation}\label{For2}\bbar{\bs{Q}}_*\big\{(\bs\pi^\infty,\bs\pi^0)\bigm|\lls G,B_\Psi\circ\widehat{D}(\bs\pi^\infty-\bs\pi^0)\rrs>\delta\big\}=0.
   	\end{equation}
   	Indeed, if this holds for all $\delta>0$ then 
   	$$\bbar{\bs{Q}}_*\big\{(\bs\pi^\infty,\bs\pi^0)\bigm|\lls G,B_\Psi\circ\widehat{D}(\bs\pi^\infty)\rrs\leq\lls G,B_\Psi\circ\widehat{D}(\bs\pi^0)\rrs\big\}=1$$
   	for all $G\in L^1(0,T;C_+(\T^d))$ and all non-decreasing $\Psi\in\bbar{C}_1(\RR_+)$. Let $\bbar{C}_{1,\uparrow}(\RR_+)$ denote the space of all non-decreasing maps $\Psi\in\bbar{C}_1(\RR_+)$. Since the spaces $L^1(0,T;C(\T^d))$ and $\bbar{C}_1(\RR_+)$ are separable, the subspaces $L^1(0,T;C(\T^d))$ and $\bbar{C}_{1,\uparrow}(\RR_+)$ are also separable and thus by arguments similar to the ones in the proof of the second claim of part (a) it follows that $\bbar{\bs{Q}}_*$ is concentrated on the set 
   	$$\W_0=\bigcap_{\Psi\in\bbar{C}_{1,\uparrow}(\RR_+)}\bigcap_{G\in L^1(0,T;C_+(\T^d))}\big\{(\bs\pi^\infty,\bs\pi^0)\bigm|\lls G,B_\Psi\circ\widehat{D}(\bs\pi^\infty)\rrs\leq\lls G,B_\Psi\circ\widehat{D}(\bs\pi^0)\rrs\big\}.$$
 Now, on the set $\bbar{\Y}_1(\T^d)$ the map $B_\Psi\circ\widehat{D}$ takes the form 
   	$$B_\Psi\circ\widehat{D}(\bs\pi)=B_{1,\Psi}(j^*\bs\pi)\equiv B_{1,\Psi}(\bs\rho_{\bs\pi})=b_\Psi(\bs\rho)\df\LL_{\T^d}\equiv\int\Psi(\lambda)\df\bs\rho_{\bs\pi}^u(\lambda)\df u,$$ 
   	where $(\bs\rho_{\bs\pi}^u)_{u\in\T^d}$ is the Lebesgue a.s.~uniquely determined disintegration of $\bs\rho_{\bs\pi}$, and thus since $\bbar{\bs{Q}}_*$ is also supported by the set $L_{w^*}^\infty(0,T;\bbar{\Y}_1(\T^d))^2$  we obtain that for $\bbar{\bs{Q}}_*$-a.s.~all pairs $(\bs\pi^\infty,\bs\pi^0)$ 
   	$$\int_0^T\int_{\T^d}G_t(u)\int_{\RR_+}\Psi(\lambda)\df\bs\rho_{\bs\pi^\infty_t}^u(\lambda)\df u\df t\leq\int_0^T\int_{\T^d}G_t(u)\int_{\RR_+}\Psi(\lambda)\df\bs\rho_{\bs\pi^0_t}^u(\lambda)\df u\df t$$
   	for all $G\in L^1(0,T;C_+(\T^d))$ and all non decreasing $\Psi\in\bbar{C}_1(\RR_+)$. Since $C([0,T]\x\T^d)\subs L^1(0,T;C(\T^d))$ this implies that  $$\int\Psi(\lambda)\df\bs\rho_{\bs\pi^\infty_t}^u(\lambda)\leq\int\Psi(\lambda)\df\bs\rho_{\bs\pi^0_t}^u(\lambda)\quad\mbox{a.s.~for all }(t,u)\in[0,T]\x\T^d,$$
   	which in turn implies claim (ii). Thus claim (ii) is reduced to proving~\eqref{For2}.

Next we note that it suffices to prove~\eqref{For2} under the additional assumption that $\Psi$ is sublinear. Indeed, let us that~\eqref{For2} holds for sublinear maps and let $\Psi$ be asymptotically linear. Then the maps $\Psi_M:=\Psi(\cdot\mn M)$, $M>0$, are sublinear and 
$$\lls G,B_{\Psi_M}\circ\widehat{D}(\bs\pi)\rrs=\lls G,B_{\Psi_M}(\bs\pi)\rrs=\lls G(U)\Psi(\Lambda\mn M),\bs\pi\rrs=\lls G(U)\Psi(\Lambda),\Pi_M^*\circ j^*\bs\pi\rrs$$ for all $M>0$. Therefore since $\Pi_M^*\circ j^*$ $w^*$-converges pointwise to $\widehat{D}$ we have that for any $\bs\pi\in L_{w^*}^\infty(0,T;\bbar{\MMM}_1(\T^d\x\RR_+))$, $G\in L^1(0,T;C(\T^d))$ 
$$\lim_{M\ra+\infty}\lls G,B_{\Psi_M}\circ\widehat{D}(\bs\pi)\rrs=\lim_{M\ra+\infty}\lls G(U)\Psi(\Lambda),\Pi_M^*\circ j^*\bs\pi\rrs=\lls G(U)\Psi(\Lambda),\widehat{D}(\bs\pi)\rrs=\lls G,B_\Psi\circ\widehat{D}(\bs\pi)\rrs.$$
Consequently 
$$\big\{(\bs\pi^\infty,\bs\pi^0)\bigm|\lls G,B_\Psi\circ\widehat{D}(\bs\pi^\infty-\bs\pi^0)\rrs>\delta\big\}\subs\bigcup_{M=1}^\infty
\big\{(\bs\pi^\infty,\bs\pi^0)\bigm|\lls G,B_{\Psi_M}\circ\widehat{D}(\bs\pi^\infty-\bs\pi^0)\rrs>\delta\big\}$$ and thus if~\eqref{For2} holds for all sublinear maps it also holds for all asymptotically linear maps. Since for sublinear maps $\Psi$ it holds that $B_\Psi=B_\Psi\circ\widehat{D}$ in order to prove (ii) it suffice to show that 
for any limit point $\bbar{\bs{Q}}_*\in\bbar{\bs{\mathcal{Q}}}_*^{\infty,\infty,0}$ for all $G\in L^1(0,T;C_+(\T^d))$ all non-decreasing maps $\Psi\in C_1(\RR_+)$ and all $\delta>0$ 
\begin{equation}\label{For2ForSublin}\bbar{\bs{Q}}_*\big\{(\bs\pi^\infty,\bs\pi^0)\bigm|\lls G,B_\Psi(\bs\pi^\infty-\bs\pi^0)\rrs>\delta\big\}=0.
\end{equation}


So let $\bbar{\bs{Q}}_*\in\bbar{\bs{\mathcal{Q}}}_*^{\infty,\infty,0}$. There exists then a diverging sequence $\{m_\ell^{(1)}\}_{l=1}^\infty$, sequences $\{\ee_i^{(1;\ell)}\}_{\ell=1}^\infty$ converging to $0$ as $i\ra\infty$ and diverging sequences $\{k_N^{(1;\ell,i)}\}_{N=1}^\infty$ such that 
$$\bbar{\bs{Q}}_*=\lim_{\ell\uparrow\infty}\lim_{i\uparrow\infty}\lim_{N\uparrow\infty}\bbar{\bs{Q}}_*^{k_N^{(1;\ell,i)},m_\ell^{(1)},\ee_i^{(1;\ell)}}$$
and then, setting $A^\delta_{G,\Psi}:=\{(\bs\pi^\infty,\bs\pi^0)\bigm|\lls G,B_\Psi(\bs\pi^\infty-\bs\pi^0)\rrs>\delta\big\}$, we have by the portmanteau theorem 
\begin{align*}
\bbar{\bs{Q}}_*(A_{G,\Psi}^\delta)&\leq\liminf_{\ell,i,N\uparrow\infty}\bbar{\bs{Q}}_*^{k_N^{(1;\ell,i)},m_\ell^{(1)},\ee_i^{(1;\ell)}}(A_{G,\Psi}^\delta)\leq\limsup_{\ell\uparrow\infty,\ee\downarrow\infty,N\uparrow\infty}\bbar{\bs{Q}}_*^{N,\ell,\ee}(A_{G,\Psi}^\delta)\\
&=\limsup_{\ell\uparrow\infty,\ee\downarrow\infty,N\uparrow\infty}P^{k_N^{(\ell)}}\big\{\lls G,B_\Psi(\bs\pi^{k_N^{(\ell)},m_\ell}-\bs\pi^{k_N^{(\ell)},\ee})\rrs>\delta\big\}.
\end{align*}
By interpolating between the processes $\bs\pi_*^{k_N^{(\ell)},m_\ell}$ and $\bs\pi^{k_N^{(\ell)},\ee}$ with the processes $\bs\pi_*^{k_N^{(\ell)},m_\ell;M;\ee}$ and $\bs\pi_*^{k_N^{(\ell)},m_\ell,\ee;M}$ and taking the limit as $M\uparrow\infty$
\begin{align*}
\bbar{\bs{Q}}_*(A_{G,\Psi}^\delta)&\leq
\limsup_{M\uparrow\infty,\ell\uparrow\infty,\ee\downarrow\infty,N\uparrow\infty}P^{k_N^{(\ell)}}\big\{|\lls G,B_\Psi(\bs\pi^{k_N^{(\ell)},m_\ell}-\bs\pi^{k_N^{(\ell)};m_\ell;M;\ee})\rrs|>\delta/3\big\}\\
&\quad+\limsup_{M\uparrow\infty,\ell\uparrow\infty,\ee\downarrow\infty,N\uparrow\infty}P^{k_N^{(\ell)}}\big\{\lls G,B_\Psi(\bs\pi^{k_N^{(\ell)},m_\ell;M;\ee}-\bs\pi^{k_N^{(\ell)},m_\ell,\ee;M})\rrs>\delta/3\big\}\\
&\quad+\limsup_{M\uparrow\infty,\ell\uparrow\infty,\ee\downarrow\infty,N\uparrow\infty}P^{k_N^{(\ell)}}\big\{|\lls G,B_\Psi(\bs\pi^{k_N^{(\ell)},m_\ell,\ee;M}-\bs\pi^{k_N^{(\ell)},\ee})\rrs|>\delta/3\big\}.
\end{align*}
By part (a), the first term in right hand side above is equal to zero and therefore 
\begin{align*}
\bbar{\bs{Q}}_*(A_{G,\Psi}^\delta)&\leq\limsup_{M\uparrow\infty,\ell\uparrow\infty,\ee\downarrow\infty,N\uparrow\infty}P^N\big\{\lls G,B_\Psi(\bs\pi^{N,\ell;M;\ee}-\bs\pi^{N,\ell,\ee;M})\rrs>\delta/3\big\}\\
&\quad+\limsup_{M\uparrow\infty,\ell\uparrow\infty,\ee\downarrow\infty,N\uparrow\infty}P^N\big\{|\lls G,B_\Psi(\bs\pi^{N,\ell,\ee;M}-\bs\pi^{N,\ee})\rrs|>\delta/3\big\}.
\end{align*}
The first term now in the right hand-side above is also equal to zero, since for any $N\in\NN$, $\ee>0$, $\ell\in\ZZ_+$, $M>0$ and any $x\in\T_N^d$ we obviously have that 
$$(\eta^\ell(x)\mn M)^{[N\ee]}\leq (\eta^\ell(x))^{[N\ee]}\mn M$$ and therefore since $G_t\geq 0$ for almost all $t\in[0,T]$ and $\Psi$ is sublinear and non-decreasing 
\begin{align*}\lls G,B_\Psi(\bs\pi^{N,\ell;M;\ee})\rrs&=\int_0^T\fr{N^d}\sum_{x\in\T_N^d}G_t\Big(\frac{x}{N}\Big)\Psi\big((\eta_t^\ell(x)\mn M)^{[N\ee]}\big)\df t\\
&\leq\int_0^T\fr{N^d}\sum_{x\in\T_N^d}G_t\Big(\frac{x}{N}\Big)\Psi\big(\eta_t^\ell(x)^{[N\ee]}\mn M\big)\df t=\lls G,B_\Psi(\bs\pi^{N,\ell,\ee;M})\rrs.\end{align*} Consequently, by Chebyshev's inequality in order to complete the proof of claim (ii) it suffices to show that 
\begin{equation}\label{ToCompleteiiFirst}\limsup_{M\uparrow\infty,\ell\uparrow\infty,\ee\downarrow\infty,N\uparrow\infty}\EE^N\big|\lls G,B_\Psi(\bs\pi^{N,\ell,\ee;M}-\bs\pi^{N,\ee})\rrs|=0.\end{equation}
By further interpolating with the process $\bs\pi^{N,\ell,\ee}$ we obtain by~\eqref{DoubleBlockFluidSolidSeparation} that 
$$\limsup_{M\uparrow\infty,\ell\uparrow\infty,\ee\downarrow\infty,N\uparrow\infty}\EE^N\big|\lls G,B_\Psi(\bs\pi^{N,\ell,\ee;M}-\bs\pi^{N,\ee})\rrs|
\leq\limsup_{\ell\uparrow\infty,\ee\downarrow\infty,N\uparrow\infty}\EE^N\big|\lls G,B_\Psi(\bs\pi^{N,\ell,\ee}-\bs\pi^{N,\ee})\rrs|$$
and thus in order to complete the proof of claim (ii) it suffices to show that 
\begin{equation}\label{ToCompleteii}
\limsup_{\ell\uparrow\infty,\ee\downarrow\infty,N\uparrow\infty}\EE^N\big|\lls G,B_\Psi(\bs\pi^{N,\ell,\ee}-\bs\pi^{N,\ee})\rrs|=0.
\end{equation}

Setting $\bbar{\bs{Q}}_2^{N,\ell,\ee}:=(\bs\pi^{N,\ell,\ee},\bs\pi^{N,\ee})_\sharp P^N$ and considering the subsequential limit set 
$$\bbar{\bs{\mathcal{Q}}}_2^{\infty,\infty,0}:=\Lim_{\ell\uparrow\infty,\ee\downarrow 0,N\uparrow\infty}\bbar{\bs{Q}}_2^{N,\ell,\ee}$$ 
we can write
$$\limsup_{\ell\uparrow\infty,\ee\downarrow\infty,N\uparrow\infty}\EE^N\big|\lls G,B_\Psi(\bs\pi^{N,\ell,\ee}-\bs\pi^{N,\ee})\rrs|=
\max_{\bs{Q}\in\bbar{\bs{\mathcal{Q}}}_2^{\infty,\infty,0}}\int\big|\lls G,B_\Psi(\bs\pi^\infty-\bs\pi^0)\rrs|\df\bs{Q}(\bs\pi^\infty,\bs\pi^0)$$ and using this equality and the Moreau-Yosida approximations $\Psi_k$ of $\Psi$ given in~\eqref{MoreauYosida}, we can reduce the proof of~\eqref{ToCompleteii} to the case that the map $\Psi$ is in addition Lipschitz, as in the proof of part (a).

 Now, since $\Psi$ is assumed Lipschitz, by the definition of the processes $\bs\pi^{N,\ell,\ee}$ and $\bs\pi^{N,\ee}$ 
\begin{align}\label{AnotherEasyTerm}
|\lls G,B_\Psi(\bs\pi^{N,\ell,\ee}-\bs\pi^{N,\ee})\rrs|&\leq\int_0^T\fr{N^d}\sum_{x\in\T_N^d}\Big|G_t\Big(\frac{x}{N}\Big)\Big|\big|\Psi\big(\eta_t^{\ell,[N\ee]}(x)\big)-\Psi(\eta_t^{[N\ee]}(x)\big)\big|\df t\nonumber\\
&\leq\Lip_\Psi\int_0^T\fr{N^d}\sum_{x\in\T_N^d}\Big|G_t\Big(\frac{x}{N}\Big)\Big|\cdot\big|\eta_t^{\ell,[N\ee]}(x)-\eta_t^{[N\ee]}(x)\big|\df t\nonumber\\
&\leq\Lip_\Psi\int_0^T\Big(\fr{N^d}\sum_{x\in\T_N^d}\big|\eta_t^{\ell,[N\ee]}(x)-\eta_t^{[N\ee]}(x)\big|\Big)\|G_t\|_\infty\df t.
\end{align}

By a standard computation on consecutive averages (see for example~\cite[(4.14)]{Faggionato2008a}) for any $\ell\leq L\in\ZZ_+$ and any family of functions $\Psi_N\colon\MM_N^d\to\RR$
\begin{equation}\label{ConsecAver}
|(\Psi_N^\ell)^L-\Psi_N^L|\leq\fr{L_\star^d}\sum_{L-\ell<|x|\leq L+\ell}\tau_x|\Psi_N|.\end{equation}
Then by inequalities~\eqref{ConsecAver} and~\eqref{AnotherEasyTerm} we obtain that
\begin{align*}
|\lls G,B_\Psi(\bs\pi^{N,\ell,\ee}-\bs\pi^{N,\ee})\rrs|&\leq\Lip_\Psi\int_0^T\Big(\fr{N^d[N\ee]_\star^d}\sum_{x\in\T_N^d}\sum_{[N\ee]-\ell<|z|\leq[N\ee]+\ell}\eta_t(x+z)\Big)\|G_t\|_\infty\df t\\
&=\Lip_\Psi\frac{([N\ee]+\ell)^d-([N\ee]-\ell)^d}{[N\ee]_\star^d}\int_0^T\ls 1,\pi_t^N\rs\|G_t\|_\infty\df t\\
&\stackrel{P^N-\mbox{a.s.}}=\Lip_\Psi\frac{([N\ee]+\ell)^d-([N\ee]-\ell)^d}{[N\ee]_\star^d}\ls 1,\pi_0^N\rs\|G\|_{\infty,1}\\
&=\Lip_\Psi\frac{2\ell}{[N\ee]_\star^d}O([N\ee]^{d-1})\ls 1,\pi_0^N\rs\|G\|_{\infty,1}.
\end{align*}
Consequently there exists a constant $C_d<+\infty$ such that
$$\limsup_{N\ra+\infty}\EE^N|\lls G,B_\Psi(\bs\pi^{N,\ell,\ee}-\bs\pi^{N,\ee})\rrs|\leq \Lip_\Psi C_d\|G\|_{\infty,1}\limsup_{N\ra+\infty}\Big(\frac{2\ell}{[N\ee]}\int\ls 1,\pi^N\rs\df\mu_0^N\Big)=0,$$
where the last limit superior is equal to $0$ by the $O(N^d)$-entropy assumption and Lemma~\ref{BoundedTotParticlInMacrLim}. This proves~\eqref{ToCompleteii} and completes the proof of (ii).

Claim (iii) is a consequence of claims (i) and (ii). Indeed, by (i) and (ii) it follows that any limit point $\bbar{\bs{Q}}_*$ as $N\ra+\infty$, $\ell\ra+\infty$ and then $M\ra+\infty$ of the family of laws $\bbar{\bs{Q}}_*^{N,\ell,\ee}$ defined in~\eqref{DBLAW} is supported on a measurable set $\W_0\subs L_{w^*}^\infty(0,T;\bbar{\Y}_1(\T^d))^2$ of trajectory pairs such that for any $(\bs\pi^\infty,\bs\pi^0)\in\W_0$ for almost all $t\in[0,T]$ it holds that $B(\bs\pi_t^\infty)=B(\bs\pi_t^0)$ and $\ls f,B_\Psi(\widehat{\bs\pi}_t^\infty)\rs\leq\ls f,B_\Psi(\widehat{\bs\pi}_t^0)\rs$ for all non-decreasing $\Psi\in\bbar{C}_1(\RR_+)$ and all $f\in C_+(\T^d)$. Therefore for any $f\in C_+(\T^d)$
	\begin{align*}
	\ls f,\rho^\perp_{\bs\pi_t^\infty}\rs-\ls f,\rho^\perp_{\bs\pi_t^0}\rs&=\ls f(U)\Lambda,\bs\pi_t^\infty\rs-\ls f(U)\Lambda,\widehat{\bs\pi}_t^\infty\rs-\ls f(U)\Lambda,\bs\pi^0_t\rs+\ls f(U)\Lambda,\widehat{\bs\pi}_t^0\rs\\
	&=\ls f,B(\bs\pi_t^\infty)\rs-\ls f,B(\bs\pi_t^0)\rs+\ls f(U)\Lambda,\widehat{\bs\pi}_t^0\rs-\ls f(U)\Lambda,\widehat{\bs\pi}_t^\infty\rs\\
	&=\ls f,B(\widehat{\bs\pi}_t^0)\rs-\ls f,B(\widehat{\bs\pi}_t^\infty)\rs\geq 0
	\end{align*} 
	for almost all $t\in[0,T]$, where the last term is non-negative by claim (ii) applied for the identity map $\Psi=\mathbbm{id}_{\RR_+}$. This proves (iii)~and completes the proof of (b).
	
\noindent(c) We recall the notation $\bbar{\bs{\mathcal{Q}}}_*^{\infty,\infty,0}:=\Lim_{\ell\uparrow\infty,\ee\downarrow 0,N\uparrow\infty}\bbar{\bs{\mathcal{Q}}}_*^{N,\ell,\ee}$ where $\bbar{\bs{Q}}_*^{N,\ell,\ee}$ are the laws defined in~\eqref{DBLAW}. We start the proof of (c) by noting that the two-blocks estimate~\eqref{TBE} is equivalent to the validity for all $\bbar{\bs{Q}}_*\in\bbar{\bs{\mathcal{Q}}}_*^{\infty,\infty,0}$, all $G\in L^\infty(0,T;C_+(\T^d))$ and all Lipschitz maps $\Psi\in\bbar{C}_{1,\uparrow}(\RR_+)$ with $\Psi(0)=0$ of the equality 
\begin{equation}\label{TBEVariant}
\int|\lls G,B_\Psi(\bs\pi^\infty-\bs\pi^0)\rrs|\df\bbar{\bs{Q}}_*(\bs\pi^\infty,\bs\pi^0)=0.
\end{equation}

 Indeed, on one-hand it is obvious that if the two-blocks estimate~\eqref{TBE} holds then~\eqref{TBEVariant} holds for all $\bbar{\bs{Q}}_*\in\bbar{\bs{\mathcal{Q}}}_*^{\infty,\infty,0}$, all $G\in L^\infty(0,T;C_+(\T^d))$ and all Lipschitz maps $\Psi\in\bbar{C}_{1,\uparrow}(\RR_+)$ with $\Psi(0)=0$. Conversely, let $\bbar{\bs{Q}}_*\in\bbar{\bs{\mathcal{Q}}}_*^{\infty,\infty,0}$ be such that~\eqref{TBEVariant} holds for all $G\in L^\infty(0,T;C_+(\T^d))$ and all Lipschitz maps $\Psi\in\bbar{C}_{1,\uparrow}(\RR_+)$ with $\Psi(0)=0$ and we will show that 
\begin{equation}\label{TBEForAFixedSubLim}
\bbar{\bs{Q}}_*\big\{B_\Psi(\bs\pi^\infty-\bs\pi^0)=0,\;\forall\Psi\in\bbar{C}_{1,\uparrow}(\RR_+)\big\}=1.
\end{equation}

 We note first for a given map $\Psi\in\bbar{C}_1(\RR_+)$, equality~\eqref{TBEVariant} holds for all $G\in L^1(0,T;C(\T^d))$ if and only if it holds for all $G\in L^\infty(0,T;C_+(\T^d))$. Indeed, for any $G\in L^1(0,T;C(\T^d))$ the maps $C_t(u)=\|G_t\|_\infty$ and $G_{+,t}(u):=G_t(u)+\|G_t\|_\infty$, $(t,u)\in[0,T]\x\T^d$, are in $L^1(0,T;C_+(\T^d))$ and 
$$\lls G_+,B_\Psi(\bs\pi)\rrs=\lls G(U)\Psi(\Lambda),\bs\pi\rrs+\int_0^T\ls C_t(U)\Psi(\Lambda),\bs\pi_t\rs\df t=\lls G,B_\Psi(\bs\pi)\rrs+\lls C,B_\Psi(\bs\pi)\rrs$$
for all $\bs\pi\in L_{w^*}^\infty(0,T;\MMM_1(\T^d\x\RR_+))$. Thus we get the estimate 
$$|\lls G,B_\Psi(\bs\pi^\infty-\bs\pi^0)\rrs|\leq|\lls G_+,B_\Psi(\bs\pi^\infty-\bs\pi^0)\rrs|+|\lls C,B_\Psi(\bs\pi^\infty-\bs\pi^0)\rrs|$$
for all $(\bs\pi^\infty,\bs\pi^0)\in L_{w^*}^\infty(0,T;\bbar{\MMM}_{1,+}(\T^d\x\RR_+))$, which shows that if equality~\eqref{TBEVariant} holds for all $G\in L^1(0,T;C_+(\T^d))$ then it also holds for all $G\in L^1(0,T;C(\T^d))$. Finally, since $L^\infty(0,T;C(\T^d))$ is dense in $L^1(0,T;C(T^d))$ it follows that equality~\eqref{TBEVariant} holds for all $L^1(0,T;C(\T^d))$ if and only if it holds for all $G\in L^\infty(0,T;C_+(\T^d))$.

Let now $\Psi\in\bbar{C}_1(\RR_+)$ be non-decreasing and let $G\in L^1(0,T;C_+(\T^d))$. By employing the approximations $\Psi_k(\lambda)=\Psi_{0,k}(\lambda)+\Psi'(\infty)\lambda$ of $\Psi$, where $\Psi_{0,k}$ are the Moreau-Yosida approximations of the sublinear part $\Psi_0(\lambda)=\Psi(\lambda)-\Psi'(\infty)\lambda$ of $\Psi$, one can reduce the proof of the two-blocks estimate~\eqref{TBE}  to the case that $\Psi$ is Lipschitz. Indeed, then $\{\Psi_k\}$ increases to $\Psi$ as $k\uparrow\infty$ and $\|\Psi_k\|_{\infty,1}\leq\|\Psi_{0,k}\|+|\Psi'(\infty)|\leq\|\Psi_0\|_{\infty,1}+|\Psi'(\infty)|$ for all large enough $k\in\NN$, and thus since $\bbar{\bs{Q}}_*$ is supported on the set $L_{w^*}^\infty(0,T;\bbar{\Y}_{1;\mathfrak{m}}(\T^d))^2$ it follows, similarly to the reduction to the case of Lipschitz maps $\Psi$ in part (a), that the maps $(\bs\pi^\infty,\bs\pi^0)\mapsto\lls G,B_{\Psi_k}(\bs\pi^\infty-\bs\pi^0)\rrs$ converge $\bbar{\bs{Q}}_*$-a.s.~pointwise to the map $\lls G,B_\Psi(\bs\pi^\infty-\bs\pi^0)\rrs$ and are dominated by an $L^\infty(\bbar{\bs{Q}}_*)$-function. Thus it follows by the dominated convergence theorem that for a given $G\in L^1(0,T;C_+(\T^d))$ equality~\eqref{TBEVariant} holds for all non-decreasing $\Psi\in\bbar{C}_1(\RR_+)$ if and only if it holds for all non-decreasing $\Psi\in\bbar{C}_1(\RR_+)\cap\Lip(\RR_+)$. Furthermore,~\eqref{TBEVariant} holds for a map $\Psi\in\bbar{C}_1(\RR_+)$ if and only if it holds for the map $\Psi+c$ for any constant $c$ and thus we can also assume that $\Psi(0)=0$ so that  $\Psi\geq 0$ in $\RR_+$ and $0\leq\Psi'(\infty)\leq\Lip_\Psi$.

Thus equality~\eqref{TBEVariant} holds for all $G\in L^1(0,T;C(\T^d))$ and all $\Psi\in\bbar{C}_{1,\uparrow}(\RR_+)$ which implies that 
$$\bbar{\bs{Q}}_*\{\lls G,B_\Psi(\bs\pi^\infty-\bs\pi^0)\rrs=0\}=1,\quad\forall G\in L^1(0,T;C(\T^d)),\;\Psi\in\bbar{C}_{1,\uparrow}(\RR_+).$$
Since the spaces $L^1(0,T;C(\T^d))$ and $\bbar{C}_{1,\uparrow}(\RR_+)$ are separable it follows then similarly to the proof of the second claim of part (a) that~\eqref{TBEForAFixedSubLim} holds. 
 
Using this equivalent characterization of the two-blocks estimate we prove first that if for any subfamily $\big\{\big(k_{k_N^{(1;\ell,i)}}^{(m_\ell^{(1)})},m_{m^{(1)}_\ell},\ee_i^{(\ell)}\big)\big\}$ of $\{(k_N^{(\ell)},m_\ell,\ee)\}$ there exists a further subfamily $\{(\bar{k}_N^{(\ell,i)},\bar{m}_\ell,\bar{\ee}_i^{(\ell)})\}$ as in~\eqref{Nightmare} such that~\eqref{ToProveReplacLemma} holds, then the two-blocks estimate~\eqref{TBE} holds. So let $\bbar{\bs{Q}}_*\in\bbar{\bs{Q}}_*^{\infty,\infty,0}$. There exists a diverging sequence $\{m_\ell^{(1)}\}_{\ell=1}^\infty$, sequences $\{\ee_i^{(1;\ell)}\}_{i=1}^\infty$, $\ell\in\NN$, converging to $0$ as $i\uparrow\infty$ for all $\ell\in\NN$ and diverging sequences $\{k_N^{(1;\ell,i)}\}_{N=1}^\infty$, $(\ell,i)\in\NN^2$, such that 
\begin{equation}\label{SomeLimit}
\bbar{\bs{Q}}_*=\lim_{\ell\ra+\infty}\lim_{i\ra+\infty}\lim_{N\ra+\infty}\bbar{\bs{Q}}_*^{k_N^{(1;\ell,i)},m_\ell^{(1)},\ee_i^{(1;\ell)}}
\end{equation}
By the assumption there exists a further subfamily $\{(\bar{k}_N^{(\ell,i)},\bar{m}_\ell,\bar{\ee}_i^{(\ell)})\}$ as in~\eqref{Nightmare} along which~\eqref{ToProveReplacLemma} holds. Recalling that $\bbar{\bs{Q}}^{N,\ell,\ee}:=(\bs\pi^{N,\ell},\bs\pi^{N,\ee})_\sharp P^N$ then with $\{(\bar{k}_N^{(\ell,i)},\bar{m}_\ell,\bar{\ee}_i^{(\ell)})\}$ being given by~\eqref{Nightmare} 
$$\bbar{\bs{Q}}_*^{k_{k_N^{(2;\ell,i)}}^{\left(1;m_\ell^{(2)},\ee_i^{(2;\ell)}\right)},m_{m_\ell^{(2)}}^{(1)},\ee_{\ee_i^{(2;\ell)}}^{\left(1;m_\ell^{(2)}\right)}}=\bbar{\bs{Q}}^{\bar{k}_N^{(\ell,i)},\bar{m}_\ell,\bar{\ee}_i^{(\ell)}}.$$
Since $\{\bbar{\bs{Q}}^{\bar{k}_N^{(\ell,i)},\bar{m}_\ell,\bar{\ee}_i^{(\ell)}}\}$ is a subfamily of $\{\bbar{\bs{Q}}_*^{k_N^{(1;\ell,i)},m_\ell^{(1)},\ee_i^{(1;\ell)}}\}$ the limit~\eqref{SomeLimit} continues to hold along this subfamily, i.e.
\begin{equation}\label{SomeSubLimit}\bbar{\bs{Q}}_*=\lim_{\ell\ra+\infty}\lim_{i\ra+\infty}\lim_{N\ra+\infty}\bbar{\bs{Q}}^{\bar{k}_N^{(\ell,i)},\bar{m}_\ell,\bar{\ee}_i^{(\ell)}}.
\end{equation}

By passing to a further subfamily which we will continue to denote by $\{(\bar{k}_N^{(\ell,i)},\bar{m}_\ell,\bar{\ee}_i^{(\ell)})\}_{(N,\ell,i)}$ we can further assume that the laws $\bs{Q}^{\bar{k}_N^{(\ell,i)},\bar{m}_\ell,\bar{\ee}_i^{(\ell)}}:=\bs\pi^{\bar{k}_N^{(\ell,i)},\bar{m}_\ell,\bar{\ee}_i^{(\ell)}}\sharp P^{\bar{k}_N^{(\ell,i)}}$ of the double block empirical process $\bs\pi^{N,\ell,\ee}$ defined in~\eqref{DoubleBlockEmpDens} converge along the subfamily $\{(\bar{k}_N^{(\ell,i)},\bar{m}_\ell,\bar{\ee}_i^{(\ell)})\}_{(N,\ell,i)}$ to some $\bs{Q}_*\in\PP L_{w^*}^\infty(0,T;\mathcal{\Y}_{1,\mathfrak{m}}(\T^d))$, i.e. 
\begin{equation}\label{SomeOtherSubLimForDoublePr}
\bs{Q}_*=\lim_{\ell,i,N\uparrow\infty}\bs{Q}^{\bar{k}_N^{(\ell,i)},\bar{m}_\ell,\bar{\ee}_i^{(\ell)}}.
\end{equation}
Of course then~\eqref{SomeSubLimit} continuous to hold and since $\{(\bar{k}_N^{(\ell,i)},\bar{m}_\ell,\bar{\ee}_i^{(\ell)})\}$ is a subfamily of $\{(k_N^{(\ell)},m_\ell,\ee\}$ we necessarily have that 
$$\bs{Q}_*\in\bs{\mathcal{Q}}_*^{\infty,\infty,0}:=\Lim_{\ell\uparrow\infty,\ee\downarrow 0,N\uparrow\infty}\bs{Q}_*^{N,\ell,\ee}$$
where here $\bs{Q}_*^{N,\ell,\ee}:=\bs{Q}^{k_N^{(\ell)},m_\ell,\ee}$ and $\bs{Q}^{N,\ell,\ee}:=\bs\pi^{N,\ell,\ee}_\sharp P^N$.

Since the map $(\bs\pi^\infty,\bs\pi^0)\mapsto|\lls G,B_\Psi(\bs\pi^\infty-\bs\pi^0)\rrs|$ is continuous it follows by~\eqref{SomeSubLimit} and the portmanteau theorem that
\begin{align*}\int|\lls G,B_\Psi(\bs\pi^\infty-\bs\pi^0)\rrs|\df\bbar{\bs{Q}}_*(\bs\pi^\infty,\bs\pi^0)&=\lim_{\ell,i,N\uparrow\infty}\int|\lls G,B_\Psi(\bs\pi^\infty-\bs\pi^0)\rrs|\df\bbar{\bs{Q}}^{\bar{k}_N^{(\ell,i)},\bar{m}_\ell,\bar{\ee}_i^{(\ell)}}(\bs\pi^\infty,\bs\pi^0)\\
&=\lim_{\ell,i,N\uparrow\infty}
\EE^{\bar{k}_N^{(\ell,i)}}\big|\lls G,B_\Psi(\bs\pi^{\bar{k}_N^{(\ell,i)},\bar{m}_\ell},\bs\pi^{\bar{k}_N^{(\ell,i)},\bar{\ee}_i^{(\ell)}})\rrs\big|.
\end{align*}

Interpolating with the process $\bs\pi_*^{\bar{k}_N^{(\ell,i)},\bar{m}_\ell;M;\bar{\ee}_i^{(\ell)}}$ and $\bs\pi_*^{\bar{k}_N^{(\ell,i)},\bar{m}_\ell,\bar{\ee}_i^{(\ell)};M}$ it follows by~\eqref{DBE} in the proof part (a) and~\eqref{ToCompleteiiFirst} that
\begin{align*}
\lim_{\ell,i,N\uparrow\infty}\EE^{\bar{k}_N^{(\ell,i)}}&|\lls G,B_\Psi(\bs\pi^{\bar{k}_N^{(\ell,i)},\bar{m}_\ell}-\bs\pi^{\bar{k}_N^{(\ell,i)},\bar{\ee}_i^{(\ell)}})\rrs|\\
&\leq \limsup_{M,\ell,i,N\uparrow\infty}\EE^{\bar{k}_N^{(\ell,i)}}|\lls G,B_\Psi(\bs\pi^{\bar{k}_N^{(\ell,i)},\bar{m}_\ell;M;\bar{\ee}_i^{(\ell)}}-\bs\pi^{\bar{k}_N^{(\ell,i)},\bar{m}_\ell,\bar{\ee}_i^{(\ell)};M})\rrs|.
 \end{align*}
 
Since we assume that $G\in L^\infty(0,T;C_+(\T^d))$ and $\Psi\geq 0$ is Lipschitz and non-decreasing, for all parameters $(N,\ell,\ee)\in\NN\x\ZZ_+\x(0,\infty)$
 \begin{align}\label{BothTRH}
 |\lls G,B_\Psi(\bs\pi^{N,\ell;M;\ee}-\bs\pi^{N,\ell,\ee;M})\rrs|&\leq\lls G,B_\Psi(\widehat{\bs\pi}^{N,\ell,\ee;M}-\widehat{\bs\pi}^{N,\ell;M;\ee})\rrs\nonumber\\
 &\qquad+\Psi'(\infty)\lls G,\rho_{\bs\pi^{N,\ell;M;\ee}}^\perp-\rho_{\bs\pi^{N,\ell,\ee;M}}^\perp\rrs
 \end{align}
and by the simple identity $a=a\mn M+(a-M)^+$ for all $a\in\RR$ we obtain that for almost all $0\leq t\leq T$  
 \begin{align*}
 \ls G_t,B_\Psi(\widehat{\bs\pi}_t^{N,\ell,\ee;M}-\widehat{\bs\pi}_t^{N,\ell;M;\ee})\rs&=\fr{N^d}\sum_{x\in\T_N^d}G_t\Big(\frac{x}{N}\Big)\big\{\Psi\big(\eta_t^{\ell,[N\ee]}(x)\mn M\big)-\Psi\big((\eta_t^\ell(x)\mn M)^{[N\ee]}\big)\big\}\\
 &\leq\Lip_\Psi\fr{N^d}\sum_{x\in\T_N^d}G_t\Big(\frac{x}{N}\Big)\big(\eta_t^{\ell,[N\ee]}(x)\mn M-(\eta_t^\ell(x)\mn M)^{[N\ee]}\big)\\
 &=\fr{N^d}\Lip_\Psi\sum_{x\in\T_N^d}G_t\Big(\frac{x}{N}\Big)\big((\eta_t^\ell(x)-M)^{+[N\ee]}-(\eta_t^{\ell,[N\ee]}(x)-M)^+\big)\\
 &=\Lip_\Psi\ls G_t,\rho_{\bs\pi_t^{N,\ell;M;\ee}}^\perp-\rho_{\bs\pi_t^{N,\ell,\ee;M}}^\perp\rs.
 \end{align*}

 Consequently, since $\Psi'(\infty)\leq\Lip_\Psi$ it follows by~\eqref{BothTRH} that
 $$|\lls G,B_\Psi(\bs\pi^{N,\ell;M;\ee}-\bs\pi^{N,\ell,\ee;M})\rrs|\leq 2\Lip_\Psi\|G\|_{\infty;\infty}\lls 1,\rho_{\bs\pi^{N,\ell;M;\ee}}^\perp-\rho_{\bs\pi^{N,\ell,\ee;M}}^\perp\rrs,$$
 where $\|G\|_{\infty;\infty}:=\|G\|_{L^\infty(0,T;C(\T^d))}<+\infty$. Therefore in order to show that the two-blocks estimate holds it suffices to show that 
 $$\limsup_{M,\ell,i,N\uparrow\infty}\EE^{\bar{k}_N^{(\ell,i)}}\lls 1,\rho_{\bs\pi^{\bar{k}_N^{(\ell,i)},\bar{m}_\ell;M;\bar{\ee}_i^{(\ell)}}}^\perp-\rho_{\bs\pi^{\bar{k}_N^{(\ell,i)},\bar{m}_\ell,\bar{\ee}_i^{(1;\ell)};M}}^\perp\rrs=0.$$
 
 
 The difference $\rho^\perp_{\bs\pi^{N,\ell;M;\ee}}-\rho^\perp_{\bs\pi^{N,\ell,\ee;M}}$ is a non-negative path-measure valued process and the weight $(\rho^\perp_{\bs\pi_t^{N,\ell;M;\ee}}-\rho^\perp_{\bs\pi_t^{N,\ell,\ee;M}})(x/N)$, $x\in\T_N^d$, can be expressed as 
 \begin{align}\label{BothPosAgain}
 (\rho^\perp_{\bs\pi_t^{N,\ell;M;\ee}}-\rho^\perp_{\bs\pi_t^{N,\ell,\ee;M}})(x/N)&=(\eta^\ell(x)-M)^{+[N\ee]}-(\eta^\ell(x)^{[N\ee]}-M)^+\nonumber\\
 &=[(\eta_t^\ell(x)-M)^+]^{[N\ee]}\1_{[0,M]}(\eta_t^\ell(x)^{[N\ee]})\nonumber\\
 &\qquad+[(\eta_t^\ell(x)-M)^-]^{[N\ee]}\1_{(M,\infty)}(\eta_t^\ell(x)^{[N\ee]})
 \end{align}
 which shows that 
 \begin{align}\label{LastSB}
 \lls 1,\rho_{\bs\pi^{N,\ell;M;\ee}}^\perp-\rho_{\bs\pi^{N,\ell,\ee;M}}^\perp\rrs&\leq 
 \EE^N\int_0^T\fr{N^d}\sum_{x\in\T_N^d}(\eta_t^\ell(x)-M)^{+[N\ee]}\1_{[0,M]}(\eta_t^\ell(x)^{[N\ee]})\df t\nonumber\\
 &\quad+\EE^N\int_0^T\fr{N^d}\sum_{x\in\T_N^d}(\eta_t^\ell(x)-M)^{-[N\ee]}\1_{(M,\infty)}(\eta_t^\ell(x)^{[N\ee]})\df t
 \end{align}

Since we have chosen the subfamily $\{(\bar{k}_N^{(\ell,i)},\bar{m}_\ell,\bar{\ee}_i^{(\ell)})\}$ so that~\eqref{SomeOtherSubLimForDoublePr} holds, it follows by the next lemma that the contribution in the limit as $N\uparrow\infty$, $i\uparrow\infty$, $\ell\uparrow\infty$ and then $M\uparrow$ of the second summand in the right hand side of~\eqref{LastSB} is zero along the subfamily $\{(\bar{k}_N^{(\ell,i)},\bar{m}_\ell,\bar{\ee}_i^{(\ell)})\}$. This shows that the two-blocks estimate holds, since the first summand in the right hand side of~\eqref{LastSB} vanishes along the subfamily $\{(\bar{k}_N^{(\ell,i)},\bar{m}_\ell,\bar{\ee}_i^{(\ell)})\}_{(N,\ell,i)}$ due to~\eqref{ToProveReplacLemma} being true for the subfamily $\{(\bar{k}_N^{(\ell,i)},\bar{m}_\ell,\bar{\ee}_i^{(\ell)})\}$. 
	 \begin{lemma} Suppose that the ZRP starts from a sequence $\{\mu_0^N\in\PP\MM_N^d\}_{N\in\NN}$ of initial profiles having asymptotically $\mathfrak{m}>0$ total mass and let $\{\bar{m}_\ell\}_{\ell=1}^\infty$ be a diverging sequence, let $\{\bar{\ee}_i^{(\ell)}\}_{i=1}^\infty$ be sequences converging to zero for each $\ell\in\NN$ and let $\{\bar{k}_N^{(\ell,i)}\}_{N=1}^\infty$ be diverging sequences for $i,N\in\NN$ such that the double-block empirical laws 
	 	\begin{equation}\label{SubFam2}
	 	\bs{Q}_*^{N,\ell,i}:=(\bs\pi^{\bar{k}_N^{(\ell,i)},\bar{m}_\ell,\bar{\ee}_i^{(\ell)}})_\sharp P^{\bar{k}_N^{(\ell,i)}}:=\bs\pi^{\bar{k}_N^{(\ell,i)},\bar{m}_\ell,\bar{\ee}_i^{(\ell)}}_\sharp P^{\bar{k}_N^{(\ell,i)}}
	 	\end{equation} converge weakly to some probability law $\bs{Q}_*\in\PP L_{w^*}^\infty(0,T;\bbar{\Y}_1(\T^d))$ as $N\uparrow\infty$, $i\uparrow\infty$ and then $\ell\uparrow\infty$. Then 
	 	$$\lim_{M\ra+\infty}\limsup_{\ell,i,N\uparrow\infty}
	 	\EE^{\bar{k}_N^{(\ell,i)}}\int_0^T\fr{(\bar{k}_N^{(\ell,i)})^d}\sum_{x\in\T_{\bar{k}_N^{(\ell,i)}}^d}[(\eta_t^{\bar{m}_\ell}(x)-M)^-]^{[\bar{k}_N^{(\ell,i)}\bar{\ee}_i^{(\ell)}]}\1_{(M,\infty)}(\eta_t^{\bar{m}_\ell}(x)^{[\bar{k}_N^{(\ell,i)}\bar{\ee}_i^{(\ell)}]})=0.$$ 
	 \end{lemma}\textbf{Proof} Let $\Psi_M\in BC(\RR_+)$ be the map $\Psi_M(\lambda)=M\cdot[(\lambda-M+1)^+\mn 1]$ and note that $$M\1_{(M,\infty)}(\lambda)\leq\Psi_M(\lambda)\leq\lambda,\quad\forall\lambda\geq 0.$$ Therefore 
 $$[(\eta^\ell(x)-M)^-]^{[N\ee]}\1_{(M,\infty)}(\eta^\ell(x)^{[N\ee]})\leq M\1_{(M,\infty)}(\eta_t^\ell(x)^{[N\ee]})\leq\Psi_M\big(\eta^\ell(x)^{[N\ee]}\big)$$
 and thus it suffices to show that 
	 $$\limsup_{M,\ell,i,N\uparrow\infty}
	 \EE^{\bar{k}_N^{(\ell,i)}}\int_0^T\fr{(\bar{k}_N^{(\ell,i)})^d}\sum_{x\in\T_{\bar{k}_N^{(\ell,i)}}^d}\Psi_M\big(\eta_t^{\bar{m}_\ell}(x)^{[\bar{k}_N^{(\ell,i)}\bar{\ee}_i^{(\ell)}]}\big)\df t=0.$$
In terms of the double-block empirical density process $\bs\pi_*^{N,\ell,i}$ of the ZRP defined in~\eqref{SubFam2} and the corresponding laws $\bs{Q}_*^{N,\ell,i}:=(\bs\pi^{\bar{k}_N^{(\ell,i)},\bar{m}_\ell,\bar{\ee}_i^{(\ell)}})_\sharp P^{\bar{k}_N^{(\ell,i)}}$ the expected value above can be written 
\begin{align*}
\EE^{\bar{k}_N^{(\ell,i)}}\int_0^T\fr{(\bar{k}_N^{(\ell,i)})^d}\sum_{x\in\T_{\bar{k}_N^{(\ell,i)}}^d}\Psi_M\big(\eta_t^{\bar{m}_\ell}(x)^{[\bar{k}_N^{(\ell,i)}\bar{\ee}_i^{(\ell)}]}\big)\df t&=\EE^{\bar{k}_N^{(\ell,i)}}\lls\Psi_M(\Lambda),\bs\pi^{\bar{k}_N^{(\ell,i)},\bar{m}_\ell,\bar{\ee}_i^{(\ell)}}\rrs\\
&=\int\lls\Psi_M(\Lambda),\bs\pi\rrs\df\bs{Q}_*^{N,\ell,i}(\bs\pi).
\end{align*}
Therefore, since by assumption $\lim_{\ell,i,N\uparrow\infty}\bs{Q}_*^{N,\ell,i}=\bs{Q}_*$ it follows that by the portmanteau theorem that
$$\limsup_{\ell,i,N\uparrow\infty}
\EE^{\bar{k}_N^{(\ell,i)}}\int_0^T\fr{(\bar{k}_N^{(\ell,i)})^d}\sum_{x\in\T_{\bar{k}_N^{(\ell,i)}}^d}\Psi_M\big(\eta_t^{\bar{m}_\ell}(x)^{[\bar{k}_N^{(\ell,i)}\bar{\ee}_i^{(\ell)}]}\big)\df t=\int\lls\Psi_M(\Lambda),\bs\pi\rrs\df \bs{Q}(\bs\pi).$$

So for the proof of the lemma it suffices to show that for any $\bs{Q}\in\bs{\mathcal{Q}}^{\infty,\infty,0}=\Lim_{\ell\uparrow,\ee\downarrow 0,N\uparrow\infty}\bs{Q}^{N,\ell,\ee}$
$$\lim_{M\ra+\infty}\int\lls\Psi_M(\Lambda),\bs\pi\rrs\df \bs{Q}(\bs\pi)=0.$$
So let $\bs{Q}\in\bs{\mathcal{Q}}^{\infty,\infty,0}$ and let us define the map $I_M\colon L_{w^*}^\infty(0,T;\bbar{\MMM}_{1,+}(\T^d\x\RR_+))\to\RR$ by 
$I_M(\bs\pi)=\lls\Psi_M(\Lambda),\bs\pi\rrs$. Since $\Psi_M\in C_1(\RR_+)$ we have that $I_M(\bs\pi)=I_M(\widehat{\bs\pi})$. Since $\Psi_M\equiv 0$ on $[0,M-1]$ we have that $\lim_{M\uparrow\infty}\Psi_M=0$ pointwise and thus, since $\Psi_M(\Lambda)\leq\Lambda$ and $\bs\rho_{\bs\pi_t}:=j^*(\widehat{\bs\pi})\in\MMM_{1,+}(\T^d\x\RR_+)$ is a measure, it follows by the dominated convergence theorem that for almost $t\in[0,T]$
$$\lim_{M\ra+\infty}\ls\Psi_M(\Lambda),\bs\pi_t\rs=\lim_{M\ra+\infty}\int\Psi_M(\Lambda)\df\bs\rho_{\bs\pi_t}=0.$$
Next, since for $\bs{Q}(L_{w^*}^\infty(0,T;\bbar{\Y}_{1,\mathfrak{m}}(\T^d))=1$ and $\Psi_M(\Lambda)\leq\Lambda$ it follows that for $\bs{Q}$-a.s.~all paths $\bs\pi$ 
	$$0\leq\ls\Psi_M(\Lambda),\bs\rho_{\bs\pi_t}\rs\leq\ls\Psi_M(\Lambda),\bs\pi_t\rs\leq\ls \Lambda,\bs\pi_t\rs =\mathfrak{m},\quad\mbox{a.s.~}\forall t\in[0,T].$$ Thus by the dominated convergence theorem once again we obtain that 
$$\lim_{M\ra+\infty}I_M(\bs\pi)=\lim_{M\ra+\infty}\int_0^T\ls\Psi_M(\Lambda),\bs\pi_t\rs\df t=0$$
for $\bs{Q}$-a.s.~all paths $\bs\pi$. Since on the $L_{w^*}^\infty(0,T;\Y_{1,\mathfrak{m}}(\TT^d))$ it holds that $I_M(\bs\pi)=\lls\Psi_M(\Lambda),\bs\pi\rrs\leq\lls\Lambda,\bs\pi\rrs\leq T\mathfrak{m}$ and this sets supports the law $\bs{Q}$, by one more application of the dominated convergence theorem we obtain that 
$$\lim_{M\ra+\infty}\int\lls\Psi_M(\Lambda),\bs\pi\rrs\df\bs{Q}(\bs\pi)=\lim_{M\ra+\infty}\int I_M(\bs\pi)\df\bs{Q}(\bs\pi)=0$$
and the proof of the lemma is complete.$\hfill\Box$\\


In order to complete the proof of (c) it remains to show that if the two-blocks estimate holds along the subfamily $\bbar{\bs{Q}}^{N,\ell,\ee}_*:=\bbar{\bs{Q}}^{k_N^{(\ell)},m_\ell,\ee}$ then for any subfamily 
\begin{equation}\label{SubFamilyForConverseProof}
\{(k_{k_N^{(1;\ell,i)}}^{(m_\ell^{(1)})},m_{m_\ell^{(1)}},\ee_i^{(1;\ell)})\}=:\{(\bar{k}_N^{(1;\ell,i)},\bar{m}_\ell^{(1)},\bar{\ee}_i^{(1;\ell)})\}
\end{equation} of $\{(k_N^{(\ell)},m_\ell,\ee)\}$ there exists a further subfamily $\{(\bar{k}_N^{(\ell,i)},\bar{m}_\ell,\bar{\ee}_i^{(\ell)})\}$ of the form~\eqref{Nightmare} such that~\eqref{ToProveReplacLemma} holds. By equality~\eqref{BothPosAgain}, for all $(N,\ell,\ee,M)\in\NN\x\ZZ_+\x(0,\infty)^2$
$$\int_0^T\fr{N^d}\sum_{x\in\T_N^d}[(\eta_t^\ell(x)-M)^+]^{[N\ee]}\1_{[0,M]}(\eta_t^\ell(x)^{[N\ee]})\df t\leq\lls 1,\rho_{\bs\pi^{N,\ell;M;\ee}}^\perp-\rho_{\bs\pi^{N,\ell,\ee;M}}^\perp\rrs.$$ Therefore if we can show that 
\begin{equation}\label{ForTBCConv}
\lim_{M,\ell,i,N\uparrow\infty}\EE^{\bar{k}_N^{(\ell,i)}}
\lls 1,\rho_{\bs\pi^{\bar{k}_N^{(\ell,i)},\bar{m}_\ell;M;\bar{\ee}_i^{(\ell)}}}^\perp-\rho_{\bs\pi^{\bar{k}_N^{(\ell,i)},\bar{m}_\ell,\bar{\ee}_i^{(\ell)};M}}^\perp\rrs=0.
\end{equation}
it follows then that~\eqref{ToProveReplacLemma} holds along the subfamily $\{(\bar{k}_N^{(\ell,i)},\bar{m}_\ell,\bar{\ee}_i^{(\ell)})\}$. For the proof of~\eqref{SubFamilyForConverseProof} we need the following variant of Proposition~\ref{LimitOfTruncatedLaw} for the case of the joint laws $\bbar{\bs{Q}}^{N,\ell,\ee}:=(\bs\pi^{N,\ell},\bs\pi^{N,\ee})_\sharp P^N$.
\begin{prop}\label{LimitOfTruncatedJointLaw} Let $\widehat{D},D^\perp\colon L_{w^*}^\infty(0,T;\bbar{\MMM}_1(\T^d\x\RR_+))\to L_{w^*}^\infty(0,T;\bbar{\MMM}_1(\T^d\x\RR_+))$ be the regular and singular decomposition operators. Let $\{m_\ell\}_{\ell=1}^\infty$ be a diverging sequence, let $\{\ee_i^{(\ell)}\}_{i=1}^\infty$ be sequences converging to $0$ for all $\ell\in\NN$ and let $\{k_N^{(\ell,i)}\}_{N\in\NN}$, $\ell,i\in\NN$ be diverging sequences such that the iterated limit
	\[\lim_{\ell\ra+\infty}\lim_{i\ra+\infty}\lim_{N\ra+\infty}\bbar{\bs{Q}}^{k_N^{(\ell,i)},m_\ell,\ee_i^{(\ell)}}=:\bbar{\bs{Q}}_*\]
	exists and set 
	\[\bbar{\bs{Q}}^{N,\ell,\ee;M}:=(\bs\pi^{N,\ell;M},\bs\pi^{N,\ee;M})_\sharp P^N\in\PP\big(L_{w^*}^\infty(0,T;\bbar{\PP}_1(\T^d\x\RR_+))^2\big).\]
Then \[\lim_{M\ra+\infty}\lim_{\ell\ra+\infty}\lim_{i\ra+\infty}\lim_{N\ra+\infty}(\widehat{D}\x\widehat{D})_\sharp\bbar{\bs{Q}}^{k_N^{(\ell,i)},m_\ell,\ee_i^{(\ell)};M}=
(\widehat{D}\x\widehat{D})_\sharp\bbar{\bs{Q}}_*\]
	and 
	\begin{equation}\label{SingApproxM}\lim_{M\ra+\infty}\lim_{\ell\ra+\infty}\lim_{i\ra+\infty}\lim_{N\ra+\infty}(D^\perp\x D^\perp)_\sharp\bbar{\bs{Q}}^{k_N^{(\ell,i)},m_\ell,\ee_i^{(\ell)};M}=
	(D^\perp\x D^\perp)_\sharp\bbar{\bs{Q}}_*.
	\end{equation}
\end{prop}\textbf{Proof} Recall that $\Pi_M\colon\bbar{C}_1(\T^d\x\RR_+)\to C_1(\T^d\x\RR_+)$, $M>0$, denotes the bounded linear operator defined by $\Pi_MF(u,\lambda)=F(u,\lambda\mn M)$. We also denote by $\Pi_M$ the induced operator on the corresponding $L^1$-spaces. Then the adjoint $\Pi_M^*\colon L_{w^*}^\infty(0,T;\MMM_1(\T^d\x\RR_+))\to L_{w^*}^\infty(0,T;\bbar{\MMM}_1(\T^d\x\RR_+))$ is bounded and $w^*$-continuous. By~\eqref{Trunc} 
\begin{equation}\label{TruncJoint} \widehat{D}\x\widehat{D}\circ(\bs\pi^{N,\ell;M},\bs\pi^{N,\ee;M})=\big((\Pi_M^*\circ j^*)\x(\Pi_M^*\circ j^*)\big)\circ(\bs\pi^{N,\ell},\bs\pi^{N,\ee})
\end{equation}
which yields 
\begin{equation}\label{ApproxFirstCoordJoint}
(\widehat{D}\x\widehat{D})_\sharp\bbar{\bs{Q}}^{N,\ell,\ee;M}=\big((\Pi_M^*\circ j^*)\x(\Pi_M^*\circ j^*)\big)_\sharp\bbar{\bs{Q}}^{N,\ell,\ee}
\end{equation}
for all $(N,\ell,\ee,M)\in\NN\x\ZZ_+\x(0,\infty)^2$. Thus since the map $\Pi_M^*\circ j^*$ is $(w^*,w^*)$-continuous we have by~\eqref{ApproxFirstCoordJoint} and the assumption that 
\begin{align*}
\lim_{\ell,i,N\uparrow\infty}(\widehat{D}\x\widehat{D})_\sharp\bs{Q}^{k_N^{(\ell,i)},m_\ell,\ee_i^{(\ell)};M}&=\big((\Pi_M^*\circ j^*)\x(\Pi_M^*\circ j^*)\big)_\sharp \lim_{\ell,i,N\uparrow\infty}\bs{Q}^{k_N^{(\ell,i)},m_\ell,\ee_i^{(\ell)}}\\
&=\big((\Pi_M^*\circ j^*)\x(\Pi_M^*\circ j^*)\big)_\sharp\bbar{\bs{Q}}_*.
\end{align*} 
Since $\Pi_M^*\circ j^*$ $w^*$-converges pointwise to $\widehat{D}$ as we have seen in the proof of Corollary~\ref{PathMeasDecomp} the map $(\Pi_M^*\circ j^*)_\sharp$ converges pointwise to $\widehat{D}_\sharp$ on $\PP L_{w^*}^\infty(0,T;\bbar{\MMM}_1(\T^d\x\RR_+))$
and therefore
$$\lim_{M\ra+\infty}\big((\Pi_M^*\circ j^*)\x(\Pi_M^*\circ j^*)\big)_\sharp\bbar{\bs{Q}}_*=(\widehat{D}\x\widehat{D})_\sharp\bbar{\bs{Q}}_*.$$

For the second limit, recalling that $T_M\colon C(\T^d)\to\bbar{C}_1(\T^d\x\RR_+)$ is the operator $T_Mf=(\Lambda-M)^+f(U)$,
\begin{equation}\label{Trunc2Joint} (D^\perp\x D^\perp)\circ(\bs\pi^{N,\ell;M},\bs\pi^{N,\ee;M})=\big((R^*\circ T_M^*)\x(R^*\circ T_M^*)\big)\circ(\bs\pi^{N,\ell},\bs\pi^{N,\ee})\end{equation} which shows that 
$$
(D^\perp\x D^\perp)_\sharp\bbar{\bs{Q}}^{N,\ell,\ee;M}=\big((R^*\circ T_M^*)\x(R^*\circ T_M^*)\big)_\sharp\bbar{\bs{Q}}^{N,\ell,\ee}$$
and as we have also seen in the proof of Corollary~\ref{PathMeasDecomp} the maps $R^*\circ T_M^*$ $w^*$-converge pointwise to $D^\perp$ as $M\ra+\infty$ and the second lmit follows as the first limit.$\hfill\Box$\\

So let $\{(\bar{k}_N^{(1;\ell,i)},\bar{m}_\ell^{(1)},\bar{\ee}_i^{(1;\ell)})\}$ be a subfamily of $\{(k_N^{(\ell)},m_\ell;\ee)\}$ as in~\eqref{SubFamilyForConverseProof}. Then since the family $\{\bbar{\bs{Q}}^{\bar{k}_N^{(1;\ell,i)},\bar{m}_\ell^{(1)},\bar{\ee}_i^{(1;\ell)}}\}$ is relatively compact there exists a further subfamily $\{(\bar{k}_N^{(\ell,i)},\bar{m}_\ell,\bar{\ee}_i^{(\ell)})\}$ of the form~\eqref{Nightmare} such that the iterated limit 
	\begin{equation}\label{ItLim}
	\bbar{\bs{Q}}_*:=\lim_{\ell,i,N\uparrow}\bbar{\bs{Q}}^{\bar{k}_N^{(\ell,i)},\bar{m}_\ell,\bar{\ee}_i^{(\ell)}}
	\end{equation} exists. Then by the assumed validity of the two-blocks estimate, for all maps $\Psi\in\bbar{C}_{1,\uparrow}(\RR_+)$ and all $G\in L^1(0,T;C(\T^d))$
$$\int|\lls G,B_\Psi(\bs\pi^\infty-\bs\pi^0)\rrs|\df\bbar{\bs{Q}}_*(\bs\pi^\infty,\bs\pi^0)=0.$$ 
By applying this equality to the maps $\Psi_M(\cdot):=\Psi(\cdot\mn M)$, $M>0$, for some $\Psi\in\bbar{C}_{1,\uparrow}(\RR_+)$
\begin{align*}
0&=\int|\lls G,B_{\Psi_M}(\bs\pi^\infty-\bs\pi^0)\rrs|\df\bbar{\bs{Q}}_*(\bs\pi^\infty,\bs\pi^0)\\
&=\int|\lls G(U)\Psi(\Lambda\mn M),\bs\pi^\infty-\bs\pi^0)\rrs|\df\bbar{\bs{Q}}_*(\bs\pi^\infty,\bs\pi^0)\\
&=\int\big|\big\lls (j\circ\Pi_M)\big(G(U)\Psi(\Lambda)\big),\bs\pi^\infty-\bs\pi^0\rrs\big|\df\bbar{\bs{Q}}_*(\bs\pi^\infty,\bs\pi^0)\\
&=\int\big|\big\lls \big(G(U)\Psi(\Lambda)\big),(j\circ\Pi_M)^*(\bs\pi^\infty-\bs\pi^0)\rrs\big|\df\bbar{\bs{Q}}_*(\bs\pi^\infty,\bs\pi^0)
\end{align*}
Since $\Pi_M^*\circ j^*$ $w^*$-converges pointwise to $\widehat{D}$ and the map $(\bs\pi^\infty,\bs\pi^0)\mapsto\big\lls \big(G(U)\Psi(\Lambda)\big),(j\circ\Pi_M)^*(\bs\pi^\infty-\bs\pi^0)\big)\rrs$ is dominated by an $L^1(\bbar{\bs{Q}}_*)$-function, by taking the limit as $M\uparrow\infty$ we obtain by the dominated convergence theorem that
\begin{align}\label{ForFluid}
0&=\lim_{M\uparrow\infty}\int\big|\big\lls \big(G(U)\Psi(\Lambda)\big),(j\circ\Pi_M)^*(\bs\pi^\infty-\bs\pi^0)\big)\rrs\big|\df\bbar{\bs{Q}}_*(\bs\pi^\infty,\bs\pi^0)\nonumber\\
&=\int\big|\big\lls\big(G(U)\Psi(\Lambda)\big),\widehat{D}(\bs\pi^\infty-\bs\pi^0)\rrs\big|\df\bbar{\bs{Q}}_*(\bs\pi^\infty,\bs\pi^0)\nonumber\\
&=\int\big|\big\lls G,B_\Psi\circ\widehat{D}(\bs\pi^\infty-\bs\pi^0)\rrs\big|\df\bbar{\bs{Q}}_*(\bs\pi^\infty,\bs\pi^0)
\end{align}
By the comparison of regular parts in statement (ii) of part (b)
$$|\lls G,B_\Psi\circ\widehat{D}(\bs\pi^\infty-\bs\pi^0)\rrs|=\lls G,B_\Psi\circ\widehat{D}(\bs\pi^0-\bs\pi^\infty)\rrs.$$
By the equation $B_\Psi=B_\Psi\circ\widehat{D}+\Psi'(\infty)B\circ D^\perp$ the equality above becomes
\begin{align*}
|\lls G,B\circ\widehat{D}(\bs\pi^\infty-\bs\pi^0)\rrs|&=\lls G,B_\Psi\circ\widehat{D}(\bs\pi^0-\bs\pi^\infty)\rrs=\Psi'(\infty)\lls G, B\circ D^\perp(\bs\pi^\infty-\bs\pi^0)\rrs\\
&=\Psi'(\infty)|\lls G, B\circ D^\perp(\bs\pi^\infty-\bs\pi^0)\rrs|
\end{align*}
and therefore if $\Psi'(\infty)\neq 0$ equality~\eqref{ForFluid} yields 
$$0=\int|\lls G,B\circ D^\perp(\bs\pi^\infty-\bs\pi^0)\rrs|\df\bbar{\bs{Q}}_*(\bs\pi^\infty,\bs\pi^0)$$
for all $G\in L^1(0,T;C_+(\T^d))$. By the linearity of $D^\perp$ this equivalent to 
$$0=\int|\lls G,B(\bs\s^\infty-\bs\s^0)\rrs|\df(D^\perp\x D^\perp)_\sharp\bbar{\bs{Q}}_*(\bs\s^\infty,\bs\s^0).$$
But by the limit~\eqref{ItLim} and Proposition~\ref{LimitOfTruncatedJointLaw}, the limit~\eqref{SingApproxM} holds and thus by the continuity of the map $(\bs\s^\infty,\bs\s^0)\mapsto\lls G,B(\bs\s^\infty-\bs\s^0)\rrs$ and the portmanteau theorem we obtain that 
\begin{align*}0&=\lim_{M,\ell,i,N\uparrow\infty}\int|\lls G,B(\bs\s^\infty-\bs\s^0)\rrs|\df(D^\perp\x D^\perp)_\sharp\bbar{\bs{Q}}_*^{\bar{k}_N^{(\ell,i)},\bar{m}_\ell,\bbar{\ee}_i^{(\ell)};M}(\bs\s^\infty,\bs\s^0)\\
&=\lim_{M,\ell,i,N\uparrow\infty}\int|\lls G,B\circ D^\perp (\bs\pi^\infty-\bs\pi^0)\rrs|\df\bbar{\bs{Q}}_*^{\bar{k}_N^{(\ell,i)},\bar{m}_\ell,\bbar{\ee}_i^{(\ell)};M}(\bs\pi^\infty,\bs\pi^0)\end{align*}
and therefore 
\begin{equation}\label{TheLimWeGetFromTBE}
\lim_{M,\ell,i,N\uparrow\infty}\EE^{\bar{k}_N^{(\ell,i)}}|\lls G,B\circ D^\perp (\bs\pi^{\bar{k}_N^{(\ell,i)},\bar{m}_\ell;M}-\bs\pi^{\bar{k}_N^{(\ell,i)},\bbar{\ee}_i^{(\ell)};M})\rrs|=0,\quad\forall G\in L^1(0,T;C_+(\T^d)).
\end{equation}
Then recalling the functional equation $B\circ D^\perp(\bs\pi)=\rho_{\bs\pi}^\perp$ from~\eqref{SingularRepresentativeBiaBerProj} we have that 
\begin{align}\label{ToComplSoonTBC}
\lls 1,\rho_{\bs\pi^{N,\ell;M;\ee}}^\perp-\rho_{\bs\pi^{N,\ell,\ee;M}}^\perp\rrs&=\lls 1,B\circ D^\perp(\bs\pi^{N,\ell;M;\ee}-\bs\pi^{N,\ell,\ee;M})\rrs\nonumber\\
&=\lls 1,B\circ D^\perp(\bs\pi^{N,\ell;M;\ee}-\bs\pi^{N,\ell;M})\rrs+\lls 1,B\circ D^\perp(\bs\pi^{N,\ell;M}-\bs\pi^{N,\ee;M})\rrs\nonumber\\
&\qquad+\lls 1,B\circ D^\perp(\bs\pi^{N,\ee;M}-\bs\pi^{N,\ell,\ee;M})\rrs.
\end{align}
The first term in the right hand side above vanishes by a change in the order of summation, i.e.
\begin{align*}
\lls 1,B\circ D^\perp(\bs\pi^{N,\ell;M;\ee}-\bs\pi^{N,\ell;M})\rrs&=\int_0^T\fr{N^d}\sum_{x\in\T_N^d}\big\{(\eta_t^\ell(x)-M)^{+[N\ee]}-(\eta^\ell_t(x)-M)^+\big\}\df t\\
&=\int_0^T\fr{N^d}\sum_{x\in\T_N^d}\Big\{\sum_{|y|\leq[N\ee]}(\eta_t^\ell(x+y)-M)^+-(\eta^\ell_t(x)-M)^+\Big\}\df t\\
&=0
\end{align*}
and the absolute value of the third term in the right hand side of~\eqref{ToComplSoonTBC} is bounded above by
\begin{align*}
|\lls 1,B\circ D^\perp(\bs\pi^{N,\ee;M}-\bs\pi^{N,\ell,\ee;M})\rrs|&=\Big|\int_0^T\fr{N^d}\sum_{x\in\T_N^d}\big\{(\eta_t^{[N\ee]}(x)-M)^+-(\eta_t^\ell(x)^{[N\ee]}-M)^+\big\}\df t\Big|\\
&\leq\int_0^T\fr{N^d}\sum_{x\in\T_N^d}\big|\eta_t^{[N\ee]}(x)-\eta_t^\ell(x)^{[N\ee]}\big|\df t
\end{align*}
and this last term converges to zero in the limit as $N\ra\infty$ by the bound~\eqref{ConsecAver} on the difference of consecutive averages. 
But since $\{(\bar{k}_N^{(\ell,i)},\bar{m}_\ell,\bar{\ee}_i^{(\ell)})\}$ is a subfamily of the original family $\{(k_N^{(\ell)},m_\ell,\ee)\}$ by the limit~\eqref{TDBE} we have that 
$$\lim_{M,\ell,i,N\uparrow\infty}\EE^{\bar{k}_N^{(\ell,i)}}|\lls G,B\circ D^\perp (\bs\pi^{\bar{k}_N^{(\ell,i)},\bar{m}_\ell;M}-\bs\pi^{\bar{k}_N^{(\ell,i)},\bar{m}_\ell;M;\bbar{\ee}_i^{(\ell)}})\rrs|=0$$
and thus it follows by the limit~\eqref{TheLimWeGetFromTBE} and inequality~\eqref{ToComplSoonTBC} that~\eqref{ForTBCConv} holds and thus the proof of Theorem~\ref{TBCTheorem} is complete.$\hfill\Box$

\subsection{On the replacement lemma}\label{RLSection}

Let us first check that for each $\Psi\in\bbar{C}_1(\RR_+)$ the map $I_\Psi$ defined in~\eqref{IPsi} is well defined and Borel measurable. In the definition of this map we consider the $L_{w^*}^\infty$-space of path-measures as the target space for the map $I_\Psi$, since due to the fact that the map $\MMM_+(\T^d)\ni\pi\mapsto\pi^{ac}\in\MMM_{+,ac}(\T^d)$ is not weakly continuous, the map $t\mapsto\Psi(\pi^{ac}_t)$ need not be a cadlag path, even when $\pi$ is cadlag. However as we will see the map $I_\Psi$ is well-defined and Borel measurable when viewed as taking values in $L_{w^*}(0,T;\MMM_+(\T^d))$. Indeed, if $\Psi\in\bbar{C}_1(\RR_+)$ then $$\|\Psi(\pi^{ac})\df\LL_{\T^d}\|_{TV}=\|\Psi(\pi^{ac})\|_{L^1(\T^d)}\leq\|\Psi\|_{\infty,1}\int_{\T^d}(1+\pi^{ac}(u))\df u\leq\|\Psi\|_{\infty,1}(1+\|\pi\|_{TV})$$
for all $\pi\in\MMM_+(\T^d)$, and thus $\Psi(\pi^{ac})\df\LL_{\T^d}$ is a finite measure. To see that the induced map $I_\Psi$ is well-defined and Borel measurable, for each $\ee\in(0,1/2)$, $M<+\infty$ we consider the map $$D(0,T;\MMM_+(\T^d))\ni\pi\mapsto I_\Psi^{\ee,M}(\pi):=\Psi\big((\pi*\iota_\ee)\mn M\big)\df\LL_{\T^d}\in D(0,T;\MMM_+(\T^d))$$
where $\pi*\iota_\ee$ is the convolution  
$$\pi*\iota_\ee(u)=\int\iota_\ee(\y-u)\df\pi(\y)$$ and $(\iota_\ee)_{0<\ee<1/2}\subs C(\T^d)\subs\PP\T^d$ is an approximation of the Dirac measure $\delta_0$, i.e.~the family $\{\iota_\ee\df\LL_{\T^d}\}\subs\PP\T^d$ converges weakly to $\delta_0$ as $\ee\ra 0$. Let us check here that the map $I_\Psi^{\ee,M}$ is continuous. Since a continuous map $f\colon M\to N$ induces a continuous map $f$ on the corresponding Skorohod spaces via $f(\mu)(t)=f(\mu(t))$ if we show that the map $\MMM_+(\T^d)\ni\pi\mapsto\Psi\big((\pi*\iota_\ee)\mn M\big)\df\LL_{\T^d}\in\MMM_{+,ac}(\T^d)$ is continuous then the induced map $I_\Psi^{\ee,M}$ on the Skorohod spaces will be continuous. The fact that $\MMM_+(\T^d)\ni\pi\mapsto\Psi\big((\pi*\iota_\ee)\mn M\big)$ is continuous follows from the presence of the convolution with the continuous function $\iota_\ee$, which strengthens weak convergence to convergence in total variation, i.e.~to convergence of the densities in $L^1(\T^d)$. Indeed, sicnce $\iota_\ee\in C(\T^d)$, for each $\ee>0$ and all $u\in\T^d$
\[\lim_{n\ra+\infty}\pi_n*\iota_\ee(u)=\lim_{n\ra+\infty}\int\iota_\ee(\y-u)\df\pi_n(\y)=\int\iota_\ee(\y-u)\df\pi(\y)=\pi*\iota_\ee(u).\]
Therefore if $\Psi\colon\RR_+\to\RR_+$ is a continuous then $\{\Psi(\pi_n*\iota_\ee)\}_{n\in\NN}$ converges pointwise in $\T^d$ to $\Psi(\pi*\iota_\ee)$. Moreover, since  $\{\pi_n\}_{n\in\NN}\subs\MMM_+(\T^d)$ converges to weakly to $\pi\in\MMM_+(\T^d)$ we have that $$C:=\sup_{n\in\NN}\|\pi_n\|_{TV}<+\infty$$ so that $|\pi_n*\iota_\ee(u)|\leq C\|\iota_\ee\|_\infty$ and
\[|\Psi(\pi_n*\iota_\ee(u))|\leq \sup_{0\leq\lambda\leq C\|\iota_\ee\|_\infty}\Psi(\lambda)<+\infty,\quad\forall u\in\T^d,\;n\in\NN\] and thus by the bounded convergence theorem we obtain that
\[\lim_{n\ra+\infty}\int|\Psi(\pi_n*\iota_\ee(u))-\Psi(\pi*\iota_\ee(u))|\df u=0.\]
This shows that for any continuous map $\Psi\colon\RR_+\to\RR_+$ the map $$\MMM_+(\T^d)\ni\pi\mapsto\Psi(\pi*\iota_\ee)\df\LL_{\T^d}\in(\MMM_{+,ac},\|\cdot\|_{TV})\cong L^1(\T^d),$$ is continuous with respect to the strong topology in the target space $\MMM_+(\T^d)$. Of course by applying this to the map $\Psi(\cdot\mn M)$ for a given map $\Psi\in C(\RR_+)$ we obtain that the map $\MMM_+(\T^d)\ni\pi\mapsto\Psi((\pi*\iota_\ee)\mn M)\df\LL_{\T^d}$ is continuous. Thus the induced map $$I_\Psi^{\ee,M}\colon D(0,T;\MMM_+(\T^d))\to D(0,T;\big(\MMM_+(\T^d),\|\cdot\|_{TV}\big)$$ on the Skorohod spaces is continuous. Since by Proposition~\ref{SkorohodToLEmbed} the natural injection $$D\big(0,T;(\MMM_+(\T^d),\|\cdot\|_{TV})\big)\hookrightarrow L_{w^*}^\infty(0,T;\MMM_+(\T^d))$$ is continuous, the map $I_\Psi^{\ee,M}$ is also continuous for the target space $L_{w^*}^\infty(0,T;\MMM_+(\T^d))$ equipped with the normed $w^*$-topology, i.e.~for any sequence $\{\pi_n\}\subs L_{w^*}^\infty(0,T;\MMM_+(\T^d))$ $w^*$-converging to $\pi\in L_{w^*}^\infty(0,T;\MMM_+(\T^d))$ we have 
$$\lim_{n\ra+\infty}\int_0^Tf(t)\big\|\Psi\big((\pi_{n,t}*\iota_\ee)\mn M\big)-\Psi\big((\pi_t*\iota_\ee)\mn M\big)\big\|_{L^1(\T^d)}\df t=0,\quad\forall f\in L^1(0,T).$$

We show next that for an appropriate choice of the approximation $(\iota_\ee)_{0<\ee<\fr{2}}$ of the convolution identity $\delta_0$, the maps $I_\Psi^{\ee,M}\colon D(0,T;\MMM_+(\T^d))\to L_{w^*}^\infty(0,T;\MMM_+(\T^d))$ converge pointwise as $\ee\ra 0$ and then $M\ra+\infty$ to the map $I_\Psi$ defined in~\eqref{IPsi}. For this we will use the following.
\begin{lemma}\label{AsToAbsDens} Let $\pi\in\MMM_+(\T^d)$, let $\LL$ be a reference measure (i.e.~the Lebesgue measure) on $\T^d$ and let $\pi=\pi^{ac}+\pi^\perp$ be the Radon-Nikodym decomposition of $\pi$ with respect to $\LL$. Then for Lebesgue almost all $u\in\T^d$ 
	$$\lim_{\ee\downarrow 0}\frac{\pi\big(u+[-\ee,\ee]^d\big)}{(2\ee)^d}=\frac{\df\pi^{ac}}{\df\LL}(u).$$
\end{lemma}\textbf{Proof} For each $\mu\in\MMM_+(\T^d)$ we define the maps 
$$D^-\mu(u)=\liminf_{\ee\downarrow 0}\frac{\mu\big(u+[-\ee,\ee]^d\big)}{(2\ee)^d},\quad D^+\mu(u)=\limsup_{\ee\downarrow 0}\frac{\mu\big(u+[-\ee,\ee]^d\big)}{(2\ee)^d},\quad u\in\T^d.$$
These maps are Borel measurable (see for example~\cite[Lemma 5.9.1]{DiBenedetto2016a}) for any $\mu\in\MMM_+(\T^d)$ and by \cite[Proposition 5.9.1]{DiBenedetto2016a} the set 
$$\mathcal{E}_\mu:=\big\{u\in\T^d\bigm|D^-\mu(u)<D^+\mu(u)\mbox{ or }D^+\mu(u)=+\infty\big\}$$ is $\LL$-null Borel set. Obviously on the complement of $\E_\mu$ the limit 
$$D\mu(u):=\lim_{\ee\ra 0}\frac{\mu\big(u+(-\ee,\ee)^d\big)}{(2\ee)^d}=D^-\mu(u)=D^+\mu(u),\quad u\in\T^d\sm\E_\mu$$ exists and is finite. By the decomposition $\pi=\pi^{ac}+\pi^\perp$ and the additivity properties of the limit inferior and limit superior
$$D^-\pi^{ac}+D^-\pi^\perp\leq D^-\pi\leq D^+\pi\leq D^+\pi^{ac}+D^+\pi^\perp,$$
and on the set $\T^d\sm(\E_{\pi^{ac}}\cup\E_{\pi^\perp})$ all the inequalities above are equalities. Since $\E_{\pi^{ac}}\cup\E_\pi^{\perp}$ is a $\LL$-null set $D\pi=D\pi^{ac}+D\pi^\perp$ $\LL$-a.s.~and thus the claim follows by~\cite[Proposition 5.10.2]{DiBenedetto2016a} according to which $D\pi^\perp=0$ $\LL$-a.s.$\hfill\Box$\\

We describe now a particular choice of an approximation $(\iota_\ee)_{0<\ee<1/2}$ of the identity $\delta_0\in\PP\T^d$ for which it is easy to verify that \begin{equation}\label{GoodConv}\liminf_{\ee\downarrow 0}\frac{\pi\big(u+[-\ee,\ee]^d\big)}{(2\ee)^d}\leq\liminf_{\ee\downarrow 0}\pi*\iota_\ee(u)\leq\limsup_{\ee\downarrow 0}\pi*\iota_\ee(u)=\limsup_{\ee\downarrow 0}\frac{\pi\big(u+[-\ee,\ee]^d\big)}{(2\ee)^d}.\end{equation}
By Lemma~\ref{AsToAbsDens} this implies that 
\begin{equation}\label{AsPointWConv}\lim_{\ee\ra 0}\pi*\iota_\ee(u)=\frac{\df\pi^{ac}}{\df\LL_{\T^d}}(u),\quad\LL_{\T^d}\mbox{a.s.~for all }u\in\T^d.\end{equation}
Of course, if we could take $\iota_\ee$ to be $k_\ee:=\fr{(2\ee)^d}\1_{[-\ee,\ee]^d}$, $\ee\in(0,\fr{2})$, then we would have
$$\pi*k_\ee(u)=\int k_\ee(\y-u)\df\pi(\y)=\fr{(2\ee)^d}\int\1_{[-\ee,\ee]^d}(\y-u)\df\pi(\y)=\frac{\pi(u+[-\ee,\ee]^d)}{(2\ee)^d}$$
and~\eqref{GoodConv} would trivially hold, but we want the want the maps $\iota_\ee$, $\ee\in(0,1/2)$, to be continuous. For this reason we consider a continuous map $j_\ee\in C_c((\T\sm\{\fr{2}\})^d))$ such that 
$$\1_{[-\frac{1-\ee}{2},\frac{1-\ee}{2}]^d}\leq j_\ee\leq\1_{(-\fr{2},\fr{2})^d},$$ we set $\ls j_\ee\rs:=\int_{\T^d}j_\ee(u)\df u$ and $\bar{j}_\ee:=\ls j_\ee\rs^{-1}j_\ee$ and define $\iota_\ee\in C_c((\T\sm\{\fr{2}\})^d)$ by 
\begin{equation}\label{ConvDef}\iota_\ee(u):=\fr{(2\ee)^d}\bar{j}_\ee\Big(\frac{u}{2\ee}\Big),\quad u\in\T^d,\;\ee\in(0,\fr{2}).\end{equation} Then $(1-\ee)^d\leq\ls j_\ee\rs\leq1$ so that $\lim_{\ee\ra 0}\ls j_\ee\rs=1$ and 
$$\fr{\ls j_\ee\rs(2\ee)^d}\1_{[-\frac{1-\ee}{2},\frac{1-\ee}{2}]^d}\Big(\frac{u}{2\ee}\Big)\leq\iota_\ee(u)\leq\fr{\ls j_\ee\rs(2\ee)^d}\1_{(-\fr{2},\fr{2})^d}\Big(\frac{u}{2\ee}\Big),$$
which shows that 
$$\frac{(2\ee(1-\ee))^d}{\ls j_\ee\rs(2\ee)^d}k_{\ee-\ee^2}(u)\leq\iota_\ee(u)\leq\fr{\ls j_\ee\rs(2\ee)^d}\1_{(-\ee,\ee)^d}(u)\leq\fr{\ls j_\ee\rs}k_\ee(u).$$
Consequently 
$$\frac{(1-\ee)^d}{\ls j_\ee\rs}\pi*k_{\ee-\ee^2}(u)\leq\pi*\iota_\ee(u)\leq\fr{\ls j_\ee\rs}\pi*k_\ee(u).$$
Since $(1-\ee)^d\leq\ls j_\ee\rs\leq 1$ it follows that  
\begin{align}\label{e}
(1-\ee)^d\pi*k_{\ee-\ee^2}(u)\leq\pi*\iota_\ee(u)\leq\fr{(1-\ee)^d}\pi*k_\ee(u).
\end{align}
Since $\pi*k_\ee(u)=\fr{(2\ee)^d}\pi(u+[-\ee,\ee]^d)$, sending $\ee$ to $0$ we obtain~\eqref{GoodConv}.

It follows that for the approximation $(\iota_\ee)$ of the identity defined in~\eqref{ConvDef}
$$\lim_{\ee\downarrow 0}\pi*\iota_\ee(u)=\frac{\df\pi^{ac}}{\df\LL}(u),\quad\mbox{a.s.-}\forall u\in\T^d.$$
Consequently, since for each fixed $M>0$ for any continuous map $\Psi\colon\RR_+\to\RR_+$ we have the trivial bound
\begin{equation}\label{MacroRepDomConvB1}
|\Psi(\pi*\iota_\ee(u)\mn M)-\Psi(\pi^{ac}(u)\mn M)|\leq2\sup_{0\leq\lambda\leq M}|\Psi(\lambda)|<+\infty
\end{equation} for all $u\in\T^d$, all $\ee\in(0,\fr{2})$ and all $\pi\in\MMM_+(\T^d)$, it follows by the bounded convergence theorem that   
$$\lim_{\ee\ra 0}\int_{\T^d}|\Psi(\pi*\iota_\ee(u)\mn M)-\Psi(\pi^{ac}(u)\mn M)|\df u=0,\quad\forall\pi\in\MMM_+(\T^d).$$
Consequently, for any path $\pi\in D(0,T;\MMM_+(\T^d))$ the paths $I_\Psi^{\ee,M}(\pi)$ converge pointwise in $[0,T]$ as $\ee\ra 0$ to the path $I_\Psi^{0,M}(\pi)(t)=\Psi(\pi^{ac}_t\mn M)\df\LL_{\T^d}$ with respect to the strong topology on the target space $\MMM_{+,ac}(\T^d)$. Since the paths $I_\Psi^{\ee,M}(\pi)\colon[0,T]\to(\MMM_{+,ac}(\T^d),\|\cdot\|_{TV})$ are cadlag and thus strongly measurable, the pointwise limit-path $I_\Psi^{0,M}(\pi):[0,T]\to\MMM_{+,ac}(\T^d)$ as $\ee\ra 0$ is also strongly measurable, and thus also $w^*$-measurable. Consequently the map 
$$D(0,T;\MMM_+(\T^d))\ni\pi\mapsto I_\Psi^{0,M}(\pi):=\Psi(\pi^{ac}\mn M)\df\LL_{\T^d}\in L_{w^*}^\infty(0,T;\MMM_{+,ac}(\T^d))$$
is well-defined, and by another application of the dominated convergence theorem it follows that $I_\Psi^{\ee,M}$ converges as $\ee\ra 0$ pointwise in $D(0,T;\MMM_+(\T^d))$ to the map $I_\Psi^{0,M}$ with respect to the normed $w^*$-convergence on the target space, i.e.~for all $\pi\in D(0,T;\MMM_+(\T^d))$
$$\lim_{\ee\ra 0}\int_0^Tf(t)\big\|\Psi\big((\pi_t*\iota_\ee)\mn M\big)-\Psi(\pi_t^{ac}\mn M)\big\|_{TV}\df t=0,\quad\forall f\in L^1(0,T).$$  
Since the maps $I_\Psi^{\ee,M}$, $\ee\in(0,1/2)$, are $w^*$-continuous and by Proposition~\ref{WeakStarAndAuxSalg} the Borel $\s$-algebra of the $w^*$-topology on $L_{w^*}^\infty(0,T;\MMM(\T^d))$ is the Borel $\s$-algebra of a separable metric space, it follows that the maps $I_\Psi^{0,M}$, $M>0$, are $w^*$-Borel measurable.

Since for each $\pi\in\MMM_+(\T^d)$ and $\Psi\in\bbar{C}_1(\T^d)$ it holds that
$$|\Psi(\pi^{ac})-\Psi(\pi^{ac}\mn M)|=|\Psi(\pi^{ac})-\Psi(M)|\1_{\{\pi^{ac}>M\}}\leq2\|\Psi\|_{\infty,1}(1+\pi^{ac})\1_{\{\pi^{ac}>M\}}$$
and $\pi^{ac}\in L^1(\T^d)$ we obtain by the dominated convergence theorem that 
\[\lim_{M\ra+\infty}\|\Psi(\pi^{ac})-\Psi(\pi^{ac}\mn M)\|_{TV}\leq2\|\Psi\|_{\infty,1}\lim_{M\ra+\infty}\int\big(1+\pi^{ac}(u)\big)\1_{\{\pi^{ac}>M\}}(u)\df u=0.\]
for each $\pi\in\MMM_+(\T^d)$. Thus for each path $\pi\in D(0,T;\MMM_+(\T^d))$ the path $I_\Psi^{0,M}(\pi)\colon[0,T]\to\MMM_{+,ac}(\T^d)$ converges as $M\uparrow +\infty$ pointwise in $[0,T]$ to the path $I_\Psi(\pi)(t)=\Psi(\pi^{ac}_t)\df\LL_{\T^d}$ in the strong topology of $\MMM_{+,ac}(\T^d)$. Since the paths $I_\Psi^{0,M}(\pi)$ are strongly Borel measurable it follows that the path $I_\Psi(\pi)$ is strongly measurable and thus also $w^*$-measurable. Since also 
\begin{equation}\|I_\Psi(\pi)\|_{TV;\infty}=\esssup_{0\leq t\leq T}\int|\Psi(\pi_t^{ac}(u))|\df u\leq\|\Psi\|_{\infty,1}(1+\|\pi\|_{TV;\infty})<+\infty
\end{equation}
it follows that the map $I_\Psi$ is well-defined with domain and target space given in~\eqref{IPsi}. Finally, since by the conservation of the total number of particles
\begin{align}\label{MacroRepDomConvB2}
\|\Psi(\pi^{ac})-\Psi(\pi^{ac}\mn M)\|_{TV;\infty}&\leq 2\|\Psi\|_{\infty,1}\esssup_{0\leq t\leq T}\int\big(1+\pi_t^{ac}(u)\big)\1_{\{\pi_t^{ac}>M\}}(u)\df u\nonumber\\
&\leq 2\|\Psi\|_{\infty,1}(1+\|\pi\|_{TV;\infty})<2\|\Psi\|_{\infty,1}(1+\mathfrak{m})<+\infty
\end{align}
$Q^{\infty}_*$-a.s.~for all $\pi\in D(0,T;\MMM_+(\T^d))$, it follows by an application of the dominated convergence theorem that for all $f\in L^1(0,T)$
$$\lim_{M\ra+\infty}\int_0^Tf(t)\|\Psi(\pi_t^{ac})-\Psi(\pi^{ac}_t\mn M)\|_{TV}\df t\leq 2\|\Psi\|_{\infty,1}\int_0^Tf(t)\int\big(1+\pi^{ac}_t(u)\big)\df u\df t=0.$$
In particular $I_\Psi=\lim_{M\ra+\infty}I_\Psi^{0,M}$ $w^*$-converges pointwise to $I_\Psi$, and since each of the maps $I_\Psi^{0,M}$ is $w^*$-Borel measurable and the $w^*$-Borel $\s$-algebra is the Borel $\s$-algebra of a separable metric space it follows that $I_\Psi$ is $w^*$-measurable.

Let $\bbar{\bs{Q}}_\Psi^\infty$ be a limit point of the sequence $\{\bbar{\bs{Q}}_\Psi^N\}$ of the laws defined in~\eqref{RLLawsInitial}. There exists then an increasing sequence $\{k_N\}_{N=1}^\infty\subs\NN$ such that $\bbar{\bs{Q}}_\Psi^\infty=\lim_{N\uparrow\infty}\bbar{\bs{Q}}_\Psi^{k_N}$. If the assumption for the validity of the full replacement lemma as stated in Theorem~\eqref{RLTheorem} holds, then we can assume that $\{k_N\}$ has been chosen so that any subfamily $\{(k_N^{(\ell)},m_\ell)\}_{(N,\ell)}$ of $\{(k_N,\ell)\}_{(N,\ell)}$ has a subfamily, still denoted by $\{(k_N^{(\ell)},m_\ell)\}_{(N,\ell)}$, such that any subfamily $\{(k_N^{(\ell,i)},m_{\ell},\ee_i^{(\ell)})\}$ of $\{(k_N^{(\ell)},m_\ell,\ee)\}$ has a further subfamily $\{(\bar{k}_N^{(\ell,i)},\bar{m}_{\ell},\bar{\ee}_i^{(\ell)})\}$ along which~\eqref{ToProveReplacLemma} holds. Since the family $\bs{Q}^{k_N,\ell}:=\bs\pi^{k_N,\ell}_\sharp P^{k_N}$, $N\in\NN$, $\ell\in\ZZ_+$ is relatively compact there exits an increasing sequence $\{m_\ell^{(0)}\}_{\ell=1}^\infty\subs\ZZ_+$ and subsequences  $\{k_N^{(0;\ell)}\}_{N=1}^\infty$, $\ell\in\NN$, of $\{k_N\}$ such that the iterated limit 
\begin{equation}\label{TBCAssumption}\bs{Q}_*^{\infty,\infty}:=\lim_{\ell\uparrow\infty}\lim_{N\uparrow\infty}\bs{Q}^{k_N^{(0;\ell)},m^{(0)}_\ell}
\end{equation} 
exists. Now, in case the assumption for the replacement lemma holds then by the choice of $\{k_N\}$ we can further assume that $\{k_N^{(0;\ell)},m^{(0)}_\ell\}$ has been chosen so that any subfamily $\{(k_N^{(\ell,i)},m_{\ell},\ee_i^{(\ell)})\}$ of $\{(k_N^{(0;\ell)},m^{(0)}_\ell,\ee)\}$ has a further subfamily $\{(\bar{k}_N^{(\ell,i)},\bar{m}_{\ell},\bar{\ee}_i^{(\ell)})\}$ along which~\eqref{ToProveReplacLemma} holds. We consider then the joint laws 
\[\bs{R}_\Psi^{k_N^{(0;\ell)},m^{(0)}_\ell,\ee}:=(\s^{k_N^{(0;\ell)},\Psi},\bs\pi^{k_N^{(0;\ell)},m^{(0)}_\ell},\bs\pi^{k_N^{(0;\ell)},\ee},\pi^{k_N^{(0;\ell)}})_\sharp P^{k_N^{(0;\ell)}},\quad(N,\ell,\ee)\in\NN^2\x(0,\infty)\]
on the product space $L_{w^*}^\infty(0,T;\MMM(\T^d))\x L_{w^*}^\infty(0,T;\bbar{\PP}_1(\T^d\x\RR_+))^2\x D(0,T;\MMM_+(\T^d)$ and we will denote by $(\s^\Psi,\bs\pi^\infty,\bs\pi^0,\pi)$ the arbitrary element of this product space and with a slight abuse of notation also the natural projections on the coordinates of this product space. The family $\{\bs{R}_\Psi^{k_N^{(0;\ell)},m^{(0)}_\ell,\ee}\}$ is relatively compact and thus there exist a subfamily $\{(k_N^{(\ell,i)},m_\ell,\ee_i^{(\ell)})\}$ of $\{(k_N^{(0;\ell)},m^{(0)}_\ell,\ee)\}$ of the form 
\[\{(k_N^{(\ell,i)},m_\ell,\ee_i^{(\ell)})\}=\left\{\left(k_{k_N^{(1;\ell,i)}}^{\left(0;m_\ell^{(1)}\right)},m^{(0)}_{m_\ell^{(1)}},\ee_i^{(1;\ell)}\right)\right\},\] where $\{m^{(1)}_\ell\}_{\ell=1}^\infty\subs\NN$ is a diverging sequence, $\{\ee_i^{(1;\ell)}\}_{i=1}^\infty\subs(0,\infty)$, $\ell\in\NN$, are sequences converging to $0$ and $\{k_N^{(1;\ell,i)}\}_{N=1}^\infty\subs\NN$, $(\ell,i)\in\NN^2$ are diverging sequences, such that the iterated limit 
	\[\bs{R}_\Psi:=\lim_{\ell,i,N\uparrow\infty}\bs{R}_\Psi^{{k_N^{(\ell,i)}},m_\ell,\ee_i^{(\ell)}}\]
exists. In the case that we assume the condition for the validity of the replacement lemma holds, then by the choice of the family $\{(k_N^{(0;\ell)},m_\ell^{(0)})\}$, the family $\{(k_N^{(0;\ell)},m^{(0)}_\ell,\ee)\}$ satisfies the assumption for the validity of the two-blocks estimate, i.e.~any subfamily $\{(k_N^{(\ell,i)},m_\ell,\ee_i^{(\ell)})\}$ of $\{(k_N^{(0;\ell)},m^{(0)}_\ell,\ee)\}$ has a further subfamily $\{(\bar{k}_N^{(\ell,i)},\bar{m}_\ell,\bar{\ee}_i^{(\ell)})\}$ along which~\eqref{ToProveReplacLemma} holds.

Since $\{k_N^{(\ell,i)}\}_{N=1}^\infty$ is a subsequence of the initial sequence $\{k_N\}$ for all $(\ell,i)\in\NN^2$ we have that 
\[\lim_{N\ra+\infty}(\s^\Psi,\pi)_\sharp\bs{R}_\Psi^{{k_N^{(\ell,i)}},m_\ell,\ee_i^{(\ell)}}=\lim_{N\ra+\infty}\bbar{\bs{Q}}_\Psi^{k_N^{(\ell,i)}}=\bbar{\bs{Q}}_\Psi^\infty.\]
Also $(\s,\bs\pi^\infty)_\sharp\bs{R}_\Psi^{{k_N^{(\ell,i)}},m_\ell,\ee_i^{(\ell)}}=\bbar{\bs{Q}}_\Psi^{k_N^{(\ell,i)},m_\ell}:=(\s^{k_N^{(\ell,i)},\Psi},\bs\pi^{k_N^{(\ell,i)},m_\ell})_\sharp P^{k_N^{(\ell,i)}}$ and thus 
\[(\s^\Psi,\bs\pi^\infty)_\sharp\bs{R}_\Psi=\lim_{\ell,i,N\ra+\infty}(\s,\bs\pi^\infty)_\sharp\bs{R}_\Psi^{{k_N^{(\ell,i)}},m_\ell,\ee_i^{(\ell)}}=\lim_{\ell,i,N\ra+\infty}\bbar{\bs{Q}}_\Psi^{k_N^{(\ell,i)},m_\ell}.\]
For each fixed $(\ell,i)\in\NN^2$ the limit $\bbar{\bs{Q}}_{\Psi,i}^{\infty,m_\ell}=\lim_{N\uparrow\infty}\bbar{\bs{Q}}_\Psi^{k_N^{(\ell,i)},m_\ell}$ exists and belongs in the closed space $\bbar{\bs{\mathcal{Q}}}_\Psi^{\infty,m_\ell}:=\Lim_{N\uparrow\infty}\bbar{\bs{Q}}^{N,m_\ell}$. Consequently $\lim_{i\uparrow\infty}\bbar{\bs{Q}}_{\Psi,i}^{\infty,m_\ell}\in\bbar{\bs{\mathcal{Q}}}_\Psi^{\infty,m_\ell}$ for all $\ell\in\NN$ and therefore
\[(\s^\Psi,\bs\pi^\infty)_\sharp\bs{R}_\Psi\in\bbar{\bs{Q}}_\Psi^{\infty,\infty}:=\Lim_{\ell,N\uparrow\infty}\bbar{\bs{Q}}_\Psi^{N,\ell}.\]
Consequently by the one-block estimate in Theorem~\ref{OBETheorem} (b) it follows that 
\[\bs{R}_\Psi\big\{(\s^\Psi,\bs\pi^\infty,\bs\pi^0,\pi)\bigm|\s^\Psi=B_{\bbar{\Psi}}(\bs\pi^\infty)\big\}=1.\]

Similarly \[(\bs\pi^\infty,\bs\pi^0)_\sharp\bs{R}_\Psi=\lim_{\ell,i,N\uparrow\infty}(\bs\pi^\infty,\bs\pi^0)_\sharp\bs{R}_\Psi^{k_N^{(\ell,i)},m_\ell,\ee_i^{(\ell)}}=\lim_{\ell,i,N\uparrow\infty}\bbar{\bs{Q}}^{k_N^{(\ell,i)},m_\ell,\ee_i^{(\ell)}},\] where $\bbar{\bs{Q}}^{N,\ell,\ee}:=(\bs\pi^{N,\ell},\bs\pi^{N,\ee})_\sharp P^N$, and since $\bbar{\bs{Q}}^{k_N^{(\ell,i)},m_\ell,\ee_i^{(\ell)}}$ is a subfamily of the family $\{\bbar{\bs{Q}}_*^{N,\ell,\ee}\}:=\{\bbar{\bs{Q}}^{k_N^{(0;\ell)},m_\ell^{(0)};\ee}\}$ we have that \[(\bs\pi^\infty,\bs\pi^0)_\sharp\bs{R}_\Psi\in\bbar{\bs{Q}}_*^{\infty,\infty,0}:=\Lim_{\ell\uparrow\infty,\ee\downarrow 0,N\uparrow\infty}\bbar{\bs{Q}}_*^{N,\ell,\ee}.\]
Therefore, since $\{(k_N^{(0;\ell)},m^{(0)}_\ell)\}$ has been chosen so that the iterated limit in~\eqref{TBCAssumption} exists, it follows by the two-blocks comparison that 
\[\bs{R}_\Psi\big\{(\s^\Psi,\bs\pi^\infty,\bs\pi^0,\pi)\bigm|B(\bs\pi^\infty)=B(\bs\pi^0)\;\mbox{ and }\;B_{\X}\circ\widehat{D}(\bs\pi^\infty)\leq B_{\X}\circ\widehat{D}(\bs\pi^0),\;\forall\X\in\bbar{C}_{1,\uparrow}(\RR_+)\big\}=1\]
and in the case that the assumption for the validity of the replacement lemma holds we have that 
\[\bs{R}_\Psi\big\{(\s^\Psi,\bs\pi^\infty,\bs\pi^0,\pi)\bigm|B(\bs\pi^\infty)=B(\bs\pi^0)\;\mbox{ and }\;B_{\X}\circ\widehat{D}(\bs\pi^\infty)= B_{\X}\circ\widehat{D}(\bs\pi^0),\;\forall\X\in\bbar{C}_{1,\uparrow}(\RR_+)\big\}=1\]
Therefore 
\begin{align*}
\bbar{\bs{Q}}_\Psi^\infty\{\s^\Psi\leq\bbar{\Psi}(\pi^{ac})\df\LL_{\T^d}\}&=\bs{R}_\Psi\big\{(\s^\Psi,\bs\pi^\infty,\bs\pi^0,\pi)\bigm|\s^\Psi\leq\bbar{\Psi}(\pi^{ac})\df\LL_{\T^d}\big\}\\
&\geq\bs{R}_\Psi\big\{(\s^\Psi,\bs\pi^\infty,\bs\pi^0,\pi)\bigm|B_{\bbar{\Psi}}\circ\widehat{D}(\bs\pi^0)=\bbar{\Psi}(\pi^{ac})\df\LL_{\T^d}\big\}
\end{align*}
and thus if we prove that 

\begin{equation}\label{ToProveRL}
\bs{R}_\Psi\big\{(\s^\Psi,\bs\pi^\infty,\bs\pi^0,\pi)\bigm|B_{\bbar{\Psi}}\circ\widehat{D}(\bs\pi^0)=\bbar{\Psi}(\pi^{ac})\df\LL_{\T^d}\big\}=1
\end{equation}
it will follow that 
\[\bbar{\bs{Q}}_\Psi^\infty\{\s^\Psi\leq\bbar{\Psi}(\pi^{ac})\df\LL_{\T^d}\}=1\]
with equality if the assumption for the validity of the replacement lemma holds.

If we consider the family of laws 
\begin{equation}\label{NotationMcroEmpStandardEmpJoint}
\bbar{\bs{Q}}^{N,\ee}:=(\bs\pi^{N,\ee},\pi^N)_\sharp P^N\in L_{w^*}^\infty(0,T;\bbar{\Y}_1(\T^d))\x D(0,T;\MMM_+(\T^d))
\end{equation} then since $\{k_N^{(\ell,i)}\}_N$ is a subsequence of $\{k_N^{(0;m^{(1)}_\ell)}\}_N$ for all $(\ell,i)\in\NN^2$, which in turn is a subsequence of $\{k_N\}$ or each $\ell\in\NN$,
\[(\bs\pi^0,\pi)_\sharp\bs{R}_\Psi=\Lim_{\ell,i,N\uparrow\infty}(\bs\pi^0,\pi)_\sharp\bs{R}_\Psi^{k_N^{(\ell,i)},m_\ell,\ee_i^{(\ell)}}=
\lim_{\ell,i,N\uparrow\infty}\bbar{\bs{Q}}^{k_N^{(\ell,i)},\ee_i^{(\ell)}}\in\Lim_{\ee\downarrow 0,N\uparrow\infty}\bbar{\bs{Q}}^{k_N,\ee}.\]
Consequently in order to prove~\eqref{ToProveRL} it suffices to show that any limit point $\bbar{\bs{Q}}_*^{\infty,0}$ of the family $\{\bbar{\bs{Q}}^{K_N,\ee}\}$ as $N\uparrow\infty$ and $\ee\downarrow 0$ is concentrated on trajectories $(\bs\pi^0,\pi)$ such that
$B_\Psi(\bs\pi^0)=\Psi(\pi^{ac})\df\LL_{\T^d}$ for any non-decreasing sublinear map $\Psi\in C_1(\RR_+)$. This follows from the next, slightly more general proposition.
\begin{prop} Let $\{k_N\}_{N=1}^\infty\subs\NN$ be a diverging sequence such that the laws $Q^{k_N}:=\pi^{k_N}P^{k_N}$ of the empirical density in the Skorohod space  converge to a law $Q_*^\infty\in\PP D(0,T;\MMM_+(\T^d))$ as $N\uparrow\infty$. Then for any subsequential limit point $\bbar{\bs{Q}}_*^{\infty,0}$ as $N\uparrow\infty$ and then $\ee\downarrow 0$ of the family of laws $\{\bbar{\bs{Q}}^{k_N,\ee}\}$, where $\bbar{\bs{Q}}^{N,\ee}$ is defined in~\eqref{NotationMcroEmpStandardEmpJoint},
	\begin{equation}\label{MacroReplacement1}
		\bbar{\bs{Q}}_*^{\infty,0}\big\{(\bs\pi^0,\pi)\bigm|B_\Psi\circ\widehat{D}(\bs\pi^0)=\Psi(\pi^{ac})\df\LL_{\T^d},\;\forall\Psi\in\bbar{C}_{1}(\RR_+)\big\}=1,
		\end{equation}
	where $\pi=\pi^{ac}+\pi^\perp$, $\pi^{ac}\ll\df\LL_{\T^d}$, $\pi^\perp\LL_{\T^d}$, is the Radon-Nikodym decomposition of $\pi\in\MMM_+(\T^d)$.
	 
	Furthermore this implies that also 
	\begin{equation}\label{MacroReplacement2}
	\bbar{\bs{Q}}_*^{\infty,0}\big\{(\bs\pi^0,\pi)\bigm|B_\Psi\circ D^\perp(\bs\pi^0)=\Psi'(\infty)\pi^\perp,\;\forall\Psi\in\bbar{C}_{1}(\RR_+)\big\}=1,
	\end{equation}
	and thus in particular the barycentric projection $B\circ D^\perp$ of the Young measures singular part operator $D^\perp$ yields the Radon-Nikodym singular part of  ordinary measures with respect to the law $\bbar{\bs{Q}}_*^{\infty,0}$.
\end{prop}\textbf{Proof} Let us check first that~\eqref{MacroReplacement1} implies~\eqref{MacroReplacement2} holds. Indeed, if~\eqref{MacroReplacement1} holds then for $\bbar{\bs{Q}}_*^{\infty,0}$-a.s.~all $(\bs\pi^0,\pi)$ it holds that $B_{\Psi_0}(\bs\pi^0)=B_{\Psi_0}\circ\widehat{D}(\bs\pi^0)=\Psi_0(\pi^{ac})\df\LL_{\T^d}$ for all sublinear maps. Thus if given $\Psi\in\bbar{C}_1(\RR_+)$ we set $\Psi_0(\lambda):=\Psi(\lambda)-\Psi'(\infty)\lambda$ we have that $\bbar{\bs{Q}}_\Psi$-a.s.
\begin{equation}\label{ForMacroSingular1}
B_\Psi\circ\widehat{D}(\bs\pi^0)=\Psi(\pi^{ac})\quad\mbox{and}\quad B_{\Psi_0}(\bs\pi^0)=\Psi_0(\pi^{ac}).
\end{equation}
By the second equality above 
\begin{equation}\label{ForMacroSingular2}B_\Psi(\bs\pi^0)-\Psi'(\infty)B(\bs\pi^0)=B_{\Psi_0}(\bs\pi^0)=\Psi_0(\pi^{ac})=\Psi(\pi^{ac})-\Psi'(\infty)\pi^{ac}\end{equation}
and since $B_\Psi(\bs\pi^0)=B_\Psi\circ\widehat{D}(\bs\pi^0)+\Psi'(\infty) B\circ D^\perp(\bs\pi^0)$ it follows by the first equality in~\eqref{ForMacroSingular1} and~\eqref{ForMacroSingular2} that
\[\Psi'(\infty) B\circ D^\perp(\bs\pi^0)-\Psi'(\infty)B(\bs\pi^0)=-\Psi'(\infty)\pi^{ac}\]
which by rearranging becomes 
\[\Psi'(\infty) B\circ D^\perp(\bs\pi^0)=\Psi'(\infty)\big(B(\bs\pi^0)-\pi^{ac}\big).\]
Now as in the proof of~Theorem~\ref{TBCTheorem}(b)(i) it easily follows that $B(\bs\pi^0)=\pi$ $\bbar{\bs{Q}}^{\infty,0}_*$-a.s.~and therefore 
\[B_\Psi\circ D^\perp(\bs\pi^0)=\Psi'(\infty) B\circ D^\perp(\bs\pi^0)=\Psi'(\infty)\pi^\perp\]
as claimed. 

Next we prove~\eqref{MacroReplacement1}. So let $\bbar{\bs{Q}}_*^{\infty,0}$ be a subsequential limit point of the family $\bbar{\bs{\mathcal{Q}}}^{k_N,\ee}$ defined in~\eqref{NotationMcroEmpStandardEmpJoint} and let $\Psi\in\bbar{C}_1(\RR_+)$. For each $\delta>0$ and $G\in L^1(0,T;C(\T^d))$ the set 
$$A_\delta^G:=\big\{(\bs\pi^0,\pi)\in L_{w^*}^\infty(0,T;\bbar{\Y}_1(\T^d))\x D(0,T;\MMM_+(\T^d))\bigm||\lls G,B_\Psi\circ\widehat{D}(\bs\pi^0)-I_\Psi(\pi)|>\delta\big\}$$
is measurable and since $L^1(0,T;C(\T^d))$ is separable in order to prove that $\bbar{\bs{Q}}_*^{\infty,0}$ is concentrated on trajectories $(\bs\pi^0,\pi)$ such that $B_\Psi(\widehat{\bs\pi}_t^0)=\Psi(\pi^{ac}_t)\df\LL_{\T^d}$ it suffices to show that $\bbar{\bs{Q}}_*^{\infty,0}(A_\delta^G)=0$ for all $\delta>0$ and all $G\in L^1(0,T;C(\T^d))$. We start by writing \[\bbar{\bs{Q}}_*^{\infty,0}(A_\delta^G)=(\widehat{D}\x\mathbbm{id}_{D(0,T;\MMM_+(\T^d))})_\sharp\bbar{\bs{Q}}_*^{\infty,0}\big\{(\bs\s^0,\pi)|\lls G,B_\Psi(\bs\s^0)-I_\Psi(\pi)\rrs|>\delta\big\}.\]
By Chebyshev's inequality it suffices to show that 
\begin{equation}\label{CheBToProvRLC}
B_{\Psi,G}:=\int|\lls G,B_\Psi(\bs\s^0)-I_\Psi(\pi)\rrs|\df(\widehat{D}\x\mathbbm{id}_{D(0,T;\MMM_+(\T^d))})_\sharp\bbar{\bs{Q}}_*^{\infty,0}(\bs\s^0,\pi)=0
\end{equation} 
for all $G\in L^1(0,T;C(\T^d))$. By using the Moreau-Yosida approximations $\{\Psi_k\}_{k\in\NN}$ of $\Psi\in\bbar{C}_1(\RR_+)$, i.e.~ $\Psi_k(\lambda):=\Psi_{0,k}(\lambda)+\Psi'(\infty)\lambda$, where $\Psi_{0,k}$ are the Moreau-Yosida approximations of $\Psi_0(\lambda):=\Psi(\lambda)-\Psi'(\infty)\lambda$ defined in~\eqref{MoreauYosida} we can see that it suffices to prove~\eqref{Term1MacroRL} under the additional assumption that $\Psi$ is Lipschitz. This follows by the dominated convergence similarly to the reduction to the case of Lipschitz maps $\Psi$ in the proof of Theorem~\ref{TBCTheorem} since the maps $\Psi_k$ converge to $\Psi$ and $\|\Psi_k\|_{\infty,1}\leq\|\Psi_0\|_{\infty;1}+|\Psi'(\infty)|$ for all large enough $k\in\NN$ and the measure $(\widehat{D}\x\mathbbm{id}_{D(0,T;\MMM_+(\T^d))})_\sharp\bbar{\bs{Q}}_*^{\infty,0}$ is concentrated on trajectories $(\bs\s^0,\pi)$ such that $\ls\Lambda,\bs\s^0_t\rs\leq m$ and $\ls 1,\pi_t\rs=m$ for almost all $t\in[0,T]$.

Since the map $I_\Psi$ is not continuous, we interpolate with the continuous map $I_\Psi^{{\ee '},M'}$, $\ee '\in(0,1/2)$, $M'\in(0,\infty)$ to obtain 
\begin{align}\label{FB}
B_{\Psi,G}&\leq\int|\lls G,B_\Psi(\bs\s^0)-I_\Psi^{\ee ',M'}(\pi)\rrs|\df(\widehat{D}\x\mathbbm{id}_{D(0,T;\MMM_+(\T^d))})_\sharp\bbar{\bs{Q}}_*^{\infty,0}(\bs\s^0,\pi)\nonumber\\
&\qquad+\int|\lls G,I_\Psi^{\ee ',M'}(\pi)-I_\Psi(\pi)\rrs|\df Q_*^\infty(\pi),
\end{align}
where $Q_*^\infty:=\lim_{N\uparrow\infty}Q^{k_N}=\pi_\sharp\bbar{\bs{Q}}_*^{\infty,0}$. Since $\bbar{\bs{Q}}_*^{\infty,0}\in\bbar{\bs{\mathcal{Q}}}_*^{\infty,0}:=\Lim_{\ee\downarrow 0,N\uparrow\infty}\bbar{\bs{\mathcal{Q}}}^{k_N,\ee}$ there exists a sequence $(\ee_i)_{i\in\NN}$ and subsequences $\{k^{(i)}_N\}_{N=1}^\infty$, $i\in\NN$, of $\{k_N\}$ such that 
\[\bbar{\bs{Q}}_*^{\infty,0}=\lim_{i,N\uparrow\infty}\bbar{\bs{Q}}^{k_N^{(i)},\ee_i}.\]
By~\eqref{ApproxFirstCoord}, $\widehat{D}\circ\bs\pi^{N,\ee;M}=\Pi_M^*\circ j^*\circ\bs\pi^{N,\ee}$ and therefore setting $\bbar{\bs{Q}}^{N,\ee;M}:=(\bs\pi^{N,\ee;M},\pi^N)_\sharp P^N$,
\begin{align*}
\big(\widehat{D}\x \mathbbm{id}\big)_\sharp\bbar{\bs{Q}}^{k_N^{(i)},\ee_i;M}&=\big(\widehat{D}\circ\bs\pi^{k_N^{(i)},\ee_i;M},\pi^{k_N^{(i}}\big)_\sharp P^{k_N^{(i)}}\\
&=(\Pi_M^*\circ j^*\circ\bs\pi^{{k_N^{(i)}},\ee_i},\pi^{k_N^{(i)}})_\sharp P^{k_N^{(i)}}
=\big((\Pi_M^*\circ j^*)\x\mathbbm{id}\big)_\sharp\bbar{\bs{Q}}^{k_N^{(i)},\ee_i}.
\end{align*}
Since the map $(\Pi_M^*\circ j^*)\x\mathbbm{id}$ is continuous
\[\lim_{i,N\uparrow\infty}\big(\widehat{D}\x\mathbbm{id}\big)_\sharp\bbar{\bs{Q}}^{k_N^{(i)},\ee_i;M}=\big((\Pi_M^*\circ j^*)\x\mathbbm{id}\big)_\sharp\bbar{\bs{Q}}_*^{\infty,0}\]
and since $\Pi_M^*\circ j^*$ $w^*$-converges pointwise to $\widehat{D}$ it follows that 
\[\lim_{M,i,N\uparrow\infty}\big(\widehat{D}\x \mathbbm{id}\big)_\sharp\bbar{\bs{Q}}^{k_N^{(i)},\ee_i;M}
=(\widehat{D}\x\mathbbm{id})_\sharp\bbar{\bs{Q}}_*^{\infty,0}.\]
Consequently, since the map $(\bs\pi^0,\pi)\mapsto|\lls G,B_\Psi(\bs\s^0)-I^{\ee',M'}(\pi)\rrs|$ is continuous in $L_{w^*}^\infty(0,T;\bbar{\Y}_1(\T^d))\x D(0,T;\MMM_+(\T^d)$, it follows by the portmanteau theorem that for each $\ee'\in(0,1/2)$, $M'<+\infty$ the first term in the right hand side of~\eqref{FB} is equal to
\[\lim_{M,i,N\uparrow\infty}\int|\lls G,B_\Psi(\bs\s^0)-I_\Psi^{\ee ',M'}(\pi)\rrs|\df(\widehat{D}\x\mathbbm{id})_\sharp\bbar{\bs{Q}}^{k_N^{(i)},\ee_i;M}(\bs\s^0,\pi)\]
and thus, since $\{k_N^{(i)}\}$ is a subsequence of $\{k_N\}$ for all $i\in\NN$, inequality~\eqref{FB} becomes
\begin{align*}
B_{\Psi,G}&\leq\lim_{M,i,N\uparrow\infty}\int|\lls G,B_\Psi(\bs\s^0)-I_\Psi^{\ee ',M'}(\pi)\rrs|\df(\widehat{D}\x\mathbbm{id})_\sharp\bbar{\bs{Q}}^{k_N^{(i)},\ee_i;M}(\bs\s^0,\pi)\\
&\qquad+\int|\lls G,I_\Psi^{\ee ',M'}(\pi)-I_\Psi(\pi)\rrs|\df Q_*^\infty(\pi)\\
&\leq
\lim_{M,i,N\uparrow\infty}\int|\lls G,B_\Psi(\bs\s^0)-I_\Psi^{\ee_i,M}(\pi)(\pi)\rrs|\df(\widehat{D}\x\mathbbm{id})_\sharp\bbar{\bs{Q}}^{k_N^{(i)},\ee_i;M}(\bs\s^0,\pi)\\
&\qquad+\limsup_{M\uparrow\infty,\ee\downarrow 0,N\uparrow\infty}\int|\lls G,I_\Psi^{\ee,M}(\pi)-I_\Psi^{\ee',M'}(\pi)\rrs|Q^{k_N}(\pi)\\
&\qquad+\int|\lls G,I_\Psi^{\ee ',M'}(\pi)-I_\Psi(\pi)\rrs|\df Q_*^\infty(\pi).
\end{align*}
By taking the limit as $\ee'\downarrow 0$ and then $M'\uparrow\infty$ we obtain
\begin{align}
B_{\Psi,G}&\leq\lim_{M,i,N\uparrow\infty}\int|\lls G,B_\Psi(\bs\s^0)-I_\Psi^{\ee_i,M}(\pi)(\pi)\rrs|\df(\widehat{D}\x\mathbbm{id})_\sharp\bbar{\bs{Q}}^{k_N^{(i)},\ee_i;M}(\bs\s^0,\pi)\label{Term1MacroRL}\\
&\qquad+\limsup_{M'\uparrow\infty,\ee'\downarrow 0,M\uparrow\infty,\ee\downarrow 0,N\uparrow\infty}\int|\lls G,I_\Psi^{\ee,M}(\pi)-I_\Psi^{\ee',M'}(\pi)\rrs|Q^{k_N}(\pi)\label{Term2MacroRL}\\
&\qquad+\limsup_{M'\uparrow\infty,\ee'\downarrow 0}\int|\lls G,I_\Psi^{\ee ',M'}(\pi)-I_\Psi(\pi)\rrs|\df Q_*^\infty(\pi).\label{Term3MacroRL}
\end{align} 

We will prove that all these terms are equal to zero. We start with the iterated limit~\eqref{Term3MacroRL}. As we have seen the maps $I^{\ee,M}_\Psi$ converge as $\ee\downarrow 0$ pointwise in $D(0,T;\MMM_+(\T^d))$ to $I^{0,M}_\Psi$ with respect to the normed $w^*$-convergence in the target space $L_{w^*}^\infty(0,T;\MMM(\T^d))$ and the maps $I^{0,M}_\Psi$ converge in the same topology to the map $I_\Psi$ as $M\uparrow\infty$. Consequently \[\lim_{\ee\downarrow 0,M\uparrow\infty}|\lls G,I^{\ee,M}(\pi)-I_\Psi(\pi)\rrs|=0\] pointwise and since for each fixed $M>0$ we have by~\eqref{MacroRepDomConvB1} the bound
\[\sup_{\ee>0}\sup_{\pi\in D(0,T;\MMM_+(\T^d))}|\lls G,I_\Psi^{\ee,M}(\pi)-I^{0,M}_\Psi(\pi)\rrs|\leq2\|G\|_{1;\infty}\sup_{0\leq\lambda\leq M}|\Psi(\lambda)|\]
and by~\eqref{MacroRepDomConvB2} the bound
\[\sup_{M>0}|\lls G,I_\Psi^{0,M}(\pi)-I_\Psi(\pi)|\leq2\|G\|_{\infty;1}\|\Psi\|_{\infty,1}(1+\mathfrak{m})\]
for $Q^\infty_*$-a.s.~all $\pi\in D(0,T;\MMM_+(\T^d))$, it follows by the dominated convergence theorem that~\eqref{Term3MacroRL} holds.

For the proof of~\eqref{Term2MacroRL} we note that the map $D(0,T;\MMM_+(\T^d))\ni\pi\mapsto|\lls G,I^{\ee,M}(\pi)-I^{\ee',M'}(\pi)\rrs|$ is continuous for each fixed $\ee,\ee',M,M'>0$ and thus since $Q_*^\infty:=\lim_{N\ra+\infty}Q^{k_N}$, it follows by the portmanteau theorem that
\[\lim_{N\ra+\infty}\int|\lls G,I_\Psi^{\ee,M}(\pi)-I_\Psi^{\ee',M'}(\pi)\rrs|\df Q^{k_N}(\pi)=\int|\lls G,I_\Psi^{\ee,M}(\pi)-I_\Psi^{\ee',M'}(\pi)\rrs|\df Q^\infty_*(\pi)\]
By further interpolating the term $I_\Psi(\pi)$ in the difference $I_\Psi^{\ee,M}(\pi)-I_\Psi^{\ee',M'}(\pi)$ in the right hand side above and taking the iterated limit as $\ee\downarrow 0$, $M\uparrow\infty$ $\ee'\downarrow 0$ and finally $M'\uparrow\infty$, we obtain that the iterated limit in~\eqref{Term2MacroRL} is bounded above by
\[2\limsup_{M\uparrow\infty,\ee\downarrow 0}\int|\lls G,I_\Psi^{\ee,M}(\pi)-I_\Psi(\pi)\rrs|\df Q^\infty_*(\pi)\]
which is equal to zero by~\eqref{Term3MacroRL}.

It remains to prove~\eqref{Term1MacroRL} under the additional assumption that $\Psi$ is Lipschitz. This follows by the next lemma/
\begin{lemma} For all $G\in L^1(0,T;C(\T^d))$ and all Lipschitz maps $\Psi\in\bbar{C}_1(\RR_+)$
\begin{equation}\label{ToCompleteRL}
\limsup_{M\uparrow\infty,\ee\downarrow 0,N\uparrow\infty}\EE^N|\lls G,B_\Psi\circ\widehat{D}(\bs\pi^{N,\ee;M})-I_\Psi^{\ee,M}(\pi^N)\rrs|=0.
\end{equation}
\end{lemma}\textbf{Proof} We recall that $k_\ee:=\fr{(2\ee)^d}\1_{[-\ee,\ee]^d}$ and note` that
 \begin{align}\label{pke}\pi^N*k_\ee(x/N)&=
\int_{\T^d}k_\ee\Big(u-\frac{x}{N}\Big)\df\pi^N(u)
=\fr{N^d}\sum_{y\in\T_N^d}\eta(y)k_\ee\Big(\frac{y-x}{N}\Big)\nonumber\\
&=\fr{(2\ee N)^d}\sum_{y\in\T_N^d}\eta(y)\1_{[-\ee,\ee]^d}\Big(\frac{y-x}{N}\Big)=\fr{(2\ee N)^d}\sum_{y:|y-x|\leq[N\ee]}\eta(y)\nonumber\\
&=\frac{(2[N\ee]+1)^d}{(2N\ee)^d}\eta^{[N\ee]}(x)=:C_{N,\ee}\cdot\eta^{[N\ee]}(x),
\end{align}
where the constant $C_{N,\ee}$ satisfies $\lim_{N\ra+\infty}C_{N,\ee}=1$ for all $\ee>0$. This implies that for each $\Psi\in\bbar{C}_1(\RR_+)$ and $G\in L^1(0,T;C(\T^d))$ 
\begin{align*}\lls G,B_\Psi\circ\widehat{D}(\bs\pi^{N,\ee;M})\rrs&=\int_0^T\fr{N^d}\sum_{x\in\T_N^d}G_t\Big(\frac{x}{N}\Big)\Psi\big(\eta_t^{[N\ee]}(x)\mn M\big)\df t\\
&=\int_0^T\fr{N^d}\sum_{x\in\T_N^d}G_t\Big(\frac{x}{N}\Big)\Psi\left(\big(C_{N,\ee}^{-1}\pi^N_t*k_\ee(x/N)\big)\mn M\right)\df t
\end{align*}
Consequently, if we define the map $K_\Psi^{N,\ee;M}\colon D(0,T;\MMM_+(\T^d)\to L_{w^*}^\infty(0,T;\MMM_+(\T^d))$ by 
\[\lls G,K^{N,\ee;M}_\Psi(\pi)\rrs=\int_0^T\fr{N^d}\sum_{x\in\T_N^d}G_t\Big(\frac{x}{N}\Big)\Psi\big(\pi_t*k_\ee(x/N)\mn M\big)\df t,\] since $\Psi$ is assumed to be in addition Lipschitz and the map $\lambda\mapsto\lambda\mn M$ is $1$-Lipschitz,
\begin{align*}|\lls G,B_\Psi\circ\widehat{D}(\bs\pi^{N,\ee;M})-K_\Psi^{N,\ee;M}(\pi^N)\rrs|
&\leq\Lip_\Psi(1-C_{N,\ee}^{-1})\int_0^T\|G_t\|_\infty\fr{N^d}\sum_{x\in\T_N^d}\pi^N_t*k_\ee(x/N)\df t\\
&=\Lip_\Psi(C_{N,\ee}-1)\int_0^T\|G_t\|_\infty\fr{N^d}\sum_{x\in\T_N^d}\eta_t^{[N\ee]}\df t\\
&\stackrel{P^N\mbox{a.s.}}=\Lip_\Psi\|G\|_{\infty;1}(C_{N,\ee}-1)\ls 1,\pi_0^N\rs.
\end{align*}
Therefore by Lemma~\ref{BoundedTotParticlInMacrLim} it follows that 
\[\limsup_{N\uparrow\infty}\EE^N|\lls G,B_\Psi\circ\widehat{D}(\bs\pi^{N,\ee;M})-K_\Psi^{N,\ee;M}(\pi^N)\rrs|=0\]
and thus the process $B_\Psi\circ\widehat{D}(\bs\pi^{N,\ee;M})$ can be replaced in the limit as $N\uparrow\infty$ by the process $K_\Psi^{N,\ee;M}(\pi^N)$.

Next we consider the process $I_\Psi^{N,\ee;M};D(0,T;\MMM_+(\T^d))\to L_{w^*}^\infty(0,T;\MMM_+(\T^d))$ defined by 
$$\lls G,I_\Psi^{N,\ee;M}(\pi)\rrs=\int_0^T\fr{N^d}\sum_{x\in\T_N^d}G_t\Big(\frac{x}{N}\Big)\Psi\big(\pi_t*\iota_\ee(x/N)\mn M\big)\df t$$ and we will show that the process $K_\Psi^{N,\ee;M}$ can be replaced by the process $I_\Psi^{N,\ee;M}$. Since for any $\theta\in(0,1)$ and $a,b,b'\geq 0$ the implication
$$\theta b'\leq a\leq \frac{b}{\theta}\quad\Lra\quad|b-a|\leq\frac{b}{\theta}-\theta^2b'\leq\Big(\fr{\theta}-\theta^2\Big)b+(b-b')$$
holds, it follows from~\eqref{e} with $\theta_\ee:=(1-\ee)^d$, $b'=\pi*k_{\ee-\ee^2}(u)$, $b=\pi*k_\ee(u)$ and $a=\pi*\iota_\ee(u)$ that 
$$|\pi*k_\ee(u)-\pi_t*\iota_\ee(u)|\leq\Big(\fr{\theta_\ee}-\theta_\ee^2\Big)\pi*k_\ee(u)+\big(\pi*k_\ee(u)-\pi*k_{\ee-\ee^2}(u)\big)$$
for all $u\in\T^d$ and thus we compute 
\begin{align*}
|\lls G,K_\Psi^{N,\ee;M}(\pi^N)-I_\Psi^{N,\ee;M}(\pi^N)\rrs|&\leq\Lip_\Psi\int_0^T\frac{\|G_t\|_\infty}{N^d}\sum_{x\in\T_N^d}|\pi_t^N*k_\ee(x/N)-\pi_t^N*\iota_\ee(x/N)|\df t\\
&\leq\Lip_\Psi\Big(\fr{\theta_\ee}-\theta_\ee^2\Big)\int_0^T\frac{\|G_t\|_\infty}{N^d}\sum_{x\in\T_N^d}\pi_t^N*k_\ee(x/N)\df t\\
&\quad+\Lip_\Psi\int_0^T\frac{\|G_t\|_\infty}{N^d}\sum_{x\in\T_N^d}(\pi_t^N*k_\ee-\pi_t^N*k_{\ee-\ee^2})(x/N)\df t
\end{align*}
By~\eqref{pke} and the conservation of the total number of particles the first term in the right hand side above is $P^N$-a.s.~equal to 
$$\Lip_\Psi\Big(\fr{\theta_\ee}-\theta_\ee^2\Big)C_{N,\ee}\|G\|_{\infty;1}\ls 1,\pi^N_0\rs.$$
For the second term we set $\bar{\eta}^{[N\ee]}(x)=\sum_{y:|y-x|\leq[N\ee]}\eta(y)$ so that by~\eqref{pke} we can write $\pi*k_\ee(x/N)=(2\ee N)^{-d}\bar{\eta}^{[N\ee]}$ and with this notation 
\begin{align*}
(\pi_t^N*k_\ee-\pi_t^N*k_{\ee-\ee^2})(x/N)&=\fr{(2\ee N)^d}\bar{\eta}^{[N\ee]}-\fr{(2\ee(1-\ee)N)^d}\bar{\eta}^{[N\ee(1-\ee)]}\\
		&=\fr{(2\ee(1-\ee)N)^d}\big((1-\ee)^d\bar{\eta}^{[N\ee]}(x)-\bar{\eta}^{[N\ee(1-\ee)]}(x)\big)\\
		&\leq\fr{(2\ee(1-\ee)N)^d}\big(\bar{\eta}^{[N\ee]}(x)-\bar{\eta}^{[N\ee(1-\ee)]}(x)\big)\\
		&=\fr{(2\ee(1-\ee)N)^d}\sum_{[N\ee(1-\ee)]<|y|\leq[N\ee]}\eta(x+y)	
		\end{align*}
		and thus since 
		$$\sum_{x\in\T_N^d}\sum_{[N\ee(1-\ee)]<|y|\leq[N\ee]}\eta(x+y)=([N\ee]_\star^d-[N\ee(1-\ee)]_\star^d)\sum_{x\in\T_N^d}\eta(x)$$
the integral in the second term can be bounded from above by  
$$\frac{[N\ee]_\star^d-[N\ee(1-\ee)]_\star^d}{(2\ee(1-\ee)N)^d}\int_0^T\frac{\|G_t\|_\infty}{N^d}\sum_{x\in\T_N^d}\eta_t(x)\df t\stackrel{P^N\mbox{-a.s.}}=
\frac{[N\ee]_\star^d-[N\ee(1-\ee)]_\star^d}{(2\ee(1-\ee)N)^d}\|G\|_{1;\infty}\ls1,\pi_0^N\rs.$$
To summarize we have the bound
\[|\lls G,K_\Psi^{N,\ee;M}(\pi^N)-I_\Psi^{N,\ee;M}(\pi^N)\rrs|\leq\Lip_\Psi\Big\{\Big(\fr{\theta_\ee}-\theta_\ee^2\Big)C_{N,\ee}+\frac{[N\ee]_\star^d-[N\ee(1-\ee)]_\star^d}{(2\ee(1-\ee)N)^d}\Big\}\|G\|_{1;\infty}\ls1,\pi_0^N\rs,\]
and thus since $\mathfrak{m}_*:=\limsup_{N\ra\infty}\int\ls 1,\pi^N\rs\df\mu_0^N<+\infty$ by Lemma~\ref{BoundedTotParticlInMacrLim}, it follows by taking the limit superior of the expected values as $N\uparrow\infty$ that 
\[\limsup_{N\uparrow\infty}\EE^N|\lls G,K_\Psi^{N,\ee;M}(\pi^N)-I_\Psi^{N,\ee;M}(\pi^N)\rrs|\leq
\Lip_\Psi\Big\{\Big(\fr{\theta_\ee}-\theta_\ee^2\Big)+\fr{(1-\ee)^d}-1\Big\}\|G\|_{1;\infty}\mathfrak{m}_*\]
and thus since $\theta_\ee=(1-\ee)^d$ converges to $1$ as $\ee\downarrow 0$ it follows that 
\[\limsup_{\ee\downarrow 0,N\uparrow\infty}\EE^N|\lls G,K_\Psi^{N,\ee;M}(\pi^N)-I_\Psi^{N,\ee;M}(\pi^N)\rrs|=0.\]

We show finally that the process $I_\Psi^{N,\ee;M}(\pi^N)$ can be replaced by the process $I^{\ee,M}(\pi^N)$. The process $I_\Psi^{N,\ee;M}(\pi^N)$ can be written as 
\[\lls G,I_\Psi^{N,\ee;M}(\pi^N)\rrs=\int_0^T\int_{\T^d}G_t\Big(\frac{[Nu]}{N}\Big)\Psi\big(\pi_t^N*\iota_\ee([Nu]/N)\mn M\big)\df u\df t\]
and therefore since $\Psi$ is assumed to be Lipschitz
\begin{align}\label{LastTerm}
|\lls G,I_\Psi^{N,\ee;M}(\pi^N)-I^{\ee,M}_\Psi(\pi^N)\rrs|&\leq\Lip_\Psi\int_0^T\|G_t\|_\infty\int_{\T^d}\big|\pi_t^N*\iota_\ee([Nu]/N)-\pi_t^N*\iota_\ee(u)\big|\df u\df t\nonumber\\
&\quad+\sup_{0\leq \lambda\leq M}\Psi(\lambda)\int_0^T\int_{\T^d}\Big|G_t\Big(\frac{[Nu]}{N}\Big)-G_t(u)\Big|\df u\df t.
\end{align}
The second term in the right hand side above is deterministic, and since $G\in L^1(0,T;C(\T^d))$, it follows by the dominated convergence theorem that it converges to $0$ as $N\uparrow\infty$ for each $\ee,M>0$. For the first term we estimate 
\begin{align*}
\big|\pi_t^N*\iota_\ee([Nu]/N)-\pi_t^N*\iota_\ee(u)\big|&=\Big|\int\iota_\ee\Big(\y-\frac{[Nu]}{N}\Big)\df\pi_t^N(\y)-\int\iota_\ee(\y-u)\df\pi_t^N(\y)\Big|\\
&=\fr{N^d}\Big|\sum_{x\in\T_N^d}\iota_\ee\Big(\frac{x-[Nu]}{N}\Big)\eta_t(x)-\sum_{x\in\T_N^d}\iota_\ee\Big(\frac{x}{N}-u\Big)\eta_t(x)\Big|\\
&\leq\fr{N^d}\sum_{x\in\T_N^d}\Big|\iota_\ee\Big(\frac{x-[Nu]}{N}\Big)-\iota_\ee\Big(\frac{x}{N}-u\Big)\Big|\eta_t(x).
\end{align*}
Now, the map $\iota_\ee$ is uniformly continuous and therefore for any $\theta>0$ there exists $\delta_{\ee,\theta}>0$ such that 
\[d_{\T^d}(u,\y)<\delta_{\ee,\theta}\quad\Lra\quad|\iota_\ee(u)-\iota_\ee(\y)|<\theta\]
and then for all $N\in\NN$ large enough so that $\fr{N}<\delta_{\ee,\theta}$
\[\big|\pi_t^N*\iota_\ee([Nu]/N)-\pi_t^N*\iota_\ee(u)\big|\leq\theta\ls 1,\pi_t^N\rs\]
and thus by the conservation of the total number of particles it follows that for large enough $N\in\NN$ the iterated integral in the first term of the right hand side of~\eqref{LastTerm} is bounded $P^N$-a.s.~from above by $\theta\|G\|_{\infty;1}\ls 1,\pi_0^N\rs$. Thus by Lemma~\eqref{BoundedTotParticlInMacrLim}
\[\limsup_{N'ra+\infty}\EE^N|\lls G,I_\Psi^{N,\ee;M}(\pi^N)-I^{\ee,M}_\Psi(\pi^N)\rrs|\leq\theta\mathfrak{m}_*,\]
where $\mathfrak{m}_*:=\limsup_{N\uparrow\infty}\EE^N\ls 1,\pi^N_0\rs<+\infty$, and thus since $\theta>0$ is arbitrary the limit superior above is equal to zero which proves that the process $I_\Psi^{N,\ee;M}(\pi^N)$ can be finally replaced by the process $I^{\ee,M}_\Psi(\pi^N)$ in the limit as $N\uparrow\infty$. This completes the proof of the lemma.$\hfill\Box$

\appendix
	\section{Functional analytic prerequisites}
	In our approach to the hydrodynamic limit of the ZRP the various empirical processes we will consider will take values in the non-separable dual $L^*$ of an appropriate separable Banach space $L$, equipped with its $w^*$-topology, denoted by $w^*_L$ or simply $w^*$. Thus each empirical processes will be a random variable $\pi\colon(\W,\F,P)\to L^*$ defined on an appropriate probability space $\W$ with the target space $L^*$ being considered as a measurable space with respect to the Borel $\s$-algebra $\B_{w^*}(L^*)$ generated by the $w^*$-topology of $L^*$ and the laws of the empirical processes will be probability measures on this Borel $\s$-algebra $\B_{w^*}(L^*)$. When the space $L^*$ is evident from the context we will just write $\B_{w^*}$ instead of $\B_{w^*}(L^*)$. 
	
	In fact the non-separable dual $L^*$ will be of the form $L^*=L_{w^*}^\infty(0,T;X^*)$ for some separable Banach space $X$ of test functions, for example $X=C(\T^d)$. Here $L_{w^*}^\infty(0,T;X^*)$ is the vector space of all \emph{$w^*$-measurable maps} $\mu\colon[0,T]\to X^*$ such $\esssup_{0\leq t\leq T}\|\mu_t\|_{X^*}<+\infty$ and a map $\mu\colon[0,T]\to X^*$ is called $w^*$-measurable if the map $[0,T]\ni t\mapsto\ls f,\mu_t\rs$ is measurable for all $f\in X$. As we will see, if $X$ is separable then the map $t\mapsto\|\mu_t\|_{X^*}$ is measurable for any $w^*$-measurable map $\mu\colon[0,T]\to X^*$ and $L_{w^*}^\infty(0,T;X^*)$ is the dual of the $L^1$-Bochner space $L^1(0,T;X)$ of all strongly measurable Bochner integrable maps $f\colon[0,T]\to X$.
	
	The main aim of this appendix is twofold. First we collect the necessary background on the definition of $L_{w^*}^\infty$-spaces of $w^*$-measurable vector valued maps and the isometry $L^1(0,T;X)^*\cong L_{w^*}^\infty(0,T;X^*)$ for separable Banach spaces $X$, without assuming that $X^*$ satisfies the Radon-Nikodym property and secondly to assert that the classical results of topological measure theory on probability measures in polish spaces are valid for the space $\PP L^*\equiv\PP(L^*,w^*)$ of probability measures on the $w^*$-dual $L^*$ of a separable Banach space $L$ when equipped with its $w^*$-topology, i.e.~for example the weak convergence is Hausdorff, satisfies the portmanteau theorem and the Prokhorov relative compactness criterion. The properties of $(L^*,w^*)$ as a topological space that the topological measure theory relies on are (Hausdorff) complete regularity, submetrizability and $\s$-compactness. The complete regularity ensures that there enough bounded and continuous functions for weak convergence of probability measures to be meaningful. In fact the standard theory of weak convergence is true under the complete regularity assumption if one restricts attention to Radon measures. The submetrizability is required for the Prokhorov theorem to be also valid for sequential compactness, and $\s$-compactness together with submetrizability ensure that all probability measures on $L^*$ are Radon, i.e.~that $(L^*,w^*)$ is a Radon space.


\subsection{Duals of separable Banach spaces and submetrizability}\label{Appendix1}
	
	Let $L$ be a Banach space. The dual space $L^*$ is a Hausdorff topological vector space when equipped with the $w^*$-topology. As such it is completely regular since any topological group that satisfies the $T_1$-separation axiom is a completely regular Hausdorff space~\cite{Husain1966a}. It is a well known fact of functional analysis that if $L$ is separable there exists a metric $d\colon L^*\x L^*\to\RR$ that metrizes the $w^*$-topology of $L^*$ on (norm) bounded subsets of $L^*$. For example one can define $d$ by 
	\begin{equation}\label{AuxMetr}
		d(\mu,\nu)=\sum_{k=1}^\infty\fr{2^k}\psi(|\nu(f_k)-\mu(f_k)|),\quad \mu,\nu\in X^*,
	\end{equation}
	where $\{f_k\}_{k=1}^\infty$ is a countable dense subset of $L$ and $\psi\colon\RR_+\to\RR_+$ is the map $\psi(t)=\frac{t}{1+t}$.
	
	This property of the $w^*$-topology is a particular example of submetrizability. A topological space $(M,\tau)$ is called \emph{submetrizable} if there exists a $\tau$-continuous metric $d$ on $M$. It is elementary to check that the metric $d$ defined in~\eqref{AuxMetr} is continuous on $L^*\x L^*$ with the product of the $w^*$-topologies and thus $L^*$ is submetrizable. It is easy to see that whenever $(M,\tau)$ is a submetrizable topological space and $d$ is a $\tau$-continuous metric on $M\x M$ then $d$ metrizes the restriction of the $\tau$-topology on each $\tau$-compact subset $K\subs M$ and thus compact sets are also sequentially compact in submetrizable spaces. Note that if $M$ is $\s$-compact then $(M,d)$ is a separable metric space for any $\tau$-continuous metric $d$ on $M$, since $K_n$ is compact metric space in the restriction of the metric $d$, and thus separable, and $M=\bigcup_{n=1}^\infty K_n$. 
	
	By the Banach-Alaoglu theorem the closed balls $B_{L^*}(0,r):=\{\mu\in L^*|\|\mu\|_{L^*}\leq r\}$ of the dual $L^*$ are compact in the $w^*$-topology and thus $L^*$ is $\s$-compact in the $w^*$-topology as the increasing union of the compact subsets $K_n:=B_{L^*}(0,n)$, $n\in\NN$. Since the $w^*$-topology is metrizable on bounded subsets, bounded subsets are also sequentially relatively compact, i.e.~for any bounded sequence $\{\mu_n\}_{n=1}^\infty\subs L^*$ there exists subsequence $\{\mu_{k_n}\}_{n=1}^\infty$ of $\{\mu_n\}$ converging to some $\mu\in L^*$ in the $w^*$-topology.  
	
	In general, unless $L^*$ is separable (an assumption to restrictive for the applications) the Borel $\s$-algebra $\B_{w^*}\equiv\B_{w^*_L}$ generated on $L^*$ by the $w^*$-topology is smaller than the strong $\s$-algebra $\B_{L^*}$ of $L^*$. For example if $L=C(\T^d)$ with the uniform norm then $L^*=(\MMM(\T^d),\|\cdot\|_{TV})$ where $\|\cdot\|_{TV}$ is the total variation norm, and if $N\subs\T^d$ is a non-Borel subset of $\T^d$ then the set $\mathcal{N}:=\{\delta_u\in\MMM(\T^d)|u\in N\}\subs\MMM(\T^d)$ is strongly closed in $\MMM(\T^d)$ and thus $\mathcal{N}$ is strongly Borel. But the Dirac map $\delta\colon\T^d\to\PP\T^d\subs\MMM(\T^d)$ is $w^*$-continuous and thus also $(\B_{\T^d},\B_{w^*_L})$-measurable. Therefore $\mathcal{N}$ can not be in $\B_{w^*_L}$ since if it were, the set $N=\delta^{-1}(\mathcal{N})$ would be a Borel subset of $\T^d$.
	
	However, since 
	$$\|\mu\|_{L^*}=\sup_{\|f\|_L\leq 1}\mu(f),$$ the norm on the dual $L^*$ of any Banach space $L$ is $w^*$-lower semicontinuous as the supremum of $w^*$-continuous functionals $\J(f)$, $f\in L$ where $\J\colon L\to L^{**}$ is the canonical injection in the double dual. Therefore the closed balls $B_{L^*}(0,r):=\{\mu\in L^*\bigm|\|\mu\|_{L^*}\leq r\}$ are $w^*$-closed, and thus in $\B_{w^*}$, and the norm $\|\cdot\|_{L^*}$ is $\B_{w^*}$-measurable (i.e.~$(\B_{w^*},\B_{\RR})$-measurable). Therefore $\B_{w^*}$ contains also all open balls 
	$$D_{L^*}(0,r):=\{\mu\in L^*\bigm|\|\mu\|_{L^*}<r\}.$$
	This implies that $\A_{L^*}\subs\B_{w^*_L}$ where $\A_{L^*}$ is the $\s$-algebra generated by the collection of all strongly open balls $D_{L^*}(\mu,r)$, $\mu\in L^*$, $r>0$. However, unless $L^*$ is separable the inclusion $\A_{L^*}\subs\B_{L^*}$ is also in general strict.

As is customary, given a collection $\Xi$ of subsets of a set $M$ we will denote by $\s(\Xi)$ the $\s$-algebra generated by $\Xi$, i.e.~the smallest $\s$-algebra $\Sigma$ on $M$ that contains $\Xi$. Then if for any $K\subs M$ and any collection of subsets $\Xi\subs M$ we set $\Xi|_{K}:=\{A\cap K|A\in\Xi\}$ then $\s(\Xi)|_{K}=\s(\Xi|_{K})$ and if $K\in\s(\Xi)$ then 
\begin{equation}\label{salgrestrict}
	\s(\Xi|_K)=\s(\Xi)|_K\subs\s(\Xi).
\end{equation}

\begin{prop}\label{wstaranddmetrsalg} If $(M,\tau)$ is a $\s$-compact submetrizable topological space and $d$ is any continuous metric then the Borel $\s$-algebra of $M$ coincides with the Borel $\s$-algebra defined by the metric $d$. 
\end{prop}\textbf{Proof} Let $\B\equiv\B_{(M,\tau)}$ be the Borel $\s$-algebra of $M$ and let $\B_d$ denote the Borel $\s$-algebra of $(M,d)$. Since the topology of $d$ is weaker than $\tau$ we have that $\B_d\subs\B$ and so we have to prove the converse. So let $B\in\B$. Since $M$ is $\s$-compact there exists an increasing sequence $\{K_n\}_{n=1}^\infty\subs M$ of compact subspaces such that $M=\bigcup_{n=1}^\infty K_n$. We set $B_n=B\cap K_n$ so that $B=\bigcup_{n=1}B_n$ and it suffices to show that $B_n\in\B_d$ for all $n\in\NN$. Note that by definition $B_n\in\B|_{K_n}$ and that by~\eqref{salgrestrict} for any topology $\tau\supseteq\tau_d$ on $M$ and any $d$-measurable subset $K\subs M$
$$\B_{(K,\tau|_K)}=\s(\tau|_K)=\s(\tau)|_K=\B_{(M,\tau)}|_K\subs\B_{(M,\tau)}.$$
Therefore $B_n\in B|_{K_n}=\B_{(K_n,\tau|_{K_n})}=\B_{\tau_d|_{K_n}}=\B_d|_{K_n}\subs\B_d$ for all $n\in\NN$ as required.$\hfill\Box$\\

 For any family $\J$ of maps $f\colon M\to(\W,\F)$ defined on $M$ and with values in a measurable space $(\W,\F)$ we will denote by $\s(\J)$ the smallest $\s$-algebra $\Sigma$ with respect to which all maps $f\in\J$ are $(\Sigma,\F)$-measurable. As we will see next if $X$ is a separable Banach spaces then the Borel $\s$-algebra $\B_{(w^*_X)}$ of the $w^*$-topology coincides with $\s(\J(X))$, i.e.~the smallest $\s$-algebra $\Sigma$ for which all functionals $\J(f)\in X^{**}$, $f\in X$, are $\Sigma$-measurable. Here $\J\colon X\hookrightarrow X^{**}$ is the natural injection on the double dual. In other words $\s(\J(X))$-is the smallest $\s$-algebra $\Sigma$ on $X^*$ for which all $w^*$-continuous linear functions on $X^*$ are measurable. We will also denote by $\mathcal{B}\alpha_{w^*}(X^*)$ ($\mathcal{B}\alpha_d(X^*)$) the \emph{Baire $\s$-algebra} of $(X^*,w^*)$ ($(X^*,d)$) i.e.~the smallest $\s$-algebra on $X^*$ with respect to which all $w^*$-continuous ($d$-continuous) functions $F\colon X^*\to\RR$ are measurable.

\begin{prop}\label{WeakStarAndAuxSalg} For any dual $(X^*,w^*)$ of a separable Banach space $X$ and any countable subset $D=\{f_k\}_{k=1}^\infty\subs X$ dense in $X$ 
	$$\s(\J(D))=\s(\J(X))=\B\alpha_{w^*}(X^*)=\B_{w^*}(X^*)=\B_d(X^*)=\B\alpha_d(X^*),$$
	 where $d$ is the metric defined on $X^*$ as in~\eqref{AuxMetr}.
\end{prop}\textbf{Proof} Since any $F\in C_{w^*}(X^*)$ is $\B_{w^*}$-measurable we obviously have $\s(\mathcal{J}(D))\subs\s(\J(X))\subs\B\alpha_{w^*}(X^*)\subs\B_{w^*}$. Since $d$ is a continuous metric on $X^*$ and $X^*$ is a $\s$-compact submetrizable space it follows by Proposition~\ref{wstaranddmetrsalg} that $\B_{w^*}=\B_d$. Furthermore the Borel and Baire $\s$-algebras of any metric space coincide. Thus $\B_d=\B\alpha_{w^*}(X^*)$ and in order to complete the proof it suffices to prove that $\B_d\subs\s(\J(D))$. 

Since $\B_d$ is the Borel $\s$-algebra of a separable metric space it suffices to show that and open $d$-ball $D(\mu_0,\ee):=\{\mu\in X|d(\mu,\mu_0)<\ee\}$ of radius $\ee>0$ around $\mu_0\in X$ is in $\s(\J(X))$. For this, for each $k\in\NN$ we consider the semimetric $d_k$ on $X^*$ given by $d_k(\mu,\nu):=\sum_{i=1}^k\fr{2^i}\psi(|\ls f_i,\mu-\nu\rs|)$. Then $\{d_k\}_{k=1}^\infty$ is non-decreasing, $d_k\leq d$ and $\lim_{k\ra+\infty}d_k=d$ pointwise on $X^*\x X^*$. A sequence $\{\mu_n\}_{n=1}^\infty\subs X^*$ converges in the metric $d$ as $n\ra+\infty$ to some $\mu\in X^*$ if and only if it converges to $\mu$ in the semimetric $d_k$ for all $k\in\NN$. We note that the semimetrics $d_k$ satisfy the following property 
\begin{equation}\label{TouchdAndDie}
d_k(\mu,\nu)=d(\mu,\nu),\quad\mbox{for some }k\in\NN\quad\Lra\quad\mu=\nu.
\end{equation} 
Indeed, if $d_k(\mu,\nu)=d(\mu,\nu)$ then also $d_{k+n}(\mu,\nu)=d(\mu,\nu)$ for all $n\in\NN$. This implies that $\psi(|\ls f_{k+n},\mu-\nu\rs|)=0$ for all $n\in\NN$. Since $\{f_{k+n}\}_{n=1}^\infty$ is dense in $X$ it separates the points of $X^*$ and therefore $\mu=\nu$.

Next we note that if for each $\ee>0$ we define $B(\mu_0,\ee):=\{\mu\in\X^*|d(\mu,\mu_0)\leq\ee\}$ the closed $d$-ball of radius $\ee>0$ and by $D_k(\mu_0,\ee)$ the open $d_k$-ball of radius $\ee$ around then 
\begin{equation}\label{ClosedDBallAsInterOpenDkBall}
D(\mu_0,\ee)=\bigcup_{n=1}^\infty B\Big(\mu_0,\ee-\fr{n}\Big)\quad\mbox{and}\quad B(\mu_0,\ee)=\bigcap_{k=1}^\infty D_k(\mu_0,\ee)
\end{equation}
for all $\ee>0$ and in particular
\begin{equation}\label{ClosedDBallAsInterOpenDkBall2}D(\mu_0,\ee)=\bigcup_{n=1}^\infty B\Big(\mu_0,\ee-\fr{n}\Big)=\bigcup_{n=1}^\infty\bigcap_{k=1}^\infty D_k\Big(\mu_0,\ee-\fr{n}\Big).\end{equation} 
The left hand side equality in~\eqref{ClosedDBallAsInterOpenDkBall} is obvious so we show the right hand side equality. On one hand if $\mu\in\bigcap_{k=1}^\infty D_k(\mu_0,\ee)$ then $d_k(\mu,\mu_0)<\ee$ for all $k\in\NN$ and therefore $d(\mu,\mu_0)=\lim_{k\ra+\infty}d_k(\mu,\mu_0)\leq\ee$ and thus $\mu\in B(\mu_0,\ee)$. This shows that $\bigcap_{k=1}^\infty D_k(\mu_0,\ee)\subs B(\mu_0,\ee)$. For the converse inclusion it suffices to show that $B(\mu_0,\ee)\subs D_k(\mu_0,\ee)$ for all $k\in\NN$. To prove this let $k\in\NN$ and note that by~\eqref{TouchdAndDie}
$$B(\mu_0,\ee)\sm D_k(\mu_0,\ee)\subs\{\mu_0\}.$$
Indeed, suppose if  $\mu_k\in B(\mu_0,\ee)\sm D_k(\mu_0,\ee)$ then $$d_k(\mu_k,\mu_0)\leq d(\mu_k,\mu_0)\leq\ee\leq d_k(\mu_k,\mu_0)$$ which according to~\eqref{TouchdAndDie} implies that $\mu_k=\mu_0$. Since $\mu_0\in B(\mu_0,\ee)\cap D_k(\mu_0,\ee)$ it follows that $B(\mu_0,\ee)\subs D_k(\mu_0,\ee)$. This proves~\eqref{ClosedDBallAsInterOpenDkBall} and thus also~\eqref{ClosedDBallAsInterOpenDkBall2} holds. Consequently in order to show that $D(\mu_0,\ee)\in\s(\J(D))$ it suffices to show that $D_k(0,\ee)\in\s(\J(D))$ for all $k\in\NN$.

But $\s(\J(D))$ contains all sets of the form \begin{equation}\label{sj}[(\ls\J(f_i),\cdot\rs)_{i=1}^k]^{-1}(A)\subs X^*,\quad A\in\B_{\RR^k},\;k\in\NN\end{equation} where $(\ls\tau(f_i),\cdot\rs)_{i=1}^k\colon X^*\to\RR^k$ is the continuous linear vector functional defined by 
$$[(\ls\J(f_i),\cdot\rs)_{i=1}^k](\mu)=(\ls f_1,\mu\rs,\ldots,\ls f_k,\mu\rs)$$
and it we can easily see that $D_k(0,\ee)$ is such a set for all $k\in\NN$. Indeed, if for each $k\in\NN$ we define the continuous map $\psi_k=\psi_k^{\mu_0}\colon\RR^k\to\RR_+$ by 
$$\psi_k(t_1,\ldots,t_k)=\sum_{i=1}^k\fr{2^i}\psi(|t_i-|\ls f_i,\mu_0\rs|$$
then we can express the map $d_k(\cdot,\mu_0)$ as $d_k(\mu,\mu_0)=\psi_k^{\mu_0}(\ls f_1,\mu\rs,\ldots,\ls f_k,\mu\rs)$ and therefore 
$$D_k(\mu_0,\ee)=\big(\psi_k\circ[(\ls\J(f_i),\cdot\rs)_{i=1}^k]\big)^{-1}\big([0,\ee)\big)=[(\ls\J(f_i),\cdot\rs)_{i=1}^k]^{-1}\big(\psi_k^{-1}\big([0,\ee)\big)\big).$$ Since $\psi_k^{-1}([0,\ee))$ is a Borel subset of $\RR^k$ we have thus expressed $D_k(0,\ee)$ as a set of the form~\eqref{sj} which proves that $D_k(0,\ee)\in\s(\J(D))$ and completes the proof.$\hfill\Box$\\

A map $T\colon X^*\to Y^*$ is called \emph{a $w^*$-measurable operator} if it is $(\s(\J(X)),\s(\J(Y)))$-measurable. Equivalently $T$ is $w^*$-measurable if it maps $w^*$-measurable curves to $w^*$-measurable curves, i.e.~if for any measurable space $(\TT,\F)$ and any $w^*$-measurable map $\mu\colon\TT\to X^*$ the map $\TT\ni t\mapsto T(\mu_t)$ is $w^*$-measurable. In the case that $X,Y$ are separable this equivalent to $T$ being $(\B_{w^*},\B_{w^*})$-measurable. 

In the next example we see how the $w^*$-lower semicontinuity of the norm $\|\cdot\|_{X^*}$ can be used in the case that $X=C(\T)$ is the space of continuous maps on a compact metric space with the uniform norm and $\MMM(\T):=C(\T^d)^*$ to show that the variation map $|\cdot|\colon\MMM(\T)\to\MMM_+(\T)$ is $w^*$-measurable.

\begin{example} \label{VarMeas} Let $\T$ be a compact metric space. Then the variation map $|\cdot|\colon\MMM(\T)\to\MMM_+(\T)$ that assigns to each $\mu\in\MMM(\T)$ its variation $|\mu|\in\MMM_+(\T)$ is $w^*$-measurable. 
\end{example}
\textbf{Proof} If we show that for each Borel set $B\subs \T$ the map $I_B\colon\MMM(\T)\to\RR_+$ given by $I_B(\mu)=|\mu|(B)$ is $w^*$-Baire then for any simple function $\phi=\sum_{i=1}^na_i\1_{B_i}$, $B_i\in\B_\T$, $a_i\in\RR$, $n\in\NN$ the map $I_f\colon\MMM(\T)\to\RR$ given by $I_\phi(\mu)=\int\phi\df|\mu|$ is $w^*$-Baire and since any bounded function $f\in B(\T)$ can be approximated pointwise by a uniformly bounded sequence $\{\phi_n\}_{n=1}^\infty$ of simple functions it follows that the map $I_f\colon\MMM(\T)\to\RR$ given by $I_f(\mu)=\ls f,|\mu|\rs$ is $w^*$-Baire as the pointwise limit of a sequence of $w^*$-Baire functions. In particular the map $\ls f,|\mu|\rs\colon\MMM(\T)\to\RR$ is $w^*$-Baire for all $f\in C(\T)$ and thus the map $|\cdot|\colon\MMM(\T)\to\MMM_+(\T)$ is $w^*$-measurable.

Let now $\A\subs\B_\T$ be the collection of all sets $B\in\B_\T$ such that the map $I_B$ is $w^*$-measurable. Obviously $\emptyset\in\A$ and since the total variation norm $\|\cdot\|_{TV}$ is measurable also $\T\in\A$. Since $I_{\T\sm B}(\mu)=|\mu|(\T\sm B)=\|\mu\|_{TV}-|\mu|(B)=\|\cdot\|_{TV}-I_B$ the collection $\A$ is closed under complements and it is obviously closed under disjoint unions. Therefore $\A$ is a $\lambda$-system and by Dynkin's $\pi$-$\lambda$ theorem in order to show that $\A=\B_\T$ it suffices to show that $\A$ contains the $\pi$-system of all open sets.

So let $U\subs \T$ be an open set. In order to show that $U\in\A$ it suffices to show that 
\begin{equation}\label{VariationOfOpenSets}|\mu|(U)=\sup_{\substack{f\in C(\T)\\|f|\leq\1_U}}\ls f,\mu\rs =:\mu_0(U).
\end{equation}
Indeed, if~\eqref{VariationOfOpenSets} holds true, then since $C(\T)$ is separable there exists a countable family $D_U\subs C(\T)$ such that $|f|\leq\1_U$ for all $f\in D_U$ and $\big\{f\in C(\T)\bigm||f|\leq\1_U\big\}\subs\bbar{D_U}$ and thus then $|\mu|(U)=\sup_{f\in D_U}\ls f,\mu\rs$ which shows that the map $\MMM(\T)\ni\mu\mapsto|\mu|(U)\in\RR_+$ is $w^*$-lower semicontinuous and $w^*$-Baire as the supremum of the linear maps $\MMM(\T)\ni\mu\mapsto\ls f,\mu\rs\in\RR$, $f\in D_U$.
Obviously $\mu_0(U)\leq|\mu|(U)$ and thus in order to show~\eqref{VariationOfOpenSets} it suffices to show the converse inequality. So let $P\cup N$ be a Hahn decomposition of $\T$ with respect to $\mu$. By the argument in the proof of Proposition~\ref{Basic}(e) we can find a sequence $\{\phi_n\}\subs C(\T^d;[0,1])$ converging  $|\mu|$-a.s.~pointwise to $\1_P-\1_N$. Furthermore there exists a sequence of compact sets $K_n\subs U$ such that $|\mu|(U\sm K_n)\leq\fr{n}$ for all $n\in\NN$ and thus if $\{\psi_n\}\subs C(\T^d;[0,1])$ is any sequence of functions such that $\1_{K_n}\leq\psi_n\leq\1_U$ then $\{\psi_n\}$ converges $|\mu|$-a.s.~pointwise to $\1_U$. Then if we set $f_n:=\phi_n\psi_n$ the sequence $\{f_n\}$ converges $|\mu|$-a.s.~pointwise to $\1_U(\1_P-\1_N)$ and $|f_n|\leq\1_U$ for all $n\in\NN$ and thus 
$$\mu_0(U)\geq\lim_{n\ra+\infty}\ls f_n,\mu\rs=\int\1_{U\cap P}\df\mu-\int\1_{U\cap N}\df\mu=\mu^+(U)+\mu^-(U)=|\mu|(U)$$
which proves~\eqref{VariationOfOpenSets} and completes the proof.$\hfill\Box$

We will say that a net $\{T_\alpha\}_{\alpha\in\A}$ of linear operators $T_\alpha\colon X^*\to Y^*$, $\alpha\in\A$, \emph{$w^*$-converges pointwise on $X^*$ to }$T\colon X^*\to Y^*$ if the net $\{T_\alpha\mu\}_{\alpha\in\A}$ to $T\mu$ for all $\mu\in X^*$, i.e.~if
$$\lim_{\alpha}\ls g,T_\alpha\mu\rs=\ls g,T\mu\rs,\quad \forall g\in Y,\;\forall\mu\in X^*.$$ We will denote by $T=w^*\mbox{-}\lim_{\alpha}T_n$ \emph{the $w^*$-limit operator} of the net $\{T_\alpha\}$. Note that by the the $w^*$-semicontinuity of the norm $\|\cdot\|_{Y^*}$ in the $w^*$-topology it follows that if $T=w^*\mbox{-}\lim_{\alpha}T_\alpha$ then the operator norm $\|\cdot\|$ on the space $B(X^*,Y^*)$ of all bounded operators $T\colon X^*\to Y^*$ satisfies
$$\|T\|\leq\liminf_{\alpha}\|T_\alpha\|$$
and thus it is lower semicontinuous with respect to pointwise $w^*$-convergence of operators. As we will see the $w^*$-measurability of operators is preserved by pointwise $w^*$-convergence of sequences of operators.

\begin{prop}\label{PointLimOp} Let $X^*,Y^*$ be the duals of the Banach spaces $X,Y$. If $\{T_n\}_{n=1}^\infty$ is a sequence of $w^*$-measurable operators $T_n\colon X^*\to Y^*$ that $w^*$-conveges to $T\colon X^*\to Y^*$ pointwise in $X^*$, then $T$ is $w^*$-measurable.
\end{prop}\textbf{Proof} Since the $\s$-algebra $\s(\J(Y))$ is generated by sets of the form $\J(g)^{-1}(A)$, $g\in Y$, $A\in\B_\RR$ it suffices to show that $T^{-1}(\J(g)^{-1}(A))$ is in $\s(\J(X))$ for all $g\in Y$, $A\in\B_\RR$. Equivalently it suffices to show that the map $\J(g)\circ T\colon X^*\to\RR$ is $\s(\J(X))$-measurable. Since the operator $T_n$, $n\in\NN$, are $w^*$-measurable the maps $\J(g)\circ T_n$ is $\s(\J(X))$-measurable for all $n\in\NN$ and by the assumption that $\{T_n\}$ $w^*$-converges pointwise to $T$ it follows that for each $g\in Y$ and each $\mu\in X^*$ 
$$\lim_{n\ra+\infty}\J(g)\circ T_n(\mu)=\lim_{n\ra+\infty}\ls g,T\mu_n\rs=\ls g,T\mu\rs=\J(g)\circ T(\mu).$$
Thus $\J(g)\circ T$ is measurable as the pointwise limit of a sequence of measurable real valued maps.$\hfill\Box$

\begin{prop}\label{ContToMeas} Any $w^*$-continuous operator $T\colon X^*\to Y^*$ between duals of Banach spaces is $w^*$-measurable and $(\B_{w^*_X},\B_{w^*_Y})$-measurable.
\end{prop}\textbf{Proof} Let $T\colon X^*\to Y^*$ be a $w^*$-continuous operator. Since continuity implies Borel measurability we only have to show that $T$ is $(\s(\J(X)),\s(\J(Y)))$-measurable. Using the fact that for any Banach space $X$ we have $(X^*,w^*)^*= X\leq X^{**}$ it follows that an operator $T\colon X^*\to Y^*$ is $w^*$-continuous if and only if $T=T_0^*$ is the adjoint of a bounded operator $T_0\colon Y\to X$, in which case $T=T_0^*$. Thus $T$ is $w^*$-measurable since if $(\TT,\F)$ is any measurable space and $\mu\colon\TT\to X^*$ a measurable map, then $\ls g,T\mu_t\rs=\ls g,T_0^*\mu_t\rs=\ls T_0g,\mu_t\rs$ for all $g\in Y$ and all $t\in\TT$ and therefore the map $\TT\ni t\mapsto\ls g,T\mu_t\rs$ is measurable since $T_0g\in X$ and the map $\TT\ni t\mapsto\mu_t$ is $w^*$-measurable.$\hfill\Box$\\

We will denote by $B(X^*,Y^*)$ the space of all bounded linear operators, by $B_{w^*}(X^*,Y^*)$ its subspace consisting of all bounded $w^*$-measurable operators and by $BC_{w^*}(X^*,Y^*)$ the subspace of $B(X^*,Y^*)$ consisting of $w^*$-continuous operators. We will always consider the space $B(X^*,Y^*)$ equipped with the topology of pointwise $w^*$-convergence of operators. In terms of these spaces, Propositions~\ref{PointLimOp} and~\ref{ContToMeas} state that the subspace $B_{w^*}(X^*,Y^*)$ is a sequentially closed subspace of $B(X^*,Y^*)$ with respect to pointwise $w^*$-convergence of operators that contains the subspace $BC_{w^*}(X^*,Y^*)$. Thus if we define \emph{the space $\A_{w^*}(X^*,Y^*)$ of $w^*$-Baire measurable operators} $T\colon X^*\to Y^*$ as the sequential closure of $BC_{w^*}(X^*,Y^*)$ in $B(X^*,Y^*)$, i.e.~
$$\A(X^*,Y^*):=\bbar{BC_{w^*}(X^*,Y^*)}^{sq}:=\bigcap\big\{F\subs B(X^*,Y^*)\bigm| F\mbox{ seq.~closed}\mbox{ and }BC_{w^*}(X^*,Y^*)\subs F\big\}$$
we have by Propositions~\ref{PointLimOp} and~\ref{ContToMeas} that $\A_{w^*}(X^*,Y^*)\subs B_{w^*}(X^*,Y^*)$. The set $\A_{w^*}(X^*,Y^*)$ is indeed a subspace of $B(X^*,Y^*)$. To see this, for any subset $A\subs B(X^*,Y^*)$ we define 
$$[A]_{sq}:=\big\{T\in B(X^*,Y^*)\bigm|\exists\{T_n\}_{n=1}^\infty\subs A\mbox{~s.t.~}T_n\lra T\;w^*\mbox{-pointwise}\big\}.$$
If $A$ is a subspace, then so is $[A]_{sq}$. If we set $\bbar{BC_{w^*}(X^*,Y^*)}^1:=[BC_{w^*}(X^*,Y^*)]_{sq}$ and via transfinite induction we define $\bbar{BC_{w^*}(X^*,Y^*)}^{\xi+1}=[\bbar{BC_{w^*}(X^*,Y^*)}^\xi]_{sq}$ for any successor ordinal $\xi>1$ smaller than the first uncountable ordinal $\w_1$ and $\bbar{BC_{w^*}(X^*,Y^*)}^\xi=\bigcup_{\zeta<\xi}\bbar{BC_{w^*}(X^*,Y^*)}^\zeta$ for any limit ordinal $\xi<\w_1$, then 
$$\bbar{BC_{w^*}(X^*,Y^*)}^{sq}=\bigcup_{\xi<\w_1}\bbar{BC_{w^*}(X^*,Y^*)}^\xi,$$
and thus $\A_{w^*}(X^*,Y^*)=\bbar{BC_{w^*}(X^*,Y^*)}^{sq}$ is a subspace as the increasing union of subspaces. 
\begin{prop}\label{ContAnalyticComp}  Let $T\colon X^*\to Y^*$ be a $w^*$-Baire measurable operator and let $R^*\colon Y^*\to Z^*$ and $j^*\colon X_0^*\to X^*$ be $w^*$-continuous operators. Then the compositions $R^*\circ T$ and $T\circ j^*$ are $w^*$-Baire measurable.\end{prop}\textbf{Proof} The claim follows by transfinite induction and the fact that whenever $\{T_n\}_{n=1}^\infty\subs B_{w^*}(X^*,Y^*)$ is a sequence of operators $w^*$-converging pointwise to $T$ and $S_1\colon Y^*\to Z^*$ and $S_2\colon Z^*\to X^*$ are $w^*$-continuous operators then $S_1\circ T_n$ $w^*$-converges pointwise to $S_1\circ T$ and $T_n\circ S_2$ $w^*$-converges pointwise to $T\circ S_2$.$\hfill\Box$

\begin{example} If we regard $B(\T^d)$ as a subspace of $\MMM(\T^d)^*\cong(\MMM(\T^d),\|\cdot\|_{TV})^*$ via the injection $I\colon B(\T^d)\to\MMM(\T^d)^*$ defined by 
	$I(f)(\mu)=\int f\df\mu$ then 
	$$C(\T^d)=(\MMM(\T^d),w^*)^*\subs B(\T^d)\subs\A_{w^*}(X^*,\RR)\subs\MMM(\T^d)^*.$$
\end{example}


\subsection{Weak convergence on completely regular spaces}\label{WCCRS}
We recall that given a Borel probability measure $\mu$ on a topological space $M$ a Borel set $B\subs M$ is called \emph{$\mu$-Radon} if it is approximated from inside by compact subsets, i.e.
$$\mu(B)=\sup_{K\Subset B}\mu(K).$$ The measure $\mu$ is called a \emph{Radon measure} if every Borel subset $B$ of $M$ is $\mu$-Radon and \emph{weakly Radon} if every open set is $\mu$-Radon. The measure $\mu$ is called \emph{tight} if the whole space $M$ is $\mu$-Radon. We will denote by $\PP_RM$, $\PP_tM$ the spaces of all Radon and tight Borel probability measures on $M$ respectively. A topological space $M$ is called \emph{a Radon space} if every Borel probability measure $\mu$ on $M$ is a Radon measure. We start by proving that $X^*$ is a Radon space.

\begin{lemma}\label{WhenTightEqualsRadon} Let $(M,\tau)$ be a topological space. Then $\PP_tM=\PP_RM$ iff every compact subspace of $M$ is a Radon space. 
\end{lemma}\textbf{Proof} We suppose first that every compact subspace of $M$ is Radon and prove that $\PP_tM\subs\PP_RM$. So let $\mu\in\PP_tM$, $B\in\B_M$ be any Borel subset of $M$ and let $\ee\in(0,1)$. Then since $\mu$ is tight, there exists a compact subset $K$ of $M$ such $\mu(M\sm K)<\frac{\ee}{2}$. Then 
$$\B_K\equiv\s(A\cap K|A\in\tau)=\B_M\cap K\equiv\{B\cap K|B\in\B_M\}$$ and by assumption the subspace $K$ is a Radon space and therefore the probability measure $\mu_K:=\fr{\mu(K)}\mu|_{\B_K}\in\PP K$ is Radon. Therefore there exists a compact subset $F$ of $K\cap B$ such that 
$\mu_{K}([K\cap B]\sm F)<\frac{\ee}{2\mu(K)}$ and for which 
$$\mu(B\sm F)\leq\mu([K\cap B]\sm F)+\mu(M\sm K)<\ee.$$
It is easy to see 
that any compact subset of the space $K$ is compact subset of $M$, which since $\mu(B\sm F)<\ee$ and $\ee>0$ was arbitrary proves that $\mu$ is Radon.\\
\indent Conversely, suppose that $\PP_tM=\PP_RM$, let $K\subs M$ be compact and let $\mu\in\PP K$. The measure $\bar{\mu}(\cdot):=\mu(K\cap\,\cdot\,)\in\PP M$ is obviously tight and therefore by assumption it is Radon. Let now $B\in\B_K=\B_M\cap K$ and $\ee>0$. Since $\bar{\mu}$ is Radon, there exists a compact subset $F$ of $M$ such that $F\subs B\subs K$ and $\bar{\mu}(B\sm F)<\ee$. Then $F$ is also compact in $K$ and $\mu(B\sm F)=\bar{\mu}(B\sm F)<\ee$, since $B\subs K$. Therefore the arbitrary measure $\mu\in\PP K$ is Radon and thus the arbitrary compact subspace $K\subs M$ is a Radon space, which completes the proof.$\hfill\Box$ 
\begin{prop}\label{TightIsRadonInSubmetrizableSpaces} Let $(M,\tau)$ be a topological space. If $M$ is $\s$-compact then $\PP_tM=\PP M$ and if $M$ is submetrizable it holds that $\PP_tM=\PP_RM$.
\end{prop}\textbf{Proof} If $M$ is $\s$-compact then $M$ can be written as the countable increasing union $M=\bigcup_{n=1}^\infty K_n$ of a family $\{K_n\}_{n=1}^\infty$ of compact sets and thus any finite measure $\mu$ is tight by continuity from below. For the second claim, if $M$ is submetrizable there exists a continuous metric $d$ on $M$, which as we have seen metrizes the restriction of $\tau$ on every compact subset $K\subs M$. Consequently, every compact subspace of $M$ is metrizable, thus polish and thus Radon. Therefore $\PP_tM=\PP_RM$ by lemma~\ref{WhenTightEqualsRadon}$\hfill\Box$
\begin{cor}\label{DualOfSepRadon} Any probability measure on a $\s$-compact submetrizable space is a Radon measure. In particular the dual $X^*$ of a separable Banach space is a Radon space with the $w^*$-topology.
\end{cor}
\begin{prop} Let $M$ be a completely regular topological space and let $\mu,\nu\in\PP M$ be weakly Radon measures, such that 
	\begin{equation}\label{BCdeterminingInCompRegSubmetr}
		\int f\df\mu=\int f\df\nu,\qquad\forall\;f\in BC(M).
	\end{equation}
	Then $\mu=\nu$. 
\end{prop}\textbf{Proof} Since $\mu,\nu$ are Borel measures, it suffices to prove that $\mu(A)=\nu(A)$ for every open set $A$. But since $\mu,\nu$ are weakly Radon, for every open $A\subs M$ we have that $\mu(A)=\sup_{K\Subset A}\mu(K)$, and likewise for $\nu$, which shows that in order to prove that $\mu=\nu$ it suffices to prove that $\mu(K)=\nu(K)$ for compact subset $K$ of $M$.\\
\indent So let $K\subs M$ be compact. Since $M$ is completely regular, for every $x\in M\sm K$ there exists a function $f_x:M\lra[0,1]$ such that $f_x(x)=1$ and $f|_K\equiv 0$. We denote by $F(M)$ the set of all finite subsets of $M$, define an upwards directed set $\A:=\big\{\alpha\in F(M)|\alpha\cap K=\emptyset\big\}$ with order the set inclusion, and define the non-decreasing net $(f_\alpha)_{\alpha\in\A}\subs C(M;[0,1])\subs BC(M)$ by $f_\alpha=\max_{x\in\alpha}f_x$. Obviously $f_\alpha|_K\equiv 0$ for every $\alpha\in\A$ and $f_\alpha(x)=1$ for all $x\in\alpha$. Consequently, $f_\alpha\stackrel{\alpha}\lra 1-\1_K$ pointwise, since given $x\in M\sm K$, $1\geq f_\alpha(x)\geq f_x(x)=1$ for every $\alpha\geq\{x\}\in\A$ 
 and for every $x\in K$ we have $f_a(x)=0$ for all $\alpha\in\A$. Furthermore, this net is obviously increasing. In other words $\1_\alpha\leq f_\alpha\leq 1-\1_K$ for all $\alpha\in\A$ and $\1_\alpha\lra 1-\1_K$ pointwise, and
\begin{equation}\label{Trivial}
	\int f_\alpha\df\rho\leq \rho(M\sm K),\quad\mbox{for }\rho=\mu,\nu\mbox{ and }\alpha\in\A.
\end{equation}
\indent On the other hand, given $\ee>0$, for each $x\in M\sm K$ we have that $f_x(x)=1>1-\ee$ and therefore 
$$M\sm K\subs\bigcup_{x\in M\sm K}\{f_x>1-\ee\}.$$ 
Then, for any compact set $F\subs M\sm K$, the family $\mathcal{U_\ee}:=(\{f_x>1-\ee\})_{x\in M\sm K}$ is an open covering of $F$, and so there exist $n=n(F,\mathcal{U}_\ee)\in\NN$ and $x_1,\dots,x_n\in M\sm K$ such that 
$F\subs\bigcup_{k=1}^n\{f_{x_k}>1-\ee\}$. Then, for $\rho=\mu,\nu$, for all $\alpha\geq\alpha_\ee:=\{x_1,\dots,x_n\}\in\A$ 
\begin{equation*}
	\rho(F)\leq\rho\Big(\bigcup_{k=1}^n\{f_{x_k}>1-\ee\}\Big)\leq\rho\big(\{f_\alpha>1-\ee\}\big)\leq\fr{1-\ee}\int f_\alpha\df\rho.
\end{equation*} Therefore, since $\ee>0$ is arbitrary, for fixed $F\Subset M\sm K$ we have that $\rho(F)\leq\liminf_{\alpha\ra\infty}\int f_\alpha\df\rho$
and since $M\sm K$ is open and $\mu,\nu$ are weakly Radon taking the supremum over all $F\Subset M\sm K$, we get that 
$\rho(M\sm K)\leq\liminf_{\alpha\ra\infty}\int f_\alpha\df\rho$. Together with (\ref{Trivial}) this proves that 
$\rho(M\sm K)=\int f_\alpha\df\rho$ for $\rho=\mu,\nu$, which by assumption\eqref{BCdeterminingInCompRegSubmetr} implies that $\mu(M\sm K)=\nu(M\sm K)$, and thus $\mu(K)=\nu(K)$.$\hfill\Box$

\begin{lemma}\label{ContImpLowContInPortmanteauLemma1}
	Let $(M,\tau)$ be a completely regular topological space and let $f\in B(M)$ be a bounded function. Then $f$ is lower semicontinuous iff 
	\begin{equation}\label{LowSemContCharAsSupofCont}
		f=\sup_{h\in BC(M),\;h\leq f}h.\end{equation}
\end{lemma}\textbf{Proof} If~\eqref{LowSemContCharAsSupofCont} holds then $f$ is lower semicontinuous and in fact the complete regularity of $M$ is not required at for this implication. So we assume that $f$ is lower semicontinuous and we will show that~\eqref{LowSemContCharAsSupofCont} holds. We note first that we can make the additional assumption that $f\geq 0$. Indeed, if $m:=\inf_{x\in M}f(x)$ and the claim holds for non-negative functions, that 
$$f=m+(f-m)=m+\sup_{h\in BC(M),\;h\leq f-m}h=\sup_{h\in BC(M),\;h\leq f}h.$$ So in the rest of the proof we assume in addition that $f\geq 0$. The $\geq$ inequality in~\eqref{LowSemContCharAsSupofCont} is obvious and thus order to complete the proof it suffices to prove that $$f(x)\leq\sup_{h\in BC(M):h\leq f}h(x)$$ for all $x\in M$. Since we assume $f$ to be $\geq 0$ we obviously have that $\sup_{h\in BC(M),\;h\leq f}h\geq 0$ and therefore if $f(x)=0$ we have nothing to prove. So we fix $x\in M$ such that $f(x)>0$ and let $\ee>0\in(0,f(x)/2)$ be arbitrary. Since $f$ is lower semicontinuous, there exists an open neighbourhood $V_x$ of $x$ such that $f(V_x)\subs(f(x)-\ee,+\infty)$, and since $M$ is completely regular, there exists a continuous function $h_x:M\lra[0,f(x)-\ee]$ such that $h_x(x)=f(x)-\ee$ and $h_x|_{V_x^c}\equiv 0$. Then, $h_x\in BC(M)$ and $0\leq h_x\leq\big[f(x)-\ee]\1_{V_x}\leq f$, where the last inequality follows from the choice of the neighbourhood $V_x$. But then
$$f(x)=\ee+h_x(x)\leq\ee+\sup_{h\in BC(M):h\leq f}h(x).$$ So letting $\ee$ tend to zero we get that $f(x)\leq\sup_{h\in BC(M):h\leq f}h(x)$ as required. $\hfill\Box$\\

A Borel probability measure $\mu$ in a topological space $(M,\tau)$ is called $\tau$-{\it{smooth}} if for any upwards directed family $\{U_\alpha\}_{\alpha\in\A}$ of open sets we have that 
$$\mu\Big(\bigcup_{\alpha\in\A}U_\alpha\Big)=\sup_{\alpha\in\A}\mu(U_\alpha).$$

It is easy to see that any weakly Radon measure on a topological space $(M,\tau)$ is $\tau$-smooth. Indeed, let $\{U_\alpha\}_{\alpha\in\A}\subs\tau$ be an upwards directed family of open sets. We obviously have that $$\mu\Big(\bigcup_{\alpha\in\A}U_\alpha\Big)\geq\sup_{\alpha\in\A}\mu(U_\alpha).$$
For the converse inequality, let $\ee>0$ be arbitrary. Then $\bigcup_{\alpha\in\A}U_\alpha$ is open and since $\mu$ is weakly Radon there exists a compact set $K\subs\bigcup_{\alpha\in\A}U_\alpha$ such that $\mu\Big(\bigcup_{\alpha\in\A}U_\alpha\Big)\leq=\mu(K)+\ee$. Now, the family $\{U_\alpha\}$ covers the compact set $K$, and therefore there exists $\alpha_1,\dots,\alpha_n\in\A$ such that $K\subs\bigcup_{k=1}^nU_{\alpha_k}$. But since $\{U_\alpha\}$ is upwards directed, there exist $\alpha_0\in\A$ such that $\bigcup_{k=1}^nU_{\alpha_k}\subs U_{\alpha_0}$, which shows that 
$$\mu\Big(\bigcup_{\alpha\in\A}U_\alpha\Big)=\mu(K)+\ee\leq\mu(U_{\alpha_0})+\ee\leq\sup_{\alpha\in\A}\mu(U_\alpha)+\ee,$$
and proves the claim.

Since any tight measure in a metric space $(M,d)$ is a Radon measure it follows that any tight measure in a metric space is $\tau_d$-smooth where $\tau_d$ is the topology defined by the metric $d$.  

\begin{lemma}\label{ContImpLowContInPortmanteauLemma2} Let $(M,\tau)$ be a topological space and let $\mu\in\PP M$ be a $\tau$-smooth measure. Then, if $f:=\sup_{u\in\mathcal{U}}u$, where $\mathcal{U}$ is any upwards directed uniformly bounded family $\mathcal{U}$ of lower semicontinuous functions $u:M\lra\RR$, we have that 
	$$\int f\df\mu=\sup_{u\in\mathcal{U}}\int u\df\mu.$$
\end{lemma}\textbf{Proof} We note first that we can assume in addition that $0\leq f(x)<1$ for all $x\in M$. Indeed, suppose this is true and set $b:=\inf_{x\in M}f(x)\leq \sup_{x\in M}f(x)=:B$. Then for any $b'<b$, we have $f-b'>0$ and the function $\bar{f}:=\frac{f-b'}{B-b'+1}$ satisfies $0<\bar{f}(x)<1$ for all $x\in M$ and $\bar{f}=\sup_{\bar{u}\in\bar{\mathcal{U}}}\bar{u}$ where $\bar{\mathcal{U}}=\{\frac{u-b'}{B-b'+1}|u\in\mathcal{U}\}$. Then,
$$\int f\df\mu=b'+(B-b'+1)\int\bar{f}\df\mu=b+(B-b'+1)\sup_{\bar{u}\in\mathcal{U}}\int\bar{u}\df\mu=\sup_{u\in\mathcal{U}}\int u\df\mu.$$ 

So in what follows we assume that $f(M)\subs(0,1)$ and let $\ee>0$ be arbitrary. We have to prove that 
$\int f\df\mu\leq\ee+\sup_{u\in\mathcal{U}}\int u\df\mu$. For each $n\in\NN$ we have 
$$\int f\df\mu\leq\sum_{k=0}^{n-1}\frac{k+1}{n}\mu\Big\{\frac{k}{n}<f\leq\frac{k+1}{n}\Big\}=\fr{n}\sum_{k=0}^{n-1}\mu\Big\{f>\frac{k}{n}\Big\}=\fr{n}+\fr{n}\sum_{k=1}^{n-1}\mu\Big\{f>\frac{k}{n}\Big\}.$$ 
We fix $n>2/\ee$. Since $f=\sup_{u\in\mathcal{U}}u$, we have that $\{f>\frac{k}{n}\}=\bigcup_{u\in\mathcal{U}}\{u>\frac{k}{n}\}$ for each $k=1,\dots,n-1$. But since each $u\in\mathcal{U}$ is lower semicontinuous, for each $k=1,\dots,n-1$ the set $U_u^k:=\{u>\frac{k}{n}\}$ is open, and the family $U^k:=\{U_u^k\}_{u\in\mathcal{U}}$ is an upwards directed family of open sets for each fixed $k=1,\dots,n-1$. Therefore since $\mu$ is $\tau$-smooth we have that 
$$\mu\Big\{f>\frac{k}{n}\Big\}=\sup_{u\in\mathcal{U}}\mu\Big\{u>\frac{k}{n}\Big\}$$ for all $k=1,\dots,n-1$, and so for each $k=1,\dots,n-1$ we can choose $u_k\in\mathcal{U}$ such that $$\mu\Big\{u_k>\frac{k}{n}\Big\}>\mu\Big\{f>\frac{k}{n}\Big\}-\frac{\ee}{2}.$$
Then, since $\mathcal{U}$ is upwards directed, there exists $u_0\in\mathcal{U}$ such that $u_0\geq u_1\mx\dots\mx u_{n-1}$, and 
\begin{align*}
	\int fd\mu&\leq\frac{\ee}{2}+\fr{n}\sum_{k=1}^{n-1}\mu\Big\{f>\frac{k}{n}\Big\}
	\leq\frac{\ee}{2}+\frac{n-1}{n}\frac{\ee}{2}+\fr{n}\sum_{k=1}^{n-1}\mu\Big\{u_k>\frac{k}{n}\Big\}\leq\ee+\fr{n}\sum_{k=1}^{n-1}\mu\Big\{u_0>\frac{k}{n}\Big\}\\
	&=
	\ee+\fr{n}\sum_{k=1}^{n-1}k\mu\Big\{\frac{k}{n}<u_0\leq\frac{k+1}{n}\Big\}\leq\ee+\int u_0d\mu\leq\ee+\sup_{u\in\mathcal{U}}\int u d\mu.
\end{align*}
Since $\ee>0$ was arbitrary this concludes the proof.$\hfill\Box$\\

Next we state the portmanteau theorem in completely regular topological spaces. It is known~\cite[Theorem 8.1]{Topsoe1970a} that the well-known characterizations of the weak convergence of nets $(\mu_\alpha)_{\alpha\in\A}\subs\PP M$ given in the portmanteau theorem in polish spaces remain valid in the more general context of completely regular topological spaces $(M,\tau)$, provided the limiting measure $\mu$ is $\tau$-smooth. Since any tight measure in a metric space is Radon and thus smooth, the portmanteau theorem is valid in any metric space under the assumption that the limiting measure $\mu$ is tight and since any measure on the $w^*$-dual of a separable Banach space $X$ is Radon by Corollary~\ref{DualOfSepRadon}, the portmanteau theorem holds in the space $\PP(X^*,w^*)$ without any assumptions on the limiting measure $\mu$.

\begin{prop}{\rm{(The portmanteau theorem)}} Let $(M,\tau)$ be a completely regular topological space, let $(\mu_\alpha)_{\alpha\in\A}$ be a net in $\PP M$, and let $\mu\in\PP M$ be a $\tau$-smooth measure. Then the following are equivalent:
	\begin{itemize}\item[{\rm{(a)}}] $\mu_\alpha\lra \mu\in\PP M$ weakly.
		\item[{\rm{(b)}}] For every closed set $F\subs M$, $\limsup_\alpha\mu_\alpha(F)\leq\mu(F)$.
		\item[{\rm{(c)}}] For every open set $U\subs M$, $\liminf_\alpha\mu_\alpha(U)\geq\mu(U)$.
		\item[{\rm{(d)}}] For every $\mu$-continuous set $A\subs M$, i.e$.$ for every Borel set $A\subs M$ such that $\mu(\partial A)=0$, it holds that $\lim_\alpha\mu_\alpha(A)=\mu(A)$.
		\item[{\rm{(b$'$)}}] For every bounded upper semicontinuous function $f:M\lra[-\infty,\infty)$, $$\limsup_\alpha\int{fd\mu_\alpha}\leq\int{fd\mu}.$$
		\item[{\rm{(c$'$)}}] For every bounded lower semicontinuous function $f:M\lra(-\infty,\infty]$, $$\liminf_\alpha\int{fd\mu_\alpha}\geq\int{fd\mu}.$$
		\item[{\rm{(d$'$)}}] For evert bounded $\mu$-a.s$.$ continuous function, $\lim_\alpha\int{fd\mu_\alpha}=\int{fd\mu}$.
	\end{itemize}
\end{prop}\textbf{Proof} Since a Borel subset $A\subs X$ is closed, open and $\mu$-continuous iff $\1_A$ is lower semicontinuous, upper semicontinuous and $\mu$-a.s$.$ continuous respectively, it follows that $(x')$ implies $(x)$, for $x=b,c,d.$ Furthermore, (b) is equivalent to (c), and (b$'$) is equivalent to (c$'$). Finally it is obvious that (d$'$) implies (a), and therefore it suffices to prove that (a)$\Lra$(c$'$), (b)$\mn$(c)$\Lra$(d), and that (d)$\Lra$(d$'$).\\
\noindent$(a)\Lra(c')$ Let $f:X\lra(-\infty,\infty]$ lower semicontinuous and bounded. By lemmas \ref{ContImpLowContInPortmanteauLemma1} and \ref{ContImpLowContInPortmanteauLemma2} $$\int fd\mu=\sup\left\{\int hd\mu\,\Big|\,h\in BC(M),\;h\leq f\right\}.$$
which as we can easily see implies that $\liminf_\alpha\int{fd\mu_\alpha}\geq\int{fd\mu}$.\\
\noindent$(b)\mn(c)\Lra(d)$ We note first that a Borel set $A\subs X$ is a $\mu$-continuous set iff $\mu(A^o)=\mu(A)=\mu(\bbar{A})$. So if $A$ is an $\mu$-continuous set, by $(b)$ and $(c)$ we have that
\begin{align*}
	\mu(A^o)\leq\liminf\mu_n(A^o)\leq\liminf\mu_n(A)\leq\limsup\mu_n(A)\leq\limsup\mu_n(\bbar{A})\leq\mu(\bbar{A}),
\end{align*}
which according to the initial remark proves (d).\\
\noindent$(d)\Lra(d')$ Let $f:X\lra\RR$ be a bounded, $\mu$-a.s$.$ continuous function and let $\ee>0$. Let $M_0\in\mathcal{B}_M$ be a full measure set, $\mu(M_0)=1$, of continuity points of $f$ and let $a,b\in\RR$ such that $a<f(x)<b$ for all $x\in M$. For each $r\in(a,b)$, we set $F_r:=\{x\in X\,|\,f(x)=r\}$. The family $\{F_r\}_{r\in(a,b)}$ is a partition of $M$, and thus for every finite subset $I$ of $(a,b)$ we have that $\sum_{r\in I}\mu(F_r)=\mu\big(\bigcup_{r\in I}F_r\big)\leq 1$ and thus
$$\sum_{r\in(a,b)}\mu(F_r)\leq 1<+\infty.$$ Consequently the set of all $r\in(a,b)$ for which $\mu(F_r)>0$ is at most countable. There exists then a partition $a=a_0<a_1<\dots<a_n=b$ of the interval $(a,b)$, such that $a_i-a_{i-1}<\ee$, $i=1,\dots,n$ and $\mu(F_{a_i})=0$, $i=0,\dots,n$. We set $E_i:=f^{-1}\big([a_{i-1},a_i)\big)$, $i=1,\dots,n$, and define the simple functions
$\phi=\sum_{i=1}^{n}{a_{i-1}\1_{E_i}}$ and $\psi=\sum_{i=1}^{n}{a_i\1_{E_i}}$. Obviously, $\phi\leq f\leq\psi$ and $\psi-\phi\leq\ee$. Also $\partial E_i\subs F_{a_{i-1}}\cup F_{a_i}\cup(M\sm M_0)$, for all $i=1,\dots,n$ and thus the $E_i$'s are $\mu$-continuous sets. Then $\lim\int{\phi d\mu_n}=\int{\phi d\mu}$ and $\lim\int{\psi d \mu_n}=\int{\psi d\mu}$ by $(d)$ and thus
\begin{align*}
	\int f\df\mu-\ee\leq\int\phi\df\mu\leq\liminf_n\int f\df\mu_n\leq\limsup_n\int f\df\mu_n\leq\int\psi\df\mu\leq \int f\df\mu+\ee.
\end{align*}
Since $\ee>0$ was arbitrary, the claim follows.$\hfill\Box$\\

Next we state the generalization of Prokhorov's relative compactness criterion on metric spaces to completely regular topological spaces, originally due to Le Cam~\cite{Cam1957a}. Prokhorov's criterion is valid in any metric space $M$ and states that a uniformly tight family $\mathcal{M}\subs\PP M$ of probability measures is relatively compact in the weak topology. It is usually stated in separable metric spaces, e.g~\cite{Billingsley1971a} but it is valid in any metric space. Indeed, if $\MMM$ is uniformly tight, then there exists a separable closed subspace $M_0$ such that $\mu(M_0)=1$ for all $\mu\in\MMM$. Then the family $$\MMM_0:=\big\{\mu|_{\B_{M_0}}\bigm|\mu\in\MMM\big\}\subs\PP M_0$$ is a uniformly tight family of probability measures in the separable space $M_0$ and thus given any sequence $\{\mu_n\}\subs\MMM$ there exists a subsequence $\{\mu_{k_n}\}$ of $\{\mu_n\}$ and $\mu_0\in\PP M_0$ such that $\mu_{k_n}|_{\B_{M_0}}\lra\mu_0$ as $n\ra+\infty$. But then for the measure $\mu\in\PP_tM$ defined by $\mu(B)=\mu_0(B\cap M_0)$ we have that $\mu_{k_n}\lra\mu$ weakly. If the metric space is complete, the converse is also true, i.e.~if a family $\MMM\subs\PP_tM$ is relatively compact then it is uniformly tight. 

\begin{theorema}\label{ProkhorovLeCam}{\rm{(Prokhorov-Le Cam)}} Let $(M,\tau)$ be a completely regular topological space. Then any uniformly tight family $\MMM\subs\PP_RM$ of probability measures is relatively compact in $\PP_RM$ in the weak topology. If $M$ is in addition submetrizable then any uniformly tight family $\MMM\subs\PP_tM$ is also sequentialy relatively compact in $\PP_tM$ in the weak topology.
\end{theorema}\textbf{Proof} For the proof of the first assertion see \cite[Ch.~3, Th.~59]{Dellacherie1978a}. The second assertion was proved in~\cite{Cam1957a}. Since the latter assertion will be frequently used in the text and the original reference is easily accesible only in Russian we will give a direct proof by using the submetrizability of $M$ and Prokhorov's theorem in metric spaces. 

So let $\MMM\subs\PP_tM$ be a uniformly tight family of probability measures and let $\{\mu_n\}_{n=1}^\infty\subs\MMM$ be a sequence in $\MMM$. We need to exhibit a subseqeunce $\{\mu_{k_n}\}_{n=1}^\infty$ of $\{\mu_n\}$ and $\mu\in\PP_tM$ such that $\mu_{k_n}\lra\mu$ in the weak topology of $\PP_t(M,\tau)$. Since $M$ is submetrizable there exists a continuous metric $d\colon M\x M\to\RR$, which necessarily metrizes the restriction of the topology $\tau$ on any $\tau$-compact subspace $K\subs M$. Since the family $\MMM$ is uniformly tight with respect to the topology $\tau$ and any $\tau$-compact set is $d$-compact, it follows that $\MMM$ is also tight in the metric space $(M,d)$. Thus by Prokhorov's theorem on metric spaces there exists a $\mu\in\PP_t(M,d)$ and a subsequence $\{\mu_{k_n}\}_{n=1}^\infty$ of $\{\mu_n\}$ such that 
$\mu_{k_n}\lra\mu$ in the weak topology of $\PP(M,d)$. In particular by the portmanteau theorem on metric spaces 
\begin{equation}\label{AuxWeak}
	\limsup_{n\ra+\infty}\mu_{k_n}(F)\leq\mu(F),\quad\forall\;d\mbox{-closed }F\subs M.
\end{equation}
If we can show that 
\begin{equation}\label{UpperWeakConvergence}
	\limsup_{n\ra+\infty}\mu_{k_n}(F)\leq\mu(F),\quad\forall\;\tau\mbox{-closed }F\subs M
\end{equation} and that 
$\mu$ is $\tau$-smooth it will follow by the portmanteau theorem for completely regular spaces that $\mu_{k_n}\lra\mu$ and the proof will be complete. 

We show first~\eqref{UpperWeakConvergence}. So let $F\subs M$ be $\tau$-closed and let $\ee>0$. Since $\MMM$ is uniformly $\tau$-tight there exists a $\tau$-compact set $K_\ee\subs M$ such that 
\begin{equation}\label{TightKEps}
	\sup_{n\in\NN}\mu_{k_n}(M\sm K_\ee)<\ee.
\end{equation}
Then $F\cap K_\ee$ is $\tau$-compact and thus it is also $d$-closed. Thus by~\eqref{AuxWeak}
\begin{equation}\label{AuxWeakCons}\limsup_{n\ra+\infty}\mu_{k_n}(F\cap K_\ee)\leq\mu(F\cap K_\ee)\leq\mu(F).\end{equation}
But $F\sm(F\cap K_\ee)=F\cap(M\sm K_\ee)$ and thus by~\eqref{TightKEps} $\mu_{k_n}(F)\leq\mu_{k_n}(F\cap K_\ee)+\ee$
for all $n\in\NN$. Taking the limit superior as $n\ra+\infty$, it follows by~\eqref{AuxWeakCons} that $\limsup_{n\ra+\infty}\mu_{k_n}(F)\leq\mu(F)+\ee$, which since $\ee>0$ was arbitrary, proves~\eqref{UpperWeakConvergence}.

It remains to check that the measure $\mu$ is $\tau$-smooth. Since the space $M$ is assumed submetrizable, by lemma~\ref{WhenTightEqualsRadon} it suffices to show that $\mu$ is $\tau$-tight, since then it is Radon and thus $\tau$-smooth. But this follows from the uniform $\tau$-tightness of $\MMM$. Indeed, given $\ee>0$ there exists a $\tau$-compact set $K_\ee\subs M$ such that $\inf_{n\in\NN}\mu_{k_n}(K_\ee)>1-\ee$ and thus $\mu(K_\ee)>1-\ee$ by~\eqref{UpperWeakConvergence}. This proves that $\mu$ is $\tau$-tight and completes the proof.$\hfill\Box$\\

In the case that $M$ is completely regular and submetrizable we do not need to assume the family $\MMM$ to consist of Radon measures due to corollary \ref{TightIsRadonInSubmetrizableSpaces}.

\subsection{The dual of $L^1(\TT,X)$}

\indent Let $(\TT,\F,m)$ be a complete finite measure space and let $X$ be a Banach space. We are mainly interested in the case where $\TT=[0,T]$ for some a finite time horizon $T>0$ where we regard the interval $[0,T]$ as a complete measure space equipped with the Lebesgue measure. The \emph{Bochner $L^p$-space} $L^p(\TT;X)$, $1\leq p\leq +\infty$, is the vector space of all strongly measurable maps $f\colon\TT\to X^*$ (i.e.~a.e$.$ pointwise limits of simple functions) equipped with the norm $$\|f\|_{L^p(\TT;X)}:=\big\|\|f_\cdot\|_{X}\|_{L^p(\TT)}.$$ Here in the right hand side we denote by $\|f_\cdot\|_{X}\in L^p(0,T)$ the map $\TT\ni t\mapsto\|f_t\|_{X}$. In the context of Banach valued-maps the duality $L^p(\TT;X)^*=L^q(\TT;X^*)$ where $1\leq p<+\infty$ and $1<q\leq+\infty$ are conjugate exponents holds if and only if $X^*$ has the Radon-Nikodym property with respect to $m$~\cite[Theorem IV.1]{Diestel1977a}. The representation of $L^\infty(\TT;X^*)\cong L^1(\TT;X)^*$ as a dual would induce a $w^*$-topology on $L^\infty(\TT;X^*)$ which is very convenient in proving that the laws of various empirical processes $\s\colon D(0,T;\MM_N^d)\to L^\infty(0,T;X^*)$ of the ZRP are relatively compact. However for the empirical processes under study the space $X$ will be such that the dual $X^*$ is not separable, for example $X=C(\T^d)$, which implies that $X^*$ does not have the Radon-Nikodym property and the duality $L^p(0,T;X)^*=L^q(0,T;X^*)$ does not hold.

A Banach space $X$ is said to have the \emph{Radon-Nikodym property with respect to $m$} if for any $X$-valued measure $\bs{\nu}\colon\F\to X$ that is \emph{$m$-absolutely continuous}, i.e.~for any $\ee>0$ there exists $\delta>0$ such that 
\begin{equation}\label{ACBVM}
m(A)<\delta\quad\Lra\quad\|\bs\nu(A)\|_X<\ee,
\end{equation}
there exists a Bochner integrable function $f\in L^1(\TT;X)$ such that $$\bs{\nu}(E)=\int_Ef\df m ,\quad\forall\;E\in\F.$$ The Banach space $X$ has the \emph{Radon-Nikodym property} if it has the Radon-Nikodym property for any finite measure space $(\TT,\F,m)$. It is known (see for instance \cite[Section 11]{Fabian2010a} where geometric characterizations of the Radon-Nikodym property are given) that a Banach space has the Radon-Nikodym property if and only if it has the Radon-Nikodym property with respect to the Lebesgue measure on $[0,1]$. As proved by Uhl~\cite{Uhl1972a} and Stegall~\cite{Stegall1975a} a dual space $X^*$ has the Radon-Nikodym property if and only if for any separable subspace $Y$ of $X$, the dual $Y^*$ is separable.

In particular if $X^*$ is not separable, as will be the case for the empirical processes of the ZRP then $X^*$ does not have the Radon-Nikodym property with respect to the Lebesgue interval $[0,T]$ and thus the inclusion $L^\infty(0,T;X^*)\hookrightarrow L^1(0,T;X)^*$ is strict. However, we can one can always describe the elements of $L^1(0,T;X)^*$ via curves taking values in $X^*$ when $X^*$ does not have the Radon-Nikodym property by relaxing strong measurability to $w^*$-measurability as described for example in \cite{PilarCembranos1997a}. This can be done since a $w^*$-measurable Radon-Nikodym derivative always exists. Our goal in this section is to give a description of the dual space $L^1(\TT;X)^*$ for general Banach spaces $X$ following \cite{Diestel1977a,PilarCembranos1997a} and \cite{Fabian2010a}.

\subsubsection{Weak-star $L^\infty$-spaces}\label{lwstarinfty}

Let $\TT=(\TT,\F,m)$ be a finite measura space and let $\LL_{w^*}(\TT;X^*)$ denote the linear space of all $w^*$-measurable maps $\mu\colon(\TT,\F)\to X^*$. We recall that $\mu$ is $w^*$-measurable if and only it is $\s(\J(X))$-measurable where $\J\colon X\hookrightarrow X^{**}$ is the canonical injection in the double dual. As is customary we use the calligraphic $\LL$ to denote that we have not identified a.e.~equal functions. In section~\ref{Appendix1} we have seen that the norm $\|\cdot\|_{X^*}$ is only $\B_{w^*}$-measurable. However, when $X$ is separable, $\s(\J(X))=\B_{w^*}$ by Proposition~\eqref{WeakStarAndAuxSalg}, and in this case the map $\|\mu_\cdot\|_{X^*}$ is measurable for all $\mu\in\LL_{w^*}(\TT;X^*)$. Thus for separable $X$ we define 
$$\LL_{w^*}^q(\TT;X^*):=\Big\{\mu\in\LL_{w^*}(\TT;X^*)\Bigm|\big\|\|\mu_\cdot\|_{X^*}\big\|_{L^q(\TT)}<+\infty\Big\}$$
for each $q\in[1,+\infty]$ and set $L^q_{w^*}(\TT;X^*)$ the quotient space modulo the relation of $m$-a.s.~equality.

  In general, unless $X^*$ is separable, the $\s$-algebra $\s(J(X))$ is smaller than the Borel $\s$-algebra $\B_{w^*}$ of the $w^*$-topology. Thus we can not conclude that $\|\mu_\cdot\|_{X^*}\colon\TT\to\RR_+$ is measurable. However one can still define the linear subspace of $\LL_{w^*}(\TT;X^*)$ consisting of $L^q$-maps as the space $\LL_{w^*}^q(\TT;X^*)$ of all $w^*$-measurable functions $\mu\in\LL_{w^*}(\TT;X^*)$ such that 
\begin{equation}\label{NonMeasSolv}\A_\mu:=\big\{g\in L^\infty(\TT)\bigm|\|\mu_t\|_{X^*}\leq g(t)\mbox{ a.s.-}\forall\;t\in\TT\big\}\neq\emptyset.
\end{equation} and the seminorm $\|\cdot\|_{\LL_{w^*}^\infty(\TT;X*)}$ on $\LL_{w^*}^q(\TT;X^*)$ given by 
$$\|\mu\|_{\LL_{w^*}^q(\TT;X^*)}=\inf_{g\in\A_\mu}\|g\|_{L^q(\TT)}.$$ The kernel $$\mathcal{N}_{X^*}:=\big\{\mu\in\LL_{w^*}^\infty(\TT;X^*)\bigm|\|\mu\|_{\LL_{w^*}^\infty(\TT;X^*)}=0\big\}$$ of the seminorm $\|\cdot\|_{\LL_{w^*}^\infty(\TT;X^*)}$ coincides with the subspace of maps $\mu\in\LL_{w^*}^\infty(\TT;X^*)$ that vanish on a measurable subset $E\in\F$ of full measure. 

Then the bilinear map $\lls\cdot,\cdot\rrs\colon L^1(\TT;X)\x L_{w^*}^\infty(\TT;X^*)\lra\RR$ given by 
$$\lls f,\mu\rrs=\int_\TT\ls f_t,\mu_t\rs\df m(t)$$ is well defined, since for all $(f,\mu)\in L^1(\TT;X)\x\LL_{w^*}^\infty(\TT;X^*)$ the map $\TT\ni t\mapsto\ls f_t,\mu_t\rs=:\ls f,\mu\rs_t$, denoted by $\ls f,\mu\rs$,  does not depend on the representatives of the $m$-a.e.~equality classes of $f$ and $\mu$, it is in $L^1(\TT)$ and 
\begin{equation}\label{CSWL1L8}
|\lls f,\mu\rrs|\leq\|f\|_{L^1(\TT;X)}\|\mu\|_{\LL^\infty_{w^*}(\TT;X^*)}.
\end{equation} 

The bilinear pairing between $L^1(\TT;X)$ and $L_{w^*}^\infty(\TT;X^*)$ induces a linear operator $S\colon L_{w^*}^\infty(\TT;X^*)\to L^1(\TT;X)^*$ via $S(\mu)(f)=\lls f,\mu\rrs$. By~\eqref{CSWL1L8} the operator $S$ is a contraction. Consider in $L^\infty_{w^*}(\TT;X^*)$ the relation $\backsim_S$ given by $\mu\backsim_S\nu$ if and only if $\mu-\nu\in\ker S$, i.e. 
\begin{equation}\label{DLinftyEqRel1}\mu\backsim_S\nu\quad\mbox{iff}\quad\lls f,\mu\rrs=\lls f,\nu\rrs\mbox{ for all }f\in L^1(\TT;X).\end{equation}
As we will see the relation $\backsim_S$ is equivalent to the relation $\backsim_{w^*\mbox{-}m}$ of \emph{$w^*$-m-a.s.~equality}, i.e.
\begin{equation}\label{DLinftyEqRel2}\mu\backsim_{w^*\mbox{-}m}\nu\quad\mbox{iff}\quad m\big\{t\in\TT\bigm|\ls f,\mu_t\rs=\ls f,\nu_t\rs\big\}=0,\forall\;f\in X.\end{equation}
Indeed, if $\mu\backsim_S\nu$ then in particular for any $A\in\F$ and any $f\in X$
$$\int_A\mu_t(f)\df m(t)=\int_\TT\mu_t(f\1_A(t))\df m(t)=\lls f\1_A,\mu\rrs=\lls f\1_A,\nu\rrs=\int_A\nu_t(f)d\df m(t)$$
which implies that $\mu_t(f)=\nu_t(f)$ for almost all $t\in\TT$ and thus $\mu\backsim_{w^*\mbox{-}m}\nu$. Conversely, if $\mu\backsim_{w^*\mbox{-}m}\nu$ holds then it is easy to see that $\lls\phi,\mu\rrs=\lls\phi,\nu\rrs$ for all all simple functions $\phi=\sum_{i=1}^nf_i\1_{A_i}$, $f_i\in X$, $A_i\in\F$ and thus $\mu\backsim_S\nu$.

By taking the quotient of the space $L_{w^*}^\infty(\TT;X^*)$ with respect to the subspace $\ker S$ the map $S$ passes to an injection $S\colon\bar{L}_{w^*}^\infty(\TT;X^*):=(\,^{L_{w^*}^\infty(\TT;X^*)}/_{w^*\mbox{-}m\mbox{-}a.s.})\to L^1(\TT;X)^*$. In the case that the Banach space $X$ is separable the operator $S$ is an injection and we do not need to take the quotient with $\ker S$. As we will show in the following two sections the map $S$ is also surjective and thus for separable $X$ the dual $L^1(\TT;X)^*$ is isometric to $L_{w^*}^\infty(\TT;X^*)$ and in the case that $X$ is not separable $L^1(\TT;X)^*$ is isometric to $\bar{L}_{w^*}^\infty(\TT;X^*)$.
\begin{prop}\label{RelationOnSeparable} If the Banach space $X$ is separable then the equivalence relation of $w^*$-$m$-a.s.~equality in $L^\infty_{w^*}(\TT;X^*)$ coincides with $m$-a.s.~equality in $L_{w^*}^\infty(\TT;X^*)$ and thus $S\colon L^\infty_{w^*}(\TT;X^*)\to L^1(\TT;X)^*$ is an injection.
\end{prop}\textbf{Proof} Let $\mu,\nu\in L_{w^*}^\infty(\TT;X^*)$ be such that $\mu\backsim_S\nu$. Since $X$ is separable there exists a countable subset $D\subs X$ dense in $X$. Then for any $f\in D$ there exists a set $E_f\in\F$ with $m(E_f)=m(\TT)$ and $\ls f,\mu_t\rs=\ls\nu_t\rs$ for all $t\in E_f$. Then set $E:=\bigcap_{f\in D}E_f$ is of full $m$-measure in $\TT$ and $\ls f,\mu_t\rs=\ls f,\nu_t\rs$ for all $f\in D$, $t\in E$. But since $D$ is dense in $X$ this implies that $\ls f,\mu_t\rs=\ls f,\nu_t\rs$ for all $f\in X$ and all $t\in E$ and thus $\mu_t=\nu_t$ for all $t\in E$ which proves that $\mu=\nu$ in $L_{w^*}^\infty(\TT;X^*)$.$\hfill\Box$
\subsubsection{Banach-valued measures}\label{BVM}

A $X$ valued set function $\bs{\nu}\colon\F\to X$ is called a \emph{Banach-valued measure} if $\bs{\nu}(\emptyset)=0\in X$ and for any disjoint sequence $\{A_n\}\subs\F$
$$\bs{\nu}\Big(\bigcup_{n=1}^\infty A_n\Big)=\sum_{n=1}^\infty\bs{\nu}(A_n)$$
where the series in the right hand side converges in the norm of $X$. Since $\bigcup_{n=1}^\infty A_n=\bigcup_{n=1}^\infty A_{\s(n)}$ for any permutation $\s\colon\NN\to\NN$ the series converges unconditionally but not necessarily absolutely when $X$ is infinite dimensional. Equivalently $\bs{\nu}$ is Banach-valued measure if and only if it is finitely additive and for any disjoint sequence $\{A_n\}_{n=1}^\infty\subs\F$
$$\lim_{n\ra+\infty}\Big\|\bs{\nu}\Big(\bigcup_{i=n}^\infty A_i\Big)\Big\|_X=0.$$
 The \emph{total variation of a Banach valued measure} $\bs{\nu}\colon\F\to X$ is the finitely additive set function 
$$|\bs{\nu}|(A)=\sup_{\mathcal{P}_A}\sum_{E\in\mathcal{P}_A}\|\bs{\nu}(E)\|_X\geq\|\bs\nu(A)\|_X,\quad A\in\F.$$
Here the supremum runs over all finite partitions $\mathcal{P}_A\subs\F$ of $A\in\F$. We say that $\bs{\nu}$ has \emph{bounded variation} if $\|\bs{\nu}\|_{TV}:=|\bs{\nu}|(\TT)<+\infty$, in which case $|\bs{\nu}|$ is a non-negative measure and we will denote by $\MMM(\TT;X)$ the space of all $X$-valued measures on $(\TT,\F)$ with bounded variation. 

A Banach-valued measure $\bs{\nu}\colon\F\to X$ on the measure space $\TT=(\TT,\F,m)$ is called $m$-\emph{absolutely continuous}, which we denote by $\bs\nu\ll m$, if for every $\ee>0$ there exists $\delta>0$ such that~\eqref{ACBVM} holds. The Banach-valued measure $\bs\nu$ is $m$-absolutely continuous if and only if its total variation $|\bs{\nu}|$ is $m$-absolutely continuous. Indeed, since $\|\bs\nu(\cdot)\|_X\leq|\bs\nu|(\cdot)$ it is obvious that if $|\bs\nu|$ is $m$-absolutely continuous then so is $\bs\nu$. For the converse, given $\ee>0$ we can choose a finite partition $\mathcal{P}_\TT=\{E_1,\ldots,E_k\}$, $i\in\NN$ of $\TT$ such that $|\bs\nu|(\TT)\leq\sum_{E\in\mathcal{P}_\TT}\|\bs\nu(E)\|_X+\ee$. Then for any $A\in\F$
\begin{align}
|\bs\nu|(A)&=|\bs\nu|(\TT)-|\bs\nu|(\TT\sm A)<\sum_{E\in\mathcal{P}_\TT}\|\bs\nu(E)\|_X+\ee-\sum_{E\in\mathcal{P}_\TT}\big\|\bs\nu\big((\TT\sm A)\cap E\big)\big\|_X\nonumber\\
&\leq\sum_{E\in\mathcal{P}_\TT}\big\|\bs\nu(E)-\bs\nu\big((\TT\sm A)\cap E\big)\big\|_X+\ee=\sum_{E\in\mathcal{P}_\TT}\|\bs\nu(A\cap E)\|_X+\ee\label{TVFromGlobToSet}
\end{align}
But since $\bs\nu$ is absolutely continuous there exists $\delta=\delta(\ee,\mathcal{P}_\TT,k)>0$ such that $m(A)<\delta$ implies $\|\bs\nu(A)\|_X\leq\frac{\ee}{k}$. So if $m(A)<\delta$ then also $m(A\cap E)<\delta$ for all $E\in\mathcal{P}_\TT$ and therefore $|\bs\nu|(A)\leq 2\ee$ by~\eqref{TVFromGlobToSet} which shows that $|\bs\nu|\ll m$.

The Banach-valued measure $\bs{\nu}\colon\F\to X$ is called $m$-Lipschitz continuous with respect to $m$ if 
$$\Lip_m(\bs\nu):=\sup_{m(A)\neq 0}\frac{\|\bs\nu(A)\|_X}{m(A)}<+\infty.$$ A Banach-valued measure $\bs\nu\in\MMM(\TT;X)$ is Lipshitz continuous if and only if the measure $|\bs\nu|\in\MMM_+(\TT)$ is $m$-Lipschitz continuous i.e.~iff $\Lip_m(|\bs\nu|)=\sup_{m(A)\neq 0}\frac{|\bs\nu|(A)}{m(A)}<+\infty$ and $\Lip_m(\bs\nu)=\Lip_m(|\bs\nu|)$. Indeed, since $\|\bs\nu(A)\|_X\leq|\nu|(A)$ for all $A\in\F$ it is obvious if $|\bs\nu|$ is $m$-Lipschitz continuous then $\bs\nu$ is Lipschitz continuous with $\Lip_m(\bs\nu)\leq\Lip_m(|\bs\nu|)$. Conversely if $\bs\nu$ is $m$-Lipschitz then for every $A\in\F$ 
$$|\bs{\nu}|(A)=\sup_{\mathcal{P}_A}\sum_{E\in\mathcal{P}_A}\|\bs{\nu}(E)\|_{X^*}\leq\Lip_m(\bs\nu)\sup_{\mathcal{P}_A}\sum_{E\in\mathcal{P}_A}m(E)=\Lip_m(\bs\nu)m(A)$$
and thus $|\bs\nu|$ is $m$-Lipschitz continuous with $\Lip_m(|\bs\nu|)\leq\Lip_m(\bs\nu)$. In particular any $m$-Lipschitz continuous measure $\bs\nu\in\MMM(\TT;X)$ is also $m$-absolutely continuous. The space of $m$-Lipschitz continuous measures $\bs\nu\in\MMM(\TT;X)$ will be denoted by $\MMM_\Lip(m;X)$. The linear space $\MMM_\Lip(m;X)$ becomes a normed space when equipped with the norm $\Lip_m\colon\MMM_\Lip(m;X)\to\RR_+$.

Given a Banach-valued measure $\bs\nu\in\MMM(\TT;X)$ we define $L^p(\bs\nu)=L^p(|\bs\nu|)$, $1\leq p\leq\infty$ and we can define the $\df\bs\nu$-integral for simple maps $\phi=\sum_{i=1}^na_i\1_{A_i}$, $a_i\in\RR$, $A_i\in\F$ by 
$$\int\phi(t)\df\bs\nu(t)=\sum_{i=1}^na_i\bs\nu(A_i).$$
Then for all simple maps $\phi$ as above
$$\Big\|\int\phi(t)\df\bs\nu(t)\Big\|_X\leq\sum_{i=1}^n|a_i|\|\bs\nu(A_i)\|_X\leq\sum_{i=1}^n|a_i||\bs\nu|(A_i)=\int|\phi(t)|\df|\bs\nu|(t)$$ 
and the $\df\bs\nu$-integral can be extended to a linear vector-valued integral $\df\bs\nu\colon L^1(\bs\nu)\to X$ by defining $\int f\df\bs\nu=\lim_{n\ra+\infty}\int\phi_n(t)\df\bs\nu(t)$ for any sequence of simple functions $\{\phi_n\}$ such that $\lim_{n\ra+\infty}\int|\phi_n(t)-f(t)|\df|\bs\nu|(t)=0$. This does not depend on the choice of the sequence $\{\phi_n\}$ of simple functions and satisfies $\|\int f\df\bs\nu\|_X\leq\int|f|\df|\bs\nu|$. 

In the case that the Banach-valued measure $\bs\nu$ takes values in a dual space, i.e.~$\bs\nu\in\MMM(\TT;X^*)$ then one can define a real valued integral $\int\ls\cdot,\df\bs\nu\rs$ on $L^1(\bs\nu;X):=L^1(|\bs\nu|;X)$. This can be done by defining for any simple function $\phi=\sum_{i=1}^nf_i\1_{A_i}$, $f_i\in X$, $A_i\in\F$ 
$$\int\ls\phi(t),\df\bs\nu(t)\rs=\sum_{i=1}^n\bs\nu(A_i)(f_i).$$
Then for any such simple map $\phi$ in canonical form so that $\|\phi\|_{X^*}=\sum_{i=1}^n\|f_i\|_X\1_{A_i}$
$$\Big|\int\ls\phi(t),\df\bs\nu(t)\rs\Big|\leq\sum_{i=1}^n\|\bs\nu(A_i)\|_{X^*}\|f_i\|_X\leq\sum_{i=1}^n|\bs\nu|(A_i)(f_i)=\int\|\phi(t)\|_X\df|\bs\nu|(t)$$
and we can extend to all maps by 
\begin{equation}\label{WeakStarIntegral} 
\int\ls f,\df\bs\nu(t)\rs=\lim_{n\ra+\infty}\int\ls\phi_n(t),\df\bs\nu(t)\rs
\end{equation}
where $\{\phi_n\}$ is any sequence of simple functions satisfying 
\begin{equation}\label{WeakStarIntegralAdmiss} \lim_{n\ra+\infty}\int\|\phi_n(t)-f(t)\|_X\df|\bs\nu|(t)=0.
\end{equation} Such a sequence $\{\phi_n\}$ exists since $\|f_\cdot\|_X\in L^1(|\bs\nu|)$ because $f$ is assumed in the Bochner space $L^1(|\bs\nu|;X)$ and the definition~\eqref{WeakStarIntegral} does not depend on the choice of sequence $\{\phi_n\}$ of simple functions satisfying~\eqref{WeakStarIntegralAdmiss}. Then 
\begin{equation}
\Big|\int\ls f(t),\df\bs\nu(t)\rs\Big|\leq\int\|f_t\|_X\df|\bs\nu|(t),\quad\forall f\in L^1(\bs\nu;X).
\end{equation}

Using this real valued integral we can see $L^1(m;X)^*$ is isometric to $\MMM_\Lip(m;X^*)$.
\begin{prop}\label{L1DualToLipMeas} The linear operator $\bs{V}\colon L^1(\TT;X)^*\to\MMM_\Lip(\TT;X^*)$ defined by assigning to each functional $J\in L^1(\TT;X)^*$ the Banach-valued measure $\bs V_J\in\MMM(\TT;X^*)$ given by $\ls f,\bs V_J(A)\rs\equiv\bs V_J(A)(f)=J(f\1_A)$ is a surjective isometry.
\end{prop}\textbf{Proof} Indeed, the set function $\bs{\nu}_J$ is obviously finitely additive and for any $f\in X$ and $\A\in\F$ satisfies $|\bs{V}_J(A)(f)|\leq\|J\|_{L^1(m;X)^*}\|f\|_Xm(A)$. Taking the supremum over all $f$ in the unit ball of $L^1(\TT;X)$ we obtain that $\big\|\bs{V}_J(A)\big\|_{X^*}\leq\|J\|_{L^1(\TT;X)^*}m(A)$ for all $A\in\F$. In particular $\bs{V}_J$ is a $m$-Lipschitz Banach-valued measure with $\Lip_m(\bs V_J)\leq\|J\|_{L^1(m;X)^*}$. Thus the map $\bs V$ is a well defined contraction and is obviously injective.
	
	Conversely, for $\bs\nu\in\MMM_\Lip(m;X^*)$ the inequality $|\bs\nu|(\cdot)\leq\Lip_m(\bs\nu) m(\cdot)$ implies that $L^1(m;X)\leq L^1(|\bs\nu|;X)$ and  
	$\int\ls f(t),\df\bs\nu(t)\rs\leq\Lip_m(\bs\nu)\|f\|_{L^1(m;X)}$ for all $f\in L^1(m;X)$. Thus the formula 
	\begin{equation}\label{Jmap}
	\mathbbm{J}_{\bs\nu}(f)=\int\ls f,\df\bs\nu\rs,\quad f\in L^1(m;X)\leq L^1(|\bs\nu|;X)
	\end{equation}
	defines a linear functional $\mathbbm{J}_{\bs\nu}\in L^1(m;X)^*$ with $\|\mathbbm{J}_{\bs\nu}\|_{L^1(m;X)}\leq\Lip_m(\bs\nu)$. Consequently the assignment $\MMM_\Lip(m;X^*)\ni\bs\nu\mapsto\mathbbm{J}_{\bs\nu}\in L^1(m;X)^*$ defines a contraction that is obviously injective.
	
	But the maps $\bs V$ and $\mathbbm{J}$ are inverse to each other since on one hand $\bs V_{\mathbbm{J}_{\bs\nu}}(A)(f)=\int\ls f\1_A(t),\df\bs\nu\rs=\ls f,\bs\nu(A)\rs$ for all $\bs\nu\in\MMM_\Lip(m;X^*)$, $A\in\F$ and $f\in X$ so that $\bs V_{\mathbbm{J}_{\bs\nu}}=\bs\nu$ for all $\bs\nu\in\MMM(m;X^*)$ and thus $\bs V\circ\mathbbm{J}={\rm{id}}_{\MMM(m;X^*)}$. On the other hand for all $J\in L^1(m;X)^*$ and all simple maps $\phi=\sum_{i=1}^nf_i\1_{A_i}\in L^1(m;X)$
	$$\mathbbm{J}_{\bs V_J}(\phi)=\int\ls \phi(t),\df\bs V_J(t)\rs=\sum_{i=1}^n\ls f_i,\bs V_J(A_i)\rs=\sum_{i=1}^nJ(f_i\1_{A_i})=J(\phi)$$
	and since this holds for all simple maps $\phi\in L^1(m;X)$ and the functionals $\mathbbm{J}_{\bs V_J}$ and $J$ are both continuous it follows that $\mathbbm{J}_{\bs V_J}=J$ for all $J\in L^1(m;X)^*$. Therefore $\mathbbm{J}\circ\bs V={\rm{id}}_{L^1(m;X)^*}$. It follows that the maps $\bs V$ and $\mathbbm{J}$ are both surjective isometries and the proof is complete.$\hfill\Box$

\subsubsection{Weak-Star Radon-Nikodym derivatives}
\begin{theorema}\label{WeakStarRNDer} Let $X$ be a Banach space and let $(\TT,\F,m)$ be a complete and finite positive measure space. There exists a linear isometric inclusion $\hat{\mu}\colon\MMM_\Lip(m;X^*)\to\LL_{w^*}^\infty(\TT;X^*)$ such that for each $\bs\nu\in\MMM_\Lip(m;X^*)$
		\begin{itemize}
			\item[(1)] The map $\hat{\mu}^{\bs\nu}\in\widehat{\LL}_{w^*}^\infty(\TT;X^*)$ for all $\bs\nu\in\MMM_\Lip(m;X^*)$
			\item[(2)] For all $f\in X$ and $A\in\F$ 
		$$\bs{\nu}(A)(f)=\int_A\ls f,\hat{\mu}^{\bs\nu}_t\rs\df m(t).$$
		\item[(3)] For all $A\in\F$ 
		$$|\bs{\nu}|(A)=\int_A\|\hat{\mu}^{\bs\nu}_t\|_{X^*}\df m(t).$$
	\end{itemize}
Furthermore, any map $\hat{\mu}\colon\MMM_\Lip(m;X^*)\to\LL_{w^*}^\infty(\TT;X^*)$ satisfying properties (1) to (3) above passes to a surjective isometry $\wt{\mu}\colon\MMM_\Lip(m;X^*)\to\bar{L}_{w^*}^\infty(\TT;X^*):=\,^{L_{w^*}^\infty(\TT;X^*)}/_{\ker S}$ when composed with the natural quotient map $[\cdot]_S\colon L_{w^*}^\infty(\TT;X^*)\to\bar{L}_{w^*}^\infty(\TT;X^*)$ of the relation of $w^*$-m-a.s.~equality.
\end{theorema}\textbf{Proof} We follow the proof based on the existence of linear liftings found in~\cite[Theorem 1.5.2]{PilarCembranos1997a}. A lifting on the space $(\TT,\F,m)$ is linear right inverse $\ell\colon L^\infty(m)\to\LL^\infty(m)$ to the quotient map $[\cdot]_m\colon\LL^\infty(m)\to L^\infty(m)$ that is unital i.e.~$\ell(1)=\1_\TT$ and monotone i.e. if $f\leq g$ in $L^\infty(m)$ then $\ell(f)(t)\leq\ell(g)(t)$ for all $t\in\TT$. If also $\ell(f\cdot g)=\ell(f)\cdot\ell(g)$ then $\ell$ is called a strong lifting. For the existence of a strong lifting on $L^\infty(m)$ on complete positive measure spaces we refer to~\cite{Dinculeanu1967a,TulceaAC1969a}. 

Using the existence of liftings it is easy to define the required isometric inclusion $\hat{\mu}\colon\MMM_\Lip(m;X^*)\to\LL_{w^*}^\infty(\TT;X^*)$. Indeed for each $\bs\nu\in\MMM_\Lip(m;X^*)$ for all $f\in X$ the signed measure $\nu_f:=\ls f,\bs\nu(\cdot)\rs$ is $m$-Lipschitz continuous since $|\nu_f(A)|\leq\Lip_m(\bs\nu)\|f\|_Xm(A)$
for all $A\in\F$. Thus $\nu_f$ has a Radon-Nikodym derivative $\frac{\df\nu_f}{\df m}\in L^\infty(m)$. Then by fixing a lifting $\ell\colon L^\infty(m)\to\LL^\infty(m)$ we define the map $\hat{\mu}^{\bs\nu}\equiv\hat{\mu}^{\bs\nu;\ell}\colon\TT\to X^*$ by $$\hat{\mu}^{\bs\nu}_t(f)=\ell\Big(\frac{\df\nu_f}{\df m}\Big)(t).$$ Let us check that indeed $\hat{\mu}_t^{\bs\nu}\in X^*$ for all $t\in\TT$. Obviously $\nu_{af+bg}=a\nu_f+b\nu_g$ for all $f,g\in X$, $a,b\in\RR$ which implies that $\frac{\df\nu_{af+bg}}{\df m}=a\frac{\df\nu_f}{\df m}+b\frac{\df\nu_g}{\df m}$ in $L^\infty(m)$ Therefore by the linearity of liftings we obtain the linearity of $\hat{\mu}_t^{\bs\nu}$ for all $t\in\TT$. To see that $\hat{\mu}_t^{\bs\nu}$ is bounded we note that since $|\bs\nu|$ is also Lipschitz $\|\frac{\df|\bs\nu|}{\df m}\|_{L^\infty(m)}\leq\Lip_m(\bs\nu)$ and thus for all $f\in X$ and $A\in\F$
$$\int_A\frac{\df\nu_f}{\df m}(t)\df m(t)=\nu_f(A)\leq\|f\|_X\|\bs\nu(A)\|_{X^*}\leq\|f\|_X|\bs\nu|(A)=\|f\|_X\int_A\frac{\df|\bs\nu|}{\df m}(t)dm(t).$$
Since this holds for all $A\in\F$ is follows that $\frac{\df\nu_f}{\df m}\leq\|f\|_X \frac{\df|\bs\nu|}{\df m}\leq\|f\|_X\cdot\Lip_m(\bs\nu)$ $m$-a.e.~in $\TT$ for all $f\in X$. Therefore by the monotonicity of liftings, for all $t\in\TT$
$$\hat{\mu}_t^{\bs\nu}(f)=\ell\Big(\frac{\df\nu_f}{\df m}\Big)(t)\leq\|f\|_X\cdot\ell\Big(\frac{\df|\bs\nu|}{\df m}\Big)(t)\leq\|f\|_X\cdot\Lip_m(\bs\nu)$$
and thus $\hat{\mu}_t^{\bs\nu}\in X^*$ with $\|\hat{\mu}_t^{\bs\nu}\|_{X^*}\leq\ell(\frac{\df|\bs\nu|}{\df m})(t)\leq\Lip_m(\bs\nu)$. Since $\hat{\mu}^{\bs\nu}$ is by definition $w^*$-measurable it follows that $\hat{\mu}^{\bs\nu}\in\LL^\infty_{w^*}(\TT;X^*)$ and $\|\hat{\mu}^{\bs\nu}\|_{\LL_{w^*}^\infty(\TT;X^*)}\leq\Lip_m(\bs\nu)$.

The map $\MMM_\Lip(m;X^*)\ni\bs\nu\to\hat{\mu}^{\bs\nu}\in\LL_{w^*}^\infty(\TT;X^*)$ is obviously linear and property (2) holds by definition. We will check that properties (1) and (3) also hold and for this it suffices to show that $\|\hat{\mu}^{\bs\nu}_t\|_{X^*}=\frac{\df|\bs\nu|}{\df m}(t)$ for almost all $t\in\TT$. Let $g$ denote the supremum of the family of all functions of the form 
$$\sum_{i=1}^n\1_{A_i}(t)\Big|\ell\Big(\frac{\df\nu_{f_i}}{\df m}\Big)(t)\Big|$$
where $\{A_i\}_{i=1}^n$ is a partition of $\TT$, $n\in\NN$ and $\{f_i\}_{i=1}^n$ is a finite sequence in the unit ball of $X$. Since $\ell(\frac{\df\nu_{f}}{\df m})(t)=\hat{\mu}^{\bs\nu}_t(f)\leq\|\hat{\mu}_t^{\bs\nu}\|_{X^*}$  for any $f\in X$ with $\|f\|_X\leq 1$ any map in this family is bounded above by $\|\hat{\mu}_\cdot^{\bs\nu}\|_{X^*}$ and thus by~\cite[Corollary IV.11.7]{Dunford1958a} the map $g$ is in $L^\infty(m)$ and $$0\leq g(t)\leq\|\hat{\mu}_t^{\bs\nu}\|_{X^*}\leq\frac{\df|\bs\nu|}{\df m}(t)\leq\Lip_m(\bs\nu).$$
By the definition of $|\bs\nu|$, given $\ee>0$ there exists a partition $\{A_i\}_{i=1}^n$ of $\TT$ and unit vectors $f_i\in X$ such that $\sum_{i=1}^n\bs\nu(A_i)(f_i)\geq|\bs\nu|(\TT)-\ee$ and thus
\begin{align*}
\int_\TT\frac{\df|\bs\nu|}{\df m}\df m(t)-\ee&=|\bs\nu|(\TT)-\ee\leq\sum_{i=1}^n\bs\nu(A_i)(f_i)=\sum_{i=1}^n\int_{A_i}\frac{\df\nu_{f_i}}{\df m}(t)\df m(t)\\
&=\int\sum_{i=1}^n\1_{A_i}(t)\ell\Big(\frac{\df\nu_{f_i}}{\df m}\Big)(t)\df m(t)\leq\int_\TT g(t)\df m(t).
\end{align*}
Since $g(t)\leq\frac{\df|\bs\nu|}{\df m}(t)$ for $m$-almost all $t\in\TT$ this implies that $g=\|\hat{\mu}^{\bs\nu}_\cdot\|_{X^*}=\frac{\df|\bs\nu|}{\df m}$ $m$-a.e.~in $\TT$ as required. Using property (3) it is now easy to see $\hat{\mu}$ is norm-preserving. Indeed,
\begin{align*}
\|\mu^{\bs\nu}\|_{\LL_{w^*}^\infty(\TT;X^*)}&=\big\|\|\mu^{\bs\nu}_\cdot\|_{X^*}\big\|_{L^\infty(m)}=\sup_{m(A)\neq 0}\fr{m(A)}\int_A\|\mu^{\bs\nu}_t\|_{X^*}\df m(t)\\
&=\sup_{m(A)\neq 0}\frac{|\bs\nu|(A)}{m(A)}=\Lip_m(\bs\nu).
\end{align*} for any $\bs\nu\in\MMM_\Lip(m;X^*)$ and the proof of the first claim is complete.

We prove next the second claim. The map $\wt{\mu}=[\cdot]_S\circ\hat{\mu}\colon\MMM_\Lip(m;X^*)\to\bar{L}_{w^*}^\infty(\TT;X^*)$ is obviously a contraction. We will show that it is also surjective. Indeed, we can define $\wt{\bs\nu}\colon\bar{L}_{w^*}^\infty(\TT;X^*)\to\MMM_\Lip(m;X^*)$ by 
\begin{equation}\label{numap}
\wt{\bs\nu}_{\mu}(A)(f)=\int_A\ls f,\mu_t\rs\df m(t),\quad \mu\in L_{w^*}^\infty(\TT;X^*),\;A\in\F,\;f\in X.\
\end{equation}
This is well-defined according to the definition of the relation $\backsim_S$ of $w^*$-m-a.s.~equality. Then for each $A\in\F$, $f\in X$ and $g\in\A_{\mu_0}$ where $\mu_0\backsim_S\mu$ is a representative of the class $\mu$
$$|\wt{\bs\nu}_\mu(A)(f)|\leq\|f\|_X\int_Ag(t)\df m(t)\leq\|f\|_X\|g\|_{L^\infty(m)}m(A)$$
and taking first the infimum over all $g\in\A_{\mu_0}$ and then the supremum over all $f\in X$ with $\|f\|_X\leq 1$ we obtain that $\|\wt{\bs\nu}(A)\|\leq\|\mu_0\|_{L^\infty_{w^*}(\TT;X^*)}m(A)$ for all $A\in\F$ and thus
$$\Lip_m(\wt{\bs\nu}_\mu)\leq\|\mu_0\|_{L^\infty_{w^*}(\TT;X^*)}.$$
Taking then the infimum over all $\mu_0\backsim_S\mu$ we obtain that the map $\wt{\bs\nu}$ is a contraction. Now by applying the map $\hat{\mu}\colon\MMM_\Lip(m;X^*)\to\LL^\infty_{w^*}(\TT;X^*)$ to $\wt{\bs\nu}_\mu$ we obtain $\hat{\mu}^{\wt{\bs\nu}_\mu}\in\LL_{w^*}^\infty(\TT;X^*)$ satisfying properties (1) to (3) of Theorem~\ref{WeakStarRNDer}. In particular by (2) for all $f\in X$ 
$$\int_A\ls f,\hat{\mu}^{\wt{\bs\nu}_\mu}_t\rs\df m(t)=\int_A\ls f,\mu_t\rs\df m(t),\quad\forall A\in\F$$
which implies that $\hat{\mu}^{\wt{\bs\nu}_\mu}\backsim_S\mu$ so that $\hat{\mu}^{\wt{\bs\nu}_\mu}=\mu$ in $\bar{L}_{w^*}^\infty(\TT;X^*)$. This proves that $\wt{\mu}\circ\wt{\bs\nu}=\mathbbm{id}_{\bar{L}_{w^*}^\infty(\TT;X^*)}$ and in particular the map $\wt{\mu}=[\cdot]_S\circ\hat{\mu}$ is surjective. Since both maps $\wt{\mu}$ and $\wt{\bs\nu}$ are contractions if we show that also $\wt{\bs\nu}\circ\wt{\mu}=\mathbbm{id}_{\MMM_\Lip(m;X^*)}$ is will follow that they are both surjective isomorphisms inverse to each other. But this is easy since for any $\bs\nu\in\MMM_\Lip(m;X^*)$ and any $A\in\F$, $f\in X$ 
$$\wt{\bs\nu}_{\wt{\mu}^{\bs\nu}}(A)(f)=\int_A\ls f,\wt{\mu}^{\bs\nu}_t\rs\df m(t)=\int_A\ls f,\hat{\mu}^{\bs\nu}_t\rs\df m(t)=\bs\nu(A)(f)$$
where the first equality is just the definition of the map $\wt{\bs\nu}$, the second is due to the fact that $\hat{\mu}^{\bs\nu}\backsim_S\wt{\mu}^{\bs\nu}$ and the last equality is by property (2) of the map $\hat{\mu}$.$\hfill\Box$\\

A few remarks are in order. First, since $\bar{L}_{w^*}^\infty(\TT;X^*)=L_{w^*}^\infty(\TT;X^*)$ when $X$ is separable, in this case the map  $\hat{\mu}\colon\MMM_\Lip(m;X^*)\to L_{w^*}^\infty(\TT;X^*)$ is an isometric isomorphism. Furthermore since the map $\|\mu_\cdot\|_{X^*}$ is $\F$-measurable for all $\mu\in L_{w^*}^\infty(\TT;X^*)$ property (1) of the map $\hat{\mu}$ in Theorem~\ref{WeakStarRNDer} is redundant and the induced map $\hat{\mu}\colon\MMM_\Lip(m;X^*)\to L_{w^*}^\infty(\TT;X^*)$ by composing with the quotient map of the relation of $m$-a.s.~equality is uniquely determined by property (2). 

In the case that $X$ is non-separable for any $\mu\in\bar{L}_{w^*}^\infty(\TT;X^*)$ any representative of its $w^*$-m-a.s~equality class satisfies property (2) of Theorem~\ref{WeakStarRNDer}, while there exists a representative $\hat{\mu}\backsim_S\mu$ in its class of $w^*$-m-a.s.~equality such that properties (1) and (3) hold namely the map $\hat{\mu}:=\hat{\mu}^{\wt{\bs\nu}_\mu}$ where $\wt{\bs{\nu}}$ is the map defined in~\eqref{numap}. For this representative $\|\mu\|_{L_{w^*}^\infty(\TT;X^*)}=\big\|\|\hat{\mu}_\cdot\|_{X^*}\big\|_{L^\infty(m)}$. Thus even when $X$ is non-separable one can always choose for each $\mu\in\bar{L}_{w^*}^\infty(\TT;X^*)$ a representative $\hat{\mu}$ from its $w^*$-m-a.s~equality class such that $\hat{\mu}\in\hat{L}_{w^*}^\infty(\TT;X^*)$.

\begin{cor} The map $S\colon\bar{L}_{w^*}^\infty(\TT;X^*)\to L^1(\TT;X)^*$ defined by $S(\mu)(f)=\lls f,\mu\rrs$ is an isometric isomorphism.
\end{cor}\textbf{Proof} It suffices to check that $S=\mathbbm{J}\circ\wt{\bs\nu}$ where $\wt{\bs\nu}\colon\bar{L}_{w^*}^\infty(\TT;X^*)\to\MMM_\Lip(m;X^*)$ and $\mathbbm{J}\colon\MMM_\Lip(m;X^*)\to L^1(\TT;X)^*$ are the isometric isomorphisms defined in~\eqref{numap} and~\eqref{Jmap} respectively. To check this let $\mu\in \bar{L}_{w^*}^\infty(\TT;X^*)$ and let $\phi=\sum_{i=1}^nf_i\1_{A_i}\in L^1(\TT;X)$, $f_i\in X$, $A_i\in\F$ be a simple function. Then
\begin{align*} \mathbbm{J}_{\wt{\bs\nu}_\mu}(\phi)&=\int\ls\phi,\df\wt{\bs\nu}_\mu\rs=\sum_{i=1}^n\wt{\bs\nu}_\mu(A_i)(f_i)=\sum_{i=1}^n\int_{A_i}\ls f_i,\mu_t\rs\df m(t)\\
&=\int_\TT \ls \phi_t,\mu_t\rs\df m(t)=S(\mu)(\phi).
\end{align*}
Since the set of simple functions is dense in $L^1(\TT;X)$ and $\mathbbm{J}_{\wt{\bs\nu}_\mu}$, $S(\mu)$ are continuous linear functionals on $L^1(\TT;X)$ it follows that $I=S(\mu)$ and thus $S=\mathbbm{J}\circ\wt{\bs\nu}$.$\hfill\Box$

\begin{prop}\label{ToApplyA2} Suppose that the Banach space $X$ is separable. If there exists a countable collection $\mathcal{A}\subs\F$ such that \begin{equation}\label{SeparableCollection}\forall\;E\in\F,\;\forall\;\ee>0,\;\exists\;A\in\A:\;m(E\triangle A)<\ee\end{equation} then the space $L^1(\TT;X)$ is separable.
\end{prop}\textbf{Proof} Let $f\in L^1(\TT;X)$, $\ee>0$. Fix a dense countable subset $D\subs X$ of $X$ and let $\A$ be a countable collection satisfying~\eqref{SeparableCollection}. Then the set $\mathcal{D}\subs L^1(\TT;X)$ consisting of all functions of the form 
$$\sum_{j=1}^nq_j\1_{A_j},\quad q_j\in D,\;A_j\in\A,\;n\in\NN$$ is obviously countable. We will show that it is also dense in $X$. Indeed, since $f\in L^1(\TT;X)$ there exists a simple function $\phi=\sum_{k=1}^nf_k\1_{E_k}\in L^1(\TT;X)$ such that $\|\phi-f\|_{L^1(\TT;X)}<\ee/2$. We set $M:=\max_{1\leq k\leq n}\|f_k\|_X$. Then for each $k=1,\dots,n$ there exists $A_k\in\A$ such that $m(E_k\triangle A_k)<\ee/4nM$ and since $D$ is dense in $X$, for each $k=1,\dots,n$ there exists $g_k\in D$ such that $\|g_k-f_k\|_X<\ee/4nm(\TT)$. Then $\psi:=\sum_{k=1}^ng_k\1_{A_k}\in\mathcal{D}$ and $\|\psi-g\|_{L^1(\TT;X)}\leq\|\psi-\phi\|_{L^1(\TT;X)}+\frac{\ee}{2}$. But 
\begin{align*}
\|\psi-\phi\|_{L^1(\TT;X)}&\leq\Big\|\psi-\sum_{k=1}^nf_k\1_{A_k}\Big\|_{L^1(\TT;X)}
+\Big\|\sum_{k=1}^nf_k\1_{A_k}-\phi\Big\|_{L^1(\TT;X)}\\
&\leq\sum_{k=1}^n\Big(\int_\TT\|g_k-f_k\|_X\1_{A_k}(t)\df m(t)+\int_\TT\|f_k\|_X\1_{A_k\triangle E_k}(t)\df m(t)\Big)\\
&\leq\frac{\ee}{4nm(\TT)}\sum_{k=1}^n\int_\TT\1_{A_k}(t)\df m(t)+M\sum_{k=1}^n m(E_k\triangle A_k)<\frac{\ee}{2},
\end{align*}
and so we have found an element of $\mathcal{D}$ that is $\ee$-close to $f\in L^1(\TT;X)$.$\hfill\Box$\\

As a consequence of Proposition~\ref{ToApplyA2}, in the case that $\TT$ is the interval $[0,T]$ for some $T>0$ equipped with the Lebesgue measure on the Lebesgue $\s$-algebra, the space $L_{w^*}^\infty(\TT;X^*)\cong L^1(\TT;X)^*$ is submetrizable and since it is completely regular as a Hausdorff topological vector space, all the results of Section~\ref{Appendix1} are applicable. In particular any probability measure on $(L_{w^*}(\TT;X^*),w^*)$ is Radon and the portmanteau and Prokhorov theorems which are well-known known in the category of polish spaces are also valid on $L_{w^*}^\infty(\TT;X^*)$ for separable $X$. 

Let us finally note that the map $\wt{\bs\nu}\colon L_{w^*}^\infty(\TT;X^*)\to\MMM_\Lip(m;X^*)$ defined in~\eqref{numap} can be viewed as the indefinite $w^*$-integral and we can equivalently use the notation 
$$\wt{\bs\nu}(E)(f)=\Big\ls f,w^*\mbox{-}\int_E\mu_t\df m(t)\Big\rs=\int_E\ls f,\mu_t\rs\df m(t),\quad f\in X.$$
\begin{prop}\label{MeasInd} For any $T\in\bbar{BC_{w^*}(X^,Y^*)}^{sq}$, $\mu\in L_{w^*}^\infty(\TT;X^*)$ and measurable $E\subs\TT$
	\begin{equation}\label{GelfandCommute}
	T\Big(w^*\mbox{-}\int_E\mu_t\df m(t)\Big)=w^*\mbox{-}\int_E T(\mu_t)\df m(t).
	\end{equation}
\end{prop}\textbf{Proof} It suffices to show that the space $\mathcal{C}$ of all bounded linear operators $T\colon X^*\to Y^*$ satisfying~\eqref{GelfandCommute} contains $BC_{w^*}(X^*,Y^*)$ and is sequentially closed with respect to pointwise $w^*$-convergence. So let $T\colon X^*\to Y^*$ be $w^*$-continuous. Then $T=S^*$ for some bounded operator $S\colon X\to Y$ and thus for any $g\in Y$
$$\Big\ls g,T\Big(w^*\mbox{-}\int_E\mu_t\df m(t)\Big)\Big\rs=\Big\ls Sg,w^*\mbox{-}\int_E\mu_t\df m(t)\Big\rs=\int_E\ls Sg,\mu_t\rs\df m(t)=\int_E\ls g,T\mu_t\rs\df m(t)$$
which shows that $BC_{w^*}(X^*,Y^*)\subs\mathcal{C}$. Let now $\{T_n\}\subs\mathcal{C}$ be a sequence of operators $w^*$-converging pointwise to an operator $T\colon X^*\to Y^*$. Then 
$$\Big\ls g,T_n\Big(w^*\mbox{-}\int_E\mu_t\df m(t)\Big)\Big\rs=\int_E\ls g,T_n\mu_t\rs\df m(t)$$
for all $n\in\NN$ and the left-hand side term converges to $\ls g,T(w^*\mbox{-}\int_E\mu_t\df m(t)\rs$ as $n\ra+\infty$. Furthermore since $\{T_n\}$ pointwise $w^*$-converges to $T$ it is norm bounded, i.e.~$C:=\sup_{n\in\NN}\|T_n\|<+\infty$. Indeed, for each fixed $\mu\in X^*$ 
$$\sup_{n\in\NN}|\ls g,T_n\mu\rs|<+\infty,\quad\forall g\in Y$$ and thus by the uniform boundedness principle we obtain that $\sup_{n\in\NN}\|T_n\mu\|_{Y^*}<+\infty$ for all $\mu\in X^*$, which by uniform boundedness principle again yields that $\sup_{n\in\NN}\|T_n\|<+\infty$. Therefore the sequence of the maps $\ls g_\cdot,T_n\mu_\cdot\rs$ is dominated by the $L^1$ function $t\mapsto C\|\mu\|_{L_{w^*}^\infty(\TT;X^*)}\|g_\cdot\|_Y$ and by the dominated convergence theorem we conclude that the right hand side term converges to $\int_E\ls g,T\mu_t\rs\df m(t)$ as $n\ra+\infty$ which completes the proof.$\hfill\Box$
\subsection{$L_{w^*}^\infty$-valued random variables}
The next proposition ensures us that the empirical processes under consideration in this article are all well defined random variables with values in the measurable space $(L^\infty_{w^*}(\TT;X^*),\B_{L_{w^*}})$ equipped with the Borel $\s$-algebra $\B_{L_{w^*}^\infty}:=\B_{(L_{w^*}^\infty(0,T;X^*),w^*)}$ of the $w^*$-topology, for an appropriate separable Banach space $X$, with $X^*$ being used to encode the empirical density of the ZRP at each time $t\in[0,T]$. These are obtained via the continuous natural inclusion of the Skorohod space $D(0,T;X^*)$ in $L_{w^*}^\infty(0,T;X^*)$ described below. For a nice survey on Skorohod spaces the reader is referred to~\cite{Ethier1986a}. 
\begin{prop}\label{SkorohodToLEmbed} For any Banach space $X$, 
	$$D(0,T;X^*)\subs L^\infty_{w^*}(0,T;X^*)$$ and the natural inclusion map is continuous with respect to the Skorohod topology and the $w^*$-topology on $L^\infty_{w^*}(0,T;X^*)\cong L^1(0,T;X)^*$. 
\end{prop}\textbf{Proof} Any cadlag path in $D(0,T;X^*)$ is strongly measurable and uniformly bounded and therefore $D(0,T;X^*)\subs L^\infty(0,T;X^*)\subs L_{w^*}^\infty(0,T;X^*)$. To show that the natural inclusion is continuous let $\{\mu_k\}\subs D(0,T;\MMM_+(\T^d))$ be a sequence converging to $\mu$ in the Skorohod metric. Then the set 
$\bigcup_{k\in\NN}\mu_k([0,T])$ is relatively compact in $X^*$ and thus 
\begin{equation}\label{SkorConvSeqUnifBound}C:=\sup_{k\in\NN}\sup_{0\leq t\leq T}\|\mu_{k,t}\|_{X^*}\mx\|\mu_t\|_{X^*}<+\infty
\end{equation} and there exists a sequence $\{\lambda_k\}_{k\in\NN}$ of Lipschitz increasing reparametrizations of $[0,T]$ such that $\g(\lambda_k):=\|\log\lambda_k'\|_{L^\infty([0,T])}\lra 0$ and
$$\lim_{k\ra\infty}\sup_{0\leq t\leq T}\|\mu_{k,t}-\mu_{\lambda_k(t)}\|_{X^*}=0.$$ 
Therefore, if we choose $k_0\in\NN$ large enough so that 
$$k\geq k_0\quad\Lra\quad\sup_{0\leq t\leq T}\|\mu_{k,t}-\mu_{\eta_{\lambda_k(t)}}\|_{X^*}\leq 1,$$ then for all $k\geq k_0$
\begin{align*}\|\mu_{k,t}-\mu_t\|_{X^*}&\leq 1+\|\mu_{\lambda_k(t)}-\mu_t\|_{X^*}\leq 1+\|\mu_{\lambda_k(t)}\|_{X^*}+\|\mu_t\|_{X^*}\\ 
&\leq 1+2\sup_{0\leq t\leq T}\|\mu_t\|_{X^*}\leq 1+2C<+\infty
\end{align*}
for almost all $t\in[0,T]$. Furthermore, since $\mu_k\lra\mu$ in the Skorohod topology we have that $\lim_{k\ra+\infty}\|\mu_{k,t}-\mu_t\|_{X^*}= 0$ for almost all $t\in[0,T]$ and therefore by~\eqref{SkorConvSeqUnifBound} we can apply the dominated convergence theorem to obtain 
$$\lim_{k\ra+\infty}\int f(t)\|\mu_{k,t}-\mu_t\|_{X^*}\df t=0,\quad\forall\;f\in L^1(0,T).$$
In other words $\{\mu_k\}$ converges to $\mu$ in the \emph{normed $w^*$-convergence}, i.e.~the sequence $\{(\|\mu_{k,\cdot}-\mu_\cdot\|_{X^*}\}_{k=1}^\infty\subs L^\infty(0,T)$ converges to $0$ in the $w^*$ topology of $L^\infty(0,T)\cong L^1(0,T)^*$. But then $\mu_k\lra\mu$ in the $w^*$-topology of $L_{w^*}^\infty(0,T;X^*)$ since
$$
\Big|\int_0^T\ls f_t,\mu_{k,t}\rs-\ls f_t,\mu_t\rs\df t\Big|\leq\int_0^T\|\mu_{k,t}-\mu_t\|_{X^*}\|f_t\|_{X}\df t\stackrel{k\ra+\infty}\lra 0$$ for any $f\in L^1(0,T;X)$. $\hfill\Box$

\subsubsection{Subspaces of $L_{w^*}^\infty(0,T;X^*)$}

For any subset $K$ of a Banach space $X$ we will use the notation
$$L^1(\TT;K):=\big\{f\in L^1(\TT;X)\bigm|f_t\in K\mbox{ for almost all }t\in\TT\big\}.$$
We also define 
$$L_{w^*}^\infty(\TT;K):=\big\{\mu\in L_{w^*}^\infty(\TT;X)\bigm|\mu_t\in K\mbox{ for almost all }t\in\TT\big\}.$$

For any positive cone $K$ in a separable Banach space $X$ (i.e.~$\lambda f+g\in K$ for any $f,g\in K$, $\lambda\geq 0$) we denote by 
\begin{equation}\label{PosCone}
K^*:=\{\mu\in X^*\bigm|\mu(f)\geq 0\mbox{ for all }f\in K\}
\end{equation} its dual cone in $X^*$. Obviously $K^*$ is $w^*$-closed subspace of $X^*$.
\begin{prop}\label{PosCon} Let $K$ be a positive closed cone in the separable Banach space $X$. Then $L^1(\TT;K)$ is a closed positive cone in $L^1(\TT;X)$ and
	$$L^\infty_{w^*}(\TT;K^*)=L^1(\TT;K)^*$$
	is a $w^*$-closed positive cone in $L_{w^*}^\infty(\TT;X^*)$. Consequently the space $\PP L_{w^*}^\infty(\TT;K^*)$ is a closed subspace of $\PP L_{w^*}^\infty(\TT;X^*)$.
\end{prop}\textbf{Proof} The subspace $L^1(\TT;K)$ is obviously a positive cone in $L^1(\TT;X)$ and $L^\infty_{w^*}(\TT;K^*)\subs L^1(\TT;K)^*$. So let $\mu=(\mu_t)_{t\in\TT}\in L^1(\TT;K)^*$ and we will show that $\mu\in L^\infty_{w^*}(\TT;K^*)$. By definition, since $\mu\in L^1(\TT;K)^*$ we have that 
$$\lls f,\mu\rrs=\int_0^T \ls f_t,\mu_t\rs\df t\geq 0,\quad\forall\;f\in L^1(\TT;K).$$ In particular, for any $g\in K$ and $A\in\F$ we have that $g\1_A\in L^1(\TT;K)$ and thus $$0\leq\lls g\1_A,\mu\rrs=\int_A\ls g,\mu_t\rs\df t.$$ Therefore for each $g\in K$ there exists a set $E_g\subs\TT$ of full Lebesgue measure such that $\ls g,\mu_t\rs\geq 0$ for all $t\in E_g$. Since $X$ is assumed separable, there exists a countable subset $D\subs K$ such that $K\subs\bbar{D}$. Then $E:=\bigcap_{g\in D}E_g$ is of full measure in $\TT$ and $\ls g,\mu_t\rs\geq 0$ for all $g\in K$ and all $t\in E$. Thus $\mu_t\in K^*$ for all $t\in E$ and $\mu=(\mu_t)_{t\in\TT}\in L_{w^*}^\infty(\TT;K^*)$. The final claim follows by the results of Section~\ref{Appendix1}.$\hfill\Box$\\

We close this section by proving that the subspace $L_{w^*}^\infty(0,T;\MMM_{ac}(\T^d))$ of $L_{w^*}^\infty(0,T;\MMM(\T^d))$ where $\MMM_{ac}(\T^d)$ is the space of absolutely continuous measures with respect to Lebesgue measure is $w^*$-measurable. We start with some terminology. Given a polish space $X$ we say that a family $\mathcal{U}_0\subs \B_X$ of Borel subsets of $X$ is \emph{absolute continuity determining class on $X$} if for all measures $\mu\in\MMM(X)$ and $\nu\in\MMM_+(X)$ if for every $\ee>0$ there exists $\delta>0$ such that
\begin{equation}\label{AbsContCrit}
U\in\mathcal{U}_0\quad\mbox{and}\quad\nu(U)<\delta\quad\Lra\quad|\mu(U)|<\ee
\end{equation}
then $\mu$ is absolutely continuous with respect to $\nu$.
\begin{lemma}\label{CountAbsContCrit} Let $\mathcal{U}$ be a base for the topology of the polish space $X$ and let $\mathcal{U}_0$ denote the collection of all finite unions of elements of $\mathcal{U}$. Then a measure $\mu\in\MMM(\T^d)$ is absolutely continuous with respect to Lebesgue measure if and only if for any $\ee>0$ there exists $\delta>0$ such that~\eqref{AbsContCrit} holds. Consequently, since $X$ has a countable base it follows that there exists a countable absolute continuity determining class on $X$ that consists of open sets.
\end{lemma}\textbf{Proof} Obviously if $\mu\ll\nu$ then~\eqref{AbsContCrit} holds and so we prove the converse. Note first that by the same argument that shows that a Banach-valued measure $\bs{\nu}$ is absolutely continuous if and only if its variation $|\bs\nu|$ is absolutely continuous also shows that~\eqref{AbsContCrit} is equivalent to requiring that for all $\ee>0$ there exists $\delta>0$ such that $|\mu|(U)<\ee$ for all $U\in\mathcal{U}_0$ with $\nu(U)<\delta$. So let $\ee>0$, choose $\delta>0$ by~\eqref{AbsContCrit} such that $|\mu|(U)<\frac{\ee}{2}$ for all $U\in\mathcal{U}_0$ with $\nu(U)<\delta$ and let $B\subs X$ be a Borel set with $\nu(B)<\delta$. We have to show that $|\mu|(B)<\ee$. Since $\mu$ is regular there exists a compact set $K\subs B$ such that $|\mu|(B\sm K)<\frac{\ee}{2}$. Then $\nu(K)\leq\nu(B)<\delta$ and thus since $\nu$ is regular there exists an open set $A\subs X$ such that $K\subs A$ and $\nu(A)<\delta$. Since $A$ is open and $\mathcal{U}$ is a basis, the set $A$ is a union $A=\bigcup_{i\in I} U_i$ of elements $U_i\in\mathcal{U}$, $i\in I$, and covers the compact $K$. Thus there exist finitely many those elements, say $U_{i_1},\ldots,U_{i_k}$, $k\in\NN$, whose union $U_0:=\bigcup_{j=1}^kU_{i_j}\in\mathcal{U}_0$ continues to cover $K$. Then $U_0\subs A$ and thus $\nu(U_0)<\delta$ which by the choice of $\delta>0$ implies that 
$|\mu|(K)\leq|\mu|(U_0)<\frac{\ee}{2}$ and thus $|\mu|(B)=|\mu|(B\sm K)+|\mu|(K)<\ee$ as required.$\hfill\Box$

\begin{prop}\label{AbsContStarMeas} Let $X$ is a compact metric space and $\nu\in\MMM_+(X)$ be a non-negative reference measure. Then the space $\MMM_{ac}(X;\nu)$ of all measures $\mu$ on $X$ that are absolutely continuous with respect to $\nu$ is a $w^*$-measurable subspace of $\MMM(X)$.
\end{prop} By Lemma~\eqref{AbsContCrit} we can express $\MMM_{ac}(X;\nu):=\{\mu\in\MMM(X)|\mu\ll\nu\}\subs\MMM(X)\cong C(X)^*$ as 
$$\MMM_{ac}(X;\nu)=
\bigcap_{n\in\NN}\bigcup_{m\in\NN}\bigcap_{\substack{U\in\mathcal{U}_0\\\nu(U)<\fr{m}}}\Big\{\mu\in\MMM(X)\Bigm|\big|\ls\1_U,\mu\rs\big|<\fr{n}\Big\}$$
for some countable absolute continuity determining class $\mathcal{U}_0$ on $X$ consisting of open sets. From this it follows that $\MMM_{ac}(X;\nu)$ is $w^*$-measurable since for any open $U\subs X$ the linear map $\ls\1_U,\cdot\rs\colon\MMM(X)\to\RR$ is $w^*$-measurable. Indeed, since $U$ is open there exists a sequence $\{f_n\}_{n=1}^\infty\subs C(X)$ such that $0\leq f_n\leq\1_U$ converging pointwise to $\1_U$. But then the sequence of maps $\J(f_n)=\ls f_n,\cdot\rs\in C(\T^d)\cong(\MMM(X),w^*)^*\leq\MMM(X)^*$, $n\in\NN$ converges to $\ls\1_U,\cdot\rs$ pointwise in $\MMM(X)$ and thus $\ls\1_U,\cdot\rs$ is $w^*$-measurable as the pointwise limit of the $w^*$-continuous maps $\J(f_n)$.$\hfill\Box$\\

The same is also true if $X$ is a locally compact polish space if we replace $C(X)$ by $C_0(X)$. In order to extend this result on the level of path-measures we further need one more lemma.

\begin{lemma}\label{InftyLipInfty} Let $X$ be a compact metric space and let $\nu\in\MMM_+(X)$ be a fixed measure. Then 
	\begin{equation}\label{LipMeasNormInTermsOfContFunct}
\Lip_\nu(\mu):=\sup_{\nu(A)\neq 0}\frac{|\mu(A)|}{\nu(A)}=\sup_{\substack{f\in C_+(X)\\\int f\df\nu\neq 0}}\frac{\int f\df|\mu|}{\int f\df\nu}=\sup_{\substack{f\in C(X)\\\int|f|\df\nu\neq 0}}\frac{|\int f\df\mu|}{\int|f|\df\nu}.
\end{equation}	
Consequently, the map $\Lip_\nu\colon\MMM(X)\to[0,+\infty]$ is $w^*$-lower semicontinuous and thus the subspace $$\MMM_{\Lip}(\nu)=\big\{\mu\in\MMM(X)\bigm|\Lip_\nu(\mu)<+\infty\big\}=\{f\df\nu|f\in L^\infty(\nu)\}=L^\infty(\nu)$$ of $\MMM(X)$ is $w^*$-measurable. Similarly the map $$L_{w^*}^\infty(0,T;\MMM(X))\ni\mu\mapsto\|\Lip_\nu(\mu_\cdot)\|_{L^\infty(0,T)}\in [0,+\infty]$$ is $w^*$-lower semicontinuous and the subspace \begin{equation}\label{inftylipinfty}
	\big\{\mu\in L_{w^*}^\infty(0,T;\MMM(X))\bigm|\|\Lip_\nu(\mu_\cdot)\|_{L^\infty(0,T)}<+\infty\big\}
\end{equation} of $L_{w^*}^\infty(0,T;\MMM(X))$ is also $w^*$-measurable.
\end{lemma}\textbf{Proof} Let $R>0$ and let $B_{L^\infty(\nu)}(0,R)$ be the closed ball of radius $R>0$ in $L^\infty(\nu)$ centred at the origin. If $\mu\in\MMM(X)$ is a measure such that $\mu\ll\nu$ with $\frac{\df\mu}{\df\nu}\in B_{L^\infty(\nu)}(0,R)$ then
\begin{equation}\label{LipMeasInTermsOfContFunct}
\Big|\int f\df\mu\Big|\leq R\int |f|\df\nu,\quad\forall f\in C(X).
\end{equation}
In particular 
$$\sup_{f\in C(X)}\frac{\big|\int f\df\mu\big|}{\int |f|\df\nu}\leq\big\|\|\mu_\cdot\|_{L^\infty(\nu)}\big\|_{L^\infty(0,T)}=\Lip_\nu(\mu).$$ We prove next the converse inequalities required for the proof of~\eqref{LipMeasNormInTermsOfContFunct}. So let $A\subs\B_X$ be such that $\nu(A)\neq 0$. Since the measures $|\mu|$ and $\nu$ are regular there exists for each $n\in\NN$ a compact set $K_n\subs X$ and an open set $U_n\subs X$ such that $|\mu|(U_n\sm K_n)\mx\nu(U_n\sm K_n)<\fr{n}$ and $\nu(K_n)>\frac{\nu(A)}{2}>0$ and we can choose a map $f_n\in C(\T^d;[0,1])$ such that $\1_{K_n}\leq f_n\leq\1_{U_n}$. Then $\int f_n\df\nu\geq\frac{\nu(A)}{2}>0$ for all $n\in\NN$ and thus for all $n\in\NN$ large enough so that $\fr{n}<\frac{\nu(A)}{2}$
	$$\frac{|\mu|(A)}{\nu(A)}\leq\frac{|\mu|(K_n)+\fr{n}}{\nu(U_n)-\fr{n}}\leq\frac{\int f_n\df|\mu|+\fr{n}}{\int f_n\df\nu-\fr{n}}.$$
	But by construction the sequence $\{f_n\}\subs C(\T^d;[0,1])$ converges pointwise $|\mu|$-a.s.~and $\nu$-a.s.~to $\1_A$ and therefore by the dominated convergence theorem the limits $\lim_{n\ra+\infty}\int f_n\df|\mu|=|\mu|(A)$ and $\lim_{n\ra+\infty}\int f_n\df\nu=\nu(A)$ exist and thus
	$$\frac{|\mu(A)|}{\nu(A)}\leq\frac{|\mu|(A)}{\nu(A)}\leq\lim_{n\ra+\infty}\frac{\int f_n\df|\mu|+\fr{n}}{\int f_n\df\nu-\fr{n}}=\lim_{n\ra+\infty}\frac{\int f_n\df|\mu|}{\int f_n\df\nu}\leq\sup_{\substack{f\in C_+(X)\\\int f\df\nu\neq 0}}\frac{\int f\df|\mu|}{\int f\df\nu}.$$
	Since this holds for any $A\in\B_X$ with $\nu(A)\neq 0$ it follows that 
	$$\Lip_\nu(\mu)\leq\sup_{\substack{f\in C_+(X)\\\int f\df\nu\neq 0}}\frac{\int f\df|\mu|}{\int f\df\nu}\leq\sup_{\substack{f\in C(X)\\\int|f|\df\nu\neq 0}}\frac{\big|\int f\df|\mu|\big|}{\int|f|\df\nu}$$
	Consequently if we can prove that 
	$$\sup_{\substack{f\in C(X)\\\int|f|\df\nu\neq 0}}\frac{\big|\int f\df|\mu|\big|}{\int|f|\df\nu}\leq\sup_{\substack{f\in C(X)\\\int|f|\df\nu\neq 0}}\frac{\big|\int f\df\mu\big|}{\int|f|\df\nu}$$
	then~\eqref{LipMeasNormInTermsOfContFunct} follows. So let $f\in C(X)$ be such that $\int|f|\df\nu\neq 0$. Let $X=P\cup N$ be a Hahn decomposition of $X$ with respect to $\mu$ and as in the proof of Proposition~\ref{Basic}(e) we can find a sequence $\{\phi_n\}\subs C(\T^d;[-1,1])$ converging $|\mu|$-a.s.~and $\nu$-a.s.~to $\1_P-\1_N$. Obviously $|\1_P-\1_N|=1$ since $P\cap N=\emptyset$ and thus
	$$\frac{\int f\df|\mu|}{\int|f|\df\nu}=\frac{\int(\1_P-\1_N)f\df\mu}{\int|(\1_P-\1_N)f|\df\nu}=\lim_{n\ra+\infty}\frac{\int\phi_nf\df\mu}{\int|\phi_nf|\df\nu}\leq\sup_{\substack{f\in C(X)\\\int f\df\nu\neq 0}}\frac{\big|\int f\df\mu\big|}{\int|f|\df\nu}.$$
Now since for each $f\in C(X)$ the map $\MMM(X)\ni\mu\mapsto\frac{|\ls f,\mu\rs|}{\ls|f|,\nu\rs}$ is $w^*$-continuous it follows that $\Lip_\nu(\cdot)\colon\MMM(X)\to[0,+\infty]$ is $w^*$-lower semicontinuous as the supremum of $w^*$-continuous linear functionals. Consequently the set
$$B_{L^\infty(\nu)}(0,R)=\{\mu\in\MMM_{ac}(X;\nu)|\|\frac{\df\mu}{\df\nu}\|_{L^\infty(\nu)}\leq R\}=\{\Lip_\nu(\cdot)\leq R\}$$
is $w^*$-closed and thus $w^*$-measurable. Thus $L^\infty(\nu)=\bigcup_{n\in\NN}B_{L^\infty(\nu)}(0,n)$ is also $w^*$-measurable.

We prove next that the map $\|\Lip_\nu(\cdot)\|_{L^\infty(0,T)}\colon L_{w^*}^\infty(0,T;\MMM(X))\to[0,+\infty]$ is also $w^*$-semicontinuous. Let us start by noting that for each $f\in C(X)$ the function $h^f\colon L_{w^*}^\infty(0,T;\MMM(X))\to[0,+\infty]$ defined by $h^f(\mu)=\|\ls f,\mu_\cdot\rs\|_{L^\infty(0,T)}$ is lower semicontinuous. Indeed, the operator 
$I_f\colon L_{w^*}^\infty(0,T;\MMM(X))\to L^\infty(0,T)$ given by 
\begin{equation}\label{IG}
I_f(\mu)(t)=\ls f,\mu_t\rs,\quad\mbox{a.s.~}\forall t\in[0,T],
\end{equation}
is $w^*$-continuous, since if $\{\mu^\alpha\}_{\alpha\in\A}\subs L_{w^*}^\infty(0,T;\MMM(\T^d))$ is a net converging to $\mu\in L_{w^*}^\infty(0,T;\MMM(\T^d))$ in the $w^*$-topology then for all $g\in L^1(0,T)$ 
\begin{align*}\lim_\alpha\int_0^Tg(t)I_f(\mu^\alpha)(t)\df t&=\lim_\alpha\int_0^T\int_Xg(t)f(x)\df\mu^\alpha_t(x)\df t\\
&=\int_0^T\int g(t)f(x)\df\mu_t(x)\df t=\int_0^Tg(t)I_f(\mu)(t)\df t,
\end{align*}
since whenever $f\in C(X)$ and $g\in L^1(0,T)$ the function given by $F(t,x)=g(t)f(x)$ is in $L^1(0,T;C(X))$. It follows then that the function $h^f$ is lower semicontinuous as it is the composition of the $w^*$-continuous function $I_f$ and the $w^*$-lower semicontinuous function $\|\cdot\|_{L^\infty(0,T)}\colon L^\infty(0,T)\to\RR$.

Consequently the map $\Lip_\nu^\infty\colon L_{w^*}^\infty(0,T;\MMM(X))\to[0,+\infty]$ given by 
$$\Lip_\nu^\infty(\mu)=\sup_{\substack{f\in C(X)\\\int|f|\df\nu\neq 0}} \frac{\|\ls f,\mu_\cdot\rs\|_{L^\infty(0,T)}}{\ls|f|,\nu\rs}$$
is $w^*$-lower semicontinuous as the supremum of $w^*$-lower semicontinuous maps and we will show that $\Lip_\nu^\infty(\mu)=\|\Lip_\nu(\mu)\|_{L^\infty(0,T)}$ for all $\mu\in L_{w^*}^\infty(0,T;\MMM(X))$. If $R:=\|\Lip_\nu(\mu)\|_{L^\infty(0,T)}<+\infty$ then $\mu_t\ll\nu$ with $\frac{\df\mu_t}{\df\nu}\in B_{L^\infty(\nu)}(0,R)$ for all $t\in E_\mu$ where $E_\mu$ is a set of full $\nu$-measure. Therefore for all $f\in C(X)$
$$|\ls f,\mu_t\rs|=\Big|\int f\frac{\df\mu_t}{\df\nu}\df\nu\Big|\leq R\int|f|\df\nu$$
for almost all $t\in E_\mu$ which yields which yields that $\|\ls f,\mu\rs\|_{L^\infty(0,T)}\leq R\ls|f|,\nu\rs$ for all $f\in C(X)$. Therefore $\Lip^\infty_\nu(\mu)\leq R=\|\Lip_\nu(\mu_\cdot)\|_{L^\infty(0,T)}$. For the converse inequality, if $\Lip_\nu^\infty(\mu)<+\infty$ then 
$$\|\ls f,\mu_\cdot\rs\|_{L^\infty(0,T)}\leq\Lip_\nu^\infty(\mu)\ls|f|,\nu\rs,\quad\forall f\in C(X).$$
Thus for each $f\in C(X)$ there exists a Borel set $E^f\subs[0,T]$ of full $\nu$-measure such that $|\ls f,\mu_t\rs|\leq\Lip^\infty_\nu(\mu)\ls|f|,\nu\rs$ for all $t\in E^f$. Since $C(X)$ is separable there exists a countable set dense set $D\subs C(X)$ in the uniform norm and then the set $E:=\bigcap_{f\in D}E^f$ is of full $\nu$-measure and $$|\ls f,\mu_t\rs|\leq\Lip_\nu^\infty(\mu)\ls|f|,\nu\rs,\quad\forall\;(t,f)\in E\x C(X).$$
Consequently $\Lip_\nu(\mu_t)\leq\Lip^\infty_\nu(\mu)$ for all $t\in E$ which shows that also $\|\Lip_\nu(\mu_\cdot)\|_{L^\infty(0,T)}\leq\Lip^\infty_\nu(\mu)$ and completes the proof.$\hfill\Box$\\

Obviously the set defined in~\eqref{inftylipinfty} is contained in the set $L_{w^*}^\infty(0,T;L^\infty(\nu))$. Note however that according to our definitions this inclusion is in general strict.

\begin{prop}\label{ACMEAS} The subspace $L_{w^*}^\infty(0,T;\MMM_{ac}(\T^d))$ of $L_{w^*}^\infty(0,T;\MMM(\T^d))$ is $w^*$-measurable.
	\end{prop}\textbf{Proof} By Lemma~\ref{CountAbsContCrit} the space $\MMM_{\Lip}([0,T])$ of Lipschitz-continuous measures on $[0,T]$ is a $w^*$-measurable subspace of $\MMM([0,T])$. Let $\mathfrak{t}\colon[0,T]\x\T^d\to[0,T]$ denote the natural projection on the first coordinate. The push forward operator $\mathfrak{t}_\sharp\colon\MMM([0,T]\x\T^d)\to\MMM([0,T])$ is $w^*$-continuous and thus the space 
$$\MMM_{\Lip(\mathfrak{t})}([0,T]\x\T^d):=(\mathfrak{t}_\sharp)^{-1}\big(\MMM_{\Lip}([0,T])\big)\leq\MMM([0,T]\x\T^d)$$
is a $w^*$-measurable subspace of $\MMM([0,T]\x\T^d)$ by Lemma~\ref{InftyLipInfty}. The space $\MMM_{ac}([0,T]\x\T^d)$ is also $w^*$-measurable by Proposition~\ref{AbsContStarMeas} and thus the space $$\MMM_{\Lip(\mathfrak{t}),ac}([0,T]\x\T^d):=\MMM_{\Lip(\mathfrak{t})}([0,T]\x\T^d)\cap\MMM_{ac}([0,T]\x\T^d)$$
is $w^*$-measurable subspace of $\MMM([0,T]\x\T^d)$.

We consider now the inclusion operator $i\colon C([0,T]\x\T^d)\to L^1(0,T;C(\T^d))$. This is a bounded injection with $\|i(f)\|_{L^1(0,T;C(\T^d))}\leq T\|f\|_\infty$. Thus its adjoint $i^*\colon L_{w^*}^\infty(0,T;\MMM(\T^d))\to \MMM([0,T]\x\T^d)$ given by 
$$\ls F,i^*\mu\rs=\lls i(F),\mu\rrs\equiv\lls F,\mu\rrs$$
is a $w^*$-continuous operator. The space $\MMM_{ac}([0,T]\x\T^d)$ is a $w^*$-measurable subspace of $\MMM([0,T]\x\T^d)$ and thus if we can show that 
\begin{equation}\label{AbsContWStMeas}
L_{w^*}^\infty(0,T;\MMM_{ac}(\T^d))=(i^*)^{-1}\big(\MMM_{\Lip(\mathfrak{t}),ac}([0,T]\x\T^d)\big)
\end{equation}
the claim follows. 

So let first $\mu\in L_{w^*}^\infty(0,T;\MMM_{ac}(\T^d))$ and we will show that $i^*(\mu)\in\MMM_{\Lip(\mathfrak{t}),ac}([0,T]\x\T^d)$. The $\mathfrak{t}$-marginal of the measure $\bs\mu:=i^*\mu$ satisfies is characterized by 
$$\int f(t)\df\mathfrak{t}_\sharp\bs\mu(t)=\int_0^T\int_{\T^d}f(t)\df\mu_t\df t=\int_0^Tf(t)\mu_t(\T^d)\df t\quad\forall\;f\in B([0,T])$$
and thus $\mathfrak{t}_\sharp\bs\mu\ll\LL_{[0,T]}$ with density $\frac{\df\mathfrak{t}_\sharp\bs\mu}{\df\LL_{(0,T)}}(t)=\mu_t(\T^d)$ for almost all $t\in[0,T]$. Thus $\mathfrak{t}_\sharp\bs\mu\in\MMM_\Lip([0,T])$ with $\|\mathfrak{t}_\sharp\bs\mu\|_{\Lip}\leq\|\mu\|_{TV;\infty}<+\infty$, i.e.~$\bs\mu\in\MMM_{\Lip(\mathfrak{t})}([0,T]\x\T^d)$. To see that also $\bs\mu\in\MMM_{ac}([0,T]\x\T^d)$ we note that since $\mu_t\in\MMM_{ac}(\T^d)$ for almost all $t\in[0,T]$ the measure $\bs\mu=i^*\mu$ is characterized by 
$$\ls F,\bs\mu\rs=\int_0^T\int_{\T^d}F(t,u)\frac{\df\mu_t}{\df\LL_{\T^d}}(u)\df u\df t.$$
It follows that $\bs\mu\ll\LL_{[0,T]\x\T^d}$ with density $$\frac{\df\bs\mu}{\df\LL_{[0,T]\x\T^d}}(t,u)=\frac{\df\mu_t}{\df\LL_{\T^d}}(u),\quad\mbox{a.s.~for all }(t,u)\in[0,T]\x\T^d,$$
which proves that $\bs\mu\in\MMM_{\Lip(\mathfrak{t}),ac}([0,T]\x\T^d)$ and thus the inclusion ``$\subs$" in~\eqref{AbsContWStMeas}.

For the converse inclusion let $\mu\in L_{w^*}^\infty(0,T;\MMM(\T^d))$ be such that $i^*(\mu)\in\MMM_{\Lip(\mathfrak{t}),ac}([0,T]\x\T^d)$ and we will show that $\mu_t\ll\LL_{\T^d}$ for almost all $t\in[0,T]$. Since $i^*(\mu)\in\MMM_{ac}([0,T]\x\T^d)$ there exists unique $\psi\in L^1([0,T]\x\T^d)$ such that $\df i^*(\mu)=\psi\df\LL_{[0,T]\x\T^d}$. Then for any $F\in C([0,T]\x\T^d)$ 
$$\lls i(F),\mu\rrs=\ls F,\bs\mu\rs=\int_0^T\int_{\T^d}F(t,u)\psi(t,u)\df u\df t$$
and thus by applying this for maps $F$ of the form $F(t,u)=f(t)g(u)$ with $f\in C([0,T])$, $g\in C(\T^d)$ we obtain that for all each fixed $g\in C(\T^d)$ 
$$\int_0^Tf(t)\int g\df\mu_t\df t=\int_0^Tf(t)\int_{\T^d}g(u)\psi(t,u)\df u\df t,\quad\forall\;f\in C([0,T]).$$
This implies that for each $g\in C(\T^d)$ 
$$\int_{\T^d}g\df\mu_t=\int_{\T^d}g(u)\psi(t,u)\df u\quad\LL_{[0,T]}\mbox{-a.s.~for all }t\in[0,T].$$
Since $C(\T^d)$ is separable this implies that $\mu_t\ll\LL_{\T^d}$ with $\df\mu_t=\psi(t,\cdot)\df\LL_{\T^d}$ and completes the proof.$\hfill\Box$

\subsubsection{Induced operators}

Let $X,Y$ be Banach spaces. Any bounded operator $S\colon X\to Y$ between Banach spaces induces a bounded operator $\bar{S}\colon L^p(\TT;X)\to L^p(\TT;Y)$, $1\leq p\leq\infty$, on the corresponding Bochner-$L^p$ spaces via $\bar{S}(f)=\big(S(f_t)\big)_{t\in\TT}$ for $f=(f_t)_{t\in\TT}\in L^p(\TT;X)$. Obviously $\|\bar{S}(f)\|_{L^p(\TT;Y)}\leq\|S\|\|f\|_{L^p(\TT;X)}$ and by checking against constant paths we see that in fact 
$$\|\bar{S}\|=\sup_{f\in L^p(\TT;X)\sm\{0\}}\frac{\|S(f)\|_{L^p(\TT;Y)}}{\|f\|_{L^p(\TT;X)}}\geq\sup_{f\in X\sm\{0\}}\frac{S(f)}{\|f\|_X}=\|S\|$$
so that the induced operator $\bar{S}$ retains the same norm. If $S$ is injective, a contraction or norm preserving then so is $\bar{S}$. If $S$ is strongly surjective in the sense that it has a bounded right inverse $T\colon Y\to X$ then so does $\bar{S}$, namely $\bar{S}\circ\bar{T}=\mathbbm{id}_{L_p(\TT;Y)}$.

\begin{prop}\label{InducedOpPropL1StrongConv} If a sequence of operators $S_n\colon X\to Y$ converges strongly to $S\colon X\to Y$, i.e.~if $\lim_{n\ra+\infty}\|S_nf-Sf\|_Y=0$ for all $f\in X$, then the sequence of induced operators $\bar{S}_n$ on the corresponding $L^1$-spaces converges strongly to $\bar{S}$. 
\end{prop}\textbf{Proof} Let $F\in L^1(0,T;X)$. Since $S_n$ converges to $S$ strongly it follows that $\lim_{n\ra+\infty}\|S_nF_t-SF_t\|_Y=0$ for almost all $0\leq t\leq T$. Furthermore $\sup_{n\in\NN}\|(S_n-S)f\|_Y<+\infty$ for all $f\in X$ and therefore by the uniform bounded principle $M:=\sup_{n\in\NN}\|S_n-S\|<+\infty$. Thus 
$\|S_nF_t-SF_t\|_Y\leq M\|F_t\|_X$ for almost all $0\leq t\leq T$ and the maps $\|(S_n-S)F_\cdot\|_Y$, $n\in\NN$, are dominated by the integrable map $M\|F_\cdot\|_X$. Therefore the claim follows by the dominated convergence theorem.$\hfill\Box$\\

The aim of this section is to provide conditions that ensure an operator $T\colon X^*\to Y^*$ induces a $w^*$-measurable operators on the respective $L_{w^*}^\infty$-spaces. First let us note that it is obvious that any bounded and $w^*$-measurable operator $T\colon X^*\to Y^*$ induces an operator $\bar{T}\colon L_{w^*}^\infty(\TT;X^*)\to L_{w^*}^\infty(\TT;Y^*)$. As we will see, this $\bar{\cdot}$ operator that maps an operator $T\in\B_{w^*}(X^*,Y^*)$ to the induced operator on the corresponding $L_{w^*}^\infty$-spaces has nice categorical properties and respects the notion of $w^*$-Baire measurability of operators. Since $w^*$-Baire measurability is known to be stronger than $w^*$-measurability we obtain a condition that ensures the $w^*$-measurability of induced operators. This stronger assumption on linear operators in order to induce $w^*$-measurable operators on the $L_{w^*}^\infty$-spaces will not pose a problem in the main text since all operators we will encounter will be $w^*$-measurable.

\begin{prop}\label{InducedOpProp} (a) Let $T\colon X^*\to Y^*$ be bounded and $w^*$-measurable. Then the formula 
	\begin{equation}\label{InducedOp}
	\bar{T}(\mu)(t)=T(\mu_t)\quad\mbox{for }\mbox{ almost all }t\in\TT,\end{equation}
	defines a bounded linear operator $\bar{T}\colon L_{w^*}^\infty(\TT;X^*)\to L_{w^*}^\infty(\TT;Y^*)$ with norm $\|\bar{T}\|=\|T\|$. If $T$ is norm-preserving then so is $\bar{T}$.\\
	\noindent(b) If $S\colon Y^*\to Z^*$ is bounded and $w^*$-measurable then $\bbar{S\circ T}=\bbar{S}\circ\bbar{T}$ and if $T_i\colon X^*\to Y^*$ are bounded and $w^*$-measurable then $\bbar{T_1+T_2}=\bbar{T_1}+\bbar{T_2}$.\\
	\noindent(c) Finally, if $T\in \A_{w^*}(X^*,Y^*)$ is $w^*$-Baire measurable then so is $\bar{T}$ and thus $\bar{T}$ is also $w^*$-measurable.
\end{prop}\textbf{Proof} (a) Since $T\colon X^*\to Y^*$ is $w^*$-measurable, for any $\mu\in\LL_{w^*}^\infty(\TT;X^*)$ the map $T\circ\mu\colon[0,T]\to Y^*$ is $w^*$-measurable and since $T$ is a bounded operator, 
$$\|T\circ\mu(t)\|_{Y^*}=\|T(\mu_t)\|_{Y^*}\leq\|T\|\|\mu_t\|_{X^*}\quad \mbox{a.s.~for all }t\in[0,T].$$ Therefore $\bar{T}(\mu)$ is an element of $L_{w^*}^\infty(\TT;Y^*)$ for all $\mu\in L_{w^*}^\infty(\TT;X^*)$ and the induced operator $\bar{T}$ is bounded with $\|\bar{T}\|\leq\|T\|$. By checking against the constant maps in $L_{w^*}(\TT;X^*)$ it follows that $\|\bar{T}\|=\|T\|$. Also, if $T$ is norm preserving then $\|T\mu_t\|_{Y^*}=\|\mu_t\|_{Y^*}$ almost everywhere for any $\mu\in L_{w^*}^\infty(\TT;X^*)$ and thus $\|\bar{T}(\mu)\|_{L_{w^*}(\TT;Y^*)}=\|\mu\|_{L_{w^*}^\infty(\TT;X^*)}$. Statement (b) is trivial.\\
\noindent(c) The fact that $\bar{T}$ is $w^*$-Baire whenever $T$ is, follows by a simple transfinite induction argument based on the following Proposition~\ref{LinftInducedwStarOpCont} according to which the application $B_{w^*}(X^*,Y^*)\ni T\mapsto \bar{T}\in B_{w^*}(L_{w^*}^\infty(\TT;X^*),L_{w^*}^\infty(\TT;Y^*))$ is sequentially continuous with respect to pointwise $w^*$-convergence of operators.$\hfill\Box$

\begin{prop}\label{LinftInducedwStarOpCont} Let $\{T_n\}_{n=1}^\infty\subs B_{w^*}(X^*,Y^*)$ be a sequence of operators $w^*$-converging to $T\in B_{w^*}(X^*,Y^*)$. Then $\bar{T}_n$ pointwise $w^*$-converges to $\bar{T}$.
\end{prop}\textbf{Proof} We have to show that $\lim_{n\ra+\infty}\lls g,\bar{T}_n\mu\rrs=\lls g,\bar{T}\mu\rrs$ for all $g\in L^1(\TT;Y)$, $\mu\in L^\infty_{w^*}(\TT;X^*)$. Since $T_n\lra T$ $w^*$-pointwise $\lim_{n\ra+\infty}\ls g_t,T_n\mu_t\rs=\ls g_t,T\mu_t\rs$ for all $t\in\TT$ and $\sup_{n\in\NN}\|T_n\|<+\infty$ and thus an application of the dominated convergence theorem concludes the proof.$\hfill\Box$\\

In the case that $X$ is non-separable, in order for a bounded and $w^*$-measurable operator $T\colon X^*\to Y^*$ to induce an operator $\bar{T}\colon\bar{L}_{w^*}^\infty(\TT;X^*)\to \bar{L}_{w^*}^\infty(\TT;Y^*)$ one has to assume in addition that $T$ respects the relation of $w^*$-a.s.~equality, i.e. that 
\begin{equation}\label{EqRelRespect}[\ls f,\mu^1_\cdot-\mu^2_\cdot\rs=0,\;m\mbox{-a.s.}]\;\forall f\in X\quad\Lra\quad[\ls g,T(\mu^1_\cdot-\mu^2_\cdot)\rs=0,\;m\mbox{-a.s.}]\;\forall g\in Y.\end{equation}
Indeed assumption~\eqref{EqRelRespect} is equivalent to $\bar{T}(\ker S_X)\leq\ker S_Y$ where $S_Z\colon L_{w^*}^\infty(\TT;Z^*)\to L^1(\TT;Z)^*$, $Z=X,Y$, is the surjective contraction that is induced by the bilinear pairing $\lls\cdot,\cdot\rrs$ between $L^1(\TT;Z)$ and $L_{w^*}^\infty(\TT;Z^*)$. This implies that $\bar{T}$ induces an operator $\bar{T}\colon\bar{L}_{w^*}^\infty(\TT;X^*)\to \bar{L}_{w^*}^\infty(\TT;Y^*)$ (denoted by the same symbol $\bar{T}$) by the formula 
\begin{equation}\label{EqRelRespect2}\bbar{T}(\mu+\ker S_X)=\bar{T}(\mu)+\ker S_Y.\end{equation}
As we will see all $w^*$-Baire operators $T\colon X^*\to Y^*$ respect the relation of $w^*$-m-a.s.~equality. This follows since $w^*$-continuous operators respect this relation and the set of operators that respect the relation of $w^*$-m-a.s.~equality is sequentially closed with respect pointwise $w^*$-convergence of operators.

\begin{prop}\label{AdjointInduced} Any adjoint operator $T=T_0^*\colon X^*\to Y^*$, where $T_0\colon Y\to X$ is bounded satisfies~\eqref{EqRelRespect} and thus induces an operator $\bar{T}\colon\bar{L}_{w^*}^\infty(\TT;X^*)\to\bar{L}_{w^*}^\infty(\TT;Y^*)$ by the formula~\eqref{EqRelRespect2} and \begin{equation}\label{InducedAdjoint}\bbar{T}=\bbar{T_0^*}=(\bbar{T_0})^*\end{equation}
where $\bbar{T_0}\colon L^1(\TT;Y)\to L^1(\TT;X)$ is the induced operator on the $L^1$-spaces. In particular $\bar{T}$ is $w^*$-continuous as the adjoint of the bounded operator $\bar{T}_0$.\end{prop}\textbf{Proof} Since $T$ is an adjoint operator it is bounded and $w^*$-measurable and thus by Proposition~\ref{InducedOpProp} it induces an operator $\wt{T}\colon\LL_{w^*}^\infty(\TT;X^*)\to\LL_{w^*}^\infty(\TT;Y^*)$ by the formula~\eqref{InducedOp}. But since $T$ is the adjoint of $T_0$, if $\mu^1,\mu^2\in\LL_{w^*}^\infty(\TT;X^*)$ are such that $\mu^1\backsim\mu^2$ then for any $g\in Y$
$$\ls g,T(\mu^1_t)-T(\mu^2_t)\rs=\ls T_0(g),\mu_t^1-\mu_t^2\rs=0\quad m\mbox{-a.s.~for all }t\in\TT.$$
Thus $T$ satisfies~\eqref{EqRelRespect} and induces an operator $\bbar{T}\colon\bar{L}_{w^*}^\infty(\TT;X^*)\to\bar{L}_{w^*}^\infty(\TT;Y^*)$ by~\eqref{EqRelRespect2}.

It remains to verify that~\eqref{InducedAdjoint} holds. So let $g\in L^1(\TT;Y)$ and $\mu\in \bar{L}_{w^*}^\infty(\TT;X^*)$. Then 
\begin{align*}
\lls g,\bbar{T}(\mu)\rrs&=\lls g,\bbar{T_0^*}(\mu)\rrs=\int\ls g_t,\wt{T_0^*}(\mu)(t)\rs\df m(t)=\int\ls g_t,T_0^*(\mu_t)\rs\df m(t)\\
&=\int\ls T_0(g_t),\mu_t\rs\df m(t)=\lls \bbar{T_0}(g),\mu\rrs=\lls g,(\bbar{T_0})^*(\mu)\rrs
\end{align*}
and since this holds for all $g\in L^1(\TT;Y)$ and all $\mu\in \bar{L}_{w^*}^\infty(\TT;X^*)$ the equality~\eqref{InducedAdjoint} follows.$\hfill\Box$

\begin{prop} The collection $\A_{w^*\mbox{-}m\mbox{-a.s.}}$ of all operators $T\in B_{w^*}(X^*,Y^*)$ that respect the relation of $w^*$-m-a.s.~equality, i.e.~that satisfy~\eqref{EqRelRespect} is sequentially closed with respect to pointwise $w^*$-convergence of operators and contains all $w^*$-continuous operators from $X^*$ to $Y^*$. Consequently 
	$$\A_{w^*}(X^*,Y^*)\subs\A_{w^*\mbox{-}m\mbox{-}a.s.}.$$ 
\end{prop}\textbf{Proof} Since an operator $T\colon X^*\to Y^*$ is $w^*$-continuous if and only if $T$ is the adjoint of a bounded operator $T_0\colon Y\to X$ is follows by Proposition~\ref{AdjointInduced} that $BC_{w^*}(X^*,Y^*)\subs\A_{w^*\mbox{-}m\mbox{-}a.s.}$. It is also easy to see that the collection $\A$ is also sequentially closed. Indeed, let $\{T_n\}\subs\A$ be a sequence $w^*$-converging pointwise to $T$, let $\mu^1,\mu^2\in L_{w^*}^\infty(\TT;X^*)$ be $w^*$-m-a.s.~equal and let $g\in Y$. Since $\{T_n\}\subs\A_{w^*\mbox{-}m\mbox{-}a.s.}$ and $\mu^1,\mu^2$ are $w^*$-m-a.s.~equal it follows that $\ls g,T_n\mu_t^1\rs=\ls g,T_n\mu^2_t\rs$ for all $t\in F_n^g$ where $F_n^g\in\F$ is a set of full measure in $\TT$. Then $F^g:=\bigcap_{n\in\NN}F_n^g$ is of full measure and $\ls g,T_n\mu^1_t\rs=\ls g,T_n\mu^2_t\rs$ for all $n\in\NN$ and all $t\in F^g$. Since $\{T_n\}$ pointwise $w^*$-converges to $T$, taking the limit as $n\ra+\infty$ we obtain that $\ls g,T\mu^1_t\rs=\ls g,T\mu^2_t\rs$ for all $t\in F^g$ and therefore $T\mu^1=T\mu^2$ $w^*$-m-a.s.~and thus $T\in\A_{w^*\mbox{-}m\mbox{-}a.s.}$.$\hfill\Box$


\begin{thebibliography}{10}

		\bibitem{AGS2000a}
			Luigi Ambrosio, Nicola Gigli and Giuseppe Savare.
			\newblock {\em Gradient Flows in Metric Spaces and in the Space of Probability Measures}, {\em Lectures in Mathematics}.
			\newblock Birkhauser, 2000.
			
		\bibitem{Armendariz2009a}
			In\'{e}s Armend\'{a}riz, Stefan Grosskinsky and Michail Loulakis.
			\newblock Zero-range condensation at criticality.
			\newblock{\em Stochastic Processes and their Applications}, 123, 2009.
			
		\bibitem{Armendariz2008a}
			In\'{e}s Armend\'{a}riz and Michail Loulakis.
			\newblock Thermodynamic Limit for the Invariant Measures in Supercritical Zero Range Processes.
			\newblock{\em Probability Theory and Related Fields}, 2008.
			
		\bibitem{Bialas1997a}
			P. Bialas, Z. Burda and D. Johnston.
			\newblock Condensation in the backgammon model.
			\newblock{\em Nuclear Physics B}, 493:505--516, 1997
		
		\bibitem{Billingsley1971a}
		Patrick Billingsley.
		\newblock{\em Weak convergence of measures: Applications in probability}, 1971.
		 
				
		\bibitem{Cam1957a}
		Le Cam.
		\newblock Convergence in distribution of Stochastic processes.
		\newblock {\em University of California Publications in Statistics}, 2(11):207--236, 1957.
		
		\bibitem{Chleboun2014a}
		    Stefan Grosskinsky and Paul Chleboun.
			\newblock Condensation in stochastic particle systems with stationary product measures.
			\newblock{\em Journal of Statistical Physics}, 154:432--465, 2014.
		
		\bibitem{Dellacherie1978a}
		Claude Dellacherie and Paul-Andr\'{e} Meyer.
		{\em Probabilities and Potential}, volume 29 of {\em Mathematics Studies}.
		\newblock North-Holland, 1978.
			
		\bibitem{Diestel1977a}
			Joseph Diestel and John Jerry Uhl.
			\newblock{\em Vector measures}, volume 15 of {\em Mathematical Surveys}.
			\newblock American Mathematical Society, Rhode Island, 1977.
			
		\bibitem{DNS2009a}
			Jean Dolbeault, Bruno Nazaret and Giuseppe Savare.
			\newblock A new class of transport distances between measures.
			\newblock{\em Calculus of Variations}, 34(2):193--231, 2009.
			
			\bibitem{DiBenedetto2016a}
			Emmanuele DiBenedetto.
			\newblock{Real Analysis Second Edition}, Birkh{\"a}user Advanced Texts   Basler Lehrb{\"u}cher
			\newblock  Birkh\"{a}user, Springer New York, 2016.
		
		\bibitem{Dinculeanu1967a}
		N.~Dinculeanu.
		\newblock{\em Vector Measures}.
		\newblock Pergamon Press, New York, 1967.
		
		\bibitem{Drouffe1998a}
			J.-M.~Drouffe and C.~Godr\'{e}che and F.~Camia,
			\newblock A simple stochastic model for the dynamics of condensation.
			\newblock{\em Journal of Physics A-Mathematical and General}, 31(1), 1998.
		
		\bibitem{Dunford1958a}
		N.~Dunford and J.~T.~Schwartz.
		\newblock{\em Linear Operators, Part I: General Theory}.
		\newblock Interscience Publishers, 1958.	
		
		\bibitem{Ethier1986a}
		N.~Ethier and Thomas G.~Kurtz.
\newblock{\em Markov Processes: Characterization and Convergence}.
\newblock Jon Wiley and Sons, 1986.
			
		\bibitem{Evans2000a}
			Martin R. Evans.
			\newblock Phase Transitions in One-dimensional Non-equilibrium Systems.
			\newblock{\em Brazilian Journal of Physics}, 30(1):196--240, 2000.
			
		\bibitem{Fabian2010a}
			Mari\'{a}n Fabian, Petr Habala, Petr H\'{a}jek, Vicente Mondesinos and V\'{a}clav Zizler.
			\newblock{\em Banach Space Theory, The Basis for Linear and Non-linear Analysis}, CMS Books in Mathematics.
			\newblock Springer, 2010
			
		\bibitem{Faggionato2008a}
		Alessandra Faggionato.
		\newblock Hydrodynamic Limit of Zero Range Processes Among Random Conductances on the Supercritical Percolation Cluster.
		\newblock{\em Electronic Journal Of Probability}, 15, 2008.
		
		\bibitem{Fornaro2012a}
		    S.~Fornaro,  S.~Lisini, G.~Toscani and G.~Savar{\'{e}}.
		    \newblock Measure valued solutions of sub-linear diffusion equations with a drift term.
		    \newblock{\em Discrete and Continuous Dynamical Systems}, 32(5):1675--1707, 2012.

			
\bibitem{Guo1988a}
M.~Z.~Guo, G.C.~Papanicolaou and S.R.S.~Varadhan.
\newblock Nonlinear diffusion limit for a system with nearest neighbor interactions.
\newblock {\em Commun.~Math.~Phys.}, 118, 31-50, 1988.
	
		\bibitem{Grosskinsky2003a}
	    	Stefan Grosskinsky, Gunter~M. Sch{\"u}tz, and Herbert Spohn.
		    \newblock Condensation in the zero range process: stationary and dynamical properties.
		    \newblock {\em J. Statist. Phys.}, 113(3-4):389--410, 2003.
		
		\bibitem{Harris2007a}
			G.~M.~Sch{\"{u}}tz and R.~J.~Harris.
            \newblock Hydrodynamics of the zero-range Process in the condensation regime.
            \newblock{\em Journal of Statistical Physics}, 127(2):419--430, 2007.
            
        \bibitem{Husain1966a}
        	Taqdir Husain.
        	\newblock{\em Introduction to Topological Groups}
        	\newblock W.B.~Saunders Company, 1966. 

		\bibitem{Kipnis1999a}
		Claude Kipnis and Claudio Landim.
		\newblock {\em Scaling limits of interacting particle systems}, volume 320 of
		{\em Grundlehren der Mathematischen Wissenschaften [Fundamental Principles of
			Mathematical Sciences]}.
		\newblock Springer-Verlag, Berlin, 1999.
		
		\bibitem{PilarCembranos1997a}
		Pilar Cembranos and Jose Mendoza.
		\newblock{\em Banach Spaces of Vector Valued Functions}, volume 1676 of {\em Lecture Notes in Mathematics}.
		\newblock Springer-Verlag, 1997.
		
		\bibitem{RevuzYor1999}
		Daniel Revuz and Marc Yor.
		\newblock{\em Continuous Martingales and Brownian Motion}, volume 293 of
		{\em Grundlehren der mathematischen Wissenschaften [Fundamental Principles of
			Mathematical Sciences]}.
		\newblock Springer-Verlag, Berlin, 1999.
		
		\bibitem{Smolyanov1976a}
		O. G. Smolyanov and S. V. Fomin.
		\newblock{Measures on Linear Topological Spaces.}
		\newblock{\em Russian Mathematical Surveys}, 31:4,1--53, 1976.
		
		\bibitem{Stamatakis2014}
		M.G.~Stamatakis.
		\newblock{Hydrodynamic Limit of Mean Zero Condensing Zero Range Processes with Sub-Critical Initial Profiles}
		\newblock{\em Journal of Statistical Physics}, 158(1), 2014.
		
		\bibitem{Stegall1975a}
		Charles Stegall.
		\newblock{The Radon-Nikodym Property in Conjugate Banach Spaces}.
		\newblock{\em Transactions of the American Mathematical Society}, 216:213--223, 1975.
		
		\bibitem{Topsoe1970a}
		Flemming Tops{\o}e. 
		\newblock{\em Topology and Measure}, volume 133, Springer, 1970.
		
		\bibitem{TulceaAC1969a}
		A.~and~C.~ Ionescu Tulcea.
		\newblock{\em Topics in the Theory of Lifting}, volume 48 of {\em Ergebnisse der Mathematik und ihrer Grenzgebiete}.
		\newblock Springer-Verlag, Berlin Heidelberg, 1969.
		
		\bibitem{Uhl1972a}
		J. J. Uhl, Jr.
		\newblock{A  note on the Radon-Nikodym property for Banach spaces.}
		\newblock{\em Rev. Roumaine	Math. Pures	Appl.},	17:113--115, 1972.
	\end{thebibliography}
\end{document}